\documentclass[a4paper, twoside, 12pt]{report}

\usepackage[utf8]{inputenc} 
\usepackage[T1]{fontenc}
\usepackage{lmodern}
\usepackage{fancyhdr}
\usepackage{fullpage}
\usepackage{geometry}
\usepackage{young}
\usepackage{tikz}
\usepackage{tikz-cd}

\usepackage{amsthm}
\usepackage{amsmath}
\usepackage{amssymb}
\usepackage{mathrsfs}
\usepackage{enumerate}

\usepackage{pstricks}
\usepackage{pstricks-add}
\usepackage{pst-plot}
\usepackage{graphicx}
\usepackage{url}

\def\op {\operatorname}
\DeclareMathOperator{\diag}{diag}
\DeclareMathOperator{\tr}{Tr}
\DeclareMathOperator{\lcm}{lcm}
\DeclareMathOperator{\Gal}{Gal}

\title{Artin groups and Iwahori-Hecke algebras over finite fields}

\newcommand{\Z}{\mathbb{Z}}
\newcommand{\N}{\mathbb{N}}

\newcommand{\C}{\mathbb{C}}
\newcommand{\F}{\mathbb{F}}
\newcommand{\Q}{\mathbb{Q}}
\newcommand{\T}{\mathbb{T}}
\newtheorem{Def}{Definition}[section]
\newtheorem{lemme}{Lemma}[section]
\newtheorem{prop}{Proposition}[section]
\newtheorem{theo}{Theorem}[section]
\newtheorem{Def2}{Definition}[chapter]
\newtheorem{lemme2}{Lemma}[chapter]
\newtheorem{prop2}{Proposition}[chapter]
\newtheorem{theo2}{Theorem}[chapter]
\newtheorem{conj2}{Conjecture}[chapter]
\newtheorem{theor}{Théorème}
\newtheorem{lemmefr}{Lemme}
\newtheorem{Deffr}{Définition}
\newtheorem{conjfr}{Conjecture}
\newtheorem{cor}[Def]{Corollary}

\newtheorem{conjecture}{Conjecture}[section]

\title{Artin groups and Iwahori-Hecke algebras over finite fields}
\author{Alexandre Esterle}
 
\begin{document}

\maketitle

\newpage

\tableofcontents

 \newpage
 
\section*{Remerciements}

Je tiens tout d'abord à remercier Olivier qui a accepté d'encadrer mon mémoire quand je lui ai demandé durant mon master 2 à Paris. Je n'étais pas encore sûr de vouloir faire une thèse à ce moment-là, j'avais décidé de ne le faire que si je trouvais un sujet qui me motivait à fond. Je me souviens encore des réunions de plusieurs heures que nous avions. J'ai ensuite commencé ma thèse avec Ivan et lui. Il m'a aidé durant toute la thèse, que ce soit pour des problèmes mathématiques, pour rencontrer des gens à des conférences ou en relisant une bonne partie de ma thèse avec minutie, ce qui m'a aidé à la rendre bien plus lisible. Je le remercie de toute son aide et son soutien durant cette thèse et en dehors. Je m'étais engagé au début du contrat de la thèse à ramener du vin de Bordeaux, je vais enfin tenir ma promesse au pot de thèse.

\smallskip

Quand j'ai demandé à Olivier s'il pouvait me proposer un sujet de thèse, il m'a dit qu'il avait un projet avec un certain Ivan Marin. J'ai pu le rencontrer durant son cours de M2 spécialisé à Amiens. Son cours était fantastique et j'ai eu droit à un entretien privé au Lutétia où il avait l'air plutôt sympa. J'en ai alors profité pour me sortir de la ville de Paris et déménager à Amiens durant ma thèse. Nous avions une équipe de choc incroyable avec Georges et Eirini, à nous quatre, nous formions le groupe de travail d'Ivan Marin et ses thésards. Au fil des années, il a tenté de contenir mes envies de tout démontrer par des calculs interminables et à trouver des preuves plus lisibles. Il m'a accueilli chez lui à plusieurs reprises, je remercie Christine pour son accueil, j'ai été honoré de devenir responsable de production de grenouilles en papier de Claire. Ivan a réussi à coencadrer ma thèse à merveille et à devenir un chef, un patron, un guide et sourtout un ami.

\smallskip

Je remercie aussi tous les membres du LAMFA. Je veux remercier en particulier Serge, le spécialiste en de nombreux domaines des mathématiques et surtout en jeu de mots du labo, il est souvent le Bouc Emiserge du laboratoire quand il s'agit de trouver quelqu'un à qui poser des questions mathématiques ultra-spécifiques, il m'a par exemple aidé en me précisant où trouver la référence exacte pour le Lemme de Goursat, ce qui est un exploit assez remarquable.

\smallskip

Je veux remercier tous les chercheurs que j'ai pu rencontrer aux cours des diverses conférences. En particulier, Jérémie qui m'a expliqué les cellules de Kazdan-Lusztig. I want to thank Kay for all the kind help he gave me during the thesis. He helped at multiple times in conferences answering the many questions which I asked him daily. He helped make some unreadable proofs much more readable and also gave me some ideas which helped me begin to write some new unreadable proofs. I also want to thank Kay and Götz for accepting to review my thesis and Cédric, David, Donna and Jérémie for accepting to be on the jury.

\smallskip

Je remercie tous les membres de l'équipe de handball de Saint-Maurice qui m'ont aidé à me défouler chaque semaine, ce qui était vraiment nécessaire pour survivre à ces trois années. Pouvoir disputer la finale de la coupe de Somme trois semaines avant ma soutenance était un moment vraiment génial.

\smallskip

Je veux maintenant remercier tous les doctorants du LAMFA de ces trois dernières années. Merci à Vianney qui m'a aidé à prendre mes marques dès le début de ma thèse et à stopper ma première crise de panique qui est survenue après avoir appuyé sur la touche inser de mon clavier. Merci aussi à Maxime qui a été très amical dès mon arrivée. Merci à Pierrot, Malal et Emna pour leur compagnie dans le bureau BC010 qui était encore assez vide au début de ma thèse, merci à Pierrot et Emna de m'avoir présenté Youcef qui est trop mignon.

 Merci à Sylvain et Anne-Sophie qui m'ont intégré à leur merveilleux duo en première année pour que je cesse de manger tout seul. Vous m'avez attendu cinq minutes avant de manger pour que je finisse le type B plus de fois que je ne peux le compter. Merci Anne-Sophie pour ta bonne humeur et la campagne terrible que j'ai eu l'honneur de mener face à toi pour l'organisation du séminaire doctorants que nous avons pu coorganiser à merveille pendant notre deuxième année. Merci Sylvain pour ta capacité à parler sans fin et ton aide quand j'envisageais encore de faire un peu de topologie dans ma thèse.
 
  Merci Ruxi pour les cours de chinois que tu m'as prodigués, je me dois te demander : "Ni Xi Huan Tu Do Ma?". Merci Valérie pour toutes tes relectures, je ne noterai malheureusement jamais les références comme tu le souhaites, c'était un plaisir de travailler aux côtés d'une personne aussi sage qui a déjà vu passer une bonne partie de sa vie devant ses yeux. Merci Clément pour tes cours de ponctuation et ton unending combat contre les anglicismes. Merci Jérémy d'avoir coorganisé à $80\%$ le séminaire doctorants la troisième année. Merci à toi et Sylvain d'avoir été à plusieurs reprises mes consultants en algorithmes avec des avis parfois quelque peu divergents mais qui m'ont bien aidés.

 Merci Clara pour tes centaines de selfies qui me permettront de garder plein de souvenirs de la thèse. Merci Gauthier, le bout-en-train du laboratoire, tes blagues sont parfois très drôles, parfois moins mais le moment de gêne ou d'hilarité qui suit est toujours mémorable. Thank you Arrianne for your company in office BC010, I hope you will keep working on your French word lists while I'm gone, here's a new one "paillasson =mat". Thank you Afaf for smiling all the time which helped creating a great atmosphere in the office everytime you are there. Sebastian, gracias por tu sumario y lo siento por la calidad de mi español. Merci à Marouan pour les petits airs d'harmonica en fin de journée. Thank you Hatice for your support during this last year, I hope you manage to realize that life is a competition and that you manage to win at life.

 Merci enfin Guillaume pour ta bonne humeur et tes histoires plus folles les unes que les autres. Je te remercie aussi pour ton niveau d'exigence très élevé en humour, j'essaierai d'atteindre les 8 ou 9 sur 10 bien plus souvent dans le futur. Merci aussi pour tes "à demain en super forme" qui ont rythmés nos journées et pour m'avoir fait connaître de grands artistes tels que Trent Reznor, François Asselineau ou Franck Syx, ça me fera bizarre de ne pas travailler avec toi au bureau l'année prochaine.
 
  Bonne chance à tous ceux qui doivent encore finir leur thèse, merci d'avoir supporté mes nombreux mails pendant les deux dernières années, je vous laisse avec deux dernières questions mathématiques.

Question mathématique du jour : Combien ai-je de frères et soeurs et combien d'entre elles (accord de proximité) pouvez-vous reconnaître et nommer dans cet amphithéâtre?

Mathematical question of the day : Which mathematical word did I use the most during the presentation and how many times did I say it?

\smallskip

Je vais remercier mes demi-frères et demi-soeurs de thèse de la moins jeune au plus jeune. 

Eirini, je te remercie pour ta gentillesse sans égale. Avoir pu te voir surmonter ta dernière année de thèse pendant ma première année m'a montré qu'il était possible de survivre à cette épreuve. Je me souviens des week-ends au laboratoire avec Georges où on travaillait dur pour tenir les rythmes infernaux imposés par le patron. Merci de m'avoir invité à Stuttgart et de m'avoir fait visité la ville alors que tu étais sur le point de démontrer la conjecture de trace pour un nouveau groupe.

Reda, tu as beau être bien plus agé, on a pu faire le master ensemble et on a pu explorer le monde de la recherche à de multiples conférences. Merci de m'avoir hébergé de nombreuses fois lors de mes séjours à Paris pour travailler avec Olivier. Merci aussi pour les discussions sur la recherche et le futur. Je te fais confiance pour assurer la descendance d'Olivier et honorer le badge de Michael Jackson que tu as récupéré à Birmingham. Merci à toi et à Edgar pour les trajets à 7h du matin gare du nord pour assister aux cours d'Ivan Marin. On a vécu l'aventure le premier jour en parcourant Amiens plusieurs fois pour mettre 1h45 à aller de la gare d'Amiens jusqu'au LAMFA.

Georges, tu m'as accueilli mon premier jour au LAMFA pour me présenter mon bureau et toute l'équipe des doctorants. Tu m'as ensuite fait assisté à des exposés uniques en leur genre au groupe de travail IMEST en parlant à un débit d'au moins 300 mots par minute et n'autorisant pas la moindre baisse d'attention. Tes histoires complotistes pendant la première année me permettaient de me concentrer à la perfection. Tu as pu te rapprocher d'I.M. à coups de Baklavas, je te félicite à nouveau pour l'obtention de son numéro. C'est un peu le comble pour un complotiste de trouver un emploi à "Bielefeld", je te souhaite une bonne continuation dans ta ville imaginaire. Merci à toi et Eirini d'avoir établi des standards très élevés pour la qualité de la thèse mais aussi et surtout pour la longueur des remerciements.

Henry, le petit dernier, ça a été un plaisir de te former cette année. Je pense que tu es déjà passé de l'état de larve à l'état de chrysalide. Tu seras à l'issue de ma soutenance le seul thésard du chef et il faudra que tu deviennes un papillon et que tu commences à voler de tes propres ailes.

\smallskip

Merci à Félix de m'avoir écouté lui expliquer que $d$ était en fait forcément inférieur ou égal à $2$ place de la République alors que ça n'avait a priori aucun intérêt pour lui qui a arrêté les mathématiques il y a bien longtemps. Merci à Planès de ne pas trop m'en vouloir de lui devoir encore 799 croissants 5 ans après le début de ma dette et d'avoir accepté de faire un peu de théorie des corps sur des serviettes à Burger King. Merci à Marc de nous avoir fait découvrir Saint-Crépin-Ibouvillers, son PMU et son kébab. Merci à Le Manach pour ses rappels sur l'ensemble vide qui m'ont remis les idées en place plus d'une fois. Merci à Sabrina pour son petit aller-retour à Amiens avant de quitter le pays qui m'a permis de dire "je me casse" à mes cobureaux un matin à 11h, merci à toi et Bastien pour votre super mariage où j'ai pu m'illustrer en faisant des maths immobile en costard et en plein soleil le lendemain matin pendant que les gens se baignaient dans la piscine. Merci aussi à Nicolas pour les billards qu'on a pu faire à défaut de babyfoot et pour le match de basket en 50 points serré jusqu'au bout qu'on a eu l'occasion de disputer. Thank you Loretta for all your support, coming to visit me in Amiens and initiating the email instead of facebook way of communication which will save the world. Thank you to both of you for hosting me in your home and letting me take care of Lucien a little bit, I hope he still remembers uncle Bop.

\smallskip

Je remercie mon oncle Alain et ma tante Maryse d'être venus à ma soutenance de thèse. J'espère que l'exposé ne sera pas trop pénible  à suivre.

I want to thank Chuck, Soraya, Taisa, Jacob and Danny for putting up with me even when I kept excluding maximal subgroups of the special linear groups of dimensions 6 and 8 even during the Olympics. I cannot express how grateful I am to Chuck, Soraya and Taisa for coming all the way from Miami to attend my defense. I still remember a talk with Chuck at Thomas' graduation where I told him I wasn't sure about starting a doctoral thesis. He told me about his experience which did not go so well but also that he never regretted beginning it. It helped me make my decision and I'm very thankful for this.

I also want to thank T-dawg for hosting me in Seattle multiple times during this thesis. His craziness helped me relax in this stressful period. Ezell's fried chicken and hiking in Washington state were also very helpful and I hope to do this many more times with you, Alice and Guthrie in the future.

I thank Skipper for coming to my defense on his trip across Europe. He was always ready to meet me in a cool restaurant in San Francisco.

\smallskip

I want to thank my mom for all her support all my life and during this thesis. Even though it was a tough couple of years, she was there all the time and gave all her support. Guillaume can testify that the phone rang multiple times a week and my answer was usually "yah". She even listened to me talk about my thesis and group theory. I'm pretty sure she still doesn't know what a perfect group is but she always asks and remembers for a few days. I know you do your best to be a great mother and sometimes feel unappreciated. I rarely say it but I'll write it down this time so you remember it, I love you mom. Also, thanks for paying for the red wine.

\smallskip

Je remercie mon père de m'avoir tenu informé avec assiduité des scores du Biarritz Olympique durant mon exil dans le nord.

\smallskip

Sophia, ma petite soeur chérie, je sais que tu ne veux plus que je t'appelle comme ça car tu n'as plus huit ans. Malheureusement pour toi, tu auras toujours huit ans pour moi. Je te remercie d'être venu  depuis le Québec pour ma soutenance. Tu ne sais pas ce qu'est le type B mais tu dois savoir que j'ai travaillé dessus toute la première année comme je te l'ai dit au téléphone à chaque fois que l'on se parlait. Tes cadeaux à base de patate et ton drapeau de la reine d'Angleterre rendent mon retour à l'appartement après le boulot plus sympatoches. Mes visites à Montréal m'ont permis de me changer les idées plus d'une fois même si tu as menti quand tu as dit qu'il allait faire en mars une température supportable pour un être humain. Merci aussi d'avoir accepté d'aller au bureau pour que je finisse vite fait le type B alors que tu venais de terminer ton baccalauréat. Heureusement, tu as pu profiter des hortillonnages qui ont dus être une expérience inoubliable pour toi. Tu es la meilleure petite soeur du monde et je te remercie d'accepter d'avoir huit ans jusqu'à la fin de ta vie.

\smallskip

Je continue mes remerciements des frères et soeurs de la plus jeune à la moins jeune avec mon frérot Thomas. Tu es en train de devenir le patriarche de la famille, merci de m'avoir aidé à organiser la soutenance et à réunir la famille. La prochaine grande occasion devrait être ton mariage. Je te remercie pour ton aide durant la thèse, notamment quand tu as été à l'initiative de l'achat de la PS4 qui m'a permis de mettre les bijections en pause plus d'une fois. Merci à Astrid pour sa bonne humeur à chaque fois que j'ai eu la chance de la voir, je ne pense pas connaître quelqu'un qui sourit autant, j'espère que Thomas a quelque chose à voir là-dedans.

\smallskip

Elsa, j'ai eu la chance de venir te voir à Annecy pendant la thèse plusieurs fois, c'était toujours avec grand plaisir. Tu m'as toujours impressionné par ta capacité de travail, j'espère avoir travaillé presque autant que toi pendant ces trois années. Tes compliments et encouragements m'ont aidé à avoir confiance en moi et à venir à bout de cette thèse. Tu m'as donné des précieux conseils pour l'enseignement, qui représente une partie très importante de ces dernières années et j'aurai encore plus besoin de tes conseils dans les années à venir.

\smallskip

Laure, je veux commencer par m'excuser pour ma tenue, tu n'es pas arrivée assez tôt pour qu'on puisse faire des courses de dernière minute cette fois-ci. Merci de m'avoir accueilli chez toi, j'ai découvert l'image des groupes Artin pour de nombreux types mais je ne pensais pas découvrir ce qu'est une figue pendant cette thèse.

\smallskip

Je peux remercier Cécile pour beaucoup de choses durant ces trois années. Je te remercie en particulier pour le voyage à Rome que nous avons fait avec Pablo pendant une des périodes les plus intenses de la thèse. Je pense que tu as apprécié le voyage aussi mais ça m'a permis de penser complétement à autre chose dans notre hôtel miteux au milieu d'une des plus belles villes du monde. J'ai hâte de pouvoir faire notre troisième voyage en Italie ensemble. Je remercie Pablo pour avoir été capable de sortir de belles phrases comme "Ces ruines, je m'en fiche comme de l'an 40" au Foro Romano à côté de ruines qui datent de l'an 40. Tu as été quand même assez agréable pour un jeune de 12 ans que l'on faisait marcher 20km par jour. Merci à Lio pour les histoires incroyables qu'elle m'a racontées au téléphone. Je remercie Gallien pour son soutien pendant la thèse, merci pour le match de rugby en loges à Bordeaux qu'il m'a offert et qui restera un très bon souvenir.

\smallskip

Je remercie enfin Anne pour son accueil pendant des phases obscures de ma thèse où mon esprit était parfois ailleurs. J'espère pouvoir être plus présent après ma thèse pour te soutenir à mon tour. Merci aussi à Zavier qui m'appelle déjà monsieur le professeur. Merci à mes nièces Afidy et Eding pour les petits exercices de mathématiques que vous m'envoyez de temps à autre et pour le jeu le logo des marques où j'ai eu l'occasion rêvée d'être "l'anim".

\smallskip

Pour finir, je veux remercier Mélyssa. Merci d'avoir changé ma vie un soir de juin où j'essayais de résoudre à la main un système à 2304 inconnues en me demandant, "ça vous dirait d'aller boire ou manger un truc un soir?". Ta rencontre est la meilleure chose qui me soit arrivée. Je t'ai rencontré le 9 Novembre 2016 après une nuit blanche où j'avais du mal à parler avant d'arriver, j'ai oublié mes soucis pendant 2 heures. Je ne savais pas encore à ce moment là à quel point tu es exceptionnelle et je n'imaginais pas une seconde que j'emménagerai avec toi quelques mois plus tard. Je ne pense pas que j'aurais pu réussir à venir à bout de ces derniers mois où nous nous levions à 7h47 tous les matins et où tu m'accompagnais à l'université que ce soit pendant la semaine ou le week-end ou les jours fériés. C'était beaucoup plus facile de mettre la tête dans la rédaction de ma thèse en me disant qu'à n'importe quel moment, la porte du bureaux des doctorants pouvait s'ouvrir avec toi derrière, un coca ou des croissants et chocolatines à la main pour me donner la force de continuer à travailler. Tes cadeaux surprises et tes petits mots que je trouve souvent dans les poches de mon jean ou mon portefeuille ou mes feuilles de calcul m'ont aidé à passer la plus belle année de ma vie alors que ça aurait probablement été une année d'épuisement physique et moral autrement. Merci aussi de m'avoir écouté plusieurs fois faire mon exposé de soutenance en notant à la seconde près le temps passé sur chaque diapositive et d'avoir noté mon étang. Cet exposé est le plus stressant que j'aurai à faire dans ma vie et il le sera beaucoup moins grâce à toi. Je t'aime et je te dédie cette thèse.

\newpage 
 
 $ $
 \newpage
 
 \section{Résumé en français de la thèse}

Nous déterminons dans cette thèse l'image des groupes de Artin associés à des groupes de Coxeter irréductibles dans leur algèbre de Iwahori-Hecke finie associée. Cela a été fait en type $A$ dans \cite{BM} et \cite{BMM}. Dans le cas générique, la clôture de Zariski de l'image a été déterminée dans tous les cas \cite{IH2}. L'approximation forte suggère que les résultats devraient être similaire dans le cas fini. Il est néanmoins impossible d'utiliser l'approximation forte sans utiliser de lourdes hypothèses et limiter l'étendue des résultats. Nous démontrons dans cette thèse que les résultats sont similaires mais que de nouveaux phénomènes interviennent de par la complexification des extensions de corps considérées. Les arguments principaux proviennent de la théorie des groupes finis. Nous utiliserons notamment un Théorème de Guralnick et Saxl \cite{GS} qui utilise la classification des groups finis simples pour les représentations de hautes dimensions. Ce théorème donne des conditions pour que des sous-groupes de groupes linéaires soient des groupes classiques dans une représentation naturelle. En petite dimension, nous utiliserons la classification des sous-groupes maximaux des groupes classiques de \cite{BHRC} pour les cas les plus compliqués.

Les résultats que nous démontrons peuvent avoir diverses applications. Par exemple, des groupes classiques finis ou des produits directs de groupes classiques finis apparaissent comme des quotients finis des groupes de Artin. Ces derniers sont des groupes fondamentaux de variétés algébriques donc cela définit des recouvrements finis de ces variétés qui peuvent être intéressant. Comme ces variétés sont définies sur le corps des rationnels $\Q$, cela peut avoir des applications au problème de Galois inverse (voir par exemple \cite{SV} pour un exemple en type $A_n$).

Ces résultats sont aussi intéressants du point de vue des groupes finis classiques. En effet, nous obtenons des générateurs explicites de ces groupes vérifiant les relations de tresses. Cela entraine des constructions intéressantes de ces groupes et de certains de leurs sous-groupes en s'intéressant à la restriction à des sous-groupes de Artin paraboliques. Nous obtenons par exemple une desciption intéressante du group $\op{Spin}_8^+(q)$ à l'aide des deux représentations de dimension $8$ en type $H_4$ (voir section \ref{sectiontriality}).

\bigskip

Dans cette thèse, nous donnerons tout d'abord dans la section \ref{DefCoxArtHec} une introduction aux groupes de Coxeter, aux groupes de Artin et aux algèbres de Iwahori-Hecke. Nous rappelerons la classification des groupes de Coxeter finis. Il y a quatre familles infinies $A_n$, $B_n$, $D_n$ et $I_2(n)$. Les groupes restant sont les groupes de Coxeter finis exceptionnels et sont noté $E_6$, $E_7$, $E_8$, $F_4$, $H_3$ et $H_4$. Ils correspondent tous à des objets géométriques et ont été classifiés en utilisant des argument géométriques par Coxeter en 1932 \cite{Cox}. Nous expliquerons ensuite comment définir les algèbres de Iwahori-Hecke dans un cadre général avant de donner des définitions plus précises sur les corps finis dans le chapitre correspondant à chaque type. Nous donnerons à la fin de cette section une idée des modèles pour les représentations irréductibles de ces algèbres.

Nous donnerons ensuite dans la section \ref{Symalg} des rappels sur les algèbres symétriques et les élements de Schur. Les éléments de Schur sont des outils qui permettent d'avoir un contrôle sur la semi-simplicité des algèbres symétriques. Les algèbres de Hecke sont des algèbres symétriques et leurs éléments de Schur dans le cas générique ont été déterminés dans tous les cas dans \cite{G-P}. Nous démontrerons après ces rappels une version du Théorème de déformation de Tits. Nous utiliserons cette version du théorème de Tits dans les différents chapitres pour montrer que sous les bonnes conditions sur les paramètres et la caractéristique du corps fini, nous pouvons spécialiser les modèles matriciels issus du cas générique aux corps finis $\F_q$. Les algèbres de Iwahori-Hecke seront alors semi-simples déployées et serons donc isomorphes à $\underset{\rho \mbox{ irr}}\oplus M_{n_\rho}(q)$. Nous pourrons alors considérer le groupe multiplicatif $A_{W_n}=<S_1,S_2,\dots,S_n>$ comme sous-groupe des éléments inversibles de l'algèbre $<S_1,S_2,\dots,S_n>$. Ce groupe sera donc un sous-groupes du produit $\underset{\rho \mbox{ irr}}\prod \op{GL_{n_\rho}(q)}$. Nous déterminerons à quoi ce groupe est isomorphe dans les différents types.

Dans la section \ref{sectionAschba}, nous rappelons le théorème de Aschbacher \cite{ASCHMAXSUBGRPS} sur les sous-groupes maximaux des groupes classiques finis. Il définit dans ce théorème $8$ classes de sous-groupes géométriques noté $\mathcal{C}_i$. Nous rappelons dans cette section une description rapide de ces classes. Nous donnerons ensuite des théorèmes qui permettent d'identifier les groupes classiques finis sous certaines conditions. Nous utiliserons ces théorèmes afin de déterminer l'image du groupe $\mathcal{A}_{W_n}$ dans les diverses représentations irréductibles des algèbres de Hecke sur les corps finis où $\mathcal{A}_{W_n}$ est le groupe dérivé de $A_{W_n}$.

\bigskip

Dans le chapitre \ref{TypeB}, nous déterminons l'image du groupe $\mathcal{A}_{B_n}$ dans son algèbre de Iwahori-Hecke associée. Nous définissons d'abord l'algèbre de Iwahori-Hecke finie sur le corps $\F_p(\alpha,\beta)$. Les représentations irréductibles en type $B_n$ sont indexées par des doubles-partitions de $n$. Les paramètres $\alpha$ et $\beta$ correspondent aux déformations de la relation d'ordre $2$ des groupes de Coxeter dans les algèbres de Iwahori-Hecke finies. Il y a deux paramètres en type $B_n$ car les générateurs ne sont pas tous conjugués. Nous devons alors considérer l'extension de corps $\F_p(\alpha,\beta)$ au-dessus de $\F_p(\alpha+\alpha^{-1},\beta+\beta^{-1})$. Les doubles-partitions avec une composante vide nous donne les mémes résultats qu'en type $A_n$. Ces résultats dépendent de l'extension de corps $\F_p(\alpha)$ au-dessus de $\F_p(\alpha+\alpha^{-1})$. Nous devrons alors distinguer les différents cas possibles pour ces extensions de corps. Cela donne $6$ possibilités et les résultats pour l'image de $\mathcal{A}_{B_n}$ sont alors différents. La preuve est dans tous les cas basée sur une récurrence où nous déterminons l'image pour $n\leq 5$ de diverses manières en regardant les représentations irréductibles une à une. Nous déterminons aussi les différentes factorisations (Proposition \ref{isomorphisme}) entre les représentations suivant les doubles-partitions qui les indexent. Ces factorisations dépendent des extensions de corps considérées et permettent de deviner quelle sera l'image de $\mathcal{A}_{B_n}$. Nous montrons ensuite par récurrence que le résultat annoncé est bien correct. Nous démontrons que les hypothèses du Théorème \ref{CGFS} sont vérifiées par $\rho(\mathcal{A}_{B_n})$ pour chaque représentation irréductible $\rho$. Nous donnons ci-dessous le résultat pour ce type dans le cas le plus simple et un des cas les plus compliqués. 

Notons $A_{1,n}=\{(\lambda_1,\emptyset),\lambda_1 \vdash n\}, A_{2,n} =\{(\emptyset,\lambda_2),\lambda_2 \vdash n\}, A_n = A_{1,n} \cup A_{2,n}$. $A\epsilon_n= \{(\lambda_1,\emptyset)\in A_{1,n},~\lambda_1 ~\mbox{pas une équerre}\}$, $ \epsilon_n=\{\lambda \vdash\vdash n, \lambda \notin A_n, \lambda~\mbox{pas une équerre}\}, \F_{\tilde{q}}=\F_p(\alpha)$. 
\begin{theor}
Si $\F_q=\F_p(\alpha,\beta)=\F_p(\alpha+\alpha^{-1}, \beta+\beta^{-1})$ et $\F_p(\alpha)=\F_p(\alpha+\alpha^{-1})$, alors le morphisme : $\mathcal{A}_{B_n} \rightarrow \mathcal{H}_{B_n,\alpha,\beta}^\times \simeq \underset{\lambda \vdash\vdash n}\prod GL(\lambda)$ se factorise à travers le morphisme surjectif
$$\Phi_{1,n} : \mathcal{A}_{B_n} \rightarrow SL_{n-1}(\tilde{q})\times \underset{(\lambda_1,\emptyset)\in A\epsilon_n, \lambda_1<\lambda_1'}\prod SL_{n_\lambda}(\tilde{q}) \times \underset{(\lambda_1,\emptyset)\in A\epsilon_n,\lambda_1=\lambda_1'}\prod OSP(\lambda)'\times$$
 $$SL_n(q)^2 \times \underset{\lambda\in \epsilon_n, \lambda < \lambda'}\prod SL_{n_\lambda}(q) \times \underset{\lambda\in \epsilon_n, \lambda=\lambda'}\prod OSP(\lambda)'.$$
\end{theor}

 \begin{theor}
Si $\F_q=\F_p(\alpha,\beta)=\F_p(\alpha+\alpha^{-1},\beta) \neq \F_p(\alpha,\beta+\beta^{-1}) = \F_p(\alpha+\alpha^{-1},\beta+\beta^{-1})$ et $\F_p(\alpha)= \F_p(\alpha+\alpha^{-1})$, alors le morphisme $\mathcal{A}_{B_n} \rightarrow \mathcal{H}_{B_n,\alpha,\beta}^\times \simeq \underset{\lambda \vdash\vdash n}\prod GL(\lambda)$ se factorise à travers le morphisme surjectif
$$\Phi_{6,n} : \mathcal{A}_{B_n} \rightarrow SL_{n-1}(\tilde{q})\times \underset{(\lambda_1,\emptyset)\in A\epsilon_n, \lambda_1<\lambda_1'}\prod SL_{n_{\lambda}}(\tilde{q}) \times \underset{(\lambda_1,\emptyset)\in A\epsilon_n,\lambda_1=\lambda_1'}\prod {OSP}(\lambda)'\times$$
 $$ SL_n(q) \times \underset{\lambda\in \epsilon_n, \lambda < \lambda', \lambda < (\lambda_1',\lambda_2'), \lambda\neq (\lambda_2,\lambda_1)}\prod SL_{n_\lambda}(q) \times \underset{\lambda \in \epsilon_n, \lambda < \lambda', \lambda=(\lambda_2,\lambda_1)}\prod SL_{n_{\lambda}}(q^{\frac{1}{2}}) \times $$
 $$ \underset{\lambda \in \epsilon_n, \lambda< \lambda',\lambda=(\lambda_1',\lambda_2')}\prod SU_{n_\lambda}(q^\frac{1}{2})\times  \underset{\lambda\in \epsilon_n, \lambda=\lambda', \lambda<(\lambda_1',\lambda_2')}\prod OSP(\lambda)' \times \underset{\lambda \in \epsilon_n, \lambda=\lambda',\lambda=(\lambda_1',\lambda_2')}\prod \widetilde{OSP}(\lambda)' .$$
 \end{theor}
 
 \bigskip
 
Dans le chapitre \ref{TypeD}, nous déterminons l'image du groupe $\mathcal{A}_{D_n}$ dans son algèbre de Iwahori-Hecke associée. Ici, l'algèbre est définie sur le corps $\F_p(\alpha)$ où $\alpha$ correspond au paramètre de déformation de la relation d'ordre $2$. Il n'y a qu'un seul paramètre car tous les générateurs sont conjugués. Cela rend les extensions de corps plus simples car il suffit de considérer $\F_p(\alpha)$ au-dessus de $\F_p(\alpha)$. Il faut utiliser la Proposition \ref{Tits} pour montrer que la spécialisation de cette algèbre est semi-simple déployée dans le Théorème \ref{Schur}.

Les représentations sont encore indexées par des doubles-partitions de $n$. Il y a dans ce cas-là un isomorphisme naturel entre la représentation indexée par la double-partition $(\lambda_1,\lambda_2$ et la représentation indexée par la représentation indexée par la double-partition $(\lambda_2,\lambda_1)$ car l'action sur le double-tableau $(\T_1,\T_2)$ correspond à l'action sur le double-tableau $(\T_2,\T_1)$. Les représentations indexées par des doubles-partitions de type $(\lambda_1,\lambda_1)$ ne sont alors plus irréductibles, elles se scindent en deux représentations irréductibles de même dimension que l'on note $(\lambda_1,\lambda_1,+)$ et $(\lambda_1,\lambda_1,-)$. La règle de branchement est ainsi plus complexe et est donnée ci-dessous.

\begin{lemmefr}
Soit $n\geq 5$ et $(\lambda,\mu)\Vdash n, \lambda>\mu$. On a alors :
\begin{enumerate}
\item Si $n_\lambda > n_\mu+1$, alors $V_{\lambda,\mu|\mathcal{H}_{D_{n-1},\alpha}}=\underset{(\tilde{\lambda},\tilde{\mu})\subset (\lambda,\mu)}\bigoplus V_{\tilde{\lambda},\tilde{\mu}}$.
\item Si $n_\lambda=n_\mu+1$ et $\mu \not\subset \lambda$, alors
$$V_{\lambda,\mu|\mathcal{H}_{D_{n-1},\alpha}}=(\underset{\tilde{\mu}\subset \mu}\bigoplus V_{\lambda,\tilde{\mu}}) \oplus (\underset{\tilde{\lambda}> \mu}{\underset{\tilde{\lambda}\subset \lambda}\bigoplus}V_{\tilde{\lambda},\mu})\oplus( \underset{\tilde{\lambda}< \mu}{\underset{\tilde{\lambda}\subset \lambda}\bigoplus}V_{\mu,\tilde{\lambda}}).$$
\item Si $n_\lambda=n_\mu+1$ et $\mu \subset \lambda$, alors
$$V_{\lambda,\mu|\mathcal{H}_{D_{n-1},\alpha}}=(\underset{\tilde{\mu}\subset \mu}\bigoplus  V_{\lambda,\tilde{\mu}}) \oplus (\underset{\tilde{\lambda}> \mu}{\underset{\tilde{\lambda}\subset \lambda}\bigoplus}V_{\tilde{\lambda},\mu})\oplus (\underset{\tilde{\lambda}< \mu}{\underset{\tilde{\lambda}\subset \lambda}\bigoplus}V_{\mu,\tilde{\lambda}})\oplus V_{\mu,\mu,+}\oplus V_{\mu,\mu,-}.$$
\item Si $n_\lambda=n_\mu$ et $\lambda>\mu$, alors $V_{\lambda,\mu|\mathcal{H}_{D_{n-1},\alpha}}=(\underset{\tilde{\mu}\subset \mu}\bigoplus V_{\lambda,\tilde{\mu}})\oplus (\underset{\tilde{\lambda}\subset \lambda}\bigoplus V_{\mu,\tilde{\lambda}}).$
\item Si $\lambda=\mu$, alors $V_{\lambda,\lambda,+|\mathcal{H}_{D_{n-1},\alpha}}=V_{\lambda,\lambda,-|\mathcal{H}_{D_{n-1},\alpha}}=\underset{\tilde{\mu} \subset \mu}\bigoplus V_{\lambda,\tilde{\mu}}.$
\end{enumerate}
\end{lemmefr}

Cette règle de branchement permet de faire un raisonnement par récurrence comme dans le chapitre \ref{TypeB}. Il faut alors traiter déterminer l'image de $\mathcal{A}_{D_n}$ dans l'algèbre de Iwahori-Hecke associée. On utilise pour cela les résultats de \cite{BMM} et les classifications des sous-groupes maximaux de certains groupes classiques sur des corps finis. La règle de branchement et le lemme de Goursat permettent d'obtenir le résultat suivant.

On écrit $A_{1,n}=\{(\lambda_1,\emptyset),\lambda_1 \vdash n\}, A_{2,n} =\{(\emptyset,\lambda_2),\lambda_2 \vdash n\}, A_n = A_{1,n} \cup A_{2,n}$ et\\
$\epsilon_n=\{\lambda \Vdash n, \lambda~\mbox{not a hook}\}$

\begin{theor}
Si $\F_q=\F_p(\alpha)=\F_p(\alpha+\alpha^{-1})$ et $n$ est impair, alors le morphisme de $\mathcal{A}_{D_n}$ dans $\mathcal{H}_{D_n,\alpha}^\times \simeq \underset{\lambda_1>\lambda_2}{\underset{\lambda \vdash\vdash n}\prod}GL_{n_\lambda}(q)$ se factorise à travers le morphisme surjectif
$$\Phi_{1',n}: \mathcal{A}_{D_n} \rightarrow SL_{n-1}(q) \times SL_n(q) \times \underset{\lambda_1>\lambda_2}{\underset{\lambda\in \epsilon_n,\lambda>\varphi(\lambda)}\prod} SL_{n_\lambda}(q)\times \underset{n_\lambda> n_\mu}{\underset{\lambda\in \epsilon_n,\lambda=\varphi(\lambda)}\prod}OSP(\lambda)'.$$ 
Si $\F_q=\F_p(\alpha)=\F_p(\alpha+\alpha^{-1})$ et $n \equiv 0~ (\bmod ~4)$, alors le morphisme de $\mathcal{A}_{D_n}$ dans $\mathcal{H}_{D_n,\alpha}^\times \simeq \underset{\lambda_1>\lambda_2}{\underset{\lambda \vdash\vdash n}\prod}GL_{n_\lambda}(\F_q) \times \underset{\lambda=(\lambda_1,\lambda_1)\vdash n}\prod GL_{n_{\lambda,+}}(q)\times GL_{n_{\lambda,-}}(q)$ se factorise à travers le morphisme surjectif
$$\Phi_{1',n}: \mathcal{A}_{D_n} \rightarrow SL_{n-1}(q) \times SL_n(q) \times \underset{\lambda_1>\lambda_2}{\underset{\lambda\in \epsilon_n,\lambda>\varphi(\lambda)}\prod} SL_{n_\lambda}(q)\times \underset{\lambda_1> \lambda_2}{\underset{\lambda\in \epsilon_n,\lambda=\varphi(\lambda)}\prod}OSP(\lambda)'\times$$ $$\underset{\lambda>\varphi(\lambda)}{\underset{\lambda=(\lambda_1,\lambda_1)\in \epsilon_n}\prod}SL_{\frac{n_\lambda}{2}}(q)^2\times \underset{\lambda=\varphi(\lambda)}{\underset{\lambda=(\lambda_1,\lambda_1)\in \epsilon_n}\prod}OSP(\lambda,+)'^2.$$ 
Si $\F_q=\F_p(\alpha)=\F_p(\alpha+\alpha^{-1})$ et $n \equiv 2~ (\bmod ~4)$ alors le morphisme de $\mathcal{A}_{D_n}$ dans $\mathcal{H}_{D_n,\alpha}^\times \simeq \underset{\lambda_1>\lambda_2}{\underset{\lambda \vdash\vdash n}\prod}GL_{n_\lambda}(\F_q) \times \underset{\lambda=(\lambda_1,\lambda_1)\vdash n}\prod GL_{n_{\lambda,+}}(q)\times GL_{n_{\lambda,-}}(q)$ se factorise à travers le morphisme surjectif
$$\Phi_{1',n}: \mathcal{A}_{D_n} \rightarrow SL_{n-1}(q) \times SL_n(q) \times \underset{\lambda_1>\lambda_2}{\underset{\lambda\in \epsilon_n,\lambda>\varphi(\lambda)}\prod} SL_{n_\lambda}(q)\times \underset{\lambda_1> \lambda_2}{\underset{\lambda\in \epsilon_n,\lambda=\varphi(\lambda)}\prod}OSP(\lambda)'\times$$ $$\underset{\lambda>\varphi(\lambda)}{\underset{\lambda=(\lambda_1,\lambda_1)\in \epsilon_n}\prod}SL_{\frac{n_\lambda}{2}}(q)^2\times \underset{\lambda=\varphi(\lambda)}{\underset{\lambda=(\lambda_1,\lambda_1)\in \epsilon_n}\prod}SL_{\frac{n_\lambda}{2}}(q).$$ 
Dans tout ce qui précéde, $OSP(\lambda)$ désigne le groupe des isométries de la forme bilinéaire définie en Proposition \ref{bilin2}.
\end{theor}

Le résultat correspondant dans le cas $\F_q=\F_p(\alpha)\neq \F_p(\alpha+\alpha^{-1})$ est similaire. Les groupes spéciaux linéaires sont alors remplacés par des groupes unitaires et les groupes symplectiques et orthogonaux sont définis sur des corps plus petits. Le résultat est donné dans le théorème \ref{result2D}. Cela conclut l'étude pour les cas classiques.

\bigskip

La seule famille infinie restante et la famille $I_2(m), m\geq 5$. Dans le cas $m$ impair, les deux générateurs du groupe de Coxeter sont conjugués. Dans le cas $m$ pair, ils ne le sont et on a donc deux paramètres pour l'algèbre de Iwahori-Hecke. On sépare donc l'étude dans le chapitre \ref{dihedralchapter} en deux sections suivant la parité de $m$. Dans les deux cas, les représentations irréductibles sont de dimension $1$ ou $2$. Nous déterminons dans les théorèmes \ref{platypodes} et \ref{bopb} alors l'image de $\rho(\mathcal{A}_{I_2(m)}$ pour les représentations irréductibles de dimension $2$ en utilisant le Théorème de Dickson (voir \cite{HUP} Théorème $8.27$) qui classifie les sous-groupes de $SL_2(q)$. La difficulté principale pour ces types provient des différentes factorisations possibles et de l'étude des extensions de corps lorsque $m$ est pair. Les différentes extensions de corps sont décrites dans les Figures \ref{fieldsIcase1} à \ref{fieldsIcase7}. Pour étudier les factorisations possibles, il faut introduire une relation d'équivalence sur les entiers. Elle est donnée dans les lemmes suivants qui dépendent de la parité de $m$.
\begin{lemmefr}
Supposons $m$ impair et $\xi_j=\theta^j+\theta^{-j}$ où $\theta$ est une racine primitive $m$-ième de l'unité dans $\overline{\F_p}$. Soit $j,l\in[\![1,\frac{m-1}{2}]\!]^2$. Il existe un automorphisme $\Psi_{l,j}$ de $\F_{q_j}=\F_p(\alpha,\xi^j+\xi^{-j})$ qui vérifie $\Psi_{l,j}(\alpha+\alpha^{-1})=\alpha+\alpha^{-1}$ et $\Psi_{l,j}(\xi^j+\xi^{-j})=\xi^l+\xi^{-l}$ si et seulement si il existe $r\in \N$ tel que $jp^r\equiv l~(\bmod ~ m)$ ou $jp^r\equiv -l~(\bmod~ m)$ et $(\alpha+\alpha^{-1})^{p^r}=\alpha+\alpha^{-1}$.

On dit que $j\sim l$ si une de ces conditions est vérifiée. Cela définit une relation d'équivalence et lorsque $j\sim l$, on a $\rho_{l|\mathcal{A}_{I_2(m)}}=\Psi_{l,j}\circ \rho_{j|\mathcal{A}_{I_2(m)}}$.
\end{lemmefr}
\begin{lemmefr}
Supposons $m$ pair et $\xi_j=\theta^j+\theta^{-j}$ où $\theta$ est une racine primitive $m$-ième de l'unité dans $\overline{\F_p}$. On dit que $j\sim l$ si $\F_p(\alpha+\alpha^{-1},\beta+\beta^{-1},\xi_j)\simeq \F_p(\alpha+\alpha^{-1},\beta+\beta^{-1},\xi_j)$ et il existe $\Phi_{j,l}\in Aut(\F_{q_j})$ tel que $\Phi_{j,l}(\alpha+\alpha^{-1})=\alpha+\alpha^{-1})$, $\Phi_{j,l}(\beta+\beta^{-1})=\beta+\beta^{-1}$ et $\Phi_{j,l}(\xi_j)=\xi_l$. Cela définit une relation d'équivalence et si $j\sim l$ alors $\Phi_{j,l}\circ \rho_{j|\mathcal{A}_{I_2(m)}} \simeq \rho_{l|\mathcal{A}_{I_2(m)}}$.
\end{lemmefr}

On a alors les théorèmes suivants

\begin{theor}
Supposons $m$ impair et que $\alpha$ vérifie les conditions données au début de la section \ref{sectionmodd}. Pour $j \in [\![1,\frac{m-1}{2}]\!]$, on pose $G_j=SL_2(q_j)$ si $\F_{q_j}=\F_p(\alpha,\theta^j+\theta^{-j})=\F_p(\alpha+\alpha^{-1},\theta^j+\theta^{-j})$ et $G_j\simeq SU_2(q_j^\frac{1}{2})$ si $\F_{q_j}=\F_p(\alpha,\theta^j+\theta^{-j})\neq\F_p(\alpha+\alpha^{-1},\theta^j+\theta^{-j})$.

On a alors que le morphisme de $\mathcal{A}_{I_2(m)}$ dans $\mathcal{H}_{I_2(m),\alpha}^\times\simeq GL_1(q_j)^2 \times \underset{j\in [\![1,\frac{m-1}{2}]\!]}\prod GL_2(q_j)$ se factorise à travers le morphisme surjectif
$$\Phi : \mathcal{A}_{I_2(m)} \rightarrow \underset{j\in [\![1,\frac{m-1}{2}]\!]/\sim}\prod G_j.$$
\end{theor}
\begin{theor}
Supposons $m$ pair et que $\alpha$ et $\beta$ vérifient les conditions données au début de la section \ref{sectionmeven}. Pour $j \in [\![1,\frac{m-2}{2}]\!]$, on pose $G_j=\rho_j([<T_t,T_s>,<T_t,T_s>])$ .

On alors que le morphisme de $\mathcal{A}_{I_2(m)}$ dans $\mathcal{H}_{I_2(m),q}^\times \simeq GL_1(q_j)^2 \times \underset{j\in [\![1,\frac{m-1}{2}]\!]}\prod GL_2(q_j)$ se factorise par le morphisme surjectif
$$\Phi : \mathcal{A}_{I_2(m)} \rightarrow \underset{j\in [\![1,\frac{m-2}{2}]\!]/\sim}\prod G_j.$$
\end{theor}

\bigskip

Cela conclut l'étude pour les familles infinies de groupes de Coxeter finis irréductibles. Il reste ensuite à traiter les groupes de Coxeter exceptionnels. C'est à dire les algèbres de Hecke associées aux groupes de Coxeter de type $E_6$, $E_7$, $E_8$, $H_3$, $H_4$ et $F_4$. Les représentations pour ces algèbres de Hecke sont données par des $W$-graphes. Nous décrivons ces objets dans le chapitre \ref{Wgraphschapter}. La définition est la suivante. Soit $W$ un

\begin{Deffr}
Soit $W$ un groupe de Coxeter, $K'$ son corps de définition, $\mathcal{H}$ son algèbre de Iwahori-Hecke de paramètres $(\alpha_s)_{s\in S}$ et $K=K'((\alpha_s)_{s\in S})$. Pour $X$ un ensemble, on note $D(X)=\{(x,x),x\in X\}$ sa diagonale. Un $W$-graphe $\Gamma$ est la donnée d'un triplet $(X,I,\mu)$ tel que
\begin{enumerate}
\item $X$ est un ensemble et $I$ est une application de $X$ dans $\mathcal{P}(S)$,
\item $\mu$ est une application de $(X\times X\setminus D(X)\times S)$ dans $K$ stable par l'involution du corps $K$ qui envoie $\sqrt{\alpha_s}$ sur $\sqrt{\alpha_s}^{-1}$.

Soit $V$ le $K'((\alpha_s)_{s\in S})$-espace vectoriel de base $(e_y)_{y\in X}$. Pour tout $s\in S$, on définit $\rho_s :V\rightarrow V$ par
\begin{eqnarray*}
e_y   \mapsto &  -e_y & \textsf{if}~ s\in I(y),\\
e_y   \mapsto & \alpha_s e_y+\underset{x\in X,s\in I(x)}\sum \sqrt{\alpha_s}\mu_{x,y}^s e_x & \textsf{if}~ s\notin I(y).
\end{eqnarray*}

\item L'application $T_s\mapsto \rho_s$ est une représenation de $\mathcal{H}$.
\end{enumerate}

\end{Deffr}

On sait qu'un tel $W$-graphe existe pour n'importe quel représentation irréductible d'une algèbre de Iwahori-Hecke \cite{Gyo}. Dans le chapitre \ref{Wgraphschapter}, nous donnons des propriétés sur les représentations qui peuvent se déduire de ces modèles. La $2$-colorabilité est une notion qui intervient dans plusieurs propositions de ce chapitre. Les $W$-graphes $2$-coloriables dans le cas des paramètres égaux ont été classifié par Gyoja (Voir la remarque après le Théorème \ref{Melyssaestparfaite}). Les notions de $W$-graphe dual de \cite{G-P} et de représentation auto-duale sont données dans la Proposition \ref{Melyssaesttropforte} et la Définition \ref{defselfdualrepresentation}. Le théorème principal de cette section est le suivant

\begin{theor}
Soit $\Gamma=(X,I,\mu)$ un $W$-graphe associé à une représentation irréductible de $\mathcal{H}$ tel que $\Gamma$ soit $2$-coloriable et tel que $\Gamma$ soit isomorphe en tant que graphe orientée pondéré au graphe $(X,\tilde{I},-\tilde{\mu})$.

\smallskip

Soit $\varphi: X\rightarrow X$ l'automorphisme de graphe de $\Gamma$ dans $(X,\tilde{I},-\tilde{\mu})$ et $x_1,x_2,\dots,x_n$ une numérotation de $X$ telle que $\varphi(x_i)=x_{n+1-i}$.

\smallskip

Soit $\langle .,.\rangle$ la forme bilinéaire définie par $<e_{x_i},e_{x_j}>=\omega(e_{x_i})\delta_{i,n+1-j}$, où $\omega$ est un coloriage de $\Gamma$ par $1$ et $-1$.

\smallskip

On a alors
$$\forall s\in S,\forall v_1,v_2\in V, \langle \rho_\Gamma(T_s)v_1,\rho_\Gamma(T_s)v_2\rangle =-\alpha\langle v_1,v_2\rangle.$$

\smallskip

Cette forme bilinéaire est non-dégénérée et elle est symétrique si $\omega(x_1)\omega(x_n)=1$ et anti-symétrique si $\omega(x_1)\omega(x_n)=-1$.

La représentation associée est alors auto-duale.
\end{theor}

On dit ensuite qu'un $W$-graphe auto-dual est un $W$-graphe vérifiant les propriétés du Théorème précédent. Cette propriété est vérifiée par certains $W$-graphes associés à des représentations irréductibles auto-duales des algèbres de Hecke mais pas par tous. Les $W$-graphes n'étant pas unique pour une représentation donnée, on établit la conjecture suivante.

\begin{conjfr}
Soit $W$ un groupe de Coxeter. Pour toute représentation irréductible auto-duale, il existe un $W$-graphe auto-dual défini sur $K$ assoicé cette représentation. Si il existe un $W$-graphe $\Gamma$ défini sur $K'$ associé à la représentation alors il existe un $W$-graphe auto-dual $\Gamma'$ associé à la représentation et une matrice $M\in GL_{\vert X\vert}(\tilde{K})$ telle que pour tout $h\in \mathcal{H}_{K}$, $M\rho_{\Gamma}(h)M^{-1}=\rho_{\Gamma'}(h)$.
\end{conjfr}

La deuxième partie de la conjecture est formulée uniquement dans l'optique de montrer la première partie de la conjecture. Afin d'utiliser cette deuxième partie, nous montrons deux conditions restrictives pour deux $W$-graphes soient associés à la même représentation irréductible dans les Propositions \ref{Unique1} et \ref{Unique2}. Nous avons ensuite démontré par des calculs avec la plateforme de calcul HPC MatriCS \cite{Matrics} la conjecture pour les types $E_6$, $E_7$, $E_8$, $H_3$ et $H_4$. Les $W$-graphes auto-duaux obtenus peuvent être téléchargés sur \cite{newgraphsEsterle}.

\bigskip

Les résultats du chapitre \ref{Wgraphschapter} permettent d'étudier les types exceptionnels. Dans le chapitre \ref{TypeEchapter}, nous déterminons l'image de $\mathcal{A}_{E_6}$, $\mathcal{A}_{E_7}$ et $\mathcal{A}_{E_8}$ dans leur algèbre de Iwahori-Hecke associée. Le groupe de Artin $A_{D_5}$ s'injecte naturellement dans le groupe de Artin $A_{E_6}$. Cela permet d'utiliser les résultats du chapitre \ref{TypeD} pour appliquer un raisonnement par récurrence.

Nous montrons dans un premier temps en utilisant la Proposition \ref{Tits} et les éléments de Schur des algèbres de Iwahori-Hecke dans le cas générique que ces algèbres de Iwahori-Hecke sont semi-simples déployées après spécialisation sous les bonnes conditions sur les paramètres de l'algèbre.

Dans la section \ref{E6section}, nous utilisons le Théorème \ref{CGFS} et les Théorème \ref{result1D} et \ref{result2D} pour déterminer l'image de $\mathcal{A}_{E_6}$ dans chaque représentation irréductible de l'algèbre de Iwahori-Hecke $\mathcal{H}_{E_6,\alpha}$. Les $W$-graphes auto-duaux obtenus et disponibles sur \cite{newgraphsEsterle} permettent de déterminer quel type de forme bilinéaire est préservée par les $\rho(\mathcal{A}_{E_6})$ pour les représentations irréductibles auto-duales. Les nouveaux $E_6$-graphes auto-duaux de dimension $10$ et $20$ sont donnés dans la section \ref{E6graphssection}. On voit sur ces figures que l'opération de symétrie sur les graphes inverse les couleurs, les formes bilinéaires associées sont donc anti-symétriques. C'est le cas pour toutes les représentations irréducibles auto-duales en type $E_6$. On utilise enfin le Lemme de Goursat (Lemme \ref{Goursat}) pour récupérer l'image totale dans l'algèbre de Iwahori-Hecke $\mathcal{A}_{E_6,\alpha}$.

\begin{theor}
Si $\F_q=\F_p(\alpha)=\F_p(\alpha+\alpha^{-1})$, alors le morphisme de $\mathcal{A}_{E_6}$ dans $\mathcal{H}_{E_6,\alpha}^\star\simeq \underset{\rho \mbox{ irr}}\prod GL_{n_\rho}(q)$ se factorise à travers le morphisme surjectif

$$\Phi : \mathcal{A}_{E_6}\rightarrow SL_6(q)\times SP_{10}(q)\times SL_{15}(q)^2\times SL_{20}(q)\times SP_{20}(q)\times SL_{24}(q)\times SL_{30}(q)$$
$$\times SP_{60}(q)\times SL_{60}(q)\times SL_{64}(q)\times SP_{80}(q)\times SL_{81}(q)\times SP_{90}(q).$$

Si $\F_q=\F_p(\alpha)\neq\F_p(\alpha+\alpha^{-1})$, alors le morphisme de $\mathcal{A}_{E_6}$ dans $\mathcal{H}_{E_6,\alpha}^\star\simeq \underset{\rho \mbox{ irr}}\prod GL_{n_\rho}(q)$ se factorise à travers le morphisme

$$\Phi : \mathcal{A}_{E_6}\rightarrow SU_6(q^{\frac{1}{2}})\times SP_{10}(q^{\frac{1}{2}})\times SU_{15}(q^{\frac{1}{2}})^2\times SU_{20}(q^{\frac{1}{2}})\times SP_{20}(q^{\frac{1}{2}})\times SU_{24}(q^{\frac{1}{2}})\times SU_{30}(q^{\frac{1}{2}})$$
$$\times SP_{60}(q^{\frac{1}{2}})\times SU_{60}(q^{\frac{1}{2}})\times SU_{64}(q^{\frac{1}{2}})\times SP_{80}(q^{\frac{1}{2}})\times SU_{81}(q^{\frac{1}{2}})\times SP_{90}(q^{\frac{1}{2}}).$$

\end{theor}

Dans la section \ref{E7section}, nous utilisons les résultats montrés dans la section \ref{E6section}. Il n'existe aucune représentation irréductible auto-duale de $\mathcal{H}_{E_7,\alpha}$. Ainsi, il n'y a aucune forme bilinéaire à considérer dans ce cas-là. Le phénomène nouveau apparaissant dans le type $E_7$ est l'existence de deux représentations irréductibles de dimension $512$ qui ne sont pas 2-coloriables. Il faut alors considérer l'extension de corps $\F_p(\sqrt{\alpha})$ au-dessus de $\F_p(\alpha)$ et l'automorphisme de corps d'ordre $2$ associé lorsque ces deux corps sont distincts. Les factorisations possibles sont alors déterminées en utilisant la Proposition \ref{Fieldfactorization} et en calculant les traces d'éléments bien choisis à l'aide du package CHEVIE \cite{CHEVIE}. Les arguments habituels permettent alors de démontrer le résultat suivant.

\begin{theor}
Si $F_q=\F_p(\sqrt{\alpha})=\F_p(\alpha+\alpha^{-1})$, alors le morphisme de $\mathcal{A}_{E_7}$ dans $\mathcal{H}_{E_7,\alpha}^\star\simeq \underset{\rho \mbox{ irr}}\prod GL_{n_\rho}(q)$ se factorise à travers le morphisme surjectif

$$\Phi : \mathcal{A}_{E_7}\rightarrow SL_7(q)\times SL_{15}(q)\times SL_{21}(q)^2\times SL_{27}(q)\times SL_{35}(q)^2\times SL_{56}(q)\times SL_{70}(q)\times SL_{84}(q)$$
$$\times SL_{105}(q)^3
\times SL_{120}(q)\times SL_{168}(q)\times SL_{189}(q)^3\times SL_{210}(q)^2\times SL_{216}(q)\times SL_{280}(q)^2$$
$$\times SL_{315}(q)\times SL_{336}(q)\times SL_{378}(q)\times SL_{405}(q) \times SL_{420}(q)\times SL_{512}(q).$$

Si $\F_p(\sqrt{\alpha})\neq \F_q=\F_p(\alpha)=\F_p(\alpha+\alpha^{-1})$ alors le morphisme de $\mathcal{A}_{E_7}$ dans $\mathcal{H}_{E_7,\alpha}^\star\simeq \underset{\rho \mbox{ irr}}\prod GL_{n_\rho}(q)$ se factorise à travers le morphisme

$$\Phi : \mathcal{A}_{E_7}\rightarrow SL_7(q)\times SL_{15}(q)\times SL_{21}(q)^2\times SL_{27}(q)\times SL_{35}(q)^2\times SL_{56}(q)\times SL_{70}(q)\times SL_{84}(q)$$
$$\times SL_{105}(q)^3
\times SL_{120}(q)\times SL_{168}(q)\times SL_{189}(q)^3\times SL_{210}(q)^2\times SL_{216}(q)\times SL_{280}(q)^2$$
$$\times SL_{315}(q)\times SL_{336}(q)\times SL_{378}(q)\times SL_{405}(q) \times SL_{420}(q)\times SU_{512}(q).$$

Si $F_q=\F_p(\alpha)\neq\F_p(\alpha+\alpha^{-1})$ alors le morphisme de $\mathcal{A}_{E_7}$ dans $\mathcal{H}_{E_7,\alpha}^\star\simeq \underset{\rho \mbox{ irr}}\prod GL_{n_\rho}(q)$ se factorise à travers le morphisme

\begin{small}
$$\Phi : \mathcal{A}_{E_7}\rightarrow SU_7(q^{\frac{1}{2}})\times SU_{15}(q^{\frac{1}{2}})\times SU_{21}(q^{\frac{1}{2}})^2\times SU_{27}(q^{\frac{1}{2}})\times SU_{35}(q^{\frac{1}{2}})^2\times SU_{56}(q^{\frac{1}{2}})\times SU_{70}(q^{\frac{1}{2}})\times SU_{84}(q^{\frac{1}{2}})$$
$$\times SU_{105}(q^{\frac{1}{2}})^3
\times SU_{120}(q^{\frac{1}{2}})\times SU_{168}(q^{\frac{1}{2}})\times SU_{189}(q^{\frac{1}{2}})^3\times SU_{210}(q^{\frac{1}{2}})^2\times SU_{216}(q^{\frac{1}{2}})\times SU_{280}(q^{\frac{1}{2}})^2$$
$$\times SU_{315}(q^{\frac{1}{2}})\times SU_{336}(q^{\frac{1}{2}})\times SU_{378}(q^{\frac{1}{2}})\times SU_{405}(q^{\frac{1}{2}}) \times SU_{420}(q^{\frac{1}{2}})\times SU_{512}(q^{\frac{1}{2}}).$$
\end{small}
\end{theor}

Dans la section \ref{E8section}, nous déterminons l'image de $\mathcal{A}_{E_8}$ dans son algèbre de Iwahori-Hecke associée. Les preuves se font par récurrence en utilisant les résultats de la section \ref{E7section}. La difficulté provient principalement des dimensions des représentations irréductibles auto-duales. La représentation irréductible auto-duale de plus haute dimension est de dimension $7168$. En utilisant la conjecture, l'obtention de la forme bilinéaire associée nécessite plus d'une semaine de calcul. Nous avons utilisé la conjecture pour obtenir toutes les formes bilinéaires associées aux représentations disponibles dans le package CHEVIE \cite{CHEVIE} de GAP. Une fois les formes bilinéaires obtenues, nous avons démontré la conjecture pour chaque représentation irréductible auto-duale. Les $E_8$-graphes auto-duaux obtenus sont téléchargeables depuis \cite{newgraphsEsterle}. Cela permet d'obtenir les formes bilinéaires correspondantes en utilisant uniquement un $2$-coloriage de la représentation. En utilisant ces résultats, nous avons démontré le résultat suivant qui conclut le chapitre \ref{TypeEchapter}.

\begin{theor}
Soit $A$ un ensemble de représentants de représentations irréductibles $2$-coloriables non auto-duales pour la relation d'équivalence $\rho \thickapprox \varphi$ si $\rho=\varphi'$ et $B$ l'ensemble des représentations irréductibles auto-duales.

Si $\F_p(\sqrt{\alpha})\neq\F_q=\F_p(\alpha)=\F_p(\alpha+\alpha^{-1})$ alors le morphisme de $\mathcal{A}_{E_8}$ dans $\mathcal{H}_{E_8,\alpha}^\star\simeq \underset{\rho \mbox{ irr}}\prod GL_{n_\rho}(q)$ se factorise à travers le morphisme surjectif
$$\Phi : \mathcal{A}_{E_8} \rightarrow \underset{\rho\in A}\prod SL_{n_{\rho}}(q) \times SL_{4096}(q^2)\times \underset{\rho\in B}\prod \Omega_{n_\rho}^+(q).$$

Si $\F_q=\F_p(\sqrt{\alpha})=\F_p(\alpha)=\F_p(\alpha+\alpha^{-1})$, alors le morphisme de $\mathcal{A}_{E_8}$ dans $\mathcal{H}_{E_8,\alpha}^\star\simeq \underset{\rho \mbox{ irr}}\prod GL_{n_\rho}(q)$ se factorise à travers le morphisme surjectif
$$\Phi : \mathcal{A}_{E_8} \rightarrow \underset{\rho\in A}\prod SL_{n_{\rho}}(q) \times SL_{4096}(q)^2\times \underset{\rho\in B}\prod \Omega_{n_\rho}^+(q).$$

Si $\F_q=\F_p(\alpha)\neq\F_p(\alpha+\alpha^{-1})$, alors le morphisme de $\mathcal{A}_{E_8}$ dans $\mathcal{H}_{E_8,\alpha}^\star\simeq \underset{\rho \mbox{ irr}}\prod GL_{n_\rho}(q)$ se factorise à travers le morphisme surjectif
$$\Phi : \mathcal{A}_{E_8} \rightarrow \underset{\rho\in A}\prod SU_{n_{\rho}}(q^{\frac{1}{2}}) \times SU_{4096}(q^{\frac{1}{2}})\times \underset{\rho\in B}\prod \Omega_{n_\rho}^+(q^\frac{1}{2}).$$
\end{theor}

\bigskip

Dans le chapitre \ref{TypeH}, nous considérons l'image des groupes $\mathcal{A}_{H_3}$ et $\mathcal{A}_{H_4}$ dans leur algèbre de Iwahori-Hecke finie associée. L'inclusion naturelle de $A_{I_2(5)}$ dans $A_{H_3}$ fait intervernir les résultats de la section \ref{sectionmodd}. Le corps de définition de $H_3$ et de $H_4$ est $\Q[\sqrt{5}]$, la relation d'équivalence définie dans la section \ref{sectionmodd} appliquée au cas $m=5$ sera alors nécessaire pour bien distinguer les différents cas. Il y a peu de représentations irréductibles en type $H_3$. En utilisant les résultats de la section \ref{sectionmodd} et les classifications de sous-groupes maximaux classiques en petite dimension \cite{BHRC}, nous obtenons le résultat suivant.

\begin{theor}
Supposons $p\notin \{2,5\}$ et que l'ordre de $\alpha$ ne divise ni $6$ ni $20$.
\begin{enumerate}
\item Supposons $1\sim 2$.
\begin{enumerate}
\item Si $\F_q=\F_p(\sqrt{\alpha})=\F_p(\alpha)=\F_p(\alpha+\alpha^{-1})$, alors le morphisme de $\mathcal{A}_{H_3}$ dans $\mathcal{H}_{H_3,\alpha}^\star \simeq GL_1(q)^2 \times GL_3(q)^2\times GL_4(q)^2\times GL_5(q)$ se factorise à travers le morphisme surjectif
$$\Phi : \mathcal{A}_{H_3} \rightarrow SL_3(q^2)\times SL_4(q)\times SL_5(q).$$
\item Si $\F_q=\F_p(\sqrt{\alpha})=\F_p(\alpha)\neq \F_p(\alpha+\alpha^{-1})$ et $\Phi_{1,2}(\sqrt{\alpha})=\sqrt{\alpha}^{-1}$ alors le morphisme de $\mathcal{A}_{H_3}$ dans $\mathcal{H}_{H_3,\alpha}^\star \simeq GL_1(q)^2 \times GL_3(q)^2\times GL_4(q)^2\times GL_5(q)$ se factorise à travers le morphisme surjectif
$$\Phi : \mathcal{A}_{H_3} \rightarrow SL_3(q)\times SU_4(q^{\frac{1}{2}})\times SU_5(q^{\frac{1}{2}}).$$
\item Si $\F_q=\F_p(\sqrt{\alpha})=\F_p(\alpha)\neq \F_p(\alpha+\alpha^{-1})$ et $\Phi_{1,2}(\sqrt{\alpha})=-\sqrt{\alpha}^{-1}$, alors le morphisme de $\mathcal{A}_{H_3}$ dans $\mathcal{H}_{H_3,\alpha}^\star \simeq GL_1(q)^2 \times GL_3(q)^2\times GL_4(q)^2\times GL_5(q)$ se factorise à travers le morphisme surjectif
$$\Phi : \mathcal{A}_{H_3} \rightarrow SL_3(q)\times SL_4(q^{\frac{1}{2}})\times SU_5(q^{\frac{1}{2}}).$$
\item Si $\F_{q^2}=\F_p(\sqrt{\alpha})\neq \F_p(\alpha)=\F_p(\alpha+\alpha^{-1})$, alors le morphisme de $\mathcal{A}_{H_3}$ dans $\mathcal{H}_{H_3,\alpha}^\star \simeq GL_1(q)^2 \times GL_3(q)^2\times GL_4(q)^2\times GL_5(q)$ se factorise à travers le morphisme surjectif
$$\Phi : \mathcal{A}_{H_3} \rightarrow SL_3(q^2)\times SU_4(q)\times SL_5(q).$$
\end{enumerate}
\item Supposons $1\nsim 2$. Lorsqu'il existe, on note $\epsilon$ l'automorphisme d'ordre $2$ de $\F_q$.
\begin{enumerate}
\item Si $\F_q=\F_p(\sqrt{\alpha})=\F_p(\alpha)\neq \F_p(\alpha+\alpha^{-1})$ et $\epsilon(\sqrt{\alpha})=\sqrt{\alpha}^{-1}$, alors le morphisme de $\mathcal{A}_{H_3}$ dans $\mathcal{H}_{H_3,\alpha}^\star \simeq GL_1(q)^2 \times GL_3(q)^2\times GL_4(q)^2\times GL_5(q)$ se factorise à travers le morphisme surjectif
$$\Phi : \mathcal{A}_{H_3} \rightarrow SU_3(q^{\frac{1}{2}})^2\times SU_4(q^{\frac{1}{2}})\times SU_5(q^{\frac{1}{2}}).$$
\item Si $\F_q=\F_p(\sqrt{\alpha})=\F_p(\alpha)\neq \F_p(\alpha+\alpha^{-1})$ et $\epsilon(\sqrt{\alpha})=-\sqrt{\alpha}^{-1}$, alors le morphisme de $\mathcal{A}_{H_3}$ dans $\mathcal{H}_{H_3,\alpha}^\star \simeq GL_1(q)^2 \times GL_3(q)^2\times GL_4(q)^2\times GL_5(q)$ se factorise à travers le morphisme surjectif
$$\Phi : \mathcal{A}_{H_3} \rightarrow SU_3(q^{\frac{1}{2}})^2\times SL_4(q^{\frac{1}{2}})\times SU_5(q^{\frac{1}{2}}).$$
\end{enumerate}
\end{enumerate} 
 
\end{theor}

\smallskip

L'étude en type $H_4$ est nettement plus compliquée. De nombreuses représentations irréductibles sont auto-duales. De plus, il y a quatre représenations irréductibles de dimension $16$. Les $H_4$-graphes sont donnés dans la section \ref{sectionnewH4graphs} de l'appendice. Il est clair en étudiant la symétrie des graphes auto-duaux munis d'un $2$-coloriage que les forme bilinéaires associées sont toutes symétriques. En particulier, les deux représentations auto-duales de dimension $8$ font intervenir le groupe $\Omega_8^+(q)$. Cela complique considérablement l'étude de l'image de $\mathcal{A}_{H_4}$ au sein de ces représentations. On sépare donc l'étude de $\mathcal{A}_{H_4}$ en $4$ sous-sections.

Dans la section \ref{H4gen}, nous démontrons que l'on a bien les propriétés voulus pour la spécialisation du modèle aux corps finis. Nous montrons aussi certaines propriétés générales qui seront utiles pour l'étude de chaque représentation.

Dans la section \ref{H4lowdim}, nous déterminons l'image de $\mathcal{A}_{H_4}$ dans les représentations de dimension inférieures à $8$. Pour déterminer $\rho_{8_r}(\mathcal{A}_{H_4}$, nous montrons que ce groupes contient un groupe suffisament grand. Pour montrer cela, nous utilisons des preuves très calculatoires qui font intervenir le Théorème de Dickson (voir \cite{HUP} Théorème $8.27$) et le Lemme de Goursat (voir Lemme \ref{Goursat}) à diverses reprises. Il faut alors montrer que des quantités sont non-nulles en effectuant des opérations qui s'apparentent à des divisions euclidiennes de polynômes en $\alpha$. Certains de ces calculs ont été mis dans la section \ref{ComputationsH4}.

Dans la section \ref{sectiontriality}, nous montrons que les représentations $\rho_{8_r}$ et $\rho_{8_{rr}}$ sont liées par la trialité. Cela fait apparaître une construction intéressante du groupe $\op{Spin}_8^+(q)$.

Enfin, dans la section \ref{H4highdim}, nous déterminons l'image de $\mathcal{A}_{H_4}$ dans les représentations de dimensions supérieures à $9$. Nous concluons la section par l'image de $\mathcal{A}_{H_4}$ dans toute l'algèbre de Iwahori-Hecke. Un des cas pour la représentation de dimension $48$ ne peut pas se traiter à l'aide des arguments habituels. Nous émettons alors la conjecture \ref{fingerscrossed} qui donne ce à quoi devrait être isomorphe $\rho_{48_rr}(\mathcal{A}_{H_4})$. Nous avons alors le résultat suivant qui conclut le chapitre \ref{TypeH}.

\begin{theor}
Supposons $p\notin \{2,3,5\}$ et $\alpha$ d'ordre ne divisant ni $20$, ni $30$, ni $48$.
\begin{enumerate}
\item Supposons $1\sim 2$ et que la conjecture \ref{fingerscrossed} est vraie.
\begin{enumerate}
\item Si $\F_q=\F_p(\sqrt{\alpha})=\F_p(\alpha)=\F_p(\alpha+\alpha^{-1})$, alors le morphisme de $\mathcal{A}_{H_4}$ dans $\mathcal{H}_{H_4,\alpha}^\star\simeq \underset{\rho~ \mbox{ irr}}\prod GL_{n_\rho}(\F_p(\sqrt{\alpha},\xi+\xi^{-1}))$ se factorise à travers le morphisme surjectif
$$\Phi :\rightarrow SL_4(q^2)\times \Omega_6^+(q^2)\times Spin_8^+(q)\times SL_9(q^2)\times \Omega_{10}^+(q)\times SL_{16}(q)^2 \times \Omega_{16}^+(q)\times \Omega_{18}^+(q)$$
$$\times \Omega_{24}^+(q)^2\times SL_{25}(q)\times \Omega_{30}^+(q)\times SL_{36}(q)\times \Omega_{40}^+(q)\times \Omega_{48}^+(q^{\frac{1}{2}}).$$
\item Si $\F_q=\F_p(\sqrt{\alpha})=\F_p(\alpha)\neq \F_p(\alpha+\alpha^{-1})$ et $\Phi_{1,2}(\sqrt{\alpha})=\sqrt{\alpha}^{-1}$, alors le morphisme de $\mathcal{A}_{H_4}$ dans $\mathcal{H}_{H_4,\alpha}^\star\simeq \underset{\rho~ \mbox{ irr}}\prod GL_{n_\rho}(\F_p(\sqrt{\alpha},\xi+\xi^{-1}))$ se factorise à travers le morphisme surjectif
$$\Phi :\rightarrow SL_4(q)\times \Omega_6^+(q)\times Spin_8^+(q^{\frac{1}{2}})\times SL_9(q^2)\times \Omega_{10}^+(q^{\frac{1}{2}})\times SU_{16}(q^{\frac{1}{2}})^2 \times \Omega_{16}^+(q^{\frac{1}{2}})\times \Omega_{18}^+(q^{\frac{1}{2}})$$
$$\times \Omega_{24}^+(q^{\frac{1}{2}})^2\times SU_{25}(q^{\frac{1}{2}})\times \Omega_{30}^+(q^{\frac{1}{2}})\times SU_{36}(q^{\frac{1}{2}})\times \Omega_{40}^+(q^{\frac{1}{2}})\times \Omega_{48}^+(q^{\frac{1}{2}}).$$
\item Si $\F_q=\F_p(\sqrt{\alpha})=\F_p(\alpha)\neq \F_p(\alpha+\alpha^{-1})$ et $\Phi_{1,2}(\sqrt{\alpha})=-\sqrt{\alpha}^{-1}$, alors le morphisme de $\mathcal{A}_{H_4}$ dans $\mathcal{H}_{H_4,\alpha}^\star\simeq \underset{\rho~ \mbox{ irr}}\prod GL_{n_\rho}(\F_p(\sqrt{\alpha},\xi+\xi^{-1}))$ se factorise à travers le morphisme surjectif
$$\Phi :\rightarrow SL_4(q)\times \Omega_6^+(q)\times Spin_8^+(q^{\frac{1}{2}})\times SL_9(q^2)\times \Omega_{10}^+(q^{\frac{1}{2}})\times SL_{16}(q) \times \Omega_{16}^+(q^{\frac{1}{2}})\times \Omega_{18}^+(q^{\frac{1}{2}})$$
$$\times \Omega_{24}^+(q^{\frac{1}{2}})^2\times SU_{25}(q^{\frac{1}{2}})\times \Omega_{30}^+(q^{\frac{1}{2}})\times SU_{36}(q^{\frac{1}{2}})\times \Omega_{40}^+(q^{\frac{1}{2}})\times \Omega_{48}^+(q^{\frac{1}{2}}).$$
\item Si $\F_{q^2}=\F_p(\sqrt{\alpha}) \neq \F_p(\alpha)=\F_p(\alpha+\alpha^{-1})$, alors le morphisme de $\mathcal{A}_{H_4}$ dans $\mathcal{H}_{H_4,\alpha}^\star\simeq \underset{\rho~ \mbox{ irr}}\prod GL_{n_\rho}(\F_p(\sqrt{\alpha},\xi+\xi^{-1}))$ se factorise à travers le morphisme surjectif
$$\Phi :\rightarrow SL_4(q^2)\times \Omega_6^+(q^2)\times Spin_8^+(q)\times SL_9(q^2)\times \Omega_{10}^+(q)\times SL_{16}(q^2) \times \Omega_{16}^+(q)\times \Omega_{18}^+(q)$$
$$\times \Omega_{24}^+(q)^2\times SL_{25}(q)\times \Omega_{30}^+(q)\times SL_{36}(q)\times \Omega_{40}^+(q)\times \Omega_{48}^+(q^{\frac{1}{2}}).$$
\end{enumerate}
\item Supposons $1\nsim 2$. Lorsqu'il existe, on note $\epsilon$ l'automorphisme d'ordre $2$ de $\F_q$.
\begin{enumerate}
\item Si $\F_q=\F_p(\sqrt{\alpha})=\F_p(\alpha)\neq \F_p(\alpha+\alpha^{-1})$ et $\epsilon(\sqrt{\alpha})=\sqrt{\alpha}^{-1}$, alors le morphisme de $\mathcal{A}_{H_4}$ dans $\mathcal{H}_{H_4,\alpha}^\star\simeq \underset{\rho~ \mbox{ irr}}\prod GL_{n_\rho}(\F_p(\sqrt{\alpha},\xi+\xi^{-1}))$ se factorise à travers le morphisme surjectif $$\Phi :\rightarrow SU_4(q^{\frac{1}{2}})^2\times \Omega_6^+(q^{\frac{1}{2}})^2\times Spin_8^+(q^{\frac{1}{2}})\times SU_9(q^{\frac{1}{2}})^2\times \Omega_{10}^+(q^{\frac{1}{2}})\times SU_{16}(q^{\frac{1}{2}})^2 \times \Omega_{16}^+(q^{\frac{1}{2}})^2$$
$$\times \Omega_{18}^+(q^{\frac{1}{2}})\times \Omega_{24}^+(q^{\frac{1}{2}})^4\times SU_{25}(q^{\frac{1}{2}})\times \Omega_{30}^+(q^{\frac{1}{2}})^2\times SU_{36}(q^{\frac{1}{2}})\times \Omega_{40}^+(q^{\frac{1}{2}})\times \Omega_{48}^+(q^{\frac{1}{2}}).$$
\item Si $\F_q=\F_p(\sqrt{\alpha})=\F_p(\alpha)\neq \F_p(\alpha+\alpha^{-1})$ et $\epsilon(\sqrt{\alpha})=-\sqrt{\alpha}^{-1}$, alors le morphisme de $\mathcal{A}_{H_4}$ dans $\mathcal{H}_{H_4,\alpha}^\star\simeq \underset{\rho~ \mbox{ irr}}\prod GL_{n_\rho}(\F_p(\sqrt{\alpha},\xi+\xi^{-1}))$ se factorise à travers le morphisme surjectif
$$\Phi :\rightarrow SU_4(q^{\frac{1}{2}})^2\times \Omega_6^+(q^{\frac{1}{2}})^2\times Spin_8^+(q^{\frac{1}{2}})\times SU_9(q^{\frac{1}{2}})^2\times \Omega_{10}^+(q^{\frac{1}{2}})\times SL_{16}(q) \times \Omega_{16}^+(q^{\frac{1}{2}})^2$$
$$\times \Omega_{18}^+(q^{\frac{1}{2}})\times \Omega_{24}^+(q^{\frac{1}{2}})^4\times SU_{25}(q^{\frac{1}{2}})\times \Omega_{30}^+(q^{\frac{1}{2}})^2\times SU_{36}(q^{\frac{1}{2}})\times \Omega_{40}^+(q^{\frac{1}{2}})\times \Omega_{48}^+(q^{\frac{1}{2}}).$$
\end{enumerate}
\end{enumerate}
\end{theor}

\bigskip

Dans le chapitre \ref{TypeF4}, nous déterminons l'image de $\mathcal{A}_{F_4}$ dans son algèbre de Iwahori-Hecke associée. Il y a deux injections naturelles de $A_{B_3}$ dans $A_{F_4}$ donc nous utiliserons les résultats du chapitre \ref{TypeB} pour démontrer les résultats. Il y a deux classes de conjugaison pour les générateurs du groupe de Coxeter de type $F_4$. L'algèbre de Iwahori-Hecke associée dépend alors de deux paramètres $\alpha$ et $\beta$. Comme les représentations sont données par des $F_4$-graphes, les modèles sont définis sur $\F_p(\sqrt{\alpha},\sqrt{\beta})$. Il faut alors considérer l'extension de corps $\F_p(\sqrt{\alpha},\sqrt{\beta})$ au-dessus de $\F_p(\alpha+\alpha^{-1},\beta+\beta^{-1})$. Une liste exhaustive des tours d'extensions possible permet de démontrer que cette extension est de degré au plus $2$. Toutes les extensions sont décrites dans la section \ref{sectionF4graphs}. Les $F_4$-graphes de \cite{G-P} sont aussi rappelés dans cette section, ils sont présentés de manière à faire apparaître les symétries pour les $F_4$-graphes associés à des graphes auto-duaux. Nous n'avons pas démontré la conjecture dans ce cas car les propriétés d'unicité démontrés dans les propositions \ref{Unique1} et \ref{Unique2} ne sont démontrées que dans le cas des paramètres égaux. Les représentations étant toutes de dimension associée inférieure à 16, on peut déterminer par le calcul la forme bilinéaire associée sans émettre d'hypothèses supplémentaires. Le groupe $\mathcal{A}_{F_4}$ n'est pas parfait et cela complique considérablement les preuves dans certains cas. De plus, les restrictions aux sous-groupes paraboliques n'apportent pas autant d'informations que dans les chapitres précédents car la clôture normale des sous-groupes paraboliques isomorphes à $A_{B_3}$ est différente de $\mathcal{A}_{F_4}$. Cela nécessite de nouveaux arguments dans ce chapitre avec des preuves parfois plus calculatoires. Le théorème suivant donne l'image de $\mathcal{A}_{F_4}$ dans son algèbre de Iwahori-Hecke associée suivant les possibles tours d'extensions de corps considérées et conclut le chapitre \ref{TypeF4}.

 \begin{theor}
On note $\F_{\tilde{q}}=\F_p(\sqrt{\alpha},\sqrt{\beta})$, $\F_{r_{\alpha}}=\F_p(\alpha+\alpha^{-1})$ et $\F_{r_{\beta}}=\F_p(\beta+\beta^{-1})$.

Dans les cas $1$, $4$, $5$ et $10$, le morphisme de $\mathcal{A}_{F_4}$ dans $\mathcal{H}_{F_4,\alpha,\beta}^\star\simeq \underset{\rho \mbox{ irr}}\prod GL_{n_\rho}(\tilde{q})$ se factorise à travers le morphisme surjectif
$$\Phi: \mathcal{A}_{F_4} \rightarrow (SL_2(r_{\alpha})\circ SL_2(r_{\beta}))\times SL_4(q)^2\times \Omega_6^+(q)^2$$
$$\times SL_8(q)^2\times SL_9(q)^2\times \Omega_{12}^+(q)\times \Omega_{16}^+(q).$$
Dans les cas $11$, $12$, $13$ et $16$, le morphisme de $\mathcal{A}_{F_4}$ dans $\mathcal{H}_{F_4,\alpha,\beta}^\star\simeq \underset{\rho \mbox{ irr}}\prod GL_{n_\rho}(\tilde{q})$ se factorise à travers le morphisme surjectif
$$\Phi: \mathcal{A}_{F_4} \rightarrow  \times (SL_2(r_{\alpha})\circ SL_2(r_{\beta}))\times SU_4(q^{\frac{1}{2}})^2\times \times \Omega_6^+(q^{\frac{1}{2}})^2$$
$$\times SU_8(q^{\frac{1}{2}})^2\times SU_9(q^{\frac{1}{2}})^2\times \Omega_{12}^+(q^{\frac{1}{2}})\times \Omega_{16}^+(q^{\frac{1}{2}}).$$
Dans les cas $2$, $6$, $9$ et $15$, le morphisme de $\mathcal{A}_{F_4}$ dans $\mathcal{H}_{F_4,\alpha,\beta}^\star\simeq \underset{\rho \mbox{ irr}}\prod GL_{n_\rho}(\tilde{q})$ se factorise à travers le morphisme surjectif
$$\Phi: \mathcal{A}_{F_4} \rightarrow  (SL_2(r_{\alpha})\circ SL_2(r_{\beta}))\times SL_4(q)\times \times \Omega_6^+(q)$$
$$\times SL_8(q^{\frac{1}{2}})\times SU_8(q^{\frac{1}{2}})\times SL_9(q)\times \Omega_{12}^+(q^{\frac{1}{2}})\times \Omega_{16}^+(q^{\frac{1}{2}}).$$
Dans les cas $3$, $7$, $8$ et $14$, le morphisme de $\mathcal{A}_{F_4}$ dans $\mathcal{H}_{F_4,\alpha,\beta}^\star\simeq \underset{\rho \mbox{ irr}}\prod GL_{n_\rho}(\tilde{q})$ se factorise à travers le morphisme surjectif
$$\Phi: \mathcal{A}_{F_4} \rightarrow  (SL_2(r_{\alpha})\circ SL_2(r_{\beta}))\times SL_4(q) \times \Omega_6^+(q)$$
$$\times SL_8(q)\times SU_8(q^{\frac{1}{2}})\times SL_9(q)\times \Omega_{12}^+(q^{\frac{1}{2}})\times \Omega_{16}^+(q^{\frac{1}{2}}).$$
 \end{theor}
 
 \bigskip
 
Pour finir, dans la section \ref{ErratBMBMM} de l'appendice, nous donnons un erratum des articles \cite{BM} et \cite{BMM}. Les résultats de ces articles sont utilisés à diverses reprises dans cette thèse, notamment dans le chapitre \ref{TypeB}. Les résultats de ces articles sont correctes mais certaines preuves sont imprécises ou contiennent des erreurs qui sont listés dans cette section.
 
\chapter{Introduction} 

\section{General Introduction}

In this doctoral thesis, we will determine the image of Artin groups associated to all finite irreducible Coxeter groups inside their associated finite Iwahori-Hecke algebra. This was done in type $A$ in \cite{BM} and \cite{BMM}. The Zariski closure of the image was determined in the generic case in \cite{IH2}. It is suggested by strong approximation that the results should be similar in the finite case. However, the conditions required to use are much too strong and would only provide a portion of the results. We show in this thesis that they are but that new phenomena arise from the different field factorizations. The techniques used in the finite case are very different from the ones in the generic case. The main arguments come from finite group theory. In high dimension, we will use a theorem by Guralnick-Saxl \cite{GS} which uses the classification of finite simple groups to give a condition for subgroups of linear groups to be classical groups in a natural representation. In low dimension, we will mainly use the classification of maximal subgroups of classical groups in \cite{BHRC} for the complicated cases.

Our results about the image of Artin groups inside the finite Hecke algebras may have various applications. For instance, finite classical groups and direct products of finite classical groups appear in this way as finite quotients of the Artin groups. Since the latter are fundamental groups of algebraic varieties, this also defines étale coverings of these varieties. Since these varieties are defined over $\Q$, this may have applications to the inverse Galois problem (see for example \cite{SV} or \cite{FriedVolklein} for applications in type $A_n$).

 They also provide new information on the finite classical groups because we get explicit generators verifying the braid relations for those groups. This can provide interesting constructions of these groups and some of their subgroups by looking at restrictions to parabolic subgroups of the Artin groups. In particular, we get a somewhat unexpected description of the $8$-dimensional spin group from the two $8$-dimensional representations in type $H_4$ (see section \ref{sectiontriality}).
 
 We find some new $W$-graphs in types $H_4$, $E_6$ and $E_8$ which provide different information from the usual ones. They are all associated in a natural way to a bilinear form which is very complicated to obtain in the previous models. In this model, the bilinear form is obtained using only the two-colorability and its matrix in a well chosen basis is anti-diagonal. The uniqueness properties can probably be extended in a more general setting and understanding which setting this is may be worth considering.
 
 \section{Outline of the thesis}

We will begin by giving an introduction of Coxeter groups, Artin groups and Iwahori-Hecke algebras. We then recall some properties of symmetric algebras and show that Iwahori-Hecke algebras are symmetric. In this section, we show a specialization result which we will use throughout the thesis to show that under certain conditions, we can specialize the models for Iwahori-Hecke algebras in the generic case to the finite case. We will then give Aschbacher's theorem \cite{ASCHMAXSUBGRPS} on maximal subgroups of classical groups and recall the different classes $\mathcal{C}_i$ and $\mathcal{S}$ which are defined in this theorem.

We will then start by determining the image of Artin groups of classical types inside their associated finite Iwahori-Hecke algebras. The matrix models in those types are given by double-partitions of an integer $n$. The general idea of the proof is then similar to the proof in type $A$. The image is first determined for small $n$. The branching rule is then used to give an inductive proof on $n$ to determine the image in the general case. It is first necessary to determine which representations factor through each other via field automorphisms or the transposed inverse automorphism. The second parameter in type $B$ will give rise to new factorizations which did not appear in the generic case. We will then have to separate the study in six different cases for the field extensions which will give different results. We get that in type $B$, if the fields extensions occuring are all trivial then the image of the representations associated to partitions of $n$ are special linear groups defined over $\F_q$ if $\lambda\neq \lambda'$ and symplectic or orthogonal groups if $\lambda=\lambda'$. In the case when the field extensions are trivial, the only factorizations appearing are between the representations labeled by hook partitions and between the representation labeled by a partition and the representation labeled by its transposed partition. This is summarized in Theorem \ref{result1}. When the field extensions a more complex, we get more factorizations and we get both unitary groups or special linear groups depending on the partition we are considering. The result for those cases are then given in Theorem \ref{result2} to Theorem \ref{result6}.

In type $D$, the matrix models are similar but there is only one conjugacy class for the generators therefore the field extensions are less complicated. However, there are additional factorizations appearing which make the branching rule more complex. The result is then similar to the one in type $B$ except there are more representations affording groups preserving non-degenerate bilinear forms. The results for this chapter are given in Theorem \ref{result1D} and Theorem \ref{result2D}.

Those results will be useful for the groups of exceptional types because of the natural inclusions of Coxeter groups of classical types inside the exceptional types. Before treating the groups of type $H$, $E$ and $F$, we will determine the image for dihedral groups and general results and $W$-graphs which will give us the matrix models for those exceptional types.

All the representations are $1$-dimensional or $2$-dimensional in dihedral type therefore determining the images inside each given representation will not be too difficult. There is a theorem by Dickson (see \cite{HUP}, Theorem $8.27$) classifying the subgroups of $PSL_2(q)$ which we will use to determine the images in those cases and for $2$-dimensional representations in all types. The difficulties will arise from the field extensions which depend on primitive roots of unity. This will lead us to define an equivalence relation between integers. The results for $I_2(5)$ will then be useful in type $H$ since there is a natural inclusion from $I_2(5)$ to $H_3$. Since $SU_2(q^{\frac{1}{2}})$ is conjugate to $SL_2(q^{\frac{1}{2}})$ inside $GL_2(q)$, we will have that the image inside every given representation is conjugate to $SL_2(q^{\frac{1}{2}})$ or $SL_2(q)$. The image inside the full Iwahori-Hecke is then given in Theorem \ref{resdihedral} and Theorem \ref{resdihedraleven}.

The matrix models for exceptional types are given by $W$-graphs. We will recall some general properties for those graphs. We will then prove some uniqueness properties and establish a conjecture for $W$-graphs associated to self-dual representations (see Definition \ref{defselfdualrepresentation} in Conjecture \ref{conjecturewgraphs}. We proved the conjecture in types $I$, $H$ and $E$ by computation using the HPC platform MatriCS \cite{Matrics}. We assume the conjecture is true in order to find non-degenerate bilinear forms preserved by the image of derived subgroups of Artin groups  under self-dual representations. It is clear that such a bilinear form exists for self-dual representations, however it is hard to determine the type of this bilinear form. The conjecture will allow us to determine explicitly those bilinear forms and to give new $W$-graphs, where a condition on the $2$-colorability of the graph is sufficient to determine the type of the form. We draw them in type $H_4$ and some of them in type $E_6$. The remaining ones can be downloaded from \cite{newgraphsEsterle}.

We then determine the image in type $H$ using the results in type $I_2(5)$. In type $H_3$, we encounter some $W$-graphs which are not $2$-colorable. It is then necessary to consider more complicated field extensions. By simple computations, we can treat those more complicated extensions and get the image inside the corresponding representations. As in the dihedral type, we need to use the equivalence relation on integers defined in Lemma \ref{Isomorphism}. In this case, we only consider fifth roots of unity therefore the equivalence relation is only on $1$ and $2$. The factorizations between the different representations will depend on whether $1\sim 2$ or $1\nsim 2$. All the groups considered in type $H_3$ are linear or unitary depending on the field extensions. In type $H_4$, the dimensions of the representations go up to $48$. There are many self-dual representations for which we have verified the conjecture and found corresponding self-dual $H_4$-graphs given in the Appendix. The study of this type is quite complicated so we separate the study in four parts. First, in section \ref{H4gen}, we prove some general properties appearing in type $H_4$ as in the other types. Then, in section \ref{H4lowdim}, we begin by considering the low-dimensional representations with the representations of dimension going up to $8$. In these cases, we cannot using Theorem \ref{CGFS} therefore we need to use some more elementary theorems and some computational proofs. The image in the $8$-dimensional representation is isomorphic to $\Omega_8^+(q)$ or $\Omega_8^+(q^{\frac{1}{2}})$. It is very difficult to exclude the case $2^{\cdot}\Omega_7(q)$. In order to do it, we provide a very computational proof using some elements in the normalizer of the image of the derived subgroup of $<S_2,S_3>$. Once this result is proven, we get a nice description of the $\op{Spin}_8^+(q)$ groups appearing via the two $8$-dimensional representations. We explain in section \ref{sectiontriality} how the triality automorphism appears between the projections of the two $8$-dimensional representations. We then prove that the image of the derived subgroup $\mathcal{A}_{H_4}$ inside the product of the two $8$-dimensional representations is the universal central extension of $P\Omega_8^+(q)$ which proves it is isomorphic to $\op{Spin}_8^+(q)$. In section \ref{H4highdim}, we determine the image of $\mathcal{A}_{H_4}$ inside the representations of dimension greater than $8$ using Theorem \ref{CGFS}. There are four $16$-dimensional representations which are not $2$-colorable which makes the study slightly more complicated but using the usual arguments, we get the image for all groups except for the image inside the $48$-dimensional representation. We did not succeed in proving that the image is what we expect for the $48$-dimensional representation but we conjecture it is. We then use Goursat's Lemma to recover the image inside the full Iwahori-Hecke algebra given in Theorem \ref{resH4}.

The image in type $E$ is determined using the results in type $D_5$. The main difficulty in this type was to determine the types of the bilinear forms. The high-dimensional representations required long computations using Matrics \cite{Matrics}. Those computations were done assuming the conjecture was true and we then proved the conjecture is true. Another difficulty arose in the proof that the specialization to finite fields still gave a split semisimple algebra. Some of the weights of the edges on the $W$-graphs considered vanish and it is then necessary to prove that the graphs remain connected once those edges vanish. Once we know if the forms are symplectic or orthogonal and the algebra is split semisimple, the usual arguments cover most of the proofs in those types. There are two representations in type $E_7$ which are not $2$-colorable, we compute the traces of well-chosen elements in order to treat those cases. The main results in this chapter are given in Theorem \ref{resultE6}, Theorem \ref{resultE7} and Theorem \ref{resultE8}.

In type $F_4$, there are two parameters and the representations are defined over $\F_p(\sqrt{\alpha},\sqrt{\beta})$. This makes the field extensions much more complex. We also have that the normal closure of parabolic subgroups does verify the properties which we can observe and use in other types. The derived subgroup $\mathcal{A}_{F_4}$ is not perfect. This complicates the use of Goursat's lemma to determine the image inside the full Iwahori-Hecke algebra. The uniqueness properties for $W$-graphs are only proven for one parameter. It seems they do not hold in type $F_4$ and the $2$-colorability of $F_4$-graphs may not be stable by isomorphism of representations. We did not prove the conjecture in this type and it may not hold. Nevertheless, the representations are of dimension at most $16$, therefore we can determine by hand the non-degenerate bilinear forms preserved by the groups. We still get the usual results depending on the field extensions except for the $4$-dimensional representation $4_1$ which is the central product of the $2$-dimensional representations $2_1$ and $2_3$. The main result in this section is given in Theorem \ref{resultF4}.

\section{Perspectives}

Conjecture \ref{conjecturewgraphs} does not need to be proved for classical types in order to conclude the doctoral thesis. After the Ph.D. is finished, one area to explore would be to prove Conjecture \ref{conjecturewgraphs} for types $A$, $B$ and $D$. A description of $W$-graphs affording the representations given by the combinatorial model associated to partitions in types $A$ and $B$ has been given by Naruse in \cite{NAR1} and \cite{NAR2}. Using those descriptions and trying to extend his work to type $D$, we hope to be able to prove the conjectures for types $A$, $B$ and $D$.

\medskip

The notion of Hecke algebras has been generalized to complex reflection groups by M.Broué, G. Malle and R.Rouquier \cite{BMR}. It has been conjectured that there exists a symmetrizing trace for all those algebras. If that conjecture is true then it will be possible to use Schur elements to determine when those algebras are semisimple and specialize those algebras to finite fields. In type $G(d,1,n)$, the representations are known and are labeled by $d$-partitions of unity and we know the Schur elements. It is then possible to study the image of the associated braid group inside those representations and this work seems like a natural extension of the work done in the Ph.D. The specialization to finite fields could become complicated due to the high number of parameters involved.

\medskip

It is still an open question whether Artin groups are linear or not. A linear group is a group for which there exists a monomorphism into some group of invertible matrices. It is also an open question whether the morphisms from the Artin groups onto the Iwahori-Hecke algebras are one-to-one in the generic case. If we could determine the kernels of those morphism in the finite case, it might be possible to find elements of the kernel in the generic case and show those morphisms are not monomorphisms. We have found some interesting elements in the kernels of some morphisms composed with projections upon irreducible representations in type $H_4$ but we have not yet had the opportunity to explore that question more.

\medskip

In order to determine the image of the Artin group inside finite Iwahori-Hecke algebras, assumptions on the parameters were made for it to be semisimple. One could consider the cases where those parameters are not semisimple and still study the same question. The Hoefsmit models can only be specialized under some assumptions for the parameters which excludes most of the semisimple cases. On the other hand, $W$-graphs afford well-defined representations for nearly all possible representations.

\medskip

The image of the Artin group inside each representation gives us a finite quotient of the Artin group. Using rigidity techniques which can be found in \cite{MM} and \cite{V}, it should be possible to get some results in inverse Galois theory using methods as in \cite{SV}. The key fact in this problem is that Artin groups can be seen as fundamental groups of certain algebraic varieties defined over $\Q$.

\section{Notations}

We will use the following notations throughout the thesis 

$\begin{array}{lcl}
\N^\star & : & \mbox{positive integers}\\
\Z & : & \mbox{ring of rational integers}\\
\Q & : & \mbox{rational field}\\
\mathfrak{S}_n & : & \mbox{symmetric group}\\
\mathfrak{A}_n & : & \mbox{alternating group}\\
\lambda\vdash n & : & \lambda\mbox{ is a partition of n}\\
\lambda \Vdash n & : & \lambda\mbox{ is a double-partition of n}\\
\T\in \lambda & : & \T \mbox{ is a tableau (resp double-tableau) associated to the partition}\\
& &  \mbox{(resp double-partition) } \lambda\\
$[\![a,b]\!]$ & : & \mbox{ set of integers in the interval }[a,b]\\
\F_q & : & \mbox{finite field with q elements}\\
M_n(q) & : & \mbox{algebra of matrices over the field }\F_q\\
\op{GL}_n(q) & : & \mbox{group of invertible matrices over the field }\F_q\\
\op{SL}_n(q) & : & \mbox{group of invertible matrices of determinant 1 over } \F_q\\
\op{SU}_n(q) & : & \mbox{group of unitary matrices over } \F_{q^2}\\
\op{SP}_n(q) & : & \mbox{group of symplectic matrices over } \F_q\\
\Omega_n^{\epsilon}(q) & : & \mbox{kernel of the spinor norm of the orthogonal group of type } \epsilon\\
A_{W_n} & : & \mbox{Artin group of type }W_n \\
\mathcal{A}_{W_n} & : & \mbox{derived subgroup of the Artin group } A_{W_n}\\
Z(G) & : & \mbox{center of a group } G\\
$[G,G]$ & : & \mbox{derived subgroup of a group } G\\
(a,b) & : & \mbox{gcd of a pair of integers (a,b)}\\
I_n & : & \mbox{identity matrix of size n x n}\\
E_{i,j} & : & \mbox{elementary matrix with non-zero coefficient in position (i,j)}\\
\lambda' & : & \mbox{transposed partition (or double-partition) of } \lambda\\
A : B & : & \mbox{split normal extension of A by B}\\
A^{\cdot} B & :  & \mbox{normal extension of A by B which is not split}\\
A.B & :  & \mbox{normal extension of A by B which may be split}\\
\diag(a_1,a_2,\dots,a_n) & : & \mbox{diagonal matrix with coefficients } a_1, a_2, \dots, a_n\\
\Phi_n & : & \mbox{n-th cyclotomic polynomial}

\end{array}$

\newpage

$ $

\newpage

\chapter{Preliminaries}

\section{Definitions and first properties of Coxeter groups, Artin groups and Iwahori-Hecke algebras}\label{DefCoxArtHec}

We here recall the definitions of Coxeter groups, Artin groups and Iwahori-Hecke algebras and give the classification of finite irreducible Coxeter groups. Those are the main objects we will use throughout the thesis. They appear in many different fields, they can be seen as groups generated by involutions verifying certain relations. They can also be considered as real reflection groups as a subclass of the complex reflection groups.

\begin{Def}
Let $W$ be a group with a generating set $S$ of elements of order $2$. For $s,s'$ in $S$, we write $m_{s,s'}$ the order of $ss'$. We say that $(W,S)$ is a Coxeter system if

 $<S | \forall s\in S,\forall s'\in S\setminus{\{s\}},s^2=1, ~ (ss')^{m_{ss'}}=1>$ is a presentation of $W$.

\smallskip

The Dynkin diagram associated to a Coxeter system $(W,S)$ is the graph with elements of $S$ as vertices. There are $m_{s,t}-2$ edges between $s$ and $t$ if $m_{s,t} \leq 4$ and one edge with weight $m_{s,t}$ otherwise.

\smallskip

A Coxeter group is said to be irreducible if there exists no non-empty disjoint subsets $S_1$ and $S_2$ of $S$ such that $S=S_1\cup S_2$ and for any pair $(s_1,s_2)\in S_1\times S_2$, we have $s_1s_2=s_2s_1$.
\end{Def}

\begin{figure}
\begin{tikzpicture}
[place/.style={circle,draw=black,
inner sep=1pt,minimum size=10mm}]
\node (1) at (0,0)[place]{$s_1$};
\node (2) at (2,0)[place]{$s_2$};
\draw (1) to (2);
\node (3) at (5.1,0)[place]{$s_n$};
\draw[dashed] (2) to (3);
\draw (-0.7,-0.7) node{$A_n$};

\node (1) at (0,-2)[place]{$t$};
\node (2) at (2,-2)[place]{$s_1$};
\draw (0.500,-2.1) -- (1.500,-2.1);
\draw (0.500,-1.9) -- (1.500,-1.9);
\node (3) at (5.1,-2)[place]{$s_{n-1}$};
\draw[dashed] (2) to (3);
\draw (-0.7,-2.7) node{$B_n$};

\node (1) at (0,-4)[place]{$s_1$};
\node (2) at (0,-6)[place]{$s_2$};
\node (3) at (2,-5)[place]{$s_3$};
\node (4) at (5.1,-5)[place]{$s_{n}$};
\draw (1) to (3);
\draw (2) to (3);
\draw [dashed] (3) to (4);
\draw (-0.7,-6.7)node{$D_n$};

\node (1) at (0,-8)[place]{$s_1$};
\node (2) at (2,-8)[place]{$s_2$};
\draw (1) to node[auto]{$n$} (2);
\draw (-0.7,-8.9)node {$I_2(n), n\geq 5$};

\node (1) at (7,0)[place]{$s_1$};
\node (2) at (9,0)[place]{$s_2$};
\node (3) at (11,0)[place]{$s_3$};
\draw (1) to node[auto]{$5$} (2);
\draw (2) to (3);
\draw (6.3,-0.7)node{$H_3$};

\node (1) at (7,-2)[place]{$s_1$};
\node (2) at (9,-2)[place]{$s_2$};
\node (3) at (11,-2)[place]{$s_3$};
\node (4) at (13,-2)[place]{$s_4$};
\draw (1) to node[auto]{$5$} (2);
\draw (2) to (3);
\draw (3) to (4);
\draw (6.3,-2.7)node{$H_4$};

\node (1) at (7,-6)[place]{$s_1$};
\node (3) at (9,-6)[place]{$s_3$};
\node (4) at (11,-6)[place]{$s_4$};
\node (2) at (11,-4)[place]{$s_2$};
\node (5) at (13,-6)[place]{$s_5$};
\node (6) at (15,-6)[place]{$s_n$};
\draw (1) to (3);
\draw (3) to (4);
\draw (2) to (4);
\draw (4) to (5);
\draw[dashed] (5) to (6);
\draw (6.3,-6.9) node{$E_n,n\in \{6,7,8\}$};

\node (1) at (7,-8)[place]{$s_1$};
\node (2) at (9,-8)[place]{$s_2$};
\node (3) at (11,-8)[place]{$s_3$};
\node (4) at (13,-8)[place]{$s_4$};
\draw (1) to (2);
\draw (9.5,-8.1) -- (10.5,-8.1);
\draw (9.5,-7.9) -- (10.5,-7.9);
\draw (3) to (4);
\draw (6.3,-8.7)node{$F_4$};
\end{tikzpicture}
\caption{Classification of finite irreducible Coxeter groups}\label{classifCox}
\end{figure}
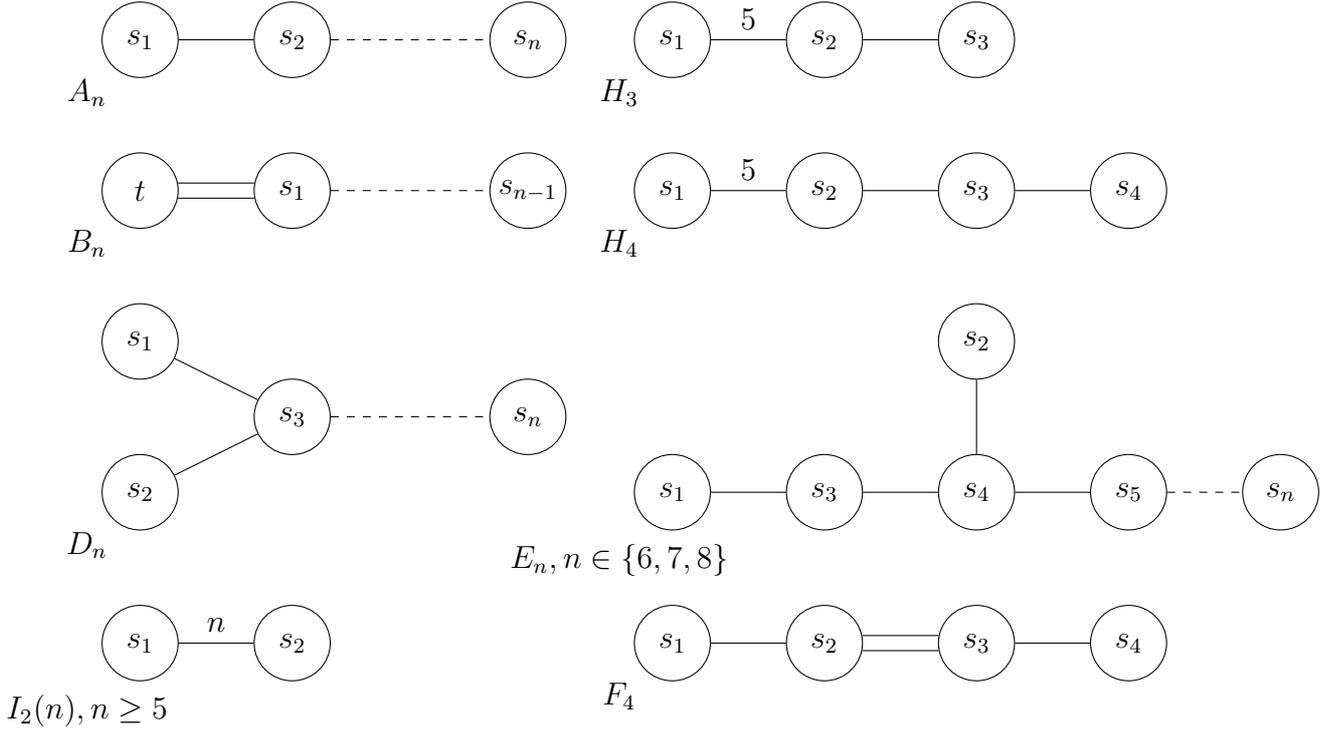

The irreducible finite Coxeter groups have been classified by Coxeter \cite{Cox} and are given in Figure \ref{classifCox}

Note that the Coxeter group of type $A_n$ is the symmetric group $\mathfrak{S}_{n+1}$ generated by the transpositions $s_i=(i ~ i+1)$. The Coxeter group of type $B_n$ is the group of symmetries of the hypercube of dimension $n$. The Coxeter group of type $I_2(n)$ is the dihedral group with $2n$ elements. 

There are natural inclusions between the groups $W_n$ and $W_{n+1}$ and between $A_{n-1}$ and $B_n$, $A_{n-1}$ and $D_n$, $I_2(5)$ and $H_3$, $A_3$ and $H_4$, $D_5$ and $E_6$, $B_3$ and $F_4$. We will use those inclusions which remain in the Artin group context which we will define below.

\begin{Def}
Let $(W,S)$ be a finite Coxeter system. To the Coxeter group $W$, we associate the following Artin group $A_W$, where the order $2$ condition has been removed 

 $$<S,\forall (s,t)\in S^2, \underset{m_{s,t}}{\underbrace{sts...}} =\underset{m_{s,t}}{\underbrace{tst...}}>$$
\end{Def}
The Artin group $A_{A_n}$ associated to the $A_n$ type Coxeter group admits the following presentation : $A_{A_n} = <(s_i)_{i\in [\![1,n]\!]}, s_is_{i+1}s_i = s_{i+1}s_is_{i+1}, \vert i-j\vert \geq 2, s_is_j=s_js_i>$.

The Artin group $A_{B_n}$ associated to the $B_n$ type Coxeter group admits the following presentation : $A_{B_n} = <t,(s_i)_{i\in [\![1,n-1]\!]}, ts_1ts_1 = s_1ts_1t, s_is_{i+1}s_i = s_{i+1}s_is_{i+1}, \vert i-j\vert \geq 2, s_is_j=s_js_i>$.

\begin{Def}
Let $A$ be an Artin group with generators $\{s_i\}_{i\in [\![1,n]\!]}$. Let $a \in A$, two expression $a=s_{i_k}..s_{i_1}$ and $a=s_{j_k}...s_{j_1}$ of $a$ are said to be equivalent if one can be deduced from the other using the braid relations.
\end{Def}

The following result is a fundamental result in the study of Iwahori-Hecke algebras, it can be found for example in \cite{G-P} (Theorem $1.2.2$).

\begin{theo}{(Matsumoto's Theorem)}\label{Matsumoto}\\
Let $(W,S)$ be a Coxeter system with $W$ a finite irreducible Coxeter group. We let $A_W^+$ be its braid monoid (the monoid of positive words in the generators of $A_W$). A reduced expression $\sigma=s_{i_k}...s_{i_1}$ of an element $\sigma$ in the braid monoid $A_W^+$ is an expression where $k$ is the minimal number of generators necessary to write the element $\sigma$ in the braid monoid.

$s_{i_k}...s_{i_1}$ and $s_{j_k}...s_{j_1}$ are reduced expressions of the same element in an Artin monoid $A_W^+$ if and only if they are equivalent.
\end{theo}

We can now define the Iwahori-Hecke algebra associated to a given Coxeter group. The Matsumoto theorem will enable us to see it is a free algebra. We define it here in a general setting, we will give more precise definitions in the finite case in each type in the following chapters.

\begin{Def}
Let $(W,S)$ be a Coxeter system, $\tilde{W}$ the associated Artin group, $R$ a ring and $(u_s)_{s\in S}$ indeterminates such that $\alpha_s=\alpha_t$ if $s$ and $t$ are conjugate in $W$.\\
The $R[u_s^{\pm 1}]$-Iwahori-Hecke Algebra $\mathcal{H}_{W,R,(u_s)_{s\in S}}$ associated to $W$ is given by the following presentation 
$$\mathcal{H}_{W,R,(u_s)_{s\in S}}=<T_1,...,T_n| \underset{m_{s_i,s_j}}{\underbrace{T_iT_jT_i...}} =\underset{m_{s_i,s_j}}{\underbrace{T_jT_iT_j...}}, (T_i -\alpha_{s_i})(T_i+1)=0>$$
\end{Def}

By Theorem \ref{Matsumoto}, we can define the element $T_{\sigma}=T_{i_1}\dots T_{i_r}$ for $\sigma=s_{i_1}\dots s_{i_r}$ in a reduced expression. The following proposition is then fairly easy to show 

\begin{prop}
The Iwahori-Hecke algebra $\mathcal{H}_{W,R,(u_s)_{s\in S}}$ is a free $R[u_s^{\pm 1}]$-module of rank $\vert W\vert$.
\end{prop}

If we consider the case $R=\Z$, we have that $\mathcal{H}_{W,\Z,(u_s)_{s\in S}}$ is a split semi-simple algebra and models for its representations are then known. In types $A_n$, $B_n$ and $D_n$, the irreducible representations are labeled by partitions or double-partitions of $n$. A partition of $n$ is a non-increasing sequence $(\lambda_i)_{i\in \N^\star}$ such that $\underset{i=1}{\overset{+\infty}\sum}\lambda_i=n$. To each partition of $n$, we can associate a Young diagram which is a diagram with $\lambda_i$ boxes in the ith row. We write $[2^2,1]$ the partition of $5$ with $\lambda_1=\lambda_2=2$ and $\lambda_3=1$. There are $5$ partitions of $4$ given by the Young diagrams
$$ \begin{Young}
\cr
 \cr
 \cr
 \cr
 \end{Young}, \begin{Young}
 & & &\cr
 \end{Young}, \begin{Young}
 & &\cr
 \cr
 \end{Young}, \begin{Young}
 & \cr
 \cr
 \cr
 \end{Young}, \begin{Young}
 & \cr
 &\cr
 \end{Young}$$ which we note respectively $[1^4],[4],[3,1],[2,1^2]$ and $[2,2]$. The irreducible representations of the Iwahori-Hecke algebra of type $A_n$ are labeled by partitions of $n$. There is a basis for each module associated to a partition of $n$ given by standard tableaux associated to that partition. A tableau associated to a partition of $n$ is numbering with the integers from $1$ to $n$. A standard tableau is a tableau such that the numbering is increasing towards the right and downwards. The representations labeled by the above partitions are therefore of respective dimensions $1$, $1$, $3$, $3$ and $2$. A basis for the representation labeled by the partition $[3,1]$ is given by the standard tableaux $$ \begin{Young}
 1& 2&3\cr
 4\cr
 \end{Young},  \begin{Young}
 1&2 &4\cr
 3\cr
 \end{Young},  \begin{Young}
 1&3 &4\cr
 2\cr
 \end{Young}. $$
 
 \smallskip
 
 In type $B_n$ and $D_n$, the irreducible representations are labeled by double-partitions of $n$. 
 A double-partition of $n$ is pair $(\lambda,\mu)$ with $\lambda$ a partition of $r$ and $\mu$ a partition of $n-r$. 
 The double-diagram associated to the double-partition $(\lambda,\mu)$ is the pair given by the Young diagrams associated to $\lambda$ and $\mu$. 
 A standard double-tableau associated to a double-partition is a numbering of the associated double-diagram with the integers from $1$ to $n$ increasingly towards the right and downwards within each component of the double-diagram. 
 There are $5$ double-partitions of $2$ given by the double-diagrams $$(\begin{Young} \cr \cr \end{Young} ,\emptyset), (\begin{Young} &\cr \end{Young} ,\emptyset), (\begin{Young} \cr \end{Young} ,\begin{Young} \cr \end{Young} ), (\emptyset, \begin{Young} \cr \cr \end{Young} ), (\emptyset, \begin{Young} \cr \cr \end{Young} )$$ which we note respectively $([1^2],\emptyset)$, $([2],\emptyset)$, $([1],[1])$, $(\emptyset,[1^2])$ and $(\emptyset,[2])$.
 The representations labeled by the above representations are of respective dimensions $1$, $1$, $2$, $1$ and $1$. A basis for the representation labeled by the double-partition $([1],[1])$ is given by the standard double-tableaux (\begin{Young} 1\cr \end{Young} ,\begin{Young}2 \cr \end{Young} ) (\begin{Young} 2\cr \end{Young} ,\begin{Young} 1\cr \end{Young} ).
 
 We will give the explicit models in the following chapters.
 
 The representations for Coxeter groups of exceptional types are afforded by $W$-graphs, they will be explained in Chapter \ref{Wgraphschapter}.
 
 In the finite field setting, we will encounter some field automorphism through which the representations will factor. We give here a condition for a field automorphism to give a representation when composed with a representation of an Iwahori-Hecke algebra.

\begin{prop}\label{Fieldfactorization}
Let $(W,S)$ be a Coxeter system and $\mathcal{H}_{(\alpha_s)_{s\in S}}$ be a finite Iwahori-Hecke $\F_p((\alpha_s)_{s\in S})$-algebra. Let $\Phi\in \op{Aut}(\F_q)$, where $\F_q=\F_p((\alpha_s)_{s\in S})$ and $\rho$ be a finite dimensional representation of $\mathcal{H}_{(\alpha_s)_{s\in S}}$. We then have that there exists a character $\eta$ such that $(\Phi\circ \rho)\otimes \eta$ is a representation of $\mathcal{H}_{(\alpha_s)_{s\in S}}$ if and only if for all $s\in S$,  $\Phi(\alpha_s)\in \{\alpha_s,\alpha_s^{-1}\}$.

Moreover, if for all $s\in S$, $\Phi(\alpha_s)\in \{\alpha_s,\alpha_s^{-1}\}$, then $\eta(T_s)=1$ if $\Phi(\alpha_s)=\alpha_s$ and $\eta(T_s)=-\alpha_s$ if $\Phi(\alpha_s)=\alpha_s^{-1}$.
\end{prop}

\begin{proof}
Assume first that for all $s\in S$, we have $\Phi(\alpha_s)\in \{\alpha_s,\alpha_s^{-1}\}$. It is clear that the braid relations are verified so we only need to check that the deformations of the relations of order $2$ are verified. Let $s\in S$, assume $\Phi(\alpha_s)=\alpha_s$. Let $\eta(T_s)=1$, we then have 
$$(\Phi(\rho(T_s))\eta(T_s)-\alpha_s)(\Phi(\rho(T_s))\eta(T_s)+1)  =  \Phi((\rho(T_s)-\alpha_s)(\rho(T_s)+1))\\
 = \Phi(0)=0.$$
 Assume $\Phi(\alpha_s)=\alpha_s^{-1}$, let $\eta(T_s)=-\alpha_s$, we then have
 \begin{eqnarray*}
 (\Phi\circ \rho)(T_s))\eta(T_s)-\alpha_s)(\Phi\circ \rho(T_s)\eta(T_s)+1) & = & \Phi((-\alpha_s^{-1}\rho(T_s)-\alpha_s^{-1})(-\alpha_s^{-1}\rho(T_s)+1))\\
 & = & \Phi(-\alpha_s^{-2}(\rho(T_s)+1)(\rho(T_s)-\alpha_s))\\
 & = & \Phi(0)\\
 & = & 0.
 \end{eqnarray*}
 This proves that $(\Phi\circ \rho)\otimes \eta$ is a representation of $\mathcal{H}_{(\alpha_s)_{s\in S}}$.
 
 \medskip
 
 Assume now that there exists a character $\eta$ such that $(\Phi\circ \rho)\otimes \eta$ is a representation of $\mathcal{H}_{(\alpha_s)_{s\in S}}$. Let $s\in S$, the eigenvalues of $\Phi(\rho(T_s))\eta(T_s)$ are then $\alpha_s$ and $-1$. We also have that the eigenvalues of $\Phi(\rho(T_s))\eta(T_s)$ are $\eta(T_s)\Phi(\alpha_s)$ and $-\eta(T_s)$. It follows that either $-\eta(T_s)=-1$ and $\eta(T_s)\Phi(\alpha_s)=\alpha_s$ or $-\eta(T_s)=\alpha_s$ and $\eta(T_s)\Phi(\alpha_s)=-1$. This implies that either $\Phi(\alpha_s)=\alpha_s$ or $\Phi(\alpha_s)=\alpha_s^{-1}$. The proof is thus concluded.

\end{proof}

\section{Symmetric algebras and specializations}\label{Symalg}

In this section, we recall the definition and some basic properties of symmetric algebras. We then give a version of Tits deformation theorem, which we will need to prove that the finite Iwahori-Hecke algebras are split semisimple under the right conditions. We then give a general result on representations of finite Iwahori-Hecke algebras, which we will use throughout the thesis.

\subsection{Definition and first properties}
In this subsection, we define symmetric algebreas and give some elementary properties of the corresponding trace. The results in this section are taken from Chapter $7$ of \cite{G-P}. Throughout this section $A$ is a commutative ring and $H$ is an $A$-algebra of finite rank.

\begin{Def}
A linear form $\tau:H\rightarrow A$ is said to be a trace if it verifies the following condition
$$\forall a,b \in H, \tau (ab) = \tau (ba).$$
If there exists a trace $\tau$ such that $(x,y) \mapsto \tau(x,y)$ is non-degenerate then we say that $H$ is a symmetric $A$-algebra and $\tau$ is a symmetrizing trace.
\end{Def}

Iwahori-Hecke algebras are symmetric algebras and the following proposition is given in \cite{G-P} (Proposition $8.1.1$).

\begin{prop}
Let $W$ be a finite irreducible Coxeter group. The Iwahori-Hecke algebra $\mathcal{H}_{W,R,(u_s)_{s\in S}}$ is then a symmetric $R[u_s^{\pm 1}]$ algebra with symmetrizing trace $\tau$ defined by $\tau(T_0)=1$ and $\tau(T_\sigma)=0$ if $\sigma\neq 1_W$.
\end{prop}

In what follows $H$ is a symmetric $A$-algebra of finite rank and $\tau$ is a trace of $H$.

\begin{Def}
Let $\mathcal{B}$ be a basis for $H$, its dual basis with regards to $\tau$ is defined to be $(\check{b})_{b\in \mathcal{B}}$, where for all $b,b'$ in $\mathcal{B}\times \mathcal{B}$, we have $\tau(\check{b}b') = \delta_{b,b'}$.

For $\sigma=s_{i_1}\dots s_{i_n}$ a reduced expression of $\sigma$, we define $u_\sigma=u_{s_{i_1}}\dots u_{s_{i_n}}$. For the Iwahori-Hecke algebra $\mathcal{H}_{W,R,(u_s)_{s\in S}}$, $(T_{\sigma})_{\sigma\in W}$ is a basis. Its dual basis with regards to $\tau$ is $(u_{\sigma^{-1}}T_{\sigma^{-1}})$. 
\end{Def}

We now consider a fixed basis $\mathcal{B}$ of $H$ in what follows.

\begin{Def}
Let $V$ and $V'$ be two right $H$-modules. For all $\varphi \in Hom_A(V,V')$, we define $I(\varphi)$ from $V$ to $V'$ by $I(\varphi)(v) = \underset{b\in \mathcal{B}}\sum \varphi(vb)\check{b}$.
\end{Def}

We consider $A$ as a subset of $H$ by considering that for $a \in A$, $a = a.1_H$.

\begin{prop}
Let $V$ and $V'$ be two right $H$-modules and $\varphi \in Hom_A(V,V')$.\\
We then have that $I(\varphi)$ is independant of the chosen basis and $I(\varphi) \in Hom_H(V,V')$.
\end{prop}

The previous proposition shows that $I$ gives us a way to obtain $H$-linear map from $A$-linear maps.

\begin{prop}
Let $V, V'$ and $V''$ be right $H$-modules and $\varphi\in Hom_H(V,V')$,\\
$ \psi \in Hom_A(V',V'')$ and $\theta \in Hom_A(V'',V)$. We have
$$ I(\psi \circ \varphi) = I(\psi) \circ \varphi, I(\varphi \circ \theta) = \varphi \circ I(\theta).$$
\end{prop}

\begin{Def}
We say an $H$-module $V$ is projective if for every surjective $H$-module morphism from $M$ to $V$ with $M$ an $H$-module, there exists $i \in Hom_H(V,M)$ such that $\pi \circ i = id_V$.
\end{Def}

\begin{lemme}[Gaschütz-Ikeda]
Let $V$ be an $H$-module which is projective as an $A$-module. $V$ is projective as an $H$-module if and only if there exists $\varphi \in End_A(V)$ such that $I(\varphi) = id_V$.
\end{lemme}

\subsection{Schur elements}

In this subsection, we recall the definition of Schur elements which will give us tools to understand the specializations to finite fields of the irreducible representations. 

Let $K$ be a field and $H$ a symmetric $K$-algebra with symmetrizing trace $\tau$. We will use the map $I$ in order to define the Schur elements. We fix a $K$-basis $\mathcal{B}$ of $H$.

\begin{Def}
An $H$-module $V$ is sais to be split semi-simple if it is simple and we have$\dim_K(\op{End}_H(V)) = 1$.

A $K$-algebra $H$ is said to be split if all its simple modules are split.
\end{Def}

\begin{theo}
Let $V$ be a split simple $H$-module. There exists a unique $c_V \in K$ such that for all $\varphi \in End_K(V)$, we have $I(\varphi) = c_VTr(\varphi)id_V$ and $c_V$ only depends on the isomorphism class of $V$.

$c_V$ is called the Schur element of $H$ associated to the isomorphism class of $V$.
\end{theo}

\begin{cor}
Let $V$ and $V'$ be two split simple $H$-modules with $n = n_V = \dim(V)$ and $m = n_{V'} = \dim(V')$. Let $\rho$(resp $\rho'$) be a representation of $H$ in $M_n(K)$(resp $M_m(K)$), we then have for all $(i,l,k,j) \in [\![ 1,n]\!]^2\times [\![ 1,m]\!]^2$
\begin{eqnarray*}
\underset{b\in \mathcal{B}}\sum \rho(b)_{i,l}\rho'(\check{b})_{k,j} & = & \delta_{i,j}\delta_{k,l} c_V ~\text{if V is isomorphic to V'}~ \rho = \rho',\\
 & = & 0~  \text{if V is not isomorphic to V'}.
\end{eqnarray*}
\end{cor}

We now recall Wedderburn's theorem on semi-simple algebras.

\begin{theo}(Wedderburn) \\
Let $H$ be a finite dimensional split semi-simple $K$-algebra. We then have $H = \underset{V}\bigoplus H(V)$, where the sum is over the isomorphism classes of simple $H$-modules. In the above, $H(V) \simeq M_{n_V}(D_V)$, where we have $D_V = End_H(V)$.
\end{theo}

We now give a semi-simplicity criteria with a condition on the Schur elements.

\begin{theo}
A split semi-simple $H$-module $V$ is projective if and only if $c_V \neq 0$.

If $H$ is split then $H$ split semi-simple if and only if all its Schur elements are non-zero and we then have 
$$\tau = \underset{V}\sum \frac{1}{c_V}\chi_V.$$
where the sum is over the isomorphism classes of $H$-modules.
\end{theo}

\subsection{Specialization theorem}
 
 We now give a proposition which we will use to show that the specializations of the Iwahori-Hecke algebras to finite fields remain split semisimple under the right conditions. This is a version from specialization theory of symmetric algebras and Schur elements which we have not found as such in the litterature.

\begin{prop}\label{Tits}
Let $A$ be a commutative integrally closed integral domain, $F$ a field containing $A$ and $H$ a free $F$-algebra of finite rank with a symetrizing trace $\tau$.\\
We write $FH=H\otimes_A F$ and consider $H$ as a subset of $FH$. If $B$ is a ring such that $A\subset B\subset F$, then we consider that $H\subset BH\subset FH$.\\
Let $\theta$ be a ring homomorphism from $A$ to a field $L$ such that $L$ is the field of fractions of $\theta(A)$.

Let $\mathcal{O}$ be a valuation ring of $F$ such that $A\subset \mathcal{O}$ and $\mathcal{J}(\mathcal{O})\cap A=\op{ker}(\theta)$, where $\mathcal{J}(\mathcal{O})$ is the unique maximal ideal of $\mathcal{O}$. Let $k=\mathcal{O}/\mathcal{J}(\mathcal{O})$ be the residue field of $\mathcal{O}$. The restriction to $A$ of the projection $\mathcal{\pi}$ from $\mathcal{O}$ to $k$ has kernel equal to $\mathcal{J}(\mathcal{O})\cap A =\op{ker}(\theta)$. We can therefore see $L$ as a subfield of $k$.

\medskip

Let $F$ be the field of fractions of $A$ and $B$ the subring of $F$ formed of elements of the form $\frac{a_1}{a_2}$, $(a_1,a_2)\in A^2$ and $\theta(a_2) \ne 0$. Assume $FH$ is split. We can define an extension $\tilde{\theta}$ of $\theta$ to $B$. Assume there exists a representation of every simple module $V$ such that for $(i,j) \in [\![ 1,n]\!]^2$, there exists $b\in B$ such that $\tilde{\theta}(\rho_V(b)_{j,i})\neq 0$. Assume there exists a basis $\mathcal{B}$ of $BH$ such that for all $b\in \mathcal{B}$, $\check{b}\in BH$ and $\tilde{\theta}(\check{b})\neq 0$, where for $b\in B$, $\check{b}$ is the unique element in $B$ such that $\tau(b\check{b})=1$.

\medskip
We then have that $LH$ is split semi-simple if and only if $\theta(c_V)\neq 0$ for every simple $FH$-module $V$
\end{prop}

\begin{proof}

Assume $\theta(c_V)\neq 0$ for every simple $FH$-module $V$. We then have that $c_V\neq 0$ for all $FH$-module simple $V$.

Let $\mathcal{B}$ be a basis as in the theorem, by Proposition $7.2.7$ of \cite{G-P}, $\left(e_{i,j,V}=\frac{1}{c_V}\underset{b\in\mathcal{B}}\sum \rho_V(b)_{j,i}\check{b}\right)_{i,j,V}$ is a basis of $BH$ corresponding to the isomorphism $BH \simeq\underset{V}\bigoplus M_{n_V}(B)$.

We have that the $\tilde{\theta}(e_{i,j,V})$ are in $LH$. Let us show that they give us a basis affording an explicit isomorphism from $LH$ to $\underset{V}\bigoplus M_{n_V}(L)$. Since $L$ is the field of fractions of $\theta(A)$, $\tilde{\theta}$ is a surjective morphism from $B$ to $L$ and the $\tilde{\theta}(e_{i,j,V})$ give us a generating family of $LH$.

 \bigskip\noindent
  We have the relations $\tilde{\theta}(e_{i,j,V})\tilde{\theta}(e_{i',j',V'}) = \delta_{V,V'}\delta_{j,i'}\tilde{\theta}(e_{i,j',V})$.
  
Let $(a_{i,j,V})\in L^{\underset{V}\sum n_V^2}$, assume $\underset{V}\sum\underset{1\leq i,j\leq n_V}\sum a_{i,j,V}\tilde{\theta}(e_{i,j,V})=0$. Let $V_0$ be a simple $FH$-module and $1\leq i_0,i_1,j_0,j_1\leq n_{V_0}$. Multiplying by $\tilde{\theta}(e_{i_0,j_0,V_0})$ on the right and $\tilde{\theta}(e_{i_1,j_1,V_0})$ on the left, we have that
$$0 = \underset{V}\sum\underset{1\leq i,j\leq n_V}\sum a_{i,j,V}\tilde{\theta}(e_{i_1,j_1,V_0})\tilde{\theta}(e_{i,j,V}) \tilde{\theta}(e_{i_0,j_0,V_0}),$$
 $$ 0 = \underset{i=1}{\overset{n_{V_0}}\sum} a_{i,i_0,V_0}\tilde{\theta}(e_{i_1,j_1,V_0})\tilde{\theta}(e_{i,j_0,V_0}) = a_{j_1,i_0,V_0}\tilde{\theta}(e_{i_1,j_0,V_0}).$$ 
 Since each $\tilde{\theta}(e_{i,j,V})$ is non-zero, these vectors are linearly independent and the relations they verify give us an isomorphism from $LH$ to $\underset{V}\bigoplus M_{n_V}(L)$. It follows that $LH$ is split semi-simple.

\bigskip

Assume now that $LH$ is split semi-simple. By Tits's deformation theorem (Theorem 7.4.6 of \cite{G-P}),  there is an isomorphism $d_L^k$ between the Grothendieck groups of finite dimensional $LH$-modules and finite dimensional $kH$-modules sending each isomorphism class of simple modules to a unique isomorphism class of simple modules. If $V$ is a simple $LH$-module, the unique $kH$-module $V'$ such that $[V']=d_L^k([V])$ verifies $\dim_k(V')=\dim_L(V)$. We have $$\dim_k(kH)=\dim_L(LH)=\underset{V}\sum \dim_L(V)^2=\underset{V'}\sum \dim_k(V')^2=\dim_k(kh/\op{rad}(kH)).$$
$kH$ is therefore semi-simple. Since $LH$ is split, $kH$ is split too. By Tits's deformation theorem applied to the specialization $\pi$ from $\mathcal{O}$ to $k$, $FH$ is semi-simple and there is an isomorphism $d_\pi$ sending isomorphism classes of $FH$-simple modules and isomorphism classes of simple $kH$-modules.

Let $V$ be a $FH$-module and $\rho_V$ a representation of $FH$ in $M_{n_V}(F)$ such that for all $h\in H$, we have $\rho_V(h)\in \mathcal{O}$. By Corollary $7.2.2$ of \cite{G-P}, we have that
$$\underset{b\in \mathcal{B}}\sum \rho_V(b)_{1,1}\rho_V(\check{b})_{1,1} = c_V, \underset{b\in \mathcal{B}}\sum \pi(\rho_V(b)_{1,1})\pi(\rho_V(\check{b})_{1,1}) = \pi(c_V) = \theta(c_V).$$

This equality corresponds to the relation verified by the Schur element of a $kH$ simple module $V'=k\tilde{V}$. We then have $c_{V'}=\theta(c_V)$, therefore $\theta(c_V)\neq 0$ because $kH$ is split semi-simple. This concludes the proof.
\end{proof}

\section{Maximal subgroups of classical groups}\label{sectionAschba}
 
In this section, we recall Aschbacher's theorem on maximal subgroups of classical groups \cite{ASCHMAXSUBGRPS} and describe the different classes. We then give some important theorems about subgroups of classical groups which we will use throughout the thesis.

We first recall the following definitions from \cite{Aschbabook}
\begin{Def}
Let $G$ be a group.

A subnormal group $H$ of $G$ is a subgroup of $G$ such that there exists a sequence $H=G_0\subset G_1\subset ...\subset G_n=G$ and for all $i\in [\![1,n]\!]$, $G_{i-1}\triangleleft G_i$

A group $G$ is quasisimple if it is perfect and $G/Z(G)$ is simple.

The components of a group $G$ are its subnormal quasisimple subgroups. Write $\op{Comp}(X)$ for the set of components of $X$ and set $E(X)=\langle \op{Comp}(X)\rangle$.

The Fitting subgroup $F(G)$ of $G$ is the largest nilpotent normal subgroup of $G$.

The generalized Fitting group of $G$ is the group $F^\star(G)=F(G)E(G)$.

\end{Def}

We begin by stating Aschbacher's theorem on maximal subgroups \cite{ASCHMAXSUBGRPS}

\begin{theo}
Let $G$ be a finite group whose generalized Fitting group is a simple classical group $G_0$ over a finite field such that $G_0\not\simeq P\Omega_8^+(q)$. Let $H$ be a proper subgroup of $G$ such that $G=HG_0$. Then either $H$ belongs to one of the geometric classes $(\mathcal{C}_i)_{i\in [\![1,8]\!]}$ or $H$ belongs to the class $\mathcal{S}$.
\end{theo}

We will not use this theorem directly but results that follow from this theorem. In particular, we will use many tables from \cite{BHRC} of the maximal subgroups of classical groups in low dimension, where all the conjugacy classes of maximal subgroups are determined. In table \ref{Aschbaclasses}, we copy table $2.1$ from \cite{BHRC} giving a rough description of the Aschbacher classes for classical groups over $\F_q$ (the unitary groups $SU_n(q)$ are defined over $\F_{q^2}$).

\begin{table}
\noindent
\begin{tabular}{ |p{0.5cm}||p{14cm}|   }
 \hline
  $\mathcal{C}_i$ & Rough description \\
 \hline
 \hline

 $\mathcal{C}_1$   & stabilizers of totally singular or non-singular subspaces\\
 \hline
 
 $\mathcal{C}_2$   & stabilizers of decompositions $V=\underset{i=1}{\overset{t}\oplus}$, $\dim(V_i)=\frac{n}{t}$\\
 \hline
 
 $\mathcal{C}_3$   & stabilizers of extension fields of $\F_q$ of prime index dividing 
 $n$\\
 \hline
 
 $\mathcal{C}_4$   & stabilizers of tensor product decompositions $V=V_1\otimes V_2$\\
 \hline
 
 $\mathcal{C}_5$   & stabilizers of subfields of $\F_q$ of prime index\\
 \hline
 
 $\mathcal{C}_6$   & normalizers of symplectic-type or extraspecial groups in absolutely irreducible representations \\
 \hline
 
 $\mathcal{C}_7$   & stabilizers of decompositions $V=\underset{i=1}{\overset{t}\otimes}V_i$, $\dim(V_i)=a$, $n=a^t$\\
 \hline
 
 $\mathcal{C}_8$   & groups of similarities of non-degenerate classical forms\\
 \hline
\end{tabular}
\smallskip

\caption{Description of the geometric Aschbacher classes}
\label{Aschbaclasses}
\end{table}

By definition of the classes, no irreducible group can be included in a group of class $\mathcal{C}_1$, no primitive group can be included in a group of class $\mathcal{C}_2$, no group containing a transvection can be included in a group of class $\mathcal{C}_3$, no tensor-indecomposable group can be included in a group of class $\mathcal{C}_4$ and no group whose trace generate the field $\F_q$ can be included in a maximal group of class $\mathcal{C}_5$. We will use this type of argument to reduce the lists of possible maximal groups containing the groups we consider for groups of low dimension in the following sections.

The main theorem we will use in high dimension is a theorem by Guralnick and Saxl \cite{GS} giving conditions for subgroups of the special linear group to be classical groups. The proof of Theorem \ref{CGFS} uses the classification of finite simple groups.

\begin{theo}[Gulralnick-Saxl]\label{CGFS}
Let $V$ be a finite-dimensional vector space of dimension $d > 8, d\neq 10$ or ($d=10$ and $p= 2$) over an algebraically closed field $\overline{\F_p}$ of characteristic $p>0$. Let $G$ be a primitive tensor-indecomposable finite irreducible subgroup of $GL(V)$. We write $v_G(V)$ the minimal dimension of $[\beta g,V] = (\beta g -1)V$, for $g\in G$ and $\beta \in \overline{\F_p}$ such that $\beta g \neq 1$. We then have either $v_G(V) > \max(2,\frac{\sqrt{d}}{2})$ or one of the following assertions.
\begin{enumerate}
\item G is a classical group in a natural representation. 
\item G is the alternating or the symmetric group of degree $c$ and $V$ is the permutation module of dimension $c-1$ or $c-2$.
\end{enumerate}
\end{theo}

The other theorems we will use are on groups generated by long-root elements, we first recall the definition of long-root elements.

\begin{Def}
If $G\simeq SL_n(q)$, $G\simeq SU_n(q)$ or $G\simeq SP_n(q)$ then a long-root element of $G$ is a transvection.

If $G\simeq \Omega_n^+(q)$, $n\geq 4$ then a long-root element $x$ of $G$ is an element of the form $x(v)=v-\langle v,a\rangle b +\langle v,b\rangle a$ for $a,b$ in a totally singular 2-space $T$ and $\langle .,.\rangle$ the non-degenerate symmetric bilinear form associated to $G$.
\end{Def}

For example, if $n=2m\geq 4$ is an integer then elements conjugate to elements of the form $\begin{pmatrix}
T_a & 0\\
0 & ^t(T_a)^{-1}
\end{pmatrix}$, where $T_a$ is a transvection of $SL_m(q)$ are long-root elements of $\Omega_4^+(q)$.

We can now state Kantor's Theorem on subgroups of orthogonal groups generated by long root elements \cite{K}. We only consider the irreducible subgroups of $\Omega_n^+(q)$ with $q$ an odd prime because it will be the only case we will need to consider in this thesis.

\begin{theo}\label{theoKantor}
Let $G$ be an irreducible subgroup of $\Omega^+_n(q)$, where $n\geq 4$ and $q=p^s$ for some prime $p$ and positive integer $s$ generated by a conjugacy class of long root elements, such that $O_p(G)\leq [G,G]\cap Z(G)$.

We then have that $G$ belongs to the following list
\begin{enumerate}
\item $\Omega^\pm(q')$ in a natural representation over $\F_{q'}$, $q'|q$,
\item $\Omega_{2m}^-(q'^{\frac{1}{2}})\leq \Omega_{2m}^+(q')$,  $q'|q$, $n=2m$ in a natural representation over $\F_{q'}$,
\item $SU_{2m}(q')\leq \Omega_{4m}^+(q')$, $n=4m$ in a natural representation over $\F_{q'}$, $q'|q$
\item $SU_{2m+1}(q')\leq \Omega_{4m+2}(q')$, $n=4m+2$ in a natural representation over $\F_{q'}$, $q'|q$,
\item $G/Z(G)\simeq P\Omega(7,q')$, $\vert Z(G)\vert=(2,q'-1)$, $G\leq \Omega_8^+(q')$ in a natural representation over $\F_q'$, $q'|q$,
\item $^3\! D_4(q')\leq \Omega_8^+(q'^3)$ in a natural representation over $\F_{q'^3}$, $q'^3|q$.
\end{enumerate}
\end{theo}

Finally we give a Theorem by Serezkin and Zaleskii for irreducible groups of the remaining classical groups generated by transvections. (First theorem of \cite{SZ}).

\begin{theo}\label{transvections}
If $G$ is an irreducible subgroup of $GL_n(q)$ generated by transvections with $q=p^r, p> 3, n> 2$, then $G$ is conjugate inside $GL_n(q)$ to $SL_n(\tilde{q}), Sp_n(\tilde{q})$ or $SU_n(\tilde{q}^{\frac{1}{2}})$ for some $\tilde{q}$ dividing $q$.
\end{theo}
 
 \chapter{Type B}\label{TypeB}
 
In this section, we will determine the image of the derived subgroup of the Artin group of type $B_n$. This is a natural continuation of the work for type $A_n$ in \cite{BM} and \cite{BMM}. The general outline of the proof is based on inductive reasoning. We first give the factorizations through field automorphisms and the transposed inverse automorphisms between irreducible representations. The main result concerning those factorizations is given in Proposition \ref{isomorphisme}. This allows us to determine what the image in the full finite Iwahori-Hecke algebra appears to be. The remainder of the section is the proof that the image appearing is indeed the image of the derived subgroup $\mathcal{A}_{B_n}$ of the Artin group. The main results of this section are given in Theorems \ref{result1} to \ref{result6}.

\bigskip

Let $p$ be a prime, $n\geq 2$ be an integer, $\alpha\in \overline{\F_p}$ of order $a$ greater than $n$ and not in $\{1,2,3,4,5,6,8,10\}$ and $\beta \in \overline{\F_p}\setminus\{-\alpha^i,-(n-1)\leq i\leq n-1\}$ different from $1$. We set $\F_q =\F_p(\alpha,\beta)$. The Artin group of type $B$ is the group generated by the elements $T=S_0,S_1,\dots,S_{n-1}$ verifying the relation $S_0S_1S_0S_1=S_1S_0S_1S_0$, for $i\in [\![1,n-2]\!]$, $S_iS_{i+1}S_i=S_{i+1}S_iS_{i+1}$ and for $(i,j)\in [\![0,n-1]\!]$ such that $\vert i-j\vert \geq 2$, $S_iS_j=S_jS_i$. The associated Iwahori-Hecke $\F_q$-algebra $\mathcal{H}_{B_n,\alpha,\beta}$ is defined by the generators indexed in the same way as for the Artin group and verifying the previous relations and deformations of the relations of order $2$ of the Coxeter groups : $(T-\beta)(T+1)=0$ and for $i\in [\![1,n-1]\!]$, $(S_i-\alpha)(S_i+1)=0$. In the sequel, we identify the Artin group with its image inside the Iwahori-Hecke algebra. We write $\ell_1,\ell_2$ for the length functions on $A_{B_n}=\langle T,S_i\rangle_{i\in [\![1,n-1]\!]}$ such that for all $i\in [\![1,n-1]\!]$, $\ell_1(S_i)=1, \ell_1(T)=0, \ell_2(S_i)=0$ and $\ell_2(T)=1$.

In Section \ref{HoefmodelB} we give in Theorem \ref{mod} the irreducible representations described by the Hoefsmit model given in \cite{G-P} and \cite{HOEF} and define a weight on standard tableaux and double-partitions of $n$ which allows us to define in Proposition \ref{bilin} a bilinear form verifying nice properties.

In Section \ref{FactorB}, we determine all the isomorphisms between different irreducible representations and then state the main results for type $B$ in Theorems \ref{result1} up to \ref{result6}.

In Section \ref{surjectivitism}, we prove the result in all possible cases depending on the properties of the field extensions $\F_q$ of $\F_p(\alpha+\alpha^{-1},\beta+\beta^{-1})$ and $\F_p(\alpha)$ of $\F_p(\alpha+\alpha^{-1})$.

\section{Hoefsmit model and first properties}\label{HoefmodelB}

In this section, we define the matrix model we will be considering throughout this section and define a natural bilinear form appearing when we consider those models. This will establish the groundwork for the rest of this section.

\begin{theo}\label{mod}
Assume $\alpha$ is of order greater than $n$ and $\beta \in \overline{\F_p}\setminus\{-\alpha^i,-(n-1)\leq i\leq n-1\}$. The following matrix model gives a list of the pairwise non-isomorphic absolutely irreducible modules $V_\lambda$ of the Iwahori-Hecke algebra $\mathcal{H}_{B_n,\alpha,\beta}$ labeled by double-partitions of $n$.\\
If $\T =(\T_1,\T_2)$ is a standard double-tableau, then
\begin{itemize}
\item if $1\in \T_1$, then $T.\T= \beta \T$,
\item if $1\in \T_2$, then $T.\T = -\T$,
\item $S_i .\T = m_i(\T)\T+(1+m_i(\T))\tilde{\mathbb{T}}$, where $\tilde{\mathbb{T}} = \T_{i\leftrightarrow i+1}$ if $\T_{i\leftrightarrow i+1}$ is standard and $0$ otherwise.
\end{itemize}
Above, we have $m_i(\T)=\frac{\alpha-1}{1-\frac{ct(\T:i)}{ct(\T:i+1)}}$, $ct(\T:j) = \alpha^{c_j(\T)-r_j(\T)}\beta$ if $j \in \T_1$ and $ct(\T:j)=-\alpha^{c_j(\T)-r_j(\T)}$ otherwise and $r_j(\T)$ (resp $c_j(\T)$) is the row (resp column) of $j$ in the tableau of $\T$ containing $j$.
\end{theo}

\begin{proof}
Let $\lambda=(\lambda_0,\lambda_1)\Vdash n$ such that $n_\lambda=\dim(V_\lambda)>1$.

We will show that $\rho_{\lambda}(\mathcal{A}_{B_n})$ contains diagonal matrices generating the algebra of diagonal matrices in $GL_n(q)$. By Proposition $5$ of \cite{M}, it is sufficient to show that for any standard double-tableau $\mathbb{T}$ associated to $\lambda$, there exist two diagonal invertible matrices $D_1$ and $D_2$ in the basis of standard double-tableaux such that $D_1 \mathbb{T}=b_1\mathbb{T}$ and $D_2 \mathbb{T}=b_2\mathbb{T}$ and $b_1\neq b_2$.

By \cite{M}, it will then be sufficient to show that for any couple of standard double-tableaux $(\T_\rho,\T_\gamma)$, there exists a matrix in $\rho_{\lambda}(\mathcal{A}_{B_n})$ such that $a_{\rho,\gamma}$ is non-zero.

In order to do that, we will use the Jucys-Murphy elements $(a_1,...,a_n)$ whose expressions in this model are given in \cite[Prop 3.16]{A-K} and are

$a_i.\T_\rho = u_{\tau(i)}\alpha^{c_i-r_i+i-1}\T_\rho$, where $i$ is in box $(r_i,c_i)$ of the $\tau(i)$-th tableau of $\T_\rho$ and $(u_1,u_2)=(\beta,-1)$.

\bigskip

Let $\T_\rho=(\T_{\rho_1},\T_{\rho_2})\in \lambda$. Since $\dim(V_\lambda)>1$, we can define $\ell=\op{min}\{j\in [\![1,n]\!], j\notin \T_{\rho,\tau(1)} ~or~ (j\in \T_{\rho,\tau(1)} ~and~  (r_j,c_j)\neq (j,1))\}$ and we have $\ell\geq 2$. 

We then have $k\in \T_{\rho,\tau(1)}$ and $(r_k,c_k)=(k,1)$ for all $k<\ell$. It follows that either $\ell\in \T_{\rho,3-\tau(1)}$ and $(r_\ell,c_\ell) =(1,1)$ or $\ell\in \T_{\rho,\tau(1)}$ and $(r_\ell,c_\ell)=(1,2)$. We have

$a_1.\T_\rho=u_{\tau(1)}\T_\rho$ and $a_\ell.\T_\rho\in \{u_{3-\tau(1)}\alpha^{\ell-1}\T_\rho,u_{\tau(1)}\alpha^\ell\T_\rho\}$, therefore $a_1.\T_\rho \neq a_\ell.\T_\rho$ since the order of $\alpha$ is strictly greater than $n$ and $\beta \notin \{-\alpha^i,-(n-1)\leq i\leq n-1\}$.

\bigskip

Let now $(\T_\rho,\T_\gamma)$ be a pair of standard double-tableaux associated to $\lambda$. There exists a permutation of $\mathfrak{S}_n$ which affords $\T_\gamma$ after permuting the numbers inside $\T_\rho$. We can decompose this permutation in a product of transpositions $(i,i+1)$ such that the path given by the successive standard double-tableaux is composed only of standard double-tableaux. It is thus sufficient to show that for any $i\in [\![1,n]\!]$ and for any standard double-tableau $\T_\rho$ such that $\T_{\rho,i\leftrightarrow i+1}$ is standard, we have $1+m_i(\T_\rho)\neq 0$.

We have $m_i(\T_\rho)=\frac{\alpha-1}{1-\frac{u_{\tau(i)}\alpha^{c_i-r_i+r_{i+1}-c_{i+1}}}{u_{\tau(i+1)}}}=-1 \Leftrightarrow \frac{u_{\tau(i)}}{u_{\tau(i+1)}}=\alpha^{1+c_{i+1}-r_{i+1}+r_i-c_i}$ and $\frac{u_{\tau(i)}}{u_{\tau(i+1)}}\in \{1,-\beta,-\beta^{-1}\}$. By the assumptions made on $\alpha$ and $\beta$, it now only remains to show that $-n\leq 1+c_{i+1}-r_{i+1}+r_i-c_i\leq n$ and $1+c_{i+1}-r_{i+1}+r_i-c_i\neq 0$ if $i$ and $i+1$ are in the same tableau and $-n+1\leq 1+c_{i+1}-r_{i+1}+r_i-c_i\leq n-1$ if $i$ and $i+1$ are in different tableaux.

\smallskip

If $i$ and $i+1$ are in the same tableau, $\vert c_{i+1}-r_{i+1}+r_i-c_i\vert +1$ represents the number of boxes to go through vertically and horizontally in order to go from $i$ to $i+1$, therefore we have the required bounds. Moreover, this quantity is non-zero by the conditions for a tableau to be standard and the fact that $\T_{\rho, i\leftrightarrow i+1}$ is standard.

\smallskip

If $i$ and $i+1$ are in different tableaux of respective sizes $j$ and $n-j$, we have $1-j \leq c_i-r_i \leq j-1$ and $1-(n-j)\leq r_{i+1}-c_{i+1} \leq n-j-1$ since $1\leq c_i,r_i\leq j$ and $1\leq c_{i+1},r_{i+1} \leq n-j$. It follows that $-n+1\leq 3-n\leq 1+c_{i+1}-r_{i+1}+r_i-c_i\leq n-1$.

\bigskip

The representations $V_\lambda$ are therefore irreducible. We have the same dimensions over $\C$ by taking $\alpha$ and $\beta$ to be irreducible parameters and the same mode. It follows that we only need to show that those modules are pairwise non-isomorphic to conclude the proof.

\bigskip 

Let us show that $R_\lambda\simeq R_\mu \Leftrightarrow \lambda =\mu$. Assume $R_\lambda \simeq R_\mu$, the eigenvalues of the Jucys-Murphy elements must be the same for both representations. We then have that for any $i\in[\![1,n]\!]$, $\{u_{\tau_{\T}(i)}\alpha^{c_{i,\T}-r_{i,\T}+i-1},\T \in \lambda\} = \{u_{\tau_{\tilde{\T}}(i)}\alpha^{c_{i,\tilde{\T}}-r_{i,\tilde{\T}}+i-1},\tilde{\T} \in  \mu\}$ and each of the elements of those sets appear with the same multiplicity. We will show that this implies that for any $i\in [\![1,n]]\!], \gamma \in \Z$, we have

$$\op{card}\{\T \in \lambda, u_{\tau_{\T}(i)}=\beta, c_{i,\T}-r_{i,\T}=\gamma\} =\op{card}\{\tilde{\T} \in \mu, u_{\tau_{\tilde{\T}}(i)}=\beta, c_{i,\tilde{\T}}-r_{i,\tilde{\T}}=\gamma\}$$

Let $\lambda=(\lambda_1,\lambda_2)$, $\mu=(\mu_1,\mu_2)$, $\lambda_1=(\lambda_{1,l})_{l\in \N^\star}$, $\mu_1=(\mu_{1,l})_{l\in \N^\star}$ and $a(i,\lambda,\beta,\gamma)$ be the above quantity.

\smallskip

We first show that $r_1=r_2$, where $\lambda_1\vdash r_1$ and $\lambda_2\vdash r_2$. We have $a(1,\lambda,\beta,0)=\binom  {n-1} {n-r_1} M_{\lambda_1}M_{\lambda_2}$, where $M_{\lambda_i}$ is the number of standard tableaux associated to the partition $\lambda_i$. This is true because counting the number of double-tableaux with $1$ in the left tableau is equivalent to choosing $n-r-1$ numbers in $[\![2,n]\!]$ and counting the number of ways the numbers can be arranged in each tableau to get a standard double-tableau. In the same way, wa have $a(1,\lambda,-1,0)=\binom{n-1} {r_1}  M_{\lambda_1}M_{\lambda_2}$. Since $1$ is either in the first row and first column of the left tableau or in the first row and first column of the right tableau for any standard double-tableau, we have  $a(1,\lambda,\beta,0)=a(1,\mu,\beta,0)$ and $a(1,\lambda,-1,0)=a(1,\mu,-1,0)$. It follows that
\begin{eqnarray*}
\binom{n-1} {n-r_1} M_{\lambda_1}M_{\lambda_2}\binom {n-1} {r_2} M_{\mu_1}M_{\mu_2} & = & \binom {n-1}{r_1} M_{\lambda_1}M_{\lambda_2} \binom  {n-1}{n-r_2} M_{\mu_1}M_{\mu_2}\\
\frac{(n-1)!(n-1)!}{(n-r_1)!(r_1-1)!(n-1-r_2)!r_2!} & = & \frac{(n-1)!(n-1)!}{r_1!(n-r_1-1)!(n-r_2)!(r_2-1)!}\\
\frac{1}{(n-r_1)r_2} & = & \frac{1}{r_1(n-r_2)}\\
nr_2-r_1r_2 & = & nr_1-r_1r_2\\
r_1 & = & r_2.
\end{eqnarray*}

\medskip

In order to show that the quantities are equal, it is sufficient to show that for any $i\in [\![1,n]\!]$, $\T\in \lambda$ and $\tilde{\T}\in \mu$, we have  $u_{\tau_{\T}(i)}\alpha^{c_{i,\T}-r_{i,\T}+i-1}=u_{\tau_{\tilde{\T}}(i)}\alpha^{c_{i,\tilde{\T}}-r_{i,\tilde{\T}}+i-1}$ if and only if $u_{\tau_{\T}(i)}=u_{\tau_{\tilde{\T}}(i)}$ and $c_{i,\T}-r_{i,\T}=c_{i,\tilde{\T}}-r_{i,\tilde{\T}}$. Let $\T\in \lambda$, $\tilde{\T}\in \mu$ such that for all  $i\in [\![1,n]\!]$,  $u_{\tau_{\T}(i)}\alpha^{c_{i,\T}-r_{i,\T}+i-1}=u_{\tau_{\tilde{\T}}(i)}\alpha^{c_{i,\tilde{\T}}-r_{i,\tilde{\T}}+i-1}$.

\smallskip

Let $i\in [\![1,n]\!]$. Assume by contradiction that $u_{\tau_{\T}(i)} \neq u_{\tau_{\tilde{\T}}(i)}$. We then have that \\ $\alpha^{c_{i,\tilde{\T}}-r_{i,\tilde{\T}}+r_{i,\T}-c_{i,\T}}=\frac{u_{\tau_{\T}(i)}}{u_{\tau_{\tilde{\T}}(i)}}\in \{-\beta,-\beta^{-1}\}$. Since $i$ is in a tableau on a different side for $\T$ and $\tilde{\T}$ and $r=r_\lambda=r_\mu$, we have that either $1-r\leq r_{i,\T}-c_{i,\T} \leq r-1$ and $1-(n-r) \leq c_{i,\tilde{\T}}-r_{i,\tilde{\T}} \leq n-r-1$ or $1-(n-r) \leq r_{i,\T}-c_{i,\T} \leq n-r-1$ and $1-r \leq c_{i,\tilde{\T}}-r_{i,\tilde{\T}} \leq r-1$. In both cases, we have $2-n\leq c_{i,\tilde{\T}}-r_{i,\tilde{\T}}+r_{i,\T}-c_{i,\T}\leq n-2$. By assumption, we have $-\beta\notin \{-\alpha^i,1-n\leq i \leq n-1\}$, therefore it follows by contradiction that $u_{\tau_{\T}(i)}=u_{\tau_{\tilde{\T}}(i)}$.

\medskip

Let us show now that $\lambda_{1,1}=\mu_{1,1}$ and $m_\lambda=m_\mu$. 

Assume by contradiction that $\mu_{1,1}<\lambda_{1,1}$. We write $i=\mu_{1,1}+1$. By the equality of the eigenvalues of the Jucys-Murphy elements, there exists $\tilde{T}\in \mu$ such that $u_{\tau_{\tilde{\T}}(i)}\alpha^{c_{i,\tilde{\T}}-r_{i,\tilde{\T}}}=\beta\alpha^{i-1}$ since there exists a standard double-tableau associated to $\lambda$ with $\mu_{1,1}+1$ in box $(1,\mu_{1,1}+1)$ of the left tableau. By the above reasonning, we get that $u_{\tau_{\tilde{\T}}(i)}=\beta$, therefore $\alpha^{\mu_{1,1}+r_{i,\tilde{\T}}-c_{i,\tilde{\T}}}=1$. For any standard double-tableau associated to $\mu$ and for any number $i$ in a box of its left tableau, we have $c_i-r_i \leq \mu_{1,1}-1< \mu_{1,1}$, therefore $0 <  \mu_{1,1}+r_{i,\tilde{\T}}-c_{i,\tilde{\T}} \leq 2n < 2a$, where $a$ is the order of $\alpha$. It follows that $\mu_{1,1}+r_{i,\tilde{\T}}-c_{i,\tilde{\T}}=a$. We have in the same way that $r_{i,\tilde{\T}}-c_{i,\tilde{\T}} \leq m_\mu-1$, where $m_\mu$ is the number of boxes in the first column of $\mu_1$. Recall that by assumption, we have $a>n$. It follows that $n < a=\mu_{1,1}+r_{i,\tilde{\T}}-c_{i,\tilde{\T}} \leq \mu_{1,1}+m_\mu-1$. We have that $\mu_{1,1}+m_\mu-1$ is equal to the number of boxes to go through horizontally and then vertically in order to get from box  $(1,\mu_{1,1})$ to box $(m_\mu,1)$, therefore $\mu_{1,1}+m_\mu-1\leq n$. This implies that $n<n$ and we get by contradiction that $\mu_{1,1}\geq \lambda_{1,1}$. In the same way, we get $\lambda_{1,1} \geq \mu_{1,1}$, therefore $\lambda_{1,1} =\mu_{1,1}$.

\smallskip

Assume by contradiction that $m_\mu < m_\lambda$. We write $i=m_\mu+1$. There exists $\tilde{T}\in \mu$ such that $u_{\tau_{\tilde{\T}}(i)}\alpha^{c_{i,\tilde{\T}}-r_{i,\tilde{\T}}}=\beta\alpha^{1-i}$ since there exists a standard double-tableau associated to $\lambda$ with $m_\mu+1$ in box $(m_\mu+1,1)$. It follows that $\alpha^{c_{i,\tilde{\T}}-r_{i,\tilde{\T}}+m_\mu}=1$. We then get that $0< c_{i,\tilde{\T}}-r_{i,\tilde{\T}}+m_\mu <2n<2a$ since $c_{i,\tilde{\T}}-r_{i,\tilde{\T}} \geq 1-m_\mu > -m_\mu$. It follows that 
$n<a=c_{i,\tilde{\T}}-r_{i,\tilde{\T}}+m_\mu\leq \mu_{1,1}+m_\mu-1\leq n$. This is a contradiction, therefore $m_\lambda\leq m_\mu$ and in the same way $m_\mu\leq m_\lambda$, therefore $m_\mu= m_\lambda$. We can now complete the proof of the equalities of the quantities.

\smallskip

Let $\T\in \lambda$, $\tilde{\T}\in \mu$ such that for all $i\in [\![1,n]\!]$,  $u_{\tau_{\T}(i)}\alpha^{c_{i,\T}-r_{i,\T}+i-1}=u_{\tau_{\tilde{\T}}(i)}\alpha^{c_{i,\tilde{\T}}-r_{i,\tilde{\T}}+i-1}$. Let  $i\in [\![1,n]\!]$, we have that $u_{\tau_{\T}(i)}=u_{\tau_{\tilde{\T}}(i)}$ and $\alpha^{c_{i,\T}-r_{i,\T}+r_{i,\tilde{\T}}-c_{i,\tilde{\T}}}=1$. We have that $1-m_\lambda \leq c_{i,\T}-r_{i,\T} \leq \lambda_{1,1}-1$ and $1-m_\mu \leq c_{i,\tilde{\T}}-r_{i,\tilde{\T}} \leq \mu_{1,1}-1$ so
$1-n\leq 2-m_\lambda-\lambda_{1,1}= 2-m_\lambda-\mu_{1,1}\leq c_{i,\T}-r_{i,\T}+r_{i,\tilde{\T}}-c_{i,\tilde{\T}} \leq \lambda_{1,1}+m_\mu
-2=\lambda_{1,1}+m_\lambda-2\leq n-1$. It follows by the assumptions on the order of $\alpha$ that $c_{i,\T}-r_{i,\T}=c_{i,\tilde{\T}}-r_{i,\tilde{\T}}$.

\bigskip

We have proved for all $i\in [\![1,n]\!]$ and $\gamma\in \Z$ that $a(i,\lambda,\beta,\gamma)=a(i,\mu,\beta,\gamma)$ and  $\lambda_{1,1}=\mu_{1,1}$. We will prove that this implies that $\lambda_1=\mu_1$.

Let $(\lambda_{1,l_i},l_i)_{1\leq i\leq s_\lambda}$ (resp $(\mu_{1,k_i},k_i)_{1\leq i \leq s_\mu}$) be the extremal boxes of $\lambda_1$ (resp $\mu_1$). For all $i,j\in [\![1,s_\lambda]\!]$ (resp $[\![1,s_\mu]\!]$), $i<j$, we have that $\lambda_{1,l_i}-l_i > \lambda_{1,l_j}-l_j$ (resp $\mu_{1,k_i}-k_i >\mu_{1,k_j}-k_j$). We have that for all $i\in [\![1,s_\lambda]\!]$, $0\neq a(n,\lambda,\beta,\lambda_{1,l_i}-l_i)=a(n,\mu,\beta,\lambda_{1,l_i}-l_i)$ and for all $i\in [\![1,s_\mu]\!]$, $0\neq a(n,\mu,\beta,\mu_{1,k_i}-k_i)=a(n,\lambda,\beta,\mu_{1,k_i}-k_i)$ since for any extremal box, there exists a tableau with $n$ in that box. It follows that $s_\lambda=s_\mu$ and for all $i\in [\![1,s_\lambda]\!], \lambda_{1,l_i}-l_i=\mu_{1,k_i}-k_i$. Since $\lambda_{1,l_1}=\mu_{1,k_1}$, we have that $k_1=l_1$.

Let us show by induction that for all $i\in [\![1,s_\lambda]\!]$, we have $l_i=k_i$ and $\lambda_{1,l_i}=\mu_{1,k_i}$. Let $j\geq 2$, assume that for all $i\leq j-1$, $(\lambda_{1,l_i},l_i)=(\mu_{1,k_i},k_i)$. Assume by contradiction that $\mu_{1,k_j}<\lambda_{1,l_j}$. We then have $k_j<l_j$ since $\mu_{1,k_j}-k_j=\lambda_{1,l_j}-l_j$. We have that $a(n-(k_j-l_{j-1}),\mu,\beta,\mu_{1,k_j}-l_{j-1})=0$ since there are $k_j-l_{j-1}$ boxes below $(l_{j-1},\mu_{1,k_j})$ and $l_j-l_{j-1}$ boxes to its right and no boxes both below and to its right. On the other hand, there exist standard double-tableaux associated to $\lambda$ such that $n-(k_j-l_{j-1})$ is in box $(l_j-k_j+l_{j-1},\lambda_{1,l_j})$ because $l_j> k_j$. It follows that $a(n-(k_j-l_{j-1}),\lambda,\beta,\lambda_{1,l_j}-(l_j-k_j+l_{j-1}))\neq 0$. We have that $\lambda_{1,l_j}-(l_j-k_j+l_{j-1}=\mu_{1,k_j}-k_j+k_j-l_{j-1}=\mu_{1,k_j}-l_{j-1}$, therefore we get a contradiction. It follows that $\mu_{1,k_j} \geq \lambda_{1,l_j}$ and in the same way, $\lambda_{1,l_j}\geq \mu_{1,k_j}$, therefore $\lambda_{1,l_j}=\mu_{1,k_j}$. We then have that $k_j=l_j$ because $\lambda_{1,l_j}-l_j=\mu_{1,k_j}-k_j$.

 We can then conclude by induction that $\lambda_1$ have $\mu_1$ the same extremal boxes. This proves that $\lambda_1=\mu_1$ because a partition is completely determined by its extremal boxes. We get in the same way that $\lambda_2=\mu_2$, therefore $\lambda=\mu$.
\end{proof}

\bigskip

We now generalize the work done in \cite{BMM} for type $A$ to type $B$, that is we define a bilinear form which is fixed by the image of the derived subgroup of the Artin group of the Iwahori-Hecke algebra. We  define a weight by $\omega(\mathbb{T}) =\omega_1(\mathbb{T})\omega_2(\mathbb{T})\omega_3(\mathbb{T})$, where $\omega_1(\mathbb{T})= \underset{ i < j, i\in \mathbb{T}_1, j \in \mathbb{T}_2}\prod (-1)$, $\omega_2(\mathbb{T})=\underset{ i<j, i,j \in \mathbb{T}_2, r_i(\mathbb{T})>r_j(\mathbb{T})}\prod (-1)$ and $\omega_3(\mathbb{T})=\underset{ i<j, i,j \in \mathbb{T}_1, r_i(\mathbb{T})>r_j(\mathbb{T})}\prod(-1)$. We then define a bilinear form $(.|.)$ by $(\T|\tilde{\T}) = \omega(\T)\delta_{\T',\tilde{\T}}$, where $\T' = (\T_2',\T_1')$ for $\T=(\T_1,\T_2)$. In the same way $\lambda'=(\lambda_2',\lambda_1')$ is the transpose of $\lambda=(\lambda_1,\lambda_2)$.

If $\mu$ is a partition of an integer $m$ with diagonal size $b(\mu) =\max\{ i, \mu_i \geq i\}$, we let $\nu(\mu) =1$ if $\frac{m-b(\mu)}{2}$ is even and $-1$ otherwise. Note that by Lemma $6$ of \cite{M2}, if $\mu=\mu'$ then $\frac{m-b(\mu)}{2}$ is an integer and $\nu(\mu)=(-1)^{\frac{m-b(\mu)}{2}}$. If $\lambda=(\lambda_1,\lambda_2)$ is a double-partition with $\lambda_1$ a partition of $r$ and $\lambda_2$ a partition of $n-r$, we let $\tilde{\nu}(\lambda) = \nu(\lambda_1)\nu(\lambda_2)(-1)^{r(n-r)}$.
 
 We now give a proposition similar to Proposition 3.1 in \cite{BMM}.

\begin{prop}\label{bilin}
For all standard double-tableaux $\T, \tilde{\T}$, we have the following properties.
\begin{enumerate}
\item $(S_i.\T|S_i.\tilde{\T}) = (-\alpha)(\T|\tilde{\T})$ and $(T.\T|T.\tilde{\T}) = (-\beta)(\T|\tilde{\T})$.
\item
For all $b\in \mathcal{A}_{B_n}, (b.\T|b.\tilde{\T}) = (\T|\tilde{\T})$.
\item
The restriction of $(.|.)$ to $V_\lambda$ when $\lambda =\lambda'$ and to $V_\lambda \oplus V_{\lambda'}$ when $\lambda \neq \lambda'$ is non-degenerate.\\
Assume that $\lambda=\lambda'$. Then $(.,.)$ is symmetric on $V_\lambda$ if $\tilde{\nu}(\lambda)=1$ and skew-symmetric otherwise.
Moreover, its Witt index is positive.
\end{enumerate}
\end{prop}

\begin{proof}
Let $\T$ and $\tilde{\T}$ be two double standard-tableaux. Recall that for $\lambda\Vdash n$, $\mathbb{T}=(\mathbb{T}_1,\mathbb{T}_2)\in \lambda$ and $i\in [\![1,n]\!]$, we write $\tau_{\T}(i)=1$ if $i\in \T_1$ and $\tau_{\T}(i)=2$ if $i\in \T_2$. We also recall that $u_1=\beta$ and $u_2=-1$.

\textbf{1.} We have $(T.\T|T.\tilde{\T})\neq 0\Leftrightarrow (\T|\tilde{\T})\neq 0 \Leftrightarrow \T'=\tilde{\T} \Rightarrow (T.\T|T.\tilde{\T})=-\beta(\T|\tilde{\T})$ because $\tau_\T(1)=3-\tau_{\T'}(1)$. Let $i\in [\![1,n-1]\!]$, if $\tau_{\T}(i)=\tau_{\T}(i+1)$ then by \cite[Prop 2.4.]{BMM}, we have $(S_i.\T|S_i.\tilde{\T})=(\T|\tilde{\T})$ because in the same way, we have $\omega(\T)=-\omega(\T_{i \leftrightarrow i+1})$ for any standard double-tableau $\T$ and $m_i(\T)=m_i(\T_{\tau_{\T}(i)})$ when $\tau_{\T}(i)=\tau_{\T}(i+1)$. 

 We now assume that $\tau_{\T}(i)\neq \tau_{\T}(i+1)$. We have that $S_i.\T = m_i(\T)\T+(1+m_i(\T))\T_{i\leftrightarrow i+1}$ and $S_i.\tilde{\T}=m_i(\tilde{\T})\tilde{\T}+(1+m_i(\tilde{\T}))\tilde{\T}_{i\leftrightarrow i+1}$. It follows that

\begin{eqnarray*}
(S_i.\T|S_i.\tilde{\T}) & = & m_i(\T)m_i(\tilde{\T})(\T|\tilde{\T})+(1+m_i(\T))m_i(\tilde{\T})(\T_{i\leftrightarrow i+1},\tilde{\T}) + \\
 &  &   m_i(\T)(1+m_i(\tilde{\T}))(\T|\tilde{\T}_{i \leftrightarrow i+1})+(1+m_i(\T))(1+m_i(\tilde{\T}))(\T_{i\leftrightarrow i+1}|\tilde{\T}_{i\leftrightarrow i+1}).
 \end{eqnarray*}
 
This is non-zero only if $\tilde{\T}=\T'$ or $\tilde{\T}=\T_{i\leftrightarrow i+1}'$. We now have two possible cases.

The first case is $\tilde{\T}=\T'$. We write $a=c_i-r_i+r_{i+1}-c_{i+1}$ with $(r_i,c_i)$ and $(r_{i+1},c_{i+1})$ the boxes in $\T$, we then have
\begin{eqnarray*}
(S_i.\T|S_i.\tilde{\T}) & = & m_i(\T)m_i(\T')\omega(\T)+(1+m_i(\T))(1+m_i(\T'))\omega(\T_{i\leftrightarrow i+1})\\
& = & -\omega(\T)(1+m_i(\T)+m_i(\T'))\\
& = & -\omega(\T)\left(1+\frac{\alpha-1}{1+\frac{u_{\tau(i)}}{u_{\tau(i+1)}}\alpha^{a}}+\frac{\alpha-1}{1+\frac{u_{\tau(i+1)}}{u_{\tau(i)}}\alpha^{-a}}\right)\\
& = & -\omega(\T)1+(\alpha-1)\left(\frac{1+\frac{u_{\tau(i)}}{u_{\tau(i+1)}}\alpha^{a}+1+\frac{u_{\tau(i+1)}}{u_{\tau(i)}}\alpha^{-a}}{1+\frac{u_{\tau(i)}}{u_{\tau(i+1)}}\alpha^{a}+\frac{u_{\tau(i+1)}}{u_{\tau(i)}}\alpha^{-a}+1}\right)\\
& = & -\alpha(\T|\tilde{\T}).
\end{eqnarray*}

The second case is $\tilde{\T}=\T_{i\leftrightarrow i+1}'$, we then have
\begin{eqnarray*}
(S_i.\T|S_i.\tilde{\T}) & = & (1+m_i(\T))(m_i(\T_{i\leftrightarrow i+1}'))\omega(\T_{i\leftrightarrow i+1})+m_i(\T)(1+m_i(\T_{i\leftrightarrow i+1}'))\omega(\T)\\
& = & -\omega(\T)(m_i(\T_{i\leftrightarrow i+1}')-m_i(\T))\\
& = & -\omega(\T)\left(\frac{\alpha-1}{1+\frac{u_{\tau(i)}\alpha^{r_{i+1}-c_{i+1}}}{u_{\tau(i+1)}\alpha^{r_i-c_i}}}-\frac{\alpha-1}{1+\frac{u_{\tau(i)}\alpha^{c_i-r_i}}{u_{\tau(i+1)}\alpha^{c_{i+1}-r_{i+1}}}}\right)\\
& = & 0\\
& = & -\alpha(\T|\tilde{\T}).
\end{eqnarray*}
This concludes the proof of \textbf{1}. \textbf{2} follows from \textbf{1}.

\textbf{3.} By definition, the bilinear form is non-degenerate. Assume $\lambda=(\lambda_1,\lambda_2)=\lambda'$. We consider $\T=(\T_1,\T_2)\in \lambda$. Since substituting $i_1<i_2<...<i_l$ in $\T_1$ by $1<2<...<l$ does not change the product and the  weight $\omega$ on $\T\in \mu$ for $\mu\Vdash n$ satisfies $\omega(\T)\omega(\T')=\nu(\mu)$ by Lemme 6 of \cite{M2}, we have that $\omega(\T)\omega(\T')=\nu(\lambda_1)\nu(\lambda_2)\underset{i<j, i\in \T_1,j\in \T_2 ~\mbox{or}~ i\in \T_2,j\in \T_1}\prod (-1)$. The cardinality of the set  $\{i<j,i\in \T_1,j\in \T_2 ~\mbox{or}~ i\in \T_2,j\in \T_1\}$ is equal to the number of pairs $(i,j)$ with $i$ in $\T_1$ and $j$ in $\T_2$, which equals $(\frac{n}{2})^2$. It follows that $\omega(\T)\omega(\T') =\tilde{\nu}(\lambda)$ for any standard double-tableau $\T$ associated to $\lambda$. For any pair $(\T,\tilde{\T})$, we have that
 $$(\tilde{\T}|\T)=\omega(\tilde{\T})\delta_{\T,\tilde{\T}'}=\tilde{\nu}(\lambda)\omega(\tilde{\T}')\delta_{\T,\tilde{\T}'}=\tilde{\nu}(\lambda)\omega(\T)\delta_{\tilde{\T},\T'}=\tilde{\nu}(\lambda)(\T|\tilde{\T}).$$
 The Witt index is positive since the basis can be partitioned in pairs $(\T|\T')$.
 \end{proof}

We remark that we have proved that $\omega(\T)\omega(\T')=\tilde{\nu}(\lambda)$ for any double-partition $\lambda$ and any standard double-tableau $\T$ in $V_\lambda$.

\section{Factorization of the image of the Artin group in the Iwahori-Hecke algebra}\label{FactorB}

In this section, we first see how the different representations are related when restricted to the derived subgroup. We see some factorisations appear between the representation associated to a given double-partition $\lambda$ and the representation associated to its transposed double-partition in Proposition \ref{transpose}. We also see that when the field extension is non-trivial, we have factorizations through field automorphism. All those results are summarized in Proposition \ref{isomorphisme}. Finally, we see that all the representations associated to hook diagrams factor through two representations as is shown in Proposition \ref{patate}.

\subsection{Isomorphisms between representations}

Let $\mathcal{L}\in End(V)$ be defined by  $\mathbb{T}\mapsto \omega(\T)\T'$. We give a generalization of Lemma 3.2. of \cite{BMM}.

\begin{prop}\label{transpose}
Let $\lambda$ be a double-partition of $n$ such that $\lambda \neq \lambda'$ (resp $\lambda=\lambda'$). The map $\mathcal{L}$ induces an endomorphism of $V_{\lambda} \oplus V_{\lambda'}$ (resp $V_\lambda$) which switches $V_{\lambda}$ and $V_{\lambda'}$ (resp leaves $V_\lambda$ stable) such that the actions of $S_r$ and $T$ satisfy
$$\mathcal{L}S_r\mathcal{L}^{-1}(-\alpha)^{-1} = {}^t\!S_r^{-1}, \mathcal{L}T\mathcal{L}^{-1}(-\beta)^{-1} = {}^t\!T^{-1}.$$
\end{prop}

\begin{proof}
This follows directly from Proposition \ref{bilin} by writing the matrix of the bilinear form and the matrix of $\mathcal{L}$. 
\end{proof}
We now suppose $\F_p(\alpha,\beta)= \F_p(\alpha+\alpha^{-1},\beta) =\F_p(\alpha,\beta+\beta^{-1})\neq \F_p(\alpha+\alpha^{-1}, \beta+\beta^{-1})$.

We then have an $\F_q$-automorphism $\epsilon$ of order $2$ such that $\epsilon(\alpha)=\alpha^{-1}$ and $\epsilon(\beta)=\beta^{-1}$. In \cite{BMM}, a function $d$ was defined on any standard tableau $\T$ associated to a partition $\lambda$ of $n$ by $d(\T)=\underset{i,j, r_i>r_j}\prod \frac{\alpha^{c_j-r_j}-\alpha^{c_i-r_i+1}}{\alpha^{c_j-r_j+1}-\alpha^{c_i-r_i}}$, where for $i\in [\![1,n]\!]$, $r_i$ denotes the row of $i$ in $\T$ and $c_i$ denotes the column of $i$ in $\T$.

 Let $\langle .,.\rangle $ be the hermitian form defined by $\langle\mathbb{T},\tilde{\mathbb{T}}\rangle = d(\mathbb{T})\delta_{\mathbb{T},\tilde{\mathbb{T}}}$, where $$d(\mathbb{T})=\tilde{d}(\mathbb{T}_1)\tilde{d}(\mathbb{T}_2)\underset{i\in \T_1, j\in \mathbb{T}_2, i<j}\prod \frac{2+\beta\alpha^{a_{i,j}-1}+\beta^{-1}\alpha^{1-a_{i,j}}}{\alpha+\alpha^{-1}+\beta\alpha^{a_{i,j}}+\beta^{-1}\alpha^{-a_{i,j}}}$$
with $a_{i,j} = c_i-r_i+r_j-c_j$ and $\tilde{d}$ induced by the $d$ defined in  \cite{BMM} applied to $\mathbb{T}_1$ and $\mathbb{T}_2$ by seeing them as standard tableaux using the ordered bijections onto $[\![1,r]\!]$ and $[\![1,n-r]\!]$.

 We now check that $d(\T)$ is well-defined and non-zero for any standard double-tableau. We prove in what follows that the big product in the expression of $d$ is indeed well-defined and non-zero for any double-tableau with no empty components.
 
 Let $\lambda \Vdash n, \T=(\T_1,\T_2)\in \lambda$ and $(i,j)$ a pair of integers such that $i<j, i\in \T_1$ and $j\in \T_2$. We set $r$ to be the number of boxes of $\T_1$, we have $1-r\leq c_i-r_i\leq r-1$ and $1-(n-r)\leq c_j-r_j \leq n-r-1$, therefore $2-n \leq a=a_{i,j}\leq n-2$. We have $\alpha+\alpha^{-1}+\beta\alpha^{a}+\beta^{-1}\alpha^{-a} = \alpha(1+\beta\alpha^{a-1})+\alpha^{-a}\beta^{-1}(1+\beta\alpha^{a-1}) =\alpha(1+\beta\alpha^{a-1})(1+\alpha^{-a-1}\beta^{-1}).$ This product never cancels because $\beta\notin \{-\alpha^i,1-n\leq i \leq n-1\}$. In the same way $2+\beta\alpha^{a-1}+\beta^{-1}\alpha^{1-a}=(1+\beta\alpha^{a-1})(1+\beta^{-1}\alpha^{1-a})$ never cancels.  
 
 Now we have defined this hermitian form, we can generalize Proposition $3.6$ from \cite{BMM}.

\begin{prop}\label{unitary}
The group $A_{B_n}$ acts in a unitary way on $V$ with respect to this hermitian form and this form is non-degenerate on $V_\lambda$ for any double-partition $\lambda$ of $n$. In particular, for any double-partition $\lambda$ of $n$, there exists a matrix $P\in GL_{n_\lambda}(q)$ such that $PR_{\lambda}(T)P^{-1}=\epsilon(R_{\lambda}^\star(T))=R_{\lambda}(T)$ and $PR_{\lambda}(S_r)P^{-1}=\epsilon(R_{\lambda}^{\star}(S_r))$.
\end{prop}

\begin{proof}
The action of $T$ is indeed unitary with regards to this hermitian form  because $\beta\epsilon(\beta)=(-1)\epsilon(-1)=1$. Let $\T$ be a standard double-tableau and $r\in [\![1,n-1]\!]$. If $\tau_{\T}(r)=\tau_{\T}(r+1)$ then the result is a consequence of Proposition 3.6. in \cite{BMM}.

We now assume that $\tau_{\T}(r)\neq\tau_{\T}(r+1)$, up to switching $\T$ and $\T_{r\leftrightarrow r+1}$ we can assume that $\tau_{\T}(r)=1$ and $\tau_{\T}(r+1)=2$.  It remains to show that $\langle\T,\T\rangle=\langle S_r.\T,S_r.\T\rangle, \langle\T_{r\leftrightarrow r+1},\T_{r\leftrightarrow r+1}\rangle = \langle S_r.\T_{r\leftrightarrow r+1},S_r.\T_{r\leftrightarrow r+1}\rangle$ and  $\langle S_r.\T,S_r.\T_{r\leftrightarrow r+1}\rangle =\langle \T,\T_{r\leftrightarrow r+1}\rangle $. In the following computation, we write $a=a_{i,i+1}$ and $\tilde{\T}=\T_{r\leftrightarrow r+1}$. We have
\begin{footnotesize}
\begin{eqnarray*}
 \langle S_r.\T,S_r.\T\rangle & = & m_r(\T)\epsilon(m_r(\T))d(\T)+(1+m_r(\T))\epsilon(1+m_r(\T))d(\tilde{\T})\\
 & = & d(\T)(m_r(\T)\epsilon(m_r(\T))+\left(\frac{\alpha+\alpha^{-1}+\beta\alpha^{a}+\beta^{-1}\alpha^{-a}}{2+\beta\alpha^{a-1}+\beta^{-1}\alpha^{1-a}}(1+m_r(\T))\epsilon(1+m_r(\T)))\right)\\
 & = & d(\T)\left(\frac{\alpha-1}{1+\beta\alpha^{a}}\frac{\alpha^{-1}-1}{1+\beta^{-1}\alpha^{-a}}+\frac{\alpha+\alpha^{-1}+\beta\alpha^{a}+\beta^{-1}\alpha^{-a}}{2+\beta\alpha^{a-1}+\beta^{-1}\alpha^{1-a}}\frac{\alpha+\beta\alpha^{a}}{1+\beta\alpha^{a}}\frac{\alpha^{-1}+\beta^{-1}\alpha^{-a}}{1+\beta^{-1}\alpha^{-a}}\right)\\
 & = & d(\T)\left(\frac{2-\alpha-\alpha^{-1}}{2+\beta\alpha^a+\beta^{-1}\alpha^{-a}}+\frac{\alpha+\alpha^{-1}+\beta\alpha^{a}+\beta^{-1}\alpha^{-a}}{2+\beta\alpha^a+\beta^{-1}\alpha^{-a}}\right)\\
 & = & d(\T)\\
 & = & \langle\T,\T\rangle.
\end{eqnarray*}
\end{footnotesize}
We also have
\begin{footnotesize}
\begin{eqnarray*}
\langle S_r.\tilde{\T},S_r.\tilde{\T}\rangle & = & m_r(\tilde{\T})\epsilon(m_r(\tilde{\T}))d(\tilde{\T})+(1+m_r(\tilde{\T}))\epsilon(1+m_r(\tilde{\T}))d(\T)\\
& = & d(\tilde{\T})\left(\frac{\alpha-1}{1+\beta^{-1}\alpha^{-a}}\frac{\alpha^{-1}-1}{1+\beta\alpha^a}+\frac{2+\beta\alpha^{a-1}+\beta^{-1}\alpha^{1-a}}{\alpha+\alpha^{-1}+\beta\alpha^a+\beta^{-1}\alpha^{-a}}\frac{\alpha+\beta^{-1}\alpha^{-a}}{1+\beta^{-1}\alpha^{-a}}\frac{\alpha^{-1}+\beta\alpha^a}{1+\beta\alpha^a}\right)\\
& = & d(\tilde{\T})\left(\frac{2-\alpha-\alpha^{-1}}{2+\beta\alpha^a+\beta^{-1}\alpha^{-a}}+\right.\\
& &\left. \frac{4+2\beta\alpha^{a+1}+2\beta^{-1}\alpha^{-a-1}+2\beta\alpha^{a-1}+\beta^2\alpha^{2a}+\alpha^{-2}+2\beta^{-1}\alpha^{1-a}+\alpha^2+\beta^{-2}\alpha^{-2a}}{(\alpha+\alpha^{-1}+\beta\alpha^a+\beta^{-1}\alpha^{-a})(2+\beta\alpha^a+\beta^{-1}\alpha^{-a})}\right)\\
& = & d(\tilde{\T})\left(\frac{2-\alpha-\alpha^{-1}}{2+\beta\alpha^a+\beta^{-1}\alpha^{-a}}+\frac{(\alpha+\alpha^{-1}+\beta\alpha^a+\beta^{-1}\alpha^{-a})^2}{(\alpha+\alpha^{-1}+\beta\alpha^a+\beta^{-1}\alpha^{-a})(2+\beta\alpha^a+\beta^{-1}\alpha^{-a})}\right)\\
& = & d(\tilde{\T})\\
& = & \langle\tilde{\T},\tilde{\T}\rangle.
\end{eqnarray*}
\end{footnotesize}
Finally, we have
\begin{footnotesize}
\begin{eqnarray*}
\langle S_r.\T,S_r.\tilde{\T}\rangle  & = & m_r(\T)\epsilon(1+m_r(\tilde{\T}))d(\T)+(1+m_r(\T))\epsilon(m_r(\tilde{\T}))d(\tilde{\T})\\
& = & d(\T)\left(\frac{\alpha-1}{1+\beta\alpha^a}\frac{\alpha^{-1}+\beta\alpha^a}{1+\beta\alpha^a}+\frac{\alpha+\alpha^{-1}+\beta\alpha^{a}+\beta^{-1}\alpha^{-a}}{2+\beta\alpha^{a-1}+\beta^{-1}\alpha^{1-a}}\frac{\alpha+\beta\alpha^a}{1+\beta\alpha^a}\frac{\alpha^{-1}-1}{1+\beta\alpha^a}\right)\\
& = &d(\T)\frac{\alpha-1}{(1+\beta\alpha^a)^2(2+\beta\alpha^{a-1}+\beta^{-1}\alpha^{1-a})}\left((\alpha^{-1}+\beta\alpha^a)(2+\beta\alpha^{a-1}+\beta^{-1}\alpha^{1-a})-\right.\\
& &\left. \alpha^{-1}(\alpha+\alpha^{-1}+\beta\alpha^a+\beta^{-1}\alpha^{-a})(\alpha+\beta\alpha^a)\right)\\
& = & d(\T)\frac{\alpha-1}{(1+\beta\alpha^a)^2(2+\beta\alpha^{a-1}+\beta^{-1}\alpha^{1-a})}\left(2\alpha^{-1}+\beta\alpha^{a-2}+\beta^{-1}\alpha^{-a}+2\beta\alpha^a+\right.\\
& & \left.\beta^2\alpha^{2a-1}+\alpha-\alpha-\alpha^{-1} -\beta\alpha^a-\beta^{-1}\alpha^{-a}-\beta\alpha^a-\beta\alpha^{a-2}-\beta^{2}\alpha^{2a-1}-\alpha^{-1}\right)\\
& = & 0\\
& = & \langle \T,\tilde{\T}\rangle.
\end{eqnarray*}
\end{footnotesize}
\end{proof}

\bigskip
We recall that $\mathcal{A}_{B_n} = [A_{B_n},A_{B_n}]$ is the derived subgroup of $A_{B_n}$. When it exists, we write $\epsilon$ for the automorphism of order $2$ of $\F_q= \F_p(\alpha,\beta)$.

\begin{lemme}\label{Lincoln}
If $\lambda$ is a double-partition of $n$, then the restriction of  $R_{\lambda}$ to $\mathcal{A}_{B_n}$ is absolutely irreducible.
\end{lemme}

\begin{proof}
Assume first it is true for $n=2$. Since $A_{B_n}$ is generated by $A_{B_{n-1}}$ and $\mathcal{A}_{B_n}$, we have the result for $n\geq 3$ by the same method as in the Lemma 3.4(i) of \cite{BMM}. 

We now show the result is true for $n=2$. We only have to show it for $([1],[1])$ since the other representations are $1$-dimensional. We will show in Section \ref{surjectivitism} (Lemmas \ref{platypus} and \ref{platypus456}) that $R_{[1],[1]}(\mathcal{A}_{B_2})\simeq SL_2(q')$ for some $q'$. The irreducibility then follows.
\end{proof}

We now show a lemma computing the normal closure of $\mathcal{A}_{B_n}$. This is a generalization to type $B$ of Lemma $2.1$ of \cite{BMM}

\begin{lemme}\label{normclosBn}
For $n\geq 4$, the normal closure $\ll \mathcal{A}_{B_{n-1}}\gg_{\mathcal{A}_{B_n}}$ of $\mathcal{A}_{B_{n-1}}$ in $\mathcal{A}_{B_n}$ is $\mathcal{A}_{B_n}$.
\end{lemme}

\begin{proof}
Let $n\geq 4$. By Lemma $2.1$ of \cite{BMM}, we have that $\mathcal{A}_{A_n}=\ll \mathcal{A}_{A_{n-1}}\gg_{\mathcal{A}_{A_n}}$, where \\$A_{A_n}=<S_1,S_2,\dots, S_{n-1}>\leq A_{B_n}$. We have that $\mathcal{A}_{B_n}$ is generated by $\mathcal{A}_{A_n}$ and $\mathcal{A}_{B_{n-1}}$ therefore the result follows.
\end{proof}

We now recall Lemma $2.2$ of \cite{BMM}.
\begin{lemme}\label{abel}
Let $G$ be a group, $k$ a field and $R_1,R_2$ two representations of $G$ in $GL_N(k)$ such that the restrictions to the derived subgroup of $G$ are equal and the restriction of at least one of them is absolutely irreducible. There exists a character $\eta : G \rightarrow k^\star$ such that $R_2 = R_1 \otimes \eta$.
\end{lemme}

\begin{prop}
If $\lambda_1$ a partition of $n$ then $R_{(\lambda_1,\emptyset)|\mathcal{A}_{B_n}}\simeq R_{(\emptyset,\lambda_1)|\mathcal{A}_{B_n}}$.
\end{prop}

\begin{proof}
The action of $T$ is diagonal and the action of $S_i$ on $(\T_1,\emptyset)$ is identical to the one on $(\emptyset,\T_1)$, therefore the proof of the result is straightforward.
\end{proof}

We now give a proposition stating the different possible factorizations between double-partitions of $n$. In Proposition $2.7$ of \cite{IH2}, a version of this proposition is given in the generic case. Note that the different possible field extensions in the finite case yield many additional factorizations.

\begin{prop}\label{isomorphisme}
Let $\lambda$ and $\mu$ be double-partitions of $n$ with no empty components. We then have the following properties.
\begin{enumerate}
\item If $\F_q = \F_p(\alpha,\beta)=\F_p(\alpha+\alpha^{-1},\beta+\beta^{-1})$, then
                \begin{enumerate}
                 \item $R_{\lambda|\mathcal{A}_{B_n}} \simeq R_{\mu|\mathcal{A}_{B_n}} \Leftrightarrow \lambda = \mu$,
                 \item $R_{\lambda|\mathcal{A}_{B_n}} \simeq R_{\mu|\mathcal{A}_{B_n}}^{\star} \Leftrightarrow \lambda = \mu'$.
                 \end{enumerate}
\item If $\F_q=\F_p(\alpha,\beta)=\F_p(\alpha+\alpha^{-1},\beta) = \F_p(\alpha,\beta+\beta^{-1}) \neq \F_p(\alpha+\alpha^{-1},\beta+\beta^{-1})$, then 
                \begin{enumerate}
               \item $R_{\lambda|\mathcal{A}_{B_n}} \simeq R_{\mu|\mathcal{A}_{B_n}} \Leftrightarrow \lambda = \mu$,
                 \item $R_{\lambda|\mathcal{A}_{B_n}} \simeq R_{\mu|\mathcal{A}_{B_n}}^{\star} \Leftrightarrow \lambda = \mu'$,
                 \item   $ R_{\lambda|\mathcal{A}_{B_n}} \simeq \epsilon \circ R_{\mu|\mathcal{A}_{B_n}} \Leftrightarrow  \lambda =\mu'$,
                 \item $R_{\lambda|\mathcal{A}_{B_n}} \simeq \epsilon\circ R_{\mu|\mathcal{A}_{B_n}}^{\star} \Leftrightarrow   \lambda=\mu$.
                 \end{enumerate}
\item If $\F_q=\F_p(\alpha,\beta)=\F_p(\alpha,\beta+\beta^{-1}) \neq \F_p(\alpha+\alpha^{-1},\beta) = \F_p(\alpha+\alpha^{-1},\beta+\beta^{-1})$, then 
                \begin{enumerate}
                 \item $R_{\lambda|\mathcal{A}_{B_n}} \simeq R_{\mu|\mathcal{A}_{B_n}} \Leftrightarrow \lambda = \mu$,
                 \item $R_{\lambda|\mathcal{A}_{B_n}} \simeq R_{\mu|\mathcal{A}_{B_n}}^{\star} \Leftrightarrow \lambda = \mu'$,
                 \item   $ R_{\lambda|\mathcal{A}_{B_n}} \simeq \epsilon \circ R_{\mu|\mathcal{A}_{B_n}} \Leftrightarrow  (\lambda_1,\lambda_2) =(\mu_1',\mu_2')$,
                 \item $R_{\lambda|\mathcal{A}_{B_n}} \simeq \epsilon\circ R_{\mu|\mathcal{A}_{B_n}}^{\star} \Leftrightarrow   (\lambda_1,\lambda_2)=(\mu_2,\mu_1)$.
                 \end{enumerate}
\item If $\F_q=\F_p(\alpha,\beta)=\F_p(\alpha+\alpha^{-1},\beta) \neq \F_p(\alpha,\beta+\beta^{-1}) = \F_p(\alpha+\alpha^{-1},\beta+\beta^{-1})$, then 
                \begin{enumerate}
               \item $R_{\lambda|\mathcal{A}_{B_n}} \simeq R_{\mu|\mathcal{A}_{B_n}} \Leftrightarrow \lambda = \mu$,
                 \item $R_{\lambda|\mathcal{A}_{B_n}} \simeq R_{\mu|\mathcal{A}_{B_n}}^{\star} \Leftrightarrow \lambda = \mu'$,
                 \item   $ R_{\lambda|\mathcal{A}_{B_n}} \simeq \epsilon \circ R_{\mu|\mathcal{A}_{B_n}} \Leftrightarrow  (\lambda_1,\lambda_2) =(\mu_2,\mu_1)$,
                 \item $R_{\lambda|\mathcal{A}_{B_n}} \simeq \epsilon\circ R_{\mu|\mathcal{A}_{B_n}}^{\star} \Leftrightarrow   (\lambda_1,\lambda_2)=(\mu_1',\mu_2')$.
                 \end{enumerate}
\end{enumerate}
\end{prop}

\begin{proof}
In all cases, $(a)$ and $(b)$ are the same and the proofs are identical.

(a) By Lemma \ref{abel} and Theorem \ref{mod}, it is sufficient to show that if there exists $\eta : A_{B_n} \rightarrow \F_q^{\star}$ such that $R_\lambda \simeq R_\mu \otimes \eta$, then $\lambda = \mu$. Assume such a character exists, since the abelianization of $A_{B_n}$ is $<\overline{T},\overline{S_1}> \simeq\Z^2$, up to conjugation we have $R_\lambda(b) = R_\mu(b)u^{\ell_1(b)}v^{\ell_2(b)}$ for some $u,v\in \F_q^\star$. Taking the eigenvalues of $S_1$ and $T$ on both sides of the equality we get that $\{\alpha,-1\} = \{u\alpha,-u\}$ and $\{\beta,-1\}=\{v\beta,-v\}$. Since $\alpha^2\neq 1\neq \beta^2$, we have $u=v=1$ which implies that $R_\lambda$ and $R_\mu$ are isomorphic representations. By theorem \ref{mod}, this implies $\lambda=\mu$.

\bigskip

(b) The implication $\lambda=\mu' \Rightarrow R_{\lambda|\mathcal{A}_{B_n}} \simeq R_{\mu|\mathcal{A}_{B_n}}^{\star}$ follows from Proposition \ref{transpose}.

Assume now $R_{\lambda|\mathcal{A}_{B_n}}\simeq R_{\mu|\mathcal{A}_{B_n}}^{\star}$, we then have $R_{\lambda'|\mathcal{A}_{B_n}}\simeq R_{\lambda|\mathcal{A}_{B_n}}^\star \simeq (R_{\mu|\mathcal{A}_{B_n}}^{\star})^{\star} = R_{\mu|\mathcal{A}_{B_n}}$. The result follow from $(a)$. 

In the same way, it is enough to show the converse implication in the remainder of the proof.

\textbf{2. d)} This result follows directly from Proposition \ref{unitary}.

\textbf{2. c)} By  \textbf{2. d)} and \textbf{2. b)}, we have 
$$\epsilon \circ R_{\lambda'|\mathcal{A}_{B_n}} \simeq \epsilon \circ (\epsilon \circ R_{\lambda'|\mathcal{A}_{B_n}}^{\star}) = R_{\lambda'|\mathcal{A}_{B_n}}^{\star} \simeq R_{\lambda|\mathcal{A}_{B_n}}.$$

\textbf{3)} In this case $\epsilon(\alpha) = \alpha^{-1}$ and $\epsilon(\beta) = \beta$.

\textbf{3. c)} For every standard double-tableau $\T = (\T_1,\T_2)$, we define $\tilde{\T}$ by $\tilde{\T}=(\T_1',\T_2')$. Let $\eta : A_{B_n} \rightarrow \F_q^{\star}$ be the character of $A_{B_n}$ defined by $\eta(T) =1$ and $\eta(S_r) =-\alpha$ for all $r$.

Let $Q : V_{(\lambda_1,\lambda_2)} \rightarrow V_{(\lambda_1',\lambda_2')},(\T_1,\T_2)\mapsto (\T_1',\T_2'), U : V_\lambda \rightarrow V_\lambda, \T \mapsto \omega(\T)\T$.

Using the same notations as in Proposition \ref{unitary}, we will show that for all $r\in [\![1,n-1]\!]$ :
$$Q^{-1}(-\alpha)\epsilon(R_{\lambda_1',\lambda_2'}(S_r))Q = U^{-1}PR_{\lambda}(S_r)P^{-1}U, Q^{-1}\epsilon(R_{\lambda_1',\lambda_2'}(T))Q = U^{-1}PR_{\lambda}(T)P^{-1}U.$$

Let $\T= (\T_1,\T_2)$ be a standard double-tableau. The second equality follows from the equalities $PR_{\lambda}(T)P^{-1} =R_{\lambda}(T)$ and $\epsilon(\beta) = \beta$. If $\T_{r\leftrightarrow r+1}$ is non-standard, the first equality is verified by $S_r$. Assume $\T_{r\leftrightarrow r+1}$  is standard, write $a= a_{r,r+1}$. If $\tau_{\T}(r)=\tau_{\T}(r+1)$ then in the basis $(\T,\T_{r\leftrightarrow r+1})$, we have :
$$R_{\lambda}(S_r) =\begin{pmatrix}
 \frac{\alpha-1}{1-\alpha^a} & \frac{\alpha-\alpha^{-a}}{1-\alpha^{-a}} \\
 \frac{\alpha+\alpha^a}{1-\alpha^a} & \frac{\alpha-1}{1-\alpha^{-a}}
\end{pmatrix}, 
 U^{-1}PR_{\lambda}(S_r)P^{-1}U = -\alpha\begin{pmatrix}
\frac{\alpha^{-1}-1}{1-\alpha^a} & \frac{\alpha^{-1}-\alpha^{-a}}{1-\alpha^{-a}}\\
\frac{\alpha^{-1}-\alpha^{a}}{1-\alpha^a} & \frac{\alpha^{-1}-1}{1-\alpha^{-a}}
\end{pmatrix},$$
$$-\alpha Q^{-1}\epsilon(R_{\lambda_1',\lambda_2'}(S_r))Q = -\alpha\epsilon(\begin{pmatrix} \frac{\alpha-1}{1-\alpha^{-a}} & \frac{\alpha-\alpha^a}{1-\alpha^a}\\
\frac{\alpha-\alpha^{-a}}{1-\alpha^{-a}} & \frac{\alpha-1}{1-\alpha^a}\end{pmatrix} = U^{-1}PR_{\lambda}(S_r)P^{-1}U.$$

If $\tau_{\T}(r)=1$ and $\tau_{\T}(r+1)=2$ then we have

$$R_\lambda(S_r) =\begin{pmatrix}
\frac{\alpha-1}{1+\beta\alpha^a} & \frac{\alpha+\beta^{-1}\alpha^{-a}}{1+\beta^{-1}\alpha^{-a}}\\
\frac{\alpha+\beta\alpha^a}{1+\beta\alpha^a} & \frac{\alpha-1}{1+\beta^{-1}\alpha^{-a}}
\end{pmatrix}, U^{-1}PR_{\lambda}(S_r)P^{-1}U = -\alpha\begin{pmatrix}
\frac{\alpha^{-1}-1}{1+\beta\alpha^a} & \frac{\beta^{-1}\alpha^{-a}+\alpha^{-1}}{1+\beta^{-1}\alpha^{-a}}\\
\frac{\beta\alpha^a+\alpha^{-1}}{1+\beta\alpha^a} & \frac{\alpha^{-1}-1}{1+\beta^{-1}\alpha^{-a}}
\end{pmatrix},$$
$$-\alpha Q^{-1}\epsilon(R_{\lambda_1',\lambda_2'})Q = -\alpha\epsilon(\begin{pmatrix} \frac{\alpha-1}{1+\beta\alpha^{-a}} & \frac{\alpha+\beta^{-1}\alpha^a}{1+\beta^{-1}\alpha^a}\\
\frac{\alpha+\beta\alpha^{-a}}{1+\beta\alpha^{-a}} & \frac{\alpha-1}{1+\beta^{-1}\alpha^a}\end{pmatrix} = U^{-1}PR_{\lambda}(S_r)P^{-1}U.$$
It follows that $ R_{(\lambda_1,\lambda_2)|\mathcal{A}_{B_n}} \simeq \epsilon \circ R_{(\lambda_1',\lambda_2')|\mathcal{A}_{B_n}}$.

\textbf{3. d)} This is a consequence of \textbf{3. c)} and \textbf{3. b)}.

\textbf{4. c)} For each standard double tableau $\T=(\T_1,\T_2)$, we set $\tilde{\T}=(\T_2,\T_1)$. Let $\eta : A_{B_n} \rightarrow \F_q^\star$ be the character defined by $\eta(T)=-\beta$ and $\eta(S_r) = 1$ for all $S_r$.

Let $Q : V_{\lambda_2,\lambda_1} \rightarrow V_{\lambda_1,\lambda_2}$ be $\T\mapsto \tilde{\T}$.

Let us show that $R_{(\lambda_1,\lambda_2)}= Q((\epsilon \circ R_{(\lambda_2,\lambda_1)})\bigotimes \eta )Q^{-1}$, where in this case $\epsilon(\alpha)=\alpha$ and $\epsilon(\beta)=\beta^{-1}$. The proof of the result is straightforward for $T$. Let $r\in[\![1,n-1]\!]$. If $r$ and $r+1$ are in the same tableau then the result is clear. Assume $r$ is in the right tableau, $r+1$ is in the left tableau and set $a=a_{r,r+1}$. In the basis $(\T,\T_{r\leftrightarrow r+1})$, we have 
$$R_{\lambda}(S_r) =\begin{pmatrix}
 \frac{\alpha-1}{1+\beta\alpha^a} & \frac{\alpha+\beta^{-1}\alpha^{-a}}{1+\beta^{-1}\alpha^{-a}} \\
 \frac{\alpha+\beta\alpha^a}{1+\beta\alpha^a} & \frac{\alpha-1}{1+\beta^{-1}\alpha^{-a}}
\end{pmatrix}$$
$$Q^{-1}\epsilon(R_{\lambda_2,\lambda_1}(S_r))Q = \epsilon(\begin{pmatrix}
\frac{\alpha-1}{1+\beta^{-1}\alpha^a} & \frac{\alpha+\beta\alpha^{-a}}{1+\beta\alpha^{-a}}\\
\frac{\alpha+\beta^{-1}\alpha^a}{1+\beta^{-1}\alpha^a} & \frac{\alpha-1}{1+\beta\alpha^{-a}} 
\end{pmatrix} = R_{\lambda}(S_r).$$
This proves that we have $ R_{(\lambda_1,\lambda_2)|\mathcal{A}_{B_n}} \simeq \epsilon \circ R_{(\lambda_2,\lambda_1)|\mathcal{A}_{B_n}}$.

\textbf{4. d)} This is a consequence of \textbf{4.c)} and \textbf{4.b)}.
\end{proof}

For $r\in [\![1,n-1]\!]$, we define the double-partitions $\lambda_{(r)}=([1^{n-r}],[r])$ and $\lambda^{(r)} = ([r],[1^{n-r}])$. The following proposition is a generalization of Proposition $3.5$ of \cite{BMM}. The analogous proposition in the generic case is given by Proposition $2.13$ of \cite{IH2}.

\begin{prop}\label{patate}
For $r\in [\![1,n-1]\!]$, there exists $\eta_{1,r},\eta_{2,r} : A_{B_n} \rightarrow \F_q^\star$ such that 

$R_{\lambda_{(r)}} \simeq (\Lambda^rR_{\lambda_{(1)}}) \otimes \eta_{1,r}$ and $R_{\lambda^{(r)}}\simeq (\Lambda^rR_{\lambda^{(1)}}) \otimes \eta_{2,r}$ .\end{prop}

\begin{proof}
Every double-tableau associated with $\lambda_{(r)}$ or $\lambda^{(r)}$ can be mapped in a one-to-one way to a set $\{i_1,i_2,...,i_r\}\subset [\![1,n]\!]$ such that $i_1<i_2<...<i_r$, where $i_k$ is the number in the $k$-th box of the right component. We write $v_I$ the corresponding double-tableau and $v_i=v_{\{i\}}$.\\
After computation, for $k\in [\![1,n-1]\!]$, we get the following properties.
\begin{enumerate}
\item If $1\in I$ then $R_{\lambda_{(r)}}(T)v_I = -v_I$.
\item If $1\notin I$ then $R_{\lambda_{(r)}}(T)v_I = \beta v_I.$
\item If $k,k+1\notin I$ then $R_{\lambda_{(r)}}(S_k)v_I = -v_I$.
\item If $k,k+1\in I$ then $R_{\lambda_{(r)}}(S_k)v_I=\alpha v_I$.
\item If $k\in I, k+1\notin I$ then $R_{\lambda_{(r)}}(S_k)v_I = \frac{\alpha-1}{1+\beta^{-1}\alpha^{k-1}}v_I+\frac{\alpha+\beta^{-1}\alpha^{k-1}}{1+\beta^{-1}\alpha^{k-1}}v_{I\Delta \{k,k+1\}}.$
\item If $k\notin I, k+1\in I$ then $R_{\lambda_{(r)}}(S_k)v_I = \frac{\alpha-1}{1+\beta\alpha^{1-k}}v_I +\frac{\alpha+\beta\alpha^{1-k}}{1+\beta\alpha^{1-k}}v_{I\Delta \{k,k+1\}}.$
\end{enumerate}
Above, $\Delta$ is the symmetric difference : $A\Delta B = (A\cup B)\setminus (A\cap B)$.

To each set $I=\{i_1,i_2,...,i_r\}$ can be given in a one-to-one way an element $u_I$ of $\Lambda^rR_{\lambda_{(1)}}$ writing $u_I = v_{i_1} \wedge v_{i_2} \wedge ... \wedge v_{i_r}$ and these $u_I$ form a basis.

 For $k\in [\![1,n-1]\!]$, we have the following properties.

\begin{enumerate}
\item If $1\in I$ then $\Lambda^rR_{\lambda_{(1)}}(T)u_I = -\beta^{r-1}u_I$.
\item If $1\notin I$ then $\Lambda^rR_{\lambda_{(1)}}(T)u_I = \beta^r u_I.$
\item If $k,k+1\notin I$ then $\Lambda^rR_{\lambda_{(1)}}(S_k)u_I = (-1)^r u_I$.
\item If $k,k+1\in I$ then $\Lambda^rR_{\lambda_{(1)}}(S_k)u_I=(-1)^{r-1}\alpha u_I$.
\item If $k\in I, k+1\notin I$ then 
$$\Lambda^rR_{\lambda_{(1)}}(S_k)u_I = (-1)^{r-1}\frac{\alpha-1}{1+\beta^{-1}\alpha^{k-1}}u_I+(-1)^{-r-1}\frac{\alpha+\beta^{-1}\alpha^{k-1}}{1+\beta^{-1}\alpha^{k-1}}u_{I\Delta \{k,k+1\}}.$$
\item If $k\notin I, k+1\in I$ then $\Lambda^rR_{\lambda_{(1)}}(S_k)u_I = (-1)^{r-1}\frac{\alpha-1}{1+\beta\alpha^{1-k}}u_I +(-1)^{r-1}\frac{\alpha+\beta\alpha^{1-k}}{1+\beta\alpha^{1-k}}u_{I\Delta \{k,k+1\}}.$
\end{enumerate}

Looking at the basis change $v_I \mapsto u_I$ and the character $\eta_{1,r}(h) =(-1)^{(r-1)\ell_1(h)}\beta^{(r-1)\ell_2(h)}$, we have the first part of the proposition. 

In the same way, writing $\eta_{2,r}(h) =(-1)^{(r-1)\ell_1(h)}(-1)^{(r-1)\ell_2(h)}$, we have the second part of the proposition.
\end{proof}

\subsection{Factorization depending on the field}\label{lalala}

In this section, we use the work from the previous sections to state the main results in type $B$. As for the factorizations shown in Proposition \ref{isomorphisme}, we have to distinguish cases depending on the field extensions appearing.

The result depends on the properties of the field extension $\F_q=\F_p(\alpha,\beta)$ of $\F_p(\alpha+\alpha^{-1},\beta+\beta^{-1})$ and the field extension $\F_{\tilde{q}}=\F_p(\alpha)$ of $\F_p(\alpha+\alpha^{-1})$. Let us consider the Hasse diagram for the field extension $\F_p(\alpha,\beta)$ of $\F_p(\alpha+\alpha^{-1},\beta+\beta^{-1})$.

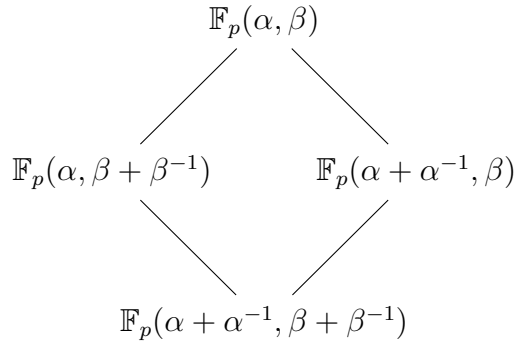
\begin{center}
\begin{tikzpicture}
\node (1) at (0,2){$\F_p(\alpha,\beta)$};
\node (2) at (-2,0){$\F_p(\alpha,\beta+\beta^{-1})$};
\node (3) at (2,0){$\F_p(\alpha+\alpha^{-1},\beta)$};
\node (4) at (0,-2){$\F_p(\alpha+\alpha^{-1},\beta+\beta^{-1})$};
\draw (1) to (2);
\draw (1) to (3);
\draw (2) to (4);
\draw (3) to (4);
\end{tikzpicture}
\end{center}

All the extensions represented by edges are of degree at most $2$ since $X^2-(\alpha+\alpha^{-1})X+1$ and $X^2-(\beta+\beta^{-1})X+1$ are the polynomials involved. Note that if $\F_p(\alpha,\beta+\beta^{-1})=\F_p(\alpha+\alpha^{-1},\beta)$ then $\beta\in \F_p(\alpha,\beta+\beta^{-1})$. Therefore $\F_p(\alpha,\beta+\beta^{-1})=\F_p(\alpha,\beta)$. This proves in particular that $\F_p(\alpha,\beta)$ cannot be an extension of degree $4$ of $\F_p(\alpha+\alpha^{-1},\beta+\beta^{-1})$, otherwise we would have that $\F_p(\alpha,\beta+\beta^{-1})$ and $\F_p(\alpha+\alpha^{-1},\beta)$ are two subfields of degree $2$ of $\F_p(\alpha,\beta)$ and would therefore be equal. This would then imply that they are equal to $\F_p(\alpha,\beta)$ and contradict the fact that the corresponding extensions are of degree $2$. Note that $\F_p(\alpha,\beta+\beta^{-1})\neq \F_p(\alpha+\alpha^{-1},\beta+\beta^{-1})$ implies that we have $\F_p(\alpha)\neq \F_p(\alpha+\alpha^{-1})$. This proves that we have the six following cases to consider

\begin{enumerate}
\item $\F_q=\F_p(\alpha,\beta)= \F_p(\alpha+\alpha^{-1},\beta+\beta^{-1})$ and $\F_p(\alpha)=\F_p(\alpha+\alpha^{-1})$.
\item $\F_q=\F_p(\alpha,\beta) = \F_p(\alpha+\alpha^{-1},\beta+\beta^{-1})$ and $\F_p(\alpha) \neq \F_p(\alpha+\alpha^{-1})$.
\item $\F_q=\F_p(\alpha,\beta)=\F_p(\alpha+\alpha^{-1},\beta)= \F_p(\alpha,\beta+\beta^{-1})\neq \F_p(\alpha+\alpha^{-1},\beta+\beta^{-1})$.
\item $\F_q = \F_p(\alpha,\beta) = \F_p(\alpha,\beta+\beta^{-1})\neq \F_p(\alpha+\alpha^{-1},\beta) = \F_p(\alpha+\alpha^{-1},\beta+\beta^{-1})$.
\item $\F_q = \F_p(\alpha,\beta) =\F_p(\alpha+\alpha^{-1},\beta) \neq \F_p(\alpha,\beta+\beta^{-1}) = \F_p(\alpha+\alpha^{-1},\beta+\beta^{-1})$ and $\F_p(\alpha) \neq \F_p(\alpha+\alpha^{-1})$.
\item $\F_q = \F_p(\alpha,\beta) = \F_p(\alpha+\alpha^{-1},\beta) \neq \F_p(\alpha,\beta+\beta^{-1})=\F_p(\alpha+\alpha^{-1},\beta+\beta^{-1})$ and $\F_p(\alpha)= \F_p(\alpha+\alpha^{-1})$.
\end{enumerate}

We remark that in the third and fourth cases, we have $\F_p(\alpha)\neq \F_p(\alpha+\alpha^{-1})$.

Before stating the main results for type $B$, we recall the two following lemmas, the first one is Lemma $2.4$ of \cite{BM} and the proof of the second one is included in the proof of Proposition $4.1$ of \cite{BMM}.

\begin{lemme}\label{Harinordoquy}
Let $\rho$ be an absolutely irreducible representation of a group $G$ in $GL_r(q)$, where $\F_q$ is a finite field such that there exists an automorphism $\epsilon$ of order $2$ of $\F_q$. If $\rho \simeq \epsilon \circ \rho^\star$, then there exists $S\in GL_r(q)$ such that $S^{-1}\rho(g)S \in GU_r(q^\frac{1}{2})$ for all $g \in G$.
 \end{lemme}

\begin{lemme}\label{Ngwenya}
Let $\rho$ and $G$ be as in the previous lemma. If $\rho \simeq \epsilon \circ \rho$, then there exists $S\in GL_r(q)$ such that $S^{-1}\rho(g)S \in GL_r(q^\frac{1}{2})$ for all $g \in G$.
 \end{lemme}

\bigskip

In certain cases, $(\lambda_1,\lambda_2)$ factorizes through $(\lambda_2,\lambda_1)$ or $(\lambda_1',\lambda_2')$, therefore we need a good order on double-partitions of $n$. We first choose for $r\leq n$ an order on partitions of $r$ such that if $r$ has $2l$ partitions different from their transpose $\{a_i,a_{i'}\}_{i\in [\![ 1,l]\!]}$ and $s$ partitions $\{a_{l+i}\}_{i\in [\![1,s]\!]}$ equal to their transpose then $a_1<a_1'<a_2<a_2'<...<a_l<a_l'<a_{l+1}=a_{l+1}'<...<a_{l+s}=a_{l+s}'$.  We also require $a_1=[r]$. This gives us that $\lambda<\mu$ implies that $\lambda'<\mu'$ whenever $\lambda \neq \mu'$. If $\lambda \Vdash n_1$ and $\mu \Vdash n_2$ then we say $\lambda >\mu$ if $n_1>n_2$ or $n_1=n_2$ and $\lambda>\mu$. We then define the order $<$ on double-partitions of $n$ in the following way, where $\lambda_1$ is a partition of $r_\lambda$ : $(\lambda_1,\lambda_2)<(\mu_1,\mu_2)$ if $r_\lambda< r_\mu$ or ($r_\lambda = r_\mu$ and $\lambda_1< \mu_1$) or ($r_\lambda =r_\mu, \lambda_1=\mu_1$ and $\lambda_2 < \mu_2$).

\begin{lemme}
If $\lambda=(\lambda_1,\lambda_2)$ is a double-partition such that $\lambda \neq \lambda', \lambda \neq (\lambda_2,\lambda_1)$ and $\lambda \neq (\lambda_1',\lambda_2')$ then exactly one element of $\{\lambda,\lambda',(\lambda_1',\lambda_2'),(\lambda_2,\lambda_1)\}$ verifies the property :

 $$ (*) ~ \lambda < \lambda' \mbox{ and } \lambda < (\lambda_1',\lambda_2').$$
\end{lemme}

\begin{proof}
Let $\lambda=(\lambda_1,\lambda_2)$ be a double-partition verifying the conditions in the lemma.
 Assume $\lambda > \lambda'$ and $\lambda < (\lambda_1',\lambda_2')$. Since $\lambda > \lambda'$, we have $r_\lambda \geq r_{\lambda'}$ and since $r_\lambda +r_\lambda' = n$, we get $r_\lambda \geq \frac{n}{2}$. Let's show that either ($\lambda' < \lambda$ and $\lambda' < (\lambda_2,\lambda_1)$) or ($(\lambda_2,\lambda_1) < (\lambda_1',\lambda_2')$ and $(\lambda_2,\lambda_1) < (\lambda_2',\lambda_1')$), i.e. either $\lambda'$ verifies $(*)$ or $(\lambda_2,\lambda_1)$ verifies $(*)$.
 Those two cases are indeed distinct because either $\lambda' < (\lambda_2,\lambda_1)$ or $(\lambda_2,\lambda_1) < \lambda'$. If $\lambda' < (\lambda_2,\lambda_1)$ then we are in the first case because we assumed $\lambda > \lambda'$. Let's now assume $\lambda' > (\lambda_2,\lambda_1)$, we must show $(\lambda_2,\lambda_1)< (\lambda_1',\lambda_2')$. This is obvious if $r_\lambda > \frac{n}{2}$. If $r_\lambda=\frac{n}{2}$, then $(\lambda_1,\lambda_2) =\lambda > \lambda'=(\lambda_2',\lambda_1')$ implies that $\lambda_1> \lambda_2'$ or ($\lambda_1=\lambda_2'$ and $\lambda_2 > \lambda_1'$), which is a contradiction. Therefore $\lambda_1>  \lambda_2'$ and since $\lambda_1 \neq \lambda_2$, this implies $\lambda_1'>\lambda_2$ by definition of our order on partitions of $r_\lambda$. This shows that $(\lambda_2,\lambda_1)< (\lambda_1',\lambda_2')$.

Assume $\lambda> \lambda'$ and $\lambda> (\lambda_1',\lambda_2')$. We then have that either $\lambda'$ verifies $(*)$ or $(\lambda_2,\lambda_1)$ verifies $(*)$ in exactly the same way as in the previous case.

Assume $\lambda < \lambda'$ and $\lambda > (\lambda_1',\lambda_2')$, let us show that $(\lambda_1',\lambda_2') < (\lambda_2,\lambda_1)$ and $(\lambda_1',\lambda_2')< (\lambda_1,\lambda_2)$, i.e. $(\lambda_1',\lambda_2')$ verifies $(*)$. It is enough to show the second inequality since we have the first one by assumption. This is obvious if $r_\lambda < \frac{n}{2}$. If $r_\lambda = \frac{n}{2}$ then $\lambda_1< \lambda_2'$ because $\lambda_1\neq \lambda_2'$, therefore $\lambda_1' < \lambda_2$ because $\lambda_1\neq \lambda_2$ and $(\lambda_1',\lambda_2') < (\lambda_2,\lambda_1)$.

Assume $\lambda < \lambda'$ and $\lambda < (\lambda_1',\lambda_2')$. To conclude the proof, it is enough to show that not one of $\lambda', (\lambda_1',\lambda_2')$ and $(\lambda_2,\lambda_1)$ verifies $(*)$ in this case. It is obvious for $\lambda'$ and $(\lambda_1',\lambda_2)$. If $r_\lambda < \frac{n}{2}$, it is also obvious for $(\lambda_2,\lambda_1)$ since $(\lambda_2,\lambda_1) > (\lambda_1',\lambda_2')$. If $r_\lambda = \frac{n}{2}$ then since $\lambda_1 < \lambda_2'$ and $\lambda_1 \neq \lambda_2'$, we have that $\lambda_2 > \lambda_1'$, therefore $(\lambda_2,\lambda_1) > (\lambda_1',\lambda_2')$.
\end{proof}

\bigskip

We are now able to state the main results for type $B$ which are a generalization of Theorem 1.1 of \cite{BMM}. The end of the proof will be in the next section. The main difference arises from the additional factorizations in the last cases for the field extensions explicited in six different cases at the beginning of Section \ref{FactorB}. The results are given in Theorems \ref{result1} up to \ref{result6}. The analogous versions of these theorems in the generic case can be seen in Theorem $2.21$ of \cite{IH2}. Note here again that in type $B$, the complexity of the field extensions involved yields a wider variety of results.

We write $A_{1,n}=\{(\lambda_1,\emptyset),\lambda_1 \vdash n\}, A_{2,n} =\{(\emptyset,\lambda_2),\lambda_2 \vdash n\}, A_n = A_{1,n} \cup A_{2,n}$. $A\epsilon_n= \{(\lambda_1,\emptyset)\in A_{1,n},~\lambda_1 ~\mbox{not a hook}\}$, $ \epsilon_n=\{\lambda \vdash\vdash n, \lambda \notin A_n, \lambda~\mbox{not a hook}\}, \F_{\tilde{q}}=\F_p(\alpha)$. 
\begin{theo}\label{result1}
If $\F_q=\F_p(\alpha,\beta)=\F_p(\alpha+\alpha^{-1}, \beta+\beta^{-1})$ and $\F_p(\alpha)=\F_p(\alpha+\alpha^{-1})$, then the morphism : $\mathcal{A}_{B_n} \rightarrow \mathcal{H}_{B_n,\alpha,\beta}^\times \simeq \underset{\lambda \vdash\vdash n}\prod GL(\lambda)$ factors through the epimorphism 
$$\Phi_{1,n} : \mathcal{A}_{B_n} \rightarrow SL_{n-1}(\tilde{q})\times \underset{(\lambda_1,\emptyset)\in A\epsilon_n, \lambda_1<\lambda_1'}\prod SL_{n_\lambda}(\tilde{q}) \times \underset{(\lambda_1,\emptyset)\in A\epsilon_n,\lambda_1=\lambda_1'}\prod OSP(\lambda)'\times$$
 $$SL_n(q)^2 \times \underset{\lambda\in \epsilon_n, \lambda < \lambda'}\prod SL_{n_\lambda}(q) \times \underset{\lambda\in \epsilon_n, \lambda=\lambda'}\prod OSP(\lambda)',$$
where $OSP(\lambda)$ is the group of isometries of the bilinear form considered in Proposition \ref{bilin}.
\end{theo}

\bigskip

\begin{theo}\label{result2}
If $\F_q=\F_p(\alpha,\beta)=\F_p(\alpha+\alpha^{-1}, \beta+\beta^{-1})$ and $\F_p(\alpha)\neq \F_p(\alpha+\alpha^{-1})$, then the morphism : $\mathcal{A}_{B_n} \rightarrow \mathcal{H}_{B_n,\alpha,\beta}^\times \simeq \underset{\lambda \vdash\vdash n}\prod GL(\lambda)$ factors through the epimorphism 
$$\Phi_{2,n} : \mathcal{A}_{B_n} \rightarrow SU_{n-1}(\tilde{q}^{\frac{1}{2}})\times \underset{(\lambda_1,\emptyset)\in A\epsilon_n, \lambda_1<\lambda_1'}\prod SU_{n_\lambda}(\tilde{q}^{\frac{1}{2}}) \times \underset{(\lambda_1,\emptyset)\in A\epsilon_n,\lambda_1=\lambda_1'}\prod \widetilde{OSP}(\lambda)'\times$$
 $$SL_n(q)^2 \times \underset{\lambda\in \epsilon_n, \lambda < \lambda'}\prod SL_{n_\lambda}(q) \times \underset{\lambda\in \epsilon_n, \lambda=\lambda'}\prod OSP(\lambda)',$$
where $OSP(\lambda)$ is the group of isometries of the bilinear form considered in Proposition \ref{bilin} and $\tilde{OSP}(\lambda)$ is the group of isometries of a form of the same type but defined over subfield of degree $2$.
\end{theo}

\begin{theo}\label{result3}
If  $\F_q=\F_p(\alpha,\beta)=\F_p(\alpha+\alpha^{-1},\beta)=\F_p(\alpha,\beta+\beta^{-1})\neq \F_p(\alpha+\alpha^{-1}, \beta+\beta^{-1})$, then the morphism $\mathcal{A}_{B_n} \rightarrow \mathcal{H}_{B_n,\alpha,\beta}^\times \simeq \underset{\lambda \vdash\vdash n}\prod GL(\lambda)$ factors through the epimorphism
$$\Phi_{3,n} : \mathcal{A}_{B_n} \rightarrow SU_{n-1}({\tilde{q}^{\frac{1}{2}}})\times \underset{(\lambda_1,\emptyset)\in A\epsilon_n, \lambda_1<\lambda_1'}\prod SU_{n_\lambda}(\tilde{q}^{\frac{1}{2}}) \times \underset{(\lambda_1,\emptyset)\in A\epsilon_n,\lambda_1=\lambda_1'}\prod \widetilde{OSP}(\lambda)'\times $$
 $$SU_n(q^{\frac{1}{2}})^2 \times \underset{\lambda\in \epsilon_n, \lambda < \lambda'}\prod SU_{n_\lambda}(q^{\frac{1}{2}}) \times \underset{\lambda\in \epsilon_n, \lambda=\lambda'}\prod \widetilde{OSP}(\lambda)',$$
where $\widetilde{OSP}(\lambda)$ is the group of isometries of a bilinear form of the same type as the one considered in proposition \ref{bilin} but defined over $\F_{q^\frac{1}{2}}$.
\end{theo}

\bigskip

\bigskip

By Proposition \ref{isomorphisme} and Lemmas \ref{Harinordoquy} and \ref{Ngwenya}, we have  the following theorem 

\begin{theo}\label{result4}
If $\F_q=\F_p(\alpha,\beta)=\F_p(\alpha,\beta+\beta^{-1}) \neq \F_p(\alpha+\alpha^{-1},\beta) = \F_p(\alpha+\alpha^{-1},\beta+\beta^{-1})$, then the morphism $\mathcal{A}_{B_n} \rightarrow \mathcal{H}_{B_n,\alpha,\beta}^\times \simeq \underset{\lambda \vdash\vdash n}\prod GL(\lambda)$ factors through the epimorphism
$$\Phi_{4,n} : \mathcal{A}_{B_n} \rightarrow SU_{n-1}({\tilde{q}^{\frac{1}{2}}})\times \underset{(\lambda_1,\emptyset)\in A\epsilon_n, \lambda_1<\lambda_1'}\prod SU_{n_{\lambda}}(\tilde{q}^{\frac{1}{2}}) \times \underset{(\lambda_1,\emptyset)\in A\epsilon_n,\lambda_1=\lambda_1'}\prod \widetilde{OSP}(\lambda)'\times$$
 $$ SL_n(q) \times \underset{\lambda\in \epsilon_n, \lambda < \lambda', \lambda < (\lambda_1',\lambda_2'), \lambda\neq (\lambda_2,\lambda_1)}\prod SL_{n_\lambda}(q) \times \underset{\lambda \in \epsilon_n, \lambda < \lambda', \lambda=(\lambda_2,\lambda_1)}\prod SU_{n_{\lambda}}(q^{\frac{1}{2}}) \times  $$
 $$  \underset{\lambda \in \epsilon_n, \lambda< \lambda',\lambda=(\lambda_1',\lambda_2')}\prod SL_{n_\lambda}(q^\frac{1}{2})\times\underset{\lambda\in \epsilon_n, \lambda=\lambda', \lambda<(\lambda_1',\lambda_2')}\prod OSP(\lambda)' \times \underset{\lambda \in \epsilon_n, \lambda=\lambda',\lambda=(\lambda_1',\lambda_2')}\prod \widetilde{OSP}(\lambda)' .$$
 \end{theo}
 
\begin{theo}\label{result5}
If $\F_q=\F_p(\alpha,\beta)=\F_p(\alpha+\alpha^{-1},\beta) \neq \F_p(\alpha,\beta+\beta^{-1}) = \F_p(\alpha+\alpha^{-1},\beta+\beta^{-1})$ and $\F_p(\alpha)\neq \F_p(\alpha+\alpha^{-1})$, then the morphism $\mathcal{A}_{B_n} \rightarrow \mathcal{H}_{B_n,\alpha,\beta}^\times \simeq \underset{\lambda \vdash\vdash n}\prod GL(\lambda)$ factors through the epimorphism
$$\Phi_{5,n} : \mathcal{A}_{B_n} \rightarrow SU_{n-1}({\tilde{q}^{\frac{1}{2}}})\times \underset{(\lambda_1,\emptyset)\in A\epsilon_n, \lambda_1<\lambda_1'}\prod SU_{n_{\lambda}}(\tilde{q}^{\frac{1}{2}}) \times \underset{(\lambda_1,\emptyset)\in A\epsilon_n,\lambda_1=\lambda_1'}\prod \widetilde{OSP}(\lambda)'\times$$
 $$ SL_n(q) \times \underset{\lambda\in \epsilon_n, \lambda < \lambda', \lambda < (\lambda_1',\lambda_2'), \lambda\neq (\lambda_2,\lambda_1)}\prod SL_{n_\lambda}(q) \times \underset{\lambda \in \epsilon_n, \lambda < \lambda', \lambda=(\lambda_2,\lambda_1)}\prod SL_{n_{\lambda}}(q^{\frac{1}{2}}) \times  $$
 $$ \underset{\lambda \in \epsilon_n, \lambda< \lambda',\lambda=(\lambda_1',\lambda_2')}\prod SU_{n_\lambda}(q^\frac{1}{2})\times \underset{\lambda\in \epsilon_n, \lambda=\lambda', \lambda<(\lambda_1',\lambda_2')}\prod OSP(\lambda)' \times \underset{\lambda \in \epsilon_n, \lambda=\lambda',\lambda=(\lambda_1',\lambda_2')}\prod \widetilde{OSP}(\lambda)' .$$
 \end{theo}
 
 \begin{theo}\label{result6}
If $\F_q=\F_p(\alpha,\beta)=\F_p(\alpha+\alpha^{-1},\beta) \neq \F_p(\alpha,\beta+\beta^{-1}) = \F_p(\alpha+\alpha^{-1},\beta+\beta^{-1})$ and $\F_p(\alpha)= \F_p(\alpha+\alpha^{-1})$, then the morphism $\mathcal{A}_{B_n} \rightarrow \mathcal{H}_{B_n,\alpha,\beta}^\times \simeq \underset{\lambda \vdash\vdash n}\prod GL(\lambda)$ factors through the epimorphism
$$\Phi_{6,n} : \mathcal{A}_{B_n} \rightarrow SL_{n-1}(\tilde{q})\times \underset{(\lambda_1,\emptyset)\in A\epsilon_n, \lambda_1<\lambda_1'}\prod SL_{n_{\lambda}}(\tilde{q}) \times \underset{(\lambda_1,\emptyset)\in A\epsilon_n,\lambda_1=\lambda_1'}\prod {OSP}(\lambda)'\times$$
 $$ SL_n(q) \times \underset{\lambda\in \epsilon_n, \lambda < \lambda', \lambda < (\lambda_1',\lambda_2'), \lambda\neq (\lambda_2,\lambda_1)}\prod SL_{n_\lambda}(q) \times \underset{\lambda \in \epsilon_n, \lambda < \lambda', \lambda=(\lambda_2,\lambda_1)}\prod SL_{n_{\lambda}}(q^{\frac{1}{2}}) \times $$
 $$ \underset{\lambda \in \epsilon_n, \lambda< \lambda',\lambda=(\lambda_1',\lambda_2')}\prod SU_{n_\lambda}(q^\frac{1}{2})\times  \underset{\lambda\in \epsilon_n, \lambda=\lambda', \lambda<(\lambda_1',\lambda_2')}\prod OSP(\lambda)' \times \underset{\lambda \in \epsilon_n, \lambda=\lambda',\lambda=(\lambda_1',\lambda_2')}\prod \widetilde{OSP}(\lambda)' .$$
 \end{theo}

 \section{Surjectivity of the morphisms $\Phi_{i,n}$}\label{surjectivitism}
 
 In this section, we conclude the proof of the theorems in the previous section by showing that for all $i\in [\![1,6]\!]$, the morphism $\Phi_{i,n}$ is surjective. The core of the proof for all $i$ will be in Section $2.3.1$. However, for small $n$, the different factorizations which appear in cases $4$, $5$ and $6$ of the field extensions described at the beginning of Section \ref{lalala} change the proofs for small cases. They require to introduce new tools and consider the maximal subgroups of some classical groups in low dimension. Those were classified in \cite{BHRC} and we will use their tables in order to treat those cases. In all cases, the proof is based on inductive reasoning. We start by proving the result for small $n$ and we then use Goursat's lemma and a Theorem by Guralnick and Saxl \cite{BHRC} which relies on the classification of finite simple groups and gives us a list of conditions for a group to be a classical group. Our induction assumptions will allow us to check that all the conditions are verified. It will be more tricky to verify those conditions in cases $4$, $5$ and $6$ because of the additional factorizations.
 
 \subsection{First case : $\F_q=\F_p(\alpha,\beta)= \F_p(\alpha+\alpha^{-1},\beta+\beta^{-1})$, $\F_p(\alpha)=\F_p(\alpha+\alpha^{-1})$}
 
 In this subsection, we prove the surjectivity of the morphism in the easiest case and establish groundwork for the other cases. This will conclude the proof of Theorem \ref{result1}. We first prove the result for $n\leq 4$ and then use induction to get the result for all $n$.

We recall Goursat's Lemma also used in \cite{BM} and \cite{BMM}, it was originally proven in \cite{Gour} and an english version of the proof can be found in \cite{Bouc} (Lemma $2.3.5$) or \cite{RIB} Lemma $5.2.1$:

\begin{lemme}[Goursat's Lemma]\label{Goursat}
Let $G_1$ and $G_2$ be two groups, $K\leq G_1 \times G_2$, and write $\pi_i: K \longrightarrow G_i$ the projection. Let $K_i=\pi_i(K)$ and $K^i = \op{ker}(\pi_{i'})$ for $\{i,i'\}=\{1,2\}$. There exists an isomorphism $\varphi : K_1/K^1 \rightarrow K_2/K^2$ such that $K=\{(k_1,k_2) \in K_1 \times K_2, \varphi(k_1K^1) =k_2K^2\}$.
\end{lemme}

We first prove that if for any $\lambda\Vdash n$, the composition of $R_\lambda$ with the projection on its corresponding quasi-simple factor is surjective, then $\Phi_{1,n}$ is surjective. We will then prove by induction that each composition is indeed surjective. In order to get the images of the hook partitions it is enough to get the images inside the representations associated with the partitions $([1^{n-1}],[1])$ and $([1],[1^{n-1}])$ . We recall Wagner's theorem which can be found for example in \cite[II, Theorem 2.3]{MM}.
\begin{theo}\label{wag}
Let $\F_r$ be a finite field, $n\in \N, n\geq 3$ and $G \subset GL_n(r)$ a primitive group generated by pseudo-reflections of order greater than or equal to $3$. Then one of the following properties is true.
\begin{enumerate}
\item $SL_n(\tilde{r}) \subset G \subset GL_n(\tilde{r})$ for some $\tilde{r}$ dividing $r$.
\item $SU_n(\tilde{r}^\frac{1}{2}) \subset G \subset GU_n(\tilde{r}^\frac{1}{2})$ for some $\tilde{r}$ dividing $r$.
\item $n\leq 4$ and the pseudo-reflections are of order $3$ and $G \simeq GU_n(2)$.
\end{enumerate}
\end{theo}
 
 \begin{prop}\label{lesgourgues}
Let $n\geq 3$ and $R_1$ (resp $R_2$)  be the representation  associated with the double-partition $([1^{n-1}],[1])$ (resp $([1],[1^{n-1}])$). We have $R_1(\mathcal{A}_{B_n})\simeq R_2(\mathcal{A}_{B_n})\simeq SL_n(q)$.
\end{prop}

\begin{proof}
Let $n\geq 3$, we will use Theorem \ref{wag}. The eigenvalues of $R_1(T)$ are $\beta$ with multiplicity $n-1$ and $-1$ with multiplicity $1$. The eigenvalues of $R_1(S_i)$ are $\alpha$ with multiplicity $1$, and $-1$ with multiplicity $n-1$. The group $G=\langle\beta^{-1}R_1(T),-R_1(S_1),...,-R_1(S_{n-1})\rangle$ is generated by pseudo-reflections. To apply Wagner's Theorem (Theorem \ref{wag}), we must show that the group is primitive. If $G$ was imprimitive, we could write $\F_q^n=V_1\oplus V_2 \oplus ... \oplus V_r$, where for all $i$ and for all $g\in G$, there exists a $j$ such that $g.V_i =V_j$. Since $R_1$ is irreducible, either $\beta^{-1}R_1(T) .V_1 \neq V_1$ or there exists $i\leq n-1$ such that $-R_1(S_i) .V_1 \neq V_1$. Assume there exists $i$ such that $-R_1(S_i).V_1 \neq V_1$. Up to reordering, we have $V_2=-R(S_i).V_1$. If $\dim(V_1)\geq 2$, then $H_{-R_1(S_i)}$ (the hyperplane fixed by $-R_1(S_i)$) has a non-empty intersection with $V_1$. It follows that $V_1\cap V_2 \neq \emptyset$, which is a contradiction, which proves that $\dim(V_1)=1$. This reasoning is valid for any $V_i$, therefore they are all one-dimensional. Let $x\in V_1$ be a non-zero vector, it can be written in a unique way as $x =x_1+x_2$ with $x_1\in \op{ker}(R_1(S_i)+\alpha)$ and $x_2\in H_{-R_1(S_i)}$. We then have that $-R_1(S_i)x = -\alpha x_1+x_2$ and $-R(S_i)(-R(S_i)x)= \alpha^2 x_1 +x_2=\alpha(x_1+x_2)+(1-\alpha)(-\alpha x_1+x_2)\in V_1\oplus V_2$. Since $\alpha\notin \{0,1\}$ this contradicts the fact that there exists $j$ such that $-R(S_i).V_2 =V_j$.
If $V_1\neq \beta^{-1}R_1(T) V_1 =V_2$, then if $x=x_1+x_2\neq 0$ with $x\in V_1,x_1\in \op{ker}(R_1(T)+\beta) x_2\in H_{\beta^{-1}R_1(T)}$, we have that $\beta^{-1}R_1(T)x = -\beta^{-1}x_1+x_2$ and $\beta^{-1}R_1(T)(\beta^{-1}R_1(T)x) = \beta^{-2}x_1+x_2=\beta^{-1}(x_1+x_2)+(1-\beta^{-1})(-\beta^{-1}x_1+x_2)\in V_1\oplus V_2$. This is absurd because $\beta^{-1}\notin \{0,1\}$. This shows that $G$ is primitive and in the same way, $\tilde{G}=\langle -R_2(T),-R_2(S_1),...,-R_2(S_{n-1})\rangle$ is primitive and generated by pseudo-reflections of order greater than or equal to 3. By Theorem \ref{wag}, we have $SL_n(\tilde{q}) \subset G,\tilde{G} \subset GL_n(\tilde{q})$ or $SU_n(\tilde{q}^\frac{1}{2}) \subset G,\tilde{G} \subset GU_n(\tilde{q}^\frac{1}{2})$ for some $\tilde{q}$ dividing $q$. If we were in the unitary case then there would exist an automorphism $\epsilon$ of order 2 of $\F_{\tilde{q}^2}$ such that $\det(M)=\epsilon(\det(M))^{-1}$ for all $M$ in $G$ or $\tilde{G}$. We also have $\det(\beta^{-1}R_1(T)) = -\beta^{-1}, \det(-R_1(S_1))= \det(-R_2(S_1)) = -\alpha$ and $\det(-R_2(T))=  -\beta$. If $G$ or $\tilde{G}$ is unitary then $\epsilon(\beta)=\beta^{-1}$ and $\epsilon(\alpha)=\alpha^{-1}$, therefore $\alpha+\alpha^{-1}$ and $\beta+\beta^{-1}$ are in $\F_{\tilde{q}}$. This contradicts the fact that $\tilde{q}^2$ divides $q$ and $\F_q=\F_p(\alpha+\alpha^{-1},\beta+\beta^{-1})$. This proves we have $SL_n(\tilde{q}) \subset G,\tilde{G} \subset GL_n(\tilde{q})$ for some $\tilde{q}$ dividing $q$. Using the determinants again, we have $\alpha$ and $\beta$ in $\F_{\tilde{q}}$, therefore $\tilde{q}=q$. We have $SL_n(q)=[G,G]=[R_1(A_{B_n}),R_1(A_{B_n})]=R_1(\mathcal{A}_{B_n})$ and $SL_n(q)=[\tilde{G},\tilde{G}]=[R_2(A_{B_n}),R_2(A_{B_n})]=R_2(\mathcal{A}_{B_n})$ which concludes the proof.
\end{proof}

By \cite{MR}, $\mathcal{A}_{B_n}$ is perfect for $n\geq 5$ but not for $n\leq 4$. Those cases must then be treated separately.

\begin{lemme}\label{platypus}
If $n\leq 4$ then $\Phi_{1,n}$ is surjective.
\end{lemme}

\begin{proof}
The representations labeled by double-partitions of $n=2$ are all one-dimensional except for $([1],[1])$. We then only need to show that $R_{[1],[1]}(\mathcal{A}_{B_2})=SL_2(q)$. We write $t=R_{[1],[1]}(T)=\begin{pmatrix} \beta & 0 \\ 0 & -1 \end{pmatrix}$ and $s=R_{[1],[1]}(S_1)= \frac{1}{\beta+1} \begin{pmatrix} \alpha -1 & \alpha +\beta\\ \alpha \beta+1 & \alpha \beta -\beta \end{pmatrix}$. First note that if $P = \begin{pmatrix} 1 & 1\\ \frac{\alpha\beta+1}{\alpha+\beta} & -1 \end{pmatrix}$ then $P^{-1}tP = \frac{1}{\alpha+1}\begin{pmatrix} \beta -1 & \alpha+\beta\\
\alpha\beta+1 & \beta\alpha -\alpha \end{pmatrix}$ and $P^{-1}sP = \begin{pmatrix}
\alpha & 0\\
0 & -1.
\end{pmatrix}$. This proves that the roles of $\alpha$ and $\beta$ are completely symmetrical in this case. It follows that up to conjugating by $P$, we can exchange the conditions on $\alpha$ and the conditions on $\beta$. We write $G=<t,s>$. We have $\det(t) = -\beta$ and $\det(s) = -\alpha$. Let $(u,v)\in \overline{\F_p}^2$ such that $u^2 = -\beta^{-1}$ and $v^2 = -\alpha^{-1}$. We set $\F_{q'} = \F_q(u,v)$. We then have $\overline{G}=\langle \overline{t},\overline{s}\rangle = \langle\overline{ut},\overline{vs}\rangle \subset PSL_2(q')$. We write $\mathfrak{S}_n$ the permutation group of $n$ elements and $\mathfrak{A}_n$ its derived subgroup. By Dickson's Theorem \cite[Chapter II, HauptSatz 8.27]{HUP}, we have that $\overline{G}$ is either abelian by abelian or isomorphic to $ \mathfrak{A}_5, \mathfrak{S}_4, PSL_2(\tilde{q})$ or $PGL_2(\tilde{q})$ for a given $\tilde{q}$ greater than or equal to $4$ and dividing $q'$. 

If $\overline{ut}^r = 1$ in $PSL_2(q')$, then $((-u)^r)^2 = 1$. Therefore $\beta^r = (-1)^r$ and by the condition on the order of $\alpha$, $\overline{G}$ cannot be isomorphic to $ \mathfrak{A}_4$ or $\mathfrak{A}_5$.

We now exclude the case $\overline{G}$ abelian by abelian. If $\overline{G}$ is abelian by abelian, then $[\overline{G},\overline{G}]$ is abelian, i.e. $\overline{ab}= \overline{ba}$ for all $a,b \in [G,G]$ or equivalently $ab=\pm ba$ for all $a,b \in [G,G]$. We have that $(tst^{-1}s^{-1}) (s^{-1}tst^{-1})-(s^{-1}tst^{-1})(tst^{-1}s^{-1}) = $
$$ \begin{pmatrix}
-\frac{(\beta-1)(\alpha-1)^2(\alpha\beta+1)(\alpha+\beta)}{\beta\alpha^2(\beta+1)} & -\frac{(\alpha^2\beta^2+\alpha\beta^3-\alpha\beta^2-\alpha^2\beta+\alpha\beta+\beta^2-\alpha-\beta)(\alpha-1)(\alpha\beta+1)}{\beta\alpha^2(\beta+1)}\\
\frac{(\alpha^2\beta^2+\alpha\beta^3-\alpha\beta^2-\alpha^2\beta+\alpha\beta+\beta^2-\alpha-\beta)(\alpha-1)(\alpha+\beta)}{\beta^2\alpha^2(\beta+1)}
& \frac{(\beta-1)(\alpha-1)^2(\alpha\beta+1)(\alpha+\beta)}{\beta\alpha^2(\beta+1)}
\end{pmatrix}.$$
This matrix is non-zero because the diagonal coefficients are non-zero by the conditions on $\beta$. This means that if $[\overline{G},\overline{G}]$ is abelian, then we have $(tst^{-1}s^{-1}) (s^{-1}tst^{-1})+(s^{-1}tst^{-1})(tst^{-1}s^{-1}) = 0$, but this matrix equals
\begin{footnotesize}
$$\begin{pmatrix}
 \frac{{\alpha}^{4}\beta+{\alpha}^{3}{\beta}^{2}-2{\alpha}^{3}\beta-2{\alpha}^{2}{\beta}^{2}+{\alpha}^{3}+4{\alpha}^{2}\beta+\alpha{\beta}^{2}-2{\alpha}^{2}-2\alpha\beta+\alpha+\beta}{{\alpha}^{2}\beta} & - {\frac { \left( {\alpha}^{2}\beta+\alpha\,{\beta}^{2}-2\,\alpha\,\beta+\alpha+\beta \right)  \left( \alpha-1 \right)  \left( \alpha\,\beta+1 \right) }{{\alpha}^{2}\beta}}\\ 
-{\frac { ( {\alpha}^{2}\beta+\alpha\,{\beta}^{2}-2\,\alpha\,\beta+\alpha+\beta ) ( \alpha+\beta )  ( \alpha-1) }{{\alpha}^{2}{\beta}^{2}}} & \frac {{\alpha}^{4}\beta+{\alpha}^{3}{\beta}^{2}-2\,{\alpha}^{3}\beta-2\,{\alpha}^{2}{\beta}^{2}+{\alpha}^{3}+4\,{\alpha}^{2}\beta+\alpha\,{\beta}^{2}-2\,{\alpha}^{2}-2\,\alpha\,\beta+\alpha+\beta}{{\alpha}^{2}\beta }
\end{pmatrix}.$$
\end{footnotesize}
The non-diagonal coefficients are non-zero if $A= {\alpha}^{2}\beta+\alpha\,{\beta}^{2}-2\,\alpha\,\beta+\alpha+\beta $ is non-zero. If $A=0$ then, the bottom right coefficient of $(tst^{-1}s^{-1})(st^{-1}s^{-1}t)+(st^{-1}s^{-1}t)(tst^{-1}s^{-1})$ is equal to $-\frac{1}{(\beta+1)\alpha^2\beta^2}$ multiplied by
$${\alpha}^{4}{\beta}^{3}+{\alpha}^{3}{\beta}^{4}-{\alpha}^{4}{\beta}^{2}-3\,{\alpha}^{3}{\beta}^{3}-2\,{\alpha}^{2}{\beta}^{4}+5\,{\alpha}^{3}{\beta}^{2}+4\,{\alpha}^{2}{\beta}^{3}+$$
$$\alpha\,{\beta}^{4}-3\,{\alpha}^{3}\beta-8\,{\alpha}^{2}{\beta}^{2}-3\,\alpha\,{\beta}^{3}+4\,{\alpha}^{2}\beta+5\,\alpha\,{\beta}^{2}+{\beta}^{3}-2\,{\alpha}^{2}-3\,\alpha\,\beta-{\beta}^{2}=$$
$$ (\alpha^2\beta^2-\alpha^2\beta+2\alpha\beta-2\alpha)A+\beta((\beta-1)A-2\alpha^2(\beta^3+1)) = -2\beta\alpha^2(\beta^3+1).$$ This is non-zero by the condition on $\beta$. The diagonal coefficients of the difference of these two commutators are identical to the ones of the difference of the previous commutators, therefore they are non-zero. This proves that $\overline{G}$ is not abelian by abelian and there exists $\tilde{q}$ greater than or equal to $4$ such that $[\overline{G},\overline{G}] \simeq PSL_2(\tilde{q})$.

For $H$ a group and $A$ an $H$-module, we write $Z^2(H,A) = \{f :H\times H \rightarrow A, \forall x,y,z \in H, z.f(x,y)f(xy,z)=f(y,z)f(xy,z)\}$ the group of cocyles and  $B^2(H,A) = \{f :H\times H \rightarrow A,\exists t : H \rightarrow A, \forall x,y \in H, f(x,y)= t(y)t(xy)^{-1}y.t(x)\}$ the group of coboundaries. We write $M(H) = H^2(H,\C^\star) = H_2(H,\Z)$ its Schur multiplier. We have $[\overline{G},\overline{G}] = \overline{[G,G]} \subset PSL_2(q)$, this inclusion gives a projective representation of $SL_2(\tilde{q})$ with associated cocyle $c\in Z^2(SL_2(\F_{\tilde{q}}),\F_q^\star)$. We will show that $H^2(SL_2(\F_{\tilde{q}}),\F_q^\star)$ is trivial which is equivalent to $Z^2(SL_2(\F_{\tilde{q}}),\F_q^\star) =B^2(SL_2(\tilde{q}),\F_q^\star)$ which implies that this cocycle is a coboundary. We write $H =SL_2(\tilde{q})$. By the Universal Coefficients Theorem  \cite[Theorem 3.2]{HAT}, we have the following exact sequence 
$$ 1 \rightarrow \op{Ext}(H_1(H,\Z),\F_q^\star)  \rightarrow H^2(H,\F_q^\star) \rightarrow \op{Hom}(H_2(H,\Z),\F_q^\star)\rightarrow 1.$$
We have $H_1(H,\Z) = H/[H,H]$ \cite{KAR} and $H_2(H,\Z) = M(H)$. Since $\tilde{q} \geq 4$, we have that $SL_2(\tilde{q})$ is perfect and the exact sequence becomes
$$ 1 \rightarrow 1 \rightarrow H^2(H,\F_q^\star) \rightarrow \op{Hom}(M(H),\F_q^\star)\rightarrow 1.$$
By \cite[Theorem 7.1.1]{KAR}, if $\tilde{q}\notin \{4,9\}$ then the Schur multiplier $M(H)$ is trivial, therefore this reduces to $ H^2(H,\F_q^\star)\simeq \{1\}$.

It remains to consider the cases $\tilde{q}= 4$ and $\tilde{q}=9$. If $\tilde{q} =4$, we have $M(H) = \Z/2\Z$ and $p = 2$, therefore  $\op{Hom}(M(H),\F_q^\star) =1$. Indeed, every morphism $\varphi$ from $M(H)$ to $\F_q^\star$ satisfies $1 = \varphi(2x)=\varphi(x)^2$ for all $x\in M(H)$. It follows that $0=\varphi(x)^2 -1= (\varphi(x)-1)^2$ for all $x\in M(H)$. This proves that $\varphi$ is trivial. If $\tilde{q} = 9$ then we have $M(H) = \Z/3\Z$ and $H^2(H,\F_q^\star)$ is trivial by the same reasoning as for $\tilde{q}=4$. In all cases, we can define a representation  $\rho$ of $SL_2(\tilde{q})$ in $SL_2(q)$.

By \cite{BN}, any representation $\sigma$ of $SL_2(q)$ in $GL_2(q)$ is up to conjugation of the form $\sigma(M)=\psi(M)$, where $\psi(M)$ is the matrix obtained from $M$ by applying $\psi \in Aut(\F_q)$ to all its coefficients. We have $\F_q = \F_{\tilde{q}}(w)$ for any $w$ generating the cyclic group $\F_q^{\star}$. There exists a homomorphism from $\F_q$ to $\F_{\tilde{q}}$ sending $1$ to $w$ and stabilizing $\F_{\tilde{q}}$. We define a representation $\tilde{\rho}$ of $SL_2(q)$ in $SL_2(q)$ such that $\tilde{\rho}(M) = \rho(M)$ for all $M$ in $SL_2(\tilde{q})$. We have $\rho(M) = \tilde{\rho}(M) = \psi(M)$ for all $M$ in $SL_2(\tilde{q})$. We have $[\overline{G},\overline{G}] \simeq PSL_2(\tilde{q})$, therefore $\psi([\overline{G},\overline{G}]) \simeq PSL_2(\tilde{q})$ is conjugate to $PSL_2(\tilde{q})$ in $GL_2(q')$. We have $\psi(tst^{-1}s^{-1}) \in \psi([\overline{G},\overline{G}])$, therefore its trace $2 - (\alpha + \alpha^{-1} + \beta + \beta^{-1})$ belongs to $\F_{\tilde{q}}$. This shows that $\alpha + \alpha^{-1} + \beta + \beta^{-1} \in \F_{\tilde{q}}$. We also have that the trace $T_1$ of $s^2t^{-1}s^{-2}t$ and the trace $T_2$ of $st^{-2}s^{-1}t^2$ are in $\F_{\tilde{q}}$. We have $T_1 = \frac{\alpha^4\beta+\alpha^3\beta^2-2\alpha^3\beta-2\alpha^2\beta^2+\alpha^3+4\alpha^2\beta+\alpha\beta^2-2\alpha^2-2\alpha\beta+\alpha+\beta}{\alpha^2\beta}$. We write $B = \alpha+\alpha^{-1}+\beta+\beta^{-1}$. We then have $T_1 = (\alpha+\alpha^{-1})B-2B+2$. The quantity $T_2$ has the same expression as $T_1$ with $\alpha$ and $\beta$ switched, therefore $T_2 = (\beta+\beta^{-1})B-2B+2$. Since $B, T_1$ and $T_2$ are in $\F_{\tilde{q}}$, we have $(\alpha+\alpha^{-1})B$ and $(\beta+\beta^{-1})B$ are in $\F_{\tilde{q}}$. We have $B = \alpha+\alpha^{-1}+\beta+ \beta^{-1}$. It follows that $B =0$ implies $\alpha\in \{-\beta,-\beta^{-1}\}$; which contradicts the assumptions on $\alpha$ and $\beta$. The quantity $B$ is therefore non-zero, therefore $\alpha+\alpha^{-1}$ and $\beta+\beta^{-1}$ are in $\F_{\tilde{q}}$, and $\F_{\tilde{q}}=\F_q$. We conclude using  Lemma 2.1 of \cite{BM}.

\bigskip

\bigskip

The double-partitions of $n=3$ to consider are $([2,1],\emptyset),([1],[1^2])$ and $([1^2],[1])$. We want to show that the image of $\Phi_3$ is equal to $SL_2(\tilde{q}) \times SL_3(q) \times SL_3(q)$. If we restrict ourselves to the image inside $SL_3(q) \times SL_3(q)$, we show that it is $SL_3(q)\times SL_3(q)$. By Proposition \ref{lesgourgues}, $SL_3(q)=R_{[1^2],[1]}(\mathcal{A}_{B_3})$ and $R_{[1],[1^2]}(\mathcal{A}_{B_3})=SL_3(q)$. We now use Goursat's Lemma : we write as in Lemma \ref{Goursat}, $K=R(\mathcal{A}_{B_3}), K_1=R_{[1^2],[1]}(\mathcal{A}_{B_3}), K_2= R_{[1],[1^2]}(\mathcal{A}_{B_3}) $ $\pi_1$ (resp $\pi_2$) the projection onto $SL_{3}(q)$(corresponding to $([1^2],[1])$) (resp $SL_{3}(q))$ (corresponding to $([1],[1^2])$)), $K^1=\op{ker}(\pi_2)$, $K^2=\op{ker}(\pi_1)$ and $\varphi$ the isomorphism given by Goursat's Lemma. We have $K=\{(x,y)\in K_1\times K_2, \varphi(xK^1)=yK^2\}$. By the same reasoning as the one in Proposition 3.1. of \cite{BM}, either  $SL_3(q)\times SL_3(q) \subset K$ or $K_1/K^1$ is non-abelian and $\varphi$ is an isomorphism of $PSL_3(q)$ and using the same notations, up to conjugation $R_2(b) = S^{\phi}(R_1(b))z(b)$ for all $b\in \mathcal{A}_{B_3}$ ($S=Id$ or $S=M\mapsto {}^t\!M^{-1}$).

Let us show that the second possibility is absurd by choosing the right elements in $\mathcal{A}_{B_3}$. For any element $b$ of $\mathcal{A}_{B_3}, \tr(R_{[1],[1^2]}(b)) = z(b)\tr(S^{\phi}(R_{[1^2],[1]}(b))$. We write $U=S_1S_2^{-1}, V=TS_1T^{-1}S_2^{-1}, W=S_2S_1S_2^{-2}$ and $X=S_2TS_1T^{-1}S_2^{-2}$, they are all elements of $\mathcal{A}_{B_3}$. By explicit computation, for both choices of $S$, we have :
$$\tr(R_{[1],[1^2]}(U))=\tr(R_{[1],[1^2]}(V))=\tr(R_{[1],[1^2]}(W))=\tr(R_{[1],[1^2]}(X))=-\frac{(\alpha-1)^2}{\alpha},$$
\begin{footnotesize}
$$\tr(S(R_{[1^2],[1]}(U)))=\tr(S(R_{[1^2],[1]}(V)))=\tr(S(R_{[1^2],[1]}(W)))=\tr(S(R_{[1^2],[1]}(X)))=-\frac{(\alpha-1)^2}{\alpha}.$$
\end{footnotesize}
This shows that $z(U)=z(V)=z(W)=z(X)$ and
$$1=z(U)z(W)^{-1}=z(UW^{-1})=\frac{-\frac{(\alpha-1)^2}{\alpha}}{\phi(-\frac{(\alpha-1)^2}{\alpha})}.$$
This proves $\phi(\alpha+\alpha^{-1})=\alpha+\alpha^{-1}$. We also have
$$1=z(UV^{-1})= \frac{3-\alpha-\alpha^{-1}-\beta-\beta^{-1}}{\phi(3-\alpha-\alpha^{-1}-\beta-\beta^{-1})}.$$
Using $\phi(\alpha+\alpha^{-1})=\alpha+\alpha^{-1}$, we have $\phi(\beta+\beta^{-1})=\beta+\beta^{-1}$, therefore $\phi=1$. We deduce that
\begin{footnotesize}
$$1=z(UX^{-1})=\frac{(\alpha-1)(\alpha\beta^2-2\alpha\beta+2\beta-1)}{\alpha\beta}\frac{-\alpha\beta}{(\alpha-1)(2\alpha\beta+\beta^2-\alpha-2\beta)}= \frac{1+2\alpha\beta-2\beta-\alpha\beta^2}{\beta^2+2\alpha\beta-2\beta-\alpha},$$
\end{footnotesize}
$$1-\alpha\beta^2=\beta^2-\alpha, (1-\beta^2)(1+\alpha)=0.$$
Since $\beta^2\neq 1$ and $\alpha^2\neq 1$, we get a contradiction. This shows that $SL_3(q)\times SL_3(q) = R(\mathcal{A}_{B_3})$.

We now set $G_1= SL_3(q) \times SL_3(q)$ and $G_2 = SL_2(\tilde{q})$, the image of $\Phi_3$ is a subgroup of $G_1 \times G_2$ for which the projections onto $G_1$ and $G_2$ are surjective. Using again Goursat's Lemma and the notation there, we have $K_1/K^1 \simeq K_2/K^2$. We have a surjective morphism $\psi$ from $K_1=G_1=SL_3(q) \times SL_3(q)$ to $K_2/K^2$, where $K_2=SL_2(\tilde{q})$. If $K_2/K^2$ was non-abelian then we would have $K_2/K^2 \simeq PSL_2(\tilde{q})$. If the restriction $\psi_1$ (resp $\psi_2$) of $\psi$ to $SL([1^2],[1])$ (resp $SL([1],[1^2])$) was not trivial then $\psi_1$ (resp $\psi_2$) would factor into an isomorphism from  $PSL_3(q)$ onto $PSL_2(q)$ since the center of $SL_3(q)$ would again be in the center of $\psi_1$ and $\psi_2$.This would lead to a contradiction, therefore their image is trivial and $\psi$ is not surjective, therefore the quotients are abelian. This shows that $K_1=[K_1,K_1]\subset K^1$ and $K_2=[K_2,K_2]\subset K^2$ then using Goursat's Lemma we conclude that the image of $\Phi_3$ is equal to $G_1 \times G_2$. This shows that $\Phi_3$ is  surjective.

\bigskip

\bigskip

The double-partitions of $4$ in our decomposition are $([1^4],\emptyset)$, $([2^2],\emptyset)$, $([2,1,1],\emptyset)$, $([1],[1^3])$, $([1^3],[1])$, $([1^2],[1^2])$ and $([2,1],[1])$ of respective dimensions $1$, $2$, $3$, $4$, $4$, $6$ and $8$ (we removed the hooks, $([3],[1])$ and $([1],[3])$). We know the restriction to the first five is surjective by \cite{BMM} and Proposition \ref{lesgourgues}, therefore we only need to show that $R_{[1^2],[1^2]}(\mathcal{A}_{B_4})=SL_6(q)$ and $R_{([2,1],[1])}(\mathcal{A}_{B_4})=SL_8(q)$.

Let us first consider the double-partition $([1^2],[1^2])$. By the branching rule and the case $n=3$ above, we have 
$$SL_3(q)\times SL_3(q)=R_{[1^2],[1]}(\mathcal{A}_{B_3})\times R_{[1],[1^2]}(\mathcal{A}_{B_3})=R_{[1^2],[1^2]}(\mathcal{A}_{B_3}) \subset R_{[1^2],[1^2]}(\mathcal{A}_{B_4})\subset SL_6(q).$$

We can now use Theorem 3 from \cite{BM}.
\begin{theo}\label{LBJ}
Let $\F_r$ be a finite field and $\Gamma < GL_N(r)$ with $N\geq 5$ and $q>3$ such that
\begin{enumerate}
\item $\Gamma$ is absolutely irreducible,
\item $\Gamma$ contains $SL_a(r)$ in a natural representation with $a\geq \frac{N}{2}$.
\end{enumerate}
 If $N\neq 2a$, then $\Gamma$ contains $SL_N(r)$. Otherwise, either $\Gamma$ contains $SL_N(r)$, or $\Gamma$ is a subgroup of $GL_{\frac{N}{2}}(r) \wr \mathfrak{S}_2$.
\end{theo} 

We use this theorem on $R_{[1^2],[1^2]}(A_{B_4})$. To get the desired result, we only need to show that $R_{[1^2],[1^2]}(A_{B_4})$ cannot be a subgroup of  $GL_{3}(r) \wr \mathfrak{S}_2$. If it were true, then we would have $R_{[1^2],[1^2]}(\mathcal{A}_{B_4})\subset SL_3(q)\times SL_3(q)$ which would contradict the irreducibility shown in Lemma \ref{Lincoln}. This shows that we have $R_{[1^2],[1^2]}(\mathcal{A}_{B_4})=SL_6(q)$.

\bigskip

We now consider the double-partition $([2,1],[1])$. Again by the branching rule and the case $n=3$, we have that the restriction to $\mathcal{A}_{B_3}$ is $ SL_3(q) \times SL_3(q)  \times SL_2(\tilde{q})$. We now use the fact that each of these groups is generated by transvections and the fact that $\mathcal{A}_{B_4}$ is normally generated by $\mathcal{A}_{B_3}$. When the characteristic is different from $2$, we can use Theorem \ref{transvections} 

We write $G= R_{([2,1],[1])}(\mathcal{A}_{B_4})$, $H = R_{([2,1],[1])}(\mathcal{A}_{B_3}) = SL_3(q) \times SL_3(q) \times SL_2(\tilde{q})$ and we pick $t_1$ (resp $t_2$ (resp $t_3$)) a transvection of $SL_3(q) \times \{I_5\}$ (resp $\{I_3\} \times SL_3(q) \times \{I_2\}$ (resp $\{I_6\} \times SL_2(\tilde{q})$). We then have $H = \langle ht_ih^{-1}, h\in H, i \in \{1,2,3\}\rangle$, therefore $G = \langle ghg^{-1}, h\in H, g\in G\rangle = \langle gt_ig^{-1}, g\in G, i\in \{1,2,3\}\rangle$ is generated by transvections and we can apply the theorem.  We also recall the following lemma \cite[Lemma 5.6]{BMM}.
\begin{lemme}\label{field}
For any prime $p$ and $m\geq 2$, the field generated over $\F_p$ by $\{\tr(g), g \in SL_m(q)\}$ is $\F_q$ and for all $m\geq 3$, the field generated over $\F_p$ by $\{\tr(g), g\in SU_m(q^{\frac{1}{2}})\}$ is $\F_q$.
\end{lemme}

By Proposition \ref{isomorphisme}, we know that  $R_{[2,1],[1]}(\mathcal{A}_{B_4})$ preserves no non-degenerate bilinear form. It also shows that it can preserve no non-degenerate hermitian form. Indeed, if it were to preserve a hermitian form then we would have $\tr(M)=\epsilon(\tr({}^t\!(M)^{-1}))$ for any $M$ in $G$ and we have $\diag([\alpha,\alpha^{-1},1,1,1,1,1,1])$ and $\diag([\beta,\beta^{-1},1,1,1,1,1,1])$ in $H \subset G$. Therefore we would have $\epsilon(\alpha+\alpha^{-1})= \alpha+\alpha^{-1}$ and $\epsilon(\beta+\beta^{-1})=\beta+\beta^{-1}$. Since $\F_q=\F_p(\alpha,\beta) =\F_p(\alpha+\alpha^{-1},\beta+\beta^{-1})$, the automorphism $\epsilon$ of order $2$ would be trivial which is a contradiction. This proves that $G$ is conjugate in $GL_8(q)$ to $SL_8(\tilde{q})$ for some $\tilde{q}$ dividing $q$. By Lemma \ref{field}, the field generated over $\F_p$ by the traces of the elements of $G$ is $\F_{\tilde{q}}$, therefore $\tilde{q}= q$ because $G$ contains $SL_3(q)$ in a natural representation. This proves that the field generated by its elements contains $\F_q$. This shows that when $p\neq 2$, $G = R_{([2,1],[1])}(\mathcal{A}_{B_4}) = SL_8(q)$.

Assume now that $p=2$, we can use the following theorem \cite[Theorem 1]{P}.
\begin{theo}
Let $V$ be a $\F_q$-vector space of dimension $n\geq 4$ with $q$ even. If $G$ is an irreductible proper subgroup of $SL(V)=SL_n(q)$ generated by a set $D$ of transvections of $G$, then $D$ is a conjugacy class of odd transpositions of $G$.
\end{theo}

Assume that $G = R_{([2,1],[1])}(\mathcal{A}_{B_4})$ is different from $SL_8(q)$. We again have that $G$ is generated by transvections and by applying the above theorem, those transvections are in a single conjugacy class of $G$. Since $O_p(G)$ is normal in $G$ and $V=\F_q^8$ is an irreducible $\F_qG$-module, we apply Clifford's Theorem \cite[Theorem 11.1]{C-R} and get that $\op{Res}^G_{O_p(G)}(V)$ is semisimple. Since $O_p(G)$ is a $p$-group, the unique irreducible $\F_qO_p(G)$-module is the trivial module, therefore $O_p(G)$ acts trivially on $V$. It follows that $O_p(G)$ is trivial. We can thus apply Kantor's Theorem \cite[Theorem II]{K} :
\begin{theo}
Let $p$ be a prime and $q=p^l$ for some $l\in \N$. Assume $G$ is an irreducible subgroup of $SL_N(q)$ generated by a conjugacy class of transvections, such that $O_p(G) \leq [G,G] \cap Z(G)$. Then $G$ is one of the following subgroups.
\begin{enumerate}
\item $G=SL_n(q')$ or $G\simeq Sp_N(q')$ in $SL_N(q')$ or $G \simeq  SU_N(q'^{\frac{1}{2}})$ in $SL_N(q')$, with $q' | q$.
\item $G\simeq  O_N^{\pm}(q') < SL_n(q')$, with $q'|q$.
\item $G \simeq  \mathfrak{S}_n < SL_N(2)$, where $N = n-d$ and $d = \op{Gcd}(n,2)$.
\item $G \simeq  \mathfrak{S}_{2n}$ in $SL_{2n-1}(2)$ or in $SL_{2n}(2)$.
\item $G \simeq  SL_2(5) < SL_2(9)$.
\item $G \simeq  3.P\Omega_6^{-,\pi} < SL_6(4)$.
\item $G \simeq  SU_4(2) < SL_5(4)$.
\item $G \simeq  A \rtimes S_N$ in $SL_N(2^i)$, where $A$ is a subgroup of diagonal matrices.
\end{enumerate}
\end{theo}

Since $\alpha$ is of order greater than $4$, we have $q \geq \tilde{q}=2^r> 8$. The group $G$ contains $H = SL_3(q) \times SL_3(q) \times SL_2(\tilde{q})$, therefore cases ($3$) to ($7$) are excluded. If we were in case ($8$), then $G$ would have at mose $(q-1)\frac{10\times 9}{2}=45(q-1)$ transvections (see proof of Theorem 1.3. page 661 of \cite{BM}). $SL_3(q)$ has $\frac{(q^3-1)(q^2-1)}{q-1} = (q-1)(q^2+q+1)(q+1)$ transvections and  $(q^2+q+1)(q+1) \geq  847(q-1)> 45(q-1)$. For the same reasons as when $p\neq 2$, $G$ is neither unitary nor symplectic nor orthogonal. The only remaining possibility is $G = SL_8(q)$ which is a contradiction since we assumed $G \neq SL_8(q)$. This proves that $G=SL_8(q)$.

\bigskip

The restriction to each double-partition of $4$ is thus surjective, it remains to show that $\Phi_{1,4}$ is surjective using Goursat's Lemma (Lemma \ref{Goursat}). This means that we have to show that the image is $SL_2(\tilde{q}) \times SL_3(\tilde{q}) \times SL_4(q) \times SL_4(q) \times SL_8(q)$.

By Theorem 1.2. of \cite{BM}, the restriction to $SL_2(\tilde{q}) \times SL_3(\tilde{q})$ is surjective. We write $G_1$ the image of this restriction and $G_2 = R_{([1],[1^3])}(\mathcal{A}_{B_4}) = SL_4(q)$. Let $K$ be the image $\mathcal{A}_{B_4}$ in $G_1 \times G_2$, using the corresponding notations in Goursat's Lemma, we have $K_1 = G_1$, $K_2 = G_2$ and $K_1/K^1 \simeq K_2/K^2$. If these quotients are abelian then the proof of $K = G_1 \times G_2$ is straightforward using Goursat's Lemma. Since the only non-abelian decomposition factor of $G_2$ is $PSL_4(q)$ and the only non-abelian decomposition factors of $G_1$ are $PSL_2(\tilde{q})$ and $PSL_3(\tilde{q})$, we have a contradiction if these quotients are non-abelian. Write now $\tilde{K} = R_{([1],[3])}(\mathcal{A}_{B_4})=SL_4(q)$ and let us consider the image $J$ of $\mathcal{A}_{B_4}$ inside $K\times \tilde{K}$. 
Using again Goursat's Lemma, this time with $K_1 = K, K_2=\tilde{K}$, we have $K_1/K^1\simeq K_2/K^2$. If the quotients are abelian then $J=K\times \tilde{K}$. If the quotients are non-abelian then there is an isomorphism $S^\phi$ from $PSL_4(q)$ to $PSL_4(q)$, where the first one corresponds to $\overline{R_{([1],[1^3])}(\mathcal{A}_{B_4})}$ and the second one to $\overline{R_{([1],[3])}(\mathcal{A}_{B_4})}$. This implies that there exists a character $z$ from $\mathcal{A}_{B_4}$ to $\F_q^{\star}$ such that up to conjugation, for every $h\in \mathcal{H}_4$, we have $R_{[1],[1^3]}(h) = S^{\phi}(R_{([1],[3])}(h))z(h)$. The isomorphism $S^\phi$ is of the form \cite[Section 3.3.4]{W} $M\mapsto \phi(M)$ or $M\mapsto \phi({}^t\!(M^{-1}))$, where $\phi$ is a field automorphism of $\F_q$. We would then have  that for all $h\in \mathcal{A}_{B_4}$, 
$$\tr(R_{([1^3],[1])}(h)) = \phi(\tr(S(R_{([1],[3])}(h)))z(h).$$

Writing $U = S_1S_2^{-1}$, $V = TS_1T^{-1}S_2^{-1}$, $X= S_2TS_1T^{-1}S_2^{-2}$, $P=S_3S_2S_3^{-2}$, $Q= TS_1T^{-1}S_3^{-1}$, $R_1=R_{([1^3],[1])}$ and $R_2=R_{([1],[3])}$, we have 
\begin{scriptsize}
$$3-\alpha-\alpha^{-1}=\tr(R_1(P)) =\tr(R_1(PQ^{-1})) = \tr(R_2(P)) = \tr(R_2(PQ^{-1}) = \tr({}^t\!(R_2(PQ^{-1})^{-1}))=\tr({}^t\!(R_2(PQ^{-1})^{-1}).
$$
\end{scriptsize}
It follows that $z(PQ^{-1})=z(P) = \frac{3-\alpha-\alpha^{-1}}{\Phi(3-\alpha-\alpha^{-1})}$.
This shows that $z(Q) = z(P)z(PQ^{-1})^{-1}=1$. We also have $\tr(R_1(Q)) = \tr(R_2(Q))=\tr({}^t\!(R_2(Q)^{-1}))= 2-\alpha-\alpha^{-1}$, therefore $1=z(Q)= \frac{2-\alpha-\alpha^{-1}}{\Phi(2-\alpha-\alpha^{-1})}$ and $\Phi(\alpha+\alpha^{-1}) = \alpha+\alpha^{-1}$. We have 
\begin{tiny}
$$\tr(R_1(U)) = \tr(R_1(V))= \tr(R_1(X)) = \tr(R_2(U)) = \tr(R_2(V))=\tr(R_2(X)) = \tr({}^t\!(R_2(X)^{-1}))= \tr({}^t\!(R_2(V)^{-1}))=\tr({}^t\!(R_2(U)^{-1}) = 3-\alpha-\alpha^{-1}.$$
\end{tiny}
This leads to $z(U) = z(V) = z(X) =1$. We have 
$$\tr(R_1(UV^{-1})) = \tr(R_2(UV^{-1})) = \tr({}^t\!(R_2(UV^{-1})^{-1})) = 4-\alpha-\alpha^{-1}-\beta-\beta^{-1}.$$

 It follows that $z(UV^{-1}) = z(U)z(V)^{-1}= 1 = \frac{4-(\alpha+\alpha^{-1})-(\beta+\beta^{-1})}{\Phi(4-(\alpha+\alpha^{-1})-(\beta+\beta^{-1}))}$, therefore $\Phi(\beta+\beta^{-1}) =\beta+\beta^{-1}$. Since $\F_q=\F_p(\alpha+\alpha^{-1},\beta+\beta^{-1})$, we have $\Phi=I_d$. 
$$\tr(R_1(UX^{-1})) = -\frac{2\alpha^2\beta+\alpha\beta^2-\alpha^2-5\alpha\beta-\beta^2+\alpha+2\beta}{\alpha\beta},$$
$$\tr(R_2(UX^{-1}))=\tr({}^t\!(R_2(UX^{-1})^{-1}))) = \frac{\alpha^2\beta^2-2\alpha^2\beta-\alpha\beta^2+5\alpha\beta-\alpha-2\beta+1}{\alpha\beta}.$$
Since $z(UX^{-1}) = z(U)z(X)^{-1}=1$ and $\Phi = I_d$, it follows that
$$\alpha^2\beta^2-2\alpha^2\beta-\alpha\beta^2+5\alpha\beta-\alpha-2\beta+1=-2\alpha^2\beta-\alpha\beta^2+\alpha^2+5\alpha\beta+\beta^2-\alpha-2\beta.$$
This shows that $\alpha^2\beta^2+1=\alpha^2+\beta^2$, therefore $(\alpha^2-1)(\beta^2-1)=0$. This contradicts the conditions on $\alpha$ and $\beta$. This contradiction shows that $J=K\times \tilde{K}$.

We conclude using Goursat's Lemma with $R_{([1^2],[1^2])}(\mathcal{A}_{B_4})=SL_6(q)$ and then with \\$R_{[2,1],[1]}(\mathcal{A}_{B_4}) =SL_8(q)$. 
\end{proof}

We now show for $n\geq 5$ that if the representation associated with each double-partition is surjective, then $\Phi_{1,n}$ is surjective.

\begin{lemme}\label{gougou}
If $n\geq 5, \F_p(\alpha,\beta)=\F_p(\alpha+\alpha^{-1},\beta+\beta^{-1})$ and $\F_p(\alpha)=\F_p(\alpha+\alpha^{-1})$ and the composition of $\Phi_{1,n}$ and the projection upon each quasi-simple group associated with each double-partition is surjective, then $\Phi_{1,n}$ is surjective.
\end{lemme}

\begin{proof}
Let $n\geq 5$, we know by \cite[Theorem 1.1]{BMM} that the restriction to double-partitions with an empty component is surjective. We first show that we can add the hook partitions. We then show by induction on the double-partitions using the order we defined that $\Phi_{1,n}$ is surjective.

\bigskip

We write $G_{0,0}=SL_{n-1}(\tilde{q})\times \underset{(\lambda_1,\emptyset)\in A\epsilon_n, \lambda_1<\lambda_1'}\prod SL_{n_\lambda}(\tilde{q}) \times \underset{(\lambda_1,\emptyset)\in A\epsilon_n,\lambda_1=\lambda_1'}\prod OSP'(\lambda)$, where $OSP'(\lambda)$ is the derived subgroup of the group of isometries of the $\F_{\tilde{q}}$-bilinear form defined in \cite{BMM}, which identifies to the one defined in this article. We then have by Theorem 1.1. of \cite{BMM} that the image of $\mathcal{A}_{B_n}$ inside $G_{0,0}$ is $G_{0,0}$.  We have $G_{0,1}= R_{[1],[n-1]}(\mathcal{A}_{B_4}) = SL_n(q)$ by Proposition \ref{lesgourgues}. We use Goursat's Lemma to show that the image of $\mathcal{A}_{B_n}$ inside $G_{0,0}\times G_{0,1}$ is equal to $G_{0,0} \times G_{0,1}$. Using the notations in Goursat's Lemma, we have $K_1=G_{0,0}, K_2=G_{0,1}$ and $K_1/K^1\simeq K_2/K^2$. If the quotients are abelian then we are done since the groups we consider are perfect. We assume that they are non-abelian and show there is a contradiction. The only non-abelian decomposition factor of  $K_2$ is $PSL_n(q)$. Since the finite classical simple groups are non-isomorphic as long as $n\geq 4$ and $q\geq 4$ \cite[Section 1.2]{W}, there would exist a decomposition factor of $K_1$ corresponding to a double-partition $\lambda$ of $n$ with its right component empty such that $PG(\lambda)=\overline{R_{\lambda}(\mathcal{A}_{B_n})} \simeq PSL_n(q)=\overline{R_{([1],[n-1])}(\mathcal{A}_{B_n})}$. Therefore, up to conjugation \cite[Section 3.3.4]{W}, we have that $R_{\lambda}(h) = S^\Phi(R_{([1],[n-1])}(h))z(h)$ for all $h\in \mathcal{A}_{B_n}$ with $z :\mathcal{A}_{B_n}\rightarrow \F_q^\star$, $\Phi$ an automorphism of $\F_q$ and $S$  being either the identity or the transpose of the inverse. Since $n\geq 5$, we have $\mathcal{A}_{B_n}$ perfect \cite[Corollary]{MR}, therefore $z$ is trivial. By Lemma \ref{abel} and since the abelianization of  $A_{B_n}$ is the group $\{\overline{T},\overline{S_1}\} \simeq \Z^2$, we have $R_{\lambda}(h) = S^\Phi(R_{([1],[n-1])}(h))u^{\ell_1(h)}v^{\ell_2(h)}$. Since the right component of $\lambda$ is empty, the only eigenvalue of $R_{\lambda}(T)$ is $\beta$. On the other hand, the eigenvalues of $S^\phi(R_{([1],[n-1])}(T))v$ are equal to $\{v\Phi(\beta), -v\}$ or $\{v\Phi(\beta^{-1}),-v\}$. Therefore we would have $-v=v\Phi(\beta)$ or $-v=v\Phi(\beta^{-1})$ which is not possible because we have $\beta\neq-1$. This contradiction shows that the image is equal to $G_{0,0}\times G_{0,1}$.

Assume now $G_{0,2} = G_{0,0}\times G_{0,1}$ and $G_{0,3}=R_{([1],[1^n-1])}(\mathcal{A}_{B_n}) =SL_n(q)$ and consider the image of $\mathcal{A}_{B_n}$ inside $G_0=G_{0,2}\times G_{0,3}$. We use Goursat's Lemma with $K_1= G_{0,2}$ and $K_2=G_{0,3}$. In the same way as before, it is sufficient to show that the quotients $K_1/K^1\simeq K_2/K^2$ are abelian. The sets of eigenvalues of $S^{\Phi}(R_{([1],[1^{n-1}])}(T)$ are again $\{\Phi(\beta),-1\}$ or $\{\Phi(\beta)^{-1},-1\}$. If the quotients were non-abelian, we would have $R_{([1],[n-1])}(h) =S^{\phi}(R_{([1],[1^{n-1}])}(h))z(h)$ for all $h\in \mathcal{A}_{B_n}$ with $S$, $
\phi$ and $z$ as before. We have $z$ trivial since $\mathcal{A}_{B_n}$ is perfect, therefore $R_{([1],[n-1])}(h) =S^{\phi}(R_{([1],[1^{n-1}])}(h))u^{\ell_1(h)}v^{\ell_2(h)}$. Let us show that $\Phi$ is trivial. We have that $R_{([1],[n-1])}(T)$ and $R_{([1],[1^{n-1}])}(T)$ both have for eigenvalues $-1$ with multiplicity $n-1$ and $\beta$ with multiplicity $1$. This shows that either $\beta=v\Phi(\beta)$ and $-1=-v$ or  $\beta=v\Phi(\beta)^{-1}$ and $-1=-v$. In both cases, $v=1$ and $\Phi(\beta+\beta^{-1})=\beta+\beta^{-1}$. The eigenvalues of $R_{([1],[n-1])}(S_1)$ are $-1$ with multiplicity $1$ and $-\alpha$ with multiplicity $n-1$. The eigenvalues of $R_{([1],[1^{n-1}])}(S_1)$ are $-1$ with multiplicity $n-1$ and $\alpha$ with multiplicity $1$. Therefore we have either $-1=u\Phi(\alpha)$ and $\alpha=-u$ or $-1=u\Phi(\alpha^{-1})$ and $\alpha = -u$. In both cases $u=-\alpha$ and $\Phi(\alpha+\alpha^{-1})=\alpha+\alpha^{-1}$. We have $\Phi$ trivial, therefore $\F_q=\F_p(\alpha+\alpha^{-1},\beta+\beta^{-1})$. This would imply $R_{([1],[n-1])|\mathcal{A}_{B_n}} \simeq S(R_{([1],[1^{n-1}])|\mathcal{A}_{B_n}})$ but $([1],[1^{n-1}])\notin \{([1],[n-1]),([1],[n-1])'\}$ when $n >2$. By Proposition \ref{isomorphisme}, this is absurd. This shows the image of $\mathcal{A}_{B_n}$ in $G_0$ is equal to $G_0$.

For $\lambda_0 \in \epsilon_n=\{\lambda \vdash\vdash n, \lambda \notin A_n, \lambda~\mbox{not a hook}\}$, we set
$$G_{\lambda_0} = SL_{n-1}(\tilde{q})\times \underset{(\lambda_1,\emptyset)\in A\epsilon_n, \lambda_1<\lambda_1'}\prod SL_{n_\lambda}(\tilde{q}) \times \underset{(\lambda_1,\emptyset)\in A\epsilon_n,\lambda_1=\lambda_1'}\prod OSP'(\lambda)\times$$ $$SL_n(q)^2\times \underset{\lambda\in \epsilon_n, \lambda < min(\lambda',\lambda_0)}\prod SL_{n_\lambda}(q) \times \underset{\lambda\in \epsilon_n, \lambda=\lambda'< \lambda_0}\prod OSP'(\lambda).$$
where $OSP(\lambda)$ is the group of isometries of the bilinear form defined before Proposition \ref{bilin}.

For the minimal element $\lambda_0$ of $\epsilon_n$, we just showed the composition $\Phi_{1,n}$ with the projection onto $G_{\lambda_0}=G_0$ is surjective. Let us show by induction (numbering the double-partitions of $n$ with the order defined previously) that for all $\lambda_0$, the composition of $\Phi_{1,n}$ with the projection onto $G_{\lambda_0}$ is surjective.

Let $\lambda_0\in \epsilon_n$. Assume that the composition is surjective onto $G_{\lambda_0}$ and let us show that the composition onto $G_{\lambda_0+1}= G_{\lambda_0}\times G(\lambda_0)$ is surjective, where $G(\lambda_0)= SL_N(q)$ if $\lambda_0\neq\lambda_0'$ and $G(\lambda_0) \simeq OSP'(\lambda_0)\in \{SP_N(q),\Omega_N^{+}(q)\}$ if $\lambda_0=\lambda_0'$. We use Goursat's Lemma with $K_1=G_{\lambda_0}$ and $K_2=G(\lambda_0)$ on the image $\Phi_{1,n}$ in $K_1\times K_2$. As before, it is sufficient to show that the quotients $K_1/K^1 \simeq K_2/K^2$ are abelian. Assume they are non-abelian. The only non-abelian Jordan-Hölder factor of $G(\lambda_0)$ is $PG(\lambda_0)$, therefore there exists $\lambda$ less than $\lambda_0$ such that up to conjugation (see \cite{W} 3.3.4., 3.5.5. and 3.7.5) $\overline{R_{\lambda}(h)}=S^{\Phi}(R_{\lambda_0}(h)z(h)$ for all $h\in \mathcal{A}_{B_n}$ (there is no triality involved since if $n\geq 5$, $\lambda=\lambda'$ and $\lambda\in \epsilon_n$ then $\dim(V_\lambda) > \dim(V_{([2,1],[2,1])})=80>8$). By the same arguments as in the induction initialization, we have that $\lambda$ has no empty components. Since $n\geq 5$, $\mathcal{A}_{B_n}$ is perfect and $z$ is trivial. We then have $R_{\lambda_0|\mathcal{A}_{B_n}}\simeq S^{\Phi}(R_{\lambda|\mathcal{A}_{B_n}})$. Let us show that $\Phi$ is trivial. By Lemma \ref{abel}, there exists $u,v\in \F_q^\star$ such that up to conjugation, for all $h\in A_{B_n}$, we have $R_ {\lambda_0}(h)= S^\Phi(R_\lambda(h))u^{\ell_1(h)}v^{\ell_2(h)}$. Comparing eigenvalues of $T$, we get either $\{\beta,-1\}=\{v\Phi(\beta),-v\}$ or $\{\beta,-1\}=\{v\Phi(\beta^{-1}),-v\}$. In the first case, either  $v=1$ and $\Phi(\beta)=\beta$ or $v=-\beta$ and $-1=v\Phi(\beta)$, therefore $\Phi(\beta+\beta^{-1})=\beta+\beta^{-1}$.
In the second case either $v=1$ and $\Phi(\beta^{-1})=\beta$ or $v=-\beta$ and $v\Phi(\beta^{-1}) = -1$, therefore $\Phi(\beta+\beta^{-1})=\beta+\beta^{-1}$. In the same way using $S_1$, we show $\Phi(\alpha+\alpha^{-1})= \alpha+\alpha^{-1}$. This shows that $\Phi$ is trivial because $\F_q=\F_p(\alpha+\alpha^{-1},\beta+\beta^{-1})$. We then have $R_{\lambda_0|\mathcal{A}_{B_n}}\simeq S(R_{\lambda|\mathcal{A}_{B_n}})$ which contradicts Proposition \ref{isomorphisme} since $\lambda< \lambda_0\leq \lambda_0'$.
\end{proof}

\bigskip

To get that $\Phi_{1,n}$ is surjective, it now only remains to show that what we assumed in Lemma \ref{gougou} is true.

\begin{theo}\label{hortillonage}
If $n\geq 5$ then for all $\lambda\vdash\vdash n$ double-partitions in our decomposition, we have $R_\lambda(\mathcal{A}_{B_n})= G(\lambda)$, where $G(\lambda)$ is the corresponding group in the following list.
\begin{enumerate}
\item $SL_{n-1}(\tilde{q})$ if $\lambda =([n-1,1],\emptyset)$.
\item $SL_N(\tilde{q})$ if $\lambda = (\lambda_1,\emptyset), \lambda_1< \lambda_1'$.
\item $SP_N(\tilde{q})$ if $\lambda = (\lambda_1,\emptyset), \lambda_1=\lambda_1'$ and ( $p=2$ or ($p\geq 3$ and $\nu(\lambda_1) = -1$)).
\item $\Omega_N^{+}(\tilde{q})$ if $\lambda=(\lambda_1,\emptyset), \lambda_1=\lambda_1'$, $p\geq 3$ and $\nu(\lambda_1)=1$.
\item $SL_n(q)$ if $\lambda \in \{([1],[n-1]),([1],[1^{n-1}])\}$.
\item $SL_N(q)$ if $\lambda \in \epsilon_n, \lambda< \lambda'$.
\item $SP_N(q)$ if $\lambda=\lambda'$ and ($p=2$ or $\tilde{\nu}(\lambda)=-1)$.
\item $\Omega_N^{+}(q)$ if $\lambda = \lambda', p\geq 3$ and $\tilde{\nu}(\lambda)= 1$.
\end{enumerate}
\end{theo}

\begin{proof}
Let $n\geq 5$. By \cite{BM} (Theorem 1.1.), it is sufficient to show it for $\lambda \in \epsilon_n$.

The result is true for $n=4$, therefore we can use induction and assume $\Phi_{n-1}$ is surjective.

The first thing to do is to take care of the double-partitions such that $n_\lambda>8$ and $n_\lambda\neq 10$. For $n=5$, the double partitions to consider are $([1^3],[1^2])$, $([1^2],[1^3])$, $([1],[2,2])$, $([1],[2,1^2])$, $([1],[3,1])$, $([2],[2,1])$ and $([1^2],[2,1])$ of respective dimensions $10$, $10$, $10$, $15$, $15$, $20$ and $20$. For $n= 6$, they are $([1^2],[1^4])$, $([1^4],[1^2])$, $([1^3],[1^3])$, $([1],[4,1])$, $([1],[2,1^3])$, $([1],[3,2])$, $([1],[2^2,1])$, $([1],[3,1,1])$, $([2],[2^2])$, $([1^2],[2^2])$, $([2],[3,1])$, $([1^2],[3,1])$, \break$([2],[2,1^2])$, $([1^2],[2,1^2])$, $([3],[2,1])$, $([1^3],[2,1])$ and $([2,1],[2,1])$ of respective dimensions $15$, $15$, $20$, $24$, $24$, $30$, $30$, $36$, $30$, $30$, $45$, $45$, $45$, $45$, $40$, $40$ and $80$. We can now note that if $n=6$, we have $n_\lambda \geq 15$, therefore by the branching rule, if $n\geq 6$ and $\lambda$ is a double-partition of $\epsilon_n$ then $n_\lambda \geq 15$. The only double-partitions $\lambda$ such that $n_\lambda\leq 8$ or $n_\lambda=10$ are $([1^3],[1^2]), ([1^2],[1^3])$ and $([1],[2,2])$ which are of dimension $10$. By Lemma \ref{platypus} and the branching rule, we have that $R_{[1],[2,2]}(\mathcal{A}_{B_4}) = SL_8(q) \times SL_2(\tilde{q})$, $R_{[1^3],[1^2]}(\mathcal{A}_{B_4})=SL_4(q)\times SL_6(q)$ and $R_{[1^3],[1^2]}(\mathcal{A}_{B_4})=SL_4(q)\times SL_6(q)$. By Theorem \ref{LBJ}, we have that $R_{([1],[2,2])}(\mathcal{A}_{B_5})=R_{[1^3],[1^2]}(\mathcal{A}_{B_5})=R_{[1^2],[1^3]}(\mathcal{A}_{B_5}) \simeq SL_{10}(q)$.

We now show in the same way as in \cite{BMM} (Part 5), that the other assumptions of Theorem \ref{CGFS} are verified. In order to do this, we use the following results shown in \cite{BMM}.

\begin{lemme}\label{tens1}
If $d\geq 6$ and $G\leq GL_d(q)$ contains an element conjugate to an element of the form $\diag(\xi,\xi^{-1},1,1,...)$ with $\xi^2\neq 1$, then $G$ is tensor-indecomposable.
\end{lemme}

\begin{lemme}\label{tens2}
If $d\geq 16$ and $G\leq GL_d(q)$ contains an element of order prime to $p$ conjugate to an element of the form $\diag(\xi,\xi,\xi^{-1}, \xi^{-1},1,1,..,1)$ with $\xi^2\neq 1$, then $G$ is tensor-indecomposable except possibly if  $G \leq G_1\otimes G_2$ with $G_1\leq GL_2(q)$.
\end{lemme}

For a block diagonal matrix with blocks $B_1,\dots,B_r$, we write $\diag(B_1,\dots,B_r)$.

\begin{lemme}
If $G$ contains a natural $SL_2(q)$ and $q \geq 8$ or $G$ contains a twisted diagonal embedding of  $SL_2(q)$ $(G \supset \{\diag(M,{}^t\!(M^{-1}),I_{N-4}), M\in SL_2(q)\})$, then case ($2$) of Theorem \ref{CGFS} is excluded.
\end{lemme}

By the proof of the imprimitivity of $G$ in \cite{BMM}, it is sufficient to show that $\mathcal{A}_{B_n}$ is normally generated by $\mathcal{A}_{B_{n-1}}$ and that $G$ contains either a transvection or an element of Jordan form $\diag(I_2+E_{1,2},I_2+E_{1,2},I_{N-4})$ to get that $G$ is imprimitive.

In order to show that we are in case ($1$) of Theorem \ref{CGFS}, we must show that for $n\geq 5$, that we have $q\geq 8$ and that for any double-partition $\lambda$ of $n$, $G=R_\lambda(\mathcal{A}_{B_n})$ contains either a natural $SL_2(q)$ and  $n_\lambda> 6$ or contains a twisted diagonal embedding of $SL_2(q)$ and $n_\lambda > 16$. We must also prove that $\mathcal{A}_{B_n}$ is normally generated $\mathcal{A}_{B_{n-1}}$ and the exceptional case of Lemma \ref{tens2} is impossible when $n_\lambda > 16$, $G$ contains a twisted diagonal embedding of $SL_2(q)$ but no natural $SL_2(q)$ in an obvious way.

Let $n\geq 5$, assume the lemma is true for all $m\leq n-1$. By Lemmas \ref{gougou} and \ref{platypus}, we have $\Phi_m$ surjective for all $m\leq n-1$.  By assumption, $\alpha$ is of order strictly greater than $5$ and not belonging to $\{1,2,3,4,5,6,8,10\}$. This implies that $\alpha$ is of order at least $7$ and that $q\geq 8$. If $\lambda$ has at most two columns then since $\lambda\in \epsilon_n, \lambda$ contains a natural $SL_2(q)$. Assume now $\lambda\in \epsilon_n$ has at least three rows or three columns.

Assume that for all $\mu \subset \lambda$ containing $([2,1],[1])$ or $([1],[2,1])$, we have $\mu' \subset \lambda$. We then have that $\lambda = \lambda'$ and $n$ is even. Since $n$ is even, we have $\mu \neq \mu'$ for any double-partition $\mu\subset \lambda$. Since $\lambda \in \epsilon_n$ and contains strictly more than two rows and two columns, there exists $\mu\subset \lambda$ containing $([1],[2,1])$ or $([2,1],[1])$. Since $\Phi_{m-1}$ is surjective and $\mu'\subset \lambda$, we have a twisted diagonal embedding of $SL_{n_\mu}(q)$ in $G= R_\lambda(\mathcal{A}_{B_n})$ and since $n_\mu \geq 8 \geq 2$, we have a twisted diagonal embedding of $SL_2(q)$. Otherwise there exists $\mu\subset \lambda$ containing $([2,1],[1])$ or $([1],[2,1])$ such that $\mu' \subset \lambda$. Since $\Phi_{m-1}$ is surjective, we get that $\lambda$ contains a natural $SL_{n_\mu}(q)$ and, therefore contains a natural $SL_2(q)$. For double-partitions which are not of dimension strictly greater than $16$, i.e. $([1],[2,1^2])$ and $([1],[3,1])$, we are in the second case.

We now show that $\mathcal{A}_{B_n}$ is normally generated by $
\mathcal{A}_{B_{n-1}}$ for $n\geq 5$. By \cite[Lemma 2.1]{BMM}, we have that $\mathcal{A}_{A_n}$ is normally generated by $\mathcal{A}_{A_{n-1}}$ for $n\geq 4$. Since $T$ commutes with $S_i$ for all $i \geq 2$, we have the same result for $\mathcal{A}_{B_n}$ for $n\geq 4$.

It now only remains to show that the exception of Lemma \ref{tens2} is impossible when there is no obvious natural $SL_2(q)$ in $G$. In order to do this, we show a proposition analogous to Proposition 2.4. of \cite{BMM}.

\begin{prop}
Let $K$ be a field. If $n\geq 7$ and $\varphi : \mathcal{A}_{B_n} \rightarrow PSL_2(K)$ is a group morphism then $\varphi =1$.
\end{prop}

\begin{proof}
Let $K$ be a field, $n\geq 7$ and $\varphi$ such a morphism. The restriction of $\varphi$ to $\mathcal{A}_{A_n} \leq \mathcal{A}_{B_n}$ is trivial by Proposition 2.4. of \cite{BMM}. By Theorem 3.9. of \cite{MR}, $\mathcal{A}_{B_n}$ is generated by $p_0=S_{n-2}S_{n-1}^{-1}, p_1=S_{n-1}S_{n-2}S_{n-1}^{-2}, q_3=S_{n-3}S_{n-1}^{-1}, b = S_{n-2}S_{n-1}^{-1}S_{n-3}S_{n-2}^{-1}, r_l=T^lS_1T^{-l}S_{n-1}, q_i = S_{n-i}S_{n-1}^{-1}, l\in \Z, 4\leq i \leq n-2$ and the following relations.
\begin{enumerate}
\item For $4\leq j\leq n-2, p_0q_j =q_jp_1$ and $p_1q_j=q_jp_0^{-1}p_1$.
\item For $l\in \Z, p_0r_l = r_lp_1$ and $p_1r_l=r_lp_0^{-1}p_1$.
\item For $3\leq i <j\leq n-2, \vert i-j\vert \geq 2, q_iq_j=q_jq_i$.
\item For $3\leq i \leq n-3, q_ir_l=r_lq_i$.
\item $p_0q_3p_0^{-1}=b, p_0bp_0^{-1}=b^2q_3^{-1}b$
\item $p_1q_3p_1^{-1} = q_3^{-1}b, p_1bp_1^{-1}=(q_3^{-1}b)^3q_3^{-2}b$.
\item For $3\leq i \leq n-3, q_iq_{i+1}q_i = q_{i+1}q_iq_{i+1}$.
\item For $l\in \Z, q_{n-2}r_lq_{n-2} = r_lq_{n-2}r_l$.
\item For $l\in \Z, r_lr_{l+1}=r_{l+1}r_{l+2}$.
\end{enumerate}

By \cite[Proposition 2.4]{BMM}, the images of all the generators except for $(r_l)_{l\in \Z}$ are trivial. By the eighth relation, we get that the images of the $r_l$ are also trivial and the desired result follows.
\end{proof}

This shows that if $n\geq 7$ and $G \leq G_1\otimes G_2$ with $G_1\leq GL_2(q)$, then $G \subset SL_{\frac{N}{2}}(q) \times SL_{\frac{N}{2}}(q)$. This contradicts the irreducibility. Since we need $n\geq 7$ to apply this reasoning, we must consider separately the cases where $n\in \{5,6\}$ and $G$ does not contain a natural $SL_2(q)$. Looking at all the cases enumerated previously, the only one to consider is $\lambda=([2,1],[2,1])$. Up to conjugation, we have $H=R_{[2,1],[2,1]}(\mathcal{A}_{B_5}) \simeq  \{\diag(M,{}^t\!(M^{-1}), N,{}^t\!(N^{-1})), M,N\in SL_{20}(q)\} \simeq SL_{20}(q) \times SL_{20}(q)$.

Assume that $G=R_{([2,1],[2,1])}(\mathcal{A}_{B_6})  \subset G_1\otimes G_2$ and that $G_1\subset GL_2(q)$. We then have a morphism $\theta$ from $G$ to $SL_2(q)$ since $\mathcal{A}_{B_n}$ is perfect for $n\geq 5$. If we consider the restriction of $\theta$ to $H$, its kernel is a subgroup of $H$ and its image is a subgroup of $SL_2(q)$. Since $PSL_{20}(q)$ is the only non-abelian composition factor of $H$, we have that if the image is non-abelian then there exists a subgroup of $SL_2(q)$ isomorphic to $PSL_{20}(q)$. This is absurd, therefore the image is abelian and the kernel contains the derived subgroup of $H$ which is equal to $H$ since $H$ is perfect. In the same way, for all $g\in G$, the restriction of $\theta$ to $gHg^{-1}$ is trivial. Since $H$ normally generates $G$, $\theta$ is trivial which contradicts the irreducibility of $G$ in the same way as in the proof of the previous proposition.

We have thus shown that we are in the first case of Theorem \ref{CGFS}. By the same reasoning as in \cite[page 16]{BMM}, we have in all cases that $q'=q$. If $\lambda=\lambda'$, we have $G \subset G(\lambda)$ by Proposition \ref{transpose}, therefore $G =G(\lambda)$. If $\lambda \neq \lambda'$, $G$ preserves no bilinear form since $R_\lambda$ is not isomorphic to $R_\lambda^\star$.  If $G$ preserves a hermitian form then there exists an automorphism $\Phi$ of order $2$ of $\F_q$ such that $M$ is conjugate to $\Phi({}^t\!(M)^{-1})$ for all $M\in G$. Since $G$ contains a natural $SL_2(q)$, we then have $\tr(\diag(\alpha,\alpha^{-1},1,1,...,1)) = \Phi(\tr(\diag(\alpha^{-1}, \alpha,1,...,1)))$ and $\tr( \diag(\beta,\beta^{-1},1,1,...,1)) = \Phi(\tr(\diag(\beta^{-1}, \beta,1,...,1)))$, therefore $\Phi(\alpha+\alpha^{-1}) =\alpha+\alpha^{-1}$ and $\Phi(\beta+\beta^{-1})= \beta+\beta^{-1}$. This implies that $\Phi=I_d$ because $\F_q=\F_p(\alpha+\alpha^{-1}, \beta+\beta^{-1})$. This is absurd and we conclude that $G =SL_{n_\lambda}(q)$.
\end{proof}

By Theorem \ref{hortillonage}, Lemma \ref{platypus} and Lemma \ref{gougou}, we have that for all $n, \Phi_{1,n}$ is surjective.

\subsection{Cases (2) and (3)}

We have shown the surjectivity of $\Phi_{1,n}$, this corresponds to the first of the six possible field extension configurations described at the beginning of subsection \ref{lalala}. The proof in cases ($2$) and ($3$) only requires small changes to the one in the first case, but the new factorizations appearing in cases ($4$) to ($6$) require more work, especially for the low dimensional representations. We treat in this section cases ($2$) and ($3$) emphasizing on the differences with the first case. This will conclude the proof of Theorem \ref{result3} and the corresponding statement in the second of the six cases listed at the beginning of subsection \ref{lalala}.

In case (2), i.e., $\F_q=\F_p(\alpha,\beta) = \F_p(\alpha+\alpha^{-1},\beta+\beta^{-1})$ and $\F_p(\alpha) \neq \F_p(\alpha+\alpha^{-1})$, the same arguments as the ones in case 1 work at every step of the proof. Indeed, $SU_2(\tilde{q}^{\frac{1}{2}})$ is also generated by a conjugacy class of transvections. Since $\tilde{q}$ is a square and $\alpha$ is order not diving $8$ by assumption, we have that $\tilde{q}\geq 16$ and $\tilde{q}^{\frac{1}{2}}> 3$. We also still have that $SL_8(q)\times SU_2(\tilde{q}^{\frac{1}{2}})$ contains $SL_8(q) \times \{I_2\}$, therefore all the arguments work in the same way. This shows that in case ($2$), $\Phi_{2,n}$ is surjective for all $n$.

In case ($3$), i.e., $\F_q=\F_p(\alpha,\beta)=\F_p(\alpha+\alpha^{-1},\beta)= \F_p(\alpha,\beta+\beta^{-1})\neq \F_p(\alpha+\alpha^{-1},\beta+\beta^{-1})$, all representations are unitary. The main differences occur in the proof that when $n=4$, the direct product of two $SU_4(q^\frac{1}{2})$ is in the image, and in the conclusion of the proof of this version of Theorem \ref{hortillonage}.
\begin{theo}
If $n\geq 5$, then for all $\lambda\vdash\vdash n$ in our decomposition, $R_\lambda(\mathcal{A}_{B_n})= G(\lambda)$, where $G(\lambda)$ is the corresponding group in the following list.
\begin{enumerate}
\item $SU_{n-1}(\tilde{q}^{\frac{1}{2}})$ if $\lambda =([n-1,1],\emptyset)$.
\item $SU_N(\tilde{q}^{\frac{1}{2}})$ if $\lambda = (\lambda_1,\emptyset), \lambda_1< \lambda_1'$.
\item $SP_N(\tilde{q}^{\frac{1}{2}})$ if $\lambda = (\lambda_1,\emptyset), \lambda_1=\lambda_1'$ and ( $p=2$ or ($p\geq 3$ and $\nu(\lambda_1) = -1$)).
\item $\Omega_N^{+}(\tilde{q}^\frac{1}{2})$ if $\lambda=(\lambda_1,\emptyset), \lambda_1=\lambda_1'$, $p\geq 3$ and $\nu(\lambda_1)=1$.
\item $SU_n(q^\frac{1}{2})$ if $\lambda \in \{([1],[n-1]),([1],[1^{n-1}])\}$.
\item $SU_N(q^\frac{1}{2})$ if $\lambda \in \epsilon_n, \lambda< \lambda'$.
\item $SP_N(q^\frac{1}{2})$ if $\lambda=\lambda'$ and ($p=2$ or ($p\geq 3$ and $\tilde{\nu}(\lambda)=-1$.
\item $\Omega_N^{+}(q^\frac{1}{2})$ if $\lambda = \lambda', p\geq 3$ and $\tilde{\nu}(\lambda)= 1$.
\end{enumerate}
\end{theo}

\begin{proof}
 We recall Proposition $4.1.$ of \cite{BMM}.
\begin{prop}\label{coolprop}
Let $q=u^2$, $\varphi$ be a non-degenerate bilinear form over $\F_q^N$, $\psi$ a non-degenerate hermitian form over $\F_q^N$. If $G\subset OSP_N(\varphi) \cap U_N(\psi)$ is absolutely irreducible, then there exists $x\in GL_N(q)$ and a non-degenerate bilinear form $\varphi'$ over $\F_u^N$ such that $^x G \subset OSP(\varphi')$ and $\varphi'$ is of the same type as $\varphi$.
\end{prop}

When $n=4$, the proof that $\Phi_4$ is surjective is the same up to the point, where we prove $\Phi$ is trivial using $\Phi(\alpha+\alpha^{-1})=\alpha+\alpha^{-1}$ and $\Phi(\beta+\beta^{-1})=\beta+\beta^{-1}$. In case $3$, $\Phi$ could also be equal to the automorphism $\epsilon$ of order $2$ of $\F_q$. It is thus necessary to show that the following is absurd :
$$\frac{\alpha^2\beta^2-2\alpha^2\beta-\alpha\beta^2+5\alpha\beta-\alpha-2\beta+1}{\alpha\beta}=\epsilon(\frac{-2\alpha^2\beta-\alpha\beta^2+\alpha^2+5\alpha\beta+\beta^2-\alpha-2\beta}{\alpha\beta}).$$
This would imply
\begin{tiny}
\begin{eqnarray*}
\alpha^2\beta^2-2\alpha^2\beta-\alpha\beta^2+5\alpha\beta-\alpha-2\beta+1 &  = &\frac{-2\alpha^{-2}\beta^{-1}-\alpha^{-1}\beta^{-2}+\alpha^{-2}+5\alpha^{-1}\beta^{-1}+\beta^{-2}-\alpha^{-1}-2\beta^{-1}}{\alpha^{-2}\beta^{-2}}\\
& = & -2\beta-\alpha+\beta^2+5\alpha\beta+\alpha^2-\alpha\beta^2-2\alpha^2\beta.
\end{eqnarray*}
\end{tiny}
This is absurd because it is the same equality we proved to be impossible in the first case.

We now adapt the end of the proof of the corresponding version of Theorem \ref{hortillonage}.  By \cite[page 18]{BMM},  we are in case ($1$) of Theorem \ref{CGFS}. If $\lambda \neq \lambda'$, $G$ contains a natural $SU_3(q^\frac{1}{2})$, therefore $q=q'$ by Lemma \ref{field}. Since $G\subset SU_{n_\lambda}(q^\frac{1}{2})$ and $G$ preserves no bilinear form by Proposition \ref{isomorphisme}, we have $G\simeq SU_{n_\lambda}(q^{\frac{1}{2}})$. If $\lambda=\lambda'$, we use Proposition \ref{coolprop} to get that $G\subset OSP(q^\frac{1}{2})$. By Lemma \ref{field}, we have that $\F_{q'}$ contains $\{x+\epsilon(x), x\in \F_q\}$. This implies that $q'=q^{\frac{1}{2}}$ because $\F_{q'}$ contains $\alpha+\alpha^{-1}$ and $\beta+\beta^{-1}$ and $q'$ divides $q^\frac{1}{2}$. We conclude that $G\simeq OSP'(q^\frac{1}{2})$.
\end{proof}

\subsection{Cases (4), (5) and (6)}

In this subsection, we finish the proof for type $B$ by considering the last three cases for the field extensions listed at the beginning of subsection \ref{lalala}. This will conclude the proofs of Theorems \ref{result4}, \ref{result5} and \ref{result6}. In these cases more factorizations appear and this complicates greatly the proof for small $n$. We will use the tables of maximal subgroups of finite classical groups in low dimension from \cite{BHRC}. This gives interesting techniques to determine if a certain subgroup $G$ of a classical group is the group itself, when given information on the subgroups of $G$.

In these cases, we can still use various arguments from the first case, but except for Proposition \ref{lesgourgues} which remains true in all these cases, all the low-dimensional cases must be done again. It is not necessary to use new arguments for Lemma \ref{gougou}. We start by studying the case $n\leq 4$.

\begin{lemme}\label{platypus456}
For $i\in \{4,5,6\}$ and $n\leq 4$, we have $\Phi_{i,n}$ surjective.
\end{lemme}

\begin{proof}
For $n=2$, using the same arguments as in the first case and  Lemma \ref{Ngwenya}, we have that $\mathrm{Im}(\Phi_2)\simeq SL_2(q^\frac{1}{2})$.

For $n=3$, we have by the factorizations in Proposition \ref{isomorphisme} that in all cases $\Phi_3$ is surjective.

The only case left to consider is $n=4$ and the double-partitions $ ([1^2],[1^2])$ and $([1],[2,1])$ of respective dimensions $6$ and $8$. We have to prove that $R_{([1^2],[1^2])}(\mathcal{A}_{B_4})\simeq SL_6(q^\frac{1}{2})$ or $SU_6(q^{\frac{1}{2}})$ and $R_{([2,1],[1])}(\mathcal{A}_{B_4}) \simeq SL_8(q^\frac{1}{2})$ or $SU_8(q^\frac{1}{2})$ depending on the case.

We start by $G=R_{([1^2],[1^2])}(\mathcal{H}_4)$ in case ($5$) or ($6$), where $\F_q = \F_p(\alpha,\beta) = \F_p(\alpha+\alpha^{-1},\beta)\neq  \F_p(\alpha,\beta+\beta^{-1})= \F_p(\alpha+\alpha^{-1},\beta+\beta^{-1})$. We then have $H=R_{([1^2],[1^2])}(\mathcal{A}_{B_3}) \simeq SL_3(q)$. Since $([1^2],[1^2])=([1^2],[1^2])$, by Proposition \ref{isomorphisme} (4.c) and Lemma \ref{Ngwenya}, up to conjugation, we have $G\subset SL_6(q^\frac{1}{2})$.

We use the classification of maximal subgroups $SL_6(q^{\frac{1}{2}})$ \cite[Tables 8.24 and 8.25]{BHRC} . Using the fact that $H$ is a subgroup of $G$, we exclude the possibility that $G$ is included in all but two of these groups, using the divisibility of the orders that would ensue.  We start by considering the sporadic maximal subgroups in Table 8.25 and get the orders of these groups using the atlas \cite{CCNPW}. We list below those groups and their order or a quantity their order divides

\begin{enumerate}
\item $2 \times 3^{.} \mathfrak{A}_6.2_3, 4320$,  
\item $2 \times 3 ^{.} \mathfrak{A}_6, 2160$,
\item $6^{.}\mathfrak{A}_6, 2160$,
\item $(q^{\frac{1}{2}}-1,6) \circ 2^{.}PSL_2(11), 12\times 660=7920$,
\item $6^{.}\mathfrak{A}_7, 15120$,
\item $6^{.}PSL_3(4)^{.}2_1^{-}, 6\times 2\times 20160=241920$,
\item $6^{.}PSL_3(4), 6\times 20160=120960$,
\item $2^{.}M_{12}, 2\times 95040=190080$,
\item $6_1^{.}PSU_4(3)^{.}2_2^{-}, 6\times 2\times 3265920=39191040$,
\item $6_1^{.}PSU_4(3), 6\times 3265920=19595520$,
\item $(q^{\frac{1}{2}}-1,6) \circ SL_3(q^\frac{1}{2}), 6\times q^{\frac{3}{2}}(q-1)(q^{\frac{3}{2}}-1)$.
\end{enumerate}

 Since $q$ is a square and $\alpha$ is of order greater than $4$, we have $q\geq 9$. This implies that $\vert SL_3(q) \vert = q^3 (q^2-1)(q^3-1) \geq 9^3(9^2-1)(9^3-1) = 42456960$, which is greater than all the orders in the list (the last one is of order $6q^\frac{1}{2}(q-1)(q^3-1)$ and $6 < q^\frac{3}{2}(q+1)(q^\frac{3}{2}+1)$). We now  look at the list in table 8.24 of the 18 geometric maximal subgroups of $SL_6(q^{\frac{1}{2}})$, which we provide below with their order or a quantity which is divisible by their order
 
\begin{enumerate}
\item $E_{q^\frac{1}{2}}^5:GL_5(q^\frac{1}{2}), q^\frac{15}{2}(q^\frac{5}{2}-1)(q^2-1)(q^\frac{3}{2}-1)(q-1)(q^\frac{1}{2}-1)$,
\item $E_{q^\frac{1}{2}}^8:(SL_4(q^\frac{1}{2})\times SL_2(q^\frac{1}{2})):(q^\frac{1}{2}-1), q^\frac{15}{2}(q^2-1)(q^\frac{3}{2}-1)(q-1)^2(q^\frac{1}{2}-1)$,
\item $E_{q^\frac{1}{2}}^9:(SL_3(q^\frac{1}{2})\times SL_3(q^\frac{1}{2})):(q^\frac{1}{2}-1), q^\frac{15}{2}(q^\frac{3}{2}-1)^2(q-1)^2(q^\frac{1}{2}-1),$
\item $E_{q^\frac{1}{2}}^{1+8}:(GL_4(q^\frac{1}{2})\times (q^\frac{1}{2}-1)), q^\frac{15}{2}(q^2-1)(q^\frac{3}{2}-1)(q-1)(q^\frac{1}{2}-1)^2$,
\item $E_{q^\frac{1}{2}}^{4+8}:SL_2(q^\frac{1}{2})^3:(q^\frac{1}{2}-1)^2, q^\frac{15}{2}(q-1)^3(q^\frac{1}{2}-1)$,
\item $GL_5(q^\frac{1}{2}), q^5(q^\frac{5}{2}-1)(q^2-1)(q^\frac{3}{2}-1)(q-1)(q^\frac{1}{2}-1)$,
\item $(SL_4(q^\frac{1}{2})\times SL_2(q^\frac{1}{2})):(q^\frac{1}{2}-1), q^\frac{7}{2}(q^2-1)(q^\frac{3}{2}-1)(q-1)^2(q^\frac{1}{2}-1)$,
\item $(q^\frac{1}{2}-1)^5\times \mathfrak{S}_6$, $6!(q^\frac{1}{2}-1)^5$,
\item $SL_2(q^\frac{1}{2})^3:(q^\frac{1}{2}-1)^2.\mathfrak{S}_3, 6q^\frac{3}{2}(q^\frac{1}{2}-1)^2(q-1)^3$,
\item $SL_3(q^\frac{1}{2})^2:(q^\frac{1}{2}-1).\mathfrak{S}_2, 2q^\frac{3}{2}(q^\frac{3}{2}-1)(q-1)(q^\frac{1}{2}-1)$,
\item $SL_3(q).(q^\frac{1}{2}+1).2, 2(q^\frac{1}{2}+1)\vert SL_3(q)\vert$,
\item $SL_2(q^\frac{3}{2}).(q+q^\frac{1}{2}+1).3, 3q^\frac{3}{2}(q^3-1)(q+q^\frac{1}{2}+1)$,
\item $SL_2(q^\frac{1}{2})\times SL_3(q^\frac{1}{2}), q^2(q^\frac{3}{2}-1)(q-1)^2$,
\item $SL_6(q_0).[(\frac{q^\frac{1}{2}-1}{q_0-1},6)]$, where $q^\frac{1}{2}=q_0^r$ and $r$ prime, $6\vert SL_6(q_0)\vert$, 
\item $(q^\frac{1}{2}-1,3)\times SO_6^{+}(q^\frac{1}{2}).2,$ $q$ odd, $ 6q^3(q^2-1)(q^\frac{3}{2}-1)(q-1)$,
\item $(q^\frac{1}{2}-1,3)\times SO_6^{-}(q^\frac{1}{2}).2$, $q$ odd, $ 6q^3(q^2-1)(q^\frac{3}{2}+1)(q-1)$,
\item $(q^\frac{1}{2}-1,3)\times SP_6(q^\frac{1}{2}), 3q^\frac{9}{2}(q-1)(q^2-1)(q^3-1)=3q^\frac{3}{2}(q-1)\vert SL_3(q)\vert$,
\item $SU_6(q^{\frac{1}{4}}).(q^{\frac{1}{4}}-1,6), 6q^{\frac{15}{4}}(q^{\frac{3}{2}}-1)(q^{\frac{5}{4}}+1)(q-1)(q^{\frac{3}{4}}+1)(q^{\frac{1}{2}}-1)$.
\end{enumerate}

  In cases $1,2,3,4,5,6,7$ and $13$, the order of the maximal subgroup divides  $q^\frac{15}{2}(q^\frac{5}{2}-1)(q^2-1)(q^\frac{3}{2}-1)^2(q-1)^2(q^\frac{1}{2}-1)$. This implies that it is sufficient to show that $\vert SL_3(q)\vert=q^3(q^3-1)(q^2-1)$ does not divide this quantity to exclude these cases.  It can be true only if $q^3-1$ divides $(q^\frac{5}{2}-1)(q^\frac{3}{2}-1)^2(q-1)^2(q^\frac{1}{2}-1)$.  The Euclidean remainder of those two quantities seen as polynomials in $q^\frac{1}{2}$ is $4q^\frac{5}{2}+2q^2-2q^\frac{3}{2}-4q-2q^\frac{1}{2}+2$. Therefore, if $q^3-1$ divides the first quantity then it divides the remainder, which is positive. Therefore it is less than or equal to it.  We have $\epsilon(\alpha)=\alpha^{-1}=\alpha^{q^\frac{1}{2}}$, therefore $\alpha^{q^\frac{1}{2}+1}=1$. Since $\alpha$ is of order strictly greater than $6$ by assumption, we have that $q^\frac{1}{2}\geq 6$ and $4q^\frac{5}{2}+2q^2-q-2q^\frac{3}{2}-4q-2q^\frac{1}{2}+2\leq 4q^\frac{5}{2}+2q^2+2\leq 4q^\frac{5}{2}+3q^2\leq 5q^\frac{5}{2}<q^3-1$. This gives us the desired contradiction.
  
  \smallskip

Cases $8,9,10$ and $12$ are excluded because $q^3$ is coprime to $q^r-1$ and $(q+q^\frac{1}{2}+1)$ for every integer $r$ and $q^3$ does not divide $6!$ or $6q^\frac{3}{2}$ since $q = (q^\frac{1}{2})^2\geq 36$.

\smallskip

In case $14$, the order of the maximal subgroup $M$ divides the quantity 
$$6q^{\frac{15}{2r}}(q^{\frac{3}{r}}-1)(q^{\frac{5}{2r}}-1)(q^{\frac{2}{r}}-1)(q^{\frac{3}{2r}}-1)(q^{\frac{1}{r}}-1),$$ where $q^{\frac{1}{2r}}=q_0$ and $r$ is a prime. If $\vert SL_3(q)\vert $ divides this quantity, then $q^3$ divides $6q^{\frac{15}{2r}}$.
If $r\geq 3$, then $6q^{\frac{15}{2r}}\leq 6q^{\frac{15}{6}}< q^3$ because $q^{\frac{1}{2}}\geq 8 > 6$ when $r\geq 3$.

It only remains to consider the case $r=2$. We then have that $(q^\frac{3}{2}-1)(q^\frac{3}{2}+1)(q-1)$ divides $6(q^{\frac{3}{2}}-1)(q^{\frac{5}{4}}-1)(q-1)(q^{\frac{3}{4}}-1)(q^{\frac{1}{2}}-1)$, therefore $q^{\frac{3}{2}}+1$ divides $6(q^{\frac{5}{4}}-1)(q^{\frac{3}{4}}-1)(q^{\frac{1}{2}}-1)$. The Euclidean remainder of the division of those two polynomials in 
$q^\frac{1}{4}$ is $-6q^{\frac{4}{4}}+6q^{\frac{3}{4}}+12q^\frac{2}{4}+6q^\frac{1}{4}-6$. This implies that $q^{\frac{3}{2}}+1$ divides the above quantity. We have $q^\frac{1}{2}\geq 6$ and it is a square when $r=2$. Therefore $q^\frac{1}{2}\geq 9$ and $q^\frac{1}{4}\geq 3$. It follows that $-6q^{\frac{4}{4}}+6q^{\frac{3}{4}}+12q^\frac{2}{4}+6q^\frac{1}{4}-6\leq -6q^\frac{4}{4}+11q^\frac{3}{4}\leq -4q^\frac{4}{4}<0$. We then have that $q^\frac{3}{2}+1$ divides $6q^\frac{4}{4}-6q^\frac{3}{4}-12q^\frac{2}{4}-6q^\frac{1}{4}+6 < 6q^\frac{4}{4}\leq 2q^\frac{5}{4}< q^\frac{3}{2}+1$, which is absurd since both quantities are positive. Case $14$ is therefore excluded.

\smallskip

The four last cases are of class $\mathcal{C}_8$, therefore they preserve a non-degenerate bilinear form. We cannot have $PGP^{-1}$ included in one of those groups because $R_{[1^2],[1^2]|\mathcal{A}_{B_4}}\not\simeq R_{[1^2],[1^2]|\mathcal{A}_{B_4}}^\star$ by Proposition \ref{isomorphisme}.

\smallskip

The only remaining case is now cases $11$ which is $SL_3(q).(q^\frac{1}{2}+1).2$. We know that $H=R_{([1^2],[1^2])}(\mathcal{A}_{B_3})\simeq SL_3(q)$ normally generates $G=R_{([1^2],[1^2])}(\mathcal{A}_{B_4})\subset P^{-1}SL_6(q^{\frac{1}{2}})P$ for a certain matrix $P$ in $GL_6(q)$. Assume $PGP^{-1}$ is a subgroup of  $M=SL_3(q).(q^\frac{1}{2}+1).2$. Since $SL_3(q)$ is perfect, $PHP^{-1}$ is perfect and the image of $PHP^{-1}$ in the quotient $\Z/2\Z$ of $M$ is trivial. The group $H$ is thus included in $SL_3(q).(q^\frac{1}{2}+1)$. Using the same argument, we have that $PHP^{-1}$ is included in the $SL_3(q)$ appearing in the expression of $M$, therefore $PHP^{-1}$ is equal to that $SL_3(q)$. For all $g\in PGP^{-1}$, we can apply the same reasoning to $gPHP^{-1}g^{-1}= SL_3(q)=PHP^{-1}$. It follows that $PGP^{-1}=SL_3(q)=PHP^{-1}$ because $H$ normally generates $G$, therefore $H=G$. This leads to a contradiction because $G$ is irreducible and $H$ is not.

This concludes the study of double-partition $([1^2],[1^2])$ in the field cases $5$ and $6$ and we have $R_{[1^2],[1^2]}(\mathcal{A}_{B_4})\simeq SL_6(q^{\frac{1}{2}})$ in those cases.

\bigskip

Assume now we are in case $4$, i.e., $\F_q=\F_p(\alpha,\beta)= \F_p(\alpha,\beta+\beta^{-1}) \neq \F_p(\alpha+\alpha^{-1},\beta)= \F_p(\alpha+\alpha^{-1},\beta+\beta^{-1})$. Then by Lemma \ref{Harinordoquy} and Proposition \ref{isomorphisme}, there exists a matrix $P$ such that $PGP^{-1} \subset SU_6(q^{\frac{1}{2}})$ and $H\simeq SL_3(q)$, writing again $G=R_{([1^2],[1^2])}(\mathcal{A}_{B_4})$ and $H = R_{([1^2],[1^2])}(\mathcal{A}_{B_3})$. The goal this time is to show that $PGP^{-1}\simeq SU_6(q^{\frac{1}{2}})$. 

We first consider the maximal subgroups of class $\mathcal{S}$ of $SU_6(q^\frac{1}{2})$ given in Table 8.27 of \cite{BHRC} and give their order or a quantity their order divides.

\begin{enumerate}
\item $2 \times 3^{.} \mathfrak{A}_6, 2160$,
\item $2\times 3^{.} \mathfrak{A}_6.2_3, 4320$,
\item $6^{.}\mathfrak{A}_6, 2160$,
\item $(q^\frac{1}{2}+1,6)\circ 2^{.}L_2(11), 12\times 660=7920$,
\item $6^{.}\mathfrak{A}_7, 15120$,
\item $6^{.}PSL_3(4), 6\times 20160=120960$,
\item $6^{.}PSL_3(4)^{.}2_1^{-}, 241920$,
\item $3^{.}M_{22}, 3*443520=1330560$,
\item $3_1^{.}U_4(3):2_2, 6\times 3265920=19595520$,
\item $6_1^{.}U_4(3), 19595520$,
\item $6_1^{.}U_4(3)^{.}2_2^{-}, 39191040$,
\item $(q^\frac{1}{2}+1,6)\circ SU_3(q^\frac{1}{2}), 6(q^\frac{3}{2}(q-1)(q^\frac{3}{2}+1)$.
\end{enumerate}
As before, we have $q\geq 9$. Therefore $\vert SL_3(q)\vert \geq 42456960$ and the last case is excluded since $q^3$ does not divide $6q^\frac{3}{2}$.

\medskip

Consider now the maximal subgroups of geometric type. We here omit the groups of class $\mathcal{C}_1$ because we know that $PGP^{-1}$ is irreducible. The remaining maximal subgroups obtained from Table 8.26. of \cite{BHRC} of $SU_6(q^{\frac{1}{2}})$ are the following

\begin{enumerate}
\item $(q^\frac{1}{2}+1)^5.\mathfrak{S}_6, 720(q^\frac{1}{2}+1)^5$,
\item $SU_2(q^\frac{1}{2})^3:(q^\frac{1}{2}+1)^2.\mathfrak{S}_3, 6q^\frac{3}{2}(q^\frac{1}{2}+1)^2(q-1)^3$,
\item $SU_3(q^\frac{1}{2})^2:(q^\frac{1}{2}+1).\mathfrak{S}_2, 2q^3(q^\frac{1}{2}+1)(q-1)^2(q^\frac{3}{2}+1)^2$,
\item $SL_3(q).(q^\frac{1}{2}-1).2, 2(q^\frac{1}{2}-1)\vert SL_3(\F_q)\vert$,
\item $SU_2(q^\frac{3}{2}).(q-q^\frac{1}{2}+1).3, 3q^\frac{3}{2}(q^3-1)(q-q^\frac{1}{2}+1)$,
\item $SU_2(q^\frac{1}{2})\times SU_3(q^\frac{1}{2}), q^2(q-1)^2(q^\frac{3}{2}+1)$,
\item $SU_6(q_0).[(\frac{q+1}{q_0+1},6)], q_0=q^\frac{1}{2r}$, $r$ odd prime, $6q_0^{15}(q_0^2-1)(q_0^3-1)(q_0^4-1)(q_0^5-1)(q_0^6-1)$,
\item $(q^\frac{1}{2}+1,3)\times SP_6(q^\frac{1}{2}), 3q^\frac{3}{2}(q-1)\vert SL_3(q)\vert$,
\item $(q^\frac{1}{2}+1,3)\times SO_6^{+}(q^\frac{1}{2}).2, 6q^3(q^2-1)(q^\frac{3}{2}-1)(q-1)$,
\item $(q^\frac{1}{2}+1,3)\times SO_6^{-}(q^\frac{1}{2}).2, 6q^3(q^2-1)(q^\frac{3}{2}+1)(q-1)$.
\end{enumerate}

In case $6$, we have that the order of the maximal subgroup is $q^2(q-1)^2(q^\frac{3}{2}+1)$. If $\vert SL_3(q)\vert$ divides this quantity, then $q^3-1$ divides $(q-1)(q^{\frac{3}{2}}+1)=q^{\frac{5}{2}}-q^{\frac{3}{2}}+q-1<q^{\frac{5}{2}}-1<q^3-1$. This is absurd, therefore this case is excluded.

\smallskip

In case $1$, we have that $q^3$ divides $720$ which is absurd since $q\geq 36$. In the same way, cases $2$ and $5$ are excluded because $q^3$ does not divide $6q^\frac{3}{2}$.

\smallskip

In case $3$, we have that $q^3-1$ divides $2(q^\frac{1}{2}+1)(q-1)^2(q^\frac{3}{2}+1)$. The  Euclidean remainder of the division of those two polynomials in $q^{\frac{1}{2}}$ is $-4q^\frac{5}{2}+8q^2-4q^\frac{3}{2}-4q+8q^\frac{1}{2}-4$n which is negative. We have
$4q^\frac{5}{2}-8q^\frac{4}{2}+4q^\frac{3}{2}+4q-8q^\frac{1}{2}+4\leq 4q^\frac{5}{2}+6q^\frac{3}{2}\leq 5q^\frac{5}{2}< q^3-1$ when $q^{\frac{1}{2}}\leq 3$.  This case is therefore also excluded.

\smallskip

In case $7$, we have that $q^3$ divides $6q_0^{15}=6q^\frac{15}{2r}\leq 6q^\frac{15}{6}<q^3$ since $q^\frac{1}{2}=q_0^r\geq 8$.

\smallskip

In cases $8$, $9$ and $10$, we have $PGP^{-1}$ included in a subgroup preserving a non-degenerate bilinear form. This would imply that $\epsilon \circ R_{[1^2],[1^2]|\mathcal{A}_{B_4}} \simeq R_{[1^2],[1^2]|\mathcal{A}_{B_4}}$ which is absurd by Proposition \ref{isomorphisme}.

\smallskip

The last remaining case is case $4$. In the same way as in case $11$ when $\F_p(\alpha,\beta)\neq \F_p(\alpha+\alpha^{-1},\beta)$, we would have that $H \neq G$ because $R_{([1^2],[1^2])}(S_2S_3^{-1})\in G\setminus H$. Since $G$ is normally generated by $H$ and $gHg^{-1}$ is perfect for all $g\in G$, we would have that $G \subset SL_3(q)$ and, therefore $G=H$, which is absurd.

\medskip

We have shown that $PGP^{-1}$ cannot be included in any maximal subgroup of $SU_6(q^{\frac{1}{2}})$. It follows that $PGP^{-1}\simeq SU_6(q^\frac{1}{2})$.

\bigskip

The only double-partition remaining for $n\leq 4$ now is $\lambda=([2,1],[1])$. It affords a representation of dimension 8 and satisfies $\lambda=(\lambda_1',\lambda_2')$.

We start by case $4$, i.e., $\F_q=\F_p(\alpha,\beta)=\F_p(\alpha,\beta+\beta^{-1})\neq \F_p(\alpha+\alpha^{-1},\beta)=\F_p(\alpha+\alpha^{-1},\beta+\beta^{-1})$ and, therefore $\F_{\tilde{q}}=\F_p(\alpha+\alpha^{-1})\neq \F_p(\alpha)$. We then have by Goursat's Lemma and the result for $n=3$ that $H=R_{([2,1],[1])}(\mathcal{A}_{B_3})\simeq SL_3(q)\times SU_2(\tilde{q}^\frac{1}{2}) \subset G=R_{([2,1],[1])}(\mathcal{A}_{B_4})$. By Proposition \ref{isomorphisme}, we know that there exists $P\in GL_8(q)$ such that for all $h\in \mathcal{H}_4$,   $PR_{([2,1],[1])}(h)P^{-1}= \epsilon(R_{([2,1],[1])})(h)$. By Lemma \ref{Ngwenya}, this implies that there exists $S\in GL_8(q)$ such that $S^{-1}R_{([2,1],[1])}(\mathcal{A}_{B_4})S\subset GL_8(q^\frac{1}{2})$ with $\gamma^{-1} P=\epsilon(S)S^{-1}$ and $\epsilon(P)P=\epsilon(\gamma)\gamma$. 

We can use the arguments used previously to see that our group is primitive, irreducible, tensor-indecomposable, preserves no symmetric, skew-symmetric or hermitian form over $\F_q^{\frac{1}{2}}$ and cannot be included in $GL_8(q')$ for $q'<q^\frac{1}{2}$. We then get that $G$ is included in no maximal subgroup of class $\mathcal{C}_1$, $\mathcal{C}_2$, $\mathcal{C}_4$, $\mathcal{C}_5$ and $\mathcal{C}_8$. It contains a transvection, therefore it cannot be included in a maximal subgroup of class $\mathcal{C}_3$. We list below the maximal subgroups remaining obtained from Tables 8.44 and 8.45 of \cite{BHRC}. We give the order of those groups or a quantity their order divides.

\begin{enumerate}
\item $((8,q^{\frac{1}{2}}+1)\circ 2^{1+6})^{.}(SP_6(2))$, $2477260800
$
\item $4_1^{.}PSL_3(4)$, $241920$
\item $((8,q^{\frac{1}{2}}+1)\circ 4_1)^{.}PSL_3(4)$, $1935360$
\item $(8\circ 4_1)^{.}PSL_3(4)$, $1935360$
\item $(8\circ 4_1)^{.}PSL_3(4).2_3$, $3870720$
\end{enumerate}

We have that $\alpha$ is of order greater than or equal to $7$ and $\epsilon(\alpha)=\alpha^{-1}$, where $\epsilon$ is the unique automorphism of order $2$ of $\F_q$. It follows that $\alpha^{q^{\frac{1}{2}}+1}=1$, therefore $q^{\frac{1}{2}}+1\geq 7$ and $q^{\frac{1}{2}}\geq 6$. We have $\F_{\tilde{q}}=\F_p(\alpha)\neq \F_p(\alpha+\alpha^{-1})$, therefore $\tilde{\epsilon}(\alpha)=\alpha^{-1}$, where $\tilde{\epsilon}$ is the unique automorphism of order $2$ of $\F_{\tilde{q}}$. This implies that $\tilde{q}^{\frac{1}{2}}\geq 6$. It follows that $\vert H\vert \geq 591963268176000$. This excludes all the maximal subgroups in the list. It follows that $S^{-1}GS$ is included in no maximal subgroup of $SL_8(q^{\frac{1}{2}})$, therefore $G\simeq SL_8(q^{\frac{1}{2}})$.

\medskip
 
We now consider cases $5$ and $6$, where our representation is now unitary by Proposition  \ref{isomorphisme}. In both cases, there exists a matrix $P$ such that $PGP^{-1}\subset SU_8(q^\frac{1}{2})$ with $G=R_{([2,1],[1])}(\mathcal{A}_{B_4})$ and we have $H=R_{([2,1],[1])}(\mathcal{A}_{B_3})\simeq SL_3(q)\times SU_2(\tilde{q}^{\frac{1}{2}})$ (resp $SL_2(\tilde{q}))$ in case $5$ (resp case $6$). This proves that $G$ contains either a natural $SL_2(\tilde{q})$ or a natural $SU_2(\tilde{q})$. We then have by Lemma \ref{tens1} that $G$ is tensor-indecomposable, therefore it is not included in any maximal subgroup of class $\mathcal{C}_4$ of $SU_8(q^{\frac{1}{2}})$. It contains a transvection, therefore it cannot be included in any group of class $\mathcal{C}_3$.  We also have that $G$ is a primitive irreducible group preserving no symmetric or skew-symmetric form over $\F_{q^\frac{1}{2}}$ si $G$ is included in no maximal subgroup of class $\mathcal{C}_1$ or $\mathcal{C}_2$ or $\mathcal{C}_5$ for $q_0=q^{\frac{1}{2}}$. Consider now the maximal subgroups of $SU_8(q^{\frac{1}{2}})$ which are not of class $\mathcal{C}_1$, $\mathcal{C}_2$, $\mathcal{C}_3$, $\mathcal{C}_4$ or $\mathcal{C}_5$ with $q_0=q^{\frac{1}{2}}$. They are given in Tables $8.46$ and $8.47$ of \cite{BHRC} and we list them below with their order or a quantity their order divides

\begin{enumerate}
\item $SU_8(q_0)$, $q_0^{28}\underset{i=2}{\overset{8}\prod}(q_0^i-(-1)^i)$ $q^{\frac{1}{2}}=q_0^r$, $r$ odd prime
\item $((8,q^{\frac{1}{2}}+1)\circ 2^{1+6})^{.}(SP_6(2))$, $2477260800
$
\item $((8,q^{\frac{1}{2}}+1)\circ 4_1)^{.}PSL_3(4)$, $1935360$
\item $(8\circ 4_1)^{.}PSL_3(4).2_3$, $3870720$
\end{enumerate}

We have that $\alpha$ of order greater than or equal to $7$ and $\epsilon(\alpha)=\alpha$, where $\alpha$ is the automorphism of order $2$ of $\F_q$. It follows that $\alpha^{q^{\frac{1}{2}}-1}=1$ and, therefore $q^{\frac{1}{2}}-1\geq 7$ and $q^{\frac{1}{2}}\geq 8$. In case $5$, we have $\F_{\tilde{q}}=\F_p(\alpha)\neq \F_p(\alpha+\alpha^{-1})$, therefore $\tilde{\epsilon}(\alpha)=\alpha^{-1}$, where $\tilde{\epsilon}$ is the automorphism of order $2$ of $\F_{\tilde{q}}$ and, therefore $\tilde{q}^{\frac{1}{2}}\geq 6$. In case $6$, we have that $\tilde{q}\geq 7$ because $\alpha$ is of order greater than or eaqual to $7$. This proves that in both cases we have $\vert H\vert \geq 6\times (6^2-1)\times 64^3\times (64^2-1)(64^3-1)=59095088588390400$. This excludes the last three cases.

\smallskip

Assume that $G$ is included in the first maximal subgroup. Let $r$ be the prime such that $q_0^r=q^{\frac{1}{2}}$. We have that $q^3$ divides $q_0^{28}=q^{\frac{28}{2r}}=q^{\frac{14}{r}}$. This implies that $3\leq \frac{14}{r}$ and, therefore $r\leq \frac{14}{3}<5$. Since $r$ is an odd prime, we have that $r=3$. We then have that $(q^3-1)$ divides $\underset{i=2}{\overset{8}\prod}(q^{\frac{i}{6}}-(-1)^i)$. The Euclidean remainder of those two polynomials in $q^{\frac{1}{6}}$ is then $2q^{\frac{17}{6}}-2q^{\frac{16}{6}}+2q^\frac{15}{6}-q^{\frac{14}{6}}+q^\frac{13}{6}-q^\frac{12}{6}-q^\frac{11}{6}+q^\frac{10}{6}-q^\frac{9}{6}-q^\frac{8}{6}+q^\frac{7}{6}-q^\frac{6}{6}-q^\frac{5}{6}+q^\frac{4}{6}-q^\frac{3}{6}+2q^\frac{2}{6}-2q^{\frac{1}{6}}+2$. The latter quantity is positive, therefore we have that $q^3-1\leq 2q^{\frac{17}{6}}-q^{\frac{16}{6}}<q^3-1$. This contradiction shows that the first cases is also excluded, therefore $G$ is included in no maximal subgroup of $SU_8(q^{\frac{1}{2}})$ and, therefore $G\simeq SU_8(q^{\frac{1}{2}})$. This concludes the proof of the lemma.
\end{proof}
 
We must now show that we can use Theorem \ref{CGFS}. The factorizations of $\lambda=(\lambda_1,\lambda_2)$ by $(\lambda_1',\lambda_2')$ and by $(\lambda_2,\lambda_1)$ change the arguments for the natural $SL_2(q)$ and twisted diagonal embeddings of $SL_3(q)$. Let $\lambda = (\lambda_1,\lambda_2)$ be a double-partition of $n\geq 5$.
 
We then have five different cases.
 
\begin{enumerate}

 \item $\lambda\neq \lambda', \lambda\neq (\lambda_2,\lambda_1)$ and $\lambda\neq (\lambda_1',\lambda_2')$. Let us show $R_{\lambda}(\mathcal{A}_{B_n})$ contains a natural $SL_3(q)$. It is sufficient to show that there exists $\mu \subset \lambda$ such that $\mu'\not\subset \lambda, (\mu_2,\mu_1) \not\subset \lambda$ and $(\mu_1',\mu_2') \not\subset \lambda$.
 
 We write $\lambda_1$ partition of $n_1$ and $\lambda_2$ partition of $n_2$ with $n=n_1+n_2\geq 5$. We only consider double-partitions with no empty component. This implies that $n_1$ and $n_2$ are greater than or equal to $1$. Since the roles of $\lambda_1$ and $\lambda_2$ are symmetrical for this, we can assume without loss of generality $n_1\geq n_2$.
 \begin{enumerate}
 \item $n_2=1$, we then have that $\lambda_2=\lambda_2'$, therefore $\lambda_1\neq \lambda_1'$. There exists $\mu_1\subset \lambda_1$ such that $\mu_1'\not\subset \lambda_1$. We then have that $\mu=(\mu_1,\lambda_2)\subset \lambda$, but $\mu'\not\subset \lambda$ and $(\lambda_2,\mu_1) \not\subset \lambda$ because $n_1-1\geq 4>1$ and $(\mu_1',\lambda_2') \not\subset \lambda$, because $\mu_1'\not\subset \lambda_1$.
 
 \item $n_1 > n_2= 2$ and $\lambda_1\neq \lambda_1'$. We set $\mu=(\lambda_1,[1])$, we have $\mu'$ and $([1],\lambda_1)\not\subset \lambda$ because $n_1>n_2$ and $(\lambda_1',[1])\not\subset\lambda$ because $\lambda_1'\neq \lambda_1$.
 
 \item $n_1>n_2=2$ and $\lambda_1=\lambda_1'$. If for all $\mu_1\subset \lambda_1$, $\mu_1\subset \lambda_2$ or $\mu_1'\subset \lambda_2$, then $n_1=3$ and $\lambda_1=[2,1]$. This implies that either $([2],[1^2])\subset \lambda$ or $([1^2],[2])\subset \lambda$. By Proposition \ref{patate}, $R_{\lambda}(\mathcal{A}_{B_n})$ contains a natural $SL_3(q)$.
 
 \item $n_1>n_2\geq 3$ and $\lambda_2\neq \lambda_2'$. There exists $\mu_2\subset \lambda_2$ such that $\mu_2'\not\subset\lambda_2$. We then set $\mu=(\lambda_1,\mu_2)$. We have that $(\mu_2,\lambda_1) \not\subset \lambda, (\mu_2',\lambda_1') \not\subset \lambda$ because $n_1>n_2$ and $(\lambda_1',\mu_2') \not\subset \lambda$ because $\mu_2'\not\subset \lambda_2$.
 
 \item $n_1>n_2\geq 3$ and $\lambda_2=\lambda_2'$, therefore $\lambda_1\neq \lambda_1'$. We know that there exists $\mu_1 \subset \lambda_1$ such that $\mu_1'\not\subset \lambda_1$. If $(\lambda_2,\mu_1)\subset \lambda$ or $(\lambda_2',\mu_1') \subset \lambda$, then $\mu_1=\lambda_2$ or $\mu_1'=\lambda_2$. We have that $\lambda_2=\lambda_2'$, therefore this contradicts $\mu_1'\not\subset \lambda_1$. This shows that $\mu_1\neq \mu_1'$. We have that $(\mu_1',\lambda_2') \not\subset \lambda$ because $\mu_1'\not\subset \lambda_1$.
 
 \item  $n_1=n_2\geq 3$. We then have that $\lambda_1\neq \lambda_1'$ or $\lambda_2\neq \lambda_2'$. If $\lambda_1\neq \lambda_1'$, we pick $\mu_1\subset \lambda_1$ such that $\mu_1'\not\subset \lambda_1$ and set $\mu=(\mu_1,\lambda_2)$, by the assumption on $\mu_1$,  $(\mu_1',\lambda_2)\not\subset \lambda, (\lambda_2,\mu_1)\not\subset \lambda$ and $(\lambda_2',\mu_1')\not\subset \lambda$ because $\lambda_2\neq \lambda_1$ and $\lambda_2'\neq \lambda_1$. If $\lambda_2\neq \lambda_2'$, we pick $\mu_2\subset \lambda_2$ such that $\mu_2'\not\subset \lambda_2$ and $\mu=(\lambda_1,\mu_2)$ verifies the required property.
 
 \end{enumerate} 
 
 \bigskip
 
 \item  $\lambda=(\lambda_1',\lambda_2'), \lambda\neq (\lambda_2,\lambda_1)$ and $\lambda\neq \lambda'$. We then have that $\mu_1'\subset \lambda$ for all $\mu_1\subset \lambda_1$ and that for all $\mu_2\subset \lambda_2$, $\mu_2'\subset \lambda_2$. We also have that $n_1+n_2\geq 5$.

\begin{enumerate}

\item $n_1\geq n_2=1$. Let $\mu_1\subset \lambda_1$, we set $\mu= (\mu_1, \lambda_2)$. We have that  $(\lambda_2,\mu_1)\not\subset \lambda$ and $\mu'\not\subset \lambda$ because $n_1-1\geq 3>1$ and $(\mu_1',\lambda_2') \subset \lambda$.
 
 \item $n_1>n_2\geq 2$. We pick $\mu_2\subset \lambda_2$ and set $\mu=(\lambda_1,\mu_2)$. We have that $(\mu_2,\lambda_1)\not\subset \lambda$ and $\mu'\not\subset \lambda$ because $n_1>n_2$ and $(\lambda_1',\mu_2') \subset \lambda$.
 
 \item $n_1=n_2\geq 2$. We pick $\mu_1\subset \lambda_1$ and set $\mu=(\mu_1,\lambda_2)$. We have that $(\lambda_2,\mu_1)\not\subset \lambda$ because  $\lambda_2\neq \lambda_1$ and $\mu'\not\subset \lambda$ because $\lambda_2'\neq \lambda_1$ and $(\mu_1',\mu_2')\subset \lambda$.
 
  In case $4$ for the fields, i.e., $\F_q=\F_p(\alpha,\beta)=\F_p(\alpha,\beta+\beta^{-1})\neq \F_p(\alpha+\alpha^{-1},\beta)$,  if $\mu\neq (\mu_1', \mu_2')$ then $R_{\lambda}(\mathcal{A}_{B_n})$ contains up to conjugation $\{\diag(
 M,\epsilon(M),I_{n_\lambda-6}), M\in SL_3(q)\}$, and a natural $SL_3(q^{\frac{1}{2}})$ if $\mu=(\mu_1',\mu_2')$ (it is possible for this to be the case for all $\mu\subset \lambda$ if we have square partitions).
 
 In cases $5$ and $6$ for the fields, i.e., $\F_q=\F_p(\alpha,\beta)=\F_p(\alpha+\alpha^{-1},\beta)\neq \F_p(\alpha,\beta+\beta^{-1})$,  if $\mu\neq (\mu_1', \mu_2')$ then $R_{\lambda}(\mathcal{A}_{B_n})$ contains up to conjugation \\$\{\diag(
 M,{}^t\!\epsilon(M^{-1}),I_{n_\lambda-6}), M\in SL_3(q)\}$, and a natural $SU_3(q^{\frac{1}{2}})$ if $\mu=(\mu_1',\mu_2')$.
 
 \end{enumerate}

 \item $\lambda=(\lambda_2,\lambda_1)\neq \lambda'$. We then have $n_1=n_2\geq 3$ and $\lambda_1=\lambda_2\neq \lambda_1'$ because $\lambda\neq \lambda'$. We can then pick $\mu_1\subset \lambda_1$ such that $\mu_1'\not\subset \lambda_1'$ and $\mu=(\mu_1,\lambda_2)$. We have $(\mu_1',\lambda_2')\not\subset \lambda$ and $\mu'\not\subset \lambda$ because $\lambda_2'\neq \lambda_1=\lambda_2$ and $(\lambda_2,\mu_1) \subset \lambda$ but $\mu\neq (\lambda_2,\mu_1)$ because $n_1-1<n_1=n_2$.
 
 In case $4$ for the fields, $R_{\lambda}(\mathcal{A}_{B_n})$ contains up to conjugation
  \\$\{\diag(M,{}^t\!\epsilon(M^{-1}),I_{n_\lambda-6}),M\in SL_3(q)\}$.
                          
  In cases $5$ and $6$ for the fields, $R_{\lambda}(\mathcal{A}_{B_n})$ contains up to conjugation \\$\{\diag(M,\epsilon(M),I_{n_\lambda-6}), M\in SL_3(q)\}$.

\item $\lambda=\lambda'\neq (\lambda_1',\lambda_2')=(\lambda_2,\lambda_1)$, we have $n_1=n_2\geq 3$ and there exists $\mu_1\subset \lambda_1$ such that $\mu_1\not\subset \lambda_2$ because $\lambda_1\neq \lambda_2$. We have $\mu'\subset \lambda$ because $\mu_1'\subset\lambda_1'=\lambda_2$, $(\lambda_2,\mu_1) \not\subset \lambda$ since $\lambda_2\neq \lambda_1$ and $(\mu_1',\lambda_2') \not\subset \lambda$ because $\lambda_2'\neq \lambda_2=\lambda_1'$. We have $\mu\neq \mu'$ because $\lambda_2'\neq \mu_1$. $R_\lambda(\mathcal{A}_{B_n})$ contains up to conjugation 
 $\{\diag(M,{}^t\!(M^{-1}),I_{n_\lambda-6}), M\in SL_3(q)\}$.

\item  $\lambda =\lambda'=(\lambda_2,\lambda_1)=(\lambda_1',\lambda_2')$. We then have $n_1=n_2\geq 3$. If $\lambda_1$ and $\lambda_2$ are square partitions, then for all $\mu \subset \lambda$, we have that $\mu=(\mu_1',\mu_2')\neq \mu'= (\mu_2,\mu_1)$, because $n_1=n_2> n_1-1=n_2-1$.
  
  In case $4$ for the fields, $R_\lambda(\mathcal{A}_{B_n})$ contains up to conjugation
  \\$\{\diag(M,{}^t\!(M^{-1}), I_{n_\lambda-6} ), M\in SU_3(q^{\frac{1}{2}})\}$.
    
    In cases $5$ and $6$ for the fields, $R_\lambda(\mathcal{A}_{B_n})$ contains up to conjugation   \\$\{\diag(M,{}^t\!(M^{-1}),I_{n_\lambda-6}), M\in SL_3(q^{\frac{1}{2}})\}$.
    
    If $\lambda_1$ or $\lambda_2$ is a square partition, then there exists $\mu\subset \lambda$ such that $\mu\neq \mu'$, $\mu\neq (\mu_2,\mu_1)$ and $\mu\neq (\mu_1',\mu_2')$. This implies that $R_\lambda(\mathcal{A}_{B_n})$ contains up to conjugation \\$\{\diag(M,{}^t\!(M^{-1}),\epsilon(M),\epsilon({}^t\!(M^{-1})),I_{n_\lambda-12}), M\in SL_3(q)\}$.
      
\end{enumerate}

We now use the notations of Theorem \ref{CGFS}. In all of the above cases except for the last one, there exists $g$ in $R_{\lambda}(\mathcal{A}_{B_n})$ such that $[g,V]\leq 2$. This implies that $v_G(V)\leq 2$, therefore $v_G(V) \leq \max(2,\frac{\sqrt{d}}{2})$. In the last case, we have in the same way an element $g$ such that $[g,V]=4$. We also have in that case that $\lambda=\lambda'=(\lambda_1,\lambda_2)$ and $n\geq 6$,.This implies that $\lambda$ contains $([2,1],[2,1])$, which is of dimension $\binom {6} {3}\times 2\times 2=80$. It follows that $d\geq 80$ and $\frac{\sqrt{d}}{2}\geq \frac{\sqrt{80}}{2}> 4$. This shows that we still have $v_G(V) \leq \max(2,\frac{\sqrt{d}}{2})$.

\bigskip

It remains to check that all the assumptions of the theorem are again verified and the classical group we get is the one we want.

The first step is to take care separately of double-partitions $\lambda$ such that $n_\lambda\leq 10$. If $n_\lambda=10$, then by the conditions of Theorem \ref{CGFS}, we can assume $p\neq 2$. The second step is to verify that the remaining double-partitions are tensor-indecomposable. The third step is to verify that they are imprimitive in the monomial case. The fourth step is to verify that they are imprimitive in the non-monomial case. The fifth step is to check that we are not in case 2. of Theorem \ref{CGFS}. The sixth and last step is to verify that we have the desired classical groups in each of the above cases.

\bigskip

\textbf{First step.} For $n=5$, it is enough to consider $([2,2],[1])$, $([2,1,1],[1])$, $([2,1],[2])$ and $([1^3],[1^2])$, for which the respective $n_\lambda$ is $10$, $15$, $20$ and $10$.\\
 We must show that $R_{([1^3],[1^2])}(\mathcal{A}_{B_5})=SL_{10}(q)$ and $R_{([2,2],[1])}(\mathcal{A}_{B_5})\simeq SL_{10}(q^\frac{1}{2})$ in case $4$ for the fields and $R_{([2,2],[1])}(\mathcal{A}_{B_5})\simeq SU_{10}(q^\frac{1}{2})$ in cases $5$ and $6$ for the fields. The other double-partitions are of dimensions greater than $10$. We know that $G=R_{([1^3],[1^2])}(\mathcal{A}_{B_5})$ contains $R_{([1^3],[1^2])}(\mathcal{A}_{B_4}) \simeq SL_4(q) \times SL_6(q^\frac{1}{2})$ and it is normally generated by this group, which is generated by transvections. Since $p\neq 2$, Theorem \ref{transvections} implies that $G$ is conjugate in $GL_{10}(q)$ to $SL_{10}(q'), SP_{10}(q')$ or $SU_{10}(q'^{\frac{1}{2}})$ for some $q'$ dividing $q$. Lemma \ref{field} implies that $q'=q$. The groups $SP_{10}(q)$ and $SU_{10}(q^\frac{1}{2})$ are excluded by Proposition \ref{isomorphisme}, because $R_{([1^3],[1^2])}$ is not isomorphic to its dual representation or its dual representation composed with the automorphism of order $2$ of $\F_q$. This shows that $G=SL_{10}(q)$. 
In case $4$, we know that $G=R_{([2,2],[1])}(\mathcal{A}_{B_5})$ is conjugate to a subgroup of $SL_{10}(q^{\frac{1}{2}})$ by Proposition \ref{isomorphisme} and Lemma \ref{Ngwenya} and that $G$ contains $R_{([2,2],[1])}(\mathcal{A}_{B_4}) \simeq  SL_8(q^\frac{1}{2})\times SL_2(\tilde{q})$. It follows that it contains a natural $SL_8(q^\frac{1}{2})$ and we can apply Theorem \ref{LBJ} to get that $G\simeq SL_{10}(q^\frac{1}{2})$.  
In cases $5$ and $6$, $G=R_{([2,2],[1])}(\mathcal{A}_{B_5})$ is conjugate to a subgroup of $SU_{10}(q^\frac{1}{2})$ by Proposition \ref{isomorphisme} and Lemma \ref{Harinordoquy} and contains $R_{([2,2],[1])}(\mathcal{A}_{B_4})$, which contains a natural $SU_{8}(q^\frac{1}{2})$ in both cases. By Theorem $1.4$ of \cite{BM}, we have indeed $G\simeq SU_{10}(q^\frac{1}{2})$.

\bigskip

\textbf{Second step.} We now show that those representations are tensor-indecomposable. Since $([2,1,1],[1])$ contains a natural $SL_3(q)$, doubles-partitions with at most two rows or at most two columns are tensor-indecomposable by Lemmas \ref{tens1} and \ref{tens2}. By the enumeration of the different cases, those lemmas cover all double-partitions of $n$ except if $\lambda=(\lambda_1,\lambda_2)=\lambda'=(\lambda_2,\lambda_1)=(\lambda_1',\lambda_2')$ and neither $\lambda_1$ nor $\lambda_2$ contains a sub-partition $\mu$ such that $\mu=\mu'$. In such a case,  $\lambda$ contains $([2,1],[2,1])$ which is of dimension $80$ and we can use the following lemma.

\begin{lemme}
If $d\geq 80$ and $G\subset GL_d(q)$ contains an element of order coprime to $p$ and conjugate in $GL_d(q)$ to the diagonal matrix $\diag(\xi,\xi,\xi,\xi,\xi^{-1},\xi^{-1},\xi^{-1},\xi^{-1},1,...,1)$ with $\xi^2\neq 1$, then $G$ is tensor-indecomposable, except possibly if $G\subset G_1\otimes G_2$ with $G_1\subset GL_a(q)$, $a\in \{2,4\}$.
\end{lemme}

\begin{proof}

Let $g=P\diag(\xi I_4,\xi^{-1} I_4,I_{d-8})P^{-1}$. Assume that $g=g_1\otimes g_2$ with $g_1\in GL_a(\F_q)$, $g_2\in GL_b(\F_q)$ with $3\leq a\leq b$ and $ab=d$. We have that $b\geq \sqrt{d}$, therefore $b\geq 9$ because $d\geq 80$. We write $\lambda_1,...,\lambda_a$ the eigenvalues of $g_1$ and $\mu_1,...,\mu_b$ the eigenvalues of  $g_2$. We then have that $\forall i \in [\![ 1,a]\!]$, $\forall j\in [\![1,b]\!]$, $\lambda_i\mu_j\in \{1,\xi,\xi^{-1}\}$. The numbers $\xi$ and $\xi^{-1}$ only appear $4$ times each. This implies the number of couples $(\lambda_1\mu_i,\lambda_2\mu_i)\in \{(1,\xi),(\xi,1),(\xi,\xi^{-1})\}$ is less than or equal to $4$ as is the number of couples $(\lambda_1\mu_i,\lambda_2\mu_i)\in \{(1,\xi^{-1}),(\xi^{-1},1),(\xi^{-1},\xi)\}$. For any $i\in [\![1,a]\!]$, the inequality $\lambda_1\mu_i\neq \lambda_2\mu_i$ implies that $(\lambda_1\mu_i,\lambda_2\mu_i)\in \{(1,\xi),(\xi,1),(1,\xi^{-1}),(\xi^{-1},1),(\xi,\xi^{-1}),(\xi^{-1},\xi)\}$. It follows that there are at most $8$ couples $(\lambda_1\mu_i, \lambda_2\mu_i)$ such that $\lambda_1\mu_i\neq \lambda_2\mu_i$. Since $b\geq 9$, there exists $i\in [\![1,a]\!]$ such that $\lambda_1\mu_i=\lambda_2\mu_i$. It follows that $\lambda_1=\lambda_2$. In the same way, we have that $\lambda_1=\lambda_j$ for all $j \in [\![1,a]\!]$.  Up to reordering, we can assume $\lambda_1\mu_1=\xi$. We then have $\lambda_2\mu_1=\lambda_3\mu_1=\xi$. Since there are exactly $4$ ways $\xi$ appears as a $\lambda_i\mu_j$, we have that $a=4$.

 By the assumptions on $\lambda$, $H=R_{\lambda}(\mathcal{A}_{B_{n-1}})$ is a direct product of groups isomorphic to  some $SL_m(q)$ with $m\geq n_{([2,1],[2])}=20$. If $G=R_{\lambda}(\mathcal{A}_{B_n})$ is not tensor-indecomposable, then $G\subset SL_2(q) \otimes SL_{\frac{d}{2}}(q)$ or $G\subset SL_4(q) \otimes SL_{\frac{d}{4}}(q)$. We then have a morphism from $G$ into $SL_2(q)$ or $SL_4(q)$. If we consider the restriction of this morphism to $H$, its kernel is a normal subgroup of $H$. The only non-abelian decomposition factors of $H$ are $PSL_m(q)$ with $m\geq 20$. If the image is non-abelian, then there exists a subgroup of $SL_2(q)$ or a subgroup of $SL_4(q)$ isomorphic to some $PSL_m(q)$. This leads to a contradiction because $m\geq 20$. It follows that the image is abelian and since $H$ is perfect, the kernel is equal to $H$.  Since $H$ normally generates $G$, the morphism is trivial on $G$ which contradicts the irreducibility of $G$. 
 \end{proof}
 
 \bigskip
 
 \textbf{Third step.} In the monomial case, the only additional case to consider is the same one as in the second step. Looking at the corresponding proof in \cite[page 14]{BMM}, we get that $(p-1)r\leq 4$ with $q=p^r$. We know that $q$ is a square, $n\geq 6$, $\alpha$ is of order greater than $n$ and $\epsilon(\alpha)\in \{\alpha,\alpha^{-1}\}$. Therefore $\alpha^{q^\frac{1}{2}-1}=1$ or $\alpha^{q^\frac{1}{2}+1}=1$. In both cases $q^\frac{1}{2}+1> 6$, and, therefore $q^\frac{1}{2}\geq 6$ and $q\geq 36$. The condition $(p-1)r\leq 4$ implies that $q\leq \max(5^1,4^1,3^2,2^4)= 16$, therefore we have a contradiction.
 
 \bigskip
 
 \textbf{Fourth step.} We know that there exists a matrix $t$ of order $p$ such as the one in \cite[page 14]{BMM} or with Jordan form $\diag(I_2+E_{1,2},I_2+E_{1,2},I_2+E_{1,2},I_2+E_{1,2},I_{n_\lambda-8})$.
 
 If $p\neq 2$, we can use the same arguments as in page $15$ of \cite{BMM} because we still have $(t-1)^2=0$.
 
 Assume now that $p=2$. Assume that $G\subset H \wr \mathfrak{S}_m = (H_1\times H_2 \times \dots \times H_m) \rtimes \mathfrak{S}_m$ with $H_1, \dots H_m$ the $m$-copies of $GL_{N/m}(q)$ permuted by $\mathfrak{S}_m$, that $V=U_1\oplus U_2 \oplus \dots \oplus U_m$ is the direct sum corresponding to the wreath product and that $t\notin H_1\times \dots \times H_m$.  Assume $t\notin H_1\times \dots \times H_m$. Up to reordering, we can assume $tU_1=U_2$.  If $\dim(U_i)\geq 5$ then we can consider linearly independent vectors $v_1,v_2,v_3,v_4,v_5$ in $U_1$ and by completing the family of vectors $(v_1,tv_1,v_2,tv_2,v_3,tv_3,v_4,tv_4,v_5,tv_5)$ which are linearly independent because $tU_1=U_2\neq U_1$, we get a basis upon which $t$ acts as a matrix of the form $M_2\oplus M_2\oplus M_2 \oplus M_2 \oplus M_2 \oplus X$ for a certain $X$ with $M_2=\begin{pmatrix}
 0 & 1\\
 1& 0
 \end{pmatrix}$. This implies that the rank of $t-1$ is greater than or equal to $5$, which is a contradiction.
 
 We can thus assume that $\dim(U_i) \leq 4$. Note that $G=R_\lambda(\mathcal{A}_{B_n})$ and $\mathcal{A}_{B_n}$ is perfect for $n\geq 5$ \cite{MR}, therefore $G$ is perfect. If $G\subset (H_1\times H_2\times \dots \times H_m) \rtimes \mathfrak{S}_m$, we get $G\subset (H_1\times H_2 \times \dots \times H_m) \rtimes \mathfrak{A}_m$ because $[\mathfrak{S}_m,\mathfrak{S}_m]\subset \mathfrak{A}_m$. 
 
  If $t$ is a transvection then by the same reasoning as above on the dimensions of $U_i$, we are in the monomial case which was done in the third step. 
 
  If $t$ is of rank $2$, then either we are in the monomial case or $\dim(U_i)=2$. The monomial case is done, therefore it is sufficient to prove that $\dim(U_i)=2$ leads to a contradiction. We take $t_1$ and $t_2$ two such elements of rank $2$. Assume $\dim(U_i)=2$, since we have $t(U_1)=U_2$ and $t_1(U_2)=t_1^2(U_1)=U_1$. If $(u_a,u_b)$ are linearly independent then $(t_1u_a-u_a,t_1u_b-u_b)$ is a basis of $\mathrm{Im}(t_1-1)$,  which is of dimension $2$ and included in $U_1\oplus U_2$ for all $i\notin\{1,2\}$, $t_i(U_i)=U_i$. It follows that the projection of $t_1$ upon $\mathfrak{S}_m$ from the semi-direct product is a transposition. This is a contradiction because the projection of $G$ upon $\mathfrak{S}_m$ is included in $\mathfrak{A}_m$.
 
 If $t$ is of rank $4$ and $R_\lambda(\mathcal{A}_ {B_{n-1}})$ does not contain in an obvious way any transvections or elements $t$ of rank $2$, then $G$ contains up to conjugation \\
 $\{\diag(M,{}^t\!(M^{-1}),\epsilon(M),{}^t\!\epsilon(M^{-1}),I_{n_\lambda-8}), M\in SL_2(q)\}$.
We consider two elements $t_1$ and $t_2$ of rank $4$. If $\dim(U_1)=4$, then if $(u_1,u_2,u_3,u_4)$ is a basis of $U_1$, $(u_1-t_1u_1,u_2-t_1u_2,u_3-t_1u_3,u_4-t_1u_4)$ is a basis of $\mathrm{Im}(t_1-1)$, which is of dimension $4$. It follows that the projection of $t_1$ upon $S_m$ is a transposition, which is absurd.

 If $\dim(U_1)=3$, then if $(u_1,u_2,u_3)$ is a basis of $U_1$, we have that $\op{Vect}\{t_1u_1-u_1,t_1u_2-u_2,t_1u_3-u_3\} \subset \mathrm{Im}(t_1-1)$. If there exists $i\notin \{1,2\}$ such that $t_i(U_i)\neq U_i$ then in the same way as before, there would exist a subspace of dimension $6$ of $\mathrm{Im}(t_1-1)$, which is of dimension $4$. This shows that the projection of $t_1$ upon $\mathfrak{S}_m$ is a transposition, which is absurd.
 
 If $\dim(U_i)=2$, then we can take $4$ distinct non-zero elements $a_1,a_2,a_3,a_4$ of $\F_q$. This is possible because $q^{\frac{1}{2}}\geq 6$. We know that $G$ contains up to conjugation the elements $t_j$ for $j\in \{1,2,3,4\}$ with $t_j=\diag(I_2+a_j E_{1,2},I_2+a_j E_{1,2},I_2+\epsilon(a_j) E_{1,2}, I_2+\epsilon(a_j)E_{1,2},I_{n_\lambda-8})$. We have that
$\mathrm{Im}(t_j-1)$ is independent of $j$. We also have that $t_1(U_1)=U_2$ and $t_1(U_2)=t_1^2(U_1)=U_1$. Since $\mathrm{Im}(t_1-1)\cap U_1\oplus U_2$ is then of dimension $2$ and the projection of $t_1$ upon $\mathfrak{S}_m$ is not a transposition, there exists $i\notin\{1,2\}$ such that $t_1(U_i)\neq U_i$. Up to reordering, we can assume $t_1(U_3)=U_4$ and $t_1(U_4)=t_1^2(U_3)=U_3$. This shows that for all $j\in \{1,2,3,4\}, \mathrm{Im}(t_j-1)=\mathrm{Im}(t_1-1)\subset U_1\oplus U_2 \oplus U_3 \oplus U_4$. Since each $t_j$ is of order $2$, it follows writing $\pi$ the projection of $G$ upon $\mathfrak{S}_m$ that we have $\{\pi(t_1),\pi(t_2),\pi(t_3),\pi(t_4)\}\subset \{I_d,(12)(34),(13)(24),(14)(23)\}$.

Let us show that $\pi(t_j)=I_d$ for all $j$. They are all conjugate in $G$. Since $H_1\times H_2\times \dots \times H_m$ is a normal subgroup of $(H_1\times H_2\times \dots \times H_m)\rtimes \mathfrak{S}_m$, it is sufficient to show it for one of them. Assume it is false for all of them. We then have $\{\pi(t_1),\pi(t_2),\pi(t_3),\pi(t_4)\}\subset \{(12)(34),(13)(24),(14)(23)\}$. Therefore, there exists a pair $(i,j), i\neq j$ such that $\pi(t_i)=\pi(t_j)$ and, therefore $\pi(t_it_j)=I_d$. But the matrix of $t_it_j$ in the basis we chose is $\diag(I_2+(a_i+a_j) E_{1,2},I_2+(a_i+a_j) E_{1,2},I_2+\epsilon(a_i+a_j) E_{1,2},I_2+\epsilon(a_i+a_j)E_{1,2},I_{n_\lambda-8})$.
We have $a_i+a_j\neq 0$ because $p=2$ and the elements $a_l$ are pairwise distinct. It follows that $t_it_j$ is conjugate to each $t_l$, therefore we have a contradiction. This shows that for all $j\in \{1,4\}, \pi(t_j)=I_d$. It follows that $t_j\in H_1\times H_2\times \dots \times H_m$, which is normal in $(H_1\times H_2\times \dots \times H_m)\rtimes \mathfrak{S}_m$. Since $G$ is normally generated by $R_\lambda(\mathcal{A}_{B_{n-1}})$, which is normally generated by elements of the form $t_j$, we have that $G\subset H_1\times H_2\times \dots \times H_m$. This contradicts the irreducibility of $G$. This is absurd and it follows that $G$ is a primitive group.
 
 \bigskip
 
\textbf{Fifth step.} If $G$ contains a natural $SL_2(q^\frac{1}{2})$ or a natural $SU_2(q^\frac{1}{2})$ then we can apply the same arguments as in \cite[page 13]{BMM}. If $G$ contains a twisted diagonal embedding or a twisted diagonal embedding composed with the automorphism of order $2$ of $\F_q$ of $SL_3(q)$, then we can apply the arguments of \cite[page 14]{BMM}. If we are not in any of the above cases, then $\lambda=\lambda'=(\lambda_2,\lambda_1)$, therefore $n\geq 6$ and we are in one of the following cases.
 \begin{enumerate}
 \item $R_\lambda(\mathcal{A}_{B_n})$ contains up to conjugation   $\{\diag(M,{}^t\!(M^{-1}),I_{n_\lambda-6} ), M\in SU_3(q^{\frac{1}{2}})\}$.
 \item $R_\lambda(\mathcal{A}_{B_n})$ contains up to conjugation   $\{\diag(M,{}^t\!(M^{-1}),I_{n_\lambda-6}), M\in SL_3(q^{\frac{1}{2}})\}$.
 \item $R_\lambda(\mathcal{A}_{B_n})$ contains up to conjugation \\$\{\diag(M,{}^t\!(M^{-1}),\epsilon(M),{}^t\!\epsilon(M^{-1}),I_{n_\lambda-12}), M\in SL_3(q)\}$.
 \end{enumerate}
 
In the first two cases, we have an element $g$ conjugate to $\diag(\xi,\xi,\xi^{-1},\xi^{-1},1,\dots,1)$ with $\xi$ of order $q^\frac{1}{2}-1$ but the order of $\alpha$ is less than or equal to $q^{\frac{1}{2}}+1$ in both cases. If $g$ is an element of $\mathfrak{S}_{n_\lambda}$ such that $[g,V]= 4$, then we have ($g=\sigma_1\sigma_2\sigma_3\sigma_4$ is the product of 4 disjoint transpositions and $g$ is of order $2$) or ($g$ is the product of 2 disjoint 3-cycles and $g$ is of order $3$) or ($g$ is a $5$-cycle and $g$ is of order $5$) or ($g$ is the disjoint product of  $2$ transpositions and a $3$-cycle and $g$ is of order $6$). Since $n_\lambda \geq 6$ and the order of $\alpha$ is greater than $n$, $q^\frac{1}{2}+1>7$, therefore $q^\frac{1}{2}-1>5$ which contradicts all the cases except for the last one. In the last case, we have that $n_\lambda\geq 7$ by the decomposition of $g$. Since $\lambda=\lambda'$, $n_\lambda$ is even and $q^\frac{1}{2}-1= q^\frac{1}{2}+1-2>n_\lambda-2> 6$, which contradicts the last case.

In the third case, we have an element $g$ conjugate to \\$\diag(\xi,\xi,\xi,\xi,\xi^{-1},\xi^{-1},\xi^{-1},\xi^{-1},1,\dots,1)$ which is of order $o(g)=q-1$. However $q^\frac{1}{2}+1>7$, therefore $q> 36$. Since $q$ is an even power of a prime number, it follows that $q>49$ and $q-1\geq 49$. We have $[g,V]= 8$. By considering the decomposition into disjoint cycles of $g$ and using the fact that the rank of $\sigma-1$ of a cycle $\sigma$ is equal to the length of the cycle minus 1, we get $o(g)\in \{\lcm(\{n_i+1\}_{i\in I}), \underset{i\in I}\sum n_i=8, n_i\in \N^\star\}$. It follows that $o(g) \leq 30 <49\leq q-1=o(g)$ which is a contradiction. 

\bigskip

\textbf{Sixth step.} We have shown that $G=R_\lambda(\mathcal{A}_{B_n})$ is a classical group in a natural representation. The last step is to show that we have the following theorem.
\begin{theo}
If $n\geq 5$, then for all double-partition $\lambda\Vdash n$ in our decomposition, $R_\lambda(\mathcal{A}_{B_n})=G(\lambda)$, where $G(\lambda)$ is given by the following list.
\begin{enumerate}
\item When $\F_q=\F_p(\alpha,\beta)=\F_p(\alpha,\beta+\beta^{-1})\neq\F_p(\alpha+\alpha^{-1},\beta)$, and $\F_{\tilde{q}}=\F_p(\alpha)\neq \F_p(\alpha+\alpha^{-1})$
\begin{enumerate}
\item $SU_{n-1}(\tilde{q}^\frac{1}{2})$ if $\lambda =([n-1,1],\emptyset)$.
\item $SU_{n_\lambda}(\tilde{q}^\frac{1}{2})$ if $\lambda = (\lambda_1,\emptyset), \lambda_1< \lambda_1'$.
\item $SP_{n_\lambda}(\tilde{q}^\frac{1}{2})$ if $\lambda = (\lambda_1,\emptyset), \lambda_1=\lambda_1'$ and ( $p=2$ or ($p\geq 3$ and $\nu(\lambda_1) = -1$)).
\item $\Omega_N^{+}(\tilde{q}^\frac{1}{2})$ if $\lambda=(\lambda_1,\emptyset), \lambda_1=\lambda_1'$, $p\geq 3$ and $\nu(\lambda_1)=1$.
\item $SL_n(q)$ if $\lambda=([1],[n-1])$.
\item  $SL_{n_\lambda}(q)$ if $\lambda\neq \lambda',\lambda\neq(\lambda_1',\lambda_2')$ and $\lambda\neq (\lambda_2,\lambda_1)$.
\item $SU_{n_\lambda}(q^{\frac{1}{2}})$, if $\lambda=(\lambda_2,\lambda_1) \neq \lambda'$.
\item $SL_{n_\lambda}(q^{\frac{1}{2}})$,if $\lambda=(\lambda_1',\lambda_2') \neq \lambda'$.
\item $SP_{n_\lambda}(q)$, if $\lambda=\lambda'\neq (\lambda_2,\lambda_1)$ and ($p=2$ or $\nu(\lambda)=-1$).
\item $\Omega^+_{n_\lambda}(q)$,if $\lambda=\lambda'\neq (\lambda_2,\lambda_1)$, $p\neq 2$ and $\nu(\lambda)=1$.
\item $SP_{n_\lambda}(q^\frac{1}{2})$, if $\lambda=\lambda'=(\lambda_2,\lambda_1)$ and ($p=2$ or $\nu(\lambda)=-1$).
\item $\Omega_{n_\lambda}^+(q^\frac{1}{2})$ if $\lambda=\lambda'= (\lambda_2,\lambda_1)$, $p\neq 2$ and $\nu(\lambda)=1$.
\end{enumerate}
\item When $\F_q=\F_p(\alpha,\beta)=\F_p(\alpha+\alpha^{-1},\beta)\neq \F_p(\alpha,\beta+\beta^{-1})$,
\begin{enumerate}
\item when $\F_{\tilde{q}}=\F_p(\alpha)=\F_p(\alpha+\alpha^{-1})$,
\begin{enumerate}
\item $SL_{n-1}(\tilde{q})$ if $\lambda =([n-1,1],\emptyset)$.
\item $SL_{n_\lambda}(\tilde{q})$ if $\lambda = (\lambda_1,\emptyset), \lambda_1< \lambda_1'$.
\item $SP_{n_\lambda}(\tilde{q})$ if $\lambda = (\lambda_1,\emptyset), \lambda_1=\lambda_1'$ and ($p=2$ or ($p\geq 3$ and $\nu(\lambda_1) = -1$)).
\item $\Omega_N^{+}(\tilde{q})$ if $\lambda=(\lambda_1,\emptyset), \lambda_1=\lambda_1'$, $p\geq 3$ and $\nu(\lambda_1)=1$.
\end{enumerate}
\item when $\F_{\tilde{q}}=\F_p(\alpha)\neq \F_p(\alpha+\alpha^{-1})$,
\begin{enumerate}
\item $SU_{n-1}(\tilde{q}^\frac{1}{2})$ if $\lambda =([n-1,1],\emptyset)$.
\item $SU_{n_\lambda}(\tilde{q}^\frac{1}{2})$ if $\lambda = (\lambda_1,\emptyset), \lambda_1< \lambda_1'$.
\item $SP_{n_\lambda}(\tilde{q}^\frac{1}{2})$ if $\lambda = (\lambda_1,\emptyset), \lambda_1=\lambda_1'$ and ( $p=2$ or ($p\geq 3$ and $\nu(\lambda_1) = -1$)).
\item $\Omega_N^{+}(\tilde{q}^\frac{1}{2})$ if $\lambda=(\lambda_1,\emptyset), \lambda_1=\lambda_1'$, $p\geq 3$ and $\nu(\lambda_1)=1$.
\end{enumerate}
\item $SL_n(q)$ if $\lambda=([1],[n-1])$,
\item $SL_{n_\lambda}(q)$ if $\lambda\neq \lambda',\lambda\neq(\lambda_1',\lambda_2')$ and $\lambda\neq (\lambda_2,\lambda_1)$, 
\item $SL_{n_\lambda}(\F_{q^\frac{1}{2}})$ if $\lambda=(\lambda_2,\lambda_1) \neq \lambda'$,
\item  $SU_{n_\lambda}(q^\frac{1}{2})$ if $\lambda=(\lambda_1',\lambda_2') \neq \lambda'$,
\item  $SP_{n_\lambda}(q)$ if $\lambda=\lambda'\neq (\lambda_2,\lambda_1)$ and ($p=2$ or $\nu(\lambda)=-1$),
\item $\Omega^+_{n_\lambda}(q)$ if $\lambda=\lambda'\neq (\lambda_2,\lambda_1)$, $p\neq 2$ and $\nu(\lambda)=1$, 
\item  $SP_{n_\lambda}(q^\frac{1}{2})$ if $\lambda=\lambda'=(\lambda_2,\lambda_1)$ and ($p=2$ or $\nu(\lambda)=-1$),
\item  $\Omega_{n_\lambda}^+(q^\frac{1}{2})$ if $\lambda=\lambda'= (\lambda_2,\lambda_1)$, $p\neq 2$ and $\nu(\lambda)=1$.
\end{enumerate}
\end{enumerate}
\end{theo}

\begin{proof}

It is sufficient to prove the result for double-partitions with no empty components which are not hooks. We know by Theorem \ref{CGFS} and the previous steps that $G(\lambda)$ is a classical group in a natural representation. The proof uses Proposition \ref{isomorphisme} and the separation of the cases made before the enumeration of the six steps. We write $\F_{q'}$ the field over which our classical group is defined. In all cases $G(\lambda) \subset SL_n(q)$, therefore $q'$ divides $q$.

Assume $\F_q=\F_p(\alpha,\beta)=\F_p(\alpha,\beta+\beta^{-1})\neq \F_p(\alpha+\alpha^{-1},\beta).$
\begin{enumerate}
\item If $\lambda\neq  \lambda', \lambda\neq (\lambda_1',\lambda_2')$ and $\lambda\neq (\lambda_2,\lambda_1)$, then $G(\lambda)$ contains a natural $SL_2(q)$. By Lemma \ref{field}, we have that $q'=q$. By Proposition \ref{isomorphisme}, $G(\lambda)$ preserves no hermitian or bilinear form, therefore $G(\lambda)=SL_{n_\lambda}(q)$.
\item If $\lambda=(\lambda_2,\lambda_1)\neq \lambda'$, then by Proposition \ref{isomorphisme} and Lemma \ref{Harinordoquy}, we have up to conjugation $G(\lambda) \subset SU_{n_\lambda}(q^\frac{1}{2})$. Up to conjugation, $G(\lambda)$ contains $\{\diag(M,{}^t\!\epsilon(M^{-1}),I_{n_\lambda-6}), M\in SL_3(q)\}$, therefore $G(\lambda)$ contains $\{\diag(M,M,I_{n_\lambda-6}),M\in SU_3(q^\frac{1}{2})\}$.
                          
Let $\varphi$ be the natural representation of $SU_3(q^\frac{1}{2})$ in $GL_3(\overline{\F_p})$ and $\rho$ the diagonal representation of $SU_3(q^\frac{1}{2})$ in $GL_{n_\lambda}(\overline{\F_p})$, given by the above subgroup of $G(\lambda)$. 

We have $\rho\simeq \varphi \oplus \varphi \oplus \mathbf{1}^{n_\lambda-6}$ with $\mathbf{1}$ the trivial representation. Let $\sigma$ be a generator of $\Gal(\F_q/\F_{q'})$. Since $G(\lambda)$ is a classical group over $\F_{q'}$, we have that $\rho\simeq \rho^\sigma$, therefore $\varphi \simeq \varphi^{\sigma}$. It follows that for every $M\in SU_3(q^\frac{1}{2})$, we have $\sigma(\tr(M))=\tr(M)$. By Lemma \ref{field}, we have that $\F_q=\F_q^{\Gal(\F_q/\F_{q'})}$ and, therefore $q'=q$. By Proposition \ref{isomorphisme}, $G(\lambda)$ preserves no bilinear form, therefore $G(\lambda)=SU_{n_\lambda}(q^\frac{1}{2})$.
\item If $\lambda=(\lambda_1',\lambda_2')\neq \lambda'$, then by Proposition \ref{isomorphisme} and Lemma \ref{Ngwenya}, up to conjugation, we have that $G(\lambda)\subset SL_{n_\lambda}(q^\frac{1}{2})$. The group $G(\lambda)$ contains either a natural $SL_3(q^\frac{1}{2})$ or a group of the form $\{\diag( M,\epsilon(M),I_{n_\lambda-6}), M\in SL_3(q)\}$.
 
 If $G(\lambda)$ contains a natural $SL_3(q^\frac{1}{2})$ then by Lemma \ref{field}, we have $q'=q^\frac{1}{2}$. We know by Proposition \ref{isomorphisme} that $G(\lambda)$ preserves no symmetric or skew-symmetric bilinear form. If we had $G(\lambda) \subset SU_{n_\lambda}(q^\frac{1}{4})$, then the natural $SL_3(q^\frac{1}{2})$ in $G(\lambda)$ would inject itself in some $SU_3(q^\frac{1}{4})$.  This is absurd because of their orders, therefore we have $G(\lambda)\simeq SL_{n_\lambda}(q^\frac{1}{2})$.
 
 If $G$ contains up to conjugation a group of the form  $\{\diag(M,\epsilon(M),I_{n_\lambda-6}), M\in SL_3(q)\}$ then it contains $\{\diag(M,M,I_{n_\lambda-6}), M\in SL_3(q^{\frac{1}{2}})\}$.  Let $\varphi$ be the natural representation of $SL_3(q^\frac{1}{2})$ in $GL_3(\overline{\F_p})$ and $\rho$ the diagonal representation of $SL_3(q^\frac{1}{2})$ in $GL_{n_\lambda}(\overline{\F_p})$ given by the above subgroup of $G(\lambda)$. We then have $\rho \simeq \varphi \oplus \varphi \oplus \mathbf{1}^{n_\lambda-6}$. Let $\sigma$ be a generator of $\Gal(\F_{q^\frac{1}{2}}/\F_{q'})$. We have $\rho \simeq \rho^\sigma$, therefore $\varphi\simeq \varphi^\sigma$. By Lemma \ref{field}, we have that $\F_{q^\frac{1}{2}}=\F_{q^\frac{1}{2}}^{\Gal(\F_{q^\frac{1}{2}}/\F_{q'})}$, therefore $q'=q^\frac{1}{2}$. We cannot have $G(\lambda)\simeq SU_{n_\lambda}(q^\frac{1}{4})$ because $SL_3(q)$ would inject itself in $SU_6(q^\frac{1}{4})$ and we know that $\frac{\vert SU_{6}(q^\frac{1}{4})\vert}{q^\frac{15}{4}} < \frac{\vert SL_3(q)\vert}{q^3}$. By Proposition \ref{isomorphisme}, $G(\lambda)$ cannot preserve any symmetric or skew-symmetric bilinear form, therefore $G(\lambda)\simeq SL_{n_\lambda}(q^\frac{1}{2})$.
 
\item Case $4$ is analogous to Case $3$.
     
\item If $\lambda=\lambda'=(\lambda_2,\lambda_1)$, then by Proposition \ref{isomorphisme} and Proposition \ref{coolprop}, $G(\lambda)$ preserves a bilinear form of the type given by Proposition \ref{bilin} defined over $\F_{q^\frac{1}{2}}$. This shows that $q'\leq q^\frac{1}{2}$ and it is enough to show that $q'=q^\frac{1}{2}$ to conclude the proof.

If $\lambda_1$ and $\lambda_2$ are square partitions, then $G(\lambda)$ contains up to conjugation the group\break $\{\diag(M,{}^t\!(M^{-1}),I_{n_\lambda-6}), M\in SU_3(q^{\frac{1}{2}})\}$.
    
Let $\varphi$ be the natural representation of $SU_3(q^\frac{1}{2})$ in $GL_3(\overline{\F_p})$, and $\rho$ the twisted diagonal representation of $SU_3(q^\frac{1}{2})$ in $GL_{n_\lambda}(\overline{\F_p})$ given by the above subgroup of $G(\lambda)$. We have $\rho\simeq \varphi \oplus \varphi^\star \oplus \mathbf{1}^{n_\lambda-6}$. Let $\sigma$ be a generator of $\Gal(\F_{q^\frac{1}{2}}/\F_{q'})$. Since $G(\lambda)$ is a classical group over $\F_{q'}$, we have $\rho\simeq \rho^\sigma$. It follows that $\varphi \simeq \varphi^{\sigma}$ or $\varphi\simeq (\varphi^\star)^\sigma$. The first possibility implies that $\F_{q^\frac{1}{2}}=\F_{q^\frac{1}{2}}^{\Gal(\F_{q^\frac{1}{2}}/\F_{q'})}$, therefore $q'=q^\frac{1}{2}$. The second possibility implies that $\varphi \simeq \varphi^{\sigma^2}$. Therefore $q'=q^\frac{1}{2}$ or $\sigma$ is of order $2$ and $SU_3(q^\frac{1}{2})$ injects into $SU_3(q^\frac{1}{4})$, which is a contradiction. In both cases, we have $q'=q^\frac{1}{2}$ and the desired result follows.

If $\lambda_1$ or $\lambda_2$ is not a square partition, then $G(\lambda)$ contains up to conjugation the group \break$\{\diag(M, {}^t\!(M^{-1}),\epsilon(M),{}^t\!\epsilon(M^{-1}),I_{n_\lambda-12}), M\in SL_3(q)\}$, and, therefore contains its subgroup $\{\diag(M,{}^t\!(M^{-1}),M,{}^t\!(M^{-1}),I_{n_\lambda-12}), M\in SL_3(q^\frac{1}{2})\}$.  Let $\varphi$ be the natural representation of $SL_3(q^\frac{1}{2})$ in $GL_3(\overline{\F_p})$ and $\rho$ be the representation of $SL_3(q^\frac{1}{2})$ in $GL_{n_\lambda}(\overline{\F_p})$ given by the above subgroup $G(\lambda)$. We have $\rho\simeq \varphi \oplus \varphi \oplus \varphi^\star\oplus \varphi^\star \oplus \mathbf{1}^{n_\lambda-6}$ . Let $\sigma$ be a generator of $\Gal(\F_{q^\frac{1}{2}}/\F_{q'})$. Since $G(\lambda)$ is a classical group defined over $\F_{q'}$, we have that $\rho\simeq \rho^\sigma$. It follows that $\varphi \simeq \varphi^{\sigma}$ or $\varphi\simeq (\varphi^\star)^\sigma$. By the same arguments as before, we have $q'=q^\frac{1}{2}$ or $SL_3(q^\frac{1}{2})$ injects itself in $SU_3(q^\frac{1}{4})$, which is not possible. This proves $q'=q^\frac{1}{2}$ and concludes the case $\F_q=\F_p(\alpha,\beta)=\F_p(\alpha,\beta+\beta^{-1})\neq \F_p(\alpha+\alpha^{-1},\beta).$ 
\end{enumerate}
 
 If $\F_q=\F_p(\alpha,\beta)=\F_p(\alpha+\alpha^{-1},\beta)\neq \F_p(\alpha,\beta+\beta^{-1})$, then all the arguments are the same up to permutation of the different cases.
\end{proof}

We have determined the image of the derived subgroup of the Artin group in all cases in type $B$. This is close to determining the image of the Artin group itself since $A_{B_n}/\mathcal{A}_{B_n}\simeq \Z^2$. We give in the following subsection an example of how to recover the group $G=R(A_{B_n})$ from $H=R(\mathcal{A}_{B_n})$ with the representation of $\mathcal{H}_{B_2,\alpha,\beta}$ labeled by the double-partition $([1],[1])$.

\section{Image of the full Artin group for the $2$-dimensional representation of $\mathcal{H}_{B_2,\alpha,\beta}$}

We first recall the results for the $2$-dimensional representation depending on the fields.
\begin{prop}\label{resultB2}
Assume the order of $\beta$ does not belong to $\{1,2,3,4,5,6,10\}$ or the order of $\alpha$ does not belong to $\{1,2,3,4,5,6,10\}$

\begin{enumerate}
\item If $\F_q=\F_p(\alpha,\beta)=\F_p(\alpha+\alpha^{-1},\beta+\beta^{-1})$ then $R_{[1],[1]}(\mathcal{A}_{B_2})=SL_2(q)$.
\item If $\F_q=\F_p(\alpha,\beta)=\F_p(\alpha,\beta+\beta^{-1})=\F_p(\alpha+\alpha^{-1},\beta)\neq\F_p(\alpha+\alpha^{-1},\beta+\beta^{-1})$ then $R_{[1^2],[1^2]}(\mathcal{A}_{B_2})\simeq SU_2(q^{\frac{1}{2}})$.
\item If $\F_q=\F_p(\alpha,\beta)=\F_p(\alpha,\beta+\beta^{-1})\neq \F_p(\alpha+\alpha^{-1},\beta)=\F_p(\alpha+\alpha^{-1},\beta+\beta^{-1})$ then $R_{[1^2],[1^2]}(\mathcal{A}_{B_2})\simeq SL_2(q^{\frac{1}{2}})$.
\item If $\F_q=\F_p(\alpha,\beta)=\F_p(\alpha+\alpha^{-1},\beta)\neq \F_p(\alpha,\beta+\beta^{-1})=\F_p(\alpha+\alpha^{-1},\beta+\beta^{-1})$ then $R_{[1^2],[1^2]}(\mathcal{A}_{B_2})\simeq SL_2(q^{\frac{1}{2}})$.
\end{enumerate}
\end{prop}

\begin{cor}
Under the same assumptions as in the previous proposition, we have that if $G=R_{[1],[1]}(A_{B_2})$, $a$ is the order of $-\alpha$ and $b$ is the order of $-\beta$ then 
\begin{enumerate}
\item If $\F_q=\F_p(\alpha,\beta)=\F_p(\alpha+\alpha^{-1},\beta+\beta^{-1})$ then
 $G \simeq SL_2(q) \rtimes \Z/\op{lcm}(a,b)\Z$.
\item $\F_p(\alpha,\beta)=\F_p(\alpha,\beta+\beta^{-1})=\F_p(\alpha+\alpha^{-1},\beta)\neq \F_p(\alpha+\alpha^{-1},\beta+\beta^{-1})$ then  $G \simeq SU_2(q^{\frac{1}{2}}) \rtimes \Z/\op{lcm}(a,b)\Z$.
 \end{enumerate}
\end{cor}

\begin{proof}

We write $s=R_{[1],[1]}(S_1)=\begin{pmatrix}
\frac{\alpha-1}{\beta+1} & \frac{\alpha\beta+1}{\beta+1}\\
\frac{\alpha+\beta}{\beta+1}& \frac{\alpha\beta-\beta}{\beta+1}
\end{pmatrix}$ and $t=R_{[1],[1]}(T)=\begin{pmatrix}
\beta & 0 \\
0 & -1
\end{pmatrix}$.

\bigskip

Assume first that $\F_q=\F_p(\alpha,\beta)=\F_p(\alpha+\alpha^{-1},\beta+\beta^{-1})$. By Proposition \ref{resultB2} and the fact that $\op{Det}(s)=-\alpha$ and $\op{Det}(t)=-\beta$, we have the following exact sequence
$$1\longrightarrow SL_2(q) \longrightarrow G \overset{det}\longrightarrow <-\alpha,-\beta> \longrightarrow 1.$$

In order to prove the result, we must show that the exact sequence is split and that $<-\alpha,-\beta>$ is isomorphic to $\Z/\op{lcm}(a,b)\Z$. We first show that the exact sequence is split.

\medskip

This is equivalent to finding a subgroup $N$ of $G$ isomorphic to $<-\alpha,-\beta>$ and such that $\op{Det}(N)= <-\alpha,-\beta>$. We have $SL_2(q)\leq G$, therefore $U=\begin{pmatrix}
0 & 1\\
-1 & 0
\end{pmatrix}\in G$. It follows that $V=stsUt^{-1}U = \begin{pmatrix}
-\alpha & 0\\
0 & -\alpha
\end{pmatrix}\in G$. We also have $sts = \begin{pmatrix} -\alpha & 0\\
0 & \alpha\beta \end{pmatrix}\in G$ and $-t =\begin{pmatrix}
-\beta & 0\\
0 & 1
\end{pmatrix}\in G$. We now distinguish five possibilities depending on the orders of $\alpha$ and the order of $\beta$. We write $c=2^\ell c'$ (resp $d=2^kd'$) the order of 
$\alpha$ (resp $\beta$) with $c'$ (resp $d'$ odd).

\smallskip

\textbf{First case} : $k=1$ or $\ell=1$. By symmetry of the roles of $\alpha$ and $\beta$, it is sufficient to show the exact sequence is split when $l=1$. We then have $(-\alpha)^{\frac{c}{2}}=(-1)^{c'}\alpha^{c'}=1$, therefore $-\alpha$ is of order $m$ for some $m$ dividing $\frac{c}{2}$. We also have $\alpha^{2m}=(-\alpha)^{2m}=1$, therefore $c$ divides $2m$. It follows that $m=\frac{c}{2}$. The order of $\alpha$ being even, we have that the order of $\alpha^2$ is also $\frac{c}{2}$. The subgroup generated by $V$ and $-t$ then verifies the desired conditions.

\smallskip

\textbf{Second case} : $k=\ell=0$. We then have that the order of $-\alpha$ is equal to the order of $-\alpha^2$. We have 
 $(sts)^{cd}=\begin{pmatrix}
(-1)^{cd}(\alpha^{c})^d & 0\\
0 & (\alpha^c)^d(\beta^d)^c
\end{pmatrix}=\begin{pmatrix} -1 & 0\\
0 & 1\end{pmatrix}$, therefore $M=(sts)^{cd}V=\begin{pmatrix}
\alpha & 0\\
0 & -\alpha
\end{pmatrix}\in G$. The subgroup $<M,-t>$ then verifies the desired conditions.

\smallskip

\textbf{Third case} : $k>1>\ell$ or $\ell>1>k$. It is sufficient to consider the case $k>1>l$. We then have that $-\alpha$ and $-\alpha^2$ have the same order since $c$ is odd. We have $t^{\frac{d}{2}}=\begin{pmatrix} \beta^{\frac{d}{2}} & 0\\
0 & (-1)^{\frac{d}{2}}\end{pmatrix} = \begin{pmatrix} -1 & 0\\
0 & 1\end{pmatrix}$, therefore $M=t^{\frac{d}{2}}V = \begin{pmatrix}
\alpha & 0\\
0 & -\alpha
\end{pmatrix}\in G$0 The subgroup $<M,-t>$ then verifies the desired conditions.

\smallskip

\textbf{Fourth case} : $k>\ell>1$ or $\ell>k>1$. It is sufficient to consider the case $k>\ell>1$. We then have that the order of $\alpha$ is equal to the order of $-\alpha$. We have $t^{\frac{d}{2^l}}=\begin{pmatrix} \beta^{2^{k-l}d'} & 0\\
0 & (-1)^{2^{k-l}d'}\end{pmatrix} = \begin{pmatrix}
\beta^{2^{k-l}d'} & 0\\
0 & 1
\end{pmatrix}$. We set $\gamma = \beta^{2^{k-l}d'}$, $\gamma$ is then of order $2^\ell$. We have that $\gamma\alpha^2$ is of order $r$, where $r$ is an integer dividing $c$. We also have that $(\gamma\alpha^2)^{\frac{c}{2}}=\gamma^{2^{\ell-1}c'}\alpha^{c}=-1$, therefore $r$ does not divide $\frac{c}{2}$. This implies that there exists an odd integer $c''$ such that $r=2^\ell c''$ and $c''$ divides $c'$. We have $1=(\gamma\alpha^2)^r= \alpha^{2r}$ and $c$ divides $2r$. Therefore $c'$ divides $c''$ and $\gamma\alpha^2$ is of order $a$. We then set $M= t^{\frac{d}{2^l}}V=\begin{pmatrix}-\gamma\alpha & 0\\
0 & -\alpha
\end{pmatrix}\in G$. The matrix $M$ is of order $c=a$ and its determinant is of order $a$ because it is equal to $\gamma\alpha^2$. This proves that we can take $<M,-t>$ as our subgroup.

\smallskip

\textbf{Fifth case} : $k=l\geq 2$. We then have that $-\alpha$ is of order $c$. We have $t^{d'}=\begin{pmatrix}
\beta^{d'} & 0\\
0 & (-1)^{d'}
\end{pmatrix}= \begin{pmatrix}
\beta^{d'} & 0\\
0 & -1
\end{pmatrix}$. We set $\gamma =\beta^{d'}$, we have that the order of $\gamma$ is $2^k$. The element $-\gamma\alpha^2$ is then of order $c$ by the same reasonning as in the fourth case. We set $M=t^{d'}V= \begin{pmatrix}
-\gamma\alpha & 0\\
0 & \alpha
\end{pmatrix}\in G$. The matrix $M$ is of order $c$, therefore the subgroup $<M,-t>$ verifies the desired conditions. 

\medskip

This shows that $G \simeq G \rtimes <-\alpha,-\beta>$. It now only remains to show that $<-\alpha,-\beta> \simeq \Z/\op{lcm}(a,b)\Z$. 

\smallskip

Let $n=q-1$, we have $<-\alpha,-\beta> \simeq <\frac{n}{a},\frac{n}{b}>$ when it is seen as a subgroup of $\Z/n\Z$. Let $\varphi :\Z/a\Z \times \Z/b\Z \rightarrow <\frac{n}{a},\frac{n}{b}>$ be the map that maps $(u,v)$ to $u\frac{n}{a}+v\frac{n}{b}$. The map $\varphi$ is a group epimorphism and its Kernel is  $\op{ker}(\varphi) = \{(u,v)\in\Z/a\Z\times \Z/b\Z,u\frac{n}{a}+v\frac{n}{b}=0\}$. Let $d=\op{Gcd}(a,b)$. We have a group isomorphism  $\psi$ from $\Z/d\Z$ to $\op{ker}(\varphi)$ defined by $\psi(k) = (ka/d,-kb/d)$. The map $\psi$ is well-defined and is clearly one-to-one. We must show that it is onto. Let $u\in\Z$, $v\in \Z$ such that $u\frac{n}{a}+v\frac{n}{b}=\ell n$ for some integer $\ell$. We have that $unb=(\ell b-v)na$ and, therefore, $ub=(\ell b-v)na$ and $u\frac{b}{d}=(\ell b-v)\frac{a}{d}$. It follows that $\frac{a}{d}$ divides $u\frac{b}{d}$ and $\frac{a}{d}$ divides $u$. This implies that there exists $k\in \Z$ such that $u=k\frac{a}{d}$ and the projection of $u$ in $\Z/a\Z$ only depends on the projection of $k$ in $\Z/d\Z$. It follows that $\psi$ is indeed onto.

It follows that $<-\alpha,-\beta>\simeq (\Z/a\Z\times \Z/b\Z)/\Z/d\Z$. Its order is therefore equal to $\frac{ab}{d}=\op{lcm}(a,b)$. Since it is a subgroup of $\Z/n\Z$, it is cyclic and therefore isomorphic $\Z/\op{lcm}(a,b)\Z$. This concludes the proof for $\F_p(\alpha,\beta)=\F_p(\alpha+\alpha^{-1},\beta+\beta^{-1})$.

\bigskip

Assume now that $\F_p(\alpha,\beta)=\F_p(\alpha,\beta+\beta^{-1})=\F_p(\alpha+\alpha^{-1},\beta)\neq \F_p(\alpha+\alpha^{-1},\beta+\beta^{-1})$. Recall that $\epsilon(\alpha)=\alpha^{-1}$ and $\epsilon(\beta)=\beta^{-1}$. Let $P=\begin{pmatrix}
\frac{\alpha+\beta}{\alpha\beta+1} & 0\\
0 & 1
\end{pmatrix}$, we then have $PsP^{-1}=\epsilon(^t\!s^{-1})$ and $PtP^{-1}=\epsilon(^t\!t^{-1})$. If we set $(X,Y)=^t\epsilon(X)PY$, then for all $M\in G$, $X,Y\in \F_q^2$, we have $(MX,MY)=(X,Y)$. This proves by Proposition \ref{Ngwenya} that there exists $Q\in GL_2(q)$ such that $\tilde{G}=QGQ^{-1}$ is a subgroup of $GU_2(q^{\frac{1}{2}})$. By Proposition \ref{resultB2}, we have the following exact sequence
$$1\longrightarrow SU_2(q^{\frac{1}{2}}) \longrightarrow \tilde{G} \overset{\op{det}}\longrightarrow <-\alpha,-\beta> \longrightarrow 1.$$
As in the previous case, it is sufficient to show that there is a splitting of $<-\alpha,-\beta>$ in $\tilde{G}$ or in $G$ since $G\simeq Q^{-1}\tilde{G}Q$. We have that $[\tilde{G},\tilde{G}]\simeq SU_2(q^{\frac{1}{2}})$ and for $M\in GL_2(\F_q)$, $M \in GU_2(q)$ if and only if $QMQ^{-1}$ stabilizes the sesquilinear form defined by $(X,Y) =^t\epsilon(X)PY$. We have that $G$ contains all the matrices of determinant $1$ preserving the above sesquilinear form. Let $M =\begin{pmatrix} \beta & 0\\
0 & \beta^{-1}\end{pmatrix}$. We have $\epsilon(^t(M^{-1}))P = MP = PM$, therefore $(MX,MY) = (X,Y)$. Since $\op{det}(M) = 1$, we have that $M\in G$. In the same way, we have $-I_2\in G$. We then have that $\begin{pmatrix}
-\alpha & 0\\
0 & -\alpha
\end{pmatrix}=Mstst^{-1}\in G$ and $-t\in G$. This proves that we can have all the matrices appearing in the previous case in $G$ and, therefore that there exists a splitting of $<-\alpha,-\beta>$ in $G$. This concludes the proof in this case.
\end{proof}

If $\F_q=\F_p(\alpha,\beta)=\F_p(\alpha,\beta+\beta^{-1})\neq \F_p(\alpha+\alpha^{-1},\beta)=\F_p(\alpha+\alpha^{-1},\beta+\beta^{-1})$ or $\F_q=\F_p(\alpha,\beta)=\F_p(\alpha+\alpha^{-1},\beta)\neq \F_p(\alpha,\beta+\beta^{-1}) =\F_p(\alpha+\alpha^{-1},\beta+\beta^{-1})$ then the exact sequence we consider is not always split, therefore the situation is slightly more complex. The $\op{Gcd}$ of the order of $-\alpha$ and the order of $-\beta$ then divides $2$ as we prove in the following lemma.

\begin{lemme}
If $\F_q=\F_p(\alpha,\beta)=\F_p(\alpha,\beta+\beta^{-1})\neq \F_p(\alpha+\alpha^{-1},\beta)=\F_p(\alpha+\alpha^{-1},\beta+\beta^{-1})$ or $\F_q=\F_p(\alpha,\beta)=\F_p(\alpha+\alpha^{-1},\beta)\neq \F_p(\alpha,\beta+\beta^{-1}) =\F_p(\alpha+\alpha^{-1},\beta+\beta^{-1})$ then $\op{Gcd}(a,b)\leq 2$, where $a$ is the order of $-\alpha$ and $b$ is the order of $-\beta$.
\end{lemme}

\begin{proof}

Since the roles of $\alpha$ and $\beta$ are symmetric, we can assume that $\F_q=\F_p(\alpha,\beta)=\F_p(\alpha,\beta+\beta^{-1})\neq \F_p(\alpha+\alpha^{-1},\beta)=\F_p(\alpha+\alpha^{-1},\beta+\beta^{-1})$. 

Let $a'=\frac{a}{d}$ and $b'=\frac{b}{d}$. The group $\F_q^\star$ is isomorphic to $\Z/n\Z$, where $n=q-1$. Since $-\alpha$ is of order $a$, its image in $\Z/n\Z$ is of the form $\overline{u\frac{n}{a}}$ with $u$ coprime to $a$. The group generated by $\beta$ is mapped to the group generated by $\overline{\frac{n}{b}}$. The image of $(-\alpha)^{a'}$ is then equal to $\overline{ua'\frac{n}{a}}=\overline{u\frac{a}{d}\frac{n}{a}}=\overline{u\frac{b}{d}\frac{n}{b}}\in <\overline{\frac{n}{b}}>$. This proves that $(-\alpha)^{a'}\in <-\beta>$, therefore $(-\alpha)^{a'}\in \F_{q^{\frac{1}{2}}}$. It follows that the polynomial $R(X)=X^{a'}-(-\alpha)^{a'}$ has its coefficients in $\F_{q^{\frac{1}{2}}}$.

 Since $\F_p(\alpha+\alpha^{-1},\beta)\neq \F_p(\alpha,\beta)$, the polynomial $X^2+(\alpha+\alpha^{-1})X+1$ is irreducible over $\F_{q^\frac{1}{2}}$ and the unique automorphism $\epsilon$ of order $2$ of $\F_q$ verifies $\epsilon(-\alpha)=-\alpha^{-1}$. Since $-\alpha$ is a root of $R\in \F_{q^{\frac{1}{2}}}[X]$, $\epsilon(-\alpha)=-\alpha^{-1}$ is also a root of $R$. If $k\in [\![0,a'-1]\!]$, we have that $((-\alpha)^{1+kd})^{a'}=(-\alpha)^{a'}$ is a root of $R$ and since $kd\in[\![0,a-d]\!]$, those roots are distinct. This proves that those are all the roots of $R$ since its degree is $a'$. It follows that there exists $k\in [\![0,a'-1]\!]$ such that $(-\alpha)^{-1}=(-\alpha)^{1+kd}$ and, therefore $(-\alpha)^{2+kd}=1$. It follows that $a$ divides $2+kd$. We have $2\leq 2+kd\leq 2+(a'-1)d=2+a-d$ and $a>1$, therefore $a=2+kd$. If $d\geq 3$ then we have $2+kd\leq 2+a-d\leq a-1$ do we cannot have $a=2+kd$. This proves by contradiction that $d\leq 2$. Note that we can have $d=2$ or $d=1$ since $2+2(a'-1)=2a'$ and $2+(a'-2)=a'$.
\end{proof}

\begin{prop}
Assume $\F_q=\F_p(\alpha,\beta)=\F_p(\alpha,\beta+\beta^{-1})\neq \F_p(\alpha+\alpha^{-1},\beta)=\F_p(\alpha+\alpha^{-1},\beta+\beta^{-1})$ or $\F_q=\F_p(\alpha,\beta)=\F_p(\alpha+\alpha^{-1},\beta)\neq \F_p(\alpha,\beta+\beta^{-1}) =\F_p(\alpha+\alpha^{-1},\beta+\beta^{-1})$. Let $d=\op{Gcd}(a,b)\leq 2$, where $a$ is the order of $-\alpha$ and $b$ is the order of $-\beta$.

If $d=1$ then we have $G\simeq (SL_2(q^{\frac{1}{2}})\rtimes \Z/ab\Z$.

If $d=2$, then we have $G \simeq J\rtimes \Z/\op{lcm}(a,b)\Z$, where $J=\{M\in SL_2(q), M\epsilon(M)^{-1} \in \{I_2,-I_2\}\}$ and $\epsilon$ is the unique automorphism of order $2$ of $\F_q$.
\end{prop}

\begin{proof}
Again, we can assume $\F_q=\F_p(\alpha,\beta)=\F_p(\alpha,\beta+\beta^{-1})\neq \F_p(\alpha+\alpha^{-1},\beta)=\F_p(\alpha+\alpha^{-1},\beta+\beta^{-1})$. Let $P =\begin{pmatrix}
-\frac{(\alpha+\beta)(\alpha-1)}{\alpha(\beta+1)} & 0\\
0 & 1
\end{pmatrix}$. We will first show that $P[G,G]P^{-1}\subset SL_2(q^{\frac{1}{2}})$.

\smallskip

Recall that we have $t=\begin{pmatrix}
\beta &0 \\
0 & -1
\end{pmatrix}$ and $s=\begin{pmatrix}
\frac{\alpha-1}{1+\beta}& \frac{\alpha+\beta^{-1}}{1+\beta^{-1}}\\
\frac{\alpha+\beta}{1+\beta} & \frac{\alpha-1}{1+\beta^{-1}}
\end{pmatrix}$. Note that $\beta\in \F_{q^{\frac{1}{2}}}$, $\alpha+\alpha^{-1}\in \F_{q^{\frac{1}{2}}}$ and for $x\in \F_q$, $x\in \F_{q^{\frac{1}{2}}}$ if and only if $\epsilon(x)=x$ with $\epsilon$ the unique automorphism of order $2$ of $\F_q$. We have $
\epsilon(\beta)=\beta$ and $\epsilon(\alpha)=\alpha^{-1}$. Let $Q=\epsilon(P^{-1})P$, we have 
$$Q=\begin{pmatrix}\frac{\alpha^{-1}(\beta+1)}{(\alpha^{-1}+\beta)(\alpha^{-1}-1)}\frac{(\alpha+\beta)(\alpha-1)}{\alpha(\beta+1)} & 0\\
0 & 1\end{pmatrix}=\begin{pmatrix} \frac{(\beta+1)(\alpha+\beta)(\alpha-1)}{(\alpha\beta+1)(1-\alpha)(\beta+1)} & 0\\
0 & 1\end{pmatrix} =\begin{pmatrix}
-\frac{\alpha+\beta}{\alpha\beta+1} & 0\\
0 & 1
\end{pmatrix}.$$

We then have $QtQ^{-1}=t=\epsilon(t)$ and 
$$QsQ^{-1}=\begin{pmatrix}
\frac{\alpha-1}{\beta+1} &-\frac{\alpha+\beta}{\beta+1} \\
-\frac{\alpha\beta+1}{\beta+1} & \frac{(\alpha-1)\beta}{\beta+1}
\end{pmatrix}=\begin{pmatrix}
\frac{-\alpha(1-\alpha^{-1})}{1+\beta} & -\frac{-\alpha(1+\beta\alpha^{-1})}{\beta+1}\\
\frac{-\alpha(\beta+\alpha^{-1})}{\beta+1} & \frac{-\alpha(1-\alpha^{-1})}{1+\beta^{-1}}
\end{pmatrix}=-\alpha\epsilon(s).$$

It follows that for all $g\in [G,G]$, $QgQ^{-1}=\epsilon(g)$ and, therefore $\epsilon(P^{-1})PgP^{-1}\epsilon(P)=\epsilon(g)$ and $\epsilon(PgP^{-1})=PgP^{-1}$. This proves that we have $P[G,G]P^{-1}\subset SL_2(q^{\frac{1}{2}})$. We write in the following $\tilde{G}=PGP^{-1}$, $\tilde{s}$ the image of $PsP^{-1}$ in $\tilde{G}/[\tilde{G},\tilde{G}]$ and $\tilde{t}$ the image of $PtP^{-1}$ in $\tilde{G}/[\tilde{G},\tilde{G}]$.

\bigskip

We have by Proposition \ref{resultB2} that $[\tilde{G},\tilde{G}]\simeq SL_2(q^{\frac{1}{2}})$. We now show that $a$ is the order of $\tilde{s}$ and $b$ is the order of $\tilde{t}$. Let $r$ be the order of $\tilde{s}$. We have $(PsP^{-1})^r\in SL_2(q^{\frac{1}{2}})$, therefore $(-\alpha)^r=\op{det}(s)^r=1$ and $a$ divides $r$. We have $(-\alpha)^a=1$, therefore $\alpha^a=(-1)^a$. The eigenvalues of $PsP^{-1}$ are $\alpha^{-1}$ and $-1$, therefore $(PsP^{-1})^a$ is conjugate to the diagonal matrix $\diag(\alpha^{-a},(-1)^a)=(-1)^aI_2$ and, therefore $(PsP^{-1})^a=(-1)^aI_2\in SL_2(q^{\frac{1}{2}})$. It follows that $\tilde{s}^a=1$ and, therefore $r$ divides $a$. This proves that $r=a$. In the same way, we have $b$ is the order of $\tilde{t}$.

\bigskip

We now determine $\tilde{G}$ in terms of $\tilde{G}\cap SL_2(q)$. The determinant gives us the following exact sequence
$$1 \longrightarrow \tilde{G} \cap SL_2(q) \longrightarrow \tilde{G}\overset{\det}\longrightarrow <-\alpha,-\beta> \longrightarrow 1.$$
The matrix $P$ commutes with all diagonal matrices, therefore $sts=\begin{pmatrix}
-\alpha & 0\\
0 & \alpha\beta
\end{pmatrix}\in \tilde G$ and $t\in \tilde{G}$.
Moreover, since $SL_2(q^{\frac{1}{2}})\subset \tilde{G}$, we have $u=\begin{pmatrix}
0 & 1\\
-1 & 0
\end{pmatrix}\in \tilde{G}$ and, therefore $ut^{-1}u =\begin{pmatrix}
1 & 0\\
0 & -\beta^{-1}
\end{pmatrix}\in \tilde{G}$. It follows that $M=stsut^{-1}u=\begin{pmatrix}
-\alpha & 0\\
0& -\alpha
\end{pmatrix}\in \tilde{G}, t\in \tilde{G}$ and $-t\in \tilde{G}$. It follows that we can use the same arguments as for $\F_q=\F_p(\alpha,\beta)=\F_p(\alpha+\alpha^{-1},\beta+\beta^{-1})$ to get that there is a splitting of $<-\alpha,-\beta>$ in $\tilde{G}$. This proves that $G\simeq (\tilde{G}\cap SL_2(q) \rtimes \Z/\op{lcm}(a,b)\Z$. It now only remains to determine $\tilde{G}\cap SL_2(q)$ depending on $d$.

\bigskip

\textbf{First case} : $d=1$. Since $d=1$, we have $\op{lcm}(a,b)=ab$. It is thus sufficient to show that $SL_2(q^{\frac{1}{2}})\simeq[\tilde{G},\tilde{G}]=\tilde{G}\cap SL_2(q)$. The group $[\tilde{G},\tilde{G}]$ is a normal subgroup of $\tilde{G}\cap SL_2(q)$, therefore we can consider the quotient of those two groups. Let $M\in (\tilde{G}\cap SL_2(\F_q))/[\tilde{G},\tilde{G}]$. Since $a$ is the order of $\tilde{s}$ and $b$ is the order of $\tilde{t}$, we can write $M$ as $\tilde{s}^k\tilde{t}^\ell$ with $0\leq k\leq a-1$ and $0\leq \ell\leq b-1$. The matrix $M$ is of determinant $1$, therefore $(-\alpha)^k(-\beta)^l=1$, which is only possible if $l=k=0$ since $a$ and $b$ are coprime. This proves that the quotient is trivial and, therefore $\tilde{G}\cap SL_2(q)\simeq SL_2(q^{\frac{1}{2}})$. It follows that $G\simeq SL_2(q^{\frac{1}{2}})\rtimes \Z/ab\Z$.

\bigskip

\textbf{Second case} : $d=2$. Since $[\tilde{G},\tilde{G}] \simeq  SL_2(q^{\frac{1}{2}})$, we have that $SL_2(q^{\frac{1}{2}})\triangleleft \tilde{G}\cap SL_2(q)$. Let $M\in \tilde{G}\cap SL_2(q)$ and $A\in SL_2(q^{\frac{1}{2}})$, we then have  $\epsilon(M)A\epsilon(M)^{-1}=\epsilon(MAM^{-1})=MAM^{-1}$ since $\epsilon$ stabilizes $SL_2(q^{\frac{1}{2}})$.

It follows that $\epsilon(M)^{-1}M$ belongs to the centralizer of $SL_2(q^{\frac{1}{2}})$ in $SL_2(q)$ which is equal to $\{I_2,-I_2\}$. It follows that $\tilde{G}\cap SL_2(q) \subset \{M\in SL_2(q), \epsilon(M)^{-1}M\in \{I_2,-I_2\}\}.$ Let us show that this inclusion is in fact an equality.

Let $\varphi : \tilde{G}\cap SL_2(q) \rightarrow \{I_2,-I_2\}$ be the map $M\mapsto \epsilon(M)^{-1}M$. The above inclusion proves that this map is a group morphism. We have $\op{ker}(\varphi)=\{M\in \tilde{G}\cap SL_2(q), \epsilon(M)=M\}\subset SL_2(q^{\frac{1}{2}})$ and $SL_2(q^{\frac{1}{2}})\subset \tilde{G}\cap SL_2(q)$, therefore $\op{ker}(\varphi) \simeq SL_2(q^{\frac{1}{2}})$. To conclude the proof, we only need to show that there exists  $M\in \tilde{G}\cap SL_2(q)$ such that $\epsilon(M)^{-1}M =-I_2$.

We have that $(-\alpha)^{\frac{a}{2}}=(-\beta)^{\frac{b}{2}}=-1$ and $\epsilon(\alpha)=\alpha^{-1}$. Let $M=Ps^{\frac{a}{2}}t^{\frac{b}{2}}P^{-1}$, we have 
\begin{eqnarray*}
M & = &\begin{pmatrix}
 \frac{(-1)^\frac{a}{2}\beta^{\frac{b}{2}+1}\alpha+\beta^\frac{b}{2}\alpha^{\frac{a}{2}+1}+\beta^{\frac{b}{2}+1}\alpha^\frac{a}{2}+\beta^\frac{b}{2}(-1)^\frac{a}{2}}{((\alpha+1)(\beta+1)} & \frac{-(\alpha+\beta)(\alpha-1)(\alpha\beta+1)((-1)^\frac{b}{2}\alpha^\frac{a}{2}-(-1)^{\frac{b}{2}+\frac{a}{2}})}{\alpha(\beta+1)(\alpha\beta+\alpha+\beta+1)}\\
\frac{ \alpha\beta^\frac{b}{2}((-1)^\frac{a}{2}-\alpha^\frac{a}{2})}{(\alpha^2-1)} & \frac{(-1)^\frac{b}{2}\alpha^{\frac{a}{2}+1}\beta+(-1)^{\frac{b}{2}+\frac{a}{2})}\alpha+(-1)^{\frac{b}{2}+\frac{a}{2})}\beta+(-1)^\frac{b}{2}\alpha^\frac{a}{2}}{\alpha\beta+\alpha+\beta+1}
\end{pmatrix}\\
& = &\begin{pmatrix}
 \frac{(-1)^{\frac{a}{2}+\frac{b}{2}}(-\beta\alpha+\alpha+\beta-1)}{\alpha\beta+\alpha+\beta+1} & \frac{(-1)^{\frac{b}{2}+\frac{a}{2}}(2(\alpha+\beta)(\alpha-1)(\alpha\beta+1))}{\alpha(\beta+1)(\alpha\beta+\alpha+\beta+1)}\\
-\frac{(-1)^{\frac{b}{2}+\frac{a}{2}} 2\alpha}{\alpha^2-1} & \frac{(-1)^{\frac{b}{2}+\frac{a}{2}}(-\alpha\beta+\alpha+\beta-1)}{\alpha\beta+\alpha+\beta+1}
\end{pmatrix}\\
& = &(-1)^{\frac{a+b}{2}}\begin{pmatrix}
 \frac{-\beta\alpha+\alpha+\beta+-1}{\alpha\beta+\alpha+\beta+1} & \frac{2(\alpha+\beta)(\alpha-1)(\alpha\beta+1)}{\alpha(\beta+1)(\alpha\beta+\alpha+\beta+1)}\\
-\frac{ 2\alpha}{\alpha^2-1} & \frac{-\alpha\beta+\alpha+\beta-1}{\alpha\beta+\alpha+\beta+1}
\end{pmatrix}\\
\epsilon(M) & = &(-1)^{\frac{a+b}{2}}\begin{pmatrix}
 \frac{-\beta\alpha^{-1}+\alpha^{-1}+\beta-1}{\alpha^{-1}\beta+\alpha^{-1}+\beta+1} & \frac{2(\alpha^{-1}+\beta)(\alpha^{-1}-1)(\alpha^{-1}\beta+1)}{\alpha^{-1}(\beta+1)(\alpha^{-1}\beta+\alpha^{-1}+\beta+1)}\\
-\frac{ 2\alpha^{-1}}{\alpha^{-2}-1} & \frac{-\alpha^{-1}\beta+\alpha^{-1}+\beta-1}{\alpha^{-1}\beta+\alpha^{-1}+\beta+1}
\end{pmatrix}\\
& = &(-1)^{\frac{a+b}{2}}\begin{pmatrix}
 \frac{-\beta+1+\beta\alpha-\alpha}{\beta+1+\alpha\beta+\alpha} & \frac{2(1+\alpha\beta)(1-\alpha)(\beta+\alpha)}{\alpha(\beta+1)(\beta+1+\alpha\beta+\alpha)}\\
-\frac{ 2\alpha}{1-\alpha^2} & \frac{-\beta+1+\alpha\beta-\alpha}{\beta+1+\alpha\beta+\alpha}
\end{pmatrix}\\
& = & -M.
\end{eqnarray*}
It follows that $\epsilon(M)^{-1}M=-I_2$ which concludes the proof.

\end{proof}

\newpage 

$ $

\newpage
\chapter{Type D}\label{TypeD}

In this section, we determine the image of the Artin group of type $D$ inside its associated finite Iwahori-Hecke algebra. The structure of this section is similar to the one in chapter \ref{TypeB} for type $B$. We first define the model and prove some basic properties on the irreducible representations. We then determine the different factorizations appearing, they will be similar to the last cases we had to consider in type $B$. We then state the main results for type $D$ in Theorem \ref{result1D} and Theorem \ref{result2D}. The last two sections prove those theorems using induction. The main differences which will arise come from the fact that we only have to consider one parameter. The branching rule (Lemma \ref{branch}) is quite a bit more complicated and we need to prove some results on orthogonal groups defined over finite fields in a less general setting as we did in type $B$.

\section{Definition of the model}

Let $n\geq 4$, $p$ a prime different from $2$, $\alpha\in \overline{\F_p}$ of order greater than $2n$. We set in this section $\F_q=\F_p(\alpha)$. As in \cite{G-P}, we take for the Iwahori-Hecke algebra of type $D$, $\mathcal{H}_{D_n,\alpha}$ the sub-algebra of $\mathcal{H}_{B_n,\alpha,1}$ generated by $U=TS_1T,S_1,\dots S_{n-1}$. More precisely, we have the following definition.

\begin{Def}
The Iwahori-Hecke algebra of $\mathcal{H}_{D_n,\alpha}$ of type $D$ is the subalgebra of $\mathcal{H}_{B_n,\alpha,1}$ generated by $U=TS_1T,S_1,\dots S_{n-1}$. We then have (see \cite{G-P} 10.4) that $\mathcal{H}_{D_n,\alpha}$ is the algebra generated by the above generators and that they verify the following relations
\begin{enumerate}
\item $U^2=(\alpha-1)U+\alpha$,
\item $\forall i \in [\![1,n-1]\!], S_i^2=(\alpha-1)S_i+\alpha$,
\item $\forall i\in [\![1,n-1]\!]\setminus\{2\}, US_i=S_iU$,
\item $US_2U=S_2US_2$,
\item $\forall i\in [\![1,n-2]\!], S_iS_{i+1}S_i=S_{i+1}S_iS_{i+1}$,
\item $\forall i,j\in [\![1,n-1]\!], \vert i-j\vert > 1, S_iS_j=S_jS_i$.
\end{enumerate}
\end{Def}

We will now give a decomposition into irreducible modules of this algebra. In order to do this, we give the action of $\mathcal{H}_{D_n,\alpha}$ on modules generated by standard double-tableaux associated to double-partitions of $n$.

\begin{prop}
Let $\lambda=(\lambda_1,\lambda_2)\Vdash n$ and $\T\in \lambda$. For $i\in [\![1,n-1]\!]$, we write $m_i(\T)=\frac{\alpha-1}{1-(-1)^{\delta_i}\alpha^{c_i-r_i+r_{i+1}-c_{i+1}}}$, where $i$ (resp $i+1$) is in box $(r_i,c_i)$ (resp $(r_{i+1},c_{i+1})$) of a component of $\T$, $\delta_i =0$ if $i$ and $i+1$ are in the same component and $\delta_i=1$ otherwise.

  The action of the generators on the standard double-tableau $\T$ is then the following
\begin{enumerate}
\item $U.\T =m_1(\T)\T-(1+m_1(\T))\tilde{\T},$ with $\tilde{\T}=0$ if $\T_{1\leftrightarrow 2}$ not standard $\tilde{\T}=\T_{1\leftrightarrow 2}$ otherwise.
\item $\forall i\in [\![1,n-1]\!]$, $S_i.\T=m_i(\T)\T+(1+m_i(\T))\tilde{\T}$ with $\tilde{\T}=0$ if $\T_{i\leftrightarrow i+1}$ is not standard and $\tilde{\T}=\T_{i\leftrightarrow i+1}$ otherwise.
\end{enumerate}

\end{prop} 

We write $V_\lambda$ for the module generated by the standard double-tableaux associated to $\lambda$. By \cite{G-P} (prop 10.4.5), we have that the action on the generators of $\mathcal{H}_{D_n,\alpha}$ commutes with the action of $\sigma$ defined by $\sigma((\T_1,\T_2))=(\T_2,\T_1)$. We then have that for any $\lambda=(\lambda_1,\lambda_2)\Vdash n$, $V_{(\lambda_1,\lambda_2)}$ is isomorphic to $V_{(\lambda_2,\lambda_1)}$ as $\mathcal{H}_{D_n,\alpha}$-module.

 If $\lambda=(\lambda_1,\lambda_1)$, then we can consider a basis $(\T_1,\T_2,\dots,\T_r,\sigma(\T_1),\sigma(\T_2),\dots,\sigma(\T_r))$ of $V_\lambda$. We then have that $V_{\lambda,+}=<\T_i+\sigma(\T_i)>_{i\in [\![1,r]\!]}$ and $V_{\lambda,-}=<\T_i-\sigma(\T_i)>_{i\in [\![1,r]\!]}$ are submodules of $V_\lambda$.
 
 We know by \cite{G-P} (10.4) that the irreducible modules of the Iwahori-Hecke algebra of type $D$ in the generic case are the modules $V_\lambda$ labeled by double-partitions $(\lambda_1,\lambda_2)$ with $\lambda_1< \lambda_2$ and the modules $V_{(\lambda_1,\lambda_1),+}$ and $V_{(\lambda_1,\lambda_1),-}$ for $\lambda_1\vdash \frac{n}{2}$. We will use Proposition \ref{Tits} in order to prove that this is also the case in the finite field setting. The Schur elements are quite complicated to write, therefore we need to introduce some new objects before giving the Schur elements.
 
 \begin{Def}\label{cacahuète}
Let $(\lambda,\mu)\Vdash n$ with $\lambda=(\lambda_1,\dots,\lambda_r,0,\dots)\vdash n_\lambda$ and \\$\mu=(\mu_1,\dots,\mu_m,0,\dots)\vdash n_\mu$.

If $r\leq m$, then we set $X_{\lambda,\mu}=\{\lambda_i+m-i\}_{i\in [\![1,r]\!]}\cup [\![0,m-r-1]\!]$ and $Y_{\lambda,\mu}=\{\mu_i+m-i\}_{i\in [\![1,m]\!]}$.

If $m > r$, then we set $X_{\lambda,\mu}=\{\lambda_i+r-i\}_{i\in [\![1,r]\!]}$ and $Y_{\lambda,\mu}=\{\mu_i+r-i\}_{i\in [\![1,m]\!]}\cup [\![0,r-m]\!]$.

For $X\subset \N$ and $u$ a given parameter, we set $\Delta(X,u)=\underset{(k,l)\in X^2, k>l}\prod (u^k -u^l)$.
\end{Def}

By \cite{G-P} (10.5.7 and 9.3.6), we have the following proposition

\begin{prop}
If $(\lambda,\mu)\Vdash n$ and $\lambda\neq \mu$, then the Schur element $c_{\lambda,\mu}$ associated to the irreducible $\mathcal{H}_{D_n,\Z[u,u^{-1}]}$-module  $V_{\lambda,\mu}$ is 
$$c_{\lambda,\mu}=\frac{2^{b-1}u^\frac{b(b-1)(4b-5)}{6}\underset{k\in X_{\lambda,\mu}}\prod \underset{h=1}{\overset{k}\prod}(u^{2h}-1)\underset{l\in Y_{\lambda,\mu}}\prod \underset{h=1}{\overset{l}\prod} (u^{2h}-1)}{(u-1)^n\Delta(X_{\lambda,\mu},u)\Delta(Y_{\lambda,\mu},u)\underset{(k,l)\in X_{\lambda,\mu}\times Y_{\lambda,\mu}}\prod (u^k+u^l)}$$ 
where $b=\vert X_{\lambda,\mu}\vert =\vert Y_{\lambda,\mu}\vert $.

If $(\lambda,\lambda)\Vdash n$, then the Schur elements associated to $V_{(\lambda,\lambda),+}$ and $V_{(\lambda,\lambda),-}$ are equal and are given by
$$c_{\lambda,\lambda}=\frac{2^{b}u^\frac{b(b-1)(4b-5)}{6}\underset{k\in X_{\lambda,\lambda}}\prod \underset{h=1}{\overset{k}\prod}(u^{2h}-1)\underset{l\in Y_{\lambda,\lambda}}\prod \underset{h=1}{\overset{l}\prod} (u^{2h}-1)}{(u-1)^n\Delta(X_{\lambda,\lambda},u)\Delta(Y_{\lambda,\lambda},u)\underset{(k,l)\in X_{\lambda,\lambda}\times Y_{\lambda,\lambda}}\prod (u^k+u^l)}.$$ 
\end{prop}

We can now state the theorem giving the semi-simple decomposition of $\mathcal{H}_{D_n,\alpha}$.

\begin{theo}\label{Schur}
Assume the order of $\alpha$ is greater than $2n$. We then have that $\mathcal{H}_{D_n,\alpha}$ is split semi-simple and its pairwise non-isomorphic irreducibles modules are $V_{\lambda,\mu}$ (with $\lambda> \mu$), $V_{\lambda,\lambda,+}$ and $V_{\lambda,\lambda,-}$.

In the following, we write $V$ for the direct sum of these irreducible modules.
\end{theo}
 
 \begin{proof}
 Let $A=\Z[u^{\pm 1}]$, $\theta:\Z[u^{\pm 1}]\rightarrow \F_q$ defined by $\theta
(k)=\overline{k}$ for $k\in \Z$ and $\theta(u)=\alpha$. By Proposition \ref{Tits}, it is sufficient to show that  $c_{\lambda,\mu}\in B$ with $B$ as in Proposition \ref{Tits}, that $\tilde{\theta}(c_{\lambda,\mu})\neq 0$ and $(1+m_i(\T))m_i(\T)\neq 0$ and that they are well-defined for $(\lambda,\mu)\Vdash n$, for $i\in [\![1,n-1]\!]$ and for $\T\in (\lambda,\mu)$ such that $\T_{i\leftrightarrow i+1}$ is also standard. The condition on the $m_i$ implies that for any pair $(i,j)$, there exists a matrix $M$ in the representation associated to $V_{\lambda,\mu}$ such that $M_{i,j}\neq 0$ since there exists a path between any pair of standard double-tableaux such that any standard double-tableau in the path is obtained by transposing a pair $(r,r+1)$.

\smallskip

Let $(\lambda,\mu)\Vdash n$ as in Definition \ref{cacahuète}, $i\in [\![1,n-1]\!]$ and $\T\in (\lambda,\mu)$.

\bigskip

We first show that $c_{\lambda,\mu}\in B$, i.e. $\theta\left((u-1)^n\Delta(X_{\lambda,\mu})\Delta(Y_{\lambda,\mu})\underset{(k,\ell)\in X_{\lambda,\mu}\times Y_{\lambda,\mu}}\prod (u^k+u^\ell)\right)\neq 0$. We have that $\theta((u-1)^n)=(\alpha-1)^n\neq 0$. 

We have $\Delta(X_{\lambda,\mu},u)=\underset{(k,\ell)\in X_{\lambda,\mu}^2, k>\ell}\prod (\alpha^k -\alpha^\ell)$. Let $(k,\ell)\in X_{\lambda,\mu}^2$ such that $k> \ell$. We have $0\leq \ell < k \leq \lambda_1+\max(m,r)-1$ and, therefore $k-\ell\leq \lambda_1+\max(m,r)-1$.
 Since $\lambda_1+r-1$ is equal to the number of boxes in the first row of $\lambda$, we have that $\lambda_1+r-1 \leq n_\lambda$. Since $m$ is the number of boxes in the first column of $\mu$, we have that $m\leq n_\mu$, therefore $\lambda_1+\max(m,r)-1=\lambda_1+\max(r-1+m-r,r-1) \leq n_\lambda+n_\mu \leq n < 2n$. The equality $\alpha^k-\alpha^\ell=0$ implies that $\alpha^{k-\ell}=1$. This quantity is therefore non-zero because the order of $\alpha$ is greater than $2n$.
 
In the same way, we have $\Delta(Y_{\lambda,\mu},u)\neq 0.$
 
Let $(k,\ell)\in X_{\lambda,\mu}\times Y_{\lambda,\mu}$. Assume by contradiction that $\theta(u^k+u^\ell)=0$. We then have $\alpha^k+\alpha^\ell=0$, therefore $\alpha^{k-\ell}=-1$ and $\alpha^{2(k-\ell)}=1$. We have $0 \leq k\leq \lambda_1+\max(r,m)-1$ and $0\leq \ell \leq \mu_1+\max(r,m)-1$. Therefore $-2n\leq 1-\max(m,r)-\mu_1 \leq 2(k-\ell)\leq 2(\lambda_1+\max(m,r)-1)\leq 2n$ and we have $\alpha^{k-\ell}\neq 1$ if $k\neq \ell$ since the order of $\alpha$ is greater than to $2n$. If $k=\ell$, then $\theta(u^k+u^\ell)=2\alpha^k\neq 0$ because $p\neq 2$.

It follows that $c_{\lambda,\mu}\in B$.
 
 \bigskip
 
We now show that $\tilde{\theta}(c_{\lambda,\mu})\neq 0$, i.e.
  $$\theta\left(2^{b}u^\frac{b(b-1)(4b-5)}{6}\underset{k\in X_{\lambda,\mu}}\prod \underset{h=1}{\overset{k}\prod}(u^{2h}-1)\underset{l\in Y_{\lambda,\mu}}\prod \underset{h=1}{\overset{l}\prod} (u^{2h}-1)\right)\neq 0.$$ We have $\theta(2^bu^\frac{b(b-1)(4b-5)}{6})=2^b\alpha^\frac{b(b-1)(4b-5)}{6}\neq 0$ because $p\neq 2$. Let $k\in \N$ such that $k\in X_{\lambda,\mu}$ or ($k\in Y_{\lambda,\mu}$ and $0\leq h\leq k$). We have shown that $k\leq n$, therefore $\theta(u^{2h}-1))= \alpha^{2h}-1\neq 0$ because $\alpha$ is of order greater than $2n$. It follows that $\tilde{\theta}(c_{\lambda,\mu})\neq 0$.
  
  \bigskip
  
  We now show that $m_i(\T)$ is well-defined and non-zero.
  
   We have that $m_i(\T)=\frac{\alpha-1}{1-(-1)^{\delta_i}\alpha^{c_i-r_i+r_{i+1}-c_{i+1}}}\neq 0$. It is sufficient to show that \\
   $1-(-1)^{\delta_i}\alpha^{c_i-r_i+r_{i+1}-c_{i+1}}\neq 0$ or $\alpha^{c_i-r_i+r_{i+1}-c_{i+1}}\neq (-1)^{\delta_i}$.
  
  If $i$ and $i+1$ are in the same component of $\T$, then $\delta_i=0$ and $\vert c_i-r_i+r_{i+1}-c_{i+1}\vert =\vert c_i-c_{i+1}+r_{i+1}-r_i\vert$ is less than the minimal number of boxes on a path within the Young diagram from the box where $i$ is and the box where $i+1$ is. It follows that $\alpha^{c_i-r_i+r_{i+1}-c_{i+1}}\neq 1$ because $\alpha$ is of order greater than $2 n > \max(n_\lambda,n_\mu) > \vert c_i-r_i+r_{i+1}-c_{i+1}\vert$ and $c_i-r_i\neq c_{i+1}-r_{i+1}$ because $\T$ is standard.
  
  If $i$ and $i+1$ are in distinct components then we have that $\delta_i=1$. It is then sufficient to show that  $-2n \leq c_i-r_i+r_{i+1}-c_{i+1} \leq 2n$ because $p\neq 2$. If $i$ is in the left tableau and $i+1$ is in the right tableau then we have $1-n\leq 1-r\leq c_i-r_i \leq \lambda_1-1\leq n-1$ and
  $1-n \leq 1-m \leq c_{i+1} -r_{i+1} \leq \mu_1-1 \leq n-1$. It follows that $2-2n \leq c_i-r_i+r_{i+1}-c_{i+1} \leq 2n-2$ which proves that $m_i(\T)$ is well-defined since $\alpha$ is of order greater than $2n$. The same reasoning shows that $m_i(\T)$ is well-defined if $i$ is in the right tableau and $i+1$ is in the left tableau.
  
  \bigskip
  
 Finally, we show that $1+m_i(\T) \neq 0$. We have
  $$1+m_i(\T) =1+\frac{\alpha-1}{1-(-1)^{\delta_i}\alpha^{c_i-r_i+r_{i+1}-c_{i+1}}}=\frac{\alpha(1-(-1)^{\delta_i}\alpha^{c_i-r_i+r_{i+1}-c_{i+1}-1})}{1-(-1)^{\delta_i}\alpha^{c_i-r_i+r_{i+1}-c_{i+1}}}.$$
  We have shown in the previous step that $-2n+2\leq c_i-r_i+r_{i+1}-c_{i+1}\leq 2n-2$, therefore $-2n+1 \leq c_i-r_i+r_{i+1}-c_{i+1}\leq 2n-3$. Since $\T_{i\leftrightarrow i+1}$ is standard if $i$ and $i+1$ are in the same component we have in that case that $\vert c_i-r_i+r_{i+1}-c_{i+1} \vert -1\geq 2$. This shows that $1+m_i(\T) \neq 0$ and concludes the proof.
  \end{proof}

  The branching rule for type $D$ is more complicated, therefore we recall it in the following proposition (a proof in a more general setting can be found in \cite{Branch}).

\begin{lemme}\label{branch}
Let $n\geq 5$ and $(\lambda,\mu)\Vdash n, \lambda>\mu$. We then have:
\begin{enumerate}
\item If $n_\lambda > n_\mu+1$, then $V_{\lambda,\mu|\mathcal{H}_{D_{n-1},\alpha}}=\underset{(\tilde{\lambda},\tilde{\mu})\subset (\lambda,\mu)}\bigoplus V_{\tilde{\lambda},\tilde{\mu}}$.
\item If $n_\lambda=n_\mu+1$ and $\mu \not\subset \lambda$, then
$$V_{\lambda,\mu|\mathcal{H}_{D_{n-1},\alpha}}=(\underset{\tilde{\mu}\subset \mu}\bigoplus V_{\lambda,\tilde{\mu}}) \oplus (\underset{\tilde{\lambda}> \mu}{\underset{\tilde{\lambda}\subset \lambda}\bigoplus}V_{\tilde{\lambda},\mu})\oplus( \underset{\tilde{\lambda}< \mu}{\underset{\tilde{\lambda}\subset \lambda}\bigoplus}V_{\mu,\tilde{\lambda}}).$$
\item If $n_\lambda=n_\mu+1$ and $\mu \subset \lambda$, then
$$V_{\lambda,\mu|\mathcal{H}_{D_{n-1},\alpha}}=(\underset{\tilde{\mu}\subset \mu}\bigoplus  V_{\lambda,\tilde{\mu}}) \oplus (\underset{\tilde{\lambda}> \mu}{\underset{\tilde{\lambda}\subset \lambda}\bigoplus}V_{\tilde{\lambda},\mu})\oplus (\underset{\tilde{\lambda}< \mu}{\underset{\tilde{\lambda}\subset \lambda}\bigoplus}V_{\mu,\tilde{\lambda}})\oplus V_{\mu,\mu,+}\oplus V_{\mu,\mu,-}.$$
\item If $n_\lambda=n_\mu$ and $\lambda>\mu$, then $V_{\lambda,\mu|\mathcal{H}_{D_{n-1},\alpha}}=(\underset{\tilde{\mu}\subset \mu}\bigoplus V_{\lambda,\tilde{\mu}})\oplus (\underset{\tilde{\lambda}\subset \lambda}\bigoplus V_{\mu,\tilde{\lambda}}).$
\item If $\lambda=\mu$, then $V_{\lambda,\lambda,+|\mathcal{H}_{D_{n-1},\alpha}}=V_{\lambda,\lambda,-|\mathcal{H}_{D_{n-1},\alpha}}=\underset{\tilde{\mu} \subset \mu}\bigoplus V_{\lambda,\tilde{\mu}}.$
\end{enumerate}
\end{lemme}

\begin{proof}

Assume first that $n_\lambda > n_\mu+1$. We then have $V_{\lambda,\mu}=\underset{\T\in (\lambda,\mu)}\bigoplus \F_q \T$. Let $r$ be the number of extremal boxes in the Young double-diagram associated to $(\lambda,\mu)$ (a box is said to be extremal if there exists a standard double-tableau containing $n$ in that box). We write $(r_\ell,c_\ell,\delta_\ell)_{l\in [\![1,r]\!]}$ the extremal box, where $r_\ell$ is the row and $c_\ell$ is the column of the box in the component the box belongs to and $\delta_\ell$ indicates which component the box belongs to. We then have $V_{\lambda,\mu}=\underset{l=1}{\overset{r} \bigoplus}\underset{j\in I_l}\bigoplus \T_{j,l}$ for $\T_{j,l}, j\in I_l$ standard tableaux associated to $(\lambda,\mu)$ such that $n$ is in box $(r_l,c_l,\delta_l)$.  

We can then define a bijection from the basis of $V_{\lambda,\mu}$ to the standard double-tableaux basis of $\underset{(\tilde{\lambda},\tilde{\mu})\subset (\lambda,\mu)}\bigoplus V_{\tilde{\lambda},\tilde{\mu}}$  by mapping a standard double-tableau $\T_{j,l}$ to the standard double-tableau $\T_{j,l}\setminus \{(r_l,c_l,\delta_l)\}$ and we have $\tilde{\lambda}>\tilde{\mu}$ for any 
$(\tilde{\lambda},\tilde{\mu}) \subset (\lambda,\mu)$ since $n_{\tilde{\lambda}}>n_{\tilde{\mu}}$. By construction, this bijection commutes with the action of $\mathcal{H}_{D_{n-1}}$ because for any $i\in [\![1,n-2]\!]$, the position of $i$ and $i+1$ is unchanged.

\bigskip

Assume now that $n_\lambda=n_\mu+1$ and $\mu \not\subset \lambda$. We can then apply a similar reasoning except that we can have $\tilde{\lambda}<\mu$ for some $\tilde{\lambda}\subset \lambda$. We then map the standard tableau $\T_{j,l}$ to $\sigma(\T_{j,l}\setminus \{(r_l,c_l,\delta_l)\})$. The action of $\mathcal{H}_{D_{n-1},\alpha}$ then still commutes with the bijection since it commutes with action of $\sigma$. 

\bigskip

Assume now that $n_\lambda=n_\mu+1$ et $\mu \subset \lambda$. We keep the same bijection as in the previous case except for the extremal box $(r_{\ell_0},c_{\ell_0},\delta_{\ell_0})$, which, when removed from the Yound double-diagram $(\lambda,\mu)$, affords the Young double diagram associated to $(\mu,\mu)$. We then map the tableau $\T_{j,\ell_0}=(\T_{j,\ell_0,1},\T_{j,\ell_0,2})$ to $\T_{j,\ell_0}\setminus\{(r_{\ell_0},c_{\ell_0},\delta_{\ell_0})\} +\sigma(\T_{j,\ell_0}\setminus\{(r_{l_0},c_{l_0},\delta_{l_0})\})$ if $\tau_{\T_{j,\ell_0}}(1)=1$ and to $\T_{j,\ell_0}\setminus\{(r_{\ell_0},c_{\ell_0},\delta_{\ell_0})\} -\sigma(\T_{j,\ell_0}\setminus\{(r_{\ell_0},c_{\ell_0},\delta_{\ell_0})\})$ otherwise. The action of $\mathcal{A}_{D_n,\alpha}$ then again commutes with the bijection because it commutes with $\sigma$.

\bigskip

Assume now that $n_\lambda=n_\mu$ and $\lambda > \mu$. The bijection is defined in the same way as before except when $\tilde{\lambda}\subset \lambda$, where we apply $\sigma$.

\bigskip

Assume now that $\lambda=\mu$. We can then number the standard double-tableaux associated to $(\lambda,\lambda)$ by $(\T_1,\T_2,\dots,\T_m,\sigma(\T_1),\sigma(\T_2),\dots,\sigma(T_m))$ such that $n$ is in left component of the $m$ first ones and in the right component of the last $m$ ones. We then have $V_{\lambda,\lambda,+}=\underset{j=1}{\overset{r}\bigoplus}\F_q(\T_j+\sigma(\T_j))$ and $V_{\lambda,\lambda,-}=\underset{j=1}{\overset{r}\bigoplus}\F_q(\T_j-\sigma(\T_j))$.

For $V_{\lambda,\lambda,+}$, we map $\T_j+\sigma(T_j)$ to the standard double-tableau obtained by removing the box of $\T_j$ containing $n$. We have to check that the action of $\mathcal{H}_{D_{n-1},\alpha}$ on $\T_j+\sigma(\T_j)$ is the same as the one on this tableau. We have $m_i(\T_j)=m_i(\sigma(\T_j)), \sigma(T_j)_{i\leftrightarrow i+1} =\sigma(T_j)_{i\leftrightarrow i+1}$ and $\T_{j,i\leftrightarrow i+1}\notin \{T_j,\sigma(\T_j)\}$, therefore the action is indeed identical.

For $V_{\lambda,\lambda,-}$, we map in the same way $\T_j-\sigma(\T_j)$ to the standard double-tableau obtained by removing the box of $\T_j$ containing $n$. In order to check that the action is the same, we can use the same arguments as for $V_{\lambda,\lambda,+}$ and check that $\T_{j,i\leftrightarrow i+1} \in \{\T_k\}_{k\in [\![1,m]\!]}$. This is true since we chose a numbering such that $n$ is in the left component only for the $m$ first tableaux and $n$ stays in the same box after permutation of $i$ and $i+1$ for any $i\leq n-2$.

\end{proof}

\textbf{Remark} : The two submodules of $V_{\lambda,\lambda}$ are not isomorphic, by the fact that $n$ goes from the left component to the right one or from the right component to the left one after applying $S_{n-1}$.

We keep the same weight on double-tableaux as for type B. Let $\lambda=(\lambda_1,\lambda_2)\Vdash n$ and $\T=(\T_1,\T_2)\in \lambda$. We define $\varphi(\T)$ to be $\T'$ if $\mu'>\lambda'$ and $\sigma(\T')$ otherwise. We define a new $\tilde{\nu}(\lambda)$ to be $\nu(\lambda_1)\nu(\lambda_2)(-1)^{n_{\lambda_1}(n-n_{\lambda_1})}$ if $\lambda_2' \geq \lambda_1'$ and $\tilde{\nu}(\lambda)=\nu(\lambda_1)\nu(\lambda_2)$ otherwise. We define the bilinear form $(\T|\tilde{\T})=\omega(\T)\delta_{\varphi(\T),\tilde{\T}}$.

\begin{prop}\label{bilin2}
For any pair of standard double-tableaux $(\T, \tilde{\T})$, we have the following properties.
\begin{enumerate}
\item $(S_i.\T|S_i.\tilde{\T}) = (-\alpha)(\T|\tilde{\T})$ and $(U.\T|U.\tilde{\T}) = (-\alpha)(\T|\tilde{\T})$.

\item
For all $d\in \mathcal{A}_{D_n}$, we have that $(d.\T|d.\tilde{\T}) = (\T|\tilde{\T})$.

Those relations stay true if we substitute one or two of the standard double-tableaux by the elements $\sigma(\T)-\T$ and $\sigma(\T)+\T$, which form bases for $V_{\lambda,+}$ and $V_{\lambda,-}$ for double-partitions $\lambda$ of the form $\lambda=(\lambda_1,\lambda_1)$.
\item
The restriction of $(.,.)$ to $V_\lambda$ if $\lambda =\varphi(\lambda)\neq (\lambda_2,\lambda_1)$ and to $V_\lambda\oplus V_{\varphi(\lambda)}$ if $\lambda\notin \{\varphi(\lambda),(\lambda_2,\lambda_1)\}$ is non-degenerate with $\varphi(\lambda)=\lambda'$ if $n_\lambda=n_\mu$ and $\mu'> \lambda'$, and $\varphi(\lambda)=(\lambda_1',\lambda_2')$ otherwise.

If $\lambda=\varphi(\lambda)\neq (\lambda_2,\lambda_1)$ then $(.,.)$ is symmetric on $V_\lambda$ if $\tilde{\nu}(\lambda)=1$ and skew-symmetric otherwise.

Moreover, its Witt index is positive.
\item If $n\equiv 0 ~(\bmod ~ 4)$ and $\lambda=(\lambda_1,\lambda_1)$, then the restriction of $(.,.)$ to $V_{\lambda,+}$ and $V_{\lambda,-}$ if $\lambda=\lambda'$ and to $V_{\lambda,+}\oplus V_{\lambda,-}$ if $\lambda\neq \lambda'$, is non-degenerate. 

If $\lambda=\lambda'$ then $(.,.)$ is symmetric on $V_{\lambda,+}$ and $V_{\lambda,-}$ if $\tilde{\mu}(\lambda)=1$ and skew-symmetric otherwise. 
Moreover, its Witt index is positive.

\item If $n\equiv 2 ~(\bmod ~ 4)$ and $\lambda=(\lambda_1,\lambda_1)$ then the restriction of $(.,.)$ to $V_{\lambda,+}\oplus V_{\lambda',-}$ is non-degenerate.
\end{enumerate}
\end{prop}

\begin{proof}

For $\mathbf{1.}$ and $\mathbf{2.}$, we have  $m_i(\sigma(\T)) =m_i(\T)$, therefore the same proof as for Proposition \ref{bilin} applies. The extension to elements of the bases of $V_{\lambda,+}$ and $V_{\lambda,-}$ follows from the bilinearity of $(.,.)$.

For $\mathbf{3.}$, the same proof also applies because $\tilde{\nu}(\lambda)=\omega(\T)\omega(\varphi(\T))$. This is true because when $\varphi(\T)=\T'$, $\tilde{\nu}(\lambda)$ does not change from the one in type $B$ and when $\varphi(\T)=\sigma(\T')$, $\tilde{\nu}(\lambda)$ is multiplied by $(-1)^{n_{\lambda_1}(n-n_{\lambda_1})}=\omega(\T)\omega(\sigma(\T))$.

$\mathbf{4.}$ We assume $n \equiv 0 ~(\bmod ~ 4)$. If $\lambda=(\lambda_1,\lambda_1)\Vdash n$ and $\T\in \lambda$ then 
$$\omega(\sigma(\T))=(-1)^{n_{\lambda_1}(n-n_{\lambda_1)})}\omega(\T)=(-1)^{(\frac{n}{2})^2}\omega(\T)=\omega(\T).$$

For any standard double-tableaux $\T,\tilde{\T}$, we have that $(\T|\tilde{\T})=\omega(\T)\delta_{\T,\varphi(\T)}$. Since $\lambda=(\lambda_1,\lambda_1)$, we have $\varphi(\lambda)=\lambda'$ and for all $\T\in \lambda$ and all $\T\in \lambda'$, we have $\varphi(\T)=\T'$.

Let $\lambda=(\lambda_1,\lambda_1)$ and $\tilde{\lambda}=(\tilde{\lambda}_1,\tilde{\lambda}_1)$ be double-partitions of $n$. If $\T\in \lambda$ and $\tilde{\T}\in \tilde{\lambda}$, then we have 

\begin{eqnarray*}
(\T+\sigma(\T)|\tilde{\T}+\sigma(\tilde{\T})) & = & (\T|\tilde{\T}) + (\T|\sigma(\tilde{\T}))+ (\sigma(\T)|\tilde{\T})) + (\sigma(\T)|\sigma(\tilde{\T}))  \\
 & = & \omega(\T)(\delta_{\T,\tilde{\T}'}+\delta_{\T,\sigma(\tilde{\T})'})+\omega(\sigma(\T))(\delta_{\sigma(\T),\tilde{\T}'}+\delta_{\sigma(\T),\sigma(\tilde{\T})'})\\
 & = & (\delta_{\T,\tilde{\T}'}+\delta_{\T,\sigma(\tilde{\T})'})(\omega(\T)+\omega(\sigma(\T))\\
 & = & \delta_{\T+\sigma(\T),\tilde{\T}'+\sigma(\tilde{\T}')}(\omega(\T)+\omega(\sigma(\T)))\\
 & = & 2\omega(\T)\delta_{\T'+\sigma(\T)',\tilde{\T}+\sigma(\tilde{\T})}.
\end{eqnarray*}
In the same way, we have that $(\T+\sigma(\T)|\tilde{\T}-\sigma(\tilde{\T}))=(\T-\sigma(\T)|\tilde{\T}+\sigma(\tilde{\T}))=0$ and\\ $(\T-\sigma(\T)|\tilde{\T}-\sigma(\tilde{\T}))=2\omega(\T)\delta_{\T'-\sigma(\T)',\tilde{\T}-\sigma(\tilde{\T})}$. The result follows.

$\mathbf{5.}$ Assume $n\equiv 2 ~(\bmod ~ 4)$.

 If $\lambda=(\lambda_1,\lambda_1)\Vdash n$ and $\T\in \lambda$, then $\omega(\sigma(\T))=(-1)^{n_{\lambda_1}(n-n_{\lambda_1)})}\omega(\T)=(-1)^{(\frac{n}{2})^2}\omega(\T)=-\omega(\T)$. It follows that if $\tilde{\lambda}=(\tilde{\lambda}_1,\tilde{\lambda}_1)\Vdash n$ and $\tilde{\T}\in \tilde{\lambda}$, then $(\T+\sigma(\T)|\tilde{\T}+\sigma(\tilde{\T}))=(\T-\sigma(\T)|\tilde{\T}-\sigma(\tilde{\T}))=0$, $(\T+\sigma(\T)|\tilde{\T}-\sigma(\tilde{\T}))=2\omega(\T)\delta_{\T'-\sigma(\T)',\tilde{\T}-\sigma(\tilde{\T})}$ and $(\T-\sigma(\T)|\tilde{\T}+\sigma(\tilde{\T}))=2\omega(\T)\delta_{\T'+\sigma(\T)',\tilde{\T}+\sigma(\tilde{\T})}$. The result follows. 
\end{proof}

\section{Factorization of the image of the Artin group inside the finite Hecke algebra}

In this section, we find the different factorizations between the irreducible representations of $\mathcal{A}_{D_n}$. Most of the factorization results are summarized in Proposition \ref{isomorphisme2}. We then state the main results for type $D$ in Theorems \ref{result1D} and \ref{result2D}.

We define the linear map $\mathcal{L}$ from $V$ to $V$ which sends $\T$ to $\mathcal{L}(\T)=\omega(\T)\varphi(\T)$.

\begin{prop}\label{transpose2}
Let $r\in [\![1,n-1]\!]$ and $\T$ a standard double-tableau, we then have
$$\mathcal{L}S_r\mathcal{L}^{-1}(\T)=(-\alpha){}^t\!(S_r^{-1})(\T), \mathcal{L}U \mathcal{L}^{-1}=(-\alpha){}^t\!(U^{-1})(\T).$$
Let $\lambda=(\lambda_1,\lambda_2)\Vdash n$. We have the following propositions.
\begin{enumerate}
\item If $\lambda\notin\{\varphi(\lambda),(\lambda_2,\lambda_1)\}$, then $\mathcal{L}$ stabilizes $V_{\lambda}\oplus V_{\varphi(\lambda)}$ and switches $V_{\lambda}$ and $V_{\varphi(\lambda)}$.
\item If $\lambda=\varphi(\lambda)\neq (\lambda_2,\lambda_1)$, then $\mathcal{L}$ stabilizes $V_\lambda$.
\item If $n \equiv 0~ (\bmod ~ 4)$ and $\lambda=(\lambda_1,\lambda_1)\neq (\lambda_1',\lambda_1')$, then $\mathcal{L}$ stabilizes $V_{\lambda,+}\oplus V_{\lambda',+}$ (resp $V_{\lambda,-} \oplus V_{\lambda',-}$) and switches $V_{\lambda,+}$ and $V_{\lambda',+}$ (resp $V_{\lambda,-}$ and $V_{\lambda',-}$).
\item If $n \equiv 0~ (\bmod ~4)$ and $\lambda=(\lambda_1,\lambda_1)=(\lambda_1',\lambda_1')$, then $\mathcal{L}$ stabilizes $V_{\lambda,+}$ and $V_{\lambda,-}$.
\item If $n\equiv 2~ (\bmod ~4)$ and $\lambda=(\lambda_1,\lambda_1)$, then $\mathcal{L}$ stabilizes $V_{\lambda,+}\oplus V_{\lambda',-}$ and switches $V_{\lambda,+}$ and $V_{\lambda',-}$.
\end{enumerate}
\end{prop}

\begin{proof}
This follows directly from Proposition \ref{bilin2} by writing the matrix of the bilinear form and the matrix of $\mathcal{L}$.
 \end{proof}

\begin{prop}\label{chaise}
For $r\in [\![1,n-1]\!]$, we write $\lambda_r$ the double-partition of $n$ defined by $\lambda_r=([r],[1^{n-r}])$ if $r\geq \frac{n}{2}$ and $\lambda_r=([1^{n-r}],r)$ if $r< \frac{n}{2}$. Taking the same notations as in Proposition \ref{patate} of type $B$, for all $d\in A_{D_n}$, we have $\eta_{1,r}(d)=\eta_{2,r}(d)$ because the length in $T$ of such an element is even and we have $\beta=1$. We define the character $\eta_r$ of $A_{D_n}$ by $\eta_r(d)=\eta_{1,r}(d)=\eta_{2,r}(d)$.

We then have by Proposition \ref{patate} that $R_{\lambda_r}\simeq (\Lambda^rR_{\lambda_1})\otimes \eta_{1,r}$ for all $r\in [\![1,n-1]\!]$.
\end{prop}

Assume $\F_q=\F_p(\alpha)\neq \F_p(\alpha+\alpha^{-1})$. We write $\epsilon$ the unique automorphism of $\F_q$ of order $2$. We have that $\epsilon(\alpha)=\alpha^{-1}$. We then define for every standard double-tableau in the same way as for type $B$, 
$$d(\T)=\tilde{d}(\T_1)\tilde{d}(\T_2)\underset{i<j}{\underset{i\in \T_1,j\in \T_2}\prod}\frac{2+\alpha^{a_{i,j}-1}+\alpha^{1-a_{i,j}}}{\alpha+\alpha^{-1}+\alpha^{a_{i,j}}+\alpha^{-a_{i,j}}}$$
 and the associated hermitian form $\langle .,.\rangle$ defined by \break$\langle \T,\tilde{\T}\rangle = d(\T)\delta_{\T,\tilde{\T}}$. We write $\Lambda$ for the set of all double-partitions $\lambda=(\lambda_1,\lambda_2)$ of $n$ such that $\lambda_1\geq \lambda_2$. 

\begin{prop}\label{unitary2}
For all $d\in A_{D_n}$, $\lambda\in \Lambda$ and $\T\in \lambda$, we have $\langle d.\T,d.\tilde{\T}\rangle=\langle \T,\tilde{\T}\rangle$. This shows that $A_{D_n}$ acts in a unitary way on those irreducible modules.  
\end{prop}

\begin{proof}
The proof of the first statement follows from Proposition \ref{unitary} and the second follows from the expression of the bases of $V_{\lambda,\pm}$ and the $\Z$-bilinearity of the hermitian form. 
\end{proof}

We now prove two lemmas which will allow us to restrict ourselves to the derived subgroup $\mathcal{A}_{D_n}$ of $A_{D_n}$.

\begin{lemme}\label{LincolnD}
If $\lambda$ is a double-partition of $n$ then the restriction of  $R_{\lambda}$ to $\mathcal{A}_{D_n}$ is absolutely irreducible.
\end{lemme}

\begin{proof}
Assume first it is true for $n=2$. Since $A_{D_n}$ is generated by $A_{D_{n-1}}$ and $\mathcal{A}_{D_n}$, we have the result for $n\geq 3$ by the same method as in the Lemma 3.4(i) of \cite{BMM}. 

We now show the result is true for $n=2$. We only have to show it for $([1],[1])$ since the other representations are $1$-dimensional. We will show in Section \ref{surjectivitism} (Lemmas \ref{platypus} and \ref{platypus456}) that $R_{[1],[1]}(\mathcal{A}_{D_2})\simeq SL_2(q')$ for some $q'$. The irreducibility then follows.
\end{proof}

We now show a lemma computing the normal closure of $\mathcal{A}_{D_n}$. This is a generalization to type $D$ of Lemma \ref{normclosBn}

\begin{lemme}\label{normclosDn}
For $n\geq 4$, the normal closure $\ll \mathcal{A}_{D_{n-1}}\gg_{\mathcal{A}_{D_n}}$ of $\mathcal{A}_{D_{n-1}}$ in $\mathcal{A}_{D_n}$ is $\mathcal{A}_{D_n}$.
\end{lemme}

\begin{proof}
Let $n\geq 4$. By Lemma $2.1$ of \cite{BMM}, we have that $\mathcal{A}_{A_n}=\ll \mathcal{A}_{A_{n-1}}\gg_{\mathcal{A}_{A_n}}$, where \\$A_{A_n}=<S_1,S_2,\dots, S_{n-1}>\leq A_{D_n}$. We have that $\mathcal{A}_{D_n}$ is generated by $\mathcal{A}_{A_n}$ and $\mathcal{A}_{D_{n-1}}$ therefore the result follows.
\end{proof}

The following proposition summarizes the results in this section :

\begin{prop}\label{isomorphisme2}
Let $\lambda$, $\mu$, $\gamma$ and $\delta$ be doubles-partitions of $n$ such that $\dim(V_\lambda)>1$, $\lambda_1>\lambda_2, \mu_1>\mu_2$, $\gamma_1=\gamma_2$ and $\delta_1=\delta_2$. We have the following properties.
\begin{enumerate}
\item The restrictions of $R_\lambda$, $R_{\gamma,+}$ and $R_{\gamma,-}$ to $\mathcal{A}_{D_n}$ are absolutely irreducible.
\item $R_{\lambda|\mathcal{A}_{D_n}}\simeq R_{\mu|\mathcal{A}_{D_n}} \Leftrightarrow \lambda=\mu$.
\item $R_{\lambda|\mathcal{A}_{D_n}}\not\simeq R_{\gamma,\pm|\mathcal{A}_{D_n}}$.
\item $R_{\gamma,\pm|\mathcal{A}_{D_n}} \simeq R_{\delta,\pm|\mathcal{A}_{D_n}} \Leftrightarrow \gamma= \delta$.
\item $R_{\gamma,\pm|\mathcal{A}_{D_n}} \not\simeq R_{\delta,\mp|\mathcal{A}_{D_n}}$.
\item $R_{\lambda|\mathcal{A}_{D_n}}\simeq R_{\mu|\mathcal{A}_{D_n}}^\star \Leftrightarrow \mu =\varphi(\lambda)$.
\item $R_{\lambda|\mathcal{A}_{D_n}}\not\simeq R_{\gamma,\pm}^\star$.
\item If $n \equiv 0~(\bmod ~ 4)$, then 
\begin{enumerate}
\item $R_{\gamma,\pm|\mathcal{A}_{D_n}} \simeq R_{\delta,\pm|\mathcal{A}_{D_n}}^\star \Leftrightarrow \gamma= \varphi(\delta)$.
\item $R_{\gamma,\pm|\mathcal{A}_{D_n}} \not\simeq R_{\delta,\mp|\mathcal{A}_{D_n}}^\star$.
\end{enumerate}
\item If $n\equiv 2~(\bmod ~ 4)$, then
\begin{enumerate}
\item $R_{\gamma,\pm|\mathcal{A}_{D_n}} \not\simeq R_{\delta,\pm|\mathcal{A}_{D_n}}^\star$.
\item $R_{\gamma,\pm|\mathcal{A}_{D_n}} \simeq R_{\delta,\mp|\mathcal{A}_{D_n}}^\star \Leftrightarrow \gamma=\varphi(\delta)$.
\end{enumerate}
\item If $\F_q=\F_p(\alpha)\neq \F_p(\alpha+\alpha^{-1})$, then 
\begin{enumerate}
\item $R_{\lambda|\mathcal{A}_{D_n}} \simeq \epsilon\circ R_{\mu|\mathcal{A}_{D_n}}^\star \Leftrightarrow \lambda=\mu$.
\item $R_{\lambda|\mathcal{A}_{D_n}}\not\simeq \epsilon \circ R_{\gamma,\pm}$.
\item $R_{\gamma,\pm|\mathcal{A}_{D_n}} \simeq \epsilon \circ R_{\delta,\pm|\mathcal{A}_{D_n}}^\star \Leftrightarrow \gamma= \delta$.
\item $R_{\gamma,\pm|\mathcal{A}_{D_n}} \not\simeq \epsilon\circ R_{\delta,\mp|\mathcal{A}_{D_n}}^\star$.
\end{enumerate}
\end{enumerate}
\end{prop}

\begin{proof}
$\textbf{1.}$ is shown in the same way as Lemma $3.4.$ of  \cite{BMM}, because $A_{D_n}$ is generated by $\mathcal{A}_{D_n}$ ans $A_{D_{n-1}}$.

Using Propositions \ref{unitary2} and \ref{transpose2}, it is sufficient to show \textbf{2,3,4} and \textbf{5} to conclude the proof. In the same way as for type $B$, we need to use Lemma \ref{abel}. If $R_{\lambda|\mathcal{A}_{D_n}}\simeq R_{\mu|\mathcal{A}_{D_n}}$ then there exists a character $\eta : A_{D_n} \rightarrow \F_q^\star$ such that $R_\lambda \simeq R_\mu \otimes \eta$. Since $A_{D_n}/\mathcal{A}_{D_n}=<\overline{S_1}>$, there exists $u\in \F_q^\star$ such that for all $d\in A_{D_n}, \eta(d)=u^{\ell(d)}$. We have $R_\lambda(S_1)=uR_\mu(S_1)$. By considering the eigenvalues, we have that $\{\alpha,-1\}=\{u\alpha,-u\}$. Therefore $u=1$ or $\alpha^2=1$. By the conditions on $\alpha$,  $u=1$ and $R_\lambda\simeq R_\mu$, therefore $\lambda=\mu$. Since the set of eigenvalues is of $R_{\gamma,\pm}(S_1)$ is also $\{\alpha,-1\}$, the rest of the proof follows.
\end{proof}

We now give a theorem for double-partitions with an empty component and then results for hook partitions.

\begin{theo}\label{empty}
Let $\lambda=(\lambda_1,\emptyset)\Vdash n$ with $\lambda_1$ not a hook and $G=R_\lambda(\mathcal{A}_{D_n})$. We then have the following properties
\begin{enumerate}
\item If $\F_q=\F_p(\alpha)=\F_p(\alpha+\alpha^{-1})$, then
\begin{enumerate}
\item if $\lambda_1\neq \lambda_1'$, then $G=SL_{n_{\lambda}}(q)$,
\item if $\lambda_1=\lambda_1'$ and $\tilde{\nu}(\lambda)=-1$, then $G\simeq SP_{n_\lambda}(q)$,
\item if $\lambda_1=\lambda_1'$ and $\tilde{\nu}(\lambda)=1$, then $G\simeq \Omega^+_{n_\lambda}(q)$.
\end{enumerate}
\item If $\F_q=\F_p(\alpha)\neq\F_p(\alpha+\alpha^{-1})$, then
\begin{enumerate}
\item if $\lambda_1\neq \lambda_1'$, then $G\simeq SU_{n_{\lambda}}(q^\frac{1}{2})$,
\item if $\lambda_1=\lambda_1'$ and $\tilde{\nu}(\lambda)=-1$, then $G\simeq SP_{n_\lambda}(q^\frac{1}{2})$,
\item if $\lambda_1=\lambda_1'$ and $\tilde{\nu}(\lambda)=1$, then $G\simeq \Omega^+_{n_\lambda}(q^\frac{1}{2})$.
\end{enumerate}
\end{enumerate}
\end{theo}

\begin{proof}
The restriction of $R_\lambda$ to $\mathcal{A}_{D_n}$ is the same as the representation $R_{\lambda_1}$ in type $A$. Since $\tilde{\nu}(\lambda)=\nu(\lambda_1)$, the result follows directly from \cite[Theorem 1.1]{BMM} after noting that $R_\lambda(\mathcal{A}_{A_n})\subset R_{\lambda}(\mathcal{A}_{D_n})$ and that we have the corresponding inclusions by Proposition \ref{bilin2}.
\end{proof}

\begin{prop}\label{timberwolves}
If $\F_q=\F_p(\alpha)=\F_p(\alpha+\alpha^{-1})$, then $R_{([1^{n-1}],[1])}(\mathcal{A}_{D_n})= SL_n(q)$ and if $\F_q=\F_p(\alpha)\neq \F_p(\alpha+\alpha^{-1})$ then $R_{([1^{n-1}],[1])}(\mathcal{A}_{D_n})\simeq SU_n(q^\frac{1}{2})$.
\end{prop}

\begin{proof}
The proof is the same one as the proof of Proposition \ref{lesgourgues}. 
\end{proof}

We write again $A_{1,n}=\{(\lambda_1,\emptyset),\lambda_1 \vdash n\}, A_{2,n} =\{(\emptyset,\lambda_2),\lambda_2 \vdash n\}, A_n = A_{1,n} \cup A_{2,n}$ and\\
$\epsilon_n=\{\lambda \Vdash n, \lambda~\mbox{not a hook}\}$

\begin{theo}\label{result1D}
If $\F_q=\F_p(\alpha)=\F_p(\alpha+\alpha^{-1})$ and $n$ is odd, then the morphism from $\mathcal{A}_{D_n}$ to $\mathcal{H}_{D_n,\alpha}^\times \simeq \underset{\lambda_1>\lambda_2}{\underset{\lambda \vdash\vdash n}\prod}GL_{n_\lambda}(q)$ factorizes through the epimorphism
$$\Phi_{1',n}: \mathcal{A}_{D_n} \rightarrow SL_{n-1}(q) \times SL_n(q) \times \underset{\lambda_1>\lambda_2}{\underset{\lambda\in \epsilon_n,\lambda>\varphi(\lambda)}\prod} SL_{n_\lambda}(q)\times \underset{n_\lambda> n_\mu}{\underset{\lambda\in \epsilon_n,\lambda=\varphi(\lambda)}\prod}OSP(\lambda)'.$$ 
If $\F_q=\F_p(\alpha)=\F_p(\alpha+\alpha^{-1})$ and $n \equiv 0~ (\bmod ~4)$, then the morphism from $\mathcal{A}_{D_n}$ to $\mathcal{H}_{D_n,\alpha}^\times \simeq \underset{\lambda_1>\lambda_2}{\underset{\lambda \vdash\vdash n}\prod}GL_{n_\lambda}(\F_q) \times \underset{\lambda=(\lambda_1,\lambda_1)\vdash n}\prod GL_{n_{\lambda,+}}(q)\times GL_{n_{\lambda,-}}(q)$ factorizes through the epimorphism
$$\Phi_{1',n}: \mathcal{A}_{D_n} \rightarrow SL_{n-1}(q) \times SL_n(q) \times \underset{\lambda_1>\lambda_2}{\underset{\lambda\in \epsilon_n,\lambda>\varphi(\lambda)}\prod} SL_{n_\lambda}(q)\times \underset{\lambda_1> \lambda_2}{\underset{\lambda\in \epsilon_n,\lambda=\varphi(\lambda)}\prod}OSP(\lambda)'\times$$ $$\underset{\lambda>\varphi(\lambda)}{\underset{\lambda=(\lambda_1,\lambda_1)\in \epsilon_n}\prod}SL_{\frac{n_\lambda}{2}}(q)^2\times \underset{\lambda=\varphi(\lambda)}{\underset{\lambda=(\lambda_1,\lambda_1)\in \epsilon_n}\prod}OSP(\lambda,+)'^2.$$ 
If $\F_q=\F_p(\alpha)=\F_p(\alpha+\alpha^{-1})$ and $n \equiv 2~ (\bmod ~4)$ then the morphism from $\mathcal{A}_{D_n}$ to $\mathcal{H}_{D_n,\alpha}^\times \simeq \underset{\lambda_1>\lambda_2}{\underset{\lambda \vdash\vdash n}\prod}GL_{n_\lambda}(\F_q) \times \underset{\lambda=(\lambda_1,\lambda_1)\vdash n}\prod GL_{n_{\lambda,+}}(q)\times GL_{n_{\lambda,-}}(q)$ factorizes through the epimorphism
$$\Phi_{1',n}: \mathcal{A}_{D_n} \rightarrow SL_{n-1}(q) \times SL_n(q) \times \underset{\lambda_1>\lambda_2}{\underset{\lambda\in \epsilon_n,\lambda>\varphi(\lambda)}\prod} SL_{n_\lambda}(q)\times \underset{\lambda_1> \lambda_2}{\underset{\lambda\in \epsilon_n,\lambda=\varphi(\lambda)}\prod}OSP(\lambda)'\times$$ $$\underset{\lambda>\varphi(\lambda)}{\underset{\lambda=(\lambda_1,\lambda_1)\in \epsilon_n}\prod}SL_{\frac{n_\lambda}{2}}(q)^2\times \underset{\lambda=\varphi(\lambda)}{\underset{\lambda=(\lambda_1,\lambda_1)\in \epsilon_n}\prod}SL_{\frac{n_\lambda}{2}}(q).$$ 
In all of the above, $OSP(\lambda)$ is the group of isometries of the bilinear form defined in Proposition \ref{bilin2}.
\end{theo}

In the unitary case, we have an analogous result.

\begin{theo}\label{result2D}
If $\F_q=\F_p(\alpha)\neq \F_p(\alpha+\alpha^{-1})$ and $n$ is odd, then the morphism from $\mathcal{A}_{D_n}$ to $\mathcal{H}_{D_n,\alpha}^\times \simeq \underset{\lambda_1>\lambda_2}{\underset{\lambda \vdash\vdash n}\prod}GL_{n_\lambda}(q)$ factorizes through the morphism
$$\Phi_{2',n}: \mathcal{A}_{D_n} \rightarrow SU_{n-1}(q^\frac{1}{2}) \times SU_n(q^\frac{1}{2}) \times \underset{\lambda_1>\lambda_2}{\underset{\lambda\in \epsilon_n,\lambda>\varphi(\lambda)}\prod} SU_{n_\lambda}(q^\frac{1}{2})\times \underset{n_\lambda> n_\mu}{\underset{\lambda\in \epsilon_n,\lambda=\varphi(\lambda)}\prod}\widetilde{OSP}(\lambda)'.$$ 
If $\F_q=\F_p(\alpha)=\F_p(\alpha+\alpha^{-1})$ and $n \equiv 0~ (\bmod ~4)$, then the morphism from $\mathcal{A}_{D_n}$ to $\mathcal{H}_{D_n,\alpha}^\times \simeq \underset{\lambda_1>\lambda_2}{\underset{\lambda \vdash\vdash n}\prod}GL_{n_\lambda}(\F_q) \times \underset{\lambda=(\lambda_1,\lambda_1)\vdash n}\prod GL_{n_{\lambda,+}}(q)\times GL_{n_{\lambda,-}}(q)$ factorizes through the morphism
$$\Phi_{2',n}: \mathcal{A}_{D_n} \rightarrow SU_{n-1}(q^\frac{1}{2}) \times SU_n(q^\frac{1}{2}) \times \underset{\lambda_1>\lambda_2}{\underset{\lambda\in \epsilon_n,\lambda>\varphi(\lambda)}\prod} SU_{n_\lambda}(q^\frac{1}{2})\times \underset{\lambda_1> \lambda_2}{\underset{\lambda\in \epsilon_n,\lambda=\varphi(\lambda)}\prod}\widetilde{OSP}(\lambda)'\times$$ $$\underset{\lambda>\varphi(\lambda)}{\underset{\lambda=(\lambda_1,\lambda_1)\in \epsilon_n}\prod}SU_{\frac{n_\lambda}{2}}(q^\frac{1}{2})^2\times \underset{\lambda=\varphi(\lambda)}{\underset{\lambda=(\lambda_1,\lambda_1)\in \epsilon_n}\prod}\widetilde{OSP}(\lambda,+)'^2.$$ 
If $\F_q=\F_p(\alpha)=\F_p(\alpha+\alpha^{-1})$ and $n \equiv 2~ (\bmod ~4)$, then the morphism from $\mathcal{A}_{D_n}$ to $\mathcal{H}_{D_n,\alpha}^\times \simeq \underset{\lambda_1>\lambda_2}{\underset{\lambda \vdash\vdash n}\prod}GL_{n_\lambda}(\F_q) \times \underset{\lambda=(\lambda_1,\lambda_1)\vdash n}\prod GL_{n_{\lambda,+}}(q)\times GL_{n_{\lambda,-}}(q)$ factorizes through the morphism
$$\Phi_{2',n}: \mathcal{A}_{D_n} \rightarrow SU_{n-1}(q^\frac{1}{2}) \times SU_n(q^\frac{1}{2}) \times \underset{\lambda_1>\lambda_2}{\underset{\lambda\in \epsilon_n,\lambda>\varphi(\lambda)}\prod} SU_{n_\lambda}(q^\frac{1}{2})\times \underset{\lambda_1> \lambda_2}{\underset{\lambda\in \epsilon_n,\lambda=\varphi(\lambda)}\prod}\widetilde{OSP}(\lambda)'\times$$ $$\underset{\lambda>\varphi(\lambda)}{\underset{\lambda=(\lambda_1,\lambda_1)\in \epsilon_n}\prod}SU_{\frac{n_\lambda}{2}}(q^\frac{1}{2})^2\times \underset{\lambda=\varphi(\lambda)}{\underset{\lambda=(\lambda_1,\lambda_1)\in \epsilon_n}\prod}SU_{\frac{n_\lambda}{2}}(q^\frac{1}{2}).$$ 
In all of the above, $\widetilde{OSP}(\lambda)$ is the group of isometries associated with the bilinear form over $\F_{q^\frac{1}{2}}$ obtained from the one in Proposition \ref{bilin2} using Proposition \ref{coolprop}.
\end{theo}

Those two theorems (except for the surjectivity) follow from Propositions \ref{bilin2}, \ref{transpose2}, \ref{chaise}, \ref{unitary2}, \ref{timberwolves}, Theorem \ref{empty} and Proposition \ref{isomorphisme2}. It now remains to check that $\Phi_{1',n}$ and $\Phi_{2',n}$ are surjective in all cases.

\section{The case $n=4$}

In this section, we prove the result for $n=4$.

The double-partitions to consider for $n=4$ are $([4],\emptyset)$, $([3,1],\emptyset)$, $([2,2],\emptyset)$, $([2,1,1],\emptyset)$, $([1^4],\emptyset)$, $([3],[1])$, $([2,1],[1])$, $([1^3],[1])$, $([2],[2])$, $([2],[1^2])$ and $([1^2],[1^2])$.

By Proposition \ref{transpose2}, if we know the image for $\lambda$, we know the image for $\varphi(\lambda)$. By Proposition \ref{empty}, we know the image for doubles-partitions with an empty component. By Proposition \ref{timberwolves}, we know the image for $([1^3],[1])$ and by Proposition \ref{chaise}, we know the image for $([2],[1^2])$ using the image of $([1^3],[1])$. The only double-partitions left to consider are $([1^2],[1^2])$ and $([2,1],[1])$. 

\begin{lemme}\label{Berkeley}
If $\F_q=\F_p(\alpha)=\F_p(\alpha+\alpha^{-1})$, then $R_{[2,1],[1]}(\mathcal{A}_{D_4}) \simeq SP_8(q)$.

If $\F_q=\F_p(\alpha)\neq \F_p(\alpha+\alpha^{-1})$, then $R_{[2,1],[1]}(\mathcal{A}_{D_4}) \simeq SP_8(q^\frac{1}{2})$.
\end{lemme}

\begin{proof}
Assume first that $\F_q=\F_p(\alpha)=\F_p(\alpha+\alpha^{-1})$. Using Proposition \ref{bilin2}, there exists $P\in GL_8(q)$ such that $G=PR_{[2,1],[1]}(\mathcal{A}_{D_4})P^{-1}\subset SP_8(q)$. Using Lemma \ref{branch}, we have that $R_{[2,1],[1]}(\mathcal{A}_{D_3})=R_{[2],[1]}\times R_{[1^2],[1]}\times R_{[2,1],\emptyset}(\mathcal{A}_{D_3})\simeq SL_3(q)\times SL_2(q)$, where $SL_3(q)$ is in a twisted diagonal embedding and $SL_2(q)$ is in a natural representation using Goursat's Lemma and the previous arguments. Using the same arguments as before and Lemma \ref{field} with the natural representation of $SL_2(q)$, we know $G$ is primitive, tensor-indecomposable, irreducible, perfect and cannot be realized in a natural representation over a proper subfield of $\F_q$. This implies that $G$ cannot be included in a maximal subgroup of class $\mathcal{C}_1, \mathcal{C}_2,\mathcal{C}_4$ or $\mathcal{C}_5$. Since it contains a transvection of $SL_2(q)$, we have that it cannot be contained in a maximal group of class $\mathcal{C}_3$. Assume that $G$ is included in a maximal subgroup of $SP_8(q)$. By the Tables 8.48. and 8.49. in \cite{BHRC}, the only possible maximal subgroups and their order or a quantity their order divides are given below

\begin{enumerate}
\item $2_{-}^{1+6.}SO_6^-(2), 51840$
\item $2_{-}^{1+6.}\Omega_6^-(2), 25920$
\item $(SP_2(q)\circ SP_2(q)\circ SP_2(q)).2^2.\mathfrak{S}_3, 24q^3(q^2-1)^3$
\item $2^{.}PSL_2(7), 336$
\item $2^{.}PSL_2(7)^{.}2, 672$
\item $2^{.}\mathfrak{A}_6, 720$
\item $2^{.}\mathfrak{A}_6.2_2, 1440$
\item $2^{.}PSL_2(17), 4896$
\item $2^{.}PSL_2(q), q(q^2-1)$
\item $2^{.}PSL_2(q^3).3, 3q^3(q^6-1)$
\end{enumerate} 

The order of $SL_3(q)\times SL_2(q)$ is $q^4(q^2-1)^2(q^3-1)$, therefore cases $3$, $9$ and $10$ are excluded. We have that $\alpha$ is of order greater than $16$ and $p\neq 2$, therefore we have that $q>17$ and, therefore $q\geq19$ and $\vert G\vert \geq 19^4(19^2-1)(19^3-1)=115828887772800$. This excludes all the remaining cases. It follows that $G$ can be included in no maximal subgroup of $SP_8(q)$, therefore $G=SP_8(q)$.

Assume now that $\F_q=\F_p(\alpha)\neq \F_p(\alpha+\alpha^{-1})$. There exists $P\in GL_8(q)$ such that $G=PR_{[2,1],[1]}(\mathcal{A}_{D_4})P^{-1}\subset SP_8(q^\frac{1}{2})$ and $G$ contains $H\simeq SU_3(q^\frac{1}{2})\times SU_2(q^\frac{1}{2})$, where $SU_3(q^\frac{1}{2})$ is in a twisted diagonal embedding and $SU_2(q^\frac{1}{2})$ is in a natural representation. We can no longer use Lemma \ref{field} in this case, but since $\epsilon(\alpha)=\alpha^{-1}$, we have up to conjugation that $\diag(I_6,\begin{pmatrix}
\alpha & 0\\
0 & \alpha^{-1}
\end{pmatrix})\in G$. It follows that $\alpha+
\alpha^{-1}$ belongs to the field generated by the traces of the elements of $G$. This shows that any field over which $G$ is realized in a natural representation contains $\F_{q^\frac{1}{2}}$. By the above, in this case, we have that $G$ is primitive, tensor-indecomposable, irreducible, perfect and cannot be realized in a natural representation over a proper subfield of $\F_{q^\frac{1}{2}}$. This implies that $G$ cannot be included in a maximal subgroup of $SP_{8}(q^\frac{1}{2})$ of class $\mathcal{C}_1, \mathcal{C}_2,\mathcal{C}_4$ or $\mathcal{C}_5$. It contains a transvection of $SU_2(q^{\frac{1}{2}})$, therefore it cannot be included in a maximal subgroup of class $\mathcal{C}_3$. Assume $G$ is included in a maximal subgroup of $SP_8(q^{\frac{1}{2}})$. We list below the possible maximal subgroups and their order or a quantity their order divides

\begin{enumerate}
\item $2_{-}^{1+6.}SO_6^-(2), 51840$
\item $2_{-}^{1+6.}\Omega_6^-(2), 25920$
\item $(SP_2(q^{\frac{1}{2}})\circ SP_2(q^{\frac{1}{2}})\circ SP_2(q^{\frac{1}{2}})).2^2.\mathfrak{S}_3, 24q^{\frac{3}{2}}(q-1)^3$
\item $2^{.}PSL_2(7), 336$
\item $2^{.}PSL_2(7)^{.}2, 672$
\item $2^{.}\mathfrak{A}_6, 720$
\item $2^{.}\mathfrak{A}_6.2_2, 1440$
\item $2^{.}PSL_2(17), 4896$
\item $2^{.}PSL_2(q^{\frac{1}{2}}), q^{\frac{1}{2}}(q-1)$
\item $2^{.}PSL_2(q^3).3, 3q^{\frac{3}{2}}(q^3-1)$
\end{enumerate} 

We have $\vert H\vert =q^2(q-1)^2(q^{\frac{3}{2}}+1)$. This excludes cases $3$, $9$ and $10$. We have $\alpha^{q^\frac{1}{2}}=\epsilon(\alpha)=\alpha^{-1}$. It follows that $q^\frac{1}{2}+1>16$ and, therefore $q^\frac{1}{2}\geq 17$ because $p\neq 2$. This implies that $\vert H\vert \geq 34042058459136$.  This proves that $G=SP_8(q^\frac{1}{2})$ and concludes the proof of the lemma. 
\end{proof}

\begin{lemme}
If $\F_q=\F_p(\alpha)=\F_p(\alpha+\alpha^{-1})$, we have $R_{([1^2],[1^2]),+}(\mathcal{A}_{D_4})=R_{([1^2],[1^2]),-}(\mathcal{A}_{D_4})=SL_3(q)$.

If $\F_q=\F_p(\alpha)\neq \F_p(\alpha+\alpha^{-1})$, we have $R_{([1^2],[1^2]),+}(\mathcal{A}_{D_4})=R_{([1^2],[1^2]),-}(\mathcal{A}_{D_4})\simeq SU_3(q^\frac{1}{2})$.
\end{lemme}

\begin{proof}
The result follows from Lemma \ref{branch} and the fact that $R_{([1^2],[1])}(\mathcal{A}_{D_3})$ is equal to the group we want in both cases. 
\end{proof}

\section{Surjectivity of $\Phi_{1',n}$ for $n\geq 5$}

In this section, we use results of the previous sections to prove by induction on $n$ the main results for type $D$. We will here conclude the proof of Theorem \ref{result1D}.

Assume first that $\F_q=\F_p(\alpha)=\F_p(\alpha+\alpha^{-1})$. Using Proposition \ref{isomorphisme2}, by the same kind of arguments as for type $B$, we can use Goursat's Lemma to show the morphism is surjective upon each component. This means it is sufficient to show the following theorem.

\begin{theo}\label{Sandburg}
Let $\lambda=(\lambda_1,\lambda_2)\Vdash n$ not a hook, such that $\lambda_1\geq \lambda_2$. We write $G(\lambda)=R_{\lambda}(\mathcal{A}_{D_n})$ if $\lambda_1> \lambda_2$, $G(\lambda,+)=R_{\lambda,+}(\mathcal{A}_{D_n})$ and $G(\lambda,-)=R_{\lambda,-}(\mathcal{A}_{D_n})$ otherwise. We then have the following possibilities.
\begin{enumerate}
\item If $\lambda=([2,1^{n-2}],\emptyset)$, then $G(\lambda)=SL_{n-1}(q)$.
\item If $\lambda=([1^{n-1}],[1])$, then $G(\lambda)=SL_n(q),$
\item If $\lambda \in \epsilon_n, \lambda_1>\lambda_2$ and $\lambda>\varphi(\lambda)$, then $G(\lambda)=SL_{n_\lambda}(q)$,
\item If $\lambda \in \epsilon_n, \lambda_1>\lambda_2$ and $\lambda=\varphi(\lambda)$, then we have the following possibilities.
\begin{enumerate}
\item If $\tilde{\nu}(\lambda)=-1$, then $G(\lambda)\simeq SP_{n_\lambda}(q)$.
\item If $\tilde{\nu}(\lambda)=1$, then $G(\lambda)\simeq \Omega_{n_\lambda}^+(q)$.
\end{enumerate}
\item If $\lambda=(\lambda_1,\lambda_1)\in \epsilon_n$, then we have the following possibilities.
\begin{enumerate}
\item If $\varphi(\lambda)>\lambda$, then $G(\lambda,+)=G(\lambda,-)\simeq SL_{\frac{n_{\lambda}}{2}}(q)$.
\item If $\varphi(\lambda)=\lambda$, then we have the following possibilities.
\begin{enumerate}
\item If $n\equiv 0~(\bmod ~ 4)$ then
\begin{enumerate}
\item  if $\tilde{\nu}(\lambda)=-1$ then $G(\lambda,+)=G(\lambda,-)\simeq SP_{n_\lambda}(q)$,
\item if $\tilde{\nu}(\lambda)=1$ then $G(\lambda,+)=G(\lambda,-)\simeq \Omega_{n_\lambda}^+(q)$.
\end{enumerate}
\item If $n\equiv 2~(\bmod ~ 4)$ then $G(\lambda,+)=G(\lambda,-)=SL_{\frac{n_\lambda}{2}}(q)$.
\end{enumerate} 
\end{enumerate}
\end{enumerate}
\end{theo}

\begin{proof}
For $n=4$, we have the result by the previous section. Theorem \ref{empty} gives us the result for double-partitions with an empty component and Proposition \ref{timberwolves} gives us the result for double-partitions with two rows or two columns and one of the components of size one. For $n\geq 5$, we proceed by induction but we must first treat the following cases separately : $([2,2],[1]), ([1^3],[1^2])$ and $(([1^3],[1^3]),\pm)$.

By Lemma \ref{Berkeley}, Theorem \ref{empty}, Lemma \ref{branch} and Goursat's Lemma, we have that \\$R_{([2,2],[1])}(\mathcal{A}_{D_4})\simeq SP_8(q)\times SP_2(q)$. By Theorem \ref{transvections} and Lemma 5.6. of \cite{BMM}, we have that $R_{([2,2],[1])}(\mathcal{A}_{D_5})\in \{SL_{10}(q),SU_{10}(q^\frac{1}{2}),SP_{10}(q)\}$. We have that $([2,2],[1])=\varphi(([2,2],[1]))$ and $\tilde{\nu}([2,2],[1])=(-1)^{\frac{4-2}{2}}(-1)^{\frac{1-1}{2}}=-1$. This implies by Proposition \ref{bilin2} that up to conjugation in $GL_{10}(q)$, we have that $R_{([2,2],[1])}(\mathcal{A}_{D_5})\subset SP_{10}(q)$, therefore we have that $R_{([2,2],[1])}(\mathcal{A}_{D_5})\simeq SP_{10}(q)$.

In the same way, we have that $R_{([1^3],[1^2])}(\mathcal{A}_{D_4})\simeq SL_4(q)\times SL_3(q)\times SL_3(q)$, therefore $G(([1^3],[1^2]))=R_{([1^3],[1^2])}(\mathcal{A}_{D_5})$ is in $\{SL_{10}(q),SU_{10}(q^\frac{1}{2}),SP_{10}(q)\}$. By Proposition \ref{isomorphisme2}, we know that $G(([1^3],[1^2]))$ preserves no bilinear form, therefore we only have to exclude the unitary case. Assume that $G(([1^3],[1^2]))$ is included up to conjugation in $SU_{10}(q^\frac{1}{2})$. There then exists an automorphism $\epsilon$ of order $2$ of $\F_q$ such that each $M$ in $G(\lambda)$ is conjugate to ${}^t\!\epsilon((M^{-1}))$. In particular $G(([1^3],[1^2]))$ contains a natural $SL_2(q)$. This implies that $\diag(I_8,\begin{pmatrix}\alpha & 0\\
0 & \alpha^{-1}
\end{pmatrix})$ is conjugate to $\diag(I_8,{}^t\!\epsilon((\begin{pmatrix}
\alpha & 0\\
0 & \alpha^{-1}
\end{pmatrix})^{-1}))$. Taking the traces of those matrices implies that $\epsilon(\alpha+\alpha^{-1}+8) =\alpha+\alpha^{-1}+8$. We have that $\F_q=\F_p(\alpha)=\F_p(\alpha+\alpha^{-1})$, therefore this shows that $\epsilon$ is trivial which is a contradiction. It follows that $G(([1^3],[1^2]))=SL_{10}(q)$.

By Lemma \ref{branch} and the fact that $R_{([1^3],[1^2])}(\mathcal{A}_{D_5})=SL_{10}(q)$, we have that $SL_{10}(q)\subset R_{([1^3],[1^3]),\pm}(\mathcal{A}_{D_6})\subset SL_{10}(q)$. It follows that $R_{([1^3],[1^3]),\pm}(\mathcal{A}_{D_6})= SL_{10}(q)$.

We now proceed to the induction on $n$ using Theorem \ref{CGFS}.

Let $n\geq 5$ and $\lambda\Vdash n$. Suppose the theorem is true for $n-1$. We use Lemma \ref{branch} for different possibilities to show that $G(\lambda)$ or $G(\lambda,\pm)$ contains a subgroup verifying the same properties as in type B.

\begin{enumerate}
\item If $\lambda=(\lambda_1,\lambda_2)$ and $\lambda_1>\lambda_2$ and $\lambda\neq \varphi(\lambda)$ then $\varphi(\lambda)=(\lambda_1',\lambda_2')$ because the order we defined for partitions of $\frac{n}{2}$ verifies that, if $\lambda_1\neq \lambda_2'$ and $\lambda_1>\lambda_2$, then $\lambda_1'>\lambda_2'$. We then have $\lambda_1\neq \lambda_1'$ or $\lambda_2\neq \lambda_2'$.
\begin{enumerate}
\item If $\lambda_2'\neq \lambda_2$, then there exists $\mu_2 \subset \lambda_2$ such that $\mu_2'\not\subset\lambda_2$. We have that $(\lambda_1',\mu_2')\not\subset (\lambda_1,\lambda_2)$ because $\mu_2'\not\subset \lambda_2$ and $(\mu_2',\lambda_1') \not\subset (\lambda_1,\lambda_2)$ because otherwise $\lambda_1'=\lambda_2$ and, therefore $\lambda_2'=\lambda_1$. This shows that $G(\lambda)$ contains a natural $SL_3(q)$.
\item If $\lambda_2=\lambda_2'$ and $\lambda_1\neq \lambda_1'$, then there exists $\mu_1\subset \lambda_1$ such that $\mu_1'\not\subset\lambda_1$. We then have $(\mu_1',\lambda_2') \not\subset (\lambda_1,\lambda_2)$ because $\mu_1'\not\subset \lambda_1$ and $(\lambda_2',\mu_1')\not\subset (\lambda_1,\lambda_2)$ because $\lambda_2'\neq \lambda_1$. This shows that $G(\lambda)$ also contains a natural $SL_3(q)$ in this case.
\end{enumerate}
\item If $\lambda=(\lambda_1,\lambda_2)=\varphi(\lambda)$ and $\lambda_1>\lambda_2$, then
\begin{enumerate}
\item If $\varphi(\lambda)=(\lambda_1',\lambda_2')$, then
\begin{enumerate}
\item If $\lambda_1$ and $\lambda_2$ are square partitions, then $R_{\lambda}(\mathcal{A}_{D_{n-1}})=G(\mu_1,\lambda_2)\times G(\lambda_1,\mu_2)$ and since $\tilde{\nu}(\mu_1,\lambda_2)=\tilde{\nu}(\lambda_1,\mu_2)=\tilde{\nu}(\lambda)$, we have that :
\begin{enumerate}
\item If $\tilde{\nu}(\lambda)=1$, then $\Omega_{n_{(\mu_1,\lambda_2)}}^+(q)\times \Omega_{n_{(\lambda_1,\mu_2)}}^+(q)\subset G(\lambda)\subset \Omega_{n_{(\lambda_1,\lambda_2)}}^+(q)$.
\item If $\tilde{\nu}(\lambda)=-1$ then $SP_{n_{(\mu,\lambda_2)}}(q)\times SP_{n_{(\lambda_1,\mu_2)}}(q) \subset G(\lambda)\subset SP_{n_\lambda}(q)$. It follows that $G(\lambda)$ is an irreducible group generated by transvections because it is normally generated by the group on the left of our inclusions, therefore by Theorem \ref{transvections}, we have that $G(\lambda)$ is equal to the group on the right and the theorem is proved in this case.
\end{enumerate}
\item If $\lambda_1$ or $\lambda_2$ is not a square partition then there exists $\mu \subset \lambda$ such that $\varphi(\mu)\neq \mu$. It follows that $\varphi(\mu)\subset \lambda$ or $\sigma(\varphi(\mu))\subset \lambda$, therefore $G(\lambda)$ contains a twisted diagonal $SL_3(q)$.
\end{enumerate}
\item If $\varphi(\lambda)=(\lambda_2',\lambda_1')$, then if $\mu \subset \lambda_2$, we have that $(\lambda_1,\mu)\subset (\lambda_1,\lambda_2), \varphi((\lambda_1,\mu))=(\lambda_1',\mu')\not\subset (\lambda_1,\lambda_2)$ because $\lambda_1\neq\lambda_1'$. We have that $(\mu',\lambda_1')\subset (\lambda_1,\lambda_2)$, therefore $G(\lambda)$ contains a twisted diagonal $SL_3(q)$.
\end{enumerate}
\item If $\lambda=(\lambda_1,\lambda_1)\neq (\lambda_1',\lambda_1')$, then there exists $\mu_1\subset \lambda_1$ such that $\mu_1'\not\subset \lambda_1$. It follows that $(\lambda_1',\mu_1')\not\subset (\lambda_1,\lambda_1)$ and $(\mu_1',\lambda_1')\not\subset (\lambda_1,\lambda_1)$. This shows that $G(\lambda,\pm)$ contains a natural $SL_3(q)$.
\item If $\lambda=(\lambda_1,\lambda_1) =(\lambda_1',\lambda_1)$ and $\lambda_1$ is not a square partition, then there exists $\mu_1\subset \lambda_1$ such that $\mu_1\neq \mu_1'$, therefore $(\lambda_1,\mu_1)\neq \varphi((\lambda_1,\mu_1))=(\lambda_1',\mu_1')$. We have that $(\mu_1',\lambda_1')\subset(\lambda_1,\lambda_1)$, therefore $G(\lambda,\pm)$ contains a twisted diagonal $SL_3(q)$.
\item If $\lambda=(\lambda_1,\lambda_1)=(\lambda_1',\lambda_1')$ and $\lambda_1$ is a square partition, then we have the two following possibilities.
\begin{enumerate}
\item If $n \equiv 0~(\bmod ~ 4)$, then for all $\mu\subset \lambda$, we have that $\tilde{\nu}(\lambda)=\nu(\lambda_1)^2(-1)^{(\frac{n}{2})^2}=1=\tilde{\nu}(\mu)$. This is because if $\lambda_1$ is a square, then the only sub-partition $\mu_1$ of $\lambda_1$ verifies $\nu(\mu_1)=\nu(\lambda_1)$. By the branching rule, we have that $\Omega^+_{\frac{n_\lambda}{2}}(q)\subset G(\lambda,\pm) \subset \Omega_{\frac{n_\lambda}{2}}^+(q)$. It follows that $G(\lambda)\simeq \Omega_{\frac{n_\lambda}{2}}^+(q)$ and the theorem is proved in this case.
\item If $n \equiv 2~(\bmod ~ 4)$, then $\tilde{\nu}(\mu)=\nu(\lambda_1)^2=1$ for all $\mu\subset \lambda$. The branching rule shows that $\Omega_{\frac{n_\lambda}{2}}^+(q)\subset G(\lambda,\pm) \subset SL_{\frac{n_\lambda}{2}}(q)$. By Proposition \ref{isomorphisme2}, $G(\lambda,\pm)$ preserves no bilinear form, therefore $G(\lambda,\pm)=SL_{\frac{n_\lambda}{2}}(q)$.
\end{enumerate}
\end{enumerate}

In all the cases where $G(\lambda)$ or $G(\lambda,\pm)$ contains a natural $SL_3(q)$ or a twisted diagonal $SL_3(q)$, we can use exactly the same arguments as in \cite{BMM} because if the morphism $A_{A_n}$ to $A_{D_n}$ defined  by $S_i \mapsto S_i$ is trivial then $A_{D_n}$ is trivial. 

The only case we need to treat separately is $(([2,1],[2,1]),\pm)$ because $n=6$. We need a separate argument to show that $G(([2,1],[2,1]),\pm)$ is tensor-indecomposable. In this case $R_{([2,1],[2,1])}(\mathcal{A}_{D_5})=G([2,1],[1^2])\times G([2,1],[2])=SL_{20}(q)\times SL_{20}(q)$. If $G(([2,1],[2,1]),\pm) \subset SL_{40}(q) \otimes SL_2(q)$, then the morphism from \break$R_{([2,1],[2,1])}(\mathcal{A}_{D_5})$ to $SL_2(q)$ is trivial. Since $R_{([2,1],[2,1])}(\mathcal{A}_{D_5})$ normally generates $G(([2,1],[2,1]),\pm)$, $G(([2,1],[2,1]),\pm)$ is included in $SL_{40}(q)\times SL_{40}(q)$. This contradicts its irreducibility.

This shows that it is sufficient to consider case $2.a.i.A$. Assume we are in case $2.a.i.A$. We then have that $G(\lambda)\subset \Omega_{n_\lambda}^+(q)$ is generated by a conjugacy class of long root elements and $G(\lambda)$ is irreducible. Since $p\neq 2$, if we check that $O_p(G(\lambda))\subset [G,G]\cap Z(G)$, then we can apply Theorem \ref{theoKantor}. Applying Clifford's Theorem (Theorem 11.1 of \cite{C-R}), we have that $Res_{O_p(G(\lambda))}^{G(\lambda)}(V)$ is semisimple and since $O_p(G(\lambda))$ is a $p$-group, its only irreducible representation over $\F_q$ is the trivial one. This shows that $Res_{O_p(G(\lambda))}^{G(\lambda)}(V)$ is trivial, therefore $O_p(G(\lambda))= 1$ and all the assumptions of Theorem I of Kantor are verified (the minimal dimension in this case is greater than or equal to the dimension of $([3,3,3],[2,2])$ and the dimension of $([4,4,4,4],[1])$ which are $42\times 2\times \binom {13} {4}\geq 5$ and $17\times 24024\geq 5$). This shows that we are in one of the following cases :
\begin{enumerate}
\item $G(\lambda)\simeq \Omega_{n_\lambda}^{+}(q')$, and $q'|q$,
\item $G(\lambda)\simeq \Omega_{n_\lambda}^{-}(q') \subset \Omega^{+}_{n_\lambda}(q'^2), q'^2 | q$ and $n_\lambda$ is even,
\item $G(\lambda)\simeq SU_{\frac{n_\lambda}{2}}(q') \subset \Omega^{+}_{n_\lambda}(q'), n_\lambda \equiv 0~ (\bmod ~ 4)$ and $q'|q$,
\item $G(\lambda)\subset \Omega_8^+(q')$ and $q'|q$,
\item $G(\lambda)\simeq [G_2(q'),G_2(q')] \subset \Omega_7(q')$ and $q'|q$,
\item $G\simeq  ~^3 D_4(q')\subset \Omega_8^+(q'^3)$ and $q'^3|q$.
\end{enumerate}
Since $n\geq 13$, $\alpha^{q-1}=1$ and $\alpha$ is of order greater than $2n$, we have $q\geq 29$ and $n_\lambda \geq \min(84\binom{13}{4},17\times 24024)$. This proves that Cases $4$, $5$ and $6$ are excluded by cardinality arguments.

Let us show that $3.$ is also excluded by cardinality arguments. We write $\vert G\vert_p$ the order of a Sylow $p$-subgroup of a group $G$, therefore that $\vert SU_{\frac{n_\lambda}{2}}(q') \vert_p=q'^{\frac{\frac{n_\lambda}{2}(\frac{n_\lambda}{2}-1)}{2}}$. We know that $G(\lambda)$ contains $\Omega_{n_1}^+(q)\times \Omega_{n_2}^+(q)$. It follows that if $\lambda_1$ is the square partition of $r$ and $\lambda_2$ is the square partition of $n-r<r$, writing $a_l$ for the number of standard tableaux associated with a square partition of $l\in \N^\star$, we have that $n_\lambda=\binom {n} {r} a_{r}a_{n-r}$, $n_1=\binom {n-1} {r-1} a_{r}a_{n-r}$ and $n_2=\binom {n-1} {r} a_{r}a_{n-r}$. Note that $a_r$ is even because $r>1$ and using the branching rule twice, we get that $a_r$ is equal to twice the dimension of the two partitions we get by removing first the only extremal node and then one of the two extremal nodes of the resulting partition. It follows that $\vert \Omega_{n_1}^+(q)\times \Omega_{n_2}^+(q)\vert_p = q^{\frac{n_1}{2}(\frac{n_1}{2}-1)+\frac{n_2}{2}(\frac{n_2}{2}-1)}$. To exclude $3$, it is sufficient to show that this quantity is strictly greater than $q^{\frac{\frac{n_\lambda}{2}(\frac{n_\lambda}{2}-1)}{2}}$. If we write $A$ the $q$-logarithm of the quotient of those two quantities, we have that :
\begin{eqnarray*}
A & = & \frac{n_1}{2}\left(\frac{n_1}{2}-1\right)+\frac{n_2}{2}\left(\frac{n_2}{2}-1\right)-\frac{\frac{n_\lambda}{2}(\frac{n_\lambda}{2}-1)}{2}\\
 & = & \frac{n_1}{2}^2+\frac{n_2}{2}^2-\frac{\frac{n_\lambda}{2}^2}{2}-\frac{n_\lambda}{4}\\
 & = & \frac{\frac{n_\lambda}{2}^2}{2}-2\frac{n_1}{2}\frac{n_\lambda-n_1}{2}-\frac{n_\lambda}{4}\\
 & = & \frac{(\frac{n_\lambda}{2}-n_1)^2}{2}-\frac{n_\lambda}{4}\\
 & = & \frac{\left(\frac{\binom {n} {r} a_ra_{n-r}}{2}-\binom {n-1} {r-1}a_r a_{n-r}\right)^2}{2}-\frac{\binom {n} {r} a_r a_{n-r}}{4}\\
 & = &  a_ra_{n-r} \left(\frac{(\binom {n} {r} (\frac{1}{2}-\frac{r}{n}))^2}{2}a_ra_{n-r}-\frac{\binom{n}{r}}{4}\right)\\
 & = & \frac{a_r a_{n-r}}{4} \binom {n} {r}\left(2a_r a_{n-r} \binom {n} {r} \left(\frac{2r-n}{2n}\right)^2-1\right). 
\end{eqnarray*} 

This shows that $A > 0$ if and only if $2a_r a_{n-r} \binom {n} {r} \frac{(2r-n)^2}{4n^2}>1$. Using the branching rule and the hook formula, we get  : $a_1=1, a_4=2, a_9=42, a_{16}=\frac{16!}{7\times 6^2\times 5^3\times 4^4\times 3^3\times 2^2}=24024>81\times 16^2, a_{25}=701149020> 81\times 25^2$ and $a_{36}>81\times 36^2$. Let $k\geq 6$, assume $a_{k^2} > 81(k^2)^2$. The branching rule shows that $a_{(k+1)^2}> 2a_{k^2} > 81(2k^4)> 81(k^4+4k^3+6k^2+4k+1)=81((k+1)^2)^2$, the last inequality being true because $k\geq 6$. It follows that for all $k\geq 4$, we have that $a_{k^2}>81\times(k^2)^2$. In our case, we have that $r\geq 16$ or $r=9$ and $n-r=4$. If $r\geq 16$ and $n-r\geq 2$ then we have $a_ra_{n_r}\geq r^2(n-r)^2\geq 4r^2\geq 2r^2+2r(n-r)\geq (r+n-r)^2\geq n^2$. It follows that $2a_r a_{n-r} \binom {n} {r} \frac{(2r-n)^2}{4n^2}\geq 2n^2\binom {n} {r} \frac{1}{4n^2}=\frac{\binom {n} {r}}{2}>1.$ If $r\geq 16$ and $n-r=1$, then $2a_r a_{n-r} \binom {n} {r} \frac{(2r-n)^2}{4n^2}=\frac{8na_{n-1}}{4n^2}=\frac{2a_{n-1}}{n}>162>1.$ If $r=9$ and $n-r=4$, then $2a_r a_{n-r} \binom {n} {r} \frac{(2r-n)^2}{4n^2}=2\times 42\times 2 \binom {13} {9} \frac{(18-13)^2}{4\times 13^2}>1.$ This shows that independently of $r$ and $n-r$, we have that $A>0$. This proves that $3.$ is excluded.

We have that $\vert \Omega_{n_\lambda}^+(q^\frac{1}{2})\vert_p=q^{\frac{\frac{n_\lambda}{2}(\frac{n_\lambda}{2}-1)}{2}}$. The previous arguments show that $2.$ is also impossible.

The only remaining possibility is $1$ and using again the same arguments, we get $q'>q^\frac{1}{2}$, therefore $q'=q$ and this concludes the proof of Theorem \ref{Sandburg}.
\end{proof}

In the unitary case, i.e. $\F_q=\F_p(\alpha)=\F_p(\alpha+\alpha^{-1})$, all the arguments are analogous.
 
\newpage

$ $

\newpage
 
 \chapter{Type $I_2(m)$}\label{dihedralchapter}
 \section{$m$ odd}\label{sectionmodd}
 
 The main difficulty in finding the image of Artin groups of dihedral type inside their finite Hecke algebras arises from the various field extensions which intervene. When $m$ is even they can be quite complex. In this section, we only consider $m$ odd. The outline of the proof is to first determine the image inside each $2$-dimensional using Dickson's Theorem and then recover the image inside the full Iwahori-Hecke algebra using Goursat's Lemma. The main difficulty will be in the use of Goursat's Lemma, we will need to introduce the equivalence relation from Lemma \ref{Isomorphism} and the the proof will be computational. The image in type $I_2(5)$ will be useful for inductive arguments in type $H_3$.
 
 \bigskip
 
 Let $m\geq 5$ be an odd integer and $p$ a prime number such that there exists $\alpha \in (\overline{\F_p})^\times$ of order not dividing $m$ and not in $\{1,2,3,4,5,6,10\}$ and $\theta\in \overline{\F_p}$ a primitive m-th root of unity. Moreover, we assume $\alpha+\alpha^{-1}\neq \theta^j+\theta^{-j}$ for all $j$ from $1$ to $\frac{m-1}{2}$. For $j\in \N$, we write $\F_{q_j}=\F_p(\alpha,\theta^j+\theta^{-j})$. We write $\F_q$ for the smallest field containing all $\F_{q_j}$,  for $j$ between $1$ and $\frac{m-1}{2}$. 
 
 \smallskip Note that we have $\F_q=\F_{q_1}$, this can be seen using Chebyshev polynomials for example.
  
  \smallskip Note also that when $\F_{q_j}=\F_p(\alpha,\theta^j+\theta^{-j})\neq \F_p(\alpha+\alpha^{-1},\theta^j+\theta^{-j})$, the extension is always of degree $2$ since $X^2-(\alpha+\alpha^{-1})X+1$ is an irreducible polynomial of degree $2$ such that $\F_p(\alpha,\theta^j+\theta^{-j})=\F_p(\alpha+\alpha^{-1},\theta^j+\theta^{-j})/(X^2-(\alpha+\alpha^{-1})X+1)$. This implies that $q_j$ is a prime and $q_j^{\frac{1}{2}}$ is well-defined.
 
 \begin{Def}
 The Iwahori-Hecke algebra of dihedral type $I_2(m)$ which we write $\mathcal{H}_{I_2(m),q}$ is the $\F_q$-algebra with the following presentation :
 
 Generators : $T_t,T_s$.
 
 Relations : 
 
$(T_s-\alpha)(T_s+1)=0$,
 
$(T_t-\alpha)(T_t+1)=0$,
 
$\underset{m}{\underbrace{T_sT_tT_s\dots}}=\underset{m}{\underbrace{T_tT_sT_t\dots}}$.
 \end{Def}
 
We can then use the Kilomoyer-Solomon model \cite{G-P} (Theorem $8.3.1$).
 
\begin{theo}\label{melyssaestlameilleure}
The following matrix model gives a decomposition into pairwise\\ non-isomorphic irreducible modules of $\mathcal{H}_{I_2(m),q}$ :
\begin{enumerate}
\item $\op{Ind} : T_s \mapsto \alpha, T_t\mapsto \alpha$.
\item $\epsilon : T_s \mapsto -1, T_t \mapsto -1$.
\item For $j\in [\![1,\frac{m-1}{2}]\!], T_s \mapsto \rho_j(T_s)=\begin{pmatrix}
-1 & 0\\
1 & \alpha
\end{pmatrix}, T_t \mapsto \rho_j(T_t)=\begin{pmatrix}
\alpha & \alpha(2+\theta^j+\theta^{-j})\\
0 & -1
\end{pmatrix}$.
\end{enumerate}
\end{theo}

\begin{proof}

We already know these models give us representations. It is thus sufficient to show that they are still irreducible and pairwise non-isomorphic in the finite field case.

The two $1$-dimensional representations are non-isomorphic because $\alpha \neq -1$ by the condition on its order.

Let us show that the $2$-dimensional representations are indeed irreducible. Let $j\in [\![1,\frac{m-1}{2}]\!]$ and $V$ the associated $\mathcal{H}_{I_2(m),q}$-module.

Let $W$ be a $\mathcal{H}_{I_2(m),q}$-submodule of $V$ associated to the representation $\rho_j$ and $x=(x_1,x_2)\in W\setminus\{(0,0)\}$.

We have $\rho_j(T_s).x+x=(0,x_1+(\alpha+1)x_2) \in W$, therefore it follows that $(0,1)\in W$ or $x_1+(\alpha+1)x_2=0$.

Assume first that $(0,1)\in W$.

We have $x=x_1(1,0)+x_2(0,1)\in W$, therefore $x_1(1,0)\in W$, which implies that $x_1=0$ or $(1,0)\in W$.

\noindent
If $(1,0)\in W$ then since $(0,1)\in W$, we have $V=W$.

\noindent
If $x_1=0$, then $(1,0)\in W$ or $x_2\alpha(2+\theta^j+\theta^{-j})=0$ since $\rho_j(T_t).x=(\alpha x_1+\alpha(2+\theta^j+\theta^{-j})x_2,-x_2)\in V$. If $(1,0)\in W$, then $V=W$ by the same reasoning as above. We have $2+\theta^j+\theta^{-j}=(1+\theta^j)(1+\theta^{-j})$ and $\theta^{2j}\neq 1$ because $m$ is odd, $j\in [\![1,\frac{m-1}{2}]\!]$ and $\theta$ is a primitive m-th root of unity. It follows that $x_2\alpha(2+\theta^j+\theta^{-j})=0$ implies $x_2=0$ and, therefore $x=(0,0)$, which is absurd because we chose $x\in  W\setminus\{(0,0)\}$.
 
Assume now $x_1+(\alpha+1)x_2=0$.
 
 \noindent
 We have $\rho_j(T_t).x+x=((\alpha+1)x_1+\alpha(2+\theta^j+\theta^{-j})x_2,0)\in W$, therefore $(\alpha+1)x_1+\alpha(2+\theta^j+\theta^{-j})x_2=0$ or $(1,0)\in W$.
 
 \noindent
 If $(\alpha+1)x_1+\alpha(2+\theta^j+\theta^{-j})x_2=0$ then by substituting $x_1$ by $-(\alpha+1)x_2$, we get $x_2(-\alpha^2-2\alpha-1+2\alpha+\alpha(\theta^j+\theta^{-j}))=0$, therefore $x_2=0$ or $\alpha^2-\alpha(\theta^j+\theta^{-j})+1=(\alpha-\theta^j)(\alpha-\theta^{-j})=0$. We cannot have $x_2=0$ because otherwise $x_1=-(\alpha+1)x_2=0$ and $x=(0,0)$. We also cannot have $\alpha=\theta^{\pm j}$ because otherwise the order of $\alpha$ would divide $m$ which contradicts the assumption on $\alpha$.
 
 \noindent
 If $(1,0)\in W$, then $x_2=0$ or $(0,1)\in W$ since $x\in W$. We therefore have $x_2=x_1=0$ or $W=V$. 
 
 This shows that in all cases $W=V$. This proves that $V$ is indeed irreducible. 
 
 It now remains to show that these representations are non-isomorphic.
 
Let  $j$ and $l$ be two integers such that $1\leq j\leq \ell\leq \frac{m-1}{2}$. We have $\tr(\rho_j(T_s)\rho_j(T_t))=\alpha(\theta^j+\theta^{-j})$. This implies that if $\rho_j$ is isomorphic to $\rho_{\ell}$ then $\alpha(\theta^j+\theta^{-j})=\alpha(\theta^{\ell}+\theta^{-\ell})$, therefore $\theta^j+\theta^{-j}=\theta^{\ell}+\theta^{-\ell}$ and, therefore $(\theta^j-\theta^l)(1-\theta^{-\ell-j})=0$. This implies that $j=l$ because $0\leq \ell-j<m$ and $2\leq j+\ell \leq m-1<m$.

\end{proof}

\begin{theo}\label{platypodes}
Let $j\in [\![1,\frac{m-1}{2}]\!]$. If $G=\mathcal{A}_{I_2(m)}$, then
\begin{enumerate}
\item If $\F_{q_j}=\F_p(\alpha,\theta^j+\theta^{-j})=\F_p(\alpha+\alpha^{-1},\theta^j+\theta^{-j})$ then we have $\rho_j(G)=SL_2(q_j)$.
\item If $\F_{q_j}=\F_p(\alpha,\theta^j+\theta^{-j})\neq \F_p(\alpha+\alpha^{-1},\theta^j+\theta^{-j})$ then we have up to conjugation in $GL_2(q_j)$ that $\rho_j(G)\simeq SU_2(q_j^\frac{1}{2})$.
\end{enumerate}
\end{theo}

\begin{proof}

The proof of this result is done in the same way as Lemma 3.5. of \cite{BM}. The proof uses Dickson's theorem (see \cite{HUP}, Theorem 8.27. chapter 2)  and shows that if we prove that $\overline{\rho_j(G)}$ is not abelian by abelian, $\overline{\rho_j(A_{I_2(m})}\notin \{\mathfrak{S}_4,\mathfrak{A}_5\}$ up to isomorphism and that $\alpha+\alpha^{-1}$ and $\theta^{j}+\theta^{-j}$ belong to the field generated by the traces of the elements of $\rho_j(G)$, then we have the result stated in the proposition. We also need to show that in the second case $\rho_j(G)\subset GU_2(q_j^\frac{1}{2})$.

We need to show that the groups considered are not abelian by abelian. In order to do this, we will show that in both cases, we have

\begin{small}
\noindent
 $A=[\rho_j(T_s),\rho_j(T_t)][\rho_j(T_t)^{-1},\rho_j(T_s)]-[\rho_j(T_t)^{-1},\rho_j(T_s)][\rho_j(T_s),\rho_j(T_t)]\neq 0$ and
 
 \noindent
 $B=[\rho_j(T_s),\rho_j(T_t)][\rho_j(T_t)^{-1},\rho_j(T_s)]+[\rho_j(T_t)^{-1},\rho_j(T_s)][\rho_j(T_s),\rho_j(T_t)]\neq 0$, where $[x,y]=xyx^{-1}y^{-1}$.
 \end{small}
 
 This will prove that $\overline{[\rho_j(T_s),\rho_j(T_t)][\rho_j(T_t)^{-1},\rho_j(T_s)]}\neq \overline{[\rho_j(T_t)^{-1},\rho_j(T_s)][\rho_j(T_s),\rho_j(T_t)]}$. If $\overline{\rho_j(G)}$ is abelian by abelian then there exists an abelian normal subgroup $H$ of $\overline{\rho_j(G)}$ such that $\overline{\rho_j(G)}/H$ is abelian. This implies that the derived subgroup of $\overline{\rho_j(G)}$ is included in $H$ and is therefore abelian which contradicts the above inequality. We now prove that $A\neq 0$ and $B\neq 0$.
 
 \noindent
 We have that the entry $A_{1,2}$ of $A$ verifies $A_{1,2}=-\frac{(2+\theta^j+\theta^{-j})^2(\alpha-1)^2(\alpha-\theta^j)(\theta^{-j}-\alpha)}{\alpha^3}$, therefore $A$ is non-zero because we assumed that the order of $\alpha$ does not divide $m$ and that $\theta$ is an $m$-th root of unity.
 
 \noindent
Assume by contradiction that $B=0$. We write $\gamma_j=\theta^j+\theta^{-j}+2$. We have $\gamma_j=(\theta^j+1)(\theta^{-j}+1)\neq 0$.
 
 \noindent
 We then have $B_{1,2}=-\frac{\gamma_j(\alpha^2-1)(-\alpha^2\gamma_j+\alpha\gamma_j^2-2\alpha\gamma_j+2\alpha-\gamma_j)}{\alpha^3}=0$, therefore $-\alpha^2\gamma_j+\alpha\gamma_j^2-2\alpha\gamma_j+2\alpha-\gamma_j=0$.

\noindent
If $p=2$ then $\gamma_j(-\alpha^2+\alpha\gamma_j-1)=0$, therefore $\alpha^2+\alpha(\theta^j+\theta^{-j}+1)=0$ and $(\alpha+\theta^j)(\alpha+\theta^{-j})=0$ This is absurd because $\alpha^m\neq 1$.

\noindent
Assume now $p\neq 2$, we set $a=\frac{\alpha^3 A_{1,2}}{\gamma_j(\alpha^2-1)}$, $b=\alpha^2A_{1,1}$ and $c=-\alpha^2A_{2,2}$.

\noindent
We then have $\frac{1}{4}(b+2(\alpha^2-1)a-c-\alpha a+a)= \alpha(\alpha^2-\alpha+1)$. However $a=b=c=0$, therefore $\alpha^2-\alpha+1=0$. It follows that $\alpha^6-1=(\alpha^3-1)(\alpha+1)(\alpha^2-\alpha+1)=0$, which is impossible by the assumption on the order of $\alpha$.

\smallskip

We need to prove that $\overline{\rho_j(G)}$ contains elements of order different from $1$, $2$, $3$ and $5$. This will show that $\overline{\rho_j(G)}$ is not isomorphic to $\mathfrak{A}_5$. We already have that it is not isomorphic to $\mathfrak{S}_4$ since $\mathfrak{S}_4$ is abelian by abelian.  The eigenvalues of $\rho_j(T_s)$ are $-1$ and $\alpha$, therefore if $\overline{\rho_j(T_s)}^r=\overline{I_2}$, we have $(\alpha^r,(-1)^r)\in \{(1,1),(-1,-1)\}$, therefore $\alpha^{2r}=1$ which implies that $r\notin \{1,2,3,5\}$ by the conditions on $\alpha$.

\smallskip

We now show that $\alpha+\alpha^{-1}$ and $\beta +\beta^{-1}$ belong to the field generated by the traces of elements of $\rho_j(G)$.

\noindent
We have $\tr(\rho_j(T_tT_s^{-1}))=\gamma_j-\alpha-\alpha^{-1}$ and $\tr([\rho_j(T_t),\rho_j(T_s)])-2=\gamma_j(\gamma_j-(\alpha+\alpha^{-1})-2)$. 

\noindent
 We have $\gamma_j-(\alpha+\alpha^{-1})-2=(\theta^j-\alpha)(1-\alpha^{-1}\theta^{-j})\neq 0$, therefore $\gamma_j$ and $\alpha+\alpha^{-1}$ are indeed in the field generated by traces of the elements of $\rho_j(G)$.

It remains to show that if $\F_{q_j}\neq \F_p(\alpha+\alpha^{-1},\theta^j+\theta^{-j})$ then we have up to conjugation in $GL_2(q_j)$ that $\rho_j(G) \subset GU_2(q_j^\frac{1}{2})$.

\noindent
By Lemma $2.4.$ of \cite{BM}, it is sufficient to show that there exists a matrix $P\in GL_2(q_j)$ such that $P\rho_j(T_s)P^{-1}=\epsilon(^t\rho_j(T_s)^{-1})$ and $P\rho_j(T_t)P^{-1}=\epsilon(^t\rho_j(T_t)^{-1})$, where $\epsilon$ is the unique automorphism of order $2$ of $\F_{q_j}$. Since $\gamma_j\in \F_{q_j^\frac{1}{2}}$, we have $\epsilon(\gamma_j)=\gamma_j$. We also have $(X-\alpha)(X-\alpha^{-1})=X^2-(\alpha+\alpha^{-1})X+1\in \F_{q_j^\frac{1}{2}}$ and $\alpha\notin \F_{q_j^\frac{1}{2}}$, therefore we have $\epsilon(\alpha)=\alpha^{-1}$.

\noindent
Set $P=\begin{pmatrix}
\frac{\alpha+1}{\gamma_j} & 1 \\
\alpha & \alpha+1
\end{pmatrix}$. We have $det(P) =\frac{(\alpha-\theta^j)(\alpha-\theta^{-j})}{\gamma_j}\neq 0$. The matrix $P$ verifies the desired property and the result follows.\end{proof}

We now provide a field-theoretic lemma to see when the representations are linked by composition with a field automorphism. This is necessary to determine the image of $\mathcal{A}_{I_2(m)}$ inside the full Iwahori-Hecke algebra.

\begin{lemme}\label{Isomorphism}
Let $j,l\in[\![1,\frac{m-1}{2}]\!]^2$. There exists an automorphism $\Psi_{l,j}$ of $\F_{q_j}=\F_p(\alpha,\xi^j+\xi^{-j})$ verifying $\Psi_{l,j}(\alpha+\alpha^{-1})=\alpha+\alpha^{-1}$ and $\Psi_{l,j}(\xi^j+\xi^{-j})=\xi^l+\xi^{-l}$ if and only if there exists $r\in \N$ such that $jp^r\equiv l~(\bmod ~ m)$ or $jp^r\equiv -l~(\bmod~ m)$ and $(\alpha+\alpha^{-1})^{p^r}=\alpha+\alpha^{-1}$.

We say that $j\sim l$ if one of those conditions is verified. This is an equivalence relation and when $j\sim l$, we have $\rho_{l|\mathcal{A}_{I_2(m)}}=\Psi_{l,j}\circ \rho_{j|\mathcal{A}_{I_2(m)}}$.
\end{lemme}

\begin{proof}
Assume there exists $r\in \N$ such that $jp^r\equiv l~(\bmod ~ m)$ or $jp^r\equiv -l~(\bmod ~ m)$ and $(\alpha+\alpha^{-1})^{p^r}=\alpha+\alpha^{-1}$. Let $\Psi$ be the automorphism of $\F_{q_j}$ defined by $\Psi(x)=x^{p^r}$ for all $x\in \F_{q_j}$. We then have $\Psi(\alpha+\alpha^{-1})=\alpha+\alpha^{-1}$ by assumption and $\Psi(\xi^j+\xi^{-j})=(\xi^j+\xi^{-j})^{p^r}=\xi^{jp^r}+\xi^{-jp^r}=\xi^l+\xi^{-l}$.

Assume now that there exists an automorphism $\Psi$ of $\F_{q_j}=\F_p(\alpha,\xi^j+\xi^{-j})$ verifying $\Psi(\alpha+\alpha^{-1})=\alpha+\alpha^{-1}$ and $\Psi(\xi^j+\xi^{-j})=\xi^l+\xi^{-l}$.

There exists $r\in [\![0,\log_p(q)-1]\!]$ such that $\Psi(x)=x^{p^r}$ for all $r\in \N$, therefore $(\alpha+\alpha^{-1})^{p^r}=\Psi(\alpha+\alpha^{-1})=\alpha+\alpha^{-1}$. The map $\Psi$ can be extended to an automorphism $\tilde{\Psi}$ of $\overline{\F_p}$ by defining $\tilde{\Psi}$ to be the automorphism sending $x$ to $x^{p^r}$ for all $x\in \overline{\F_p}$. We then have $\xi^l+\xi^{-l}=\Psi(\xi^j+\xi^{-j})=\tilde{\Psi}(\xi^j+\xi^{-j})=\tilde{\Psi}(\xi^j)+\tilde{\Psi}(\xi^{-j})$. It follows that $(\tilde{\Psi}(\xi^j)\xi^l-1)(\xi^{-l}-\tilde{\Psi}(\xi^{-j}))=0$, therefore $\tilde{\Psi}(\xi^j)\in \{\xi^l,\xi^{-l}\}$. This proves that $\xi^{jp^r}\in \{\xi^l,\xi^{-l}\}$, therefore $jp^r\equiv l~(\bmod ~ m)$ or $jp^r\equiv -l~(\bmod ~ m)$.

The fact that $\sim$ is an equivalence relation follows from the fact that for all $r\in \N$, $\op{Gcd}(m,p^r)=1$. 
\end{proof}

\begin{theo}\label{resdihedral}
Assume $m$ odd and $\alpha$ satisfies the conditions given at the beginning of this section. For $j \in [\![1,\frac{m-1}{2}]\!]$, we set $G_j=SL_2(q_j)$ if $\F_{q_j}=\F_p(\alpha,\theta^j+\theta^{-j})=\F_p(\alpha+\alpha^{-1},\theta^j+\theta^{-j})$ and $G_j=SU_2(q_j^\frac{1}{2})$ if $\F_{q_j}=\F_p(\alpha,\theta^j+\theta^{-j})\neq\F_p(\alpha+\alpha^{-1},\theta^j+\theta^{-j})$.

We then have that the morphism from $\mathcal{A}_{I_2(m)}$ to $\mathcal{H}_{I_2(m),\alpha}^\times\simeq GL_1(q_j)^2 \times \underset{j\in [\![1,\frac{m-1}{2}]\!]}\prod GL_2(q_j)$ factorizes through the surjective morphism :
$$\Phi : \mathcal{A}_{I_2(m)} \rightarrow \underset{j\in [\![1,\frac{m-1}{2}]\!]/\sim}\prod G_j.$$
\end{theo}

\begin{proof}

We know by Theorem \ref{platypodes} that the composition of the morphism from $\mathcal{A}_{I_2(m)}$ to $\mathcal{H}_{I_2(m),\alpha}^\times$ with the projection upon each representation is surjective. We also know by lemma \ref{Isomorphism} that the morphism to $\mathcal{H}_{I_2(m),\alpha}^\times$ factorizes through the morphism $\Phi$.
We will now use Goursat's lemma and induction on $j\in [\![1,\frac{m-1}{2}]\!]$ in order to conclude the proof of this theorem.

\noindent
For $j_0\in [\![1,\frac{m-1}{2}]\!]$. We define $\Phi_{j_0}(\mathcal{A}_{I_2(m)})$ to be the image of $\mathcal{A}_{I_2(m)}$ inside $\underset{j\in [\![1,j_0]\!]/\sim}\prod GL_2(q_j)$.

We know that $\Phi_1(\mathcal{A}_{I_2(m)})=G_1$. Let $j_0\in [\![1,\frac{m-3}{2}]\!]$, assume that $\Phi_{j_0}(\mathcal{A}_{I_2(m)})=\underset{j\in [\![1,j_0]\!]/\sim}\prod G_{j}$.

\noindent
Consider $\Phi_{j_0+l}(\mathcal{A}_{I_2(m)})\subset \underset{j\in [\![1,j_0]\!]}\prod G_{j} \times G_{j_l}$ for $\ell$ the smallest positive integer such that $j_\ell\nsim j$ for all $j\in [\![1,j_0]\!]$. We know that the projection upon each factor is surjective. Let $K_1= \underset{j\in [\![1,j_0]\!]}\prod G_{j}$ and $K_2=G_{j_0+\ell}$ as in Goursat's Lemma.
 We then have $K_1/K^1\simeq K_2/K^2$. If the quotients are abelian then we are done since both groups are perfect. Assume that those quotients are non-abelian. There is only one non-abelian decomposition factor of $K_2$ and it is isomorphic to $PSL_2(q_{j_0+\ell})$ or $PSU_2(q_{j_0+\ell}^\frac{1}{2})$ depending on the field $\F_{q_{j_0+\ell}}$. We write that decomposition factor $PG_{j_0+\ell}$. The isomorphism then implies that there exists $j_1\in [\![1,j_0]\!]$ such that $\overline{\rho_{j_1}(\mathcal{A}_{I_2(m})}\simeq PG_{j_1} \simeq \overline{\rho_{j_0+\ell}(\mathcal{A}_{I_2(m})}$. We then have that there exists $z:\mathcal{A}_{I_2(m)} \rightarrow \overline{\F_p}^\times$ and $\Psi\in Aut(\F_{q_{j_0+\ell}})$ such that up to conjugation, for all $h\in \mathcal{A}_{I_2(m)}, \rho_{j_0+\ell}(h)=\Psi(\rho_{j_1}(h))z(h)$. We will prove this is absurd by considering traces of some elements in $\mathcal{A}_{I_2(m)}$ under these representations. We may first note that for all $M\in SL_2(\overline{\F_q})$, we have $1=\op{det}(z(h)M)=z(h)^2\op{det}(M)=z(h)^2$, therefore for all $h \in \mathcal{A}_{I_2(m)}$, $z(h)\in \{-1,1\}$. We write as before in the sequel $\gamma_{j_0+\ell}=\theta^{j_0+\ell}+\theta^{-(j_0+\ell)}+2$ and $\gamma_{j_1}=\theta^{j_1}+\theta^{-j_1}+2$.

\bigskip

Assume first that $z(T_tT_s^{-1})=z([T_t,T_s])=1$.
We then consider the traces of those two elements inside each representation, $\tr(\rho_{j_0+\ell}(T_tT_s^{-1})==\Psi(\tr(\rho_{j_1}(T_tT_s^{-1})))$ and $\tr(\rho_{j_0+\ell}([T_t,T_s]))=\Psi(\tr(\rho_{j_1}([T_t,T_s])))$. This implies that
 $$-(\alpha+\alpha^{-1})+\gamma_{j_0+\ell}=-\Psi(\alpha+\alpha^{-1})+\Psi(\gamma_{j_1})$$ 
  $$\gamma_{j_0+\ell}(-(\alpha+\alpha^{-1})+\gamma_{j_0+l}-2)+2=\Psi(\gamma_{j_1}(-(\alpha+\alpha^{-1})+\gamma_{j_1}-2)+2).$$ Since $\alpha+\alpha^{-1}+\gamma_{j_0+\ell}\neq 0$, it follows that $\Psi(\theta^{j_1}+\theta^{-j_1})=\gamma_{j_0+\ell}$ and then $\Psi(\alpha+\alpha^{-1})=\alpha+\alpha^{-1}$. This implies $j_1\sim j_0+\ell$, which is absurd by assumption.
 
\medskip 
 
Assume now $z(T_tT_s^{-1})=-z([T_t,T_s])=1$. By considering the traces of the same elements, we get
\begin{small}
$$A=-(\alpha+\alpha^{-1})+\gamma_{j_0+\ell}=-\Psi(\alpha+\alpha^{-1})+\Psi(\gamma_{j_1})=\Psi(A_1)$$ 
  $$B=\gamma_{j_0+\ell}(-(\alpha+\alpha^{-1})+\gamma_{j_0+\ell}-2)+2=-\Psi(\gamma_{j_1}(-(\alpha+\alpha^{-1})+\gamma_{j_1}-2)+2)=-\Psi(B_1).$$ 
  \end{small}
  We have $z(T_tT_s^{-1}[T_t,T_s])=-1$, therefore
   $$C=\tr(\rho_{j_0+\ell}(T_tT_1^{-1}[T_t,T_s]))=-\Psi(\tr(\rho_{j_1}(T_tT_s^{-1}[T_t,T_s])))-\Psi(C_1)$$ but $C=AB-A$ and $C_1=A_1B_1-A_1$, therefore we have $AB-A=C=-\Psi(C_1)=-\Psi(A_1B_1-A_1)=AB+A$ and it follows that $A=0$, which is absurd as was proven before.

\medskip
In the same way as above if $z(T_tT_s^{-1})=z([T_t,T_s])=1$, we have $A=-\Psi(A_1), B=-\Psi(B_1)$ and $AB-A=\Psi(A_1B_1-A_1)=A_1B_1+A_1$, therefore $A=A_1=0$, which is absurd.

\medskip
The last case remaining is $-z(T_tT_s^{-1})=z([T_t,T_s])=1$. We then have 
$$-(\alpha+\alpha^{-1})+\gamma_{j_0+\ell}=\Psi(\alpha+\alpha^{-1})-\Psi(\gamma_{j_1})$$ 
  $$\gamma_{j_0+\ell}(-(\alpha+\alpha^{-1})+\gamma_{j_0+\ell}-2)+2=\Psi(\gamma_{j_1}(-(\alpha+\alpha^{-1})+\gamma_{j_1}-2)+2).$$
  so
  $$\gamma_{j_0+\ell}(-(\alpha+\alpha^{-1})+\gamma_{j_0+\ell}-2)=\Psi(\gamma_{j_1})(\alpha+\alpha^{-1}-\gamma_{j_0+l}-2).$$
  
 We  consider a third element in this case and distinguish two possibilities.
 
 \smallskip
 
  Assume first that $z([T_t^2,T_s])=1$. We then have that $z([T_t,T_s][T_t^2,T_s]T_tT_s^{-1})=-1$, therefore $E=\tr(\rho_{j_0+\ell}([T_t,T_s][T_t^2,T_s]T_tT_s^{-1}))=-\Psi(\tr(\rho_{j_1}([T_t,T_s][T_t^2,T_s]T_tT_s^{-1})))=-\Psi(E_1)$ and $D=\tr(\rho_{j_0+\ell}([T_t^2,T_s]))=\Psi(\tr(\rho_{j_1}([T_t^2,T_s])))=\Psi(D_1)$.
 
\noindent
 We also have $ABD-AD-AB+D+A-2=E=-\Psi(A_1B_1D_1)+\Psi(A_1D_1)+\Psi(A_1B_1)+\Psi(D_1)-\Psi(A_1)-2$ and $-\Psi(E_1)=-\Psi(A_1B_1D_1-A_1D_1-A_1B_1+D_1+A_1-2)$. It follows that $2(\Psi(D_1)-2)=0$, therefore $\Psi(D_1-2)=0$. We then have $0=D_1-2=\frac{\gamma_{j_1}(\alpha-1)^2(\alpha^2-\alpha\gamma_{j_1}+2\alpha+1}{\alpha^2}$, therefore $\alpha+\alpha^{-1}=\gamma_{j_1}-2=\theta^{j_1}+\theta^{-j_1}$. This would imply that $\alpha\in\{\theta^{j_1},\theta^{-j_1}\}$ and, therefore $\alpha^m=1$ which contradicts our assumptions on $\alpha$.
 
 \smallskip Assume now that $z([T_t^2,T_s])=-1$. We then have $z([T_t^2,T_s][T_t,T_s])=-1$, therefore $F=\tr(\rho_{j_0+\ell}([T_t^2,T_s][T_t,T_s]))=-\Psi(\tr(\rho_{j_1}([T_t^2,T_s][T_t,T_s])))=-\Psi(F_1)$. We also have $F=BD-B=-\Psi(B_1D_1)-\Psi(B_1)$ and $-\Psi(F_1)=-\Psi(B_1D_1-B_1)=-\Psi(B_1D_1)+\Psi(B_1)$. It follows that $\Psi(B_1)=0$, therefore $B_1=0$ and $B=\Psi(B_1)=0$.
 
 \noindent We also have $z([T_t,T_s][T_t^2,T_s](T_t T_s^{-1}))=1$ and, therefore $E=\Psi(E_1)$. We have $E=ABD-AD-AB+D+A-2=-AD+D+A-2=-\Psi(A_1D_1)-\Psi(D_1)-\Psi(A_1)-2$ and $\Psi(E_1)=\Psi(A_1B_1D_1-A_1D_1-A_1B_1+D_1+A_1-2)=-\Psi(A_1D_1)+\Psi(D_1)+\Psi(A_1)-2$. It follows that $\Psi(A_1+D_1)=0$, therefore $A_1+D_1=0$. We then have $0=\alpha^{-1}(\alpha-1)^2B_1+(A_1+D_1)=\alpha-2+\alpha^{-1}+\gamma_{j_1}=0$, therefore $\alpha+\alpha^{-1}=2-\gamma_{j_1}$. We then have $0=B_2=-(\alpha+\alpha^{-1})\gamma_{j_1}+\gamma_{j_1}^2-2\gamma_{j_1}+2=-(2-\gamma_{j_1})\gamma_{j_1}+\gamma_{j_1}^2-2\gamma_{j_1}+2=2\gamma_{j_1}^2-4\gamma_{j_1}^2+2=2(\gamma_{j_1}-1)^2$. It follows that $\gamma_{j_1}=1$, therefore $\alpha+\alpha^{-1}=2-1=1$ and $A_1=-\alpha-\alpha^{-1}+\gamma_{j_1}=0$. We then have $A=-\Psi(A_1)=0$, therefore $-\alpha-\alpha^{-1}+\gamma_{j_0+\ell}=0$ and $\gamma_{j_0+\ell}=\alpha+\alpha^{-1}=\gamma_{j_1}$. This implies that $j_0+\ell=j_1$, which is absurd. This concludes the proof.
\end{proof}

\section{$m$ even}\label{sectionmeven}
 
Let $m\geq 5$ even, $p$ a prime, $\alpha, \beta \in \overline{\F_p}$ of orders not belonging to $\{1,2,3,4,5,6,10\}$ and $\theta\in \overline{\F_p}$ a primitive $m$-th root of unity. Note that $p\neq 2$ because there exists a primitive $m-th$ root of unity with $m$ even. We assume that $j\in [\![1,\frac{m-2}{2}]\!]$, we have $(\alpha+\sqrt{\alpha\beta}(\theta^j+\theta^{-j})+\beta)(\alpha\beta-\sqrt{\alpha\beta}(\theta^j+\theta^{-j})+1)\neq 0$. For $j\in \N$, we write $\F_{q_j}=\F_p(\alpha,\beta,\sqrt{\alpha\beta}(\theta^j+\theta^{-j}))$. We write $\F_q$ the smallest field containing $\F_{q_j}$ for $j$ in $[\![1,\frac{m-2}{2}]\!]$. Note that we have $\F_q=\F_{q_1}$ as in section \ref{sectionmodd}.

 \begin{Def}
The Iwahori-Hecke algebra of dihedral type $I_2(m)$ which we write $\mathcal{H}_{I_2(m),q}$ is the $\F_q$-algebra with the following presentation
 
 Generators : $T_t,T_s$.
 
 Relations : 
 
$(T_s-\alpha)(T_s+1)=0$,
 
$(T_t-\beta)(T_t+1)=0$,
 
$\underset{m}{\underbrace{T_sT_tT_s\cdots}}=\underset{m}{\underbrace{T_tT_sT_t\cdots}}$.
 \end{Def}
 
We then give the Kilmoyer-Solomon matrix model given in \cite{G-P} (Theorem $8.3.1$) in the finite field case.
 
\begin{theo}
Under the assumptions made on $\alpha$, $\beta$ and $\theta$, the following matrix model gives us a decomposition into pairwise non-isomorphic irreducible $\mathcal{H}_{I_2(m),q}$-modules 
\begin{enumerate}
\item $Ind : T_s \mapsto \alpha, T_t\mapsto \beta$.
\item $\epsilon : T_s \mapsto -1, T_t \mapsto -1$.
\item For $j\in [\![1,\frac{m-2}{2}]\!], T_s \mapsto \rho_j(T_s)=\begin{pmatrix}
-1 & 0\\
1 & \alpha
\end{pmatrix}, T_t \mapsto \rho_j(T_t)=\begin{pmatrix}
\beta & \alpha+\sqrt{\alpha\beta}(\theta^j+\theta^{-j})+\beta\\
0 & -1
\end{pmatrix}$.
\end{enumerate}
\end{theo}

\begin{proof}

We know that these models give us representations of $\mathcal{H}_{I_2(m),q}$. It is thus sufficient to show that they are irreducible and pairwise non-isomorphic. The two $1$-dimensional representations are non isomorphic by the conditions on $\alpha$ and $\beta$. Let us show that the $2$-dimensional representations are irreducible. Let $j\in [\![1,\frac{m-1}{2}]\!]$. Let $W$ be a non-trivial $\mathcal{H}_{I_2(m),q}$-submodule of the module $V$ associated to the representation $\rho_j$. Let $x=(x_1,x_2)\in W\setminus \{(0,0)\}$. We have that $\rho_j(T_s)\cdot x+x=(0,x_1+(\alpha+1)x_2) \in W$, therefore $(0,1)\in W$ or $x_1+(\alpha+1)x_2=0$.

\smallskip

Assume first that $(0,1)\in W$. Since $x\in W$, we have that $x_1=0$ or $(1,0)\in W$. 

If $(1,0)\in W$ then $W=V$.

Assume now $x_1=0$. We then have that $\rho_j(T_t)\cdot x+x=((\beta+1)x_1+(\alpha+\sqrt{\alpha\beta}(\theta^j+\theta^{-j})+\beta)x_2,0)\in W$. It follows that $(1,0)\in W$ or $(\alpha+\sqrt{\alpha\beta}(\theta^j+\theta^{-j})+\beta)x_2=0$. By the assumptions on $\alpha$, $\beta$ and $\theta$, we have that $\alpha+\sqrt{\alpha\beta}(\theta^j+\theta^{-j})+\beta\neq 0$, therefore the latter possibility would imply $x=(0,0)$, which contradicts our assumptions. It follows that $(1,0)\in W$ and $W=V$.

\smallskip

Assume now that $x_1+(\alpha+1)x_2=0$. We consider again $\rho_j(T_t)\cdot x+x$ and we get that $(1,0)\in W$ or $(\beta+1)x_1+(\alpha+\sqrt{\alpha\beta}(\theta^j+\theta^{-j})+\beta)x_2=0$. 

Assume $(1,0)\in W$. Since $x\in W$, we then have $(0,1)\in W$ or $x_2=0$. The latter would imply that $x_2=x_1=0$ since $x_1+(\alpha+1)x_2=0$. It follows that $(0,1)\in W$ and, therefore $V=W$. 

Assume now by contradiction that $(\beta+1)x_1+(\alpha+\sqrt{\alpha\beta}(\theta^j+\theta^{-j})+\beta)x_2=0$. We have $x_1=-(\alpha+1)x_2$, it follows that $(-\alpha\beta+\sqrt{\alpha\beta}(\theta^j+\theta^{-j})-1)x_2=0$, therefore $x_2=0$ and $x_1=x_2=0$. This is absurd, therefore we are in the first case and $V=W$.

This proves that in all cases, we have that $W=V$ and $\rho_j$ is therefore irreducible.

\medskip

It remains to show that those representations are pairwise non-isomorphic. Let $j$ and $\ell$ be two integers such that $1\leq j\leq \ell\leq \frac{m-2}{2}$. We have $\tr(\rho_j(T_t)\rho_j(T_t))=\sqrt{\alpha\beta}(\theta^j+\theta^{-j})$ and $\tr(\rho_\ell(T_t)\rho_\ell(T_t))=\sqrt{\alpha\beta}(\theta^\ell+\theta^{-\ell})$. Assume now that those representations are isomorphic. We then have that $\sqrt{\alpha\beta}(\theta^j+\theta^{-j})=\sqrt{\alpha\beta}(\theta^\ell+\theta^{-\ell})$. It follows that $\theta^j+\theta^{-j}=\theta^\ell+\theta^{-\ell}$ and, therefore $(\theta^j-\theta^\ell)(1-\theta^{-\ell-j})=0$. This implies that $j=\ell$ since $0\leq \ell-j\leq m$ and $2\leq j+\ell\leq m-2<m$. \end{proof}

The main difference between $m$ even and $m$ odd arises in the field extensions we have to consider. We describe the different cases we encounter in what follows. First note that for $j\in [\![1,\frac{m-2}{2}]\!]$, if we set $P_j=\begin{pmatrix} 0 & \alpha+\sqrt{\alpha\beta}(\theta^j+\theta^{-j})\\
1 & 0\end{pmatrix}$ then by the assumptions made on $\alpha$, $\beta$ and $\theta$, the matrix $P_j$ is invertible. Moreover, we have that $P_j\rho_j(T_s)P_j^{-1}=\begin{pmatrix}
 \alpha & \alpha+\sqrt{\alpha\beta}(\theta^j+\theta^{-j})+\beta\\
 0 & -1
\end{pmatrix}$ and $P_j\rho_j(T_t)P_j^{-1}=\begin{pmatrix}
 -1 & 0\\
 1 & \beta
 \end{pmatrix}$. This shows that the roles of $\alpha$ and $\beta$ are perfectly symmetric.
We now fix a $j\in [\![1,\frac{m-2}{2}]\!]$ and we write $\gamma=\sqrt{\alpha\beta}$, $u_j=\theta^j+\theta^{-j}$ and $\xi=\frac{\gamma u_j(\alpha\beta-\alpha-\beta+1)}{\alpha\beta}$. We set $L=\F_p(\alpha,\beta,\gamma u_j)$, $L_1=\F_p(\alpha+\alpha^{-1},\beta,\gamma u_j)$, $L_2=\F_p(\alpha,\beta,\xi)$, $L_3=\F_p(\alpha,\beta+\beta^{-1},\xi)$, $K_1=\F_p(\alpha+\alpha^{-1},\beta,\xi)$, $K_2=\F_p(\alpha+\alpha^{-1},\beta+\beta^{-1},\gamma u_j)$, $K_3=\F_p(\alpha,\beta+\beta^{-1},\theta)$ and $K=\F_p(\alpha+\alpha^{-1},\beta+\beta^{-1},\xi)$.
 
Note that $L=L_2$ since $\alpha\beta-\alpha-\beta+1=(\alpha-1)(\beta-1)\neq 0$ and $\frac{\alpha\beta\xi}{\alpha\beta-\alpha-\beta+1}=\gamma u_j$. If we consider the polynomials $X^2-(\alpha+\alpha^{-1})X+1$ and $X^2-(\beta+\beta^{-1})X+1$, we get that $([L:L_1],[L:L_3],[L_2,K_1],[L_2,K_3],[K_1:K],[K_3:K])\in \{1,2\}^6$. We then have the following Hasse diagram

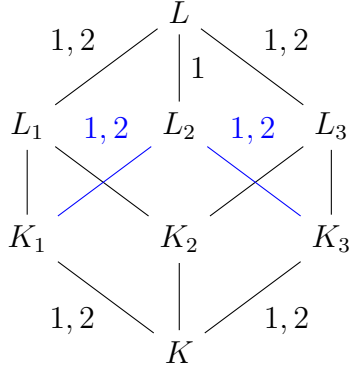
\begin{center}
\begin{tikzpicture}
\node (1) at (0,3){$L$};
\node (2) at (-2,1.5){$L_1$};
\node (3) at (0,1.5){$L_2$};
\node (4) at (2,1.5){$L_3$};
\node (5) at (-2,0){$K_1$};
\node (6) at (0,0){$K_2$};
\node (7) at (2,0){$K_3$};
\node (8) at (0,-1.5){$K$};
\draw (1) to node[auto]{$1$} (3);
\draw (1) to node[auto,swap]{$1,2$} (2);
\draw (1) to node[auto]{$1,2$} (4);
\draw (2) to (5);
\draw (2) to (6);
\draw (4) to (6);
\draw (4) to (7);
\draw[blue] (3) to node[auto,swap,very near start]{$1,2$} (5);
\draw[blue] (3) to node[auto,very near start]{$1,2$} (7);
\draw (6) to (8);
\draw (5) to node[auto,swap]{$1,2$} (8);
\draw (7) to node[auto]{$1,2$} (8);
\end{tikzpicture}
\end{center}

We then see that $[L:K]=[L:L_2][L_2:K_1][K_1:K]\in \{1,2,4\}$. Assume that $[L:K]=4$, we then have that $[L_2:K_1]=[L_2:K_3]=2$. By unicity of the subfields of degree $2$, we then have $K_1=K_3$ and, therefore $\alpha\in K_1$. It follows that $K_1=L_2$ and, therefore $[K_1:K_1]=2$ which is a contradiction. This proves that $[L:K]\in \{1,2\}$. By uniqueness of the subfields of a given degree we cannot have $[L:L_1]=[L:L_3]=2$ or $[L_1:K_1]=[L_1:K_2]=2$ or $[L_2:K_1]=[L_2:K_3]=2$ or $[L_3:K_3]=[L_3:K_2]=2$. It follows that the possible ways $L$ is an extension of $K$ corresponds to Hasse diagrams described in Figures \ref{fieldsIcase1} to \ref{fieldsIcase7}. We write in dashed red lines the extensions of degree $2$, in red the subfields of degree $2$ of $L$, in dotted black lines the extensions of degree $1$ and in black the fields equal to $L$.

\begin{figure}
\centering
\begin{tikzpicture}
\node (1) at (0,3){$L$};
\node (2) at (-2,1.5){$L_1$};
\node (3) at (0,1.5){$L_2$};
\node (4) at (2,1.5){$L_3$};
\node (5) at (-2,0){$K_1$};
\node (6) at (0,0){$K_2$};
\node (7) at (2,0){$K_3$};
\node (8) at (0,-1.5){$K$};
\draw[dotted] (1) to  (3);
\draw[dotted] (1) to  (2);
\draw[dotted] (1) to  (4);
\draw[dotted] (2) to (5);
\draw[dotted](2) to (6);
\draw[dotted] (4) to (6);
\draw[dotted](4) to (7);
\draw[dotted] (3) to  (5);
\draw[dotted] (3) to  (7);
\draw[dotted] (6) to (8);
\draw[dotted] (5) to  (8);
\draw[dotted] (7) to  (8);
\end{tikzpicture}
\caption{Field extensions in type $I_2(m)$, $m$ even in case $1$}\label{fieldsIcase1}
\end{figure}

\begin{figure}
\centering
\begin{tikzpicture}
\node (1) at (0,3){$L$};
\node (2) at (-2,1.5){$L_1$};
\node (3) at (0,1.5){$L_2$};
\node (4) at (2,1.5){$L_3$};
\node (5) at (-2,0){$K_1$};
\node (6) at (0,0){$K_2$};
\node (7) at (2,0){$K_3$};
\node[red] (8) at (0,-1.5){$K$};
\draw[dotted] (1) to  (3);
\draw[dotted] (1) to  (2);
\draw[dotted] (1) to  (4);
\draw[dotted] (2) to (5);
\draw[dotted](2) to (6);
\draw[dotted] (4) to (6);
\draw[dotted](4) to (7);
\draw[dotted] (3) to  (5);
\draw[dotted] (3) to  (7);
\draw[dashed][red] (6) to (8);
\draw[dashed][red] (5) to  (8);
\draw[dashed][red] (7) to  (8);
\end{tikzpicture}
\caption{Field extensions in type $I_2(m)$, $m$ even in case $2$}\label{fieldsIcase2}
\end{figure}

\begin{figure}
\centering
\begin{tikzpicture}
\node (1) at (0,3){$L$};
\node (2) at (-2,1.5){$L_1$};
\node (3) at (0,1.5){$L_2$};
\node (4) at (2,1.5){$L_3$};
\node (5) at (-2,0){$K_1$};
\node[red] (6) at (0,0){$K_2$};
\node (7) at (2,0){$K_3$};
\node[red] (8) at (0,-1.5){$K$};
\draw[dotted] (1) to  (3);
\draw[dotted] (1) to  (2);
\draw[dotted] (1) to  (4);
\draw[dotted] (2) to (5);
\draw[dashed][red] (2) to (6);
\draw[dashed][red] (4) to (6);
\draw[dotted](4) to (7);
\draw[dotted] (3) to  (5);
\draw[dotted] (3) to  (7);
\draw[dotted] (6) to (8);
\draw[dashed][red] (5) to  (8);
\draw[dashed][red] (7) to  (8);
\end{tikzpicture}
\caption{Field extensions in type $I_2(m)$, $m$ even in case $3$}\label{fieldsIcase3}
\end{figure}

\begin{figure}
\centering
\begin{tikzpicture}
\node (1) at (0,3){$L$};
\node[red] (2) at (-2,1.5){$L_1$};
\node (3) at (0,1.5){$L_2$};
\node (4) at (2,1.5){$L_3$};
\node[red] (5) at (-2,0){$K_1$};
\node[red] (6) at (0,0){$K_2$};
\node (7) at (2,0){$K_3$};
\node[red] (8) at (0,-1.5){$K$};
\draw[dotted] (1) to  (3);
\draw[dashed][red] (1) to  (2);
\draw[dotted] (1) to  (4);
\draw[dotted] (2) to (5);
\draw[dotted](2) to (6);
\draw[dashed][red] (4) to (6);
\draw[dotted](4) to (7);
\draw[dashed][red] (3) to  (5);
\draw[dotted] (3) to  (7);
\draw[dotted] (6) to (8);
\draw[dotted] (5) to  (8);
\draw[dashed][red] (7) to  (8);
\end{tikzpicture}
\caption{Field extensions in type $I_2(m)$, $m$ even in case $4$}\label{fieldsIcase4}
\end{figure}

\begin{figure}
\centering
\begin{tikzpicture}
\node (1) at (0,3){$L$};
\node (2) at (-2,1.5){$L_1$};
\node (3) at (0,1.5){$L_2$};
\node[red] (4) at (2,1.5){$L_3$};
\node (5) at (-2,0){$K_1$};
\node[red] (6) at (0,0){$K_2$};
\node[red] (7) at (2,0){$K_3$};
\node[red] (8) at (0,-1.5){$K$};
\draw[dotted] (1) to  (3);
\draw[dotted] (1) to  (2);
\draw[dashed][red] (1) to  (4);
\draw[dotted] (2) to (5);
\draw[dashed][red](2) to (6);
\draw[dotted] (4) to (6);
\draw[dotted](4) to (7);
\draw[dotted] (3) to  (5);
\draw[dashed][red] (3) to  (7);
\draw[dotted] (6) to (8);
\draw[dashed][red] (5) to  (8);
\draw[dotted] (7) to  (8);
\end{tikzpicture}
\caption{Field extensions in type $I_2(m)$, $m$ even in case $5$}\label{fieldsIcase5}
\end{figure}

\begin{figure}
\centering
\begin{tikzpicture}
\node (1) at (0,3){$L$};
\node (2) at (-2,1.5){$L_1$};
\node (3) at (0,1.5){$L_2$};
\node (4) at (2,1.5){$L_3$};
\node (5) at (-2,0){$K_1$};
\node (6) at (0,0){$K_2$};
\node[red] (7) at (2,0){$K_3$};
\node[red] (8) at (0,-1.5){$K$};
\draw[dotted] (1) to  (3);
\draw[dotted] (1) to  (2);
\draw[dotted] (1) to  (4);
\draw[dotted] (2) to (5);
\draw[dotted](2) to (6);
\draw[dotted] (4) to (6);
\draw[dashed][red](4) to (7);
\draw[dotted] (3) to  (5);
\draw[dashed][red] (3) to  (7);
\draw[dashed][red] (6) to (8);
\draw[dashed][red] (5) to  (8);
\draw[dotted] (7) to  (8);
\end{tikzpicture}
\caption{Field extensions in type $I_2(m)$, $m$ even in case $6$}\label{fieldsIcase6}
\end{figure}

\begin{figure}
\centering
\begin{tikzpicture}
\node (1) at (0,3){$L$};
\node (2) at (-2,1.5){$L_1$};
\node (3) at (0,1.5){$L_2$};
\node (4) at (2,1.5){$L_3$};
\node[red] (5) at (-2,0){$K_1$};
\node (6) at (0,0){$K_2$};
\node (7) at (2,0){$K_3$};
\node[red] (8) at (0,-1.5){$K$};
\draw[dotted] (1) to  (3);
\draw[dotted] (1) to  (2);
\draw[dotted] (1) to  (4);
\draw[dashed][red] (2) to (5);
\draw[dotted](2) to (6);
\draw[dotted] (4) to (6);
\draw[dotted](4) to (7);
\draw[dashed][red] (3) to  (5);
\draw[dotted] (3) to  (7);
\draw[dashed][red] (6) to (8);
\draw[dotted] (5) to  (8);
\draw[dashed][red] (7) to  (8);
\end{tikzpicture}
\caption{Field extensions in type $I_2(m)$, $m$ even in case $7$}\label{fieldsIcase7}
\end{figure}

\begin{theo}\label{bopb}
For $j\in [\![1,\frac{m-2}{2}]\!]$, let $G=[<T_t,T_s>,<T_t,T_s>]$. We then have 
\begin{enumerate}
\item In case $1$, $\rho_j(G)=SL_2(q_j)$.
\item In case $2$, $\rho_j(G)\simeq SU_2(q_j^\frac{1}{2})$.
\item In case $3$, $\rho_j(G)\simeq SU_2(q_j^\frac{1}{2})$.
\item In cases $4$ and $5$, $\rho_j(G)\simeq SL_2(q_j^\frac{1}{2})$.
\item In cases $6$ and $7$, $\rho_j(G)\simeq SU_2(q_j^\frac{1}{2})$.
\end{enumerate}
\end{theo}

\begin{proof}

We write $\epsilon$ the unique automorphism of order $2$ of $\F_{q_j}$ when it exists. By symmetry of the roles of $\alpha$ and $\beta$, it is sufficient to consider cases $1$, $2$, $3$, $4$ and $6$. As for the proof of Theorem \ref{platypodes}, it is sufficient to show that
\begin{enumerate}
\item $A=[\rho_j(T_s),\rho_j(T_t)][\rho_j(T_t)^{-1},\rho_j(T_s)]-[\rho_j(T_t)^{-1},\rho_j(T_s)][\rho_j(T_s),\rho_j(T_t)]\neq 0$,
\item  $B=[\rho_j(T_s),\rho_j(T_t)][\rho_j(T_t)^{-1},\rho_j(T_s)]+[\rho_j(T_t)^{-1},\rho_j(T_s)][\rho_j(T_s),\rho_j(T_t)]\neq 0$,
\item $\overline{\rho_j(G)}$ contains elements of order different from $1$, $2$, $3$ and $5$.
\item $\xi, \alpha+\alpha^{-1}$ and $\beta+\beta^{-1}$ are in the field $L_G$ generated by the traces of the elements of $\rho_j([G,G])$,
\item In case $2$ and $3$, $\rho_j\simeq \epsilon \circ ^t\!\rho_j^{-1}$,
\item In case $4$, $\rho_{j|[G,G]}\simeq \epsilon \circ \rho_{j|[G,G]}$,
\item In case $6$, $\rho_{j|[G,G]}\simeq \epsilon \circ ^t\!\rho_{j|[G,G]}^{-1}$.
\end{enumerate} 

We now prove these assertions.

\begin{enumerate}
\item We have that $A_{1,2}=\frac{(\beta-1)(\alpha-1)(\alpha\beta-\gamma u_j+1)(\alpha+\gamma u_j+\beta)^2}{\alpha^2\beta^2}\neq 0$ by the assumptions on $\alpha$, $\beta$ and $\theta$.
\item We have that $B_{1,2}=-\frac{(\beta-1)(\alpha+1)(\alpha+\gamma u_j+\beta)(\alpha^2\beta+\alpha\beta^2+\alpha+\beta-2\alpha\beta+(\alpha\beta-\alpha-\beta+1-\gamma u_j)\gamma u_j)}{\alpha^2\beta^2}$. It follows that $B_{1,2}=0$ implies
$t=\alpha^2\beta+\alpha\beta^2+\alpha+\beta-2\alpha\beta+(\alpha\beta-\alpha-\beta+1-\gamma u_j)\gamma u_j=0$.

If $p\neq 2$ then we have $0=\frac{1}{2\alpha}\left(\frac{1}{\beta}(\frac{\alpha^2\beta^2(B_{1,1}+B_{2,2})}{2\alpha}-t)+(2-\beta)t\right)=\Phi_6(\beta)$ .This is absurd by the conditions on $\beta$.

If $p=2$ then we have $0=t=\alpha^2\beta+\alpha\beta^2+\alpha+\beta+(\alpha\beta+\alpha+\beta+1+\gamma u_j)\gamma u_j=(\alpha+\gamma u_j+\beta)(\alpha\beta+\gamma u_j +1)$ which contradicts our assumptions.

It follows that $B\neq 0$.
\item We have $(2-tr([\rho_j(T_s),\rho_j(T_t)])(\alpha+\alpha^{-1}-2)=tr([\rho_j(T_s)^2,\rho_j(T_t)])-2\in L_G$. We also have $C=2-tr([rho_j(T_s),\rho_j(T_t)])\in L_G$. We have $C=\frac{(\alpha\beta-\gamma u_j+1)(\alpha+\gamma u_j+\beta)}{\alpha\beta}\neq 0$. It follows that $\alpha+\alpha^{-1}\in L_G$.

We have that $C(\beta+\beta^{-1}-2)=tr([\rho_j(T_t)^2,\rho_j(T_s)])-2$, therefore we also have $\beta+\beta^{-1}\in L_G$.

We have $\xi=C-(\alpha+\alpha^{-1}+\beta+\beta^{-1})+\frac{2-tr([\rho_j(T_s)\rho_j(T_t)\rho_j(T_s),\rho_j(T_t)]}{C}\in L_G$. It follows that $K\subset L_G$ as required.

\item We can use the same arguments as for the proof of Theorem \ref{platypodes} to get that $\overline{\rho_j(G)}$ contains elements of order different from $1$, $2$, $3$ and $5$.

\item In case $2$, we have $L=L_1=L_2=L_3=K_1=K_2=K_3\neq K$, therefore $\epsilon(\alpha)=\alpha^{-1}$, $\epsilon(\beta)=\beta^{-1}$ and $\epsilon(\xi)=\xi$. It follows that $\epsilon(\gamma u_j)\frac{\alpha^{-1}\beta^{-1}-\alpha^{-1}-\beta^{-1}+1}{\alpha^{-1}\beta^{-1}}=\gamma u_j\frac{\alpha\beta-\alpha-\beta+1}{\alpha\beta}$, therefore $\epsilon(\gamma u_j)=\alpha^{-1}\beta^{-1}\gamma u_j$. We then set $P=\begin{pmatrix}
\frac{\alpha(\beta+1)}{\alpha+\gamma u_j+\beta} & \alpha\\
1 & \alpha+1
\end{pmatrix}$. We have that $\op{det}(P)=\frac{\alpha(\alpha\beta-\gamma u_j+1)}{\alpha \gamma u_j+\beta}\neq 0$ and 
$$P\rho_j(T_s)P^{-1}=\begin{pmatrix}
-1 & \alpha \\
0 & \alpha
\end{pmatrix} =\epsilon(\begin{pmatrix}
-1 & \alpha^{-1} \\
0 & \alpha^{-1}
\end{pmatrix})=\epsilon(^t\!\rho_j(T_s)^{-1})$$ and
$$P\rho_j(T_t)P^{-1}=\begin{pmatrix}
\beta & 0\\
\frac{\alpha+\gamma u_j+\beta}{\alpha} & -1
\end{pmatrix}=\epsilon(\begin{pmatrix}
\beta^{-1} & 0\\
\frac{\beta+\gamma u_j+\alpha}{\beta} & -1
\end{pmatrix})=\epsilon(^t\!\rho_j(T_t)^{-1}).$$
It follows that $\rho_j\simeq \epsilon \circ ^t\!\rho_j^{-1}$.

In case $3$, we have $L=L_1=L_2=L_3=K_1=K_3\neq K=K_2$, therefore $\epsilon(\alpha)=\alpha^{-1}$, 
$\epsilon(\beta)=\beta^{-1}$, $\epsilon(\gamma u_j)=\gamma u_j$ et $\epsilon(\xi)=\xi$. It follows that $\gamma u_j=\alpha\beta \gamma u_j$, therefore $u_j=0$ or $\alpha\beta=1$.

Assume first that $u_j=0$. We then have $\theta^j+\theta^{-j}=0$, therefore $\theta^{2j}=-1$ and $\theta^{4j}=1$. This implies that  $j=\frac{m}{4}$. Let $P=\begin{pmatrix} \frac{\alpha(\beta+1)}{\alpha+\beta} & \alpha\\
1 & \alpha+1\end{pmatrix}$. We have $det(P)=\frac{\alpha(\alpha\beta+1)}{\alpha+\beta}$ and $(P\rho_j(T_s)P^{-1},P\rho_j(T_t)P^{-1})=(\epsilon(^t\rho_j(T_s)^{-1}),\epsilon(^t(\rho_j(T_t)^{-1}))$.

Assume now that $\alpha\beta=1$. We then have $\gamma=1$, it follows that $\alpha^2+\alpha u_j+1=\alpha(\alpha+u_j+\alpha^{-1})=\alpha(\alpha+\gamma u_j+\beta)\neq 0$. Moreover, we have that $u_j-2=\theta^j+\theta^{-j}-2=(\theta^j-1)(1-\theta^{-j})\neq 0$. Let $P=\begin{pmatrix}
\frac{\alpha(\alpha+1)}{\alpha^2+\alpha  u_j+1} & \alpha\\
1 & \alpha+1 \end{pmatrix}$, we then have $\op{det}(P)=-\frac{\alpha^2( u_j-2)}{\alpha^2+\alpha u_j+1}\neq 0$ and $(P\rho_j(T_s)P^{-1},P\rho_j(T_t)P^{-1})=(\epsilon(^t\!\rho_j(T_s)^{-1}),\epsilon( ^t \!\rho_j(T_t)^{-1}))$. 

We then have that in both cases $\rho_j\simeq \epsilon \circ ^t\!\rho_j^{-1}$.

\item In case $4$, we have $L=L_2=L_3=K_3\neq L_1=K_1=K_2=K$, therefore $\epsilon(\alpha)=\alpha^{-1}$, $\epsilon(\beta)=\beta$, $\epsilon(\gamma u_j)=\gamma u_j$ and $\epsilon(\xi)=\xi$. We have that $\gamma u_j\frac{\alpha^{-1}\beta-\alpha^{-1}-\beta+1}{\alpha^{-1}\beta}=\gamma u_j \frac{\alpha\beta-\alpha-\beta+1}{\alpha\beta}$. It follows that $\gamma u_j(\alpha-1)^2(1-\beta)=0$, therefore $u_j=0$ and $j=\frac{m}{4}$. Let $P=\begin{pmatrix}
 1 & \alpha+1\\
 0 & -\alpha
 \end{pmatrix}$, we have $det(P)=-\alpha\neq 0$. We have $P\rho_j(T_s)P^{-1}=\begin{pmatrix}
 \alpha & 0\\
 -\alpha & -1
 \end{pmatrix} =-\alpha\epsilon(\rho_j(T_s))$  and $P\rho_j(T_t)P^{-1}=\begin{pmatrix}
 \beta & \beta+\alpha^{-1}\\
 0 & -1
 \end{pmatrix}=\epsilon(\rho_j(T_t))$. We then have that $\rho_{j|[G,G]}\simeq \epsilon\circ \rho_{j|[G,G]}$.
 
 \item In case $6$, we have that $L=L_1=L_2=L_3=K_1=K_2\neq K=K_3$, therefore $\epsilon(\alpha)=\alpha$, $\epsilon(\beta=\beta^{-1})$ and $\epsilon(\xi)=\xi$. It follows that $\epsilon(\gamma u_j)=\frac{\alpha\beta^{-1}(\alpha\beta-\alpha-\beta+1)}{\alpha\beta(\alpha\beta^{-1}-\alpha-\beta^{-1}+1)}\gamma u_j=-\beta^{-1}\gamma u_j$. Let $P=\begin{pmatrix}
 \beta+1 & \alpha+\gamma u_j+\beta\\
 \alpha\beta-\gamma u_j+1 & 0
 \end{pmatrix}$, we have that $det(P)=-(\alpha\beta-\gamma u_j+1)(\alpha+\gamma u_j+\beta)\neq 0$. We also have $P\rho_j(T_s)P^{-1}=-\alpha\epsilon(^t(\rho_j(T_s)^{-1})$ and $P\rho_j(T_t)P^{-1}=\epsilon(^t(\rho_j(T_t)^{-1})$. It follows that $\rho_{j|[G,G]}\simeq \epsilon \circ ^t\rho_{j|[G,G]}^{-1}$ and the proof is concluded.
\end{enumerate}

\end{proof}

\begin{lemme}\label{Isomorphismeven}
We say that $j\sim l$ if $\F_p(\alpha+\alpha^{-1},\beta+\beta^{-1},\xi_j)\simeq \F_p(\alpha+\alpha^{-1},\beta+\beta^{-1},\xi_j)$ and there exists $\Phi_{j,l}\in Aut(\F_{q_j})$ such that $\Phi_{j,l}(\alpha+\alpha^{-1})=\alpha+\alpha^{-1})$, $\Phi_{j,l}(\beta+\beta^{-1})=\beta+\beta^{-1}$ and $\Phi_{j,l}(\xi_j)=\xi_l$. This defines an equivalence relation and if $j\sim l$ then $\Phi_{j,l}\circ \rho_{j|\mathcal{A}_{I_2(m)}} \simeq \rho_{l|\mathcal{A}_{I_2(m)}}$.
\end{lemme}

\begin{proof}
Let us show this is an equivalence relation. Let $j,l,k\in [\![1,\frac{m-2}{2}]\!]$. We have $j\sim j$. If $j\sim l$ and $l\sim k$ then it it clear that $j\sim k$. This relation is symmetric because $\Phi_{l,j}=\Phi_{j,l}^{-1}$ verifies the desired conditions if $j\sim l$.

\smallskip

By Theorem \ref{bopb}, the second part of the statement is also true because $SU_2(q^{\frac{1}{2}})\simeq SL_2(q^{\frac{1}{2}})$.
\end{proof}

We now give the image of the derived subgroup of the Artin group in the full Iwahori-Hecke algebra.

\begin{theo}\label{resdihedraleven}
Assume $m$ even and $\alpha$ and $\beta$ verify the conditions given at the beginning of this section. For $j \in [\![1,\frac{m-2}{2}]\!]$, we set $G_j=\rho_j([<T_t,T_s>,<T_t,T_s>])$ .

We then have that the morphism from $\mathcal{A}_{I_2(m)}$ to $\mathcal{H}_{I_2(m),q}^\times \simeq GL_1(q_j)^2 \times \underset{j\in [\![1,\frac{m-1}{2}]\!]}\prod GL_2(q_j)$ factorizes through the surjective morphism
$$\Phi : \mathcal{A}_{I_2(m)} \rightarrow \underset{j\in [\![1,\frac{m-2}{2}]\!]/\sim}\prod G_j.$$
\end{theo}

\begin{proof}
We know by Theorem \ref{bopb} that the composition of the morphism from $\mathcal{A}_{I_2(m)}$ to $\mathcal{H}_{I_2(m),q}^\times$ with the projection upon each representation is surjective. We know by Lemma \ref{Isomorphismeven} that it factorizes through the morphism. We will now use Goursat's lemma and induction on $j\in [\![1,\frac{m-2}{2}]\!]$ in order to conclude the proof of this theorem. For $j_0\in [\![1,\frac{m-2}{2}]\!]$, we define $\Phi_{j_0}(\mathcal{A}_{I_2(m)})$ to be the image of $\mathcal{A}_{I_2(m)}$ inside $\underset{j\in [\![1,j_0]\!]/\sim}\prod GL_2(q_j)$. We know that $\Phi_1(\mathcal{A}_{I_2(m)})=G_1$. Let $j_0\in [\![1,\frac{m-4}{2}]\!]$, assume $\Phi_{j_0}(\mathcal{A}_{I_2(m)})=\underset{j\in [\![1,j_0]\!]/\sim}\prod G_{j}$.

\noindent
Consider $\Phi_{j_0+\ell}(\mathcal{A}_{I_2(m)})\subset \underset{j\in [\![1,j_0]\!]}\prod G_{j} \times G_{j_l}$ for $\ell$ the smallest positive integer such that $j_\ell\nsim j$ for all $j\in [\![1,j_0]\!]$. We know that the projection upon each factor is surjective. Let $K_1= \underset{j\in [\![1,j_0]\!]}\prod G_{j}$ and $K_2=G_{j_0+\ell}$ as in Goursat's Lemma.
 We then have $K_1/K^1\simeq K_2/K^2$. If the quotients are abelian then we are done since both groups are perfect. Assume that those quotients are non-abelian. There is only one decomposition factor of $K_2$ and it is equal to $PSL_2(q_{j_0+\ell})$, $PSU_2(q_{j_0+\ell}^\frac{1}{2})$ or $PSL_2(q_{j_0+\ell}^\frac{1}{2})$  depending on the field $\F_{q_{j_0+\ell}}$. We write that decomposition factor $PG_{j_0+\ell}$. The isomorphism then implies that there exists $j_1\in [\![1,j_0]\!]$ such that $\F_p(\alpha+\alpha^{-1},\beta+\beta^{-1},\xi_{j_1})\simeq \F_p(\alpha+\alpha^{-1},\beta+\beta^{-1},\xi_{j_0+\ell})$ and $\overline{\rho_{j_1}(\mathcal{A}_{I_2(m})}\simeq PG_{j_1} \simeq \overline{\rho_{j_0+\ell}(\mathcal{A}_{I_2(m})}$. We then have that there exists $z:\mathcal{A}_{I_2(m)} \rightarrow \overline{\F_p}^\times$ and $\Psi\in Aut(\F_{q_{j_0+\ell}})$ such that up to conjugation, for all $h\in \mathcal{A}_{I_2(m)}, \rho_{j_0+\ell}(h)=\Psi(\rho_{j_1}(h))z(h)$. We will prove this is absurd by considering traces of some elements in $\mathcal{A}_{I_2(m)}$ under these representations. We may first note that for all $M\in SL_2(\overline{\F_q})$, we have $1=\op{det}(z(h)M)=z(h)^2\op{det}(M)=z(h)^2$, therefore for all $h \in \mathcal{A}_{I_2(m)}$, $z(h)\in \{-1,1\}$. We write as before in what follows $\gamma=\sqrt{\alpha\beta}$, $u_{j_0+\ell}=\theta^{j_0+\ell}+\theta^{-(j_0+\ell)}$, $u_{j_1}=\theta^{j_1}+\theta^{-j_1}$, $\xi_{j_1}=\frac{\gamma u_{j_1}(\alpha\beta-\alpha-\beta+1)}{\alpha\beta}$ and $\xi_{j_0+\ell}=\frac{\gamma u_{j_0+\ell}(\alpha\beta-\alpha-\beta+1)}{\alpha\beta}$.
 
\medskip 
 
\textbf{1.} Assume first that $z([T_t,T_s])=1$ and $\tr(\rho_{j_1}([T_t,T_s]))\neq 0$. We then have that 
$$A_1=\tr(\rho_{j_0+\ell}([T_t,T_s]))=\Phi(\tr(\rho_{j_1}([T_t,T_s])))=\Phi(B_1).$$
 We have $\tr(\rho_{j_0+\ell}([T_t^{-1},T_s]))=\tr(\rho_{j_0+\ell}([T_t,T_s]))$ and $\tr(\rho_{j_1}([T_t^{-1},T_s]))=\tr(\rho_{j_1}([T_t,T_s]))$, therefore $z([T_t^{-1},T_s])=1$ because $\tr(\rho_{j_1}([T_t,T_s]))\neq 0$. It follows that $z([T_t,T_s][T_t^{-1},T_s])=1$. We then have
$$A_2=\tr(\rho_{j_0+\ell}([T_t,T_s][T_t^{-1},T_s]))=\Phi(\tr(\rho_{j_0+\ell}([T_t,T_s][T_t^{-1},T_s])))=\Phi(B_2).$$ We have $A_2=A_1^2-(\beta+\beta^{-1}-2)(A_1-2)-2$ and $B_2=B_1^2-(\beta+\beta^{-1}-2)(B_1-2)-2$. It follows that 
 $$A_1^2-(\beta+\beta^{-1}-2)(A_1-2)-2=\Phi(B_1^2-(\beta+\beta^{-1}-2)(B_1-2)-2)=A_1^2-\Phi(\beta+\beta^{-1}-2)(A_1-2)-2.$$
 $$((\beta+\beta^{-1})-\Phi(\beta+\beta^{-1}))(A_1-2)=0.$$
 We have $A_1-2=-\frac{(\alpha\beta-\gamma u_j+1)(\alpha+\gamma u_j+\beta)}{\alpha\beta}\neq 0$ by assumption, therefore $\Phi(\beta+\beta^{-1})=\beta+\beta^{-1}$.
 Note that we also have $A_1=\tr(\rho_{j_0+\ell}([T_s,T_t]))=\tr(\rho_{j_0+\ell}([T_s^{-1},T_t]))$ and $B_1=\tr(\rho_{j_1}([T_s,T_t]))=\tr(\rho_{j_1}([T_s^{-1},T_t]))$. It follows that $z([T_s,T_t][T_s^{-1},T_t])=1$ and 
 $$A_3=\tr(\rho_{j_0+\ell}([T_s,T_t][T_s^{-1},T_t]))=\Phi(\tr(\rho_{j_1}([T_s,T_t][T_s^{-1},T_t])))=\Phi(B_3).$$ We have $A_3=A_1^2-(\alpha+\alpha^{-1}-2)(A_1-2)-2$ and $B_3=B_1^2-(\alpha+\alpha^{-1}-2)(B_1-2)-2$. It follows that by the same reasoning as before $\alpha+\alpha^{-1}=\Phi(\alpha+\alpha^{-1})$.
 
 \smallskip
 
\textbf{1.1.} Assume now that $z([T_tT_sT_t,T_s])=1$. We then have $A_4=\tr(\rho_{j_0+\ell}([T_tT_sT_t,T_s]))=\Phi(tr(\rho_{j_1}([T_tT_sT_t,T_s])))=\Phi(B_4)$. We also have $\xi_{j_0+\ell}=\frac{2-A_4}{2-A_1}+2-A_1-(\alpha+\alpha^{-1}+\beta+\beta^{-1})$ and 
$\xi_{j_1}=\frac{2-B_4}{2-B_1}+2-B_1-(\alpha+\alpha^{-1}+\beta+\beta^{-1})$. It follows that $\xi_{j_0+\ell}=\Phi(\xi_{j_1})$ and, therefore $j_0 +\ell\simeq j_1$, which contradicts our assumptions.

\smallskip

\textbf{1.2.} Assume now that $z([T_tT_sT_t,T_s])=-1$. We then have $A_4=\tr(\rho_{j_0+\ell}([T_tT_sT_t,T_s]))=\Phi(tr(\rho_{j_1}([T_tT_sT_t,T_s])))=-\Phi(B_4)$. We have $z([T_tT_sT_t,T_s][T_t,T_s])=-1$, therefore $A_5=\tr(\rho_{j_0+\ell}([T_tT_sT_t,T_s][T_t,T_s]))=-\Phi(tr(\rho_{j_1}([T_tT_sT_t,T_s][T_t,T_s])))=-\Phi(B_5)$. We have $A_5=A_4A_1-A_1$ and $B_5=B_4B_1-B_1$, therefore $A_4A_1-A_1=-\Phi(B_4)\Phi(B_1)+\Phi(B_1)=A_4A_1+A_1$. It follows that $2A_1=0$ which contradicts our assumption since $p\neq 2$.

\medskip

\textbf{2.} Assume now that $z([T_t,T_s])=-1$ and $\tr(\rho_{j_1}([T_t,T_s]))\neq 0$. We have that $A_1=\tr(\rho_{j_1}([T_t,T_s]))=-\Phi(\tr(\rho_{j_1}([T_t,T_s])))=-\Phi(B_1)$.

 We then have that $z([T_t,T_s])=z([T_t^{-1},T_s])=z([T_s^{-1},T_t])$.
 
  It follows that $z([T_t,T_s][T_t^{-1},T_s])=z([T_s,T_t][T_s^{-1},T_t])=1$. We then have 
 $$A_2=\tr(\rho_{j_0+\ell}([T_t,T_s][T_t^{-1},T_s]))=\Phi(\tr(\rho_{j_0+\ell}([T_t,T_s][T_t^{-1},T_s])))=\Phi(B_2),$$
$$A_3=\tr(\rho_{j_0+\ell}([T_s,T_t][T_s^{-1},T_t]))=\Phi(\tr(\rho_{j_1}([T_s,T_t][T_s^{-1},T_t])))=\Phi(B_3).$$
 We have $A_2=A_1^2-(\beta+\beta^{-1}-2)(A_1-2)-2$, $B_2=B_1^2-(\beta+\beta^{-1}-2)(B_1-2)-2$, $A_3=A_1^2-(\alpha+\alpha^{-1}-2)(A_1-2)-2$ and $B_3=B_1^2-(\alpha+\alpha^{-1}-2)(B_1-2)-2$. It follows that 
\begin{footnotesize}
$$A_1^2-(\beta+\beta^{-1}-2)(A_1-2)-2=\Phi(B_1^2)-\Phi(\beta+\beta^{-1}-2)\Phi(B_1-2)-\Phi(2)=A_1^2-\Phi(\beta+\beta^{-1}-2)(-A_1-2)-2$$
\end{footnotesize}
and 
\begin{footnotesize}
$$A_1^2-(\alpha+\alpha^{-1}-2)(A_1-2)-2=\Phi(B_1^2)-\Phi(\alpha+\alpha^{-1}-2)\Phi(B_1-2)-\Phi(2)=A_1^2-\Phi(\alpha+\alpha^{-1}-2)(-A_1-2)-2$$
\end{footnotesize}
It follows that
$$(\beta+\beta^{-1}-2)=\Phi(\beta+\beta^{-1}-2)\frac{-A_1-2}{A_1-2}$$
 and 
 $$(\alpha+\alpha^{-1}-2)=\Phi(\alpha+\alpha^{-1}-2)\frac{-A_1-2}{A_1-2}$$

\smallskip

\textbf{2.1.} Assume now that $z([T_t^2,T_s^2])=1$. We then have that $A_8=\tr(\rho_{j_0+\ell}([T_t^2,T_s^2]))=\Phi(\tr(\rho_{j_1}([T_t^2,T_s^2])))=\Phi(B_8)$. We have $A_8=(\beta+\beta^{-1}-2)(\alpha+\alpha^{-1}-2)(A_1-2)+2$ and $B_8=(\beta+\beta^{-1}-2)(\alpha+\alpha^{-1}-2)(B_1-2)+2$. It follows that
\begin{eqnarray*}
(\beta+\beta^{-1}-2)(\alpha+\alpha^{-1}-2)(A_1-2)+2 & = & \Phi(\beta+\beta^{-1}-2)\Phi(\alpha+\alpha^{-1}-2)(\Phi(B_1)-2)+2\\
(\beta+\beta^{-1}-2)(\alpha+\alpha^{-1}-2)(A_1-2) & = & \frac{(A_1-2)^2(-A_1-2)}{(-A_1-2)^5}(\beta+\beta^{-1}-2)(\alpha+\alpha^{-1}-2)\\
(A_1-2)(-A_1-2) & = & (A_1-2)^2\\
-A_1^2+4& = & A_1^2-4A_1+4\\
2(A_1^2-2A_1) & = & 0\\
2A_1(A_1-2) & = & 0.
\end{eqnarray*}
This is a contradiction since $p\neq 2$, $A_1\neq 2$ and, by assumption, $A_1\neq 0$.

\medskip

\textbf{2.2.} Assume now that $z([T_t^2,T_s^2])=-1$. We then have that 
$$A_8=\tr(\rho_{j_0+\ell}([T_t^2,T_s^2]))=-\Phi(\tr(\rho_{j_1}([T_t^2,T_s^2])))=-\Phi(B_8).$$ We have $A_8=(\beta+\beta^{-1}-2)(\alpha+\alpha^{-1}-2)(A_1-2)+2$ and $B_8=(\beta+\beta^{-1}-2)(\alpha+\alpha^{-1}-2)(B_1-2)+2$. It follows that
\begin{footnotesize}
\begin{eqnarray*}
(\beta+\beta^{-1}-2)(\alpha+\alpha^{-1}-2)(A_1-2)+2 & = & -\Phi(\beta+\beta^{-1}-2)\Phi(\alpha+\alpha^{-1}-2)(\Phi(B_1)-2)-2\\
(\beta+\beta^{-1}-2)(\alpha+\alpha^{-1}-2)(A_1-2)+2 & = & \frac{(A_1+2)(A_1-2)^2}{(A_1+2)^2}(\beta+\beta^{-1}-2)(\alpha+\alpha^{-1}-2)-2\\
(\beta+\beta^{-1}-2)(\alpha+\alpha^{-1}-2)(A_1^2-4-A_1^2+4A_1-4)& = & -4(A_1+2)\\
(\beta+\beta^{-1}-2)(\alpha+\alpha^{-1}-2)(A_1-2) & = &-A_1-2\\
A_8-2 & = & -A_1-2\\
A_8 & = & -A_1.
\end{eqnarray*}
\end{footnotesize}

\smallskip

\textbf{2.2.1.} Assume now that $z([T_tT_sT_t,T_s])=-1$. We then have $z([T_tT_sT_t,T_s][T_t,T_s])=1$, 
$$A_4=\tr(\rho_{j_0+\ell}([T_tT_sT_t,T_s]))=\Phi(tr(\rho_{j_1}([T_tT_sT_t,T_s])))=-\Phi(B_4),$$ 
$$A_5=\tr(\rho_{j_0+\ell}([T_tT_sT_t,T_s][T_t,T_s]))=\Phi(tr(\rho_{j_1}([T_tT_sT_t,T_s][T_t,T_s])))=\Phi(B_5).$$
 We have $A_5=A_4A_1-A_1$ and $B_5=B_4B_1-B_1$, therefore $A_4A_1-A_1=\Phi(B_4)\Phi(B_1)-\Phi(B_1)=A_4A_1+A_1$. It follows that $A_4=A_1$ which contradicts our assumptions.

\smallskip

\textbf{2.2.2.} Assume now that $z([T_tT_sT_t,T_s])=1$. We then have $A_4=\tr(\rho_{j_0+\ell}([T_tT_sT_t,T_s]))=\Phi(tr(\rho_{j_1}([T_tT_sT_t,T_s])))=\Phi(B_4)$.

\smallskip

\textbf{2.2.2.1.} Assume that $z([T_t^2,T_s])=1$. We then have
 $$A_6=\tr(\rho_{j_0+\ell}([T_t^2,T_s]))=\Phi(tr(\rho_{j_1}([T_t^2,T_s])))=\Phi(B_6).$$
  We then have $z([T_t,T_s][T_t^2,T_s][T_tT_sT_t,T_s])=-1$, therefore 
  $$A_7=\tr(\rho_{j_0+\ell}([T_t,T_s][T_t^2,T_s][T_tT_sT_t,T_s]))=-\Phi(\tr(\rho_{j_1}([T_t,T_s][T_t^2,T_s][T_tT_sT_t,T_s])))=-\Phi(B_7).$$
   We have $A_7=A_1A_4A_6-A_1A_4-A_1A_6-A_1-\xi_{j_0+\ell}(A_1-2)$ and $B_7=B_1B_4B_6-B_1B_4-B_1B_6-B_1-\xi_{j_1}(B_1-2)$. It follows that $\xi_{j_0+\ell}=\Phi(\xi_{j_1})\frac{A_1+2}{A_1-2}$.

We have $z([T_t^2,T_s^2][T_t,T_s])=1$, therefore 
$$A_9=\tr(\rho_{j_0+\ell}([T_t^2,T_s^2][T_t,T_s]))=\Phi(\tr(\rho_{j_1}([T_t^2,T_s^2][T_t,T_s])))=\Phi(B_9).$$
 We also have $A_9=A_8A_1-A_8+(\xi_{j_0+\ell}-1)(A_1-2)$ and $B_9=B_8B_1-B_8+(\xi_{j_1}-1)(B_1-2)$. It follows that 
\begin{eqnarray*}
A_8A_1-A_8+(\xi_{j_0+\ell}-1)(A_1-2) & = & \Phi(B_8B_1-B_8+(\xi_{j_1}-1)(B_1-2))\\
& = & A_8A_1+A_8+(\Phi(\xi_{j_1})-1)(-A_1-2)\\
& = & A_8A_1-A_1+(\xi_{j_0+\ell}\frac{A_1-2}{A_1+2}-1)(-A_1-2)\\
A_1 +(\xi_{j_0+\ell}-1)(A_1-2)& = & -A_1-(\xi_{j_0+\ell}(A_1-2)-A_1-2)\\
A_1+\xi_{j_0+\ell}(A_1-2)-A_1+2 & = & -A_1-\xi_{j_0+\ell}(A_1-2)+A_1+2\\
2\xi_{j_0+\ell}(A_1-2) & = & 0.
\end{eqnarray*}
This implies that $\xi_{j_0+ \ell}=0=\xi_{j_1}$. It then follows that $u_{j_0+\ell}=0=u_{j_1}$. It follows that $\theta^{j_0+\ell}+\theta^{-(j_0+\ell)}=\theta^{j_1}+\theta^{-j_1}$ and, therefore $(\theta^{j_0+\ell+j_1}-1)(\theta^{-j_1}-\theta^{-(j_0+l)})=0$. Since we have $1\leq j_0+\ell,j_1\leq \frac{m-2}{2}$, this implies that $j_0+\ell=j_1$ which contradicts our assumptions.

\smallskip

\textbf{2.2.2.2.} Assume now that $z([T_t^2,T_s])=-1$. We then have 
$$A_6=\tr(\rho_{j_0+\ell}([T_t^2,T_s]))=-\Phi(tr(\rho_{j_1}([T_t^2,T_s])))=-\Phi(B_6).$$ We then have $z([T_t^2,T_s][T_t,T_s])=1$, therefore
$$A_{10}=\tr(\rho_{j_0+\ell}([T_t^2,T_s][T_t,T_s]))=\Phi(\tr(\rho_{j_1}([T_t^2,T_s][T_t,T_s])))=\Phi(B_{10}).$$ We have $A_{10}=A_1A_6-A_1$ and $B_{10}=B_1B_6-B_1$. We then have $A_1A_6-A_1=A_{10}=\Phi(B_{10})=A_1A_6+A_1$, therefore $2A_1=0$ and $A_1=0$. This contradicts our assumptions.

\medskip

\textbf{3.} Assume now that $A_1=\tr(\rho_{j_0+\ell}([T_t,T_s]))=0$. We then have $B_1=\tr(\rho_{j_1}([T_t,T_s]))=0$.

\smallskip

\textbf{3.1.} Assume $z([T_t^2,T_s])=1$. We then have $A_6=\tr(\rho_{j_0+\ell}([T_t^2,T_s]))=\Phi(\tr(\rho_{j_1}([T_t^2,T_s])))=\Phi(B_6)$. We also have $A_6=2(\beta+\beta^{-1}-1)-A_1(\beta+\beta^{-1}-2)=2(\beta+\beta^{-1}-1)$ and $B_6=2(\beta+\beta^{-1}-1)-B_1(\beta+\beta^{-1}-2)=2(\beta+\beta^{-1}-1)$. It follows that $\Phi(\beta+\beta^{-1})=\beta+\beta^{-1}$. 

\smallskip

\textbf{3.1.1.} Assume $z([T_t^2,T_s^2])=1$. We then have 
$$A_8=\tr(\rho_{j_0+\ell}([T_t^2,T_s^2]))=\Phi(\tr(\rho_{j_1}([T_t^2,T_s^2])))=\Phi(B_8).$$
 We have 
 $$A_8=(\beta+\beta^{-1}-2)(\alpha+\alpha^{-1}-2)(A_1-2)+2=2-2(\beta+\beta^{-1}-2)(\alpha+\alpha^{-1}-2),$$
 $$B_8=(\beta+\beta^{-1}-2)(\alpha+\alpha^{-1}-2)(B_1-2)+2=2-2(\beta+\beta^{-1}-2)(\alpha+\alpha^{-1}-2).$$
  It follows that $\alpha+\alpha^{-1}=\Phi(\alpha+\alpha^{-1})$. 

\smallskip

\textbf{3.1.1.1.} Assume $z([T_tT_sT_t,T_s])=1$. We then have 
$$A_4=\tr(\rho_{j_0+\ell}([T_tT_sT_t,T_s]))=\Phi(tr(\rho_{j_1}([T_tT_sT_t,T_s])))=\Phi(B_4).$$ We have 
$\xi_{j_0+\ell}=\frac{2-A_4}{2-A_1}+2-A_1-(\alpha+\alpha^{-1}+\beta+\beta^{-1})$ and 
$\xi_{j_1}=\frac{2-B_4}{2-B_1}+2-B_1-(\alpha+\alpha^{-1}+\beta+\beta^{-1})$.
It follows that $\xi_{j_0+\ell}=\Phi(\xi_{j_1})$ and, therefore, $j_0+\ell\sim j_1$. This contradicts our assumptions.

\smallskip

\textbf{3.1.1.2.} Assume now $z([T_tT_sT_t,T_s])=-1$. We then have $A_4=-\Phi(B_4)$. By the same computations as in \textbf{3.1.1.1}, we get that 
$$\xi_{j_0+\ell}=\frac{2-A_4}{2}+2-(\alpha+\alpha^{-1}+\beta+\beta^{-1})=\Phi(\frac{2+B_4}{2}+2-(\alpha+\alpha^{-1}+\beta+\beta^{-1}))=\Phi(\xi_{j_1})+\Phi(B_4).$$
We also have $z([T_t^2,T_s][T_tT_sT_t,T_s])=-1$, therefore $A_{11}=\tr(\rho_{j_0+\ell}([T_t^2,T_s][T_tT_sT_t,T_s]))=-\Phi(\tr(\rho_{j_1}([T_t^2,T_s][T_tT_sT_t,T_s])))=-\Phi(B_{11})$. We also have $A_{11}  =  A_4A_6-A_4-A_6+\xi_{j_0+\ell}(A_1-2)+2$ and $B_{11}  =  B_4B_6-B_4-B_6+\xi_{j_1}(B_1-2)+2$. It follows that $\xi_{j_0+\ell}=-A_6-\Phi(\xi_{j_1})+2$. We then have $\Phi(\xi_{j_1})=-A_6-\xi_{j_0+\ell}+2=\xi_{j_0+\ell}+A_4$ and, therefore $2\xi_{j_0+\ell}=-A_6-A_4-2$. We have 
$$-A_6-A_4+2=-A_1^2-\xi_{j_0+\ell}A_1-(\alpha+\alpha^{-1}-2)A_1+2\xi_{j_0+\ell}+2(\alpha+\alpha^{-1}-1).$$
It follows that $\alpha+\alpha^{-1}-1=0$ and, therefore, $0=\alpha^2-\alpha+1=\Phi_6(\alpha)$, which contradicts our assumptions.

\smallskip

\textbf{3.1.2.} Assume now $z([T_t^2,T_s^2])=-1$. We then have 
$$A_8=\tr(\rho_{j_0+\ell}([T_t^2,T_s^2]))=-\Phi(\tr(\rho_{j_1}([T_t^2,T_s^2])))=-\Phi(B_8).$$ We have $z([T_t^2,T_s][T_t^2,T_s^2])=-1$, therefore 
$$A_{12}=\tr(\rho_{j_0+\ell}([T_t^2,T_s^2][T_t^2,T_s]))=-\Phi(\tr(\rho_{j_1}([T_t^2,T_s^2][T_t^2,T_s])))=-\Phi(B_{12}).$$
 We also have $A_{12}=A_8A_6-A_6$ and $B_{12}=B_8B_6-B_6$, therefore $A_8A_6-A_6=A_8A_6+A_6$ and $A_6=0$. We then have $$0=A_6=-(\beta+\beta^{-1}-2)A_1+2(\beta+\beta^{-1}-1)=\frac{2}{\beta}\Phi_6(\beta).$$
This contradicts our assumptions.

\smallskip

\textbf{3.2.} Assume now $z([T_t^2,T_s])=-1$.  We then have 
$$A_6=\tr(\rho_{j_0+\ell}([T_t^2,T_s]))=-\Phi(\tr(\rho_{j_1}([T_t^2,T_s]))=-\Phi(B_6).$$
 We have $A_6=2(\beta+\beta^{-1}-1)=B_6$, therefore $\Phi(\beta+\beta^{-1}-1)=-(\beta+\beta^{-1}-1)$. We have $A_6\neq 0$, $\tr(\rho_{j_0+\ell}(T_t^{-2},T_s]))=A_6$ and $\tr(\rho_{j_1}(T_t^{-2},T_s]))=B_6$. It follows that $z([T_t^{-2},T_s])=z([T_t^2,T_s])=-1$ and $z([T_t^2,T_s][T_t^{-2},T_s])=1$. We then have $A_{12}=\tr(\rho_{j_0+\ell}([T_t^2,T_s][T_t^{-2},T_s]))=\Phi(B_{12})$. We also have $A_{12}=A_6^2-(\beta+\beta^{-1})^2(A_6-2)-2=B_{12}$. It follows that 
\begin{eqnarray*}
(\frac{A_6}{2}+1)^2(A_6-2) & = & \Phi\left((\frac{A_6}{2}+1)^2(A_6-2)\right)\\
(\frac{A_6^2}{4}+A_6+1)(A_6-2) &  = & (\frac{A_6^2}{4}-A_6+1)(-A_6-2)\\
\frac{A_6^3}{4}+A_6^2+A_6-\frac{A_6^2}{2}-2A_6-2 & = & -\frac{A_6^3}{4}+A_6^2-A_6-\frac{A_6^2}{2}+2A_6-2\\
\frac{A_6^3}{2}-2A_6 & = & 0\\
\frac{A_6}{2}(A_6^2-4) & = & 0\\
\frac{A_6}{2}(A_6-2)(A_6+2) & = & 0\\
4(\beta+\beta^{-1}-1)(\beta+\beta^{-1}-2)(\beta+\beta^{-1}) & = & 0\\
4\beta^{-3}\Phi_6(\beta)\Phi_1(\beta)^2\Phi_4(\beta) & = & 0.
\end{eqnarray*}
This is a contradiction by the assumptions on the order of $\beta$. This concludes the proof.
\end{proof}

\chapter{W-graphs}\label{Wgraphschapter}

Before extending our study to types $H_3$ and $H_4$, we need to introduce the notion of $W$-graphs and give some properties which they verify.  We also prove some new properties and propose a new conjecture. In this section, $(W,S)$ is a Coxeter system with $W$ a finite Coxeter group, $K'$ is the splitting field of $W$, $K=K'((\sqrt{\alpha_s})_{s\in S}), \tilde{K}=K'((\alpha_s)_{s\in S})$, $C'$ is the ring of integers of $K'$ and $C=C'((\sqrt{\alpha_s})_{s\in S})$. We consider the Iwahori-Hecke algebra given by the presentation $\mathcal{H}=\mathcal{H}_{W,(\alpha_s)_{s\in S}}=<T_1,...,T_n| \underset{m_{s_i,s_j}}{\underbrace{T_iT_jT_i...}} =\underset{m_{s_i,s_j}}{\underbrace{T_jT_iT_j...}}, (T_i -\alpha_{s_i})(T_i+1)=0>$, where $\alpha_{s_i}=\alpha_{s_j}$ if $s_i$ and $s_j$ are in the same conjugacy class of $W$.

 W-graphs were introduced in 1979 by Kazhdan-Lusztig \cite{KL} for one-parameter families and the definition was extended to all Coxeter groups in \cite{G-P}. We here give the definition from \cite{G-P}. We will prove some uniqueness properties (Proposition \ref{Unique1} and Proposition \ref{Unique2}) in the one-parameter case and establish a conjecture (Conjecture \ref{conjecturewgraphs}) for certain $W$-graphs.
 
 \bigskip
 
 We here give the definition of $W$-graphs which can be found in \cite{G-P}

\begin{Def2}
For $X$ a set, we write $D(X)=\{(x,x),x\in X\}$ its diagonal. A $W$-graph $\Gamma$ is given by a triple $(X,I,\mu)$ such that
\begin{enumerate}
\item $X$ is a set and $I$ is a map from $X$ to $\mathcal{P}(S)$,
\item $\mu$ is a map from $(X\times X\setminus D(X)\times S)$ to $K$ stable by the field involution sending $\sqrt{\alpha_s}$ to $\sqrt{\alpha_s}^{-1}$.
Let $V$ be a $\tilde{K}$-vector space with basis $(e_y)_{y\in X}$. For all $s\in S$, we define $\rho_s :V\rightarrow V$ by
\begin{eqnarray*}
e_y   \mapsto &  -e_y & \mbox{if}~ s\in I(y),\\
e_y   \mapsto & \alpha_s e_y+\underset{x\in X,s\in I(x)}\sum \sqrt{\alpha_s}\mu_{x,y}^s e_x & \mbox{if}~ s\notin I(y).
\end{eqnarray*}

\item The map $T_s\mapsto \rho_s$ affords a representation of $\mathcal{H}$.
\end{enumerate}

\end{Def2}

For $\Gamma$ a $W$-graph, we write $\rho_\Gamma$ its associated representation and $V_\Gamma$ the corresponding $\mathcal{H}_{K,(\alpha_s)_{s\in S}}$-module.

\begin{Def2}
A $W$-graph $(X,I,\mu)$ is said to be $2$-colorable whenever there exists a map $\omega :X\rightarrow \{-1,1\}$ such that for any $(s,x,y)\in S\times X^2$ verifying $\mu^s(x,y)\neq 0$, we have $\omega(x)=-\omega(y)$.
\end{Def2}

The data given by the triple $(X,I,\mu)$ can be represented by a weighted oriented graph for which the set $X$ represents the vertices and the map $I$ represents the weight of the vertices. The non-zero values of the map $\mu$ represent the oriented weighted edges of the graph. After they were introduced, it was shown by Alvis and Lusztig \cite{A-L} that there exist $W$-graphs affording all irreducible representations of Coxeter groups of types $H_3$ and $H_4$. Using those results, Gyoja \cite{Gyo} showed that any irreducible representation of an Iwahori-Hecke algebra associated to a finite non-crystalographic Coxeter group in the equal parameters case could be afforded by a $W$-graph. Moreover, he showed the following result

\begin{theo2}\label{Melyssaestparfaite}
If all the parameters of the Hecke algebra are equal, then the following statements hold
\begin{enumerate}
\item Every irreducible $\mathcal{H}_{K,\alpha}$-module is afforded by $V_\Gamma$ for a $W$-graph $\Gamma$ over $C$.
\item An irreducible $H_{q,K}(W)$-module is afforded by $V_\Gamma$ for a $2$-colorable $W$-graph $\Gamma$ over $K(\sqrt{\alpha})$ if and only it admits a form over $K(\sqrt{\alpha})$. The representation is then said to be non-exceptional.
\end{enumerate}
\end{theo2} 

\textbf{Remark} : He also classified all the exceptional representations for Iwahori-Hecke algebras, that is, the ones which are not non-exceptional, and obtained

\begin{enumerate}
\item $2$ irreducible representations of dimension 512 of $E_7$,
\item $4$ irreducible representations of dimension 4096 of $E_8$,
\item $2$ irreducible representations of dimension 4 of $H_3$,
\item $4$ irreducible representations of dimension 16 of $H_4$,
\end{enumerate}

The link between $2$-colorable representations and non-exceptionality can be seen through the following proposition.

\begin{prop2}\label{color}
If $(X,I,\mu)$ is a $2$-colorable $W$-graph then

$$\rho_{(X,I,\mu)}\simeq \rho_{(X,I,-\mu)}.$$

If $\sigma: K((\sqrt{\alpha_s})_{s\in S}) \rightarrow K((\sqrt{\alpha_s})_{s\in S})$ leaves $K$ stable and maps $\sqrt{\alpha_s}$ to $-\sqrt{\alpha_s}$ for all $s\in S$ then $\rho_{(X,I,-\mu)}=\sigma \circ \rho_{(X,I,\mu)}$.
\end{prop2}

\begin{proof}
Let $(X,I,\mu)$ be a $2$-colorable $W$-graph and $\omega:X\rightarrow \{-1,1\}$ be an associated coloring. Let $\mathcal{L}:V\rightarrow V$ be the linear map defined by $\mathcal{L}(e_x)=\omega(x)e_x$ for all $x\in X$. We have $\mathcal{L}^{-1}=\mathcal{L}$. Let $s\in S$ and $y\in X$. If $s\in I(y)$, then
\begin{eqnarray*}
\mathcal{L}(\rho_{(X,I,\mu)}(\mathcal{L}(e_y))) & = &  -\omega(y)^2e_y \\
& = &  -e_y \\
 & = &  \rho_{(X,I,-\mu)}(e_y).
\end{eqnarray*}

If $s\notin I(y)$ then 
\begin{eqnarray*}
\mathcal{L}(\rho_{(X,I,\mu)(T_s)}(\mathcal{L}(e_y)))  & = & \mathcal{L}(\omega(y)\alpha_s e_y+\underset{x\in X,s\in I(x)}\sum \sqrt{\alpha_s}\mu_{x,y}^s\omega(y) e_x )\\
& =& \omega(y)^2\alpha_s e_y+\underset{x\in X,s\in I(x)}\sum \sqrt{\alpha_s}\mu_{x,y}^s\omega(x)\omega(y) e_x \\
& =& \alpha_s e_y+\underset{x\in X,s\in I(x)}\sum \sqrt{\alpha_s}(-\mu_{x,y}^s) e_x\\
 & = &  \rho_{(X,I,-\mu)}(T_s)(e_y).
\end{eqnarray*}

The second part of the proposition is proved in the same way.\end{proof}

Proposition \ref{color} shows that two different $W$-graphs can give isomorphic representations. We now provide some uniqueness conditions with the following propositions.

\begin{prop2}\label{Unique1}
If $\Gamma_1=(X,I_1,\mu_1)$ and $\Gamma_2=(Y,I_2,\mu_2)$ are $W$-graphs such that $\rho_{\Gamma_1}\simeq \rho_{\Gamma_2}$ and $\rho_{\Gamma_1}$ is irreducible, then there exists a bijection $\varphi:X\rightarrow Y$ such that for all $x\in X$ we have $I_2(\varphi(x))=I_1(x)$.
\end{prop2}

\begin{proof}
Let $\Gamma_1$ and $\Gamma_2$ be as above. If $n=1$, then the result is straightforward. Let us assume $n\geq 2$. Let us show that for all $S'\subset S$, we have $$\dim(\underset{s\in S'}\cap \op{ker}(\rho_{\Gamma_1}(T_s)+1))= \vert \{x\in X,S'\subset I(x)\}\vert.$$

If there exists $x_0\in X$ such that $I(x_0)=\emptyset$, then by definition of a $W$-graph, we would have that $\tilde{V}=\underset{x\in X\setminus\{x_0\}}\bigoplus Ke_x$ is a stable subvector-space of dimension $n-1$ which contradicts the irreducibility of $\rho_{\Gamma_1}$.

Let $S'=\{s_{i_l}\}_{l\in [\![1,k]\!]}\subset S$ for some $k\in \N^\star$. We label the vertices $\{e_{x_i}\}_{i\in[\![1,n_{\Gamma_1}]\!]}$ in a such a way that there exist $r_{S'}\in \N^\star$ verifying $S'\subset I(x_{h})$ for all $h\in [\![1,r_{S'}]\!]$ and $S'\not\subset I(x_{h})$ for all $h\in [\![r_{S'}+1,n_{\Gamma_1}]\!]$.

\medskip

Let $x\in \underset{s\in S'}\cap \op{ker}(\rho_{\Gamma_1}(T_s)+1)$. There exists a unique family $(\lambda_i)_{i\in [\![1,n_{\Gamma_1}]\!]}\in \tilde{K}^{n_{\Gamma_1}}$ such that $x=\underset{i\in [\![1,n_{\Gamma_1}]\!]}\sum \lambda_i e_{x_i}$.

Let $h_0\in [\![r_{S'}+1,n_{\Gamma_1}]\!]$, we will show that $\lambda_{h_0}=0$. There exists $j\in [\![1,k]\!]$ such that $s_{i_j}\notin I(x_{h_0})$. We then have $\rho_{\Gamma_1}(T_{s_{i_j}})(x)=\alpha_{s_{i_j}}\lambda_{h_0}e_{x_{h_0}}+\underset{h\neq h_0}\sum a_{h}e_{x_h}$ for some $a_h\in \tilde{K}$. 

Since $x\in \op{ker}(\rho_{\Gamma_1}(T_{s_{i_j}})+1)$, we also have $\rho_{\Gamma_1}(T_{s_{i_j}})(x)=-x=-\underset{i\in [\![1,n_{\Gamma_1}]\!]}\sum \lambda_i e_{x_i}$. This implies that $\alpha_{s_{i_j}}\lambda_{h_0}=-\lambda_{h_0}$, therefore $\lambda_{h_0}=0$. 

We then conclude that $x\in \op{Vect}_{\tilde{K}}((e_{x_i})_{i\in [\![1,r_{S'}]\!]})$, therefore 
$$\underset{s\in S'}\cap \op{ker}(\rho_{\Gamma_1}(T_s)+1)\subset \op{Vect}_{\tilde{K}}((e_{x_i})_{i\in [\![1,r_{S'}]\!]}).$$ The reverse inclusion follows from the definition of $r_{S'}$. 

This proves that $\dim(\underset{s\in S'}\cap \op{ker}(\rho_{\Gamma_1}(T_s)+1))= \vert \{x\in X,S'\subset I(x)\}\vert$.

\bigskip

We now have that for all $S'\subset S$, $\vert \{x\in X,S'\subset I_1(x)\}\vert=\vert \{x\in X,S'\subset I_2(x)\}\vert$ because the representations are isomorphic. It follows by induction on $\vert S\vert -\vert S'\vert $ that for all $S'\subset S$,  $\vert \{x\in X,S'= I_1(x)\}\vert=\vert \{x\in X,S'= I_2(x)\}\vert$ which concludes the proof. \end{proof}

The next proposition gives a fairly good uniqueness property which will be used for the computations on $W$-graphs.

\begin{prop2}\label{Unique2}
Let $\Gamma_1=(X,I,\mu)$ and $\Gamma_2=(X',I',\mu')$ be $W$-graphs as above and $\varphi:X\rightarrow X'$ be a bijection such that for all $x\in X$, $I'(\varphi(x))=I(x)$. We assume that $X$ and $X'$ are labeled in such a way that if $I(x_i)\subsetneq I(x_j)$ then $i<j$ and $\varphi(x_i)=x_i'$. We also assume that the images of $\mu$ and $\mu'$ are included in $\tilde{K}$ and that $\mu$ and $\mu'$ are independent of $S$. We then write $\mu_{x,y}$ instead of $\mu_{x,y}^s$.

If there exists $M\in GL_n(\tilde{K})$ such that for all $T\in \mathcal{H}_{\tilde{K}}$,
$M\rho_{\Gamma_1}(T)M^{-1}=\rho_{\Gamma_2}(T)$ then $M$ is block diagonal with blocks of length $\vert \{x\in X,I_1(x)=I_1(x_i)\}\vert$.

\end{prop2}

\begin{proof}
Note first that by Proposition \ref{Unique1}, we can choose a labeling as required.

We can choose numberings verifying desired by Proposition \ref{Unique1}. We write $(e_i)_{i\in [\![1,m]\!]}$ the basis corresponding to $\Gamma_1$ and $(e_i')_{i\in [\!1,m]\!]}$ the basis corresponding to $\Gamma_2$. Assume there exists a matrix $M\in GL_n(\tilde{K})$ such that for all $T\in \mathcal{H}_{\tilde{K}}$, $M\rho_{\Gamma_1}(T)M^{-1}=\rho_{\Gamma_2}(T)$. We will first show that the matrix $M$ is block lower-triangular and then that it is diagonal in the basis indexed by the numbering we chose on the vertices. We write $\rho_1=\rho_{\Gamma_1}$, $\rho_2=\rho_{\Gamma_2}$ and $m$ the dimension of both representations for the remainder of the proof.

Let $x_j\in X$. Let $s\in I(x_j)$. We then have $M\rho_1(T_s)(e_j)=-Me_j=-\underset{i=1}{\overset{m}\sum}m_{i,j}e_i'$. On the other hand, we have $\rho_2(T_s)Me_j=\underset{i=1}{\overset{m}\sum}m_{i,j}\rho_2(T_s)e_i'=-\underset{i\in [\![1,m]\!], s\in I(x_i)}\sum a_ie_i'+\underset{i\in [\![1,m]\!], s\notin I(x_i)}\sum\alpha_s m_{i,j}e_i'$. Since those two quantities are assumed to be equal, we have  $m_{i,j}=0$ for all $i\in [\![1,m]\!]$ such that $s\notin I(x_i)$.

Let now $(i,j)\in [\![1,m]\!]^2$ be such that $I(x_i)\neq I(x_j)$ and $i<j$. By the assumption on the numbering, we have that there exists $s\in I(x_j)$ such that $s\notin I(x_i)$. The above computation implies that $m_{i,j}=0$, therefore $M$ is block lower-triangular.

\bigskip

We have proven that for all $i\in [\![1,m]\!]$, $Me_i=\underset{j,I(x_i)\subset I(x_j)}\sum m_{i,j}e_j'$. Let $i\in [\![1,m]\!]$ and $j_0\in[\![1,m]\!]$ such that $I(x_i)\subsetneq I(x_{j_0})$ and $s\in I(x_{j_0})\setminus I(x_i)$. We have

\begin{eqnarray*}
M\rho_1(T_{s})(e_i) & =  & M(\alpha_s e_i+\sqrt{\alpha_s}\underset{s \in I(x_\ell)}{\underset{\ell}\sum} \mu_{\ell,i}e_i) \\
& = &  \alpha_s \underset{I(x_i)\subset I(x_j)}{\underset{j}\sum} m_{i,j}e_j'+\sqrt{\alpha_s}\underset{s\in I(x_\ell)}{\underset{\ell}\sum} \mu_{\ell,i}\underset{I(x_\ell)\subset I(x_k)}{\underset{k}\sum} m_{\ell,k}e_k'\\
& = & \alpha_s m_{i,j_0}e_{j_0}'+\alpha_s \underset{I(x_i)\subset I(x_j)}{\underset{j\neq j_0 }\sum} m_{i,j}e_j'+\sqrt{\alpha_s}\underset{s\in I(x_\ell)}{\underset{\ell}\sum} \underset{I(x_k)\subset I(x_\ell)}{\underset{k}\sum}\mu_{\ell,i}m_{\ell,k}e_k'\\
& =& \left(\alpha_s m_{i,j_0}+\sqrt{\alpha_s}\underset{s\in I(x_\ell)\subset I(x_{j_0})}{\underset{\ell}\sum}\mu_{\ell,i}m_{\ell,j_0}\right)e_{j_0}'+\underset{l\neq j_0}\sum b_\ell e_\ell'  .
\end{eqnarray*}
where the coefficients $b_\ell$ are elements of $K$ which can be deduced from the above equalities.

We also have 
\begin{eqnarray*}
\rho_2(s)M(e_i) & = & -\underset{s\in I(x_j)}{\underset{I(x_i)\subset I(x_j)}{\underset{j}\sum}} m_{i,j}e_j'+\underset{j,I(x_i)\subset I(x_j),s\notin I(x_j)}\sum m_{i,j}\left(\alpha_s e_j'+\sqrt{\alpha_s}\underset{\ell,s\in I(x_\ell)}\sum \mu'_{\ell,j}e_\ell'\right)\\
& =& \left(-m_{i,j_0}+\sqrt{\alpha_s}\underset{s\notin I(x_j)}{\underset{j}\sum}m_{i,j}\mu_{j_0,j}'\right)e_{j_0}'+\underset{\ell\neq j_0}\sum c_\ell e_\ell'.
\end{eqnarray*}
where the coefficients $c_\ell$ are elements of $K$ which can be deduced from the above equalities.

Since we have assumed those two quantities to be equal, we get that 
$$\alpha_s m_{i,j_0}+\sqrt{\alpha_s}\underset{s\in I(x_l)\subset I(x_{j_0})}{\underset{l}\sum}\mu_{l,i}m_{l,j_0}=-m_{i,j_0}+\sqrt{\alpha_s}\underset{s\notin I(x_j)}{\underset{j}\sum}m_{i,j}\mu_{j_0,j}'.$$

Since $1+\alpha_s$ and $\sqrt{\alpha_s}$ are $\tilde{K}$-linearly independant, we get that $m_{i,j_0}=0$. This is true for all $(i,j)$ such that $I(x_i)\subsetneq I(x_{j_0})$, therefore we have that the matrix $M$ is block-diagonal. \end{proof}

\textbf{Remark} : To prove that the matrix is block lower-triangular, we don't need to assume anything on the images of $\mu$ or $\mu'$ or on the independence with regards to $S$ and we can take $M\in GL_m(K)$. If we don't assume that the image of $\mu'$ is in $\tilde{K}$ then the result does not hold. The two following $W$-graphs give us a counter-example in the more general setting
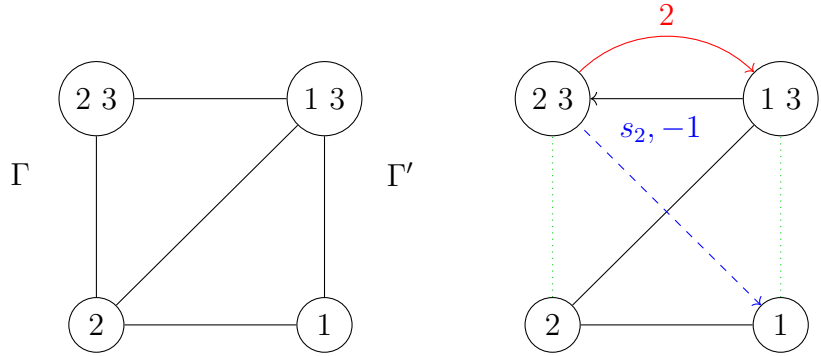
\begin{center}
\begin{tikzpicture}

\node[circle, draw =black] (4) at (0,0) {$2 ~ 3$};
\node[circle, draw =black] (3) at (3,0) {$1 ~ 3$};
\node[circle, draw =black] (2) at (0,-3) {$2$};
\node[circle, draw =black] (1) at (3,-3) {$1$};
\draw (1) to (2);
\draw (1) to (3);
\draw (2) to (3);
\draw (2) to (4);
\draw (3) to (4);
\draw (-1,-1) node{$\Gamma$};

\node[circle, draw =black] (4) at (6,0) {$2 ~ 3$};
\node[circle, draw =black] (3) at (9,0) {$1 ~ 3$};
\node[circle, draw =black] (2) at (6,-3) {$2$};
\node[circle, draw =black] (1) at (9,-3) {$1$};
\draw (1) to (2);
\draw[dotted] [-] [black!18!green](1) to (3);
\draw (2) to (3);
\draw[dotted] [-] [black!18!green](2) to (4);
\draw [<-] [red] (3) to [bend right =45] node[auto,swap] {$2$} (4);
\draw [->] (3) to (4);
\draw[dashed]  [<-] [blue] (1) to  node[auto,swap,pos=0.84] {$s_2,-1$} (4);
\draw (4,-1) node{$\Gamma'$};
\end{tikzpicture} 
\end{center}

The green dotted edges are of weight $1+\sqrt{\alpha}+\sqrt{\alpha}^{-1}$. The blue dashed edge indicates that $\mu^{s_2}_{4,1}=-1$ and $\mu^{s_1}_{4,1}=\mu^{s_3}_{4,1}=0$.

Then the matrix $M=\begin{pmatrix}
1 & 0 & 0 & 0\\
0 & 1 & 0 & 0\\
1 & 0 & 1& 0\\
0 & 1 & 0& 1\end{pmatrix}$ verifies for all $s\in S$, 
$$M\rho_{\Gamma}(T_s)M^{-1}=\rho_{\Gamma'}(T_s).$$

\bigskip

We now give a proposition which can be found in \cite{G-P} (11.1.7) providing a construction of a $W$-graph associated to the dual representation of a $W$-graph.

\begin{prop2}\label{Melyssaesttropforte}
Let $\Gamma=(X,I,\mu)$ be a $W$-graph. We define its dual $W$-graph $\Gamma^\star=(X,\tilde{I},\tilde{\mu})$ by 
\begin{enumerate}
\item $\forall x\in X, \tilde{I}(x)=S\setminus I(x)$,
\item $\forall (x,y,s)\in X\times X\setminus D(X)\times S, \tilde{\mu}_{x,y}^s=-\mu_{y,x}^s$.
\end{enumerate}
We then have that $\rho_{\Gamma}\simeq \rho_{\Gamma^\star}^\star$, where $\rho^\star(T_s) =\rho(-\alpha_s T_s^{-1})$.
\end{prop2}

\begin{Def2}\label{defselfdualrepresentation}
We say a representation $\rho$ of $\mathcal{H}_K$ is self-dual if $\rho\simeq \rho^\star$.
\end{Def2}

\begin{prop2}\label{existbilin}
If $\rho:\mathcal{H}_K\rightarrow GL_{n_\rho}(K)$ is an irreducible self-dual representation of a Iwahori-Hecke algebra then there is either a symmetric or skew-symmetric non-degenerate bilinear form $\langle .,.\rangle$ associated to $\rho$ in the following way
$$ \forall s\in S,\forall (u,v)\in V^2, \langle \rho(T_s)u,\rho(T_s)v\rangle =\langle u,v\rangle.$$

\end{prop2}

\begin{proof}
Let $\rho:\mathcal{H}_K\rightarrow GL_{n_\rho}(K)$ be a self-dual representation. Then there exists $P\in GL_n(k)$ such that for all $s\in S$, 
$$P\rho(T_s)P^{-1}=-\alpha_s~^t\rho(T_s)^{-1}.$$ This implies that $P\rho(T_s)P^{-1}=-\alpha_s~ ^t\! (-\alpha_s P^{-1}~ ^t\!\rho(T_s)^{-1}P)^{-1}= ^t\!P\rho(T_s)^t\!P^{-1}$. It follows that for all $s\in S$, $^t\! P^{-1} P\rho(T_s)(^t\! P^{-1} P)^{-1}=\rho(T_s)$. Hence, for all $h\in \mathcal{H}_K$, $^t\! P^{-1} P\rho(h)(^t\! P^{-1} P)^{-1}=\rho(h)$. By Schur's lemma, there exists $\lambda\in \overline{\F_q}$ such that $^t\! P^{-1} P=\lambda$. Thus $P=\lambda ^t\! P= \lambda ^t (\lambda ^t \! P)=\lambda^2 P$. It follows that $\lambda\in \{-1,1\}$.

This bilinear form $\langle . \vert . \rangle$ associated to $P$ is non-degenerate. It is symmetric when $\lambda=1$ and skew-symmetric when $\lambda=-1$. For all $s\in S, (u,v)\in V^2$, we have $\langle \rho(T_s)u,\rho(T_s)v\rangle =^t\! (\rho(T_s)u)P\rho(T_s)v=^t\! u (^t \rho(T_s) P\rho(T_s)) v= ^t\! u (-\alpha_s P)v=-\alpha_s \langle u,v\rangle$. This concludes the proof.
\end{proof}

Although such a bilinear form always exists for any given self-dual representation, it is difficult to obtain the bilinear form explicitely and to determine whether it is symmetric or anti-symmetric. Since there is a combinatorial way to define a dual $W$-graph, it is natural to expect $W$-graphs affording self-dual representations to be isomorphic to their dual $W$-graph. This seems to never be the case. Nevertheless, we have $(X,I,\mu)\simeq (X,\tilde{I},-\tilde{\mu})$ in many cases. We know by Proposition \ref{color} that if $(X,\tilde{I},-\tilde{\mu})$ is $2$-colorable, then we have $\rho_{(X,\tilde{I},-\tilde{\mu})}\simeq \rho_{(X,\tilde{I},\tilde{\mu})}$. When this is the case, we can define a bilinear form using only the $2$-coloring of the graph. The construction is given in the following theorem.

\begin{theo2}\label{bilinwgraphs}
Let $\Gamma=(X,I,\mu)$ be a $W$-graph affording an irreducible representation of $\mathcal{H}$ such that $\Gamma$ is $2$-colorable and $\Gamma$ is isomorphic as an oriented weighted graph to $(X,\tilde{I},-\tilde{\mu})$.

\smallskip

Let $\varphi:X\rightarrow X$ be the graph automorphism from $\Gamma$ to $(X,\tilde{I},-\tilde{\mu})$ and $x_1,x_2,\dots,x_n$ be a numbering of $X$ such that $\varphi(x_i)=x_{n+1-i}$.

\smallskip

Let $\langle .,.\rangle$ be the bilinear form defined by $\langle e_{x_i},e_{x_j}\rangle=\omega(e_{x_i})\delta_{i,n+1-j}$, where $\omega$ corresponds to a coloring of $\Gamma$ with $1$ and $-1$.

\smallskip

We then have 
$$\forall s\in S,\forall v_1,v_2\in V, \langle \rho_\Gamma(T_s)v_1,\rho_\Gamma(T_s)v_2\rangle =-\alpha\langle v_1,v_2\rangle.$$

\smallskip

This bilinear form is non-degenerate and it is symmetric if $\omega(x_1)\omega(x_n)=1$ and skew-symmetric if $\omega(x_1)\omega(x_n)=-1$.

The associated representation is then self-dual.
\end{theo2}

\begin{proof}
Let $\Gamma=(X,I,\mu)$ be a $W$-graph as above and $\omega:X\rightarrow \{-1,1\}$ be the corresponding 2-coloring. Let $s\in S$ and $x,y\in V$. First note that since $I(\varphi(x))=S\setminus I(x)$, we have $\langle e_x,e_z\rangle =0$ for all $z\in X$ such that $I(z)\neq S\setminus I(x)$. We now have four different cases to consider.

If $s\in I(x)\cap I(y)$ then $\langle \rho_\Gamma(T_s) e_x, \rho_\Gamma(T_s)e_y\rangle = \langle -e_x,-e_y\rangle = \langle e_x,e_y\rangle=0=-\alpha_s\langle e_x,e_y\rangle$.

If $s\in I(x)$, $s\notin I(y)$ then 
\begin{eqnarray*}
\langle \rho_\Gamma(T_s)e_x,\rho_\Gamma(T_s)e_y\rangle & =& \langle -e_x,\alpha_s e_y+\sqrt{\alpha_s}\underset{z\in X,s\in I(z)}\sum \mu_{z,y} e_z\rangle\\
& = & -\alpha_s\langle e_x,e_y\rangle -\sqrt{\alpha_s}\underset{z\in X,s\in I(z)}\sum \mu_{z,y} \langle e_x,e_z\rangle\\
& = & -\alpha_s\langle e_x,e_y\rangle .
\end{eqnarray*}

If $s\notin I(x)$, $s\in I(y)$ then 
\begin{eqnarray*}
\langle \rho_\Gamma(T_s)e_x,\rho_\Gamma(T_s)e_y\rangle & =& \langle \alpha_s e_x+\sqrt{\alpha_s}\underset{z\in X,s\in I(z)}\sum \mu_{z,x} e_z,- e_y\rangle\\
& = & -\alpha_s\langle e_x,e_y\rangle -\sqrt{\alpha_s}\underset{z\in X,s\in I(z)}\sum \langle e_z,e_y\rangle\\
& = & -\alpha_s\langle e_x,e_y\rangle .
\end{eqnarray*}

If $s\notin I(x)\cup I(y)$ then we have $\langle e_x,e_y\rangle =0$ and
\begin{eqnarray*}
\langle \rho_\Gamma(T_s)e_x,\rho_\Gamma(T_s)e_y\rangle & =& \langle \alpha_s e_x+\sqrt{\alpha_s}\underset{z\in X,s\in I(z)}\sum \mu_{z,x} e_z,\alpha_s e_y+\sqrt{\alpha_s}\underset{z'\in X,s\in I(z')}\sum \mu_{z',x} e_{z'}\rangle\\
& = & \sqrt{\alpha_s}\underset{z\in X,s\in I(z)}\sum \mu_{z,x} \langle e_z,e_y\rangle+\sqrt{\alpha_s}\underset{z'\in X,s\in I(z')}\sum \mu_{z',y}\langle e_x,e_{z'}\rangle \\
& = & \sqrt{\alpha_s} \mu_{\varphi(y),x}\omega(\varphi(y))+\sqrt{\alpha_s}\mu_{\varphi(x),y}\omega(x)\\
& =& \sqrt{\alpha_s}(\mu_{\varphi(y),\varphi(\varphi(x))}\omega(\varphi(y))+\mu_{\varphi(x),y}\omega(x))\\
& =& \sqrt{\alpha_s}(\mu_{\varphi(x),y}\omega(\varphi(y))+\mu_{\varphi(x),y}\omega(x)\\
& =& \sqrt{\alpha_s}\mu_{\varphi(x),y}(\omega(\varphi(y))+\omega(x)).
\end{eqnarray*}
If $\mu_{\varphi(x),y}=0$ then the above quantity is equal to zero. If $\mu_{\varphi(x),y}\neq 0$ then $\mu_{\varphi(y),x}\neq 0$, therefore $\omega(\varphi(y))=-\omega(x)$ and the above quantity is again zero.

This proves that in all cases, $\langle \rho_\Gamma(T_s)e_x,\rho_\Gamma(T_s)e_y\rangle=-\alpha_s\langle e_x,e_y\rangle$. It only remains to show that this bilinear form is symmetric if $\omega(x_1)\omega(x_n)=1$ and skew-symmetric if $\omega(x_1)\omega(x_n)=-1$. To show this, we first prove that for all $x\in X$, we have $\omega(x)\omega(\varphi(x))=\omega(x_1)\omega(x_n)$. If $\Gamma$ was a disconnected graph, then $\rho_\gamma$ would be reducible, therefore $\Gamma$ is connected. This means we only need to show that for a given $x$, we have $\omega(y)\omega(\varphi(y))=\omega(x)\omega(\varphi(x))$ for all $y\in X$ such that $\mu_{x,y}\neq 0$ or $\mu_{y,x}\neq 0$, . Let $x\in X$ and $y\in X$ such that $\mu_{x,y}\neq 0$ or $\mu_{y,x}\neq 0$. Then $\mu_{\varphi(y),\varphi(x)}\neq 0$ or $\mu_{\varphi(x),\varphi(y)}\neq 0$. This implies that $\omega(y)=-\omega(x)$ and $\omega(\varphi(y))=-\omega(\varphi(x))$, therefore $\omega(y)\omega(\varphi(y))=(-\omega(x))(-\omega(\varphi(x))=\omega(x)\omega(\varphi(x))$. This shows that for all $x\in X$, we have $\omega(x)\omega(\varphi(x))=\omega(x_1)\omega(x_n)$. Let now $(x,y)\in X^2$. We have $\langle e_x,e_y\rangle =\omega(x)\delta_{y,\varphi(x)}=\omega(x_1)\omega(x_n)\omega(\varphi(x))\delta_{\varphi(y),x}=\omega(x_1)\omega(x_n)\omega(y)\delta_{x,\varphi(y)}=\omega(x_1)\omega(x_n)\langle e_y,e_x\rangle$ as required.
\end{proof}

\begin{Def2}\label{defselfdualwgraph}
A $W$-graph satisfying the conditions in Theorem \ref{bilinwgraphs} is said to be self-dual.
\end{Def2}

Not all $W$-graphs affording self-dual representations are self-dual. For example, in type $D_4$, there is an irreducible self-dual eight-dimensional representation which corresponds to the double Young diagram $\lambda=([2,1],[1])$. By Proposition \ref{bilin2}, we know that the associated bilinear form must be skew-symmetric since $\tilde{\nu}(\lambda)=\nu([2,1])\nu([1])=(-1)^\frac{3-1}{2}(-1)^\frac{1-1}{2}=-1$. We can associate to this representation the following $W$-graph which we determined using the restrictions to the three different $A_{A_3}$ and $A_{A_1}^3$

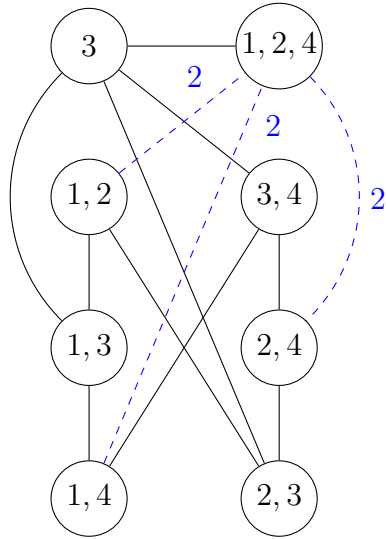
\begin{center}
\begin{tikzpicture}
[place/.style={circle,draw=black,
inner sep=1pt,minimum size=10mm}]
\node (1) at (1,0)[place] {$3$};
\node (2) at (1,-2)[place]{$1,2$};
\node (3) at (1,-4) [place]{$1,3$};
\node (4) at (1,-6) [place]{$1,4$};
\node (8) at (3.5,0) [place]{$1,2,4$};
\node (7) at (3.5,-2) [place]{$3,4$};
\node (6) at (3.5,-4) [place]{$2,4$};
\node (5) at (3.5,-6) [place]{$2,3$};
\draw [-] (1) to [bend right =45] (3);
\draw[dashed] [-] [blue](8) to [bend left =45]node [auto] {2} (6);
\draw [-] (1) to (5);
\draw[dashed] [-][blue] (4) to node [auto,swap,pos =0.95] {2} (8);
\draw [-] (1) to (7);
\draw[dashed] [-] [blue] (2) to  node [auto,near end] {2} (8);
\draw [-] (2) to (5);
\draw [-] (7) to (4);
\draw [-] (2) to (3);
\draw [-] (7) to (6);
\draw [-] (3) to (4);
\draw [-] (6) to (5);
\draw [-] (1) to (8);
\end{tikzpicture}
\end{center}

Here the $W$-graphs are presented as in \cite{G-P}. If $I(x)=\{s_3\}$, then we write $3$ inside the vertex $x$. We don't write the edges when $\mu_{x,y}=0$. If $\mu_{x,y}=\mu_{y,x}$, then we write non-oriented edges. If $\mu_{x,y}=1$, then we do not write the weight on the edge and we write $\mu_{x,y}$ as a weight on the edge otherwise. The $W$-graph does not satisfy the desired properties because the edge between nodes $3$ and $3,4$ is of weight $1$ whereas the edge between notes $1,2$ and $1,2,4$ is of weight $2$.

\smallskip

However, since the properties appear to be natural and are verified by some $W$-graphs, we propose the following conjecture.

\begin{conj2}\label{conjecturewgraphs}
Let $W$ be a Coxeter group. For any irreducible self-dual representation, there exists a self-dual $W$-graph defined over $K$. If there exists a $W$-graph $\Gamma=(X,I,\mu)$ defined over $K'$ then there exist a self-dual $W$-graph $\Gamma'=(X',I',\mu')$ defined over $K'$ and a matrix $M\in GL_{\vert X\vert}(\tilde{K})$ such that for all $h\in \mathcal{H}_{K}$, $M\rho_{\Gamma}(h)M^{-1}=\rho_{\Gamma'}(h)$.
\end{conj2}

In order to prove the conjecture for exceptional groups, we effectively find the $W$-graphs verifying the right properties by assuming they exist and finding the bilinear form preserving the $W$-graph we are working with. This can be seen in the following proposition.

\begin{prop2}
Assume the conjecture holds and that $\alpha_s=\alpha_{s'}$ for all $s,s'\in S$. We write $\alpha$ for the unique parameter of the Iwahori-Hecke algebra. Then any bilinear form associated to a $W$-graph affording an irreducible self-dual representation defined over $K'$ is represented by a block anti-diagonal matrix.
\end{prop2}

\begin{proof}
Assume the conjecture is true and that $\Gamma=(X,I,\mu)$ is a self-dual $W$-graph defined over $K'$. There exists $\Gamma'=(X',I',\mu')$ such that $\mu'(X\times X\setminus D(X)\times S)\subset K'$ and $\rho_{\Gamma}\simeq \rho_{\Gamma'}$. We order the vertices of $X$ in order to have $I(x_i)=I(x'_i)$ for all $i$, and if $I(x_i)\subsetneq I(x_j)$ then $i< j$. We consider the matrices with respect to the bases corresponding to those orders. 

There exists $M\in GL_{\vert X\vert}(\tilde{K})$ such that for all $h\in \mathcal{H}_{K}$, $M\rho_{\Gamma}(h)M^{-1}=\rho_{\Gamma'}(h)$. By Propositions \ref{Unique1} and \ref{Unique2}, $M$ is block-diagonal. Since $\Gamma'$ is self-dual, there exists an anti-diagonal matrix $\tilde{P}$ corresponding to its $2$-coloring such that for all $s\in S$, $\tilde{P}\rho_{\Gamma'}(T_s)\tilde{P}^{-1}=-\alpha~ ^t\!\rho_{\Gamma'}(T_s)^{-1}$. We also have by Proposition \ref{existbilin} that there exists $P\in GL_{\vert X\vert}(K)$ such that for all $s\in S$, $P\rho_{\Gamma}(T_s)P^{-1}=-\alpha ^t\! \rho_{\Gamma}(T_s)^{-1}$.

By substituting $M\rho_{\Gamma}(T_s)M^{-1}=\rho_{\Gamma'}(T_s)$ in the first expression, we get $\tilde{P}M\rho_{\Gamma}(T_s)M^{-1}\tilde{P}^{-1}=-\alpha ^t\!(M\rho_{\Gamma}(T_s)M^{-1})^{-1}=-\alpha ^t\!M^{-1}~ ^t\!\rho_{\Gamma}(T_s)^{-1}~^t\! M$.
 It follows that $^t\! M\tilde{P}M\rho_{\Gamma}(T_s)(^t\!M\tilde{P}M)^{-1}=-\alpha ^t\! \rho_{\Gamma}(T_s)^{-1}=P\rho_{\Gamma}(T_s)P^{-1}$. Hence, by Schur's lemma, there exists $\lambda\in \overline{\F_p}$ such that $P=\lambda ^t\! M\tilde{P}M$. Since $M$ is block-diagonal and $\tilde{P}$ is anti-diagonal, we get that $P$ is block anti-diagonal.
\end{proof}

\textbf{Remark} : This proposition can be useful to to find the matrix $P$ after assuming the conjecture is true. For the $W$-graph associated to $([2,1],[1])$, all the blocks are of size one, therefore we only need to look for an anti-diagonal matrix. This means we have only $8$ unknowns, assuming the $2$-coloring is preserved up to permutation of vertices having the same image under $I$, we get that $P$ must be anti-symmetric, therefore we only have four unknowns.

 After solving the equations afforded by each generator, we get $$P=\begin{pmatrix} 0 & 0 & 0& 0& 0 & 0 & 0 & 1\\
0 & 0 & 0& 0& 0 & 0 & 2 & 0\\
0 & 0 & 0& 0& 0 & -2 & 0 & 0\\
0 & 0 & 0& 0& 2 & 0 & 0 & 0\\
0 & 0 & 0& -2& 0 & 0 & 0 & 0\\
0 & 0 & 2& 0& 0 & 0 & 0 & 0\\
0 & -2 & 0& 0& 0 & 0 & 0 & 0\\
-1 & 0 & 0& 0& 0 & 0 & 0 & 0\\
\end{pmatrix}$$

We then look for a matrix $M$ such that $^t\!M \tilde{P} M=P$.

Using Gaussian reduction for quadratic forms, we get
$$M=\begin{pmatrix}1 & 0 & 0& 0& 0 & 0 & 0 & 0\\
0 & 2 & 0& 0& 0 & 0 & 0 & 0\\
0 & 0 & 2& 0& 0 & 0 & 0 & 0\\
0 & 0 & 0& 2& 0 & 0 & 0 & 0\\
0 & 0 & 0& 0& 1 & 0 & 0 & 0\\
0 & 0 & 0& 0& 0 & 1 & 0 & 0\\
0 & 0 & 0& 0& 0 & 0 & 1 & 0\\
0 & 0 & 0& 0& 0 & 0 & 0 & 1\\
\end{pmatrix}$$

By computing $M\rho(T_s)M^{-1}$ for all $s$, we obtain the following $W$-graph.

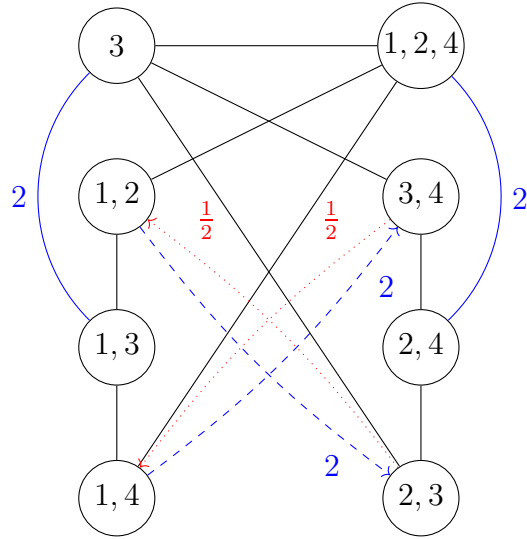
\begin{center}
\begin{tikzpicture}
[place/.style={circle,draw=black,
inner sep=1pt,minimum size=10mm}]
\node (9) at (6,0)[place] {$3$};
\node (10) at (6,-2)[place]{$1,2$};
\node (11) at (6,-4) [place]{$1,3$};
\node (12) at (6,-6) [place]{$1,4$};
\node (16) at (10,0) [place]{$1,2,4$};
\node (15) at (10,-2) [place]{$3,4$};
\node (14) at (10,-4) [place]{$2,4$};
\node (13) at (10,-6) [place]{$2,3$};
\draw [-] [blue](9) to [bend right =45]node [auto,swap] {2} (11);
\draw [-] [blue](16) to [bend left =45]node [auto] {2} (14);
\draw [-] (9) to (13);
\draw [-] (12) to (16);
\draw [-] (9) to (15);
\draw [-] (10) to  (16);
\draw[dashed] [->] [blue](10) to [bend right =8] node [auto,swap, very near end] {$2$}(13);
\draw[dotted] [->] [red](13) to [bend right =8] node[auto,swap, very near end] {$\frac{1}{2}$}(10);
\draw[dotted] [->] [red](15) to[bend right =8] node[auto,swap, very near start] {$\frac{1}{2}$} (12);
\draw[dashed] [->] [blue](12) to [bend right =8] node[auto,swap,very near end] {$2$}(15);
\draw [-] (10) to (11);
\draw [-] (15) to (14);
\draw [-] (11) to (12);
\draw [-] (14) to (13);
\draw [-] (9) to (16);
\end{tikzpicture}
\end{center}

We obtain in the same way self-dual $W$-graphs for all self-dual representations in types $I_2(m),H_3,H_4,E_6,E_7$ and $E_8$. We give the new $W$-graphs in type $H_4$ in the Appendix. All the new $W$-graphs can be downloaded from \cite{newgraphsEsterle}.

\chapter{Type E}\label{TypeEchapter}

In this section we determine the images of the Artin groups of type $E$ inside their associated Iwahori-Hecke algebras. Their representations are given by $W$-graphs, therefore we will use the results from Chapter \ref{Wgraphschapter}. We have defined in that chapter self-dual representations in Definition \ref{defselfdualrepresentation} and self-dual $W$-graphs in Definition \ref{defselfdualwgraph}. We have established the Conjecture \ref{conjecturewgraphs} which states that there exists a self-dual $W$-graph affording any self-dual representation of a $W$-graph. There are no self-dual representations in type $E_7$. In types $E_6$ and $E_8$, we have proved the conjecture and obtained new $W$-graphs which are self-dual for each of the self-dual representations. We only provide in this Appendix the graphs for the $10$-dimensional self-dual representation of $E_6$ and its $20$-dimensional self-dual representation. All the remaining ones are of dimension greater than $60$ and can be downloaded from \cite{newgraphsEsterle}.

\section{Type $E_6$}\label{E6section}
Let $p\notin \{2,3\}$ be a prime and $\alpha\in \overline{\F_p}$ be of order not dividing $5,8,9$ or $12$. We write $\F_q=\F_p(\alpha)$. There are $25$ irreducible representations of the Iwahori-Hecke algebra $\mathcal{H}_{E_6,\alpha}$ which we define below. They are all of dimension less than or equal to $90$ and there are $5$ self-dual representations, they associated to the $E_6$-graphs $10_s$, $20_s$, $60_s$, $80_s$ and $90_s$. We have found self-dual $E_6$-graphs \cite{newgraphsEsterle} for each of these representations. Using the $2$-coloring, the bilinear form defined in Theorem \ref{bilinwgraphs} is skew-symmetric for each of those self-dual representations.

\begin{Def}
The Iwahori-Hecke algebra $\mathcal{H}_{E_6,\alpha}$ of type $E_6$ is the $\F_q$-algebra generated by $S_1,S_2,S_3,S_4,S_5,S_6$ and the following relations 
\begin{enumerate}
\item $\forall i\in \{1,2,3,4,5,6\}$, $(S_i-\alpha)(S_i+1)=0$.
\item $S_1S_3S_1=S_3S_1S_3$.
\item $\forall i \in \{2,4,5,6\}$, $S_1S_i=S_iS_1$.
\item $S_2S_4S_2=S_4S_2S_4$.
\item $\forall i\in \{3,5,6\}$, $S_2S_i=S_iS_2$.
\item $\forall i\in \{3,4,5\}$, $S_{i}S_{i+1}S_{i}=S_{i+1}S_iS_{i+1}$.
\item $S_3S_5=S_5S_3$.
\item $S_3S_6=S_6S_3$.
\item $S_4S_6=S_6S_4$.
\end{enumerate}
For $\sigma$ in the Coxeter group $E_6$, we set $T_{\sigma}=S_{i_1}\dots S_{i_k}$ whenever $\sigma=s_{i_1}\dots s_{i_k}$ is a reduced expression.
\end{Def}

This means we consider $E_6$ as in the CHEVIE package of GAP3 \cite{CHEVIE} with the following Dynkin diagram

\begin{center}
\begin{tikzpicture}
[place/.style={circle,draw=black,
inner sep=1pt,minimum size=10mm}]
\node (1) at (0,0)[place]{$S_1$};
\node (2) at (2,0)[place]{$S_3$};
\node (3) at (4,0)[place]{$S_4$};
\node (4) at (4,2)[place]{$S_2$};
\node (5) at (6,0)[place]{$S_5$};
\node (6) at (8,0)[place]{$S_6$};
\draw (1) to (2);
\draw (2) to (3);
\draw (3) to (4);
\draw (3) to (5);
\draw (5) to (6);
\end{tikzpicture}
\end{center}

\begin{prop}
Under our assumptions on $p$ and $\alpha$, $\mathcal{H}_{E_6,\alpha}$ is split semisimple, the representations afforded by the $W$-graphs are irreducible and pairwise non-isomorphic over $\F_q$. The restrictions of the irreducible representations of $\mathcal{H}_{E_6,\alpha}$ to $\mathcal{H}_{D_5,\alpha}$ are the same as in the generic case.
\end{prop}

\begin{proof}
We will apply Proposition \ref{Tits}. Let $A=\Z[\sqrt{u}^{\pm 1}]$ and $F=\Q(\sqrt{u})$. We have a symetrizing trace defined by $\tau(T_0)=1$ and $\tau(T_{\sigma})=0$ for all $\sigma\in E_6\setminus \{1_{E_6}\}$. By \cite{G-P}, $\mathcal{H}_{E_6,u}$ is then a free symmetric $F$-algebra of rank $51840$. By \cite{Bourb} V.3. Corollary 1, $A$ is integrally closed. Let $\theta$ be the ring homomorphism from $A$ to $L=\F_q$ defined by $\theta(u)=\alpha$ and $\theta(k)=\overline{k}$. We know $FH$ is split. The basis formed by the elements $T_\sigma$, $\sigma\in E_6$ verifies the conditions of the Proposition \ref{Tits}. The $E_6$-graphs are still connected after specialization since all the weights are in $\{-6,-4,-3,-2,-3/2,-1,-1/2,-1/3,-1/6,1/3,1/2,1,3/2,2,3\}$.

It now only remains to check that the Schur elements associated to the specialized representations are in $B$ and do not vanish under $\theta$ with $B$ as in Proposition \ref{Tits}. The Schur elements are given in Table \ref{SchurelmtsE6}. For a pair $(\rho,\rho^\star)$ of representations, we only give the Schur element of one of the representations since the other is obtained by applying the involution $\sqrt{u}\mapsto \sqrt{u}^{-1}$. The conditions on $\alpha$ and $p$ imply that the Schur elements verify the right conditions and the proof is concluded.

\begin{table}
\centering
$\begin{array}{lcr}
 1_p & : & (\Phi_2^4\Phi_3^3\Phi_4^2\Phi_5\Phi_6^2\Phi_8\Phi_9\Phi_{12})(u).\\
 6_p & : & \frac{1}{u}(\Phi_2^4\Phi_3^3\Phi_4^2\Phi_5\Phi_6^2\Phi_{12})(u).\\
 10_s & : & \frac{3}{u^7}(\Phi_2^4\Phi_3^3\Phi_4^2)(u).\\
 15_p & : & \frac{2}{u^3}(\Phi_2^4\Phi_3^3\Phi_4^2\Phi_{12})(u).\\
 15_q & : &  \frac{2}{u^3}(\Phi_2^4\Phi_3^3\Phi_4^2\Phi_6^2)(u).\\
 20_p & : & \frac{1}{u^2}(\Phi_2^4\Phi_3^3\Phi_4\Phi_6^2\Phi_9)(u).\\
 20_s & : & \frac{6}{u^7}(\Phi_2^4\Phi_3^3\Phi_{12})(u).\\
 24_p & : & \frac{1}{u^6}(\Phi_2^4\Phi_3^3\Phi_5\Phi_6^2)(u).\\
 30_p & : & \frac{2}{u^3}(\Phi_2^4\Phi_3^3\Phi_6^2\Phi_8)(u).\\
 60_s & : & \frac{2}{u^7}(\Phi_2^4\Phi_3^3\Phi_6^2)(u).\\
 60_p & : & \frac{1}{u^5}(\Phi_2^4\Phi_3^3\Phi_4\Phi_6^2)(u).\\
 64_p & : & \frac{1}{u^4}(\Phi_2\Phi_3^3\Phi_5\Phi_9)(u).\\
 80_s & : & \frac{6}{u^7}(\Phi_3^3\Phi_4^2\Phi_6^2)(u).\\
 81_p & : & \frac{1}{u^6}(\Phi_2^4\Phi_4^2\Phi_5\Phi_8)(u).\\
 90_s & : & \frac{3}{u^7}(\Phi_2^4\Phi_4^2\Phi_9)(u).
\end{array}$
\caption{Schur elements in type $E_6$}
\label{SchurelmtsE6}
\end{table}

\end{proof}

Using the CHEVIE package of GAP3\cite{CHEVIE}, we give in Table \ref{resE6D5} the restriction table from $\mathcal{H}_{E_6,\alpha}$ to its subalgebra $\mathcal{H}_{D_5,\alpha}$ generated by $S_1,S_2,S_3,S_4$ and $S_5$ which is naturally isomorphic to the Iwahori-Hecke of type $D_5$ with parameter $\alpha$. They correspond in the generic case to the induction/restriction tables of the corresponding Coxeter groups.

\begin{table}
\centering
\tabcolsep=0.1cm
\begin{tabular}{ |p{0.65cm}||p{0.65cm}|p{0.65cm}|p{0.65cm}|p{0.65cm}|p{0.65cm}|p{0.65cm}|p{0.65cm}|p{0.65cm}|p{0.65cm}|p{0.65cm}| p{0.65cm}|p{0.65cm}|p{0.65cm}|p{0.65cm}|p{0.65cm}|p{0.65cm}|p{0.65cm}|p{0.65cm}|  }
 \hline
 & $1$ & $1'$ & $4$ & $4'$ & $5_1$ & $5_1'$ & $5_2$ & $5_2'$  & $6$ & $10_1$ & $10_1'$ & $10_2$ & $10_2'$ & $10_3$ & $15$ & $15'$ & $20$ & $20'$ \\
 \hline

 $1_p$   & $1$    &  &    &  &  &  &  &  &  & &  &    &  &  &  &  &  &  \\
 \hline
 
 $6_p$   & $1$    &  &    &  &  &  & $1$ &  &  & &  &    &  &  &  &  &  &  \\
 \hline
 $10_s$   &     &  &    &  &  &  &  &  &  & &  &    &  & $1$ &  &  &  &  \\
 \hline
 $15_p$   &     &  &    &  &  &  & $1$ &  &  & $1$&  &    &  &  &  &  &  &  \\
 \hline
 $15_q$   &     &  &    &  & $1$ &  &  &  &  & &  &   $1$ &  &  &  &  &  &  \\
 \hline
 $20_p$   & $1$    &  &   $1$ &  &  &  & $1$ &  &  & &  &   $1$ &  &  &  &  &  &  \\
 \hline
 $20_s$   &     &  &    &  &  &  &  &  &  & $1$& $1$ &    &  &  &  &  &  &  \\
 \hline
 $24_p$   &     &  &   $1$ &  &  &  &  &  &  & &  &    &  &  &  &  & $1$ &  \\
 \hline
 $30_p$   &     &  &    &  &  &  & $1$ &  &  & &  &  $1$  &  &  & $1$ &  &  &  \\
 \hline
 
 $60_s$   &     &  &    &  &  $1$& $1$ &  &  &  & &  &    &  & $1$ &  &  &  $1$&  $1$\\
 \hline
 $60_p$   &     &  &    &  & $1$  &  &  &  &  & &  & $1$   &  & $1$ & $1$ &  & $1$ &  \\
 \hline
 $64_p$   &     &  &  $1$  &  &  &  & $1$ &  &  & $1$ &  &  $1$  &  &  & $1$ &  & $1$ &  \\
 \hline
 $80_s$   &     &  &    &  &  &  &  &  &  & &  &    &  & $1$ & $1$ & $1$ & $1$ & $1$ \\
 \hline
 $81_p$   &     &  &    &  &  &  &  &  & $1$ & $1$ &  & $1$   &  &  & $1$ &  & $1$ & $1$ \\
 \hline
 $90_s$   &     &  &    &  &  &  &  &  &  & $1$ & $1$  &    &  &  & $1$ & $1$ & $1$ & $1$ \\
 \hline
 
 \end{tabular}
 \smallskip
 
 $1:([5],\emptyset)$, $4:([4,1],\emptyset)$, $5_1:([3,2],\emptyset)$, $5_2:([4],[1])$, $6:([3,1^2])$, $10_1:([3],[1^2])$, $10_2:([3],[2])$, $10_3:([2,2],[1])$, $15:([3,1],[1])$ and $20 :[2,1],[2])$.
 \caption{Restriction table from $\mathcal{H}_{E_6,\alpha}$ to $\mathcal{H}_{D_5,\alpha}$.}
\label{resE6D5}
\end{table}

\begin{prop}\label{resE6derivedsubgroup}
The restrictions to $\mathcal{A}_{E_6}$ of the representations afforded by those $E_6$-graphs are absolutely irreducible and the representations of dimension greater than one are pairwise non-isomorphic.
\end{prop}

\begin{proof}
As in \cite{BMM} Lemma $3.4$, we only need to prove that $A_{E_6}$ is generated by $A_{E_5}$ and $\mathcal{A}_{E6}$. This true because $s_6=s_6s_1^{-1}s_1$, $s_6s_1^{-1}\in \mathcal{A}_{E_6}$ and $s_1\in A_{D_5}$.

\smallskip

Let $\rho_1$ and $\rho_2$ be two irreducible representations of $\mathcal{H}_{E_6,\alpha}$ such that $\rho_{1|\mathcal{A}_{E_6}}\simeq \rho_{2|\mathcal{A}_{E_6}}$. We have $A_{E_6}/\mathcal{A}_{E_6}\simeq <\overline{S_1}>\simeq \Z$. It follows by Lemma \ref{abel} that there exists $x\in \F_q^\star$ such that for all $i\in [\![1,6]\!]$, we have $\rho_1(S_i)$ conjugate to $x\rho_2(S_i)$. Since the representations are irreducible we have that the set of eigenvalues of $\rho_{1}(S_i)$ is $\{-1,\alpha\}$ and the set of eigenvalues of $\rho_2(S_i)$ is $\{-x,x\alpha\}$. This implies that $x=1$ or ($x=-\alpha$ and $\alpha^2=1$). The latter contradicts our assumptions on $\alpha$ therefore $x=1$ and $\rho_1\simeq \rho_2$.
\end{proof}

Before determining the image of the Artin groups inside this Iwahori-Hecke algebra, we need as in the other cases a Lemma on Artin groups which will allow us to use the restriction from $\mathcal{H}_{E_6,\alpha}$ to $\mathcal{H}_{D_5,\alpha}$.

\begin{lemme}\label{normclosD5E6}
The normal closure $\ll \mathcal{A}_{D_5} \gg_{\mathcal{A}_{E_6}}$ of $\mathcal{A}_{D_5}$ inside $\mathcal{A}_{E_6}$ is equal to $\mathcal{A}_{E_6}$, where we identify $A_{D_5}$ as a subgroup of $A_{E_6}$ using the natural isomorphism from $A_{D_5}$ to $<S_i,i\in [\![1,5]\!]>$.
\end{lemme}

\begin{proof}
By \cite{MR}, we have $\mathcal{A}_{D_5}=<S_3S_1^{-1},S_1S_3S_1^{-1},S_4S_1^{-1},S_2S_1^{-1}S_5S_1^{-1},S_3S_1^{-1},S_4S_3^{-1}>$ and $\mathcal{A}_{E_6}=<S_3S_1^{-1},S_1S_3S_1^{-1},S_4S_1^{-1},S_2S_1^{-1}S_5S_1^{-1},S_3S_1^{-1},S_4S_3^{-1},S_6S_1^{-1}>$. This proves we only need to show that $S_6S_1^{-1}\in \ll \mathcal{A}_{D_5} \gg_{\mathcal{A}_{E_6}}$.
We have $S_6S_5S_6=S_5S_6S_5$, therefore $S_6=S_5S_6S_5(S_5S_6)^{-1}$ and 
$$S_6S_1^{-1}=S_5S_6S_5S_1^{-1}(S_5S_6)^{-1}=((S_5S_1^{-1})(S_6S_1^{-1}))(S_5S_1^{-1})((S_5S_1^{-1})(S_6S_1^{-1}))^{-1}.$$
This concludes the proof.
\end{proof}

We will now use the above information and the usual techniques to determine the image of the Artin group inside $\mathcal{H}_{E_6,\alpha}$.

\begin{prop}\label{ImageinsiderepsE6}
If $\F_q=\F_p(\alpha)=\F_p(\alpha+\alpha^{-1})$ then we have
\begin{enumerate}
\item $\rho_{6_p}(\mathcal{A}_{E_6})\simeq SL_6(q)$,
\item $\rho_{10_s}(\mathcal{A}_{E_6})\simeq SP_{10}(q)$,
\item $\rho_{15_p}(\mathcal{A}_{E_6})\simeq SL_{15}(q)$,
\item $\rho_{15_q}(\mathcal{A}_{E_6})\simeq SL_{15}(q)$,
\item $\rho_{20_p}(\mathcal{A}_{E_6})\simeq SL_{20}(q)$,
\item $\rho_{20_s}(\mathcal{A}_{E_6})\simeq SP_{20}(q)$,
\item $\rho_{24_p}(\mathcal{A}_{E_6})\simeq SL_{24}(q)$,
\item $\rho_{30_p}(\mathcal{A}_{E_6})\simeq SL_{30}(q)$,
\item $\rho_{60_s}(\mathcal{A}_{E_6})\simeq SP_{60}(q)$,
\item $\rho_{60_p}(\mathcal{A}_{E_6})\simeq SL_{60}(q)$,
\item $\rho_{64_p}(\mathcal{A}_{E_6})\simeq SL_{64}(q)$,
\item $\rho_{80_s}(\mathcal{A}_{E_6})\simeq SP_{80}(q)$,
\item $\rho_{81_p}(\mathcal{A}_{E_6})\simeq SL_{81}(q)$,
\item $\rho_{90_s}(\mathcal{A}_{E_6})\simeq SP_{90}(q)$.
\end{enumerate}

If $\F_q=\F_p(\alpha)\neq\F_p(\alpha+\alpha^{-1})$ then we have
\begin{enumerate}
\item $\rho_{6_p}(\mathcal{A}_{E_6})\simeq SU_6(q^{\frac{1}{2}})$,
\item $\rho_{10_s}(\mathcal{A}_{E_6})\simeq SP_{10}(q^{\frac{1}{2}})$,
\item $\rho_{15_p}(\mathcal{A}_{E_6})\simeq SU_{15}(q^{\frac{1}{2}})$,
\item $\rho_{15_q}(\mathcal{A}_{E_6})\simeq SU_{15}(q^{\frac{1}{2}})$,
\item $\rho_{20_p}(\mathcal{A}_{E_6})\simeq SU_{20}(q^{\frac{1}{2}})$,
\item $\rho_{20_s}(\mathcal{A}_{E_6})\simeq SP_{20}(q^{\frac{1}{2}})$,
\item $\rho_{24_p}(\mathcal{A}_{E_6})\simeq SU_{24}(q^{\frac{1}{2}})$,
\item $\rho_{30_p}(\mathcal{A}_{E_6})\simeq SU_{30}(q^{\frac{1}{2}})$,
\item $\rho_{60_s}(\mathcal{A}_{E_6})\simeq SP_{60}(q^{\frac{1}{2}})$,
\item $\rho_{60_p}(\mathcal{A}_{E_6})\simeq SU_{60}(q^{\frac{1}{2}})$,
\item $\rho_{64_p}(\mathcal{A}_{E_6})\simeq SU_{64}(q^{\frac{1}{2}})$,
\item $\rho_{80_s}(\mathcal{A}_{E_6})\simeq SP_{80}(q^{\frac{1}{2}})$,
\item $\rho_{81_p}(\mathcal{A}_{E_6})\simeq SU_{81}(q^{\frac{1}{2}})$,
\item $\rho_{90_s}(\mathcal{A}_{E_6})\simeq SP_{90}(q^{\frac{1}{2}})$.
\end{enumerate}
\end{prop}

\begin{proof}
Assume first $\F_q=\F_p(\alpha)=\F_p(\alpha+\alpha^{-1})$. The representations considered are all defined over $\F_p(\sqrt{\alpha})$. If $\F_p(\sqrt{\alpha})\neq\F_p(\alpha)$ then $X^2-\alpha$ is an irreducible polynomial over $\F_p(\alpha)$. The unique automorphism $\sigma$ of order $2$ of $\F_p(\sqrt{\alpha})$ fixes $\F_p$ pointwise and verifies $\sigma(\sqrt{\alpha})=-\sqrt{\alpha}$. Hence, Proposition \ref{color} and Lemma \ref{Ngwenya} imply that the representations can be considered to be defined over $\F_q$.

We have found non-degenerate skew-symmetric bilinear forms defined over $\F_p$ associated to the representation $\rho_{\widetilde{10_s}},\rho_{\widetilde{20_s}},\rho_{\widetilde{60_s}},\rho_{\widetilde{80_s}}$ and $\rho_{\widetilde{90_s}}$, therefore we have all the corresponding inclusions for the images of $\mathcal{A}_{E_6}$. 

It now only remains to prove that all the inclusions are isomorphisms. We prove it separately for $6_p$ and $10_s$ and we will use Theorem \ref{CGFS} for the remaining ones. 

By Table \ref{resE6D5} and Theorem \ref{Sandburg}, we have that $\rho_{6_p}(\mathcal{A}_{D_5})\simeq SL_5(q)\times \{1\}$. It follows by Lemma \ref{normclosD5E6} that $\rho_{6_p}(\mathcal{A}_{E_6})$ is generated by transvections. Theorem \ref{transvections} then shows that there exists $q'$ such that up to conjugation in $GL_6(q)$, we have $\rho_{6_p}(\mathcal{A}_{E_6})\in \{SL_6(q'),SP_6(q'),SU_6(q'^{\frac{1}{2}})\}$. Since it contains a natural $SL_5(q)$, we have that $q'=q$ and $\rho_{6_p}(\mathcal{A}_{E_6})\simeq SL_6(q)$.

By Table \ref{resE6D5}, we have that $\rho_{10_s}(\mathcal{A}_{D_5})\simeq \rho_{[2,2],[1]}(\mathcal{A}_{D_5})$. We have $\varphi([2,2],[1])=([2,2],[1])$ and $\tilde{\nu}([2,2],[1])=\nu([2,2])\nu([1])$ since $[2,2]'>[1]'$. It follows that $\tilde{\nu}([2,2],[1])=(-1)^{\frac{4-2}{2}}(-1)^{\frac{1-1}{2}}=-1$ and by Theorem \ref{Sandburg}, $\rho_{[2,2],[1]}(\mathcal{A}_{D_5})\simeq SP_{10}(q)$. It follows that $\rho_{10_s}(\mathcal{A}_{E_6})\simeq SP_{10}(q)$. 

By Table \ref{resE6D5}, the non self-dual representations contain a natural $SL_2(q)$ and the self-dual ones contain a  twisted diagonal $SL_3(q)$. We can therefore apply Theorem \ref{CGFS} and Lemmas \ref{tens1}, \ref{tens2}, \ref{exceptiontens2} to conclude the proof.

\smallskip

Assume now $\F_q=\F_p(\alpha)\neq\F_p(\alpha+\alpha^{-1})$. The representations can again be considered to be defined over $\F_q$. We have that $X^2-(\alpha+\alpha^{-1})X+1$ is an irreducible $\F_p(\alpha+\alpha^{-1})$-polynomial over $\F_p$, therefore we have an automorphism $\epsilon$ of order $2$ of $\F_q$ mapping $\alpha$ to $\alpha^{-1}$. To conclude the proof, we only need to use previous arguments if we prove that all the representations considered verify $\rho \simeq \epsilon \circ \rho^\star$.

\smallskip

First consider the non self-dual representations of dimension different from $15,20$ and $60$. Let $\rho$ be such a representation, we know that $\epsilon\circ \rho^\star$ is an irreducible representation, therefore we have $\epsilon\circ \rho^\star \simeq \rho$ or $\epsilon \circ \rho^\star\simeq \rho^\star$ since those are the only irreducible representations of the same dimension. Assume $\epsilon \circ \rho^\star \simeq \rho^\star$. Then $\rho\simeq \epsilon \circ \rho$ and Lemma \ref{Harinordoquy} implies that up to conjugation $\rho\leq SL_{n_\rho}(q^{\frac{1}{2}})$. By Table \ref{resE6D5}, those representations all contain a natural $SU_{a}(q^{\frac{1}{2}})$ with multiplicity $1$ for some $a\geq 5$. This implies that $SU_{a}(q^{\frac{1}{2}})$ is a subgroup of $SL_{a}(q^{\frac{1}{2}})$ which is absurd by simple cardinality arguments since $a\geq 3$. It follows by contradiction that $\rho\simeq \epsilon\circ \rho^\star$. The result is obvious for the self-dual representations of dimension different from $15$, $20$ and $60$ because there is only one possibility when $\rho\simeq \rho^\star$.

\smallskip

Consider now the $60$-dimensional representations, we have $\rho_{60_s}\simeq \rho_{60_s}^\star$ and $\rho_{60_p}\not\simeq \rho_{60_p}^\star$. We know that $\epsilon \circ \rho_{60_s}^\star\simeq \rho_{60_s}$ or $\epsilon \circ \rho_{60_s}^\star\simeq \rho_{60_p}$ or $\epsilon \circ \rho_{60_s}^\star\simeq \rho_{60_p}^\star$. We have $\epsilon \circ \rho_{60_s}^\star\simeq \epsilon \circ\rho_{60_s} \simeq (\epsilon\circ \rho_{60_s}^\star)^\star$. This proves that the only possibility is the first one, therefore $\epsilon \circ \rho_{60_s}^\star \simeq \rho_{60_s}$. We have $\epsilon \circ \rho_{60_p}^\star \simeq \rho_{60_p}$ or $\epsilon \circ \rho_{60_p}^\star \simeq \rho_{60_p}^\star$ or $\epsilon \circ \rho_{60_p}^\star \simeq \rho_{60_s}$. The second possibility is excluded by the same reasonning as for the representations of dimension different from $15$ or $20$. The third possibility would imply that $\epsilon \circ \rho_{60_s}\simeq \rho_{60_p}^\star$. By the above, this would imply $\rho_{60_s}\simeq \rho_{60_p}^\star$ which is absurd. It follows that $\epsilon \circ \rho_{60_p}^\star \simeq \rho_{60_p}$. The arguments are identical for the $20$-dimensional representations. 

\smallskip

It only remains to consider the $15$-dimensional representations. We have $\rho_{15_p}\not\simeq \rho_{15_p}^\star$, $\rho_{15_q}\not\simeq \rho_{15_q}^\star$. The are therefore four $15$-dimensional representations. We have $\epsilon \circ \rho_{15_p}^\star\simeq \rho_{15_p}$ or $\epsilon \circ \rho_{15_p}^\star\simeq \rho_{15_p}^\star$ or $\epsilon \circ \rho_{15_p}^\star\simeq \rho_{15_q}$ or $\epsilon \circ \rho_{15_p}^\star\simeq \rho_{15_q}^\star$. Using the same arguments as before we have $\epsilon \circ \rho_{15_p}^\star \not\simeq \rho_{15_p}^\star$. Assume now $\epsilon \circ \rho_{15_p}^\star\simeq \rho_{15_q}$. We have that $\rho_{15_p|\mathcal{H}_{D_5,\alpha}}\simeq \rho_{[4],[1]}\oplus \rho_{[3],[1^2]}$ and $\rho_{15_q|\mathcal{H}_{D_5,\alpha}}\simeq \rho_{[3,2],\emptyset}\oplus \rho_{[3],[2]}$. We know by Proposition \ref{unitary2} that $\epsilon \circ \rho^\star\simeq \rho$ for every representation $\rho$ of $\mathcal{H}_{D_5}$. This implies that under the assumption $\epsilon \circ \rho_{15_p}^\star\simeq \rho_{15_q}$, we would have $\rho_{[4],[1]}\oplus \rho_{[3],[1^2]}\simeq \rho_{[3,2],\emptyset}\oplus \rho_{[3],[2]}$ which is absurd. We exclude in the same way $\epsilon \circ \rho_{15_p}^\star\simeq \rho_{15_q}^\star$. It follows that $\epsilon \circ \rho_{15_p}^\star \simeq \rho_{15_p}$. In the same way $\epsilon \circ \rho_{15_q}^\star \simeq \rho_{15_q}$. This concludes the proof.
\end{proof}

We now state the main theorem for type $E_6$.

\begin{theo}\label{resultE6}
Assume $\F_q=\F_p(\alpha)=\F_p(\alpha+\alpha^{-1})$, we then have that the morphism from $\mathcal{A}_{E_6}$ to $\mathcal{H}_{E_6,\alpha}^\star\simeq \underset{\rho \mbox{ irr}}\prod GL_{n_\rho}(q)$ factorizes through the surjective morphism

$$\Phi : \mathcal{A}_{E_6}\rightarrow SL_6(q)\times SP_{10}(q)\times SL_{15}(q)^2\times SL_{20}(q)\times SP_{20}(q)\times SL_{24}(q)\times SL_{30}(q)$$
$$\times SP_{60}(q)\times SL_{60}(q)\times SL_{64}(q)\times SP_{80}(q)\times SL_{81}(q)\times SP_{90}(q).$$

Assume $\F_q=\F_p(\alpha)\neq\F_p(\alpha+\alpha^{-1})$, we then have that the morphism from $\mathcal{A}_{E_6}$ to $\mathcal{H}_{E_6,\alpha}^\star\simeq \underset{\rho \mbox{ irr}}\prod GL_{n_\rho}(q)$ factorizes through the surjective morphism

$$\Phi : \mathcal{A}_{E_6}\rightarrow SU_6(q^{\frac{1}{2}})\times SP_{10}(q^{\frac{1}{2}})\times SU_{15}(q^{\frac{1}{2}})^2\times SU_{20}(q^{\frac{1}{2}})\times SP_{20}(q^{\frac{1}{2}})\times SU_{24}(q^{\frac{1}{2}})\times SU_{30}(q^{\frac{1}{2}})$$
$$\times SP_{60}(q^{\frac{1}{2}})\times SU_{60}(q^{\frac{1}{2}})\times SU_{64}(q^{\frac{1}{2}})\times SP_{80}(q^{\frac{1}{2}})\times SU_{81}(q^{\frac{1}{2}})\times SP_{90}(q^{\frac{1}{2}}).$$

\end{theo}

\begin{proof}

By \cite{MR}, $\mathcal{A}_{E_6}$ is perfect. Furthermore, Lemma \ref{Goursat} gives that the morphism is surjective unless there exists two different representations $\rho_1$ and $\rho_2$ in the decomposition such that $\Psi\circ \rho_{1|\mathcal{A}_{H_4}}\simeq \rho_{2|\mathcal{A}_{H_4}}$ for some field automorphism $\Psi$. By Proposition \ref{Fieldfactorization}, we have that $\Psi(\alpha+\alpha^{-1})=\alpha+\alpha^{-1}$. This shows that $\Psi$ must be trivial over $\F_p(\alpha+\alpha^{-1})$. It follows by the previous propositions that there are no such representations in the decompositions and the proof is concluded.
\end{proof}

\section{Type $E_7$}\label{E7section}
Let $p$ be a prime different from $2$ and $3$ and $\alpha\in \overline{\F_p}$ of order not dividing $8,10,12,14$ and $18$. We write $\F_q=\F_p(\alpha)$. There are $60$ irreducible representations of $\mathcal{H}_{E_7,\alpha}$, none of them are self-dual. The highest dimensional representation is of dimension $512$. They are all $2$-colorable except for the two $512$-dimensional representations.

\begin{Def}
The Iwahori-Hecke algebra $\mathcal{H}_{E_7,\alpha}$ of type $E_7$ is the $\F_q$-algebra generated by $S_1,S_2,S_3,S_4,S_5,S_6,S_7$ and the following relations 
\begin{enumerate}
\item $\forall i\in \{1,2,3,4,5,6,7\}$, $(S_i-\alpha)(S_i+1)=0$.
\item $S_1S_3S_1=S_3S_1S_3$.
\item $\forall i \in \{2,4,5,6,7\}$, $S_1S_i=S_iS_1$.
\item $S_2S_4S_2=S_4S_2S_4$.
\item $\forall i\in \{3,5,6,7\}$, $S_2S_i=S_iS_2$.
\item $\forall i\in \{3,4,5,6\}$, $S_{i}S_{i+1}S_{i}=S_{i+1}S_iS_{i+1}$.
\item $S_3S_5=S_5S_3$.
\item $S_3S_6=S_6S_3$.
\item $S_3S_7=S_7S_3$.
\item $S_4S_6=S_6S_4$.
\item $S_4S_7=S_7S_4$.
For $\sigma$ in the Coxeter group $E_7$, we set $T_{\sigma}=S_{i_1}\dots S_{i_k}$ for any reduced expression $\sigma=s_{i_1}\dots s_{i_k}$ is a reduced expression.
\end{enumerate}
\end{Def}

This means we consider $E_7$ as in the CHEVIE package of GAP3 \cite{CHEVIE} with the following Dynkin diagram

\begin{center}
\begin{tikzpicture}
[place/.style={circle,draw=black,
inner sep=1pt,minimum size=10mm}]
\node (1) at (0,0)[place]{$S_1$};
\node (2) at (2,0)[place]{$S_3$};
\node (3) at (4,0)[place]{$S_4$};
\node (4) at (4,2)[place]{$S_2$};
\node (5) at (6,0)[place]{$S_5$};
\node (6) at (8,0)[place]{$S_6$};
\node (7) at (10,0)[place]{$S_7$};
\draw (1) to (2);
\draw (2) to (3);
\draw (3) to (4);
\draw (3) to (5);
\draw (5) to (6);
\draw (6) to (7);
\end{tikzpicture}
\end{center}

\begin{prop}
Under our assumptions on $p$ and $\alpha$, $\mathcal{H}_{E_7,\alpha}$ is split semisimple, the representations afforded by the $W$-graphs are irreducible and pairwise non-isomorphic over $\F_q$. The restrictions of the irreducible representations of $\mathcal{H}_{E_7,\alpha}$ to $\mathcal{H}_{E_6,\alpha}$ are the same as in the generic case.
\end{prop}

\begin{proof}
We will apply Proposition \ref{Tits}. Let $A=\Z[\sqrt{u}^{\pm 1}]$ and $F=\Q(\sqrt{u})$. We have a symetrizing trace defined by $\tau(T_0)=1$ and $\tau(T_{\sigma})=0$ for all $\sigma\in E_7\setminus \{1_{E_7}\}$. By \cite{G-P}, $\mathcal{H}_{E_7,u}$ is then a free symmetric $F$-algebra of rank $2903040$. By \cite{Bourb} V.3. Corollary 1, $A$ is integrally closed. Let $\theta$ be the ring homomorphism from $A$ to $L=\F_q$ defined by $\theta(u)=\alpha$ and $\theta(k)=\overline{k}$. We know $FH$ is split. The basis formed by the elements $T_\sigma$, $\sigma\in E_7$ verifies the conditions of the Proposition \ref{Tits}. The $E_7$-graphs remains connected since all the weights lie in $\{-3,-2,-1,1,2,3\}$. 

It now only remains to check that the Schur elements associated to the specialized representations are in $B$ and do not vanish under $\theta$ with $B$ as in Proposition \ref{Tits}. The Schur elements are given in Table \ref{SchurelmtsE7}. They were obtained using Proposition $9.3.6$ and Table $E.6$ of \cite{G-P}. For a pair $(\rho,\rho^\star)$ of representations, we only give the Schur element of one of the representations since the other is obtained by applying the involution $\sqrt{u}\mapsto \sqrt{u}^{-1}$. The conditions on $\alpha$ and $p$ imply that the Schur elements verify the right conditions and the proof is concluded.
\end{proof}

\begin{table}
\centering
$\begin{array}{lcr}
 1_a & : & (\Phi_2^7\Phi_3^3\Phi_4^2\Phi_5\Phi_6^3\Phi_7\Phi_8\Phi_9\Phi_{10}\Phi_{12}\Phi_{14}\Phi_{18})(u).\\
 7_a' & : & \frac{1}{u}((\Phi_2^7\Phi_3^3\Phi_4^2\Phi_5\Phi_6^3\Phi_8\Phi_9\Phi_{10}\Phi_{18})(u).\\
 15_a' & : & \frac{2}{u^4}(\Phi_2^7\Phi_3^3\Phi_4^2\Phi_6^3\Phi_7)(u).\\
 21_a & : & \frac{2}{u^3}(\Phi_2^7\Phi_3^3\Phi_4^2\Phi_5\Phi_6^3\Phi_{18})(u).\\
 21_b' & : &  \frac{1}{u^3}(\Phi_2^7\Phi_3^3\Phi_4^2\Phi_5\Phi_6^3\Phi_8\Phi_{10}\Phi_{12})(u).\\
 27_a & : & \frac{1}{u^2}(\Phi_2^7\Phi_3\Phi_4^2\Phi_5\Phi_6\Phi_7\Phi_8\Phi_{10}\Phi_{14})(u).\\
 35_a' & : & \frac{6}{u^7}(\Phi_2^7\Phi_3^3\Phi_4^2\Phi_{18})(u).\\
 35_b & : & \frac{2}{u^3}(\Phi_2^7\Phi_3^3\Phi_4^2\Phi_6^3\Phi_9\Phi_{10})(u).\\
 56_a' & : & \frac{2}{u^3}(\Phi_2^3\Phi_3^3\Phi_4^2\Phi_5\Phi_6\Phi_8\Phi_9\Phi_{12})(u).\\
 70_a' & : & \frac{3}{u^7}(\Phi_2^7\Phi_3^3\Phi_4^2\Phi_6^3)(u).\\
 84_a & : & \frac{2}{u^{10}}(\Phi_2^7\Phi_3^3\Phi_5\Phi_6^3)(u).\\
 105_a' & : & \frac{2}{u^4}(\Phi_2^7\Phi_3^3\Phi_4^2\Phi_6^3\Phi_{14})(u).\\
 105_b & : & \frac{1}{u^6}(\Phi_2^7\Phi_3^3\Phi_4^2\Phi_6^3\Phi_8)(u).\\
 105_c & : & \frac{1}{u^{12}}(\Phi_2^7\Phi_3^3\Phi_4^2\Phi_6^3\Phi_8)(u).\\
 120_a & : & \frac{2}{u^4}(\Phi_2^3\Phi_3^3\Phi_4^2\Phi_6\Phi_7\Phi_8\Phi_{12})(u).\\
 168_a & : & \frac{1}{u^6}(\Phi_2^7\Phi_3^3\Phi_5\Phi_6^3\Phi_{10})(u).\\
 189_a & : & \frac{2}{u^8}(\Phi_2^7\Phi_3\Phi_4^2\Phi_5\Phi_{14})(u).\\
 189_b' & : & \frac{1}{u^5}(\Phi_2^7\Phi_3\Phi_4^2\Phi_5\Phi_6\Phi_8\Phi_{10})(u).\\ 
 189_c' & : & \frac{1}{u^7}(\Phi_2^7\Phi_3\Phi_4^2\Phi_5\Phi_6\Phi_8\Phi_{10})(u).\\
 210_a & : & \frac{1}{u^6}(\Phi_2^7\Phi_3^3\Phi_4^2\Phi_6^3\Phi_{12})(u).\\
 210_b & : & \frac{1}{u^{10}}(\Phi_2^7\Phi_3^3\Phi_4^2\Phi_6^3)(u).\\
 216_a' & : & \frac{2}{u^8}(\Phi_2^3\Phi_3\Phi_4^2\Phi_5\Phi_7\Phi_8)(u).\\
 280_a' & : & \frac{3}{u^7}(\Phi_2^7\Phi_3^3\Phi_6^3\Phi_{12})(u).\\
 280_b & : & \frac{2}{u^7}(\Phi_2^3\Phi_3^3\Phi_4^2\Phi_8\Phi_9)(u).\\
 315_a' & : & \frac{6}{u^7}(\Phi_2^7\Phi_4^2\Phi_6^3\Phi_9)(u).\\
 336_a' & : & \frac{2}{u^{10}}(\Phi_2^3\Phi_3^3\Phi_4^2\Phi_5\Phi_6\Phi_{12})(u).\\
 378_a' & : & \frac{1}{u^9}(\Phi_2^7\Phi_3\Phi_4^2\Phi_5\Phi_6\Phi_{10})(u).\\
 405_a & : & \frac{2}{u^8}(\Phi_2^7\Phi_4^2\Phi_6\Phi_7\Phi_{10})(u).\\
 420_a & : & \frac{2}{u^{10}}(\Phi_2^7\Phi_3^3\Phi_6^3\Phi_{10})(u).\\
 512_a' & : & \frac{2}{u^{11}}(\Phi_3^3\Phi_5\Phi_7\Phi_9)(u).
\end{array}$
\caption{Schur elements in type $E_7$}
\label{SchurelmtsE7}
\end{table}

The restriction table from $\mathcal{H}_{E_7,\alpha}$ to its subalgebra $\mathcal{H}_{E_6,\alpha}$ generated by $S_1,S_2,S_3,S_4,S_5$ and $S_6$ which is naturally isomorphic to the Iwahori-Hecke of type $E_6$ with parameter $\alpha$ is then given by Table \ref{resE7E6}. It is obtained using the CHEVIE package of GAP3 \cite{CHEVIE}. They correspond in the generic case to the induction/restriction tables of the corresponding Coxeter groups.

\begin{table}
\noindent
\tabcolsep=0.1cm
\begin{tabular}{ |p{0.675cm}||p{0.33cm}|p{0.33cm}|p{0.33cm}|p{0.33cm}|p{0.51cm}|p{0.51cm}|p{0.51cm}|p{0.51cm}|p{0.51cm}|p{0.51cm}| p{0.51cm}|p{0.51cm}|p{0.51cm}|p{0.51cm}|p{0.51cm}|p{0.51cm}|p{0.51cm}|p{0.51cm}|p{0.51cm}|p{0.51cm}|p{0.51cm}|p{0.51cm}|p{0.51cm}|p{0.51cm}|p{0.51cm}|  }
 \hline
 & $1_p$ & $1_p'$ & $6_p$ & $6_p'$ & $10_s$ & $15_p$ & $15_p'$ & $15_q$  & $15_q'$ & $20_p$ & $20_p'$ & $20_s$ & $24_p$ & $24_p'$ & $30_p$ & $30_p'$ & $60_s$ & $60_p$ & $60_p'$ & $64_p$ & $64_p'$ & $80_s$ & $81_p$ & $81_p'$ & $90_s$  \\
 \hline

 $1_a$   & $ 1$   &  &    &  &  &  &  &  &  & &  &    &  &  &  &  &  &  &  &  &  &  &  & & \\
 \hline
 $7_a'$   &  $1$   &  &  $1$  &  &  &  &  &  &  & &  &    &  &  &  &  &  &  &  &  &  &  &  & &\\
 \hline
 $15_a'$   &     &  &    &  &  &  &  & $1$ &  & &  &    &  &  &  &  &  &   &  &  &  &  &  & &\\
 \hline
 $21_a$   &     &  &   $1$ &  &  & $1$ &  &  &  & &  &    &  &  &  &  &  &   &  &  &  &  &  & &\\
 \hline
 $21_b'$   & $1$    &  &    &  &  &  &  &  &  & $1$&  &    &  &  &  &  &  &  &  &  &  &  &  & & \\
 \hline
 $27_a$   & $1$    &  & $1$   &  &  &  &  &  &  & $1$&  &    &  &  &  &  &  &  &  &  &  &  &  & & \\
 \hline
 $35_a'$   &     &  &    &  &  & $1$ &  &  &  &  &   &  $1$  &  &  &  &  &  &   &  &  &  &  &  & & \\
\hline 
 $35_b$   &     &  &    &  &  &  &  & $1$ &  & $1$&  &    &  &  &  &  &  &  &  &  &  &  &  & & \\
 \hline
 $56_a'$   &     &  & $1$   &  &  &  &  &  &  & $1$&  &    &  &  & $1$ &  &  &  &  &  &  &  &  & & \\
 \hline
 $70_a'$   &     &  &    &  & $1$ &  &  &  &  &  &   &    &  &  &  &  &  &  $1$ &  &  &  &  &  & & \\
\hline 
 $84_a$   &     &  &    &  &  &  &  &  &  &  &   &    & $1$ &  &  &  & $1$ &   &  &  &  &  &  & & \\
\hline
 $105_a'$   &     &  & $1$   &  &  & $1$ &  &  &  & $1$& &    &  &  &  &  &  &  &  & $1$ &  &  &  &  &\\
 \hline
 $105_b$   &     &  &    &  &  &  &  & $1$ &  & &  &    &  &  &  $1$&  &  & $1$  &  &  &  &  &  & & \\
 \hline
 $105_c$   &     &  &    &  &  &  &  &  &  &  &   &    &  &  $1$&  &  &  &   &  &  &  &  &  $1$& & \\
 \hline
  $120_a$   &     &  &  $1$  &  &  &  &  &  &  &$1$ &  &    &  &  & $1$ &  &  &  &  & $1$ &  &  &  & & \\
 \hline
 $168_a$   &     &  &    &  &  &  &  &  &  & $1$ &  &    &  $1$&  &  &  &  &  $1$&  &  $1$&  &  &  &  &\\
 \hline
 $189_a$   &     &  &    &  &  & $1$ &  &  &  &  &   & $1$   &  &  &  &  &  &   &  & $1$ &  &  &  & & $1$\\
\hline 
 $189_b'$   &     &  &    &  &   &  &  & $1$ &  &$1$ &  &    &  &  & $1$ &  &  & $1$ &  & $1$ &  &  &  & & \\
 \hline
 $189_c'$   &     &  &    &  &  &  &  &  &  & $1$ &   &    & $1$&  &  &  &  &   &  & $1$ &  &  & $1$ & & \\
\hline 
 $210_a$   &     &  &    &  &  &  $1$&  &  &  & $1$ &  &    &  &  &  $1$&  &  &  &  & $1$  &  &  & $1$ & & \\
 \hline
 $210_b$   &     &  &    &  & $1$  &  &  &  &  &  &   &    &  &  &  &  & $1$ &  $1$ &  &  &  &  $1$&  & & \\
\hline 
 $216_a'$   &     &  &    &  &  &  &  & $1$  &  &  &   &    &  &  &  &  & $1$ &  $1$ &  &  &  &  & $1$ & & \\
\hline 
 $280_a'$   &     &  &    &  &  & $1$ &  &  &  &  &   &    &  &  & $1$ &  &  &   &  & $1$ &  &  & $1$ & & $1$\\
\hline 
 $280_b$   &     &  &    &  &  &  &  & $1$ &  &  &   &    &  &  &  &  & $1$ & $1$  &  & $1$ &  &  & $1$ & & \\
\hline 
 $315_a'$   &     &  &    &  &  &  &  &  &  &  &   &    &  &  & $1$ &  &  & $1$  &  & $1$ &  & $1$ & $1$ & & \\
\hline 
 $336_a'$   &     &  &    &  &  &  &  &  &  &  &   & $1$ &  &  &  &  &  &   &  & $1$ &  &  & $1$ & $1$ &$1$ \\
 \hline
 $378_a'$   &     &  &    &  &  &  &  &  &  &  &   &    & $1$ &  &  &  & $1$ &  $1$ &  & $1$ &  & $1$ &  & & $1$\\
\hline 
 $405_a$   &     &  &    &  &  &  &  &  &  &  &   &    &  &  & $1$ &  &  & $1$  &  & $1$ &  & $1$ & $1$ & & $1$\\
\hline 
 $420_a$   &     &  &    &  &  &  &  &  &  &  &   &    & $1$ &  &  &  &  &   &  & $1$ &  &$1$  & $1$ &$1$ & $1$\\
\hline
 $512_a'$   &     &  &    &  &  &  &  &  &  &  &   &    &  &  &  &  & $1$ &  $1$ & $1$ &  &  & $1$ &$1$  & $1$& $1$\\
\hline
 
 \end{tabular}
 \smallskip
 
 \caption{Restriction table from $\mathcal{H}_{E_7,\alpha}$ to $\mathcal{H}_{E_6,\alpha}$.}
\label{resE7E6}
\end{table}

\begin{prop}\label{resE7derivedsubgroup}
The restrictions to $\mathcal{A}_{E_7}$ of the representations afforded by those $E_7$-graphs are absolutely irreducible and the representations of dimension greater than $1$ are pairwise non-isomorphic.
\end{prop}

\begin{proof}
As in \cite{BMM} Lemma $3.4$, we only need to prove that $A_{E_7}$ is generated by $A_{E_6}$ and $\mathcal{A}_{E_7}$. This true because $s_7=s_7s_1^{-1}s_1$, $s_7s_1^{-1}\in \mathcal{A}_{E_6}$ and $s_1\in A_{E_6}$.

We now prove the second part of the statement. Let $\rho_1$ and $\rho_2$ be two irreducible representations of $\mathcal{H}_{E_7,\alpha}$ such that $\rho_{1|\mathcal{A}_{E_7}}\simeq \rho_{1|\mathcal{A}_{E_7}}$. By Lemma \ref{abel}, there exists a character $\xi:A_{E_7}\mapsto \F_q^{\star}$ such that $\rho_1\simeq \rho_2\otimes \xi$. This means there exists $x\in \F_q$ such that for all $i\in [\![1,7]\!]$, $\rho_1(S_i)$ is conjugate to $x\rho_2(S_i)$. We know for any representation $\rho$ of dimension greater than $1$, the set of eigenvalues of $\rho$ is equal to $\{\alpha,-1\}$. This implies that $\{\alpha,-1\}=\{x\alpha,-x\}$. We then have $x=1$ or ($x=-\alpha$ and $-\alpha^2=-1$). It follows that $x=1$ and $\rho_1\simeq \rho_2$.
\end{proof}

We now prove the usual lemma computing the normal closure of $\mathcal{A}_{E_6}$ inside $\mathcal{A}_{E_7}$.

\begin{lemme}\label{normclosE6E7}
The normal closure $\ll \mathcal{A}_{E_6} \gg_{\mathcal{A}_{E_7}}$ of $\mathcal{A}_{E_6}$ inside $\mathcal{A}_{E_7}$ is equal to $\mathcal{A}_{E_7}$, where we identify $A_{E_6}$ as a subgroup of $A_{E_7}$ using the natural isomorphism from $A_{E_6}$ to $<S_i,i \in [\![1,6]\!]>$.
\end{lemme}

\begin{proof}
By \cite{MR}, we have $\mathcal{A}_{E_7}=<S_1S_3S_1^{-1},S_4S_3^{-1},S_iS_1^{-1}, i\in [\![1,7]\!]>$ and 
$$\mathcal{A}_{E_6}=<S_1S_3S_1^{-1},S_4S_3^{-1},S_iS_1^{-1}, i\in [\![1,6]\!]>.$$
 This proves we only need to show that $S_7S_1^{-1}\in \ll \mathcal{A}_{E_6} \gg_{\mathcal{A}_{E_7}}$.
We have $S_7S_6S_7=S_6S_7S_6$, therefore $S_7=S_6S_7S_6(S_6S_7)^{-1}$ and 
$$S_7S_1^{-1}=S_6S_7S_6S_1^{-1}(S_6S_7)^{-1}=((S_6S_1^{-1})(S_7S_1^{-1}))(S_6S_1^{-1})((S_6S_1^{-1})(S_7S_1^{-1}))^{-1}.$$
This concludes the proof.
\end{proof}

Note now that there are no self-dual representations in type $E_7$, they are all $2$-colorable except for the two representations of dimension $512$. We then have the following proposition.

\begin{prop}\label{ImageinsiderepsE7}
If $\F_q=\F_p(\sqrt{\alpha})=\F_p(\alpha)=\F_p(\alpha+\alpha^{-1})$ then for any irreducible representation $\rho$ of $\mathcal{H}_{E_7,\alpha}$, $\rho(\mathcal{A}_{E_7})\simeq SL_{n_\rho}(q)$, where $n_\rho=\dim(\rho)$.

If $\F_p(\sqrt{\alpha})\neq \F_q=\F_p(\alpha)=\F_p(\alpha+\alpha^{-1})$ then for any irreducible representation  $\rho$ of $\mathcal{H}_{E_7,\alpha}$ such that $n_\rho\neq 512$, we have $\rho(\mathcal{A}_{E_7})\simeq SL_{n_\rho}(q)$. We have $\rho_{512_a}(\mathcal{A}_{E_7})\simeq SU_{512}(q)$.

If $\F_q=\F_p(\alpha)\neq \F_p(\alpha+\alpha^{-1})$ then for any irreducible representation $\rho$ of $\mathcal{H}_{E_7,\alpha}$, we have $\rho(\mathcal{A}_{E_7})\simeq SU_{n_\rho}(q^{\frac{1}{2}})$.
\end{prop}

\begin{proof}
Assume first that $\F_q=\F_p(\alpha)=\F_p(\alpha+\alpha^{-1})$. Let $\rho$ be an irreducible representation of $\mathcal{H}_{E_7,\alpha}$ of dimension different from $7$ and $512$. The associated $E_7$-graph is then $2$-colorable, therefore by Proposition \ref{color}, the image of $A_{E_7}$ under $\rho$ is inluded up to conjugation in $GL_{n_\rho}(q)$ even when $\F_p(\sqrt{\alpha})\neq \F_p(\alpha)$. By Table \ref{resE6D5} and Proposition \ref{ImageinsiderepsE6}, $\rho(\mathcal{A}_{E_7})$ contains a natural $SL_r(q)$ for some $r\geq 6$ whenever $n_\rho>1$.
 We can therefore apply Theorem \ref{CGFS}. We get that $\rho(\mathcal{A}_{E_7})$ is a classical group over $\F_{q'}$ in a natural representation for some $q'$ dividing $q$. Since it contains a natural $SL_r(q)$, we get $q'=q$ and the representation is not unitary.
  By Corollary \ref{resE7derivedsubgroup}, it cannot preserve a bilinear form, therefore we get $\rho(\mathcal{A}_{E_7})\simeq SL_{n_\rho}(q)$. We now have to consider the $7$-dimensional representations $\rho_{7_a'}$ and $\rho_{7_a}=\rho_{7_a'}^{\star}$. By Lemma \ref{normclosE6E7} and Proposition \ref{ImageinsiderepsE6}, $\rho_{7_a}(\mathcal{A}_{E_7})$ is normally generated by a natural $SL_6(q)$. It follows by Corollary \ref{resE7derivedsubgroup} that it is an irreducible subgroup of $SL_7(q)$ generated by transvections. We can then apply Theorem \ref{transvections} and the same arguments as above give $\rho_{7_a'}(q)\simeq SL_{7}(q)$. This concludes the proof for the representations of dimension different from $512$.

  \smallskip
  
  Consider now the representation $\rho_{512_a'}$. If $\F_p(\sqrt{\alpha})=\F_p(\alpha)$ then $\rho(\mathcal{A}_{E_7})$ is included in $SL_{512}(q)$. By Proposition \ref{ImageinsiderepsE6}, it contains a natural $SP_{60}(q)$, therefore we can apply the above reasoning to get that $\rho_{512_a'}(\mathcal{A}_{E_7})\simeq SL_{512}(q)$. Assume now $\F_p(\sqrt{\alpha})\neq \F_p(\alpha)$. We then have that $X^2-\alpha$ is an irreducible polynomial of degree $2$ of $\F_p(\alpha)$ such that $\F_p(\sqrt{\alpha})=\F_p(\alpha)/(X^2-\alpha)$. It follows that $\F_p(\sqrt{\alpha})$ is an extension of degree $2$ of $\F_p(\alpha)$ and there exists a unique automorphism $\Psi$ of $\F_p(\sqrt{\alpha})$ which fixes $\F_q$ pointwise and such that $\Psi(\sqrt{\alpha})=-\sqrt{\alpha}$. Note that the restriction of $\rho_{512_a'}$ and $\rho_{512_a}$ to $\mathcal{A}_{E_6}$ are identical, therefore we cannot get any information from this restriction. Using CHEVIE \cite{CHEVIE}, we have
   $$\tr(\rho_{512_a'}(S_1S_7^{-1}S_6S_2^{-1}S_3S_5S_4^{-2}) )=228+5\alpha^{-4}-34\alpha^{-3}+104\alpha^{-2}-189\alpha^{-1}$$
   $$+\sqrt{\alpha}^{-1}-\sqrt{\alpha}-189\alpha+104\alpha^2-34\alpha^{3}+5\alpha^4.$$
   
  It follows that $\rho_{512_a'}\simeq \Psi\circ \rho_{512_a'}$ then $\sqrt{\alpha}-\sqrt{\alpha}^{-1}=\Psi(\sqrt{\alpha}^{-1}-\sqrt{\alpha})=\sqrt{\alpha}^{-1}-\sqrt{\alpha}$. Hence $2\sqrt{\alpha}=2\sqrt{\alpha}^{-1}$ and $\alpha=1$ which contradicts our assumptions on $\alpha$. This proves that $\Psi\circ\rho_{512_a'}\simeq \rho_{512_a}=\rho_{512_a'}^\star$. It follows that $\rho_{512_a'}(\mathcal{A}_{E_7})$ is included in $SU_{512}(q)$ up to conjugation in $GL_{512}(q^2)$. Since $\rho_{512_a}$ contains a natural $SP_{60}(q)$, we can apply Theorem \ref{CGFS} and $\rho_{512_a}$ is a classical group over $\F_{q'}$ for some $q'$ dividing $q^2$. However $q$ divides $q'$ because $\rho_{512_a'}(\mathcal{A}_{E_7})$ contains a natural $SP_{60}(q)$, therefore $q'\in \{q,q^2\}$. We know $\rho_{512_a}$ does not preserve any non-degenerate bilinear form because $\rho_{512_a|\mathcal{A}_{E_7}}\not\simeq \rho_{512_a'|\mathcal{A}_{E_7}}$. It follows that up conjugation in $GL_{512}(q^2)$, we have $\rho_{512_a}(\mathcal{A}_{E_7})\in \{SL_{512}(q),SU_{512}(q^{\frac{1}{2}}),SL_{512}(q^2),SU_{512}(q)\}$. Assume we are in one the first two cases, we would then have using the trace of the same element as above that $\sqrt{\alpha}^{-1}-\sqrt{\alpha}\in \F_q$. This would imply that $\frac{1-\alpha}{\sqrt{\alpha}}\in \F_q$, therefore $\sqrt{\alpha}\in \F_q$ since $\alpha\neq 1$ and $\alpha\in \F_q$. The third possibility is excluded since $SL_{512}(q^2)$ cannot be injected inside $SU_{512}(q)$. We can then conclude that $\rho_{512_a}(\mathcal{A}_{E_7})\simeq SU_{512}(q)$.
  
  \medskip
  
  Assume now that $\F_q=\F_p(\alpha)\neq \F_p(\alpha+\alpha^{-1})$. Then $X^2-(\alpha+\alpha^{-1})X+1$ is a $\F_p(\alpha+\alpha^{-1})$-irreducible polynomial of degree $2$, and$\F_q=\F_p(\alpha+\alpha^{-1})/(X^2 -(\alpha+\alpha^{-1})X+1)$ is an extension of degree $2$. There is a unique automorphism $\epsilon$ of degree $2$ of $\F_q$. It fixes $\F_p(\alpha+\alpha^{-1})$ pointwise and $\epsilon(\alpha)=\alpha^{-1}$. We can then consider the extension $\F_p(\sqrt{\alpha})$ of $\F_p(\alpha)$.  We have $\epsilon(\sqrt{\alpha})^2=\epsilon(\alpha)=\epsilon(\alpha)=\alpha^{-1}$, therefore $\epsilon(\sqrt{\alpha})\in \{-\sqrt{\alpha}^{-1},\sqrt{\alpha}^{-1}\}$. It follows that $\epsilon^2(\sqrt{\alpha})=\sqrt{\alpha}$ which implies that $\sqrt{\alpha}\in  \F_q$, therefore $\F_p(\alpha)=\F_p(\sqrt{\alpha})$.
  
  \smallskip
  
    Let $\rho$ be an irreducible representation of $\mathcal{H}_{E_7,\alpha}$ of degree greater than $1$. We have that $\epsilon \circ \rho^\star$ is an irreducible representation. Assume that the only representations of dimension $n_\rho$ are $\rho$ and $\rho^{\star}$. We then have $\epsilon \circ \rho^\star \simeq \rho$ or $\epsilon \circ \rho^\star \simeq \rho^\star$. The set of eigenvalues of $\rho_{S_1}$ is $\{\alpha,-1\}$. In the second case we would have $\epsilon(^t\!\rho(S_i)^{-1})$ conjugate to $^t\!\rho(S_i)^{-1}$, therefore $\{\epsilon(-1,\epsilon(\alpha^{-1})\}=\{-1,\alpha\}=\{-1,\alpha^{-1}\}$, therefore $\alpha^2=1$ which contradicts our assumptions. It follows that $\rho\simeq \epsilon \circ\rho^{\star}$, and Lemma \ref{Ngwenya} gives $\rho(\mathcal{A}_{E_7})\subset SU_{n_{\rho}}(q^{\frac{1}{2}})$. We can then apply the same reasoning as bove to conclude that $\rho( \mathcal{A}_{E_7})\simeq SU_{n_\rho}(q^{\frac{1}{2}})$.
    
    \smallskip
    
    It only remains to consider the representations of dimension $21$, $35$, $105$, $189$, $210$ and $280$. By Proposition \ref{ImageinsiderepsE6}, we only need to check that the restrictions to $\mathcal{H}_{E_6,\alpha}$ are different for the other representations of the same dimension. This is true by Table \ref{resE7E6}.
\end{proof}

\begin{theo}\label{resultE7}
Assume $F_q=\F_p(\sqrt{\alpha})=\F_p(\alpha+\alpha^{-1})$. Then the morphism from $\mathcal{A}_{E_7}$ to $\mathcal{H}_{E_7,\alpha}^\star\simeq \underset{\rho \mbox{ irr}}\prod GL_{n_\rho}(q)$ factorizes through the surjective morphism

$$\Phi : \mathcal{A}_{E_7}\rightarrow SL_7(q)\times SL_{15}(q)\times SL_{21}(q)^2\times SL_{27}(q)\times SL_{35}(q)^2\times SL_{56}(q)\times SL_{70}(q)\times SL_{84}(q)$$
$$\times SL_{105}(q)^3
\times SL_{120}(q)\times SL_{168}(q)\times SL_{189}(q)^3\times SL_{210}(q)^2\times SL_{216}(q)\times SL_{280}(q)^2$$
$$\times SL_{315}(q)\times SL_{336}(q)\times SL_{378}(q)\times SL_{405}(q) \times SL_{420}(q)\times SL_{512}(q).$$

Assume $\F_p(\sqrt{\alpha})\neq \F_q=\F_p(\alpha)=\F_p(\alpha+\alpha^{-1})$. Then the morphism from $\mathcal{A}_{E_7}$ to $\mathcal{H}_{E_7,\alpha}^\star\simeq \underset{\rho \mbox{ irr}}\prod GL_{n_\rho}(q)$ factorizes through the surjective morphism

$$\Phi : \mathcal{A}_{E_7}\rightarrow SL_7(q)\times SL_{15}(q)\times SL_{21}(q)^2\times SL_{27}(q)\times SL_{35}(q)^2\times SL_{56}(q)\times SL_{70}(q)\times SL_{84}(q)$$
$$\times SL_{105}(q)^3
\times SL_{120}(q)\times SL_{168}(q)\times SL_{189}(q)^3\times SL_{210}(q)^2\times SL_{216}(q)\times SL_{280}(q)^2$$
$$\times SL_{315}(q)\times SL_{336}(q)\times SL_{378}(q)\times SL_{405}(q) \times SL_{420}(q)\times SU_{512}(q).$$

Assume $F_q=\F_p(\alpha)\neq\F_p(\alpha+\alpha^{-1})$. Then the morphism from $\mathcal{A}_{E_7}$ to \\$\mathcal{H}_{E_7,\alpha}^\star\simeq \underset{\rho \mbox{ irr}}\prod GL_{n_\rho}(q)$ factorizes through the surjective morphism

\begin{small}
$$\Phi : \mathcal{A}_{E_7}\rightarrow SU_7(q^{\frac{1}{2}})\times SU_{15}(q^{\frac{1}{2}})\times SU_{21}(q^{\frac{1}{2}})^2\times SU_{27}(q^{\frac{1}{2}})\times SU_{35}(q^{\frac{1}{2}})^2\times SU_{56}(q^{\frac{1}{2}})\times SU_{70}(q^{\frac{1}{2}})\times SU_{84}(q^{\frac{1}{2}})$$
$$\times SU_{105}(q^{\frac{1}{2}})^3
\times SU_{120}(q^{\frac{1}{2}})\times SU_{168}(q^{\frac{1}{2}})\times SU_{189}(q^{\frac{1}{2}})^3\times SU_{210}(q^{\frac{1}{2}})^2\times SU_{216}(q^{\frac{1}{2}})\times SU_{280}(q^{\frac{1}{2}})^2$$
$$\times SU_{315}(q^{\frac{1}{2}})\times SU_{336}(q^{\frac{1}{2}})\times SU_{378}(q^{\frac{1}{2}})\times SU_{405}(q^{\frac{1}{2}}) \times SU_{420}(q^{\frac{1}{2}})\times SU_{512}(q^{\frac{1}{2}}).$$
\end{small}
\end{theo}

\begin{proof}

The proof of this theorem is very similar to the proof of Theorem \ref{resultE6}.

By \cite{MR}, $\mathcal{A}_{E_7}$ is perfect. We have by Goursat's Lemma \ref{Goursat} that the morphism is surjective unless there exists two different representations $\rho_1$ and $\rho_2$ in the decomposition such that there exists a field automorphism $\Psi$ verifying $\Psi\circ \rho_{1|\mathcal{A}_{E_6}}\simeq \rho_{2|\mathcal{A}_{E_6}}$. By Proposition \ref{Fieldfactorization}, we have that $\Psi(\alpha+\alpha^{-1})=\alpha+\alpha^{-1}$. This shows that $\Psi$ must be trivial over $\F_p(\alpha+\alpha^{-1})$. It follows by the previous propositions that there are no such representations in the decompositions and the proof is concluded.
\end{proof}

\section{Type $E_8$}\label{E8section}
Let $p$ be a prime different from $2$, $3$ and $5$ and $\alpha\in \overline{\F_p}$ of order not dividing $14$, $18$, $20$, $24$ and $30$. We write $\F_q=\F_p(\alpha)$. There are $112$ irreducible representations of $\mathcal{H}_{E_8,\alpha}$, $18$ of them are self-dual. For each self-dual representation, we have found a self dual-W-graph \cite{newgraphsEsterle}. Using the $2$-coloring of those graphs, we have that all of the associated bilinear forms are symmetric. The highest dimensional representation is of dimension $7168$.

\begin{Def}
The Iwahori-Hecke algebra $\mathcal{H}_{E_8,\alpha}$ of type $E_8$ is the $\F_q$-algebra generated by $S_1,S_2,S_3,S_4,S_5,S_6,S_7,S_8$ and the following relations 
\begin{enumerate}
\item $\forall i\in \{1,2,3,4,5,6,7,8\}$, $(S_i-\alpha)(S_i+1)=0$.
\item $S_1S_3S_1=S_3S_1S_3$.
\item $\forall i \in \{2,4,5,6,7,8\}$, $S_1S_i=S_iS_1$.
\item $S_2S_4S_2=S_4S_2S_4$.
\item $\forall i\in \{3,5,6,7,8\}$, $S_2S_i=S_iS_2$.
\item $\forall i\in \{3,4,5,6,7\}$, $S_{i}S_{i+1}S_{i}=S_{i+1}S_iS_{i+1}$.
\item $\forall i \in \{5,6,7,8\}$, $S_3S_i=S_iS_3$.
\item $\forall i \in \{6,7,8\}$, $S_4S_i=S_iS_4$.
For $\sigma$ in the Coxeter group $E_8$, if $\sigma=s_{i_1}\dots s_{i_k}$ is a reduced expression we set $T_{\sigma}=S_{i_1}\dots S_{i_k}$.
\end{enumerate}
\end{Def}

This means we consider $E_8$ as in the CHEVIE package of GAP3 \cite{CHEVIE} with the following Dynkin diagram

\begin{center}
\begin{tikzpicture}
[place/.style={circle,draw=black,
inner sep=1pt,minimum size=10mm}]
\node (1) at (0,0)[place]{$S_1$};
\node (2) at (2,0)[place]{$S_3$};
\node (3) at (4,0)[place]{$S_4$};
\node (4) at (4,2)[place]{$S_2$};
\node (5) at (6,0)[place]{$S_5$};
\node (6) at (8,0)[place]{$S_6$};
\node (7) at (10,0)[place]{$S_7$};
\node (8) at (12,0)[place]{$S_8$};
\draw (1) to (2);
\draw (2) to (3);
\draw (3) to (4);
\draw (3) to (5);
\draw (5) to (6);
\draw (6) to (7);
\draw (7) to (8);
\end{tikzpicture}
\end{center}

\begin{prop}
Under our assumptions on $p$ and $\alpha$, $\mathcal{H}_{E_8,\alpha}$ is split semisimple, the representations afforded by the $W$-graphs are irreducible and pairwise non-isomorphic over $\F_q$. The restrictions of the irreducible representations of $\mathcal{H}_{E_8,\alpha}$ to $\mathcal{H}_{E_7,\alpha}$ are the same as in the generic case.
\end{prop}

\begin{proof}
We will apply Proposition \ref{Tits}. Let $A=\Z[\sqrt{u}^{\pm 1}]$ and $F=\Q(\sqrt{u})$. We have a symetrizing trace defined by $\tau(T_0)=1$ and $\tau(T_{\sigma})=0$ for all $\sigma\in E_8\setminus \{1_{E_8}\}$. By \cite{G-P}, $\mathcal{H}_{E_8,u}$ is then a free symmetric $F$-algebra of rank $696729600$. By \cite{Bourb} V.3. Corollary 1, $A$ is integrally closed. Let $\theta$ be the ring homomorphism from $A$ to $L=\F_q$ defined by $\theta(u)=\alpha$ and $\theta(k)=\overline{k}$. We know $FH$ is split. The basis formed by the elements $T_\sigma$, $\sigma\in E_8$ satisfies the conditions of Proposition \ref{Tits}. We work with the $E_8$-graphs implemented in CHEVIE \cite{CHEVIE} for the non-self dual $E_8$-graphs. For the self-dual ones, we use the new ones we found \cite{newgraphsEsterle} which we denote by $\tilde{\rho}$ if $\rho$ is the initial representation, they verify the properties of Theorem \ref{bilinwgraphs} and they all verify $\omega(e_{x_1})\omega(e_{x_n})=-1$ for any $2$-coloring $\omega$. For the $4536$-dimensional representation, we need to define two different $E_8$-graphs over $\Q(\sqrt{u})$ because the bilinear form is different for $p=11$, therefore we work with the one in CHEVIE. Going through the weights of the $E_8$-graphs in CHEVIE \cite{CHEVIE} and \cite{newgraphsEsterle}, we see that most of the $E_8$-graphs considered remain connected since the weights of the $E_8$-graphs which do not belong to $\{2240_x,2240_x',4200_y,4200_y',\widetilde{4480_y},\widetilde{5670_y},\widetilde{7168_w}\}$ belong to the set composed of $-24$, $-16$, $-10$, $-8$, $-20/3$, $-6$, $-5$, $-4$, $-10/3$, $-3$, $-8/3$, $-5/2$, $-2$, $-5/3$,  $-3/2$, $-4/3$, $-5/4$, $-1$, $-8/9$, $-5/6$, $-3/4$, $-2/3$, $-5/8$, $-3/5$, $-1/2$, $-5/12$, $-2/5$, $-3/8$, $-1/3$, $-1/4$, $-1/8$, $1/9$, $1/4$, $1/3$, $3/8$, $2/5$, $1/2$, $5/8$, $2/3$, 
  $3/4$, $5/6$, $8/9$, $1$, $4/3$, $3/2$, $5/3$, $2$, $5/2$, $8/3$, $3$, $10/3$, $4$, $5$, $6$, $8$.
  
  \medskip
  
  We now consider the remaining $E_8$-graphs separately. Consider first $2240_x$. We have four weighted edges which vanish when $p=7$ and none otherwise. Using CHEVIE \cite{CHEVIE},
  when $p=7$, the edges $1572\rightarrow 12$, $2183\rightarrow 1$, $2229\rightarrow 1$ and $2229\rightarrow 121$ all vanish, therefore we need to prove the $E_8$-graph remains connected when those edges vanish. The edge $1572\rightarrow 12$ can be replaced by the path $1572\rightarrow 24\rightarrow 100\rightarrow 12$. The edge $2183\rightarrow 1$ can be replaced by the path $2183\rightarrow 7\rightarrow 627 \rightarrow 1$. The edge $2229\rightarrow 1$ can be replaced by $2229\rightarrow 22\rightarrow 59\rightarrow 1$. The edge $2229\rightarrow 121$ can be replaced by the path $2229\rightarrow 2230 \rightarrow 2234 \rightarrow 121$. This proves $2240_x$ remains connected, this is true as well for $2240_x'$ because it is its dual $E_8$-graph.
  
  \smallskip
  
  Consider now $4200_y$. When $p\neq 7$, none of the weights vanish, therefore all the edges remain. For $p=7$, the edges $ 2700\rightarrow 20$, $ 3465\rightarrow 1 $, $ 3465\rightarrow 2$, $ 4075\rightarrow 894 $, $4172\rightarrow 1$, $ 4172\rightarrow 399$, $4190\rightarrow 2$, $4190\rightarrow  12$ and $4190\rightarrow 399$ are the only ones disappearing. The edge $2700\rightarrow 20$ can be replaced by the path $2700 \rightarrow 563 \rightarrow 630 \rightarrow 20$. The edge $3465\rightarrow 1$ can be replaced by the path $3465 \rightarrow 103 \rightarrow 119 \rightarrow 1$. The edge $3465\rightarrow 2$ can be replaced by the path $3465 \rightarrow 217 \rightarrow 287 \rightarrow 2$. The edge $4075\rightarrow 894$ can be replaced by the path $4075 \rightarrow 908 \rightarrow 2900 \rightarrow 894$. The edge $4172\rightarrow 1$ can be replaced by the path $4172 \rightarrow 141 \rightarrow 178\rightarrow 1$. The edge $4172\rightarrow 399$ can be replaced by the path $4172 \rightarrow 432 \rightarrow 2150\rightarrow 399$. The edge $4190\rightarrow 2$ can be replaced by the path $4190 \rightarrow 68 \rightarrow 197\rightarrow 2$. The edge $4190\rightarrow 12$ can be replaced by the path $4190 \rightarrow 31 \rightarrow 136\rightarrow 12$. The edge $4190\rightarrow 399$ can be replaced by the path $4190 \rightarrow 385 \rightarrow 20\rightarrow 399$. This proves $4200_y$ and $4200_y'$ remain connected.
  
  \smallskip
  
  Consider now $\widetilde{4480_y}$. When $p\notin \{7,11\}$, none of the weights vanish, therefore all the edges remain and the specialization is still connected. We give the proof of the connectedness for $p\in \{7,11\}$ in subsection \ref{conec4480} of the Appendix. The information is obtained using \cite{newgraphsEsterle}.

  \smallskip
  
Consider now $\widetilde{5670_y}$. When $p\neq 7$, none of the weights vanish, therefore all the edges remain. We give the proof of the connectedness for $p=7$ in subsection \ref{connec5670} of the Appendix.
  
\smallskip

Consider now $\widetilde{7168_w}$. When $p\neq 7$, none of the weights vanish, therefore the $E_8$-graph remains connected. We give the proof of the connectedness for $p=7$ in subsection \ref{connec7168} of the Appendix.  

\medskip

It now only remains to check that the Schur elements associated to the specialized representations are in $B$ and do not vanish under $\theta$ with $B$ as in Proposition \ref{Tits}. The Schur elements are given in Tables \ref{Schurelmts1E8}, \ref{Schurelmts11E8} and \ref{Schurelmts2E8}. They were obtained using Proposition $9.3.6$ and Table $E.7$ of \cite{G-P}. For a pair $(\rho,\rho^\star)$ of representations, we only give the Schur element of one of the representations since the other is obtained by applying the involution $\sqrt{u}\mapsto \sqrt{u}^{-1}$.  The conditions on $\alpha$ and $p$ imply that the Schur elements verify the right conditions and the proof is concluded.
\end{proof}

\begin{table}
\begin{enumerate}
\item $1_x$ : $(\Phi_2^8\Phi_3^4\Phi_4^4\Phi_5^2\Phi_6^4\Phi_7\Phi_8^2\Phi_9\Phi_{10}^2\Phi_{12}^2\Phi_{14}\Phi_{15}\Phi_{18}\Phi_{20}\Phi_{24}\Phi_{30})(u)$.
\item $8_z$ : $\frac{1}{u}(\Phi_2^8\Phi_3^4\Phi_4^2\Phi_5^2\Phi_6^4\Phi_7\Phi_8\Phi_9\Phi_{10}^2\Phi_{12}\Phi_{14}\Phi_{15}\Phi_{18}\Phi_{30})(u)$.
\item $28_x$ : $\frac{2}{u^3}(\Phi_2^8\Phi_3^4\Phi_4^2\Phi_5^2\Phi_6^4\Phi_8\Phi_9\Phi_{10}^2\Phi_{12}\Phi_{30})(u)$.
\item $35_x$ : $\frac{1}{u^2}(\Phi_2^8\Phi_3^4\Phi_4^4\Phi_5\Phi_6^4\Phi_8^2\Phi_9\Phi_{10}\Phi_{12}^2\Phi_{18}\Phi_{24})(u)$.
\item $50_x$ : $\frac{2}{u^4}(\Phi_2^8\Phi_3^4\Phi_4^4\Phi_6^4\Phi_7\Phi_9\Phi_{10}\Phi_{12}^2)(u)$.
\item $56_z$ : $\frac{6}{u^7}(\Phi_2^8\Phi_3^4\Phi_4^2\Phi_5^2\Phi_6^4\Phi_{30})(u)$.
\item $84_x$ : $\frac{2}{u^3}(\Phi_2^8\Phi_3^4\Phi_4^2\Phi_5^2\Phi_6^4\Phi_8\Phi_{10}^2\Phi_{12}\Phi_{15}\Phi_{18})(u)$.
\item $112_z$ : $\frac{2}{u^3}(\Phi_2^4\Phi_3^4\Phi_4^4\Phi_5^2\Phi_6^2\Phi_8\Phi_9\Phi_{12}^2\Phi_{15}\Phi_{20})(u)$.
\item $160_z$ : $\frac{2}{u^4}(\Phi_2^4\Phi_3^4\Phi_4^2\Phi_5\Phi_6^2\Phi_7\Phi_8^2\Phi_9\Phi_{12}\Phi_{24})(u)$.
\item $175_x$ : $\frac{3}{u^8}(\Phi_2^8\Phi_3^4\Phi_4^4\Phi_6^4\Phi_8^2)(u)$.
\item $210_x$ : $\frac{2}{u^4}(\Phi_2^8\Phi_3^4\Phi_4^4\Phi_5\Phi_6^4\Phi_{12}^2\Phi_{14}\Phi_{18})(u)$.
\item $300_x$ : $\frac{2}{u^6}(\Phi_2^8\Phi_3^4\Phi_4^2\Phi_6^4\Phi_7\Phi_8\Phi_{12}\Phi_{18})(u)$.
\item $350_x$ : $\frac{6}{u^8}(\Phi_2^8\Phi_3^4\Phi_4^4\Phi_6^4\Phi_{24})(u)$.
\item $400_z$ : $\frac{2}{u^6}(\Phi_2^4\Phi_3^4\Phi_4^4\Phi_6^2\Phi_7\Phi_8\Phi_9\Phi_{12}^2)(u)$.
\item $448_z$ : $\frac{3}{u^7}(\Phi_2^8\Phi_3^4\Phi_5^2\Phi_6^4\Phi_{10}^2\Phi_{12})(u)$.
\item $525_x$ : $\frac{1}{u^{12}}(\Phi_2^8\Phi_3^4\Phi_4^4\Phi_6^4\Phi_8^2\Phi_{12}^2)(u)$.
\item $560_z$ : $\frac{1}{u^5}(\Phi_2^8\Phi_3^4\Phi_4^2\Phi_5\Phi_6^4\Phi_9\Phi_{10}\Phi_{12}\Phi_{18})(u)$.
\item $567_x$ : $\frac{1}{u^6}(\Phi_2^8\Phi_3\Phi_4^4\Phi_5^2\Phi_6\Phi_8^2\Phi_{10}^2\Phi_{20})(u)$.
\item $700_x$ : $\frac{2}{u^6}(\Phi_2^8\Phi_3^4\Phi_4^2\Phi_6^4\Phi_8\Phi_9\Phi_{12}\Phi_{14})(u)$.
\item $700_{xx}$ : $(\Phi_2^8\Phi_3^4\Phi_4^2\Phi_8\Phi_9)(u)$.
\item $840_x$ : $\frac{2}{u^{12}}(\Phi_2^8\Phi_3^4\Phi_4^2\Phi_5\Phi_6^4\Phi_{12})(u)$.
\item $840_z$ : $\frac{2}{u^{10}}(\Phi_2^8\Phi_3^4\Phi_4^2\Phi_5\Phi_6^4\Phi_{18})(u)$.
\item $972_x$ : $\frac{2}{u^{10}}(\Phi_2^8\Phi_4^2\Phi_5^2\Phi_6\Phi_7\Phi_8\Phi_{10}^2)(u)$.
\item $1008_z$ : $\frac{3}{u^7}(\Phi_2^8\Phi_3\Phi_4^2\Phi_5^2\Phi_6\Phi_9\Phi_{10}^2\Phi_{18})(u)$.
\item $1050_x$ : $\frac{2}{u^8}(\Phi_2^8\Phi_3^4\Phi_4^4\Phi_6^4\Phi_{12}^2)(u)$.
\item $1296_z$ : $\frac{2}{u^{10}}(\Phi_2^4\Phi_3\Phi_4^4\Phi_5^2\Phi_7\Phi_8\Phi_{20})(u)$.
\item $1344_x$ : $\frac{2}{u^7}(\Phi_2^4\Phi_3^4\Phi_4^2\Phi_5^2\Phi_6^2\Phi_8\Phi_{12}\Phi_{15})(u)$.
\end{enumerate}
\caption{Schur elements for non self-dual representations in type $E_8$}
\label{Schurelmts1E8}
\end{table}

\begin{table}
\begin{enumerate}
\item $1400_x$ : $\frac{6}{u^8}(\Phi_2^8\Phi_3^4\Phi_6^4\Phi_8^2\Phi_{12}^2)(u)$.
\item $1400_z$ : $\frac{6}{u^7}(\Phi_2^8\Phi_3^4\Phi_4^2\Phi_6^4\Phi_{10}^2\Phi_{15})(u)$.
\item $1400_{zz}$ : $\frac{2}{u^{10}}(\Phi_2^8\Phi_3^4\Phi_4^2\Phi_6^4\Phi_9\Phi_{10})(u)$.
\item $1575_x$ : $\frac{3}{u^8}(\Phi_2^8\Phi_3\Phi_4^4\Phi_6\Phi_8^2\Phi_9\Phi_{18})(u)$.
\item $2100_x$ : $\frac{2}{u^{13}}(\Phi_2^8\Phi_3^4\Phi_4^2\Phi_6^4\Phi_8\Phi_{18})(u)$.
\item $2240_x$ : $\frac{2}{u^{10}}(\Phi_2^4\Phi_3^4\Phi_4^2\Phi_5\Phi_6^2\Phi_8\Phi_9\Phi_{12})(u)$.
\item $2268_x$ : $\frac{2}{u^{10}}(\Phi_2^8\Phi_3\Phi_4^2\Phi_5^2\Phi_8\Phi_{10}^2\Phi_{14})(u)$.
\item $2400_z$ : $\frac{2}{u^{15}}(\Phi_2^8\Phi_3^4\Phi_6^4\Phi_7\Phi_{18})(u)$.
\item $2800_z$ : $\frac{2}{u^{13}}(\Phi_2^4\Phi_3^4\Phi_4^4\Phi_8\Phi_9\Phi_{12}^2)(u)$.
\item $2835_x$ : $\frac{1}{u^{14}}(\Phi_2^8\Phi_3\Phi_4^4\Phi_5\Phi_6\Phi_8^2\Phi_{10})(u)$.
\item $3200_x$ : $\frac{2}{u^{15}}(\Phi_2^4\Phi_3^4\Phi_4^2\Phi_6^2\Phi_7\Phi_9\Phi_{12})(u)$.
\item $3240_z$ : $\frac{1}{u^9}(\Phi_2^8\Phi_3\Phi_4^2\Phi_5\Phi_6\Phi_7\Phi_8\Phi_{10}\Phi_{14})(u)$.
\item $3360_z$ : $\frac{2}{u^{12}}(\Phi_2^4\Phi_3^4\Phi_4^4\Phi_5\Phi_6^2\Phi_{12}^2)(u)$.
\item $4096_x$ : $\frac{2}{u^{11}}(\Phi_2\Phi_3^4\Phi_5^2\Phi_7\Phi_9\Phi_{15})(u)$.
\item $4096_z$ : $\frac{2}{u^{11}}(\Phi_2\Phi_3^4\Phi_5^2\Phi_7\Phi_9\Phi_{15})(u)$.
\item $4200_x$ : $\frac{2}{u^{12}}(\Phi_2^8\Phi_3^4\Phi_4^2\Phi_6^4\Phi_{10}\Phi_{12})(u)$.
\item $4200_z$ : $\frac{1}{u^{15}}(\Phi_2^8\Phi_3^4\Phi_4^2\Phi_6^4\Phi_8\Phi_{12})(u)$.
\item $4536_z$ : $\frac{1}{u^{13}}(\Phi_2^8\Phi_3\Phi_4^2\Phi_5^2\Phi_6\Phi_8\Phi_{10}^2)(u)$.
\item $5600_z$ : $\frac{2}{u^{15}}(\Phi_2^8\Phi_3^4\Phi_6^4\Phi_9\Phi_{14})(u)$.
\item $6075_x$ : $\frac{1}{u^{14}}(\Phi_2^8\Phi_4^4\Phi_7\Phi_8^2\Phi_{14})(u)$.
\end{enumerate}
\caption{Schur elements for non self-dual representations in type $E_8$}
\label{Schurelmts11E8}
\end{table}

\begin{table}
\begin{enumerate}
\item $70_y$ : $\frac{30}{u^{16}}(\Phi_2^8\Phi_3^4\Phi_4^4\Phi_{30})(u)$.
\item $168_y$ : $\frac{8}{u^{16}}(\Phi_2^8\Phi_3^4\Phi_5^2\Phi_6^4)(u)$.
\item $420_y$ : $\frac{5}{u^{16}}(\Phi_2^8\Phi_3^4\Phi_4^4\Phi_6^4)(u)$.
\item $448_w$ : $\frac{12}{u^{16}}(\Phi_2^4\Phi_3^4\Phi_4^2\Phi_5^2\Phi_{24})(u)$.
\item $1134_y$ : $\frac{6}{u^{16}}(\Phi_2^8\Phi_4^4\Phi_5^2\Phi_6\Phi_{18})(u)$.
\item $1344_w$ : $\frac{4}{u^{16}}(\Phi_2^4\Phi_3^4\Phi_4^2\Phi_5^2\Phi_6^2\Phi_{12})(u)$.
\item $1400_y$ : $\frac{24}{u^{16}}(\Phi_2^8\Phi_3^4\Phi_{10}^2\Phi_{12}^2)(u)$.
\item $1680_y$ : $\frac{20}{u^{16}}(\Phi_2^8\Phi_3^4\Phi_6^4\Phi_{20})(u)$.
\item $2016_w$ : $\frac{6}{u^{16}}(\Phi_2^4\Phi_3\Phi_4^2\Phi_5^2\Phi_8^2\Phi_9)(u)$.
\item $2100_y$ : $\frac{1}{u^{20}}(\Phi_2^8\Phi_3^4\Phi_4^4\Phi_6^4\Phi_{12}^2)(u)$.
\item $2688_y$ : $\frac{8}{u^{16}}(\Phi_3^4\Phi_4^4\Phi_5^2\Phi_{12}^2)(u)$.
\item $3150_y$ : $\frac{6}{u^{16}}(\Phi_2^8\Phi_3\Phi_4^4\Phi_9\Phi_{10}^2)(u)$.
\item $4200_y$ : $\frac{8}{u^{16}}(\Phi_2^8\Phi_3^4\Phi_6^4\Phi_{10}^2)(u)$.
\item $4480_y$ : $\frac{120}{u^{16}}(\Phi_3^4\Phi_4^4\Phi_6^4\Phi_9\Phi_{10}^2)(u)$.
\item $4536_y$ : $\frac{24}{u^{16}}(\Phi_2^8\Phi_5^2\Phi_6^4\Phi_{12}^2)(u)$.
\item $5600_w$ : $\frac{6}{u^{16}}(\Phi_2^4\Phi_3^4\Phi_4^2\Phi_8^2\Phi_{15})(u)$.
\item $5670_y$ : $\frac{30}{u^{16}}(\Phi_2^8\Phi_4^4\Phi_6^4\Phi_{15})(u)$.
\item $7168_w$ : $\frac{12}{u^{16}}(\Phi_3^4\Phi_5^2\Phi_6^2\Phi_8^2\Phi_{12})(u)$.
\end{enumerate}
\caption{Schur elements for self-dual representations in type $E_8$}
\label{Schurelmts2E8}
\end{table}

The restriction table from $\mathcal{H}_{E_8,\alpha}$ to its subalgebra $\mathcal{H}_{E_7,\alpha}$ generated by $S_1,S_2,S_3,S_4,S_5,S_6$ and $S_7$ which is naturally isomorphic to the Iwahori-Hecke of type $E_7$ with parameter $\alpha$ is then given by Tables \ref{resE8E71} and \ref{resE8E72}. It is obtained using the CHEVIE package of GAP3 \cite{CHEVIE}. They correspond in the generic case to the induction/restriction tables of the corresponding Coxeter groups.

\begin{table}
$\begin{array}{lll}
1_x & = & 1_a \\
8_z &= & 1_a+7_a'\\
28_x & = & 7_a'+21_a\\
35_x & = & 1_a+7_a'+27_a\\
50_x & = & 15_a'+35_b\\
56_z & = & 21_a+35_a'\\
70_y & = & 35_a+35_a'\\
84_x & = & 1_a+21_b'+27_a+35_b\\
112_z & = & 1_a+7_a'+21_b'+27_a+56_a'\\
160_z & = & 7_a'+21_a+27_a+105_a'\\
168_y & = & 84_a+84_a'\\
175_x & = & 70_a'+105_b\\
210_x & = & 7_a'+27_a+56_a'+120_a\\
300_x & = & 27_a+105_a'+168_a\\
350_x & = & 21_a+35_a'+105_a'+189_a\\
400_z & = & 15_a'+35_b+56_a'+105_b+189_b'\\
420_y & = & 210_b+210_b'\\
448_w  & = & 35_a+35_a'+189_a+189_a'\\
448_z & = & 21_b'+70_a'+168_a+189_b'\\
525_x & =& 21_b'+105_c+189_c'+210_a\\
560_z & = & 7_a'+21_a'+27_a+35_b+56_a'+105_a'+120_a+189_b'\\
567_x & = & 7_a'+21_a+21_b'+27_a+56_a'+105_a'+120_a+210_a\\
700_x & =& 27_a+35_b+56_a'+105_b+120_a+168_a+189_b'\\
700_{xx} & = & 15_a'+84_a'+105_c+216_a'+280_b\\
840_x & = & 84_a+168_a+210_b+378_a'\\
840_z & = & 105_a'+168_a+189_a+378_a'\\
972_x & = & 35_b+84_a+168_a+189_c'+216_a'+280_b\\
1008_z & = & 21_a+27_a+56_a'+105_a'+120_a+189_c'+210_a+280_a'\\
1050_x & = & 15_a'+35_b+105_b+189_b'+210_a+216_a'+280_b\\
1134_y & = & 189_a+189_a'+378_a+378_a'\\
1296_z & = & 21_a+35_a'+105_a'+120_a+189_a+210_a+280_a'+336_a'\\
1344_x & = & 21_b'+27_a+35_b+105_a'+120_a+168_a+189_b'+189_c'+210_a+280_b\\
1344_w & =& 84_a+84_a'+210_b+210_b'+378_a+378_a'\\
1400_x & = & 56_a'+105_b+120_a+189_b'+210_a+315_a'+405_a\\
1400_y & = & 280_a+280_a'+420_a+420_a'\\
1400_z & = & 21_a'+27_a+56_a'+105_a'+120_a+168_a+189_b'+189_c'+210_a+315_a'\\
1400_{zz} & = & 15_a'+70_a'+105_b+189_b'+210_b+216_a'+280_b+315_a'\\
1575_x &  = & 21_a+56_a'+105_a'+120_a+189_a+189_b'+210_a+280_a'+405_a\\
1680_y & = & 35_a+35_a'+189_a+189_a'+280_a+280_a'+336_a+336_a'\\
2016_w & = & 70_a+70_a'+210_b+210_b'+216_a+216_a'+512_a+512_a'\\
2100_x & = & 35_a'+105_a'+189_a+189_c'+210_a+280_a'+336_a+336_a'+420_a\\
2100_y & = & 105_c+105_c'+189_c+189_c'+336_a+336_a'+420_a+420_a'
\end{array}$\smallskip
\caption{Restriction table from $\mathcal{H}_{E_8,\alpha}$ to $\mathcal{H}_{E_7,\alpha}$.}
\label{resE8E71}
\end{table}

\begin{table}
$\begin{array}{lll}
2240_x & =&  70_a'+105_b+120_a+168_a+189_b'+210_b+280_b+315_a'+378_a'+405_a\\
2268_x &  = & 56_a'+105_a'+120_a+168_a+189_c'+210_a+280_a'+315_a'+405_a+420_a\\
2400_z & = & 35_a'+189_a+189_a'+210_a+280_a'+336_a+336_a'+405_a+420_a'\\
2688_y & = & 216_a+216_a'+280_b+280_b'+336_a+336_a'+512_a+512_a'\\
2800_z & =& 105_c+120_a+189_c'+210_a+280_a'+315_a'+336_a'+405_a+420_a+420_a' \\
2835_x & = & 70_a'+105_b+210_b+210_b'+216_a'+280_b+315_a'+405_a+512_a+512_a' \\ 
3150_y & = & 70_a+70_a'+210_b+210_b'+378_a+378_a'+405_a+405_a'+512_a+512_a'\\
3200_x & = & 84_a+105_c'+168_a+189_c'+216_a+280_b+336_a'+378_a'+420_a+512_a+512_a'\\
3240_z & =& 35_b+56_a'+105_a'+105_b+120_a+168_a+2\times 189_b'+189_c'+210_a\\
& & +216_a'+280_a'+280_b+315_a'+378_a'+405_a\\
3360_z & = & 105_b+189_b'+210_a+216_a'+280_a'+280_b+315_a'+336_a'+405_a+512_a+512_a'\\
4096_x & = & 105_a'+120_a+168_a+189_a+189_b' +189_c'+210_a\\
& &+280_a'+280_b+315_a'+336_a'+378_a'+405_a+420_a+512_a\\
4096_z & = & 105_a'+120_a+168_a+189_a+189_b' +189_c'+210_a\\
& &+280_a'+280_b+315_a'+336_a'+378_a'+405_a+420_a+512_a'\\
4200_x & = & 105_b+168_a+189_b'+210_a+210_b+216_a'+280_a'+280_b+315_a'+378_a'+405_a\\
& & +420_a+512_a+512_a'\\
4200_y & = & 84_a+84_a'+210_b+210_b'+216_a+216_a'+280_b+280_b'+378_a+378_a'\\
& & +420_a+420_a'+512_a+512_a'\\
4200_z & = & 84_a'+105_b+105_c+210_a+216_a'+280_a'+280_b+315_a'+378_a+378_a'\\
& & +405_a+420_a'+512_a+512_a'\\
4480_y & = & 210_b+210_b'+315_a+315_a'+378_a+378_a'+405_a+405_a'+420_a+420_a'+512_a+512_a' \\
4536_y & = & 280_a+280_a'+315_a+315_a'+336_a+336_a'+405_a+405_a'+420_a+420_a'+512_a+512_a'\\
4536_z & = & 70_a'+84_a+168_a+189_b'+189_c'+210_b+210_b'+216_a'+280_b+315_a'\\
& & +2\times 378_a'+405_a+420_a+512_a+512_a'\\
5600_w & = & 189_a+189_a'+280_a+280_a'+280_b+280_b'+336_a+336_a'+378_a+378_a'\\
& & +405_a+405_a'+420_a+420_a'+512_a+512_a'\\
5600_z & = & 105_c'+168_a+189_a+189_c'+210_a+280_a'+315_a'+336_a+336_a'+378_a'+405_a+405_a'\\
& & +2\times 420_a+420_a'+512_a+512_a'\\
5670_y & = & 189_a+189_a'+280_a+280_a'+315_a+315_a'+336_a+336_a'+378_a+378_a'\\
& & +405_a+405_a'+420_a+420_a'+512_a+512_a'\\
6075_x & = & 105_c+189_a+189_b'+189_c'+210_a+216_a'+280_a'+280_b+315_a'+336_a+336_a'\\
& & + 378_a+378_a'+2\times 405_a+420_a+420_a'+512_a+512_a'\\
7168_w & = & 210_b+210_b'+216_a+216_a'+280_b+280_b'+315_a+315_a'+336_a+336_a'+378_a+378_a'\\
 & &+ 405_a+405_a'+420_a+420_a'+2\times 512_a+2\times 512_a' 
\end{array}$
\smallskip
\caption{Restriction table from $\mathcal{H}_{E_8,\alpha}$ to $\mathcal{H}_{E_7,\alpha}$}
\label{resE8E72}
\end{table}

\begin{prop}\label{resE8derivedsubgroup}
The restrictions to $\mathcal{A}_{E_8}$ of the representations afforded by those $E_8$-graphs are absolutely irreducible and the representations of dimension greater than $1$ are pairwise non-isomorphic.
\end{prop}

\begin{proof}
As in \cite{BMM} Lemma $3.4$, we only need to prove that $A_{E_8}$ is generated by $A_{E_7}$ and $\mathcal{A}_{E_8}$. This true because $s_8=s_8s_1^{-1}s_1$, $s_8s_1^{-1}\in \mathcal{A}_{E_7}$ and $s_1\in A_{E_7}$.

We now prove the second part of the statement. Let $\rho_1$ and $\rho_2$ be two irreducible representations of $\mathcal{H}_{E_8,\alpha}$ such that $\rho_{1|\mathcal{A}_{E_8}}\simeq \rho_{1|\mathcal{A}_{E_8}}$. By Lemma \ref{abel},there exists a character $\xi:A_{E_8}\mapsto \F_q^{\star}$ such that $\rho_1\simeq \rho_2\otimes \xi$. This means there exists $x\in \F_q$ such that for all $i\in [\![1,8]\!]$, $\rho_1(S_i)$ is conjugate to $x\rho_2(S_i)$. We know for any representation $\rho$ of dimension greater than $1$, the set of eigenvalues of $\rho$ is equal to $\{\alpha,-1\}$. This implies that $\{\alpha,-1\}=\{x\alpha,-x\}$. We then have $x=1$ or ($x=-\alpha$ and $-\alpha^2=-1$). It follows that $x=1$ and $\rho_1\simeq \rho_2$.
\end{proof}

We now prove a lemma computing the normal closure of $\mathcal{A}_{E_7}$ inside $\mathcal{A}_{E_8}$ as we did in the other types.

\begin{lemme}\label{normclosE7E8}
The normal closure $\ll \mathcal{A}_{E_7} \gg_{\mathcal{A}_{E_8}}$ of $\mathcal{A}_{E_7}$ inside $\mathcal{A}_{E_8}$ is equal to $\mathcal{A}_{E_8}$, where we identify $A_{E_7}$ as a subgroup of $A_{E_8}$ using the natural isomorphism from $A_{E_7}$ to $<S_i, i\in [\![1,7]\!]>$.
\end{lemme}

\begin{proof}
By \cite{MR}, we have $\mathcal{A}_{E_8}=<S_1S_3S_1^{-1},S_4S_3^{-1},S_iS_1^{-1}, i\in [\!
[2,8]\!]>$ and 
$$\mathcal{A}_{E_7}=<S_1S_3S_1^{-1},S_4S_3^{-1},S_iS_1^{-1}, i\in [\![2,7]\!]>.$$
 This proves we only need to show that $S_8S_1^{-1}\in \ll \mathcal{A}_{E_7} \gg_{\mathcal{A}_{E_8}}$.
We have $S_8S_7S_8=S_7S_8S_7$, therefore $S_8=S_7S_8S_7(S_7S_8)^{-1}$ and 
$$S_8S_1^{-1}=S_7S_8S_7S_1^{-1}(S_7S_8)^{-1}=((S_7S_1^{-1})(S_8S_1^{-1}))(S_7S_1^{-1})((S_7S_1^{-1})(S_8S_1^{-1}))^{-1}.$$
This concludes the proof.
\end{proof}

We now determine the image of $\mathcal{A}_{E_8}$ inside each given representation. Note that the new self-dual $W$-graphs we found all verify $\omega(e_1)\omega(e_n)=1$ for any $2$-coloring $\omega$, therefore the bilinear forms appearing are self-dual. For the self-dual representation of dimension 4536, since we do not work with the new $W$-graph, we use the bilinear form which is availiable at \cite{newgraphsEsterle} and see that it is symmetric. The $E_8$-graphs are all $2$-colorable except for $4096_x$, $4096_x'$, $4096_z$ and $4096_z'$.

\begin{prop}\label{ImageinsiderepsE8}
Assume $\F_p(\sqrt{\alpha})\neq \F_q=\F_p(\alpha)=\F_p(\alpha+\alpha^{-1})$. Let $\rho$ be a representation of dimension $n_\rho$ associated to a $W$-graph. If $\rho$ is not self-dual then $\rho(\mathcal{A}_{E_8})\simeq SL_{n_\rho}(q)$ if $n_\rho\neq 4096$ and $\rho(\mathcal{A}_{E_8})=\simeq SL_{n_\rho}(q^2)$ if $n_\rho=4096$.  If $\rho$ is self-dual then $\rho(\mathcal{A}_{E_8})\simeq \Omega_{n_\rho}^+(q)$. 

Assume $\F_q=\F_p(\sqrt{\alpha})=\F_p(\alpha)=\F_p(\alpha+\alpha^{-1})$. Let $\rho$ be a representation of dimension $n_\rho$ associated to a $W$-graph. If $\rho$ is not self-dual then $\rho(\mathcal{A}_{E_8})\simeq SL_{n_\rho}(q)$. If $\rho$ is self-dual then $\rho(\mathcal{A}_{E_8})\simeq \Omega_{n_\rho}^+(q)$. 

Assume $\F_q=\F_p(\alpha)\neq \F_p(\alpha+\alpha^{-1})$. Let $\rho$ be a representation of dimension $n_\rho$ associated to a $W$-graph. If $\rho$ is not self-dual then $\rho(\mathcal{A}_{E_8})\simeq SU_{n_\rho}(q^{\frac{1}{2}})$.  If $\rho$ is self-dual then $\rho(\mathcal{A}_{E_8})\simeq \Omega_{n_\rho}^+(q^{\frac{1}{2}})$. 
\end{prop}

\begin{proof}
Assume first $\F_q=\F_p(\alpha)=\F_p(\alpha+\alpha^{-1})$. Let $\rho$ be a non-self-dual representation of dimension $n_{\rho}$ associated to a $W$-graph. 
\noindent
Assume first that the corresponding $E_8$-graph is $2$-colorable. We can then consider the representation is defined over $\F_q$. If $n_\rho\geq 28$ then we can apply Theorem \ref{CGFS} since by Tables \ref{resE8E71} and \ref{resE8E72} and Proposition \ref{ImageinsiderepsE7}, $\rho(\mathcal{A}_{E_7})$ contains a natural $SL_7(q)$. It follows that $\rho(\mathcal{A}_{E_8})$ is a classical group over $\F_q$. We know that no non-degenerate bilinear form is preserved by this group because the representation is not self-dual. It cannot be unitary because it contains a natural $SL_7(q)$. We can then conclude that $\rho(\mathcal{A}_{E_8})=SL_{n_\rho}(q)$. We know $\rho_{8_z}(\mathcal{A}_{E_8})$ is an irreducible group normally generated by $\rho_{8_z}(\mathcal{A}_{E_7})$. Since $\rho_{8_z}(\mathcal{A}_{E_7})$ is a natural $SL_7(q)$, we have that $\rho_{8_z}(\mathcal{A}_{E_8})$ is an irreducible subgroup of $GL_8(q)$ generated by transvections, therefore by Theorem \ref{transvections}, $\rho_{8_z}(\mathcal{A}_{E_8})$ is isormorphic to $SL_8(q')$, $SU_8(q'^{\frac{1}{2}})$ or $SP_8(q')$ for some $q'$ dividing $q$. It contains a natural $SL_7(q)$, therefore $q'=q$ and $\rho_{8_z}(\mathcal{A}_{E_8})=SL_8(q)$. 
\noindent
Assume now that the corresponding $W$-graph is not $2$-colorable. The representations we have to consider are then the representations of dimension $4096$. If $\F_q=\F_p(\sqrt{\alpha})=\F_p(\alpha)$ then we can apply the previous reasoning. If $\F_p(\sqrt{\alpha})\neq \F_q=\F_p(\alpha)=\F_p(\alpha+\alpha^{-1})$ then $X^2-\alpha$ is an irreducible polynomial over $\F_q$. We then have $\F_{q^2}=\F_p(\sqrt{\alpha})=\F_q/(X^2-\alpha)$ and there is a unique field automorphism $\Psi$ of degree $2$ of $\F_{q^2}$, ot fixes $\F_q$ pointiwise and $\Psi(\sqrt{\alpha})=-\sqrt{\alpha}$. By Proposition \ref{ImageinsiderepsE7}, for any representation $\varphi$ of dimension different from $512$ of $\mathcal{A}_{E_7}$, we have $\Psi \circ \varphi \simeq \varphi$. We also have $\Psi \circ \rho_{512_a}\simeq \rho_{512_a'}$. By Table \ref{resE8E72}, we have that $\Psi\circ \rho_{4096_x}$ is not isomorphic to $\rho_{4096_x}$ or $\rho_{4096_z'}$ because otherwise we would have $\Psi\circ \rho_{512_a}\simeq \rho_{512_a}$, therefore $\rho_{512_a}\simeq \rho_{512_a'}$. We also have that $\Psi \circ \rho_{4096_x}$ is not isomorphic $\rho_{4096_x'}$ because otherwise, we would have $\Psi\circ \rho_{420_a}\circ \rho_{420_a'}$, therefore $\rho_{420_a}\simeq \rho_{420_a'}$. It follows that $\Psi \circ \rho_{4096_x}\simeq \rho_{4096_z}$. We know by the above reasonning that $\rho_{4096_x}(\mathcal{A}_{E_8})$ is a classical group over $\F_{q'}$ for some $q'$ dividing $q^2$. Furthermore, $q$ divides $q'$ since it contains a natural $SL_{420}(q)$. It does not preserve any non-degenerate bilinear form since $\rho$ is not self-dual. We cannot have $\rho_{4096_x}(\mathcal{A}_{E_8})\simeq SL_{4096}(q)$ because $\Psi\circ \rho_{4096_x}$ is not isomorphic to $\rho_{4096_x}$. We also cannot have $\rho_{4096_x}(\mathcal{A}_{E_8})\simeq SU_{4096}(q)$ because $\rho_{4096_x}$ is not isomorphic to $\Psi\circ \rho_{4096_x'}$. It follows that $\rho_{4096_x}(\mathcal{A}_{E_8})\simeq SL_{4096}(q^2)$. We get the same result for the remaining representations of dimension $4096$ using $\rho_{4096_x'}\simeq \rho_{4096_x}^{\star}$ and $\rho_{4096_z}\simeq \Psi\circ\rho_{4096_x}$.

\medskip

Assume now that $\F_q=\F_p(\alpha)\neq \F_p(\alpha+\alpha^{-1})$. Then the same argument as in the proof of Proposition \ref{ImageinsiderepsE7} gives $\F_q=\F_p(\sqrt{\alpha})$. Using the same reasoning as above, we only need to prove that for any representation $\rho$, we have $\epsilon \circ \rho\simeq \rho$, where $\epsilon$ is the unique field automorphism of order $2$ of $\F_q$. The polynomial $X^2-(\alpha+\alpha^{-1})X+1$ is an irreducible $\F_p(\alpha+\alpha^{-1})$-polynomial of degree $2$, therefore we have $\epsilon(\alpha)=\alpha^{-1}$. We know by Proposition \ref{ImageinsiderepsE7} that $\epsilon \circ \varphi^\star \simeq \varphi$ for any representation $\varphi$ of $\mathcal{A}_{E_7}$. Finally, by Table \ref{resE8E71} and Table \ref{resE8E72}, no pair of distinct representations have the same restriction to $\mathcal{A}_{E_7}$. The result follows.
\end{proof}

We can now state the main theorem for type $E_8$.

\begin{theo}\label{resultE8}
We let $A$ be a set of representatives of the irreducible $2$-colorable non self-dual representations for the equivalence relation $\rho \thickapprox \varphi$ if $\rho=\varphi'$ and $B$ be the set of irreducible self-dual representations.

Assume $\F_p(\sqrt{\alpha})\neq\F_q=\F_p(\alpha)=\F_p(\alpha+\alpha^{-1})$. Then the morphism from $\mathcal{A}_{E_8}$ to $\mathcal{H}_{E_8,\alpha}^\star\simeq \underset{\rho \mbox{ irr}}\prod GL_{n_\rho}(q)$ factorizes through the surjective morphism
$$\Phi : \mathcal{A}_{E_8} \rightarrow \underset{\rho\in A}\prod SL_{n_{\rho}}(q) \times SL_{4096}(q^2)\times \underset{\rho\in B}\prod \Omega_{n_\rho}^+(q).$$

Assume $\F_q=\F_p(\sqrt{\alpha})=\F_p(\alpha)=\F_p(\alpha+\alpha^{-1})$. Then the morphism from $\mathcal{A}_{E_8}$ to $\mathcal{H}_{E_8,\alpha}^\star\simeq \underset{\rho \mbox{ irr}}\prod GL_{n_\rho}(q)$ factorizes through the surjective morphism
$$\Phi : \mathcal{A}_{E_8} \rightarrow \underset{\rho\in A}\prod SL_{n_{\rho}}(q) \times SL_{4096}(q)^2\times \underset{\rho\in B}\prod \Omega_{n_\rho}^+(q).$$

Assume $\F_q=\F_p(\alpha)\neq\F_p(\alpha+\alpha^{-1})$. Then the morphism from $\mathcal{A}_{E_8}$ to $\mathcal{H}_{E_8,\alpha}^\star\simeq \underset{\rho \mbox{ irr}}\prod GL_{n_\rho}(q)$ factorizes through the surjective morphism
$$\Phi : \mathcal{A}_{E_8} \rightarrow \underset{\rho\in A}\prod SU_{n_{\rho}}(q^{\frac{1}{2}}) \times SU_{4096}(q^{\frac{1}{2}})\times \underset{\rho\in B}\prod \Omega_{n_\rho}^+(q^\frac{1}{2}).$$
\end{theo}

\begin{proof}

By \cite{MR}, $\mathcal{A}_{E_8}$ is perfect. We can then apply exctly the same arguments as in the proof of of Proposition \ref{resultE6}.
\end{proof}

\chapter{Type $H$}\label{TypeH}

In this section, we determine the images of the Artin groups in types $H_3$ and $H_4$. The main difference with the other types is that in the generic case, some of the irreducible representations are defined over $\Q[\sqrt{5}](\sqrt{u})$ but not $\Q(\sqrt{u})$. We therefore need arguments as in the dihedral cases to understand the field extensions. Lemma \ref{IsomorphismH3} summarizes the information we need about the field extensions. The proof is then similar to the proof in the other types. We first determine the image inside each irreducible representation and then recover the full image using Goursat's Lemma. In order to determine the image inside each representation, we will use inductive arguments using the image in type $I_2(5)$ for type $H_3$ and the image in type $H_3$ for type $H_4$. The main results in this section are in Theorems \ref{resicos} and \ref{resH4}. The $8$-dimensional irreducible representations in type $H_4$ use additional assumptions on the order of $\alpha$ which might not be necessaray. The proof is highly computational, therefore it might be complicated to ommit those assumptions on the order of $\alpha$. The image in the product of those representations gives rise to a nice description of the $Spin_8^+(q)$ group as can be seen in Proposition \ref{Spin}. This section uses many results from chapter \ref{Wgraphschapter}. We proved Conjecture \ref{conjecturewgraphs} in type $H_4$ and the self-dual $H_4$-graphs we obtained are available in section \ref{sectionnewH4graphs} of the Appendix and at \cite{newgraphsEsterle}.

\section{Type $H_3$}

Let $p$ be a prime number, $p\notin \{2,5\}$ and $\alpha\in (\overline{\F_p})^\times$ such that the order of $\alpha$ does not divide $20$ and does not belong to $\{1,2,3,4,5,6,10\}$. Let $\xi \in \overline{\F_p}$ be a primitive fifth-root of unity. We set $\F_q=\F_p(\alpha)$ and $\F_r=\F_p(\alpha,\xi+\xi^{-1})$.

\begin{Def}
The Iwahori-Hecke $\mathcal{H}_{H_3,\alpha}$ is the $\F_q$-algebra generated by the generators $S_1,S_2,S_3$ and the following relations :
\begin{enumerate}
\item $\forall i \in \{1,2,3\}, (S_i-\alpha)(S_i+1)=0.$
\item $S_1S_2S_1S_2S_1=S_2S_1S_2S_1S_2$,
\item $S_1S_3=S_3S_1$,
\item $S_2S_3S_2=S_3S_2S_3$.
For $\sigma$ in the Coxeter group $H_3$, if $\sigma=s_{i_1}\dots s_{i_k}$ is a reduced expression we set $T_{\sigma}=S_{i_1}\dots S_{i_k}$.
\end{enumerate}
\end{Def}

\begin{prop}
Under the conditions we assumed on $\alpha$, the Iwahori-Hecke algebra is split semi-simple and the models given by specialization of the $W$-graphs are irreducible.
\end{prop}

\begin{proof}
We want to apply Proposition \ref{Tits}. Let $A=\Z[\frac{1+\sqrt{5}}{2}][\sqrt{u}^{\pm 1}]$, $B=\op{Frac}(A)$ and $F=\Q[\sqrt{5}](\sqrt{u})$. We have a symetrizing trace defined by $\tau(T_0)=1$ and $\tau(T_\sigma)=0$ for all $\sigma\in H_3\setminus \{1_{H_3}\}$. $\mathcal{H}_{H_3,u}$ is then a free $F$-algebra of rank $120$. 
$A$ is an integrally closed integral domain because $\Z[\frac{1+\sqrt{5}}{2}]$ is integrally closed (see \cite{Bourb} V.3. Corollary 1 and \cite{NivZucMont} Theorem 9.20.).

Let $\theta$ be the ring homomorphism from $A$ to $L=\F_q$ defined by $\theta(\frac{1+\sqrt{5}}{2})=\xi+\xi^{-1}+1$, $\theta(u)=\alpha$ and $\theta(k)=\overline{k}$. We know $FH$ is split. The basis formed by the elements $T_\sigma, \sigma \in H_3$ verifies the conditions of the Proposition \ref{Tits}.

All the $W$-graphs are connected and remain connected after we specialize the weights because none of them vanishes under $\theta$. Indeed, if they were not connected then they would afford reducible representations over $\Q(\sqrt{5})(\sqrt{u})$.

We now only need to check that the Schur elements associated to these irreducible representations are in $B$ and don't vanish under $\theta$. For $n\in \N$, $\Phi_n$ is the $n$-th cyclotomic polynomial and $\Phi_{5,a}(u)=u^2+\frac{1+\sqrt{5}}{2}u+1, \Phi_{5,b}(u)=u^2+\frac{1-\sqrt{5}}{2}u+1, \Phi_{10,a}(u)=u^2-\frac{1+\sqrt{5}}{2}u+1$ and $\Phi_{10,b}(u)=u^2+\frac{\sqrt{5}-1}{2}u+1$. If $\chi$ is an irreducible character then the character $\chi^\star$ associated to the dual representation of $\chi$ has a Schur element $c_{\chi^\star}=a(c_{\chi})$, where $a$ is the involution of $\Q(\sqrt{5})(\sqrt{u})$ sending $\sqrt{u}$ to $\sqrt{u}^{-1}$. We define the field automorphism of $\Q(\sqrt{5})$, written $x\mapsto \overline{x}$ by $\overline{\frac{\sqrt{5}-1}{2}}=\frac{-1-\sqrt{5}}{2}$ and $\overline{q}=q$ for all $q\in \Q$. The Schur elements are given in Table \ref{SchurelmtsH3} (obtained using Table E.2. of the Appendix  and corollary 9.3.6 of \cite{G-P}) 

\begin{table}
\centering
$\begin{array}{lcr}
1_r & : & (\Phi_2^3\Phi_3\Phi_5\Phi_6\Phi_{10})(u).\\
3_s' & :&  \dfrac{5+\sqrt{5}}{2}\dfrac{\Phi_2^3(u)\Phi_5(u)\Phi_{10}(u)}{u\Phi_{5,b}(u)\Phi_{10,b}(u)}.\\
 \overline{3_{s}}' & : & \dfrac{5-\sqrt{5}}{2}\dfrac{\Phi_2^3(u)\Phi_5(u)\Phi_{10}(u)}{u\Phi_{5,b}(u)\Phi_{10,b}(u)}.\\
 5_r & : & \dfrac{1}{u^2}\Phi_2^3(u)\Phi_3(u)\Phi_6(u).\\
4_{r'} &  : &  \dfrac{2}{u^3}\Phi_3(u)\Phi_5(u).
\end{array}$
\caption{Schur elements in type $H_3$}
\label{SchurelmtsH3}
\end{table}

Since $\alpha$ is of order not dividing $20$ and different from $6$ and $p\neq 2$, we only need to check that $\theta(\frac{5+\sqrt{5}}{2})\neq 0$, $\theta(\frac{5-\sqrt{5}}{2})\neq 0$, $\alpha^2\pm (\xi+\xi^{-1})\alpha+1\neq 0$, $\alpha^2\pm (\xi+\xi^{-1}+1)\alpha+1\neq 0$, $\alpha^2\mp (\xi^2+\xi^{-2})\alpha+1\neq 0$.

We have $\theta(\frac{5+\sqrt{5}}{2})=\xi+\xi^{-1}+3=-\xi^2-\xi^{-2}+2=-(\xi-\xi^{-1})^2\neq 0$ and $\theta(\frac{5-\sqrt{5}}{2})=-\xi-\xi^{-1}+2=-(\xi^2-\xi^{-2})^2\neq 0$.

We also have $\alpha^2+(\xi+\xi^{-1})\alpha+1=(\alpha+\xi)(\alpha+\xi^{-1})\neq 0$ because $\alpha$ is of order not dividing $20$. In the same way, $\alpha^2-(\xi+\xi^{-1})\alpha+1=(\alpha-\xi)(\alpha-\xi^{-1})\neq 0$. The last two inequalities are shown in exactly the same way because the order of $\xi$ is the same as the order of $\xi^2$.

\end{proof}

 We now show a lemma on the Artin groups $A_{H_3}$ and $A_{I_2(5)}$ which will be useful later on.

\begin{lemme}\label{normclosdih}
 We write $A_{H_3}=<t,s_1,s_2, ts_1ts_1t=s_1ts_1ts_1, s_1s_2s_1=s_2s_1s_2, ts_2=s_2t>$ and $A_{I_2(5)}$ its subgroup generated by $t$ and $s_1$. It is naturally isomorphic to the Artin group of type $I_2(5)$.

The derived subgroup $\mathcal{A}_{H_3}$ is equal to the normal closure of $\mathcal{A}_{I_2(5)}$ in $A_{H_3}$.
\end{lemme}

\begin{proof}

By \cite{MR}, we know that $\mathcal{A}_{I_2(5)}$ is generated by $s_1t^{-1}, ts_1t^{-2}$ and $t^2s_1t^{-3}$ and that $\mathcal{A}_{H_3}$ is generated by $s_1t^{-1}, ts_1t^{-2}$, $t^2s_1t^{-3}$ and $s_2t^{-1}$.

It is thus sufficient to show that $s_2t^{-1}$ can be written as a product of conjugates of the generators of $\mathcal{A}_{I_2(5)}$ or of their inverses. We have 
$$s_2t^{-1}=s_1s_2s_1s_2^{-1}s_1^{-1}t^{-1}=s_1(t^{-1}t)s_2s_1s_2^{-1}(t^{-2}t^2)s_1^{-1}t^{-1}=(s_1t^{-1})s_2(ts_1t^{-2})s_2^{-1}(ts_1t^{-2})^{-1}.$$

\end{proof}

\begin{prop}\label{resH3derivedsubgroup}
The restrictions to $\mathcal{A}_{H_3}$ of the representations of dimension greater than $1$ afforded by those $W$-graphs are absolutely irreducible and pairwise non-isomorphic.
\end{prop}

\begin{proof}
As in \cite{BMM} Lemma $3.4$, we only need to prove that $A_{H_3}$ is generated by $A_{I_2(5)}$ and $\mathcal{A}_{H_3}$ in order to prove the absolute irreducibility. This is true because $s_3=s_3s_1^{-1}s_1$ and $s_3s_1^{-1}\in \mathcal{A}_{H_3}$ and $s_1\in A_{I_2(5)}$.

\smallskip

Let $\rho_1$ and $\rho_2$ be irreducible representations afforded by $H_3$-graphs. Assume $\rho_{1|\mathcal{A}_{H_3}}\simeq \rho_{2|\mathcal{A}_{H_3}}$. By Lemma \ref{abel}, we have that there exists a character $\eta$ such that $\rho_{1}\simeq \rho_2\otimes \eta$. We have $A_{H_3}/\mathcal{A}_{H_3}=<\overline{S_1}>\simeq \Z$. It follows that there exists $u\in \overline{\F_p}^\star$ such that $\rho_1(S_1)$ is conjugate to $u\rho_2(S_1)$. The eigenvalues of $\rho_1(S_1)$ and $\rho_2(S_1)$ are $-1$ and $\alpha$. It follows that $\{-1,\alpha\}=\{-u,u\alpha\}$. We then have either $u=1$ or ($u=-\alpha$ and $\alpha^2=1$). The latter contradicts our assumptions on the order of $\alpha$ therefore $u=1$ and $\eta$ is the trivial morphism. It follows that $\rho_1\simeq \rho_2$.
\end{proof}

\begin{figure}
\begin{tikzpicture}

\node[circle, draw =black] (4) at (0,0) {$2 ~ 3$};
\node[circle, draw =black] (3) at (2.25,0) {$1 ~ 3$};
\node[circle, draw =black] (2) at (0,-2) {$2$};
\node[circle, draw =black] (1) at (2.25,-2) {$1$};
\draw (1) to (2);
\draw (1) to (3);
\draw (2) to (3);
\draw (2) to (4);
\draw (3) to (4);
\draw (-1,-1) node{$4_{r}'$};

\draw (4,0) node{$3_s'$};
\node[circle, draw=black] (5) at (5,0) {$1$};
\node[circle, draw=black](6) at (6.75,0) {$2$};
\node[circle, draw=black](7) at (8.5,0) {$3$};
\draw (5) to (6);
\draw (6) to (7);
\draw (5.85,0.3) node{$\lambda$};

\draw (11,0) node{$1_r$};
\draw (12,0) node[circle,draw = black]{$\emptyset$};

\draw (4,-2) node{$5_r$};
\node[circle,draw=black] (8) at (5,-2) {$1~3$};
\node[circle, draw=black] (9) at (6.75,-2) {$2$};
\node[circle, draw=black] (10) at (8.5,-2) {$1$};
\node[circle, draw=black] (11) at (10.25,-2) {$2$};
\node[circle, draw=black] (12) at (12,-2) {$3$};
\draw (8) to (9);
\draw (9) to (10);
\draw (10) to (11);
\draw (11) to (12);
\end{tikzpicture} 
\caption{$W$-graphs in type $H_3$}\label{WgraphsH3}
\end{figure}
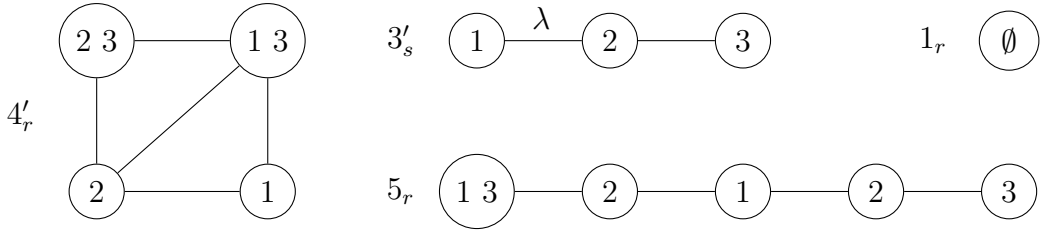

The four $W$-graphs provided in \cite{G-P} (see Figure \ref{WgraphsH3}, where $\lambda=\xi+\xi^{-1}+1$)  then determine all the irreducible representations of the Iwahori-Hecke algebra over $\F_q$ up to taking the dual $W$-graph or the algebraic conjugate by the involution $x\mapsto \overline{x}$.

\bigskip

Before determining the image of the Artin groups inside each representation, we show that we cannot always have $1\sim 2$ as in Lemma \ref{Isomorphism}. Recall that $1\sim 2$ if there exists a field automorphism $\Phi$ of $\F_p(\alpha,\xi+\xi^{-1})$ such that $\Phi(\alpha+\alpha^{-1})=\alpha+\alpha^{-1}$ and
$\Phi(\xi+\xi^{-1})=\xi^2+\xi^{-2}$.

\begin{lemme}\label{IsomorphismH3}
We have $1\sim 2$ as in Lemma \ref{Isomorphism} if and only if $\F_p(\alpha,\xi+\xi^{-1})=\F_p(\alpha+\alpha^{-1},\xi+\xi^{-1})$.

If $1\sim 2$ and $\F_p(\alpha)=\F_p(\alpha+\alpha^{-1})$ then $[\F_p(\alpha,\xi+\xi^{-1}):\F_p(\alpha)]=2$ and $\Phi_{1,2}(\alpha)=\alpha$.

If $1\sim 2$ and $\F_p(\alpha)\neq \F_p(\alpha+\alpha^{-1})$ then $[\F_p(\alpha,\xi+\xi^{-1}):\F_p(\alpha)]=1$ and $\Phi_{1,2}(\alpha)=\alpha^{-1}$.

If $1\nsim 2$ then $\F_p(\alpha,\xi+\xi^{-1})=\F_p(\alpha)\neq \F_p(\alpha+\alpha^{-1})$.
\end{lemme}

\begin{proof}

The Hasse diagram representing the inclusions between fields in this case can be found in Figure \ref{HasseH3H41}. Note that when $1\sim 2$, we have $\F_p(\alpha+\alpha^{-1},\xi+\xi^{-1})\neq \F_p(\alpha+\alpha^{-1})$. 

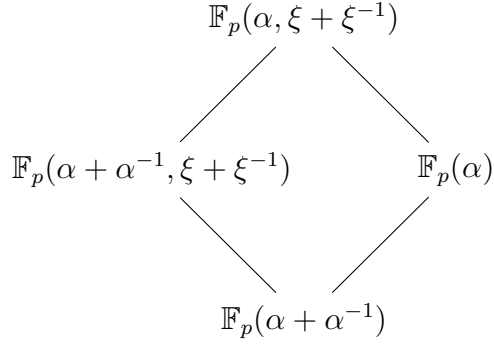
\begin{figure}
\centering
\begin{tikzpicture}
\node (1) at (0,2){$\F_p(\alpha,\xi+\xi^{-1})$};
\node (2) at (-2,0){$\F_p(\alpha+\alpha^{-1},\xi+\xi^{-1})$};
\node (3) at (2,0){$\F_p(\alpha)$};
\node (4) at (0,-2){$\F_p(\alpha+\alpha^{-1})$};
\draw (1) to (2);
\draw (1) to (3);
\draw (2) to (4);
\draw (3) to (4);
\end{tikzpicture}
\caption{Field extensions in types $H_3$ and $H_4$.1}\label{HasseH3H41}
\end{figure}

Assume by contradiction that $1\sim 2$ and $L_1=\F_p(\alpha,\xi+\xi^{-1})\neq\F_p(\alpha+\alpha^{-1},\xi+\xi^{-1})=L_2$. Since $1\sim 2$, we have $L_2=\F_p(\alpha+\alpha^{-1},\xi+\xi^{-1})\neq \F_p(\alpha+\alpha^{-1})=L_4$, because otherwise, $\Phi_{1,2}$ would stabilize $\F_p(\alpha+\alpha^{-1},\xi+\xi^{-1})$. It follows that $[L_1:L_2]=[L_2:L_4]=2$ and $[L_1:L_4]=4$. Let $L_3=\F_p(\alpha)$. We have $[L_1:L_3]\leq 2$, $[L_3:L_4]\leq 2$ and $4=[L_1:L_4]=[L_1:L_3][L_3:L_4]$. This implies that $[L_1:L_3]=[L_3:L_4]=2$. By uniqueness of the subfield of degree $2$ of $L_1$, we have $L_2=L_3$. It follows that $\alpha\in L_2=\F_p(\alpha+\alpha^{-1},\xi+\xi^{-1})$, therefore we have that $L_2=L_1$ which is a contradiction.

\medskip

Assume now that $\F_p(\alpha+\alpha^{-1},\xi+\xi^{-1})\neq \F_p(\alpha+\alpha^{-1})$. We then have that there exists an automorphism $\Phi$ of order $2$ of $\F_q$ permuting the roots of $X^2+X-1$. We have $(\xi^2+\xi^{-2})^2+\xi^2+\xi^{-2}-1=\xi+\xi^{-1}+2+\xi^2+\xi^{-2}-1=0$ so $\Phi(\xi^2+\xi^{-2})=\xi+\xi^{-1}$. This proves by definition of $\sim$ that $1\sim 2$

\medskip

Assume now $1\sim 2$ and $\F_p(\alpha)=\F_p(\alpha+\alpha^{-1})$. We have that $\Phi_{1,2}$ stabilizes $\F_p(\alpha+\alpha^{-1})$ and $\alpha\in \F_p(\alpha+\alpha^{-1}$, therefore we have that $\Phi_{1,2}(\alpha)=\alpha$.

\smallskip

Assume now $1\sim 2$ and $\F_p(\alpha)\neq \F_p(\alpha+\alpha^{-1})$. Using the notations of the first part of the proof, we here have that $[L_2:L_4]=[L_3:L_4]=2$. This implies that the unique automorphism of order $2$ of $L_3$ is $\Phi_{1,2}$ and $X^2-(\alpha+\alpha^{-1})X+1$ is an irreducible polynomial of $L_4[X]$. It follows that $\Phi_{1,2}(\alpha)=\alpha^{-1}$.

\smallskip

Assume $1\nsim 2$. We then have $\F_p(\alpha,\xi+\xi^{-1})\neq \F_p(\alpha+\alpha^{-1},\xi+\xi^{-1})$. This implies that $\alpha\notin \F_p(\alpha+\alpha^{-1},\xi+\xi^{-1})$. It follows that $\F_p(\alpha)=\F_p(\alpha,\xi+\xi^{-1})$. We then have by the Hasse diagram in Figure \ref{HasseH3H41} that $\F_p(\alpha)\neq \F_p(\alpha+\alpha^{-1})$.
\end{proof}

\bigskip

We now determine the image in each of those representations before determining the image in the full Iwahori-Hecke algebra.

\begin{prop}\label{resrepreflicos}

If $\F_r=\F_p(\alpha,\xi+\xi^{-1})=\F_p(\alpha+\alpha^{-1},\xi+\xi^{-1})$, then we have\\
 $\rho_{3_s'}(\mathcal{A}_{H_3})\simeq SL_3(r)$ and $\rho_{\overline{3_s'}}(\mathcal{A}_{H_3})=SL_3(r).$

If $\F_r=\F_p(\alpha,\xi+\xi^{-1})\neq \F_p(\alpha+\alpha^{-1},\xi+\xi^{-1})$, then we have $\rho_{3_s'}(\mathcal{A}_{H_3})\simeq SU_3(r^{\frac{1}{2}})$ and $\rho_{\overline{3_s'}}(\mathcal{A}_{H_3})\simeq SU_3(r^{\frac{1}{2}})$.

If $1\sim 2$ and $\F_p(\alpha)=\F_p(\alpha+\alpha^{-1})$, then $\Phi_{1,2}\circ \rho_{\overline{3_s'}|\mathcal{A}_{H_3}}\simeq \rho_{3_s'|\mathcal{A}_{H_3}}$.

If $1\sim 2$ and $\F_p(\alpha)\neq \F_p(\alpha+\alpha^{-1})$, then $\Phi_{1,2}\circ \rho_{\overline{3_s'}|\mathcal{A}_{H_3}}\simeq \rho_{3_s|\mathcal{A}_{H_3}}$.

\smallskip

If $\F_p(\alpha)=\F_p(\alpha+\alpha^{-1})$, then $r=q^2$.

If $\F_p(\alpha)\neq \F_p(\alpha+\alpha^{-1})$, then $r=q$.
\end{prop}

\begin{proof}

First note that by Proposition \ref{color} above and Proposition $4.1.$ of \cite{BMM}, we can see this representation as a representation over $\F_q$ even if $\F_p(\sqrt{\alpha},\xi+\xi^{-1})\neq \F_r=\F_p(\alpha,\xi+\xi^{-1})$.
\medskip

Let $\tilde{G}=<\alpha^{-1}\rho_{3_s'}(S_1),\alpha^{-1}\rho_{3_s'}(S_2),\alpha^{-1}\rho_{3_s'}(S_3)>$. Let us show that it is primitive. 
If $\tilde{G}$ was imprimitive, we could write $\F_r^n=V_1\oplus V_2 \oplus ... \oplus V_r$, where for all $i$ and for all $g\in G$, there exists a $j$ such that $g.V_i =V_j$. Since $R_1$ is irreducible, there exists $i\in [\![1,3]\!]$ such that $-R_1(S_i) .V_1 \neq V_1$. Assume there exists $i$ such that $-R_1(S_i).V_1 \neq V_1$. Up to reordering, we have $V_2=-R(S_i).V_1$. If $\dim(V_1)\geq 2$ then $H_{-R_1(S_i)}$ (the hyperplane fixed by $-R_1(S_i)$) has a non-empty intersection with $V_1$, therefore $V_1\cap V_2 \neq \emptyset$ which is a contradiction, therefore $\dim(V_1)=1$. This reasoning is valid for any $V_i$, therefore they are all one-dimensional. Let $x\in V_1$ be a non-zero vector, it can be written in a unique way as $x =x_1+x_2$ with $x_1\in \op{ker}(R_1(S_i)+\alpha)$ and $x_2\in H_{-R_1(S_i)}$. We then have that $-R_1(S_i)x = -\alpha x_1+x_2$ and $-R(S_i)(-R(S_i)x)= \alpha^2 x_1 +x_2=\alpha(x_1+x_2)+(1-\alpha)(-\alpha x_1+x_2)\in V_1\oplus V_2$. Since $\alpha\notin \{0,1\}$ this contradicts the fact that there exists $j$ such that $-R(S_i).V_2 =V_j$.
This shows that $\tilde{G}$ is primitive.

\medskip

By Wagner's theorem on groups generated by reflections (\cite{Wa} and Theorem 2.3. of \cite{MM}), since $\tilde{G}$ is primitive and is generated by pseudo-reflections, there exists $r'$ dividing $r$ such that $SL_3(r')\leq \tilde{G}\leq GL_3(r')$ or $SU_3(r'^{\frac{1}{2}})\leq \tilde{G}\leq GU_3(r'^{\frac{1}{2}})$.

We now show that $r'$ must be equal to $r$. We have $\op{det}(\alpha^{-1}\rho_{3_s'}(S_1))=-\alpha^{-1}$, therefore $\alpha\in \F_{r'}$. We also have $\tr((\alpha^{-1}\rho_{3_s'}(S_1)\alpha^{-1}\rho_{3_s'}(S_2))=1-2\alpha^{-1}+\alpha^{-1}((\xi+\xi^{-1}+1)^2)\in \F_{r'}$, therefore $(\xi+\xi^{-1}+1)^2\in \F_{r'}$. We have $(\xi+\xi^{-1}+1)^2=\xi^2+\xi^{-2}+1+2\xi+2\xi^{-1}+2=\xi+\xi^{-1}+2$, therefore it follows that $\xi+\xi^{-1}\in \F_{r'}$. Since $\F_r=\F_p(\alpha,\xi+\xi^{-1})$, we have $r'=r$.

\bigskip

Assume now $\F_r=\F_p(\alpha,\xi+\xi^{-1})=\F_p(\alpha+\alpha^{-1},\xi+\xi^{-1})$. If $\tilde{G}\leq GU_3(r^{\frac{1}{2}})$, we let $\epsilon$ be the automorphism of order $2$ of $\F_r$. We then have $\epsilon(\op{det}(\alpha^{-1}\rho_{3_s'}(S_1))=\op{det}(^t\!(\alpha^{-1}\rho_{3_s'}(S_1))^{-1})$, therefore $\epsilon(-\alpha^{-1})=-\alpha$ and $\epsilon(\alpha+\alpha^{-1})=\alpha+\alpha^{-1}$.

\medskip

We also have $\epsilon(\tr((\alpha^{-1}\rho_{3_s'}(S_1)\alpha^{-1}\rho_{3_s'}(S_2))))=\tr(^t\!(\alpha^{-1}\rho_{3_s'}(S_1)\alpha^{-1}\rho_{3_s'}(S_2))^{-1})$, therefore $\epsilon(1-2\alpha^{-1}+\alpha^{-1}(\xi+\xi^{-1}+1)^2)=1-2\alpha+\alpha(\xi+\xi^{-1}+1)^2$. This implies that $\epsilon((\xi+\xi^{-1}+1)^2)=(\xi+\xi^{-1}+1)^2$, therefore $\epsilon(\xi+\xi^{-1})=\xi+\xi^{-1}$. This would imply that $\xi+\xi^{-1}\in \F_{r^{\frac{1}{2}}}$, therefore $\F_r=\F_{r^{\frac{1}{2}}}$ which is absurd.

It follows that $SL_3(r)\leq \tilde{G}\leq GL_3(r)$, therefore $\rho_{3_s'}(\mathcal{A}_{H_3})=[\rho_{3_s'}(A_{H_3}),\rho_{3_s'}(A_{H_3})]=[\tilde{G},\tilde{G}]=SL_3(r)$.

\bigskip

Assume now $\F_r=\F_p(\alpha,\xi+\xi^{-1})\neq \F_p(\alpha+\alpha^{-1},\xi+\xi^{-1})$. The unique automorphism of order $2$ of $\F_r$ stabilizes $\F_p$ and verifies $\epsilon(\alpha)=\alpha^{-1}$ and $\epsilon(\xi+\xi^{-1})=\xi+\xi^{-1}$. By Lemma $2.4.$ of \cite{BM}, we only need to show that $\rho_{3_s'}\simeq \epsilon \circ \rho_{3_s}^\star$ in order to prove that $\rho_{3_s'}(A_{H_3})\leq GU_3(r^{\frac{1}{2}})$. Once this is shown, we will have that there exists $P\in GL_3(r)$ such that for all $i\in \{1,2,3\}$, $P\rho_{3_s'}(S_i)P^{-1}=\epsilon(^t\!\rho_{3_s'}(S_i)^{-1})$. We then have that for all $i\in\{1,2,3\}$, $P(\alpha^{-1}\rho_{3_s'}(T_s))P^{-1}=\alpha^{-1}\epsilon(^t\!\rho_{3_s'}(S_i)^{-1})=\epsilon(\alpha^{-1})^{-1}\epsilon(^t\!\rho_{3_s'}(S_i)^{-1})=\epsilon(^t\!(\alpha^{-1}\rho_{3_s'}(S_i))^{-1})$, therefore $\tilde{G}\leq GU_3(r^{\frac{1}{2}})$.

We have $\epsilon(\alpha)=\alpha^{-1}$, therefore $\epsilon(\sqrt{\alpha}^2)=\epsilon(\sqrt{\alpha})^2=\sqrt{\alpha}^{-2}$ and $(\epsilon(\sqrt{\alpha})-\sqrt{\alpha}^{-1})(\epsilon(\sqrt{\alpha})+\sqrt{\alpha}^{-1})=0$. It follows that $\epsilon(\sqrt{\alpha})\in \{\pm\sqrt{\alpha}^{-1}\}$. We will show that $\rho_{3_s'}\simeq \epsilon\circ \rho_{3_s'}^\star$ in both cases and the proof will be completed.

\medskip

Assume $\epsilon(\sqrt{\alpha})=\sqrt{\alpha}^{-1}$. Let $P=\begin{pmatrix}
 -\frac{\alpha+1}{\sqrt{\alpha}} & \lambda & 0\\
 \lambda & -\frac{\alpha+1}{\sqrt{\alpha}} & 1\\
 0 & 1 & -\frac{\alpha+1}{\sqrt{\alpha}}
 \end{pmatrix}$, we have for all $i\in \{1,2,3\}$, $P\rho_{3_s'}(S_i)P^{-1}=\epsilon(^t\!\rho_{3_s'}(S_i)^{-1})$. Note that the matrix is invertible in $\F_r$ because $\op{det}(P)=-\frac{(\alpha+1)(-\alpha\lambda^2+\alpha^2+\alpha+1)}{\alpha\sqrt{\alpha}}=0$ would imply $-\xi^2-\xi^{-2}+1=\lambda^2=\alpha+\alpha^{-1}+1$ and then $\alpha\in \{-\xi^2,-\xi^{-2}\}$ which is absurd.
 
 \medskip
 
Assume $\epsilon(\sqrt{\alpha})=-\sqrt{\alpha}^{-1}$. Let $P=\begin{pmatrix}
 1 & -\frac{\lambda\sqrt{\alpha}}{\alpha+1} & 0\\
 \frac{\lambda\sqrt{\alpha}}{\alpha+1} & -1 & \frac{\sqrt{\alpha}}{\alpha+1}\\
 0 & -\frac{\sqrt{\alpha}}{\alpha+1} & 1
 \end{pmatrix}$. It verifies the same conditions as in the previous case and the determinant is again non-zero, therefore the proof is completed for the $W$-graph $3_s'$. 
 
 \smallskip
 
 The proof is identical for $\overline{3_s'}$. Assume now $1\sim 2$.
 
 \bigskip
 
 Assume first that $\F_p(\alpha)=\F_p(\alpha+\alpha^{-1})$. By Lemma \ref{IsomorphismH3}, we have that $\Phi_{1,2}(\alpha)=\alpha$ and $\Phi_{1,2}(\lambda)=\Phi_{1,2}(\xi+\xi^{-1}+1)=\xi^2+\xi^{-2}+1=-\xi-\xi^{-1}=1-\lambda$. 
 It follows by Proposition \ref{Fieldfactorization} that there exists a character $\eta:A_{H_3}\rightarrow \F_r^\times$ such that $\Phi_{1,2}\circ \rho_{\overline{3_s'}}$ is an irreducible representation of $\mathcal{H}_{H_3,\alpha}$.
   We have $\rho_{3_s'|\mathcal{A}_{I_2(5)}}\simeq \rho_{3_s|\mathcal{A}_{I_2(5)}}\simeq \rho_1\times 1$, $\rho_{\overline{3_s'}|\mathcal{A}_{I_2(5)}}\simeq \rho_{\overline{3_s}|\mathcal{A}_{I_2(5)}}\simeq \rho_2\times 1$, where $1$ represents the trivial representation and $\Phi_{1,2}\circ \rho_2 \simeq \rho_1$.
    It follows that either $\Phi_{1,2}\circ \rho_{\overline{3_s'}|\mathcal{A}_{H_3}}\simeq \rho_{3_s'|\mathcal{A}_{H_3}}$ or $\Phi_{1,2}\circ \rho_{\overline{3_s'}|\mathcal{A}_{H_3}}\simeq \rho_{3_s|\mathcal{A}_{H_3}}$. We have $\tr(\rho_{\overline{3_s'}}(S_1S_2S_3^{-2}))=\alpha^2-\alpha+1-\lambda\alpha^{-1}$. It follows that
   $$\Phi_{1,2}(\tr(\rho_{\overline{3_s'}}(S_1S_2S_3^{-2})))=\alpha^2-\alpha+1+(\lambda-1)\alpha^{-1}.$$
   We have $\tr(\rho_{3_s'}(S_1S_2S_3^{-2}))=\alpha^2-\alpha+1+(\lambda-1)\alpha^{-1}$ and $\tr(\rho_{3_s}(S_1S_2S_3^{-2}))=(\lambda-1)\alpha+1-\alpha^{-1}+\alpha^{-2}$. Assume by contradiction that $\Phi_{1,2}\circ \rho_{\overline{3_s'}|\mathcal{A}_{H_3}}\simeq \rho_{3_s|\mathcal{A}_{H_3}}$. We then have 
  \begin{eqnarray*}
  \alpha^2-\alpha+1+(\lambda-1)\alpha^{-1} & = & (\lambda-1)\alpha+1-\alpha^{-1}+\alpha^{-2}\\
  \alpha^2-\alpha^{-2} & =& \lambda\alpha-\lambda\alpha^{-1}\\
  (\alpha-\alpha^{-1})(\alpha+\alpha^{-1}) & = & \lambda (\alpha-\alpha^{-1})\\
\alpha+\alpha^{-1} & = & -\xi^2-\xi^{-2}\\
(\alpha\xi^2+1)(\xi^{-2}+\alpha^{-1})& = & 0.  
  \end{eqnarray*}
  It then follows that $\alpha^{10}=1$ which contradicts our assumptions on the order of $\alpha$. This proves that $\Phi_{1,2}\circ \rho_{\overline{3_s'}|\mathcal{A}_{H_3}}\simeq \rho_{3_s'|\mathcal{A}_{H_3}}$.
  
  \smallskip 
  
  Assume now that $\F_p(\alpha)\neq \F_p(\alpha+\alpha^{-1})$. By Lemma \ref{IsomorphismH3}, we have that $\Phi_{1,2}(\alpha)=\alpha^{-1}$ and $\Phi_{1,2}(\lambda)=1-\lambda$. By the same arguments as in the previous case, we have that either $\Phi_{1,2}\circ \rho_{\overline{3_s'}|\mathcal{A}_{H_3}}\simeq \rho_{3_s'|\mathcal{A}_{H_3}}$ or $\Phi_{1,2}\circ \rho_{\overline{3_s'}|\mathcal{A}_{H_3}}\simeq \rho_{3_s|\mathcal{A}_{H_3}}$. Assume by contradiction that $\Phi_{1,2}\circ \rho_{\overline{3_s'}|\mathcal{A}_{H_3}}\simeq \rho_{3_s'|\mathcal{A}_{H_3}}$. We then have 
  \begin{small}
 $$ \alpha^{-2}-\alpha^{-1}+1+(\lambda-1)\alpha=\Phi_{1,2}(\tr(\rho_{\overline{3_s'}}(S_1S_2S_3^{-2})))=\tr(\rho_{3_s'}(S_1S_2S_3^{-2}))=\alpha^2-\alpha+1+(\lambda-1)\alpha^{-1}.$$
 \end{small}
 The same computation as in the previous case shows this is a contradiction and the proof is thus completed.
 
 \medskip
 
 We now show that $r=q^2$ if $\F_p(\alpha)=\F_p(\alpha+\alpha^{-1})$ and $r=q$ if $\F_p(\alpha)\neq \F_p(\alpha+\alpha^{-1})$. If $1\sim 2$, then the result follows by Lemma \ref{IsomorphismH3}. Assume now that $1\nsim 2$. Note that $\F_p(\alpha)=\F_p(\alpha+\alpha^{-1},\xi+\xi^{-1})$ then $\F_p(\alpha)=\F_p(\alpha,\xi+\xi^{-1})$ therefore the two former fields cannot be simultaneously subfields of degree $2$ of the latter. This implies that $[\F_p(\alpha,\xi+\xi^{-1}):\F_p(\alpha)]\leq 2$. By Lemma \ref{IsomorphismH3}, we have $\F_p(\alpha,\xi+\xi^{-1})\neq \F_p(\alpha+\alpha^{-1},\xi+\xi^{-1})$. The result then follows from the fact that $[\F_p(\alpha,\xi+\xi^{-1}):\F_p(\alpha)]= \dfrac{2}{[\F_p(\alpha):\F_p(\alpha+\alpha^{-1})]}$.
\end{proof}

\begin{prop}\label{H3dim4}

If $\F_q=\F_p(\sqrt{\alpha})=\F_p(\alpha)=\F_p(\alpha+\alpha^{-1})$, then we have $\rho_{4_r'}(\mathcal{A}_{H_3})=SL_4(q)$.

\smallskip

Assume that $\F_q=\F_p(\sqrt{\alpha})=\F_p(\alpha)\neq\F_p(\alpha+\alpha^{-1})$ and $\epsilon$ is the unique automorphism of order $2$ of $\F_q$.

If $\epsilon(\sqrt{\alpha})=\sqrt{\alpha}^{-1}$, then we have $\rho_{4_r'}(\mathcal{A}_{H_3})\simeq SU_4(q^{\frac{1}{2}})$.

If $\epsilon(\sqrt{\alpha})=-\sqrt{\alpha}^{-1}$, then we have $\rho_{4_r'}(\mathcal{A}_{H_3})\simeq SL_4(q^{\frac{1}{2}})$.

\smallskip

If $\F_{q'}=\F_p(\sqrt{\alpha})\neq \F_q=\F_p(\alpha)=\F_p(\alpha+\alpha^{-1})$, then we have $\rho_{4_r'}(\mathcal{A}_{H_3})=SU_4(q)$.
\end{prop}
   
\begin{proof}

Let $\beta=\xi+\xi^{-1}$ and $M=\begin{pmatrix}
1 & -\frac{\sqrt{\alpha}(3+2\beta)}{(\alpha+1)(\beta+1)} & \frac{\beta+2}{\beta+1} & -\frac{\sqrt{\alpha}(2+\beta)}{(\alpha+1)(\beta+1)}\\
-\frac{1}{\alpha+1} & \frac{1}{\sqrt{\alpha}} & -\frac{2+\beta}{(\alpha+1)(1+\beta)} & \frac{1}{\sqrt{\alpha}(\beta+1)}\\
1 & -\frac{\sqrt{\alpha}\beta^2}{\alpha+1} & -\beta & \frac{\sqrt{\alpha}\beta}{(\alpha+1)}\\
-\frac{1}{\alpha+1} & \frac{1}{\sqrt{\alpha}} & \frac{\beta}{\alpha+1} & -\frac{1}{\sqrt{\alpha}\beta}
\end{pmatrix}$, we have $\det(M)=\dfrac{5(\beta+1)\Phi_5(\alpha)}{\alpha(\alpha+1)^4(1+\beta)^2\beta}$.

Let $\rho_1$ and $\rho_2$ be the two $2$-dimensional irreducible representations of $\mathcal{H}_{I_2(5)}=<S_1,S_2>$ defined by  $\rho_1(S_1)=\begin{pmatrix}
-1 & 0\\
1 & \alpha
\end{pmatrix}$, $\rho_1(S_2)=\begin{pmatrix}
-1 & 0\\
1 & \alpha
\end{pmatrix}$, $\rho_2(S_1)=\begin{pmatrix}
\alpha & \alpha(2+\beta)\\
0 & -1
\end{pmatrix}$ and $\rho_2(S_2)=\begin{pmatrix}
\alpha & \alpha\beta^2\\
0 & -1
\end{pmatrix}$. We then have 
$$M(\rho_{4_r'}(S_1))M^{-1}=\begin{pmatrix}
\rho_1(S_1) & 0\\
0 & \rho_2(S_1)
\end{pmatrix}, M(\rho_{4_r'}(S_2))M^{-1}=\begin{pmatrix}
\rho_1(S_2) & 0\\
0 & \rho_2(S_2)
\end{pmatrix}.$$
We set $H=(\rho_1\times \rho_2)(\mathcal{A}_{I_2(5)})$ and $G=\rho_{4_r'}(\mathcal{A}_{H_3})$.

\bigskip

Assume first that $\F_q=\F_p(\sqrt{\alpha})=\F_p(\alpha)=\F_p(\alpha+\alpha^{-1})$ and $1\sim 2$. By Theorem \ref{resdihedral}, we have $\rho(\mathcal{A}_{I_2(5)})\simeq SL_2(q^2)$. More precisely, we have that $\rho_{4_r'}(\mathcal{A}_{I_2(5)})=\{\begin{pmatrix} N & 0\\
0 & \Phi_{1,2}(N) \end{pmatrix}, N\in SL_2(q^2)\}\subset  MGM^{-1}$.
Consider now the maximal subgroups of $SL_4(q)$. They are given in Tables $8.8$ and $8.9$ of \cite{BHRC}. Since our group is irreducible, we can remove the groups of class $\mathcal{C}_1$. We can also remove the groups of class $\mathcal{C}_5$ because the field generated by the traces of the elements of our group if $\F_q$. Note that $\xi+\xi^{-1}\notin \F_q$ so $q$ is an odd power of $p$ and the representation is not self-dual. It follows that $MGM^{-1}$ is included in no maximal subgroup of $\mathcal{C}_8$. The only remaining maximal subgroups are listed below with their order or a quantity their order divides.
\begin{enumerate}
\item  $(q-1)^3.\mathfrak{S}_4$, $24(q-1)^3$
\item  $SL_2(q)^2:(q-1).2$, $2q^2(q^2-1)^2(q-1)$
\item  $SL_2(q^2).(q+1).2$, $2(q+1)q^2(q^4-1)$
\item  $(4\circ 2^{1+4})^{.}  \mathfrak{S}_6$, $92160$
\item  $(4\circ 2^{1+4})^{.} \mathfrak{A}_6$, $46080$
\item  $(q-1,4) \circ 2^{.}PSL_2(7)$, $1344$  
\item  $(q-1,4) \circ 2^{.}\mathfrak{A}_7$, $20160$
\item  $(q-1,4) \circ 2^{.}PSU_4(2)$, $103680$.
\end{enumerate}

By the conditions on the order of $\alpha$, we have $q\geq 8$ and $q\neq 11$. Since $p\neq 2$ and $q$ is an odd power of $p$, we have $q\geq 13$. It follows that $\vert H\vert =q^2(q^4-1)\geq 4826640$. This implies that cases $\mathbf{4}$, $\mathbf{5}$, $\mathbf{6}$, $\mathbf{7}$ and $\mathbf{8}$ are excluded. If we were in case $\mathbf{1}$, we would have that $q^2$ divides $24$ which is absurd since $p\neq 2$. In case $\mathbf{2}$, we would have that $q^2+1$ divides $2(q^2-1)(q-1)=2(q^2+1)(q-1)-4(q-1)$, therefore we would have that $q^2+1$ divides $4(q-1)<q^2+1$ since $q\geq 13$.

\smallskip

The only remaining case is case $\mathbf{3}$. By \cite{MR}, $\mathcal{A}_{H_3}$ is perfect, it would then follow that $G\simeq SL_2(q^2)$. This would imply that $G=\rho_{4_r'}(\mathcal{A}_{I_2(5)})$. The coefficient on the first row and third column of $M\rho_{4_r'}(S_1S_3^{-1})M^{-1}$ is equal to $\frac{1}{5}(-\alpha-3\sqrt{\alpha}-1)\beta-(3\alpha-4\sqrt{\alpha}-3)$. 

If this is equal to zero and $-\alpha-3\sqrt{\alpha}-1\neq 0$, then $\beta=-\frac{3\alpha+4\sqrt{\alpha}+3}{\alpha+3\sqrt{\alpha}+1}\in \F_p$ which is false.

 If the coefficient is equal to zero and $-\alpha-3\sqrt{\alpha}-1=0$ then $\alpha+1=-3\sqrt{\alpha}$, therefore the coefficient is equal to $5\sqrt{\alpha}\neq 0$ which is absurd. This proves that $M\rho_{4_r'}(S_1S_3^{-1})M^{-1}\notin \rho_{4_r'}(\mathcal{A}_{I_2(5)})$. This proves that if $\F_q=\F_p(\sqrt{\alpha})=\F_p(\alpha)=\F_p(\alpha+\alpha^{-1})$ and $1\sim 2$ then $G=SL_4(q)$.

\medskip

Assume now that $\F_q=\F_p(\sqrt{\alpha})=\F_p(\alpha)=\F_p(\alpha+\alpha^{-1})$ and $1\nsim 2$. By Theorem \ref{resdihedral}, we have $\rho(\mathcal{A}_{I_2(5)})\simeq SL_2(q')\times SL_2(q')$, where $\F_{q'}=\F_p(\alpha,\xi+\xi^{-1})=\F_p(\alpha+\alpha^{-1},\xi+\xi^{-1})$. By Lemma \ref{IsomorphismH3}, we have that $\F_q=\F_q'$. By Lemma \ref{normclosdih}, $G$ is generated by transvections. We can then apply Theorem \ref{transvections} to get that $G\in \{SL_4(\tilde{q}),SU_4(\tilde{q}^\frac{1}{2}),SP_4(\tilde{q})\}$ for some $\tilde{q}$ dividing $q$. We have $\tilde{q}=q$ because $G$ contains a natural $SL_2(q)$. We have that $G$ cannot be preserved by any non-degenerate bilinear form because $\rho_{4_r'}$ is not self-dual. Assume by contradiction that $G\simeq SU_4(q^{\frac{1}{2}})$. There exists an automorphism $\epsilon$ of order $2$ of $\F_q$ such that $\epsilon \circ \rho_{4_r'}^\star\simeq \rho_{4_r'}$. By Lemma \ref{abel}, we can apply Proposition \ref{Fieldfactorization} to $\epsilon$. It follows that we have that $\epsilon(\alpha)\in \{\alpha,\alpha^{-1}\}$, therefore $\epsilon(\alpha+\alpha^{-1})=\alpha+\alpha^{-1}$. This implies that $\epsilon$ is of order $1$ since $\F_q=\F_p(\alpha+\alpha^{-1})$. This is a contradiction, therefore $G=SL_4(q)$.

\bigskip

Assume now $\F_q=\F_p(\sqrt{\alpha})=\F_p(\alpha)\neq\F_p(\alpha+\alpha^{-1})$ and $1\sim 2$. We have $\Phi_{1,2}(\alpha)=\alpha^{-1}$, therefore $\Phi_{1,2}(\sqrt{\alpha})\in \{\pm\sqrt{\alpha}^{-1}\}$. By Theorem \ref{resdihedral}, we have $H\simeq SL_2(q)$.

\medskip

 Assume first $\Phi_{1,2}(\sqrt{\alpha})=\sqrt{\alpha}^{-1}$. Let $R=\begin{pmatrix}
\frac{\alpha-\sqrt{\alpha}+1}{\sqrt{\alpha}} & 1 & 1 & \frac{\sqrt{\alpha}}{\alpha+1}\\
1 & \frac{\alpha+1}{\sqrt{\alpha}}& \frac{\alpha+\sqrt{\alpha}+1}{\alpha+1} & 1\\
1 & \frac{\alpha+\sqrt{\alpha}+1}{\alpha+1} & \frac{\alpha+1}{\sqrt{\alpha}} & 1\\
\frac{\sqrt{\alpha}}{\alpha+1} & 1 & 1  & \frac{\alpha-\sqrt{\alpha}+1}{\sqrt{\alpha}} 
\end{pmatrix}$. We have $\op{det}(R)=\frac{\Phi_5(\sqrt{\alpha})\Phi_{10}(\sqrt{\alpha})^3}{\alpha^2(\alpha+1)^4}\neq 0$ and
 $$(R^{-1}\rho_{4_r'}(S_1)R,R^{-1}\rho_{4_r'}(S_2)R,R^{-1}\rho_{4_r'}(S_3)R)=(\epsilon(^t\rho_{4_r'}(S_1)^{-1}), \epsilon(^t\rho_{4_r'}(S_2)^{-1}),\epsilon(^t\rho_{4_r'}(S_3)^{-1})).$$
 It follows by Lemma \ref{Harinordoquy} that up to conjugation in $GL_4(q)$, we have $G\leq SU_4(q^{\frac{1}{2}})$. We know $G$ is irreducible so it is included in no subgroup of class $\mathcal{C}_1$. We list below the remaining maximal subgroups of $SU_4(q^{\frac{1}{2}})$ with their order or a quantity their order divides. The tables are obtained using Table $8.10$ and $8.11$ of \cite{BHRC}.
 
 \begin{enumerate}
 \item $(q^{\frac{1}{2}}+1)^3.\mathfrak{S}_4$, $24(q^{\frac{1}{2}}+1)^3$
 \item $SU_2(q^{\frac{1}{2}})^2:(q^{\frac{1}{2}}+1).2$, $2q(q^{\frac{1}{2}}+1)(q-1)^2$
 \item $SL_2(q).(q^{\frac{1}{2}}-1).2$, $2(q^{\frac{1}{2}}-1)\vert SL_2(q)\vert$
 \item $SU_4(q_0)$, $q^{\frac{1}{2}}=q_0^r$, $r$ odd prime, $q_0^6(q_0^2-1)(q_0^3+1)(q_0^4-1)$
 \item $SO_4^+(q^{\frac{1}{2}}).[(q^{\frac{1}{2}}+1,4)]$, $4q(q-1)^2$
  \item $SO_4^-(q^{\frac{1}{2}}).[(q^{\frac{1}{2}}+1,4)]$, $4q(q^2-1)$
  \item $(4\circ 2^{1+4})^{.} S_6$, $92160$
  \item $(4\circ 2^{1+4})^{.} S_6$, $46080$
  \item $(q^{\frac{1}{2}}+1,4)\circ 2^{.}PSL_2(7)$, $1344$
  \item $(q^{\frac{1}{2}}+1,4)\circ 2^{.}\mathfrak{A}_7$, $20160$
  \item $4_2^{.}PSL_3(4)$, $5040$
  \item $(q^{\frac{1}{2}}+1,4)\circ 2^{.}PSU_4(2)$, $207360$
 \end{enumerate}
 We have $\vert H\vert =q(q^2-1)$. We have $\Phi_{1,2}(\alpha)=\alpha^{-1}$ so $\alpha^{q^{\frac{1}{2}}+1}=1$. This implies that $q^{\frac{1}{2}}+1\geq 7$, therefore we have that $q\geq 49$. It follows that $\vert H\vert \geq 117600$. This excludes cases $\mathbf{7}$ to $\mathbf{11}$. We have that $q$ is a square and $p$ is odd, therefore $q\geq 49$ implies that $q=49$ or $q\geq 81$. If $q=49$ then $\vert H\vert =117600$ does not divide $207360$, therefore case $\mathbf{12}$ is excluded. If $q\geq 81$ then $\vert H\vert \geq 531360$. Therefore, case $\mathbf{12}$ is excluded for any $q$.
 
 \smallskip
 
Assume by contradiction that we are in case $\mathbf{1}$. We then have that $q$ divides $24$ which is a contradiction. Case $\mathbf{1}$ is therefore excluded.

\smallskip

Assume by contradiction that we are in case $\mathbf{2}$. We then have that $q+1$ divides $(q^{\frac{1}{2}}+1)(q-1)$. It follows that $q+1$ divides $(q^{\frac{1}{2}}+1)(q+1)-(q^{\frac{1}{2}}+1)(q-1)=2q^{\frac{1}{2}}<q+1$. This contradiction proves that case $\mathbf{2}$ is excluded.

\smallskip

Assume by contradiction that we are in case $\mathbf{4}$. We then have that $q$ divides $q_0^{6}=q^{\frac{3}{r}}$. This implies that $r=3$, therefore $(q^2-1)$ divides $(q^{\frac{2}{3}}-1)(q+1)(q^{\frac{4}{3}}-1)$. It follows that $q-1$ divides $(q^{\frac{2}{3}}-1)(q^{\frac{4}{3}}-1)=q^2-q^{\frac{4}{3}}-q^{\frac{2}{3}}+1=q^2-1-q^{\frac{1}{3}}(q-1)-q^{\frac{2}{3}}-q^{\frac{1}{3}}+2$. This implies that $q-1$ divides $q^{\frac{2}{3}}+q^{\frac{1}{3}}-2<q-1$. This contradiction proves that case $\mathbf{4}$ is also excluded.

\smallskip

Assume now by contradiction that we are in case $\mathbf{3}$. Since $G$ is perfect, we have that $G\leq SL_2(q)$. It follows that $G=H$ but this is absurd by the computations made when $\F_q=\F_p(\sqrt{\alpha})=\F_p(\alpha)=\F_p(\alpha+\alpha^{-1})$ and $1\sim 2$.

\smallskip

Assume by contradiction that we are in case $\mathbf{5}$ or $\mathbf{6}$. By \cite{MR}, we have that $\mathcal{A}_{H_3}$ is perfect. It follows that we would have $G\leq SO_4^{\pm}(q^{\frac{1}{2}})$. This contradicts the fact that $\rho_{4_r'}$ is not self-dual.

\smallskip

It follows that $G$ is included in no maximal subgroup of $SU_4(q^{\frac{1}{2}})$. This implies that $G\simeq SU_4(q^{\frac{1}{2}})$.

\medskip

Assume now $\Phi_{1,2}(\sqrt{\alpha})=-\sqrt{\alpha}^{-1}$. Let $R'=\begin{pmatrix}\frac{\sqrt{\alpha}}{\alpha+1} &  1 & 1 &\frac{\alpha-\sqrt{\alpha}+1}{\sqrt{\alpha}} \\
1 & \frac{\alpha+\sqrt{\alpha}+1}{\alpha+1} & \frac{\alpha+1}{\sqrt{\alpha}} & 1\\
 1 & \frac{\alpha+1}{\sqrt{\alpha}} & \frac{\alpha+\sqrt{\alpha}+1}{\alpha+1}  & 1\\
 \frac{\alpha-\sqrt{\alpha}+1}{\sqrt{\alpha}}  & 1 & 1  & \frac{\sqrt{\alpha}}{\alpha+1}
 
 \end{pmatrix}$. We have $\op{det}(R')=\op{det}(R)\neq 0$ and
 $$(R'^{-1}\rho_{4_r'}(T)R',R'^{-1}\rho_{4_r'}(S_1)R',R'^{-1}\rho_{4_r'}(S_2)R')=(-\alpha\epsilon(\rho_{4_r'}(T)), -\alpha\epsilon(\rho_{4_r'}(S_1)),-\alpha\epsilon(\rho_{4_r'}(S_2))).$$
 It follows by Lemma \ref{Harinordoquy} that up to conjugation in $GL_4(q)$, we have $G\leq SL_4(q^{\frac{1}{2}})$. We have here $q\geq 49$ so $\vert H\vert\geq  117600$. We can then apply the same reasoning as for $\F_q=\F_p(\sqrt{\alpha})=\F_p(\alpha)=\F_p(\alpha+\alpha^{-1})$ and $1\sim 2$ to prove that $G\simeq SL_4(q^{\frac{1}{2}})$.

\bigskip

Assume $\F_q=\F_p(\sqrt{\alpha})=\F_p(\alpha)\neq\F_p(\alpha+\alpha^{-1})$ and $1\nsim 2$. Let $\epsilon$ be the unique automorphism of order $2$ of $\F_q$. By the same arguments as in the previous case, if $\epsilon(\sqrt{\alpha})=\sqrt{\alpha}^{-1}$ then we have that up to conjugation, $G\leq SU_4(q^{\frac{1}{2}})$. If $\epsilon(\sqrt{\alpha})=-\sqrt{\alpha}^{-1}$, then we have that, up to conjugation, $G\leq SL_4(q^{\frac{1}{2}})$. We have $H\simeq SU_2(q^{\frac{1}{2}})\times SU_2(q^{\frac{1}{2}})$. By Lemma \ref{normclosdih}, we have that $G$ is normally generated by $H$. It follows that $G$ is an irreducible subgroup of $GL_4(q)$ generated by transvections.

\smallskip

Assume $\epsilon(\sqrt{\alpha})=\sqrt{\alpha}^{-1}$. By Theorem \ref{transvections}, we have that $G$ is conjugate in $GL_4(q)$ to $SU_4(q'^{\frac{1}{2}})$, $SP_4(q')$ or $SL_4(q')$ for some $q'$ dividing $q$. We know that $G$ contains a natural $SU_2(q^{\frac{1}{2}})$. By Lemma \ref{field}, we have that $q^{\frac{1}{2}}$ divides $q$. $\rho_{4_r'}$ is not self-dual so the symplectic case is excluded. The groups $SL_4(q)$ and $SL_4(q^{\frac{1}{2}})$ are not included in $SU_4(q^{\frac{1}{2}})$, therefore we have $G\simeq SU_4(q^{\frac{1}{4}})$ or $G\simeq SU_4(q^{\frac{1}{2}})$. The natural $SU_2(q^{\frac{1}{2}})$ cannot be included in $SU_4(q^{\frac{1}{4}})$, therefore we have that $G\simeq SU_4(q^{\frac{1}{2}})$.

\smallskip

Assume now that $\epsilon(\sqrt{\alpha})=-\sqrt{\alpha}^{-1}$. By Theorem \ref{transvections} and the facts that $G$ contains a natural $SU_2(q^{\frac{1}{2}})$ and $\rho_{4_r'}$ not self-dual, we have that $G$ is conjugate in $GL_4(q)$ to $SU_4(q^{\frac{1}{2}})$, $SL_4(q^{\frac{1}{2}})$ or $SL_4(q)$. $G$ is conjugate to a subgroup of $SL_4(q^{\frac{1}{2}})$ therefore we have that $G\simeq SL_4(q^{\frac{1}{2}})$. 

\bigskip

Assume now that $\F_{q^2}=\F_p(\sqrt{\alpha})\neq \F_p(\alpha)=\F_p(\alpha+\alpha^{-1})$ and $1\sim 2$. By Lemma \ref{IsomorphismH3}, we have that the unique automorphism of order $2$ of $\F_{q^2}$ is $\Phi_{1,2}$ and $\Phi_{1,2}(\alpha)=\alpha$. The polynomial $X^2-\alpha$ is irreducible in $\F_q[X]$, therefore we have that $\Phi_{1,2}(\sqrt{\alpha})=-\sqrt{\alpha}$. Let $Q=E_{1,4}+E_{2,3}+E_{3,2}+E_{4,1}$, we have that
\begin{footnotesize}
$$(Q\rho_{4_r'}(S_1)Q^{-1},Q\rho_{4_r'}(S_2)Q^{-1},Q\rho_{4_r'}(S_3)Q^{-1})=(-\alpha \epsilon(^t(\rho_{4_r'}(S_1)^{-1})),-\alpha \epsilon(^t(\rho_{4_r'}(S_2)^{-1})),-\alpha \epsilon(^t(\rho_{4_r'}(S_3)^{-1})).$$
\end{footnotesize}
It follows by Lemma \ref{Harinordoquy} that up to conjugation in $GL_4(q^2)$, we have that $G\leq SU_4(q)$. We have $H\simeq SL_2(q^2)$. We have here $\Phi_{1,2}(\alpha)=\alpha$, therefore $\alpha^{q-1}=1$. This implies that $q\geq 8$ and $q\neq 11$. We also have that $q$ is not a square since $\xi+\xi^{-1}\notin \F_q$. It follows that $q\geq 13$ and we can apply the same reasonning as when $\F_q=\F_p(\sqrt{\alpha})=\F_p(\alpha)\neq\F_p(\alpha+\alpha^{-1})$, $1\sim 2$ and $\Phi_{1,2}(\sqrt{\alpha})=\sqrt{\alpha}^{-1}$ to conclude that $G\simeq SU_4(q)$.

\smallskip

Assume now that $\F_{q^2}=\F_p(\sqrt{\alpha})\neq \F_p(\alpha)=\F_p(\alpha+\alpha^{-1})$ and $1\nsim 2$. The unique automorphism $\epsilon$ of order $2$ of $\F_q$ still verifies $\epsilon(\sqrt{\alpha})=-\sqrt{\alpha}$, therefore we have that $G\leq SU_4(q)$ up to conjugation in $GL_4(q^2)$. By Lemma \ref{IsomorphismH3} and Theorem \ref{resdihedral}, we have that $H\simeq SU_2(q)\times SU_2(q)$. By Lemma \ref{normclosiscos}, we have that $\mathcal{A}_{H_3}$ is normally generated by $\mathcal{A}_{I_2(5)}$. It follows that $\rho_{4_r'}(\mathcal{A}_{H_3})$ is generated by transvections. The group $G$ contains a natural $SU_2(q)$, therefore by Lemma \ref{field}, we have that the field generated by the traces of the elements of $G$ contains $\F_q^{\frac{1}{2}}$. We also have that $G$ preserves no non-degenerate bilinear form because $\rho_{4_r'}$ is not self-dual. By Theorem \ref{transvections}, we have that $G$ is conjugate in $GL_4(q^2)$ to $SL_4(q)$, $SU_4(q)$ or $SL_4(q^2)$. We can then conclude that $G\simeq SU_4(q)$ conclude from the fact that $G\leq SU_4(q)$.
\end{proof}

 \begin{prop}\label{H3dim5}
 If $\F_q=\F_p(\alpha)=\F_p(\alpha+\alpha^{-1})$, then we have $\rho_{5_r}(\mathcal{A}_{H_3})\simeq SL_5(q)$.
 If $\F_q=\F_p(\alpha)\neq \F_p(\alpha+\alpha^{-1})$, then we have $\rho_{5_r}(\mathcal{A}_{H_3})\simeq SU_5(q^{\frac{1}{2}})$
 \end{prop}
 
 \begin{proof}
 Let $G=\rho_{5_r}(\mathcal{A}_{H_3})$ and $H=\rho_{5_r}(\mathcal{A}_{I_2(5)})$. First note that by Proposition \ref{color}, we can assume $G\leq SL_5(q)$. Let us now consider the restriction to $\mathcal{A}_{I_2(5)}$. Let $M$ be the following matrix
$$\begin{pmatrix}
\sqrt{\alpha}(\alpha+1)(\beta+1)& -\alpha(\beta+2)(\beta+1) & \sqrt{\alpha}(\alpha+1)(\beta+2) & -\alpha(\beta+2) & 0\\
-\sqrt{\alpha}(\beta+1) & (\alpha+1)(\beta+1) & -\sqrt{\alpha}(\beta+2) & \alpha+1 & -\sqrt{\alpha}\\
\beta(\alpha+1) & -\sqrt{\alpha}\beta^3 & \beta(\beta^2-1)(\alpha+1) & -\sqrt{\alpha}\beta(\beta^2-1) & 0\\
-\beta & \beta\frac{\alpha+1}{\sqrt{\alpha}} & -\beta(\beta^2-1) & \frac{(\beta^2-1)(\alpha+1)}{\sqrt{\alpha}\beta} & -\frac{\beta^2-1}{\beta}\\
0 & 0 &0 & 0 & 1
\end{pmatrix}.$$
We then have $\op{det}(M)=5\Phi_5(\alpha)\neq 0$ and for all $h\in \mathcal{H}_{I_2(5),\alpha}$, 
$$M\rho_{5_s}(h)M^{-1}=\begin{pmatrix}
\rho_1(h) & 0 & 0\\
 0 & \rho_2(h) & 0\\
 0 & 0 & \op{Ind}(h)
\end{pmatrix}.$$
The representations above are given in Theorem \ref{melyssaestlameilleure}.

\bigskip
 
 Assume now $\F_q=\F_p(\alpha)=\F_p(\alpha+\alpha^{-1})$ and $1\sim2$. We then have $q$ is an odd power of $p$. By Lemma \ref{IsomorphismH3} and Theorem \ref{resdihedral}, we have that $H\simeq SL_2(q^2)$.
 \smallskip
We now consider the maximal subgroups of $SL_5(q)$ given in Tables $8.18$ and $8.19$ of \cite{BHRC}. We know that $G$ is irreducible, therefore it cannot be contained in any group of class $\mathcal{C}_1$. Since $G$ contains $H$, we have that the field generated by the elements of $G$ contains $\F_q$, therefore we have that $G$ cannot be included in any group of class $\mathcal{C}_5$. The group $G$ contains a $SL_2(q^2)$ twisted by the field automorphism of degree $2$ of $\F_{q^2}$, therefore it cannot be a subgroup of $SU_5(q^{\frac{1}{2}})$. We know that $\rho_{5r}$ is not self-dual, therefore $G$ cannot preserve any non-degenerate bilinear form. It follows that $G$ cannot be included in any group of class $\mathcal{C}_8$. We list below the remaining maximal subgroups remaining and their order or a quantity their order divides.
 
 \begin{enumerate}
 \item  $(q-1)^4 : \mathfrak{S}_5$, $120(q-1)^4$
 \item  $\frac{q^5-1}{q-1} : 5$, $5(q^4+q^3+q^2+q+1)$
   \item  $5_+^{1+2} : Sp_2(5)$, $15000$
      \item $(q-1,5)\times PSL_2(11)$, $3300$
        \item $M_{11}$, $7920$
         \item  $(q-1,5)\times PSU_4(2)$, $129600$
    \end{enumerate}
    
  We have $\alpha^{q-1}=1$ so $q-1\geq 7$. It follows that $q\geq 9$ and $\vert H\vert \geq 
  531360$. This excludes cases $\mathbf{3}$ to $\mathbf{6}$. 
  
  Assume we are in case $\mathbf{1}$ or $\mathbf{2}$. We have $\vert H\vert =q^2(q^4-1)$. It follows that $q^2$ divides $120$. This is absurd because $p\neq 2$.

\smallskip  
  
  It follows that $G\simeq SL_5(q)$.
  
  \medskip
  
  Assume now $\F_q=\F_p(\alpha)=\F_p(\alpha+\alpha^{-1})$ and $1\nsim 2$. We then have that  $\rho_{5_s}(\mathcal{A}_{I_2(5)})\simeq SU_2(q)\times SU_2(q)$. The group $G$ is then irreducible and generated by transvections. Therefore, by Theorem \ref{transvections},
  we have  that $G\in \{SL_5(q')Sp_5(q'),SU_5(q'^{\frac{1}{2}}),q'|q\}$. Since $G$ contains a natural $SL_2(q)$, we have that $q'=q$. We have $G\neq SU_5(q^{\frac{1}{2}})$ because $G$ contains a natural $SU_2(q)$ and $G\neq Sp_5(q)$ because $\rho_{5_s} \not\simeq \rho_{5_s}^{\star}$. This proves that $G\simeq SL_5(q)$.
  
  \bigskip
  
  Assume now that $\F_q=\F_p(\alpha)\neq\F_p(\alpha+\alpha^{-1})$ and $1\sim2$. By Lemma \ref{IsomorphismH3} and Theorem \ref{resdihedral}, we have that $H\simeq SL_2(q)$ and $\Phi_{1,2}$ is the unique automorphism of order $2$ of $\F_q$. We have $\Phi_{1,2}(\alpha)=\alpha^{-1}$, therefore by Proposition \ref{Fieldfactorization}, we have that $\Phi_{1,2} \circ \rho_{5_r|\mathcal{A}_{H_3}}\simeq \rho_{5_r|\mathcal{A}_{H_3}}$ or $\Phi_{1,2} \circ \rho_{5_r|\mathcal{A}_{H_3}}\simeq \rho_{5_r'|\mathcal{A}_{H_3}}$. We have 
\begin{eqnarray*}
\tr(\rho_{5_r}(S_1S_2S_3^{-2})) & = & \alpha^2-2\alpha+2-\alpha^{-1}\\
\Phi_{1,2}(\tr(\rho_{5_r}(S_1S_2S_3^{-2}))) & = & \alpha^{-2}-2\alpha^{-1}+2-\alpha\\
\tr(\rho_{5_r'}(S_1S_2S_3^{-2})) & = & \alpha^{-2}-2\alpha^{-1}+2-\alpha
\end{eqnarray*}
If $\alpha^2-2\alpha+2-\alpha^{-1}= \alpha^{-2}-2\alpha^{-1}+2-\alpha$, then 
$$0=\alpha^2-\alpha^{-2}-(\alpha-\alpha^{-1})=(\alpha-\alpha^{-1})(\alpha+\alpha^{-1}-1)=\alpha^{-2}(\alpha^2-1)\Phi_6(\alpha).$$
This proves that $\Phi_{1,2} \circ \rho_{5_r|\mathcal{A}_{H_3}}\simeq \rho_{5_r'|\mathcal{A}_{H_3}}$, therefore we have that  $\Phi_{1,2} \circ \rho_{5_r|\mathcal{A}_{H_3}}^\star\simeq \rho_{5_r|\mathcal{A}_{H_3}}$. By Lemma \ref{Harinordoquy}, we have that $G$ is conjugate in $GL_5(q)$ to a subgroup of $SU_5(q^{\frac{1}{2}})$.

We now consider the maximal subgroups of $SU_5(q^{\frac{1}{2}})$ given in Tables $8.20$ and $8.21$ of \cite{BHRC}. $G$ is irreducible so it cannot be included in a maximal subgroup of class $\mathcal{C}_1$. It contains a diagonal $SL_2(q)$ twisted by $\epsilon$, therefore it cannot be included in a maximal subgroup of class $\mathcal{C}_5$. We list below the remaining maximal subgroups with their order or a quantity their order divides.

\begin{enumerate}
\item $(q^{\frac{1}{2}}+1)^4:\mathfrak{S}_5$, $120(q^{\frac{1}{2}}+1)^4$
\item $\frac{q^5+1}{q+1}:5$, $5(q^4-q^3+q^2-q+1)$
\item $5_+^{1+2}:SP_2(5)$, $15000$
\item $(q^{\frac{1}{2}}+1,5)\times PSL_2(11)$, $3300$
\item $(q^{\frac{1}{2}}+1,5)\times PSU_4(2)$, $129600$
\end{enumerate}

We have $\vert H\vert =q(q^2-1)$. We have $\alpha^{q^{\frac{1}{2}}+1}=1$, therefore $q^{\frac{1}{2}}\geq 6$. This implies that $q^{\frac{1}{2}}\geq 7$ and $q\geq 49$. This implies that $q=49$ or $q\geq 81$ since $q$ is as square and $p\neq 2$. We then have that $\vert H\vert =117600$ or $\vert H\vert \geq 531360$. Cases $\mathbf{3}$, $\mathbf{4}$ and $\mathbf{5}$ are therefore excluded. We have that $q$ is a square and $p\neq 2$, therefore $q$ does not divide $120$ and cases $\mathbf{1}$ and $\mathbf{2}$ are excluded.

It follows that $G\simeq SU_5(q^{\frac{1}{2}})$.
  
 \medskip
 
  Assume now that $\F_q=\F_p(\alpha)\neq\F_p(\alpha+\alpha^{-1})$ and $1\nsim2$. We then have by the same arguments as in the previous case that $G$ is conjugate in $GL_5(q)$ to a subgroup of $SU_5(q^{\frac{1}{2}})$. We have $H\simeq SU_2(q^{\frac{1}{2}})\times SU_2(q^{\frac{1}{2}})$. By Lemma \ref{normclosdih}, we have that $G$ is normally generated by $H$. This implies that $G$ is an irreducible group generated by transvections. We also have that $G$ contains a natural $SU_2(q^{\frac{1}{2}})$ and $\rho_{5_r}$ is not self-dual. It follows by Theorem \ref{transvections} that $G$ is conjugate in $GL_5(q)$ to $SU_5(q^{\frac{1}{2}})$, $SL_5(q^{\frac{1}{2}})$ or $SL_5(q)$. $G$ is conjugate to a subgroup of $SU_5(q^{\frac{1}{2}})$ so $G\simeq SU_5(q^{\frac{1}{2}})$ and the proof is concluded.
 \end{proof}
 
 \begin{theo}\label{resicos}
 Under the assumptions on $\alpha$ and $p$, we have the following results.
\begin{enumerate}
\item Assume $1\sim 2$.
\begin{enumerate}
\item If $\F_q=\F_p(\sqrt{\alpha})=\F_p(\alpha)=\F_p(\alpha+\alpha^{-1})$, then the morphism from $\mathcal{A}_{H_3}$ to $\mathcal{H}_{H_3,\alpha}^\star \simeq GL_1(q)^2 \times GL_3(q)^2\times GL_4(q)^2\times GL_5(q)$ factorizes through the surjective morphism 
$$\Phi : \mathcal{A}_{H_3} \rightarrow SL_3(q^2)\times SL_4(q)\times SL_5(q).$$
\item If $\F_q=\F_p(\sqrt{\alpha})=\F_p(\alpha)\neq \F_p(\alpha+\alpha^{-1})$ and $\Phi_{1,2}(\sqrt{\alpha})=\sqrt{\alpha}^{-1}$, then the morphism from $\mathcal{A}_{H_3}$ to $\mathcal{H}_{H_3,\alpha}^\star \simeq GL_1(q)^2 \times GL_3(q)^2\times GL_4(q)^2\times GL_5(q)$ factorizes through the surjective morphism 
$$\Phi : \mathcal{A}_{H_3} \rightarrow SL_3(q)\times SU_4(q^{\frac{1}{2}})\times SU_5(q^{\frac{1}{2}}).$$
\item If $\F_q=\F_p(\sqrt{\alpha})=\F_p(\alpha)\neq \F_p(\alpha+\alpha^{-1})$ and $\Phi_{1,2}(\sqrt{\alpha})=-\sqrt{\alpha}^{-1}$, then the morphism from $\mathcal{A}_{H_3}$ to $\mathcal{H}_{H_3,\alpha}^\star \simeq GL_1(q)^2 \times GL_3(q)^2\times GL_4(q)^2\times GL_5(q)$ factorizes through the surjective morphism 
$$\Phi : \mathcal{A}_{H_3} \rightarrow SL_3(q)\times SL_4(q^{\frac{1}{2}})\times SU_5(q^{\frac{1}{2}}).$$
\item If $\F_{q^2}=\F_p(\sqrt{\alpha})\neq \F_p(\alpha)=\F_p(\alpha+\alpha^{-1})$, then the morphism from $\mathcal{A}_{H_3}$ to $\mathcal{H}_{H_3,\alpha}^\star \simeq GL_1(q)^2 \times GL_3(q)^2\times GL_4(q)^2\times GL_5(q)$ factorizes through the surjective morphism 
$$\Phi : \mathcal{A}_{H_3} \rightarrow SL_3(q^2)\times SU_4(q)\times SL_5(q).$$
\end{enumerate}
\item Assume $1\nsim 2$. When it exists, we write $\epsilon$ the automorphism of order $2$ of $\F_q$.
\begin{enumerate}
\item If $\F_q=\F_p(\sqrt{\alpha})=\F_p(\alpha)\neq \F_p(\alpha+\alpha^{-1})$ and $\epsilon(\sqrt{\alpha})=\sqrt{\alpha}^{-1}$, then the morphism from $\mathcal{A}_{H_3}$ to $\mathcal{H}_{H_3,\alpha}^\star \simeq GL_1(q)^2 \times GL_3(q)^2\times GL_4(q)^2\times GL_5(q)$ factorizes through the surjective morphism 
$$\Phi : \mathcal{A}_{H_3} \rightarrow SU_3(q^{\frac{1}{2}})^2\times SU_4(q^{\frac{1}{2}})\times SU_5(q^{\frac{1}{2}}).$$
\item If $\F_q=\F_p(\sqrt{\alpha})=\F_p(\alpha)\neq \F_p(\alpha+\alpha^{-1})$ and $\epsilon(\sqrt{\alpha})=-\sqrt{\alpha}^{-1}$, then the morphism from $\mathcal{A}_{H_3}$ to $\mathcal{H}_{H_3,\alpha}^\star \simeq GL_1(q)^2 \times GL_3(q)^2\times GL_4(q)^2\times GL_5(q)$ factorizes through the surjective morphism 
$$\Phi : \mathcal{A}_{H_3} \rightarrow SU_3(q^{\frac{1}{2}})^2\times SL_4(q^{\frac{1}{2}})\times SU_5(q^{\frac{1}{2}}).$$
\end{enumerate}
\end{enumerate} 
 
\end{theo}

\begin{proof}
Assume first that $1\sim 2$. We then have the result by Lemmas \ref{Goursat} and \ref{IsomorphismH3} and Propositions \ref{resrepreflicos}, \ref{H3dim4} and \ref{H3dim5}.

\smallskip

Assume now that $1\nsim 2$. By \cite{MR}, $\mathcal{A}_{H_3}$ is perfect, it follows that by Lemmas \ref{Goursat} and \ref{IsomorphismH3} and Propositions \ref{resrepreflicos}, \ref{H3dim4} and \ref{H3dim5}, we only need to prove that in all cases there exists no non-trivial field automorphism $\varphi$ of $\F_p(\alpha)$ such that $\varphi \circ \rho_{\overline{3_s}|\mathcal{A}_{H_3}}\simeq \rho_{3_s|\mathcal{A}_{H_3}}$ or $\varphi \circ \rho_{\overline{3_s}|\mathcal{A}_{H_3}}\simeq \rho_{3_s'|\mathcal{A}_{H_3}}$ to conclude the proof. Assume there exists such an automorphism $\varphi$. By Proposition \ref{Fieldfactorization} and Lemma \ref{abel}, we have that $\varphi(\alpha)\in \{\alpha,\alpha^{-1}\}$. It follows that $\varphi(\alpha+\alpha^{-1})=\alpha+\alpha^{-1}$.  If $\F_p(\alpha)=\F_p(\alpha+\alpha^{-1})$, this proves that $\varphi$ is trivial. If $\F_p(\alpha)\neq \F_p(\alpha+\alpha^{-1})$, then $\varphi=\epsilon$ and $\varphi \circ \rho_{\overline{3_s}|\mathcal{A}_{H_3}}\simeq \rho_{\overline{3_s'}|\mathcal{A}_{H_3}}$. The result thus follows from Proposition \ref{resH3derivedsubgroup}.
\end{proof}
\section{Type $H_4$, general aspects}\label{H4gen}

Let $p$ be a prime number, $p\notin \{2,3,5\}$ and $\alpha\in \overline{\F_p}$ such that the order of does not divide $20$, $30$ or $48$. Let $\xi \in \overline{\F_p}$ be a primitive fifth-root of unity. We set $\F_q=\F_p(\alpha)$. The irreducible representations are given by the $H_4$-graphs in subsection \ref{sectionnewH4graphs} of the Appendix. The highest dimensional representation is of dimension $48$ and our results for this representation are still conjectural in some cases.

\begin{Def}
The Iwahori-Hecke $\mathcal{H}_{H_4,\alpha}$ is the $\F_q$-algebra generated by the generators $S_1,S_2,S_3,S_4$ and the following relations :
\begin{enumerate}
\item $\forall i \in \{1,2,3,4\}, (S_i-\alpha)(S_i+1)=0.$
\item $S_1S_2S_1S_2S_1=S_2S_1S_2S_1S_2$,
\item $S_1S_3=S_3S_1$,
\item $S_1S_4=S_4S_1$,
\item $S_2S_3S_2=S_3S_2S_4$,
\item $S_2S_4=S_4S_2$,
\item $S_3S_4S_3=S_4S_3S_4$.
For $\sigma$ in the Coxeter group $H_4$, if $\sigma=s_{i_1}\dots s_{i_k}$ is a reduced expression we set $T_{\sigma}=S_{i_1}\dots S_{i_k}$.
\end{enumerate}
\end{Def}
 
We then use the $W$-graphs given in the Appendix which are the $W$-graphs given in \cite{G-P} together with the new ones verifying the conditions in Theorem \ref{bilinwgraphs}.
 
 \begin{prop}
If the order of $\alpha$ does not divide $12, 20$ or $30$ then the Iwahori-Hecke algebra is split semi-simple and the models given by specialization of the $W$-graphs are irreducible and pairwise non-isomorphic. The restriction rules from the generic case then apply to the specialized case.
\end{prop}

\begin{proof}
As in type $H_3$, we let $A=\Z[\frac{1+\sqrt{5}}{2}][\sqrt{u}^{\pm 1}]$ and $F=\Q[\sqrt{5}](\sqrt{u})$. We have a symetrizing trace defined by $\tau(T_0)=1$ and $\tau(T_\sigma)=0$ for all $\sigma\in H_4\setminus\{1_{H_4}\}$. The algebra $\mathcal{H}_{H_4,u}$ is then a free $F$-algebra of rank $14400$. 
We have that $A$ is an integrally closed integral domain because $\Z[\frac{1+\sqrt{5}}{2}]$ is integrally closed (see \cite{Bourb} V.3. Corollary 1 and \cite{NivZucMont} Theorem 9.20.).

Let $\theta$ be the ring homomorphism from $A$ to $L=\F_q$ defined by $\theta(\frac{1+\sqrt{5}}{2})=\xi+\xi^{-1}+1$, $\theta(u)=\alpha$ and $\theta(k)=\overline{k}$. We know $FH$ is split. The basis formed by the elements $T_\sigma, \sigma \in H_4$ verifies the conditions of the Proposition \ref{Tits}. 

We now only need to check that the Schur elements associated to these irreducible representations are in $B$ and don't vanish when specialized under $\theta$. For $n\in \N$, $\Phi_n$ is the $n$-th cyclotomic polynomial and $\Phi_{5,a}(u)=u^2+\frac{1+\sqrt{5}}{2}u+1$, $\Phi_{5,b}(u)=u^2+\frac{1-\sqrt{5}}{2}u+1$, $\Phi_{10,a}(u)=u^2-\frac{1+\sqrt{5}}{2}u+1$, $\Phi_{10,b}(u)=u^2+\frac{\sqrt{5}-1}{2}u+1$, $\Phi_{15,a}(u)=u^4-\frac{1+\sqrt{5}}{2}u^3+\frac{1+\sqrt{5}}{2}u^2-\frac{1+\sqrt{5}}{2}u+1$, $\Phi_{15,b}(u)=u^4+\frac{\sqrt{5}-1}{2}u^3+\frac{1-\sqrt{5}}{2}u^2+\frac{\sqrt{5}-1}{2}u+1$, $\Phi_{20,a}(u)=u^4-\frac{1+\sqrt{5}}{2}u^2+1$, $\Phi_{20,b}(u)=u^4+\frac{\sqrt{5}-1}{2}u^2+1$, $\Phi_{30,a}(u)=u^4+\frac{1+\sqrt{5}}{2}u^3+\frac{1+\sqrt{5}}{2}u^2+\frac{1+\sqrt{5}}{2}u+1$ and $\Phi_{30,b}(u)=u^4+\frac{1-\sqrt{5}}{2}u^3+\frac{1-\sqrt{5}}{2}u^2+\frac{1-\sqrt{5}}{2}u+1$. If $\chi$ is an irreducible character then the character $\chi^\star$ associated to the dual representation of $\chi$ has a Schur element $c_{\chi^\star}=a(c_{\chi})$, where $a$ is the involution of $\Q(\sqrt{5})(\sqrt{u}) $ sending $\sqrt{u}$ to $\sqrt{u}^{-1}$. We define the field automorphism of $\Q(\sqrt{5})$, written $x\mapsto \overline{x}$ by $\overline{\frac{\sqrt{5}-1}{2}}=\frac{-1-\sqrt{5}}{2}$ and $\overline{k}=k$ for all $k\in \Q$. They are given Table \ref{SchurelmtsH4} (obtained using Table E.2. of the Appendix  and Corollary 9.3.6 of \cite{G-P}). 

It then only remains to check that the $H_4$-graphs are still connected after specialization. We must then check which weights vanish under $\theta$ since they are connected in the generic case. The $48$-dimensional $H_4$-graph we found verifying the conditions of Theorem \ref{bilinwgraphs} is not defined for $p=29$. We therefore only use the bilinear form obtained using its existence. We consider the one found in \cite{A-L} which is also available in \cite{G-P} and the CHEVIE Package of GAP3 \cite{CHEVIE}. Recall that $\beta=\xi+\xi^{-1}+1$ and $p\notin \{2,3,5\}$. We first prove that for $p\neq 19$, we have
$$(1-\beta)(\beta+1)(2\beta+3)(\beta+2)(3\beta+4)(2\beta+1)(\beta+3)(\beta+5)\neq 0.$$
We have $1-\beta=1-\xi-\xi^{-1}=-\xi^{-1}(\xi^2-\xi+1)=-\xi^{-1}(\Phi_6(\xi))\neq 0$.\\
$\lambda+1=\xi^{-1}\Phi_3(\xi)\neq 0$, $\lambda+2=\xi^{-1}\Phi_2(\xi)\neq 0$ and $\lambda=\xi^{-1}\Phi_4(\xi)\neq 0$.

Assume by contradiction that $2\beta+3=0$. We then have $\beta=-\frac{3}{2}$. We have $\beta^2+\beta-1=0$. It follows that $0=\frac{9}{4}-\frac{3}{2}-1=\frac{9-6-4}{4}=-\frac{1}{4}$, which is absurd.

Assume by contradiction that $3\beta+4=0$. We then have $\beta=-\frac{4}{3}$. It follows that $0=\frac{16}{9}-\frac{4}{3}-1=-\frac{5}{9}$, which is absurd.

Assume by contradiction that $2\beta+1=0$. We then have $\beta=-\frac{1}{2}$. It follows that $0=\frac{1}{4}-\frac{1}{2}-1=\frac{-5}{4}$, which is absurd.

Assume by contradiction that $\beta+3=0$. We then have $\xi+\xi^{-1}=-3$, therefore $0=9-3-1=5$, which is absurd.

Assume by contradiction that $\beta+5=0$. We then have $\xi+\xi^{-1}=-5$ therefore $0=25-5-1=19$. This implies that $p=19$, which contradicts our assumption on $p$.

\smallskip

For the $40$-dimensional representation, the only weight which can vanish is $\frac{7}{3}$. We have that none of the weights vanish for $p\notin \{7,19\}$.

 For $p=7$, only the brown edges in the figure in the Appendix vanish for the $40$-dimensional representation. There are only $2$ such edges in the $H_4$-graph $\tilde{40}_r$. Since the graph is symmetric, we only need to check that the path represented by one of the brown edges can be replaced. One of the brown edges connects the right vertex $x$ with $I(x)=\{s_1,s_2,s_4\}$ to a vertex $y$ with $I(y)=\{s_2,s_4\}$. It can be replaced by the path going through yellow then green then black edges with vertices $x_0=x$, $x_1$, $x_2$ and $x_3=y$ such that $I(x_1)=\{s_2,s_4\}$ and $I(x_2)=\{s_3\}$.
 
 For $p=19$, we only need to consider the $H_4$-graph $\tilde{30_s}$. Only the blue edges vanish and it is clear from the figure in the Appendix that the graph remains connected without those edges.

\begin{table}
\centering
$\begin{array}{lcr}
 1_r & : &  (\Phi_2^4\Phi_3^2\Phi_4^2\Phi_5^2\Phi_6^2\Phi_{10}^2\Phi_{12}\Phi_{15}\Phi_{20}\Phi_{30})(u).\\
 4_t & : &  \dfrac{5+\sqrt{5}}{2}\dfrac{(\Phi_2^4\Phi_3^2\Phi_5^2\Phi_6^2\Phi_{10}^2\Phi_{15}\Phi_{30})(u)}{u\Phi_{5,a}(u)\Phi_{10,a}(u)\Phi_{15,a}(u)\Phi_{30,a}(u)}\\
 \overline{4_t} & : &  \dfrac{5-\sqrt{5}}{2}\dfrac{(\Phi_2^4\Phi_3^2\Phi_5^2\Phi_6^2\Phi_{10}^2\Phi_{15}\Phi_{30})(u)}{u\Phi_{5,b}(u)\Phi_{10,b}(u)\Phi_{15,b}(u)\Phi_{30,b}(u)}\\
 9_s & : &  \dfrac{5+\sqrt{5}}{2}\dfrac{(\Phi_2^4\Phi_4^2\Phi_5^2\Phi_{10}^2\Phi_{20})(u)}{u^2\Phi_{5,a}(u)\Phi_{10,a)(u)}\Phi_{20,b}(u)}\\
 \overline{9_s} & : &  \dfrac{5-\sqrt{5}}{2}\dfrac{(\Phi_2^4\Phi_4^2\Phi_5^2\Phi_{10}^2\Phi_{20})(u)}{u^2\Phi_{5,b}(u)\Phi_{10,b}(u)\Phi_{20,a}(u)}\\
 16_{rr} & : &  \dfrac{2(\Phi_2\Phi_3^2\Phi_5^2\Phi_{15})(u)}{u^3}\\
 16_r & : &  \dfrac{2(\Phi_2\Phi_3^2\Phi_5^2\Phi_{15})(u)}{u^3}\\
 25_r & : &  \dfrac{(\Phi_2^4\Phi_3^2\Phi_4^2\Phi_6^2\Phi_{12})(u)}{u^4}\\
 36_{rr} & : &  \dfrac{(\Phi_2^4\Phi_5^2\Phi_{10}^2)(u)}{u^5}\\
 24_s & : &   120(\dfrac{18+8\sqrt{5}}{2})\dfrac{(\Phi_5^2\Phi_6^2\Phi_{10}^2)(u)}{u^6\Phi_{5,a}^2(u)\Phi_{10,b}^2(u)}\\
 \overline{24_s} & : &  120(\dfrac{18-8\sqrt{5}}{2})\dfrac{(\Phi_5^2\Phi_6^2\Phi_{10}^2)(u)}{u^6\Phi_{5,b}^2(u)\Phi_{10,a}^2(u)}\\
 24_t & : &  30(\dfrac{7-3\sqrt{5}}{2})\dfrac{(\Phi_4^2\Phi_5^2\Phi_{15})(u)}{u^6\Phi_{5,b}(u)^2\Phi_{15,a}(u)}\\
 \overline{24_t} & : &   30(\dfrac{7+3\sqrt{5}}{2})\dfrac{(\Phi_4^2\Phi_5^2\Phi_{15})(u)}{u^6\Phi_{5,a}(u)^2\Phi_{15,b}(u)}\\
 40_r & : &  \dfrac{40}{u^6}(\Phi_3^2\Phi_{10}^2)(u)\\
 48_{rr} & : &  \dfrac{12}{u^6}(\Phi_5^2\Phi_{12})(u)\\
 18_r & : &  \dfrac{10}{u^6}(\Phi_2^4\Phi_5\Phi_6^2)(u)\\
 30_s & : &  30(\dfrac{3-\sqrt{5}}{2})\dfrac{(\Phi_2^4\Phi_{10}^2\Phi_{15})(u)}{u^6\Phi_{10,a}^2(u)\Phi_{15,b}(u)}\\
 \overline{30_s} & : &  30(\dfrac{3+\sqrt{5}}{2})\dfrac{(\Phi_2^4\Phi_{10}^2\Phi_{15})(u)}{u^6\Phi_{10,b}^2(u)\Phi_{15,a}(u)}\\
 16_t & : &  \dfrac{20(\Phi_3^2\Phi_5^2\Phi_{20})(u)}{u^6\Phi_{5,b}(u)^2\Phi_{20,b}(u)}\\
 \overline{16_t} & : &  \dfrac{20(\Phi_3^2\Phi_5^2\Phi_{20})(u)}{u^6\Phi_{5,a}(u)^2\Phi_{20,a}(u)}\\
 6_s & : &  30(\dfrac{3+\sqrt{5}}{2})\dfrac{(\Phi_2^4\Phi_5^2\Phi_{30})(u)}{u^6\Phi_{5,b}(u)^2\Phi_{30,a}(u)}\\
 \overline{6_s} & : &  30(\dfrac{3-\sqrt{5}}{2})\dfrac{(\Phi_2^4\Phi_5^2\Phi_{30})(u)}{u^6\Phi_{5,a}(u)^2\Phi_{30,b}(u)}\\
 8_r & : &  \dfrac{8}{u^6}(\Phi_3^2\Phi_4^2\Phi_5)(u)\\
 8_{rr} & : &  \dfrac{8}{u^6}(\Phi_3^2\Phi_5^2)(u)\\
 10_r & : &  \dfrac{10}{u^6}(\Phi_2^4\Phi_3^2\Phi_{10})(u)
\end{array}$
\caption{Schur elements in type $H_4$}
\label{SchurelmtsH4}
\end{table}

\bigskip

We now prove that the Schur elements do not vanish under $\theta$. The first $9$ representations are not self-dual and all the remaining ones are. The representations $\rho$ are stable by the field automorphism $\overline{ }$ sending $\sqrt{5}$ to $-\sqrt{5}$ if and only if the representation $\overline{\rho}$ does not appear in the list in Table \ref{SchurelmtsH4}.

By the conditions on $\alpha$, none of the cyclotomic polynomials appearing in \ref{SchurelmtsH4} cancel $\alpha$. We have $\Phi_{5,a}\Phi_{5,b}=\Phi_5$, $\Phi_{10,a}\Phi_{10,b}=\Phi_{10}$, $\Phi_{15,a}\Phi_{15,b}=\Phi_{15}$, $\Phi_{20,a}\Phi_{20,b}=\Phi_{20}$ and $\Phi_{30,a}\Phi_{30,b}=\Phi_{30}$, therefore none these polynomials cancel $\alpha$. 

The roles of $\xi+\xi^{-1}$ and $\xi^2+\xi^{-2}$ are symmetric, therefore it remains to check that $3+\xi+\xi^{-1}\neq 0$, $5-8(\xi+\xi^{-1})\neq 0$, $2-3(\xi+\xi^{-1})\neq 0$ and $1-\xi-\xi^{-1}\neq 0$.

\smallskip 

If $3+\xi+\xi^{-1}=0$, then $2-\xi^2-\xi^{-2}=0$. Therefore $(1-\xi)(1-\xi^{-1})=0$, which is absurd.

\smallskip 

If $5-8(\xi+\xi^{-1})=0$, then $\xi+\xi^{-1}=\frac{5}{8}$. Therefore $(\frac{5}{8})^2+\frac{5}{8}-1=0$ and $25+40-8=0$, which implies that $57=0$. We then have $19=0$ because $p\neq 3$. In $\F_{19}$, we have $4^2+4-1=0$ and $(-5)^2+(-5)-1=0$, therefore $\xi+\xi^{-1}\in \{4,14\}$. It follows that $5-8(\xi+\xi^{-1})\in \{5-8\times 4,5-8\times (-5)\}=\{-27,45\}$. Therefore $27=0$ or $45=0$ in $\F_{19}$, which is absurd.

\smallskip

 If $2-3(\xi+\xi^{-1})= 0$, then $\xi+\xi^{-1}=\frac{2}{3}$. Therefore $(\frac{2}{3})^2+\frac{2}{3}-1=0$ and $4+6-3=0$. This would imply $7=0$, therefore $p=7$. Note that $X^2+X-1$ has no roots in $\F_7$ and $\xi+\xi^{-1}$ cancels $X^2+X-1$, therefore $\xi+\xi^{-1}\notin \F_7$ and $\xi+\xi^{-1}=\frac{2}{3}=10=3$. This implies $3\notin \F_7$ which is absurd.

\smallskip

If $1-\xi-\xi^{-1}=0$, then $2+\xi^2+\xi^{-2}=0$. Therefore $(1+\xi)(1+\xi^{-1})=0$, which is absurd. This concludes the proof.

\end{proof}

The restrictions of the representations to $\mathcal{H}_{H_3}$ available in the CHEVIE package of GAP3 \cite{CHEVIE} are given in Table \ref{resH4H3} (we don't write the restriction of a dual representation or a representation obtained after applying the field automorphism sending $\xi+\xi^{-1}$ to $\xi^2+\xi^{-2}$ if the restriction of the representation is already given), where for a representation $R$, the dual representation is denoted by $R'$. In the generic case, they are the same as the induction/restriction tables of the corresponding finite Coxeter groups, therefore they can be computed easily.

\begin{table}
\centering
\begin{tabular}{ |p{1cm}||p{1cm}|p{1cm}|p{1cm}|p{1cm}|p{1cm}|p{1cm}|p{1cm}|p{1cm}|p{1cm}|p{1cm}|   }
 \hline
 & $1_r$ & $1_r'$ & $3_s$ & $3_s'$ & $\overline{3_s}$ & $\overline{3_s'}$ & $4_r$ & $4_r'$ & $5_r$ & $5_r'$ \\
 \hline

 $1_r$   & $1$    & $0$ &   $0$ & $0$ & $0$ & $0$ & $0$ & $0$ & $0$ & $0$\\
 $4_t$   & $1$    & $0$ &   $0$ & $1$ & $0$ & $0$ & $0$ & $0$ & $0$ & $0$\\
 $9_s$   & $1$    & $0$ &   $0$ & $1$ & $0$ & $0$ & $0$ & $0$ & $1$ & $0$\\
 $16_{rr}$   & $1$    & $0$ &   $0$ & $1$ & $0$ & $1$ & $0$ & $1$ & $1$ & $0$\\
 $16_r$   & $1$    & $0$ &   $0$ & $1$ & $0$ & $1$ & $1$ & $0$ & $1$ & $0$\\
 $25_r$   & $1$    & $0$ &   $0$ & $1$ & $0$ & $1$ & $1$ & $1$ & $2$ & $0$\\
 $36_{rr}$   & $1$    & $0$ &   $0$ & $2$ & $0$ & $2$ & $1$ & $1$ & $2$ & $1$\\
 $24_s$   & $0$    & $0$ &   $0$ & $0$ & $1$ & $1$ & $1$ & $1$ & $1$ & $1$\\
 $24_t$   & $0$    & $0$ &   $0$ & $0$ & $1$ & $1$ & $1$ & $1$ & $1$ & $1$\\
 $40_r$   & $0$    & $0$ &   $1$ & $1$ & $1$ & $1$ & $1$ & $1$ & $2$ & $2$\\
 $48_{rr}$   & $0$    & $0$ &   $1$ & $1$ & $1$ & $1$ & $2$ & $2$ & $2$ & $2$\\
 $18_r$   & $0$    & $0$ &   $0$ & $0$ & $0$ & $0$ & $1$ & $1$ & $1$ & $1$\\
 $30_s$   & $0$    & $0$ &   $1$ & $1$ & $1$ & $1$ & $1$ & $1$ & $1$ & $1$\\
 $16_t$   & $0$    & $0$ &   $1$ & $1$ & $0$ & $0$ & $0$ & $0$ & $1$ & $1$\\
 $6_s$   & $0$    & $0$ &   $1$ & $1$ & $0$ & $0$ & $0$ & $0$ & $0$ & $0$\\
 $8_r$   & $0$    & $0$ &   $0$ & $0$ & $0$ & $0$ & $1$ & $1$ & $0$ & $0$\\
 $8_{rr}$   & $0$    & $0$ &   $0$ & $0$ & $0$ & $0$ & $1$ & $1$ & $0$ & $0$\\
 $10_r$   & $0$    & $0$ &   $0$ & $0$ & $0$ & $0$ & $0$ & $0$ & $1$ & $1$\\
 \hline
\end{tabular}
\smallskip

\caption{Restriction from $\mathcal{H}_{H_4,\alpha}$ to $\mathcal{H}_{H_3,\alpha}$.}
\label{resH4H3}
\end{table}

\begin{lemme}\label{normclosiscos}
The normal closure $\ll \mathcal{A}_{H_3}\gg_{\mathcal{A}_{H_4}}$ of $\mathcal{A}_{H_3}$ inside $\mathcal{A}_{H_4}$ is equal to $\mathcal{A}_{H_4}$.
\end{lemme}

\begin{proof}
We note $s_1,s_2,s_3$ and $s_4$ the Coxeter generators of $A_{H_4}$ such that $s_1,s_2$ and $s_3$ generate $A_{H_3}$ and $s_1s_2s_1s_2s_1=s_2s_1s_2s_1s_2$. By \cite{MR}, we have $$\mathcal{A}_{H_3}=<s_1s_2^{-1},s_2s_1s_2^{-2},s_2^2s_1s_2^{-3},s_2^3s_1
s_2^{-4},s_3s_1^{-1}>,$$
$$\mathcal{A}_{H_4}=<s_1s_2^{-1},s_2s_1s_2^{-2},s_2^2s_1s_2^{-3},s_2^3s_1
s_2^{-4},s_3s_1^{-1},s_4s_1^{-1}>.$$
 It follows that we only need to show that $s_4s_1^{-1}\in \ll \mathcal{A}_{H_3}\gg_{\mathcal{A}_{H_4}}$. This is true because
$$s_4s_1^{-1}=s_3s_4s_3(s_3s_4)^{-1}s_1^{-1}=(s_3s_4)(s_3s_1^{-1})(s_3s_4)^{-1}=(s_3s_1^{-1}s_4s_1^{-1})(s_3s_1^{-1})(s_3s_1^{-1}s_4s_1^{-1})^{-1}.$$

\end{proof}

\begin{prop}\label{resH4derivedsubgroup}
The restrictions of the irreducible representations of dimension greater than $1$ of $\mathcal{H}_{H_4,\alpha}$ to $\mathcal{A}_{H_4}$ are absolutely irreducible and pairwise non-isomorphic.
\end{prop}

\begin{proof}
As in \cite{BMM} Lemma $3.4$, we only need to prove that $A_{H_3}$ is generated by $A_{H_3}$ and $\mathcal{A}_{H_4}$. This true because $s_4=s_4s_1^{-1}s_1$ and $s_4s_1^{-1}\in \mathcal{A}_{H_4}$ and $s_1\in A_{H_3}$.

The proof of the second part of the statement is the same as for Proposition \ref{resH3derivedsubgroup}.
\end{proof}

\section{Type $H_4$, low dimensional representations}\label{H4lowdim}

We now determine the image of the Artin group of type $H_4$ inside the low-dimensional representations of the Iwahori-Hecke algebra.

\begin{prop}\label{dim4H4}
\noindent

\begin{enumerate}
\item Assume $1\sim 2$.
\begin{enumerate}
\item If $\F_q=\F_p(\alpha)=\F_p(\alpha+\alpha^{-1})$, then $\rho_{4_t}(\mathcal{A}_{H_4})\simeq SL_4(q^2)$.
\item If $\F_q=\F_p(\alpha)\neq \F_p(\alpha+\alpha^{-1})$, then $\rho_{4_t}(\mathcal{A}_{H_4})\simeq SL_4(q)$.
\end{enumerate}
\item If $1\nsim 2$, then $\rho_{4_t}(\mathcal{A}_{H_4})\simeq SU_4(q^{\frac{1}{2}})$.
\end{enumerate}

If $1\sim 2$ and $\F_p(\alpha)=\F_p(\alpha+\alpha^{-1})$, then $\Phi_{1,2}\circ \rho_{\overline{4_t}|\mathcal{A}_{H_3}}\simeq \rho_{4_t|\mathcal{A}_{H_3}}$.

If $1\sim 2$ and $\F_p(\alpha)\neq \F_p(\alpha+\alpha^{-1})$, then $\Phi_{1,2}\circ \rho_{\overline{4_t}|\mathcal{A}_{H_3}}\simeq \rho_{4_t'|\mathcal{A}_{H_3}}$.
\end{prop}

\begin{proof}
The proof of the first part of the statement is identical to the proof of Proposition \ref{resrepreflicos}.

\smallskip

The proof of the second part of the statement follows from Proposition \ref{resrepreflicos} and Table \ref{resH4H3}.
\end{proof}

\begin{prop}
\noindent

\begin{enumerate}
\item Assume $1\sim 2$.
\begin{enumerate}
\item If $\F_q=\F_p(\alpha)=\F_p(\alpha+\alpha^{-1})$, then $\rho_{6_s}(\mathcal{A}_{H_4})\simeq \Omega_6^+(q^2)$.
\item If $\F_q=\F_p(\alpha)\neq \F_p(\alpha+\alpha^{-1})$, then $\rho_{6_s}(\mathcal{A}_{H_4})\simeq \Omega_6^+(q)$.
\item We have $\Phi_{1,2}\circ \rho_{6_s|\mathcal{A}_{H_4}}\simeq \rho_{\overline{6_s}|\mathcal{A}_{H_4}}$.
\end{enumerate}
\item Assume $1\nsim 2$. We then have $\rho_{6_s}(\mathcal{A}_{H_4})\simeq \Omega_6^+(q^{\frac{1}{2}})$.
\end{enumerate}
\end{prop}

\begin{proof}

Let $\F_r=\F_p(\alpha,\xi+\xi^{-1})$. By Proposition \ref{bilinwgraphs}, we have $\rho_{6_s}(\mathcal{A}_{H_4})\leq \Omega_6^+(q')$, where $\F_{q'}=\F_p(\sqrt{\alpha},\xi+\xi^{-1})$. By Proposition \ref{color}, we also have $\rho_{6_s}\simeq \sigma \circ \rho_{6_s}$, where $\sigma=\op{I_d}_{\F_r}$ if $\F_q'=\F_r$ and $\sigma$ is the automorphism of order $2$ of $\F_{q'}$ otherwise. We then have by Proposition 4.1 of \cite{BMM} that up to conjugation, we have $\rho_{6_s}(\mathcal{A}_{H_4})\leq \Omega_6^+(r)$.

\medskip

Assume first that $1\sim 2$. By Proposition \ref{IsomorphismH3}, we have that $\F_p(\alpha,\xi+\xi^{-1})=\F_p(\alpha+\alpha^{-1},\xi+\xi^{-1})$.
By Table \ref{resH4H3} and Theorem \ref{resicos}, we have that $\rho_{6_s}(\mathcal{A}_{H_3})$ is up to conjugation a twisted $SL_3(r)$. It follows that $\rho_{6_s}(\mathcal{A}_{H_4})$ is an irreducible subgroup of $\Omega_6^+(r)$ generated by long root elements. By Kantor's Theorem \ref{theoKantor}, we have that $\rho_{6_s}(\mathcal{A}_{H_4})$ is conjugate in $GL_6(r)$ to one of the following groups 
\begin{enumerate}
\item $\Omega_6^{\pm}(r)$,
\item $\Omega_6^{-}(\sqrt{r})\leq \Omega_6^{+}(r)$.
\end{enumerate}

\medskip

We have $\vert \Omega_6^{-}(\sqrt{r})\vert=\sqrt{r}^6(\sqrt{r}^3+1)(\sqrt{r}^4-1)(\sqrt{r}^2-1)/2 $. This means that if $\rho_{6_s}(\mathcal{A}_{H_4}$ was isomorphic to 
$\Omega_6^{-}(\sqrt{r})$ then we would have that $r^3(r^2-1)(r^3-1)$ would divide $r^3(\sqrt{r}^3+1)(r^2-1)(r-1)$, therefore $r^3-1$ would divide $(r^{\frac{3}{2}}+1)(r-1)=r^{\frac{5}{2}}-r^{\frac{3}{2}}+r-1<r^3-1$ which is absurd because both quantities are positive.

\medskip

We cannot have $\rho_{6_s}(\mathcal{A}_{H_4})\simeq \Omega_6^{-}(r)$ because otherwise we would have $\Omega_6^{-}(r)\leq \Omega_6^{+}(r)$, therefore $r^6(r^3+1)(r^4-1)(r^2-1)$ would divide $r^6(r^3-1)(r^4-1)(r^2-1)$. This would imply that $r^3+1$ divides $r^3-1$ which is absurd.

\medskip

This implies that $\rho_{6_s}(\mathcal{A}_{H_4})\simeq \Omega_6^{+}(r)$. 

\medskip

We know that $\Phi_{1,2}\circ \rho_{6_s|\mathcal{A}_{H_4}}\simeq \rho_{6_s|\mathcal{A}_{H_4}}$ or $\Phi_{1,2}\circ \rho_{6_s|\mathcal{A}_{H_4}}\simeq \rho_{\overline{6_s}|\mathcal{A}_{H_4}}$. The Table \ref{resH4H3} and Proposition \ref{resrepreflicos} imply that the latter possibility is the only one possible.

\bigskip

Assume now $1\nsim 2$. By Proposition \ref{IsomorphismH3}, we have that $\F_p(\alpha,\xi+\xi^{-1})\neq \F_p(\alpha+\alpha^{-1},\xi+\xi^{-1})$. There exists then a unique automorphism $\epsilon$ of order $2$ of $\F_r$. We have $\epsilon(\alpha)=\alpha^{-1}$. It follows by Proposition \ref{Fieldfactorization} that $\epsilon \circ \rho_{6_s|\mathcal{A}_{H_4}}\simeq \rho_{6_s|\mathcal{A}_{H_4}}$ or $\epsilon \circ \rho_{6_s|\mathcal{A}_{H_4}}\simeq \rho_{\overline{6_s}|\mathcal{A}_{H_4}}$. By Table \ref{resH4H3} and Proposition \ref{resrepreflicos}, we have that $\epsilon \circ \rho_{6_s|\mathcal{A}_{H_4}}\simeq \rho_{6_s|\mathcal{A}_{H_4}}$. This implies by Proposition $4.1$ of \cite{BMM} that, up to conjugation in $GL_4(r)$, we have $\rho_{6_s|\mathcal{A}_{H_4}}\leq \Omega_6^+(r^{\frac{1}{2}})$. By Proposition \ref{resrepreflicos} and Table \ref{resH4H3}, we have that $\rho_{6_s}(\mathcal{A}_{H_3}\simeq SU_3(r^{\frac{1}{2}})$. We can again apply Theorem \ref{theoKantor}, and we have that  $\rho_{6_s}(\mathcal{A}_{H_4})$ is conjugate in $GL_6(r)$ to one of the following groups 
\begin{enumerate}
\item $\Omega_6^{+}(r^{\frac{1}{2}})$,
\item $\Omega_6^{-}(r^{\frac{1}{4}})\leq \Omega_6^{+}(r)$.
\end{enumerate}

We have $\vert SU_3(r^{\frac{1}{2}})\vert =r^{\frac{3}{2}}(r-1)(r^{\frac{3}{2}}+1)$ and $\vert \Omega_6^{-}(r^{\frac{1}{4}})\vert=r^{\frac{3}{2}}(r^{\frac{3}{4}}+1)(r-1)(r^{\frac{1}{2}}-1)/2 $. Assume by contradiction that we are in the second case. We then have that $(r^{\frac{3}{2}}+1)$ divides $(r^{\frac{3}{4}}+1)(r^{\frac{1}{2}}-1)/2=\frac{1}{2}(r^{\frac{5}{4}}-r^{\frac{3}{4}}+r^{\frac{1}{2}}-1)<r^{\frac{3}{2}}+1$. This is a contradiction, therefore we have that $\rho_{6_s}(\mathcal{A}_{H_4})\simeq \Omega_6^{+}(r^{\frac{1}{2}})$.

\medskip

The proof is then concluded by Proposition \ref{resrepreflicos}.
\end{proof}

\begin{prop}
If $\F_q=\F_p(\alpha)=\F_p(\alpha+\alpha^{-1})$, then we have $\rho_{8_r}(\mathcal{A}_{H_4})\simeq\Omega_8^+(q)$.

If $\F_q=\F_p(\alpha)\neq \F_p(\alpha+\alpha^{-1})$, then we have $\rho_{8_r}(\mathcal{A}_{H_4})\simeq\Omega_8^+(q^{\frac{1}{2}})$.
\end{prop}

\begin{proof}
Assume first that $\F_p(\alpha)=\F_p(\alpha+\alpha^{-1})$.
By Proposition \ref{bilinwgraphs}, we have that $G=\rho_{8_r}(\mathcal{A}_{H_4})\leq \Omega_8^+(q)$ up to conjugation in $GL_8(q)$. We will use the same theorem as for the $6$-dimensional representation but one of the cases will be much more technical to exclude. We have by Table \ref{resH4H3} and Theorem \ref{resicos} that $\rho_{8_r}(\mathcal{A}_{H_3})$ is conjugate in $GL_8(q)$ to a twisted diagonal $SL_4(q)$. It follows that $G$ is irreducible and generated by long root elements. By Theorem \ref{theoKantor}, we have that $G$ belongs to the following list
\begin{enumerate}
\item $\Omega_8^{+}(q)$
\item $\Omega_8^{-}(\sqrt{q})$
\item $SU_4(q)$
\item $G/Z(G)=P\Omega_7(q)$, $Z(G)=2$
\item $^3\!D_4(\sqrt[3]{q})$
\end{enumerate}

We first exclude cases $\mathbf{2}$, $\mathbf{3}$ and $\mathbf{5}$ because they cannot occur by simple cardinality arguments. We will then exclude case $\mathbf{4}$ showing that $G$ contains a group which cannot be contained in a group of the same order as in case $\mathbf{4}$. We have $\vert SL_4(q)\vert =q^6(q^2-1)(q^3-1)(q^4-1)$.

\medskip

We have $\vert \Omega_8^{-}(\sqrt{q})\vert =\sqrt{q}^{12}(\sqrt{q}^4+1)(\sqrt{q}^2-1)(\sqrt{q}^4-1)(\sqrt{q}^6-1)=q^6(q^2+1)(q-1)(q^2-1)(q^3-1)=\vert SL_4(q)\vert \frac{(q^2+1)(q-1)}{q^4-1}<\vert SL_4(q)\vert$, therefore the second case is excluded.

\medskip

We have $\vert SU_4(q)\vert =q^6(q^2-1)(q^3+1)(q^4-1)$ and $q^3-1$ cannot divide $q^3+1$ because it is greater than $2$, therefore the third case is excluded.

\medskip

We have $\vert ^3\!D_4(\sqrt[3]{q})\vert = \sqrt[3]{q}^{12}(\sqrt[3]{q}^{8}+\sqrt[3]{q}^{4}+1)(\sqrt[3]{q}^{6}-1)(\sqrt[3]{q}^{2}-1)$, therefore it cannot contain a group isomorphic to $SL_4(q)$ because $q^6$ does not divide $q^4$.

\medskip

We now want to show that the fourth case is also excluded but it can contain a twisted diagonal $SL_4(q)$. We therefore have to construct a different subgroup of our group $G$.

We order the vertices with the graded lexicographic order ($I(x_1)=\{s_1\}$, $I(x_2)=\{s_2\}$, $I(x_3)=\{s_1,s_3\}$, $I(x_4)=\{s_1,s_4\}$, $I(x_5)=\{s_2,s_3\}$, $I(x_6)=\{s_2,s_4\}$, $I(x_7)=\{s_1,s_3,s_4\}$ and $I(x_8)=\{s_2,s_3,s_4\}$) and consider the matrices with respect to the associated basis ordered the same way.

Let $P=\begin{pmatrix}
1 & 0 & 0 & 1 & 0 & 0 & 0 & 0\\
0 & 1 & 0 & 0 & 0 & 1 & 0 & 0\\
0 & 0 & 1 & 0 & 0 & 0 & 1 & 0\\
0 & 0 & 0 & 0 & 1 & 0 & 0 & 1\\
0 & 0 & 0 & 0 & 1 & 0 & 0 & -1\\
0 & 0 & -1 & 0 & 0 & 0 & 1 & 0\\
0 & 1 & 0 & 0 & 0 & -1 & 0 & 0\\
-1 & 0 & 0 & 1 & 0 & 0 & 0 & 0
\end{pmatrix}$ and\\ $X =\begin{pmatrix}
1 & 0 & -\frac{v^2-v+1}{v} & 0 & 0 & 0 & 0 & 0\\
1 & -\frac{v^2-v+1}{v} & 0 & 0 & 0 & 0 & 0 & 0\\
1 & 0 & 0 & 0 & 0 & 0 & 0 & 0\\
1 & -\frac{v^2+1}{v} & -\frac{v^2+1}{v} & \frac{v^4+v^3+2v^2+v+1}{v^2} & 0 & 0 & 0 & 0\\
0 & 0 & 0 & 0 & 0 & 0 & 1 & \frac{v}{v^2+v+1}\\
0 & 0 & 0 & 0 & 0 & 1 & 0 & 0\frac{v}{v^2+v+1}\\
0 & 0 & 0 & 0 & 1 & \frac{v}{v^2-v+1} & \frac{v}{v^2-v+1} & \frac{v^2}{(v^2+1)(v^2-v+1)}\\
0 & 0 & 0 & 0 & 0 & 0 & 0 & 1\\
\end{pmatrix}$.

We then have for $i\in \{1,2,3\}$, $P\rho_{8_r}(S_i)P^{-1}=\begin{pmatrix} \rho_{4_r'}(S_i) & 0\\
0 & -\alpha~^t\!\rho_{4_r'}(S_i)^{-1}\end{pmatrix}$.

We also have for $i\in \{2,3\}$, 
$(XP)\rho_{8_r}(S_i)(XP)^{-1}=$
$$\begin{pmatrix}
\begin{pmatrix} \rho_{2_r}(S_i) & 0\\
0 & \begin{pmatrix} \rho_{1_r}(S_i) & 0\\
0 & -\alpha~^t\!\rho_{1_r}(S_i)^{-1} \end{pmatrix}\end{pmatrix} & 0\\
0 & \begin{pmatrix} -\alpha~^t\!\rho_{2_r}(S_i)^{-1} & 0\\
0 & \begin{pmatrix} -\alpha~^t\!\rho_{1_r}(S_i)^{-1} & 0\\
0 & \rho_{1_r}(S_i)\end{pmatrix}\end{pmatrix}
\end{pmatrix}$$
where $2_r$ and $1_r$ are given by the following W-graphs and the bases are ordered in the anti-lexicographic way for $2_r$ ($I(e_{x_1})=\{s_3\}$ and $I(e_{x_2})=\{s_2\}$)

\begin{center}
\begin{tikzpicture}
[place/.style={circle,draw=black,
inner sep=1pt,minimum size=10mm}]

\node at (0.5,0.8)[place]{$\emptyset$};
\draw (0,0) node{$1_r$};

\node (1) at (6.5,0.8)[place]{$2$};
\node (2) at (8.4,0.8)[place]{$3$};
\draw (1) to (2);
\draw (6,0) node{$2_r$};

\end{tikzpicture}
\end{center}

We have $\rho_{4_r'}(\mathcal{A}_{H_4})=SL_4(q)$ and by  Lemma $3.5$ of \cite{BM}, $\rho_{2_r}(\mathcal{A}_{A_2})=SL_2(q)$, where $A_{A_2}=<S_2,S_3>$. We note $H=(XP) \rho_{8_r}(\mathcal{A}_{A_2})(XP)^{-1}$. We will consider a large subgroup of the normalizer of $H$ and show it cannot be contained in a group corresponding to the fourth case of our list. Let
$$n=(XP)\rho_{8_r}((S_1S_3S_2)^5(S_2S_3S_2)^{-5})(XP)^{-1},$$
$$m=(XP)\rho_{8_r}((S_3S_2S_4)^4(S_2S_3S_2)^{-4})(XP)^{-1}.$$
 Let then $u=[n,m]$. We have for $i\in \{2,3\}$, 
$$(XP)^{-1}u(XP)\rho_{8_r}(S_i)=\rho_{8_r}(S_i)(XP)^{-1}u(XP), (XP)^{-1}m(XP)\rho_{8_r}(S_i)=\rho_{8_r}(S_i)(XP)^{-1}m(XP).$$

This shows that $N_1=<u,m>\leq C_G(H)$. We will now determine what this group $N_1$ is. If we let $$R=\begin{pmatrix}
1 & 1 & 0 & 0 & 0 & 0 & 0 & 0\\
0 & 0 & 0 & 0 & 1 & 1 & 0 & 0\\
0 & 0 & 1 & 1 & 0 & 0 & 0 & 0\\
0 & 0 & 0 & 0 & 0 & 0 & 1 & 1\\
0 & 0 & 0 & 0 & 1 & -1 & 0 & 0\\
1 & -1 & 0 & 0 & 0 & 0 & 0 & 0\\
0 & 0 & 0 & 0 & 0 & 0 & 1 & -1\\
0 & 0 & 1 & -1& 0 & 0 & 0 & 0\\
\end{pmatrix},$$
we have 
$$R^{-1}N_1R\subset \{\begin{pmatrix}
M_1 & 0 & 0 & 0\\
0 & M_2 & 0 & 0\\
0 & 0 & M_3 & 0\\
0 & 0 & 0 & M_4
\end{pmatrix}, M_1,M_2,M_3,M_4\in GL_2(q)\}.$$
Let $\pi_1,\pi_2,\pi_3$ and $\pi_4$ be the corresponding projections and $u'=R^{-1}uR, m'=R^{-1}mR$ and $N_1'=<u',m'>$. We will first show that $\pi_1(N_1')=SL_2(q)$ then determine the images under the other projections and then determine fully $N_1'$ using Goursat's Lemma.

\bigskip

Let $u_1=\pi_1(u')$ and $n_1=\pi_1(m')$. We have $\op{det}(u_1)=\op{det}(n_1)=1$. As in the proof of Theorem \ref{platypodes}, if we prove that $G_1=<u_1,m_1>$  contains elements whose traces generate $\F_q$, that $\overline{G_1}\notin\{\mathfrak{A}_5,\mathfrak{S}_4\}$ and that $\overline{G_1}$ is not abelian by abelian, then we have $G_1=SL_2(q)$.

\medskip

Let us first show that $G_1$ is not abelian by abelian. Assume by contradiction that $G_1$ is abelian by abelian, we would then have $[u_1,n_1][u_1,n_1^{-1}]\in \{\pm [u_1,n_1^{-1}][u_1,n_1]\}$. Let $A_1=[u_1,n_1][u_1,n_1^{-1}]+[u_1,n_1^{-1}][u_1,n_1]$ and $A_2=[u_1,n_1][u_1,n_1^{-1}]-[u_1,n_1^{-1}][u_1,n_1]$.

\smallskip

Assume $A_1=0$, we then have $a_1=\frac{4 A_1[1,1]\alpha^{15}}{(\alpha-1)^{5}(\alpha^2+1)^5(\alpha^4+1)\Phi_3(\sqrt{\alpha})\Phi_6(\sqrt{\alpha})}=0$ and \\$a_2=-\frac{4 A_1[1,2]\alpha^{15}}{(\alpha-1)^{6}(\alpha^2+1)^5(\alpha^4+1)\Phi_3(\sqrt{\alpha})\Phi_6(\sqrt{\alpha})}=0$, where $A_1[i,j]$ is the coefficient in row $i$ and column $j$ of $A_1$. If we let $v=\sqrt{\alpha}$, we have 
$$a_1=v^{18}+v^{17}-v^{16}+2v^{14}-2v^{12}+2v^{11}+2v^{10}+2v^{9}-2v^8+2v^7+2v^6-2v^4+v^2+v-1,$$
$$a_2=v^{16}-v^{15}+2v^{14}-v^{13}+2v^{12}-v^{11}+2v^{10}+v^9+2v^8-v^7+2v^6+v^5+2v^4+v^3+2v^2+v+1.$$

We will show that for any prime $p$, those two polynomials in $\F_p[v]$ are coprime. We write $\op{Rem}$ the Euclidean remainder in $\F_p[v]$. We have 

\begin{eqnarray*}
a_3  = & \frac{1}{4v^3}\op{Rem}(a_1,a_2) & =-v^{12}+v^{11}-v^{10}-v^8-v^4-v^2-v-1,\\
a_4  = &\op{Rem}(a_2,a_3) & =-v^{11}+v^{10}+v^9+v^8-v^7+v^2+v+1,\\
a_5  = &\op{Rem}(a_3,a_4) & =-2v^{10}-v^9-v^4-v^3-2v^2-2v-1,\\
a_6  = &4\op{Rem}(a_4,a_5) & = v^9+4v^8-4v^7+2v^5-v^4+v^3+2v^2+1,\\
a_7  = &\op{Rem}(a_5,a_6) & =  -36v^8+28v^7+4v^6-16v^5+8v^4-4v^3-16v^2-8,\\
a_8  = &81\times \op{Rem}(a_6,a_7) & = -14v^7+7v^6+8v^5-4v^4+2v^3-10v^2-18v-5,\\
a_9  = &\frac{7}{81}\op{Rem}(a_7,a_8) & = -v^6+2v^3+2v^2-1,\\
a_{10}  = & \frac{1}{4}\op{Rem}(a_8,a_9) & = 2v^5-8v^4-3v^3+v^2-v-3,\\
a_{11}  = & \frac{2}{7}\op{Rem}(a_9,a_{10}) & = -5v^4-v^3+v^2-v-2,\\
a_{12}  = & 25\times \op{Rem}(a_{10},a_{11}) & = -23v^3-27v^2-3v+9,\\
a_{13}  = & \frac{23^2}{50}\op{Rem}(a_{11},a_{12}) & = -43v^2-38v-1,\\
a_{14}  = & \frac{43^2}{4\times 23^2}\op{Rem}(a_{12},a_{13}) & = 3v+8,\\
a_{15}  = & -\frac{9}{43^2}\op{Rem}(a_{13},a_{14}) & =1.
\end{eqnarray*}

These computations are only correct if $p\notin \{2,3,5,7,23,43\}$ because the computations are made in $\Q[v]$ and they can only be specialized if the polynomials considered are of same degree in $\F_p$ and in $\Q$. We have $p\notin \{2,3,5\}$ by assumption. We check the remaining cases using the \verb+Gcd+ function in GAP4 with $a_1$ and $a_2$ polynomials in the indeterminate \verb+v:=Indeterminate(GF(p))+. This shows that $A_1$ cannot be equal to zero.

\smallskip

Assume now $A_2=0$, we then have $b_1=\frac{4\alpha^{15}A_2[1,1]}{\alpha^2-\alpha+1}=0$, $b_2=-\frac{4\alpha^{15}A_2[1,2]}{(\alpha^2-\alpha+1)(\alpha-1)^3(\alpha^2+1)^3}=0$ and $b_3=\frac{4\alpha^{15}A_2[1,2]}{(\alpha^2-\alpha+1)(\alpha-1)^3(\alpha^2+1)^3(\alpha+1)^2}=0$. We then have that \\$b_4=\op{Rem}(\frac{\op{Rem}(b_2,b_3)}{8},9\op{Rem}(b_3,\frac{\op{Rem}(b_2,b_3)}{8}))\times \frac{1}{-9v(v^2+v+1)(v^4+1)}=0$.
 
 We then show using the same techniques as for $A_1$ that $b_1$ and $b_4$ are coprime for \\
 $p\notin\{ 2, 3, 5, 7, 11, 13, 17, 47, 167, 233, 293, 449, 5303, 13649, 15797,25913, 245071\}$. Again, by assumption, we have $p\notin\{2,3,5\}$ and using GAP4, we get that $b_1$ and $b_4$ are coprime except for $p=11$. For $p=11$, we have $\op{Gcd}(b_1,b_4)=v^4+6v^2+1=(v^2-6v-1)(v^2+6v-1)$. Again, using GAP4, we see that this polynomial divides $v^{24}-1$ in $\F_{11}[v]$, therefore if $\sqrt{\alpha}$ is a root of this polynomial then $\alpha^{12}=1$ which contradicts our assumptions. This concludes the proof that $G_1$ is not abelian by abelian.
 
 \medskip
 
 Assume now that $\overline{G_1}\in\{\mathfrak{A}_5,\mathfrak{S}_4\}$. We then have that its elements are of order less than or equal to $5$.
 The eigenvalues of $n_1$ are $\alpha^2$ and $\alpha^{-2}$, therefore if $\overline{n_1}^r=\overline{I_2}$ then $\alpha^{2r}\in \{-1,1\}$, therefore $\alpha^{4r}=1$ and $r\geq 7$ by the assumptions on $\alpha$. This is absurd since the elements in $\overline{G_1}$ are of order less than or equal to $5$. It also implies that $\overline{G_1}\simeq \mathfrak{S}_4$ because $\mathfrak{A}_5$ contains no element of order $4$. We have shown above that the commutator subgroup of $\overline{G_1}$ was not abelian. This leads to a contradiction because the commutator subgroup of $\mathfrak{S}_4$ is the Klein group of order $4$, which is abelian.
 
 \medskip
 
  It only remains to show that $\F_q$ is generated by traces of elements in $G_1$. We have $\op{T_r}(n_1)=\alpha^2+\alpha^{-2}$ and $\tr(u_1)=\alpha^4+\alpha^{-4}-(\alpha^3+\alpha^{-3})+2(\alpha^2+\alpha^{-2})-3(\alpha^2+\alpha^{-2})+4$. We then have $\tr(n_1)+2=(\alpha+\alpha^{-1})^2\neq 0$ and 
 $$\frac{\tr{n_1}^2+2\tr(n_1)+2-\tr(u_1)}{\tr(n_1)+2}=\alpha+\alpha^{-1}.$$
 
 This shows that $\alpha+\alpha^{-1}$ belongs to the field generated by traces of the elements of $G_1$, therefore $\F_q$ is generated by traces of the elements of $G_1$. It follows that $G_1=SL_2(q)$.
 
 \bigskip If $u_3=\pi_3(u')$, $n_3=\pi_3(m')$ and $C_{13}=\begin{pmatrix}
 0 & 1\\
 1 & 0
\end{pmatrix}$, then we have that $C_{13} u_1 C_{13}^{-1}=u_3$ and $C_{13} n_1 C_{13}^{-1}=n_3$. We therefore get that $G_3=<u_3,n_3>=SL_2(q)$ and the representation $\pi_3$ of $N_1'$ factors through $\pi_1$.

\bigskip

We now show that if $u_2=\pi_2(u')$ and $n_2=\pi_2(m')$ then $G_2=<u_2,n_2>\simeq SL_2(q).<\alpha^4>$, where $SL_2(q).<\alpha^4>$ denotes a group having $SL_2(q)$ as a normal subgroup and the subgroup of $\F_q^\star$ generated by $\alpha^4$ as its quotient group by $SL_2(q)$. Again as in the proof of Theorem \ref{platypodes}, to prove that $G_2$ contains $SL_2(q)$ as a normal subgroup, we only need to show that $\overline{G_2}$ is not abelian by abelian, that $\overline{G_2}$ contains elements of order greater than $6$ and that the field generated by the traces of elements of $[G_2,G_2]$ is $\F_q$.

\medskip

We first show that $\overline{G_2}$ is not abelian by abelian. Assume $G_2$ is abelian by abelian. We then have $B_1=[u_2,n_2][u_2,n_2^{-1}]-[u_2,n_2^{-1}][u_2,n_2]=0$ or $B_2=[u_2,n_2][u_2,n_2^{-1}]+[u_2,n_2^{-1}][u_2,n_2]=0$.

\smallskip

Assume $B_1=0$, we then have $0=B_1[1,2]=\frac{(\alpha^3-1)(\alpha^5-1)(\alpha-1)^3(\alpha^4+1)^2}{\sqrt{\alpha}^{15}}$, which is absurd.

\smallskip

Assume $B_2=0$. We then have $0=b_2=\frac{\sqrt{\alpha}^{15}B_2[1,2]}{(\alpha^2-\alpha+1)(\alpha-1)^3(\alpha^4+1)}$ and $0=b_3=\frac{\sqrt{\alpha}^{23}B_2[2,1]}{(\alpha^2-\alpha+1)(\alpha^8-1)(\alpha-1)^2}$. It follows that 

$$0=b_2=(\alpha^4-\alpha^3-\alpha^2-\alpha+1)(\alpha^6+2\alpha^5+3\alpha^4+2\alpha^3+3\alpha^2+2\alpha+1).$$

$$0=b_3=(\alpha^4+\alpha^2+\alpha+1)(\alpha^4-\alpha^3-\alpha^2-\alpha+1).$$

We then have $0=b_4=\frac{1}{(\alpha+1)(\alpha-\sqrt{\alpha}+1)(\alpha+\sqrt{\alpha}+1)}\op{Rem}(b_2,b_3)=\alpha^4-\alpha^3-\alpha^2-\alpha+1$.

We also have $b_1=\frac{\alpha^9B_2[1,1]}{\alpha^2-\alpha+1}=0$, therefore

$$0= \alpha^{17}-2\alpha^{16}-\alpha^{14}+\alpha^{13}-\alpha^{12}+2\alpha^{11}+\alpha^{10}-2\alpha^9+3\alpha^8-3\alpha^7+\alpha^6-\alpha^5+5\alpha^8-2\alpha^3-\alpha^2-\alpha+1.$$

We then have
\begin{eqnarray*}
b_5 = & \frac{1}{2}\op{Rem}(b_1,b_4) & = 83\alpha^3+60\alpha^2+20\alpha-48\\
b_6 = & 83^2\op{Rem}(b_4,b_5) & = 31\alpha^2-45\alpha+25\\
b_7 = & 31^2\op{Rem}(b_5,b_6) & =2\times 3\times 5\times 83^2\alpha-186003\\
b_8 = & \frac{2^2\times 5^2}{31^2}\op{Rem}(b_6,b_7) & = 1.
\end{eqnarray*}

Those computations show that $1=0$ for $p\notin\{31,83\}$. For $p\in  \{31,83\}$, we check using GAP4 that $b_1$ and $b_4$ are coprime in $\F_p[v]$. It follows that $B_2\neq 0$ and $\overline{G_2}$ is not abelian by abelian.

\medskip

We must now show that $\overline{G_2}\notin \{\mathfrak{A}_5,\mathfrak{S}_4\}$. Assume $\overline{G_2}\in\{\mathfrak{A}_5,\mathfrak{S}_4\}$. We have that the order of $n_2$ must belong to $\{1,2,3,5\}$. The eigenvalues of $n_2$ are $1$ and $\alpha^{-4}$. Therefore, if $\overline{n_2}^r=I_2$, then we have $\alpha^{4r}=1^r\in \{-1,1\}$. It follows that $\alpha^{4r}=1$. By the assumptions on the order of $\alpha$, we cannot have $r\in \{1,2,3,4,5\}$, therefore we have a contradiction.

\medskip

We now show that the traces of the elements of $[G_2,G_2]$ generate $\F_q$. Let $\F$ be the field generated by those traces. We have $C_1=-(\tr([u_2,n_2])-4)=\alpha^5+\alpha^{-5}-2(\alpha^4+\alpha^{-4})+(\alpha^3+\alpha^{-3})-(\alpha^2+\alpha^{-2})+2(\alpha+\alpha^{-1})\in \F$. We then have $C_1-2=\frac{(\alpha^3-1)(\alpha^5-1)(\alpha-1)^2}{\alpha^5}\neq 0$, therefore 

$$C_2=\frac{\tr([[u_2,n_2],n_2])-C_1^2+6C_1-10}{C_1-2}=\alpha^4+\alpha^{-4}\in \F.$$

It follows by induction that $\alpha^{4r}+\alpha^{-4r}\in \F$ for all $r\in \N$.

We then have $C_3=C_1+2C_2=\alpha^5+\alpha^{-5}+\alpha^3+\alpha^{-3}-(\alpha^2+\alpha^{-2})+2(\alpha+\alpha^{-1})\in \F$.

We also have $C_4=\tr([n_2,u_2][u_2^{-1},n_2]+(C_1-2)(\alpha^8+\alpha^{-8})-2C_1^2+11C_1-(\alpha^4+\alpha^{-4})-16\in \F$, therefore

$$C_4=\alpha^7+\alpha^{-7}-2(\alpha^6+\alpha^{-6})+2(\alpha^3+\alpha^{-3})-3(\alpha^2+\alpha^{-2})+\alpha+\alpha^{-1}\in \F.$$

We then have 

$$C_5=\frac{1}{2^4\times 3}(\tr([[[u_2,n_2],n_2],n_2])-(\alpha^{20}+\alpha^{-20})-52(\alpha^{16}+\alpha^{-16})+(6C_4+34C_3-427)(\alpha^{12}+\alpha^{-12}))+$$
$$(60C_4-5C_3^2+255C_3-1601)(\alpha^8+\alpha^{-8})+(210C_4-13C_3^2+726C_3-16\times 211)(\alpha^4+\alpha^{-4})+$$

$$336C_4+1034C_3-16\times 269)\in \F.$$

Since $C_5=\alpha^2+\alpha^{-2}$, we get that $\alpha^2+\alpha^{-2}\in \F$. We then have $\frac{C_3+C_5}{C_5^2}=\alpha+\alpha^{-1}\in \F$. It follows that $\F=\F_q$ and we have $G_2=SL_2(q)$.

\bigskip

Let $u_4=\pi_4(u')$ and $n_4=\pi_4(m')$.

 Let $C_{24}=\begin{pmatrix} v^4+v^3+3v^2+v+1 & v^4+v^3+v^2+v+1\\
-v^4-v^3-v^2-v-1 & -v^4-v^3-3v^2-v-1
\end{pmatrix}$, we have $\op{det}(C_{24})=-4\alpha(\alpha+1)(v^2+v+1)$.

 We also have $C_{24}u_2C_{24}^{-1}=^t\!u_4^{-1}$ and $C_{24}n_2C_{24}^{-1}=^t\!n_4^{-1}$.
 
 \bigskip
 
 By what has been done above, we have that $N_1'=<u',m'>\simeq \pi_1\times \pi_2(N_1')\leq SL_2(q)\times GL_2(q)$, $\pi_1(N_1')=SL_2(q)$ and $SL_2(q)\subset \pi_2(N_1')$.  We will show that $[N_1',N_1']\simeq SL_2(q)\times SL_2(q)$. By Goursat's lemma, we only need to show that we cannot have $\overline{\pi_1}\simeq \Phi(\overline{\pi_2})$, which is equivalent to $\pi_1\simeq \Phi(\pi_2)\otimes \chi$, where $\chi:N_1'\rightarrow \F_q^\times$ is a character of $N_1'$ and $\Phi\in \op{Aut}(\F_q)$. Assume such a character exists, the following proof will be quite computational. We have that there exists $M\in GL_2(q)$ such that for all $h\in [N_1',N_1']$, $\tr(\pi_1(h))=\Phi(\tr(\pi_2(h)))\chi(h)$ and $\chi(h)\in\{-1,1\}$.

\smallskip

Let $D_1=\tr([u_1,n_1])$ and $E_1=\tr([u_2,n_2])$.
We have 

$$-\frac{v^{14}}{\Phi_6(\alpha)}D_1=v^{24}-2v^{22}+4v^{20}-8v^{18}+9v^{16}-12v^{14}+14v^{12}-12v^{10}+9v^8
-8v^6+4v^4-2v^2+1,$$

$$-\frac{v^{10}}{\Phi_6(\alpha)\Phi_8(\alpha)}E_1=v^8-v^6-v^4-v^2+1.$$

We have $D_1=\Phi(E_1)$ or $D_1=-\Phi(E_1)$. If one of them vanishes, then they both vanish because $\Phi$ is an automorphism. This proves it is sufficient to show that they cannot vanish at the same time to prove that neither of them vanishes. Assume $D_1=E_1=0$. We then have that the following quantities vanish

\begin{eqnarray*}
v_1 = & \op{Rem}(-\frac{v^{14}}{\Phi_6(\alpha)}D_1,-\frac{v^{10}}{\Phi_6(\alpha)\Phi_8(\alpha)}E_1) & = 20v^6+10v^4+5v^2-10\\
v_2 = & 2^3\op{Rem}(-\frac{v^{10}}{\Phi_6(\alpha)\Phi_8(\alpha)}E_1,v_1) & = -4v^4-v^2+2\\
v_3 = & \frac{2^2}{5}\op{Rem}(v_1,v_2) & = 11v^2-6\\
v_4 = & \frac{11^2}{2^5}\op{Rem}(v_2,v_3) & =1.
\end{eqnarray*}

For any $p\neq 2$, this proves that $0=1$ which is absurd. It follows that neither of them vanishes, therefore $\chi([u,m'])=\frac{\Phi(E_1)}{D_1}$. We have $D_1=\tr([u_1^{-1},n_1])$ and $E_1=\tr([u_2^{-1},n_2])$, therefore $\chi([u_1,n_1])=\chi([u_1^{-1},n_1])$ and $\chi([u_1,n_1][u_1^{-1}n_1])=1$.
We let $D_2=\tr(\chi([u_1,n_1][u_1^{-1}n_1])$ and $E_2=\tr(\chi([u_2,n_2][u_2^{-1}n_2])$. We have $D_2=\Phi(E_2)$.

Let $D_3=\tr([u_1,n_1^2])$ and $E_3=[u_2,n_2^2]$. We have $D_3=\Phi(E_3)$ or $D_3=-\Phi(E_3)$, therefore if we show that they cannot vanish simultaneously, then we have that neither of them vanishes. Assume by contradiction that $D_3=E_3=0$. We then use the Euclidean algorithm to prove that $1=0$ modifying it slightly using the conditions on $\alpha$ to dividing by non-zero quantities. In order to complete the algorithm, we need to invert the primes $7$, $13$, $17$, $19$, $23$, $47$, $53$, $193$, $599$, $881$, $1471$, $2503$, $3559$, $13967$, $44101$, $180811$, $382843$ and $981391$. See subsection \ref{computationD3E3} for the detailed computation. We check using GAP4 that those polynomials remain coprime when $p$ is a prime in the previous list. This is true except for $p=599$.

For $p=599$, we have $\op{Gcd}(-\frac{v^{22}}{\Phi_6(\alpha)}D_3,-\frac{v^{18}}{\Phi_6(\alpha)}E_3)=v^4+394v^2+1=0$. Let $v_{59901}=\frac{v^4+394v^2+1}{v^2}=0$. We then have

\begin{tiny}
\begin{eqnarray*}
D_1  & = & D_1-(-(v^{12}+v^{-12})+397(v^{10}+v^{-10})-156424(v^8+v^{-8})+61630673(v^6+v^{-6})-24282328759(v^4+v^{-4})\\
& & +9567175900402(v^2+v^{-2})-3769443022429664)v_{59901}-2479367974099310\times 599)\\
& = & 160\\
&\in & \pm \Phi(E_1)\\
& \in & \pm \Phi(E_1-(-v^8-v^{-8}+396(v^6+v^{-6})-156024(v^4+v^{-4})+61473061(v^2+v^{-2})-24220230012)v_{59901}-15930964405\times 599)\\
& \in & \pm \Phi(15)\\
& \in & \pm 15. 
\end{eqnarray*}
\end{tiny}

We then have $145=0$ or $175=0$. This concludes the proof that $D_3$ and $E_3$ are non-zero.

\smallskip

Case $\mathbf{1}$ : $D_1=\Phi(E_1)$.

\smallskip

Case $\mathbf{1.1}$ : $D_3=\tr([u_1,n_1^2])=\Phi(\tr([u_2,n_2^2]))=\Phi(E_3)$.

We then have $D_4 = D_1^3+D_1D_2-4D_1D_3+D_3^2-2D_2=\Phi(E_4)$, where $E_4=E_1^3+E_1E_2-4E_1E_3+E_3^2-2E_2$. We have

\begin{tiny}
$$0=\frac{v^{46}(E_1^5-4E_1^3E_3+36E_1^3+7E_1^2E_3+6E_1E_3^2-58E_1^2+18E_1E_2
-116E_1E_3-6E_1E_4+5E_2E_3)}{(v^2-1)^8(v^4+1)^{12}(v^6-1)^4\Phi_6(v^2)}$$
$$+\frac{v^{46}(E_2E_4-4E_3E_4+76E_1-42E_2+50E_3-10E_4-36)}{(v^2-1)^8(v^4+1)^{12}(v^6-1)^4\Phi_6(v^2)}.$$

\end{tiny}

It follows that the same expression, where the $E$ is replaced by $D$, is equal to zero but this is equal to 
\smallskip

\centerline{\fontsize{12}{12}\selectfont $\frac{(v^{12}-2v^{10}+4v^8-2v^6+4v^4-2v^2+1)(v^{16}-v^{14}+2v^{12}-3v^{10}+6v^8-3v^6+2v^4-v^2+1)(v^{16}-4v^{14}+8v^{12}-8v^{10}+7v^8-8v^6+8v^4-4v^2+1)}{v^{22}}$ }

This means that one of those three factors must vanish, we therefore treat those three cases separately.

\smallskip

Cases $\mathbf{1.1.1}$ : $D_5=\frac{v^{12}-2v^{10}+4v^8-2v^6+4v^4-2v^2+1}{v^6}=0$. We then define the following elements in order to get another polynomial in $v$ which will vanish and we will then prove that the two polynomials obtained this way cannot simultaneously vanish. Let us assume first that $p\neq 7$.

\begin{tiny}
\begin{eqnarray*}
D_{401} & = &  -\frac{1}{7}(D_1-(-(v^8+v^{-8})+(v^6+v^{-6})-(v^4+v^{-4})+6(v^2+v^{-2})+1)D_5)\\
& = &\frac{v^8+v^6+2v^4+v^2+1}{v^4}\\
  & = &-\frac{1}{7}\Phi(E_1)\\
D_{402}  & = &\frac{1}{2}(D_2-(v^{24}+v^{-24}-2(v^{22}+v^{-22}+4(v^{20}+v^{-20})-12(v^{18}+v^{-18}) +14(v^{16}+v^{-16})-22(v^{14}+v^{-14})\\
& & +52(v^{12}+v^{-12})-30(v^{10}+v^{-10})+55(v^8+v^{-8})-158(v^6+v^{-6})-56(v^4+v^{-4})-118(v^2+v^{-2})+508)D_5)\\
& = & \frac{479v^8-598v^6+311v^4-598v^2+479}{v^4}\\
& = & \frac{1}{2}\Phi(E_2)\\
D_{403} & = &  D_3-(\frac{-(v^{32}+1)+(v^{30}+v^{2})-(v^{28}+v^{4})+6(v^{26}+v^{6})-(v^{24}+v^{8})+v^{22}+v^{10}-24(v^{20}+v^{12})-17(v^{18}+v^{14})+v^{16}}{v^{16}}D_5)\\
& =& \frac{121v^8-55v^6+146v^4-55v^2+121}{v^4}\\
& = & \Phi(E_3)\\
D_{201} & = & 479D_{401}-D_{402}\\
& = & \frac{1077v^4+647v^2+1077}{v^2 }\\
& = & \Phi(-\frac{479}{7}E_1-\frac{1}{2}E_2)\\
D_{202}  & = & \frac{1}{2^4}(121D_{401}-D_{403})\\
 & = & \frac{11v^4+6v^2+11}{v^2}\\
 & = & \Phi(-\frac{121}{112}E_1-\frac{1}{16}E_3)\\
 D_{101} & = & 7^1 2^4 v^{26} (1077D_{202}-11D_{201}+655)\\
  & =& 0\\
  & =& -v^{26}(46013 E_1+7539E_3-616E_2-73360)\\
  & = & 616(v^{52}+1)-1232(v^{50}+v^2)+616(v^{48}+v^4)-1232(v^{46}+v^6)+11235(v^{44}+v^8)-20006(v^{42}+v^{10})+11851(v^{40}+v^{12})\\
  & & -12467(v^{38}+v^{14})+84793(v^{36}+v^{16}) -148348(v^{34}+v^{18})+86025(v^{32}+v^{20})-77254(v^{30}+v^{22})+148964(v^{28}+v^{24})\\
  & & -199634v^{26}
\end{eqnarray*}
\end{tiny}

We then let $F_1=v^6D_5$ and $F_2=-v^{26}(46013 E_1+7539E_3-616E_2-73360)$. Those two polynomials in $v$ vanish. Let us prove that this is absurd using the Euclidean algorithm.
\begin{tiny}
\begin{eqnarray*}
F_3  & = & \frac{1}{2^4}\op{Rem}(F_1,F_2) \\
& =&  -1560220v^{10}-984236v^{8}-1419140v^{6}-109343v^{4}-108563v^{2}-272023\\
F_4 & = & 2^2 5^2 181^2 431^2 \op{Rem}(F_2, F_3) \\
&  = & 2890736267584v^8+196484466595v^6+2504145304402v^4
-1211843170318v^2+887713181987\\
F_5&  = & \frac{2^{10} 11^2 37^2 41^2 53^2 51071^2}{5^2 181^2 431^2} \op{Rem}(F_3, F_4)\\
 & =& -108261500375317071v^6-36647807575597274v^4 +33134576659475654v^2-32149953533188943\\
F_6 & = & \frac{3^2 101^2 349^2 1023778455893^2}{2^{10} 11^2 37^2 41^2 53^2 51071^2}\op{Rem}(F_4,F_5)\\
& =&  20498163516987363083353v^4-12957989741793097471291v^2+6283381428888280265305\\
F_7 & = & \frac{20498163516987363083353^2}{2^3 3^2 7^2 101^2 139^1 349^2 262400449^1 1023778455893^2}\op{Rem}(F5, F6) \\
&= & -275631295456v^2+156402674859\\
F_8 & = & \frac{2^10 11^2 47^2 16660499^2}{20498163516987363083353^2}\op{Rem}(F6, F7) \\
& = & 1.
\end{eqnarray*}
\end{tiny}

This leads to the contradiction $1=0$ if $p$ is different from $7$, $11$, $37$, $41$, $47$, $53$, $101$, $139$, $181$, $349$, $431$, $51071$, $16660499$, $262400449$, $1023778455893$ or $20498163516987363083353$. We are under the assumption that $p\neq 7$, therefore we first check the other primes. Using GAP4, we see that $F_1$ and $F_2$ are coprime unless $p\in \{139,262400449\}$.

For $p=139$, we have $\op{Gcd}(F_1,F_2)=v^4+82v^2+1$. Let then $F_{13901}=v^2+v^{-2}+82=0$. We have

\begin{tiny}
\begin{eqnarray*}
D_1 & = & D_1-(-(v^{12}+v^{-12})+85(v^{10}+v^{-10})-6976(v^8+v^{-8})+571961(v^6+v^{-6})-46893847(v^4+v^{-4})\\
 & & +3844723522(v^2+v^{-2})-315220434992)F_{13901}-185902059153\times 139\\
& = & 71\\
& = & \Phi(E_1)\\
& = & \Phi(E_1-(-(v^8+v^{-8})+84(v^6+v^{-6})-6888(v^4+v^{-4})+564733(v^2+v^{-2})-46301220)F_{13901}-27306263\times 139)\\
& = & \Phi(21)\\
& = & 21.
\end{eqnarray*}
\end{tiny}

This implies $71=21$, therefore $50=0$ which is a contradiction since $p=139$.

For $p=262400449$, we have $\op{Gcd}(F_1,F_2)=v^4+144711873 v^2+1$. Let then $F_{26240044901}=v^2+v^{-2}+144711873=0$. We have

\begin{tiny}
\begin{eqnarray*}
D_1 & = & D_1-(-(v^{12}+v^{-12}+144711876(v^{10}+v^{-10})-20941526621303754(v^8+v^{-8})+3030487540848227798559380(v^6+v^{-6})\\
& & -438547528139311032518667960215007(v^4+v^{-4})+63463034196559901414852807139575347118760(v^2+v^{-2})\\
& & -9183854544847233051891171601164688387512952822508)F_{26240044901}
\\
& & -5064826670873590871335015289738121630701687311347\times 262400449\\
& = & 19405199\\
& = & \Phi(E_1)\\
& = & \Phi(E_1-(-(v^8+v^{-8})+144711875(v^6+v^{-6})-20941526476591875(v^4+v^{-4})
+3030487519906700743120001(v^2+v^{-2})\\
& & -
438547525108823428845860731880000)F_{26240044901}
-241855659926909504463506824998245\times 262400449)\\
& = & \Phi(212787997)\\
& = & 212787997.
\end{eqnarray*}
\end{tiny}

This implies $212787997=19405199$, therefore $193382798=0$ which is a contradiction since $p=262400449$.

The only remaining prime is $p=7$. We have 

\begin{tiny}
\begin{eqnarray*}
D_1 & = & D_1-(-(v^8+v^{-8})+(v^6+v^{-6})-(v^4+v^{-4})+6(v^2+v^{-2})+1)D_5\\
& = & -\frac{7(v^8+v^6+2v^4+v^2+1)}{v^4}\\
& =& 0.
\end{eqnarray*}
\end{tiny}

This  leads to a contradiction because we have proven that $D_1\neq 0$. This concludes the proof of case \textbf{1.1.1}.

Case $\mathbf{1.1.2}$ : $D_6=\frac{v^{16}-v^{14}+2v^{12}-3v^{10}+6v^8-3v^6+2v^4-v^2+1}{v^8}=0$.

We then have 

\begin{tiny}
\begin{eqnarray*}
D_{601} & = & D_1-(-(v^6+v^{-6})+2(v^4+v^{-4})-3(v^2+v^{-2})+4)D_6\\
& = & \frac{2(2v^{12}-4v^{10}+6v^8-7v^6+6v^4-4v^2+2)}{v^6}\\
& = & \Phi(E_1)\\
D_{602} & = & D_2-(v^{22}+v^{-22}-3(v^{20}+v^{-20})+7(v^{18}+v^{-18})-14(v^{16}+v^{-16})+19(v^{14}+v^{-14})
-27(v^{12}+v^{-12})+37(v^{10}+v^{-10})\\
& & -48(v^8+v^{-8})+74(v^6+v^{-6})-82(v^4+v^{-4})+70(v^2+v^{-2})-52)D_6\\
& = & -\frac{2(24v^{12}-16v^{10}-24v^8+63v^6-24v^4-16v^2+24)}{v^6}\\
& = & \Phi(E_2)\\
D_{603} & = & D_3-(-(v^{14}+v^{-14})+2(v^{12}+v^{-12})-3(v^{10}+v^{-10})+4(v^8+v^{-8})-(v^6+v^{-6})+2(v^4+v^{-4})-7(v^2+v^{-2})+12)D_6\\
& = & -\frac{2(6v^{12}-4v^{10}-6v^8+15v^6-6v^4-4v^2+6)}{v^6}\\
& = & \Phi(E_3)\\
D_{404} &= &  -12D_{601}-D_{602}\\
& = & \frac{2(32v^8-96v^6+147v^4-96v^2+32)}{v^4}\\
& = & \Phi(-12E_1-E_2)\\
D_{405} & = & -3D_{601}-D_{603}\\
& =&  \frac{8(2v^8-6v^6+9v^4-6v^2+2)}{v^4}\\
& = & \Phi(-3E_1-E_3)\\
D_{102} & = & D_{404}-4D_{405}-6\\
& = & 0\\
& = & \Phi(-E_2+4E_3-6)\\
& =& -\Phi(\frac{(v^6-1)(v^{10}-1)(v^8+1)^2(v^2-1)^2(v^8-v^5+v^3+1)(v^8+v^5-v^3+1)}{v^{26}})
\end{eqnarray*}
\end{tiny}

Since $\Phi$ is an automorphism, we have that $G_2 = (v^8-v^5+v^3+1)(v^8+v^5-v^3+1)=0$. We also have by assumption $G_1=v^6D_6=0$. The Euclidean remainder of the division of $G_1$ by $G_2$ as polynomials in $v$ is equal to $-v^2(v^6-1)\Phi_6(v^2)(v^2-1)^2$. This leads to a contradiction by the assumptions on $\alpha$ since $\alpha=v^2$.

\smallskip

Case $\mathbf{1.1.3}$ : $D_7=\frac{v^{16}-4v^{14}+8v^{12}-8v^{10}+7v^8-8v^6+8v^4-4v^2+1}{v^8}=0$.

Assume first $p\notin \{11,13\}$. We now define new quantities which will permit us to find polynomials which vanish in $v$.

\begin{tiny}
\begin{eqnarray*}
D_{604}  & = & D_1-(-v^6-v^{-6}-v^4-v^{-4}-3(v^2+v^{-2})+2)D_7\\
&= & -\frac{13v^{12}-23v^{10}+23v^8-16v^6+23v^4-23v^2+13}{v^6}\\
& = & \Phi(E_1)\\
D_{605} & = & D_2-(v^{22}+v^{-22}+4(v^{18}+v^{-18})-6(v^{16}+v^{-16})-v^{14}-v^{-14}-32(v^{12}+v^{-12})-7(v^{10}+v^{-10})-10(v^8+v^{-8})\\
& & +137(v^6+v^{-6})+232(v^4+v^{-4})+331(v^2+v^{-2})-158)D-7\\
& = & -\frac{476v^{12}-751v^{10}+838v^8-645v^6+838v^4-751v^2+476}{v^6}\\
& = & \Phi(E_2)\\
D_{606} & = & D_3-(-v^{14}-v^{-14}-v^{12}-v^{-12}-3(v^{10}+v^{-10})+2(v^8+v^{-8})+8(v^6+v^{-6})+26(v^4+v^{-4})+28(v^2+v^{-2})-2)D_7\\
& = & -\frac{149v^{12}-271v^{10}+283v^8-208v^6+283v^4-271v^2+149}{v^6}\\
& = & \Phi(E_3)\\
D_{607} & = & D_1^2-(v^{20}+v^{-20}-2(v^{18}+v^{-18})+7(v^{16}+v^{-16})-18(v^{14}+v^{-14})+24(v^{12}+v^{-12})-62(v^{10}+v^{-10})+77(v^8+v^{-8})\\
& & -82(v^6+v^{-6})+224(v^4+v^{-4})-30(v^2+v^{-2})+207)D_7\\
& = & -\frac{476v^{12}-751v^{10}+838v^8-645v^6+838v^4-751v^2+476}{v^6}\\
& = & \Phi(E_1^2)\\
D_{406} & = & 1697D_{604}+13D_{605}\\
& = & \frac{3245v^8-4025v^6+3736v^4-4025v^2+3245}{v^4}\\
& = & \Phi(1697 E_1+13E_2)\\
D_{407} & = & \frac{1}{4}(149 D_{604}+13D_{606})\\
& = & \frac{24v^8-63v^6+80v^4-63v^2+24}{v^4}\\
& = & \frac{1}{4}\Phi(149 E_1+13 E_3)\\
D_{408} & = & -13D_{605}-476D_{604}\\
& = & \frac{1185v^8-54v^6-769v^4-54v^2+1185}{v^4}\\
& = & \Phi(-13 E_1^2-476 E_1)\\
D_{203} & = & \frac{1}{13}(24D_{406}-3245 D_{407})\\
& = & \frac{8295v^4-13072v^2+8295}{v^2}\\
& = & \Phi(-\frac{24661}{4}E_1+24E_2-\frac{3245}{4}E_3)\\
D_{204} & = & \frac{1}{3^2 11^1 13^1}(-1185 D_{407}+24D_{408})\\
& = & \frac{57v^4-88v^2+57}{v^2}\\
& = & \Phi(-\frac{395}{2^2 3^1 11^1}E_3-\frac{5699}{2^2 3^1 11^1}E_1-\frac{8}{3^1 11^1}E_1^2)\\
D_{103}& =& 57D_{203}-8295D_{204}+15144\\
 & = & 0\\
 & = & \Phi(\frac{73822}{11}E_1+1368E_2-\frac{235610}{11}E_3+\frac{22120}{11}E_1^2+15144)\\
 & = & \Phi(\frac{2}{11v^{26}}(7524v^{36}-45144v^{34}+127908v^{32}-214660v^{30}+262569v^{28}-329886v^{26}
 +445989v^{24}-513465v^{22}+494686v^{20}\\
 & & -471434v^{18}+494686v^{16}
 -513465v^{14}+445989v^{12}-329886v^{10}+262569v^8-214660v^6
 +127908v^4-45144v^2+7524)\\
 & & (v^{16}+4v^{14}+8v^{12}+8v^{10}
 +7v^8+8v^6+8v^4+4v^2+1)).
\end{eqnarray*}
\end{tiny}

It follows that one of the two factors inside the last expression vanishes. We separate the two possibilities in order to make the computations easier.

\smallskip

Case $\mathbf{1.1.3.1}$ : $H_2=v^{16}+4v^{14}+8v^{12}+8v^{10}
 +7v^8+8v^6+8v^4+4v^2+1=0$. We set $H_1=v^6D_7=0$. We then use the Euclidean algorithm to obtain a contradiction.

 \begin{eqnarray*}
 H_3 & = -\frac{1}{8}\op{Rem}(H_1, H_2) & = v^{14}+2v^{10}+2v^6+v^2\\
H_4 & = \op{Rem}(H_2, H_3) & = 6v^{12}+5v^8+7v^4+1\\          
H_5 & = 6\op{Rem}(H_3, H_4) &=  7v^{10}+5v^6+5v^2\\
H_6 & =   7\op{Rem}(H_4, H_5) & = 5v^8+19v^4+7\\
H_7 & =  5\op{Rem}(H_5, H_6) & = -108v^6-24v^2\\
H_8 & =  9\op{Rem}(H_6, H_7, v) & = 161v^4+63\\
H_9 & =  23\op{Rem}(H_7, H_8, v) & = 420v^2\\
H_{10} & =  \frac{1}{3^2 7^1}\op{Rem}(H_8, H_9, v) & =  1.
 \end{eqnarray*}

 For $p=23$, we get $H_8=63=17\neq 0$, therefore for any prime different from $7$, we get a contradiction. Using GAP4, we get for $p=7$ that $\op{Gcd}(H_1,H_2)=v^4+1$, which is also a contradiction.
 
 \smallskip
 
 Case $\mathbf{1.1.3.2}$ : \begin{tiny}$I_2=7524(v^{36}+1)-45144(v^{34}+v^2)+127908(v^{32}+v^4)-214660(v^{30}+v^6)+262569(v^{28}+v^8)
 -329886(v^{26}+v^{10})+445989(v^{24}+v^{12})-513465(v^{22}+v^{14})+494686(v^{20}+v^{16})-471434v^{18}
=0$. \end{tiny} We set again $I_1=v^6 D_7$. We then prove that $0=1$ for most primes using the Euclidean algorithm
\begin{tiny}
\begin{eqnarray*}
I_3 & = &  \op{Rem}(I_2, I_1)\\
& = & 214002v^{14}-916685v^{12}+1035136v^{10}-754561v^8+879809v^6
-1139783v^4+640260v^2-178015\\
I_4 & = &  2^2 3^8 1321^2\op{Rem}(I_1, I_3)\\
& = & 200475369505v^{12}-267706231982v^{10}+178081604207v^8
-175843177159v^6+298516540603v^4-183940914006v^2+56598272159\\
I_5 & = &  \frac{5^2 137^2 90023^2 3251^2}{2^2 3^8 1321^2}\op{Rem}(I_3, I_4)\\
& = & 2229625335992657v^{10}-5627663796526377v^8+6805441741717849v^6
-3483903150258508v^4+844734325177536v^2\\
& & +92479019896356\\
I_6 & = &  \frac{31^2 389471963^2 184669^2}{5^2 137^2 90023^2 3251^2}\op{Rem}(I_4, I_5)\\
& = & 20737893552210092217v^8-72972566540716005022v^6+73587018705359365791v^4
-34948076592200136378v^2\\
& & +5778167564925374963\\
I_7 & = &  \frac{3^2 23^2 29^2 10363764893658217^2}{7^1 11^1 31^2 389471963^2 184669^2}\op{Rem}(I_5, I_6)\\
& = & 7525561044044520343864v^6-8534947218256775108025v^4
+4450478132504667710979v^2-590406861955024796687\\
I_8 & = &  \frac{2^6 1497344053^2628242472811^2}{3^2 23^2 29^2 10363764893658217^2}\op{Rem}(I_6, I_7)\\
& = & 689630866150855463473281v^4-536688784588813135219275v^2
+249996942573872461979743\\
I_9 & = &  \frac{3^1 79^2 599418320761723^2 4854431^2}{ 2^7 103^1 159073628242472811^2 1497344053^290023}\op{Rem}(I_7, I_8)\\
& = & -343470144464v^2+361073206339\\
I_{10} & = & \frac{2^8 21466884029^2}{79^2 599418320761723^2 4854431^2}\op{Rem}(I_8,I_9)\\
& = & 1.
\end{eqnarray*}
\end{tiny}

This leads to a contradiction if $p$ is different from $7$, $23$, $29$, $31$, $79$, $103$, $137$, $1321$, $3251$, $90023$, $159073$, $184669$, $4854431$, $10363764893658217$, $1497344053$, $21466884029$, 
$389471963$, $599418320761723$ and $628242472811$. We check using GAP4 that $I_1$ and $I_2$ remain coprime except for $p\in \{7,103,90023,159073\}$.

\smallskip

For $p=7$, we have $\op{Gcd}(I_1,I_2)=v^8+v^6+4v^4+v^2+1=(v^2+v-1)^2(v^2-v-1)^2$. 

If $v^2 =1-v$, then $v^4=1+v^2-2v=2-3v$, $v^8=4+9v^2-12v=13-21v=13=-1$ and $v^{16}=1$, which contradicts the fact that $\alpha$ is of order not dividing $8$. 

If $v^2=v+1$, then $v^4=v^2+2v+1=2+3v$, $v^8=4+9v^2+12v=13+21v=13=-1$, $v^{16}=1$, which is a contradiction.

\smallskip

For $p=103$, we have $\op{Gcd}(I_1,I_2)=v^4+60v^2+1$. Let then $I_{10301}=\frac{v^4+60v^2+1}{v^2}$. We then have

\begin{tiny}
\begin{eqnarray*}
D_1 & = & D_1-(-v^{12}-v^{-12}+63(v^{10}+v^{-10})-3786(v^8+v^{-8})+227111(v^6+v^{-6})-13622895(v^4+v^{-4})+817146618(v^2+v^{-2})\\
& & -49015174220)I_{10301}-28536661747\times 103\\
& = & 61\\
& = & \Phi(E_1)\\
& = & \Phi(E_1-(-v^8-v^{-8}+62(v^6+v^{-6})-3720(v^4+v^{-4})+223139(v^2+v^{-2})-26692311)I_{10301}-7792534\times 103)\\
& = &\Phi(44)\\
& = & 44.
\end{eqnarray*}
\end{tiny}

This would imply $17=0$ which is absurd since $p=103$.

\smallskip

For $p=90023$, we have $\op{Gcd}(I_1,I_2)=v^4+18030v^2+1$. Let then $I_{9002301}=\frac{v^4+18030v^2+1}{v^2}$. We then have

\begin{tiny}
\begin{eqnarray*}
D_1 & = & D_1-(-v^{12}-v^{-12}+18033(v^{10}+v^{-10})-325134996(v^8+v^{-8})+5862183959861(v^6+v^{-6})-105695176471158855(v^4+v^{-4})\\
& & +1905684025912810195818(v^2+v^{-2})
-34359482881512791359439720)I_{9002301}-6881591066086528735823931\times 90023\\
& = & 19589\\
& = & \Phi(E_1)\\
& = & \Phi(E_1-(-v^8-v^{-8}+18032(v^6+v^{-6})-325116960(v^4+v^{-4})+5861858770769(v^2+v^{-2})-105689313311848112)I_{9002301}\\
& & -21167682784276284\times 90023)\\
& = & \Phi(3294)\\
 & = & 3294.
\end{eqnarray*}
\end{tiny}

This would imply that $16295=0$ which is absurd since $p=90023$.

\smallskip

For $p=159073$, we have $\op{Gcd}(I_1,I_2)=v^4+69018v^2+1$. Let then $I_{15907301}=\frac{v^4+69018v^2+1}{v^2}$. We then have

\begin{tiny}
\begin{eqnarray*}
D_1 & = & D_1-(-v^{12}-v^{-12}+69021(v^{10}+v^{-10})-4763691384(v^8+v^{-8})+328780451871905(v^6+v^{-6})\\
& & -22691769222531447927(v^4+v^{-4})+1566140527871895021153810(v^2+v^{-2})\\
& & -108091886929970681347462210688)I_{15907301}
-46898504768253791840823752710\times 159073\\
 & = & 118972\\
 & = & \Phi(E_1)\\
 & = & \Phi(E_1-(-v^8-v^{-8}+69020(v^6+v^{-6})-4763622360(v^4+v^{-4})+328775687973461(v^2+v^{-2})\\
 & & -22691440427788708940)I_{15907301}
 -9845277544193984759\times 159073)\\
 & = & \Phi(105595)\\
 & =& 105595.
\end{eqnarray*}
\end{tiny}

This would imply that $13377=0$, which is absurd since $p=159073$.

\smallskip

It now only remains to consider $p\in \{11,13\}$ to conclude Case \textbf{1.1.3}.

Assume first $p=13$. We have, using the same notations as before, that $24D_{406}-3245D_{407}=0=\Phi(24 E_1-3245 E_2)$. Therefore $H_{1301}=v^{26}(24 E_1-3245 E_2)=0$. We let $H_{1302}=v^8 D_7=0$. We then have using GAP4 that $\op{Gcd}(H_{1301},H_{1302})=1$ which leads to a contradiction.

Assume now $p=11$. We have, using the same notations as before, that $-1185D_{407}+24D_{408})=0=\Phi(-1185 E_2+24E_3)$. Therefore $H_{1101}=v^{26}(-1185E_2+24E_3)=0$. We let $H_{1102}=v^8 D_7=0$. We then have using GAP4 that $\op{Gcd}(H_{1101},H_{1102})=1$ which leads to a contradiction.

\smallskip

Case $\mathbf{1.2}$ : $D_3=\tr([u_1,n_1^2])=-\Phi(\tr([u_2,n_2^2]))=-\Phi(E_3)$.

We have $D_3 E_3 D_1 E_1\neq 0$, therefore in this case, we have $\chi([u',m'^2])=-1$ and $\chi([u',m'])=1$. Therefore $\chi([u',m'][u',m'^2])=-1$. It follows that 
$$D_8=\tr([u_1,n_1][u_1,n_1^2])=-\Phi(\tr([u_2,n_2][u_2,n_2^2]))=-\Phi(E_8).$$

We then have $0=D_1 D_3-D_8-D_1=\Phi(-E_1E_3+E_8-E_1)=2\Phi_(E_1)$. It follows that $E_1=0$, which is absurd by what was proven before case $\mathbf{1}$.

\smallskip

Case $\mathbf{2}$ : $D_1=-\Phi(E_1)$.

\smallskip

Case $\mathbf{2.1}$ : $D_3=\Phi(E_3)$.

This is the worst case in terms of computations, we will not write all the polynomials appearing because their degree can be very high. We will give all the necessary elements to consider and the final results of the computations. We have $2-E_1=\frac{(v^{10}-1)(v^6-1)(v^2-1)^2}{v^{10}}\neq 0$, therefore $D_1+2=\Phi^{-1}(2-E_1)\neq 0$. The elements defined in subsection \ref{computationscase21} of the Appendix and the corresponding computations then lead us to the contradiction $1=0$ except for the primes which we needed to invert in order to do the computations. We check using GAP4 that $v^{86}D_{13}$ and $v^{70}D_{12}$ are coprime for the primes $p$ we inverted in the computations. It is true except for $p\in \{7,17,43,79,1013,1747\}$.
     
     \smallskip
     
     For $p=7$, we have $\op{Gcd}(v^{86}D_{13},v^{70}D_{12})=(v^2+v-1)^3(v^2-v-1)^3(v^8+v^6-v^4+v^2+1)$, therefore $v^8+v^6-v^4+v^2+1=0$, because $(v^2+v-1)(v^2-v-1)=0$ would imply $v^{16 }=1$ as we have seen in case \textbf{1.1.3.2}. Let now $K_{701}=\frac{v^8+v^6-v^4+v^2+1}{v^4}$. We then have
     
     \begin{tiny}
    $$ D_1 = D_1-(-(v^{10}+v^{-10}+4(v^8+v^{-8})-12(v^6+v^{-6})+31(v^4+v^{-4})-67(v^2+v^{-2})+135)K_{701} =-\frac{7(27v^4-35v^2+27)}{v^2}=0.$$
    \end{tiny}
    
This is absurd because we proved that we cannot have $D_1=0$.

\smallskip

For $p=17$, we have $\op{Gcd}(v^{86}D_{13},v^{70}D_{12})=\Phi_6(\alpha)\neq 0$, which is a contradiction.

\smallskip

For $p=43$, we have $\op{Gcd}(v^{86}D_{13},v^{70}D_{12})=v^4+10v^2+1.$ We then set $K_{4301}=\frac{v^4+10v^2+1}{v^2}$. We have

\begin{tiny}
\begin{eqnarray*}
D_1 & = &  D_1-(-(v^{12}+v^{-12}+13(v^{10}+v^{-10})-136(v^8+v^{-8})+1361(v^6+v^{-6})-13495(v^4+v^{-4})\\
& & + 133618(v^2+v^{-2})-1322720)K_{4301}-301395\times 43\\
& = & 17\\
& =& -\Phi(E_1)\\
 & = & -\Phi(E_1-(-(v^8+v^{-8})+12(v^6+v^{-6})-120(v^4+v^{-4})+1189(v^2+v^{-2})-11772)K_{4301}-2682\times 43)\\
 & = & -\Phi(20)\\
 & = & -20.
\end{eqnarray*}
\end{tiny}

This implies that $37=0$ which is a contradiction because $p=43$.

\smallskip

For $p=79$, we have $\op{Gcd}(v^{86}D_{13},v^{70}D_{12})=v^4+4v^2+1$. Let $K_{7901}=\frac{v^4+4v^2+1}{v^2}$, we then have a contradiction because

\begin{tiny}
$$D_1= D_1-(-(v^{12}+v^{-12})+7(v^{10}+v^{-10})-34(v^8+v^{-8})+143(v^6+v^{-6})-559(v^4+v^{-4})+2122(v^2+v^{-2})-7964)K_{7901}-350\times 79=0.$$
\end{tiny}

For $p=1013$, we have $\op{Gcd}(v^{86}D_{13},v^{70}D_{12})=v^4+179v^2+1$. Let $K_{101301}=\frac{v^4+179v^2+1}{v^2}$, we have 

\begin{tiny}
\begin{eqnarray*}
D_1 & = & D_1-(-(v^{12}+v^{-12})+182(v^{10}+v^{-10})-32584(v^8+v^{-8})+5832368(v^6+v^{-6})-1043961309(v^4+v^{-4})\\
& & +186863241972(v^2+v^{-2})-33447476351714)K_{101301}-5909895893852\times 1013\\
& = & 824\\
& = & -\Phi(E_1)\\
& = & -\Phi(E_1-(-(v^8+v^{-8})+181(v^6+v^{-6})-32399(v^4+v^{-4})+5799241(v^2+v^{-2})-1038031742)K_{101301}-183411730\times 1013)\\
& = & -\Phi(850)\\
& = & 163.
\end{eqnarray*}
\end{tiny}

This would imply that $661=0$ which is absurd because $p=1013$.

\smallskip

For $p=1747$, we have $\op{Gcd}(v^{86}D_{13},v^{70}D_{12})=v^4+482v^2+1$. Let $K_{174701}=\frac{v^4+482v^2+1}{v^2}$, we have 

\begin{tiny}
\begin{eqnarray*}
D_1 & =& D_1-(-v^{12}-v^{-12}+485(v^{10}+v^{-10})-233776(v^8+v^{-8})+112679561(v^6+v^{-6})-54311314647(v^4+v^{-4})\\
& & +26177940980322(v^2+v^{-2})-12617713241200592)K_{174701}-3481216614983814\times 1747\\
&= & 1680\\
& = & -\Phi(E_1)\\
& = & -\Phi(E_1-(-(v^8+v^{-8})+484(v^6+v^{-6})-233288(v^4+v^{-4})+112444333(v^2+v^{-2})-54197935220)K_{174701}-14953165361\times 1747)\\
& = & -\Phi(1711)\\
& = & 36
\end{eqnarray*}
\end{tiny}

This would imply $1644=0$ which is absurd because $p=1747$. This shows that Case $\mathbf{2.1}$ is absurd.

\smallskip

Case $\mathbf{2.2}$ : $D_3=-\Phi(E_3)$.

We then have $D_8=\tr([u_1,n_1][u_1,n_1^2])=\Phi(E_8)=\Phi(\tr([u_2,n_2][u_2,n_2^2])$. We have as in case $\mathbf{1.2}$, $D_1D_3-D_8-D_1=0=\Phi(E_1E_3-E_8+E_1)$ and $E_1E_3-E_8-E_1=0$, therefore $\Phi(2E_1)=0$ and $E_1=0$ which is absurd.

\medskip

We have proven $(\pi_1\times \pi_2)([N_1',N_1'])= SL_2(q)\times SL_2(q)$. We now want to determine $(\pi_1\times \pi_2)(N_1')$. Since $\op{det}(u_1)=\op{det}(u_2)=\op{det}(n_1)=1$ and $\op{det}(n_2)=\frac{1}{v^8}$, we have that $det( (\pi_1\times \pi_2)(N_1'))\simeq <v^8>\simeq \Z/\frac{o(\alpha)}{\op{Gcd}(o(\alpha),4)}\Z$, where $o(\alpha)$ is the order of $\alpha$. Since $(\pi_1\times \pi_2)(N_1')\leq SL_2(q)\times GL_2(q)$, we have that the kernel of the determinant is contained in $SL_2(q)\times SL_2(q)$ and since $(\pi_1 \times \pi_2)([N_1',N_1'])=SL_2(q)\times SL_2(q)$, we get that the kernel is equal to $SL_2(q)\times SL_2(q)$. This shows that it is a normal subgroup of $(\pi_1\times \pi_2)(N_1')$, the resulting quotient is a cyclic group of order $\frac{o(\alpha)}{\op{Gcd}(o(\alpha),4)}$. It follows that $N_1'\simeq (SL_2(q)\times SL_2(q)). \Z/\frac{o(\alpha)}{\op{Gcd}(o(\alpha),4)}\Z$, where this denotes an extension which may be split.

Recall now that $H=(XP)\rho_{8_r}(\mathcal{A}_{A_2})(XP)^{-1}\simeq SL_2(q)$, $N_1\simeq N_1'$ and $N_1\leq C_G(H)$. It follows that the order $\vert HN_1\vert $ is equal to $\frac{\vert H\vert \vert N_1\vert}{\vert H\cap N_1\vert}$. The intersection of $H$ and $N_1$ is at most $\Z/2\Z$. By noting that $R^{-1}\begin{pmatrix}-I_2 & 0 & 0 & 0\\
0 & I_2 & 0 & 0\\
0 & 0 & -I_2 & 0\\
0 & 0 &0 & I_2
\end{pmatrix}R=\begin{pmatrix}-I_2 & 0 & 0 & 0\\
0 & I_2 & 0 & 0\\
0 & 0 & -I_2 & 0\\
0 & 0 &0 & I_2
\end{pmatrix}\in N_1'$, we get that the intersection is exactly $\Z/2\Z$. It follows that the order of $HN_1$ is equal to $\frac{1}{2}\vert SL_2(q)^3\vert \frac{o(\alpha)}{\op{Gcd}(o(\alpha),4)}=\frac{1}{2\op{Gcd}(o(\alpha),4)}q^3(q^2-1)^3o(\alpha)$. Assume now that we are in the fourth case of Theorem \ref{theoKantor}, we would have $G\simeq 2^{\cdot}\Omega_7(q)$. It would follow that $\frac{1}{2\op{Gcd}(o(\alpha),4)}q^3(q^2-1)^3o(\alpha)$ divides $\frac{2q^9(q^2-1)(q^4-1)(q^6-1)}{\op{Gcd}(2,q-1)}=q^9(q^2-1)^3(q^2+1)(q^4+q^2+1)$. We would then have that $o(\alpha)$ divides $2\op{Gcd}(o(\alpha),4)(q^2+1)(q^4+q^2+1)$ which divides $8(q^2+1)(q^4+q^2+1)$. Since $o(\alpha)$ divides $q^2-1$, $q^2+1=q^2-1+2$ and $q^4+q^2+1=q^4-1+q^2-1+3$, we would have $o(\alpha)$ divides $48$ which contradicts our assumptions on $\alpha$. 

\bigskip

\bigskip

Assume now that $\F_{q^2}=\F_p(\sqrt{\alpha})\neq \F_p(\alpha)=\F_p(\alpha+\alpha^{-1})$. 
 We have, by Table \ref{resH4H3} and Proposition \ref{H3dim4}, that $\rho_{8_r}(\mathcal{A}_{H_3})\simeq SU_4(q)$. It follows that $G$ is again an irreducible group generated by long root elements. By Propositions \ref{color} and \ref{bilinwgraphs}, we have that $G\leq \Omega_8^+(q)$, up to conjugation in $GL_8(q^2)$. It follows by Theorem \ref{theoKantor} that $G$ is isomorphic to one of the following groups
 
 \begin{enumerate}
\item $\Omega_8^{+}(q)$
\item $\Omega_8^{-}(\sqrt{q})$
\item $SU_4(q)$
\item $G/Z(G)=P\Omega_7(q)$, $Z(G)=2$
\item $^3\!D_4(\sqrt[3]{q})$
\end{enumerate}

By Lemma 3.5 of \cite{BM}, we again have $\rho_{8_r}(\mathcal{A}_{A_2,1})\simeq SL_2(q)$, where $A_{A_2,1}=<S_2,S_3>$. This proves that we again have a subgroup of $G$ of order $\frac{1}{2\op{Gcd}(o(\alpha),4)}q^3(q^2-1)^3o(\alpha)$. This excludes the fourth case by the same arguments as before. We also have that $G$ contains a subgroup of order equal to $q^6(q^2-1)(q^3+1)(q^4-1)$.

\smallskip

Assume now by contradiction that $G\simeq \Omega_8^{-}(\sqrt{q})$. We then have $\vert G\vert =q^6(q^2+1)(q-1)(q^2-1)(q^3-1)$. It then follows that $(q^2-1)$ divides $8(q^2+1)(q^2-q+1)$. We then have that $q^2-1$ divides $8(q^2+1)(q^2-q+1)-8(q^2-1)(q^2-q+1)=16(q^2-q+1)$, therefore $q^2-1$ divides $-16(q^2-q+1)+16(q^2-1)=16(q-2)=16q-32$. We have $16q-32<q^2-1$ if $q\geq 16$, therefore $q\leq 15$. By the conditions on $\alpha$ and $p$, this is absurd.

\smallskip

Assume now by contradiction that $G\simeq ^3\!D_4(\sqrt[3]{q})$. We then have $\vert G\vert = \sqrt[3]{q}^{12}(\sqrt[3]{q}^{8}+\sqrt[3]{q}^{4}+1)(\sqrt[3]{q}^{6}-1)(\sqrt[3]{q}^{2}-1).$ This is absurd because $q^6$ does not divide $q^4$.

\smallskip

Assume now by contradiction that $G\simeq SU_4(q)$. We then have $\vert G\vert=q^6(q^2-1)(q^3+1)(q^4-1)$. It follows that $(q^2-1)$ divides $8(q^3+1)(q^2+1)$. This implies that $(q^2-1)$ divides $8(q^3+1)(q^2+1)-8(q^3+1)(q^2-1)=16(q^3+1)$. We then have that $q^2-1$ divides $16(q^3+1)-16(q^3-q)=16(q+1)$. It follows that $q-1$ divides $16$ which is absurd by the conditions on $\alpha$ and $p$.

\smallskip

This proves that we also have $G\simeq \Omega_8^+(q)$ when $\F_{q^2}=\F_p(\sqrt{\alpha})\neq \F_p(\alpha)=\F_p(\alpha+\alpha^{-1})$.

\bigskip

Assume now that $\F_q=\F_p(\sqrt{\alpha})=\F_p(\alpha)\neq \F_p(\alpha+\alpha^{-1})$. Let $\epsilon$ be the unique automorphism of order $2$ of $\F_q$. We have $\epsilon(\alpha)=\alpha^{-1}$. It follows by Proposition \ref{Fieldfactorization} that $\epsilon \circ \rho_{8_r|\mathcal{A}_{H_4}}\simeq \rho_{8_r|\mathcal{A}_{H_4}}$ or $\epsilon \circ \rho_{8_r|\mathcal{A}_{H_4}}\simeq \rho_{8_{rr}|\mathcal{A}_{H_4}}$. We have 
\begin{eqnarray}
\tr(\rho_{8_r}(S_1S_3S_2^{-1}S_4^{-1})) & = &  \alpha^2+\alpha^{-2}-5(\alpha+\alpha^{-1})+8\\
\epsilon(\tr(\rho_{8_r}(S_1S_3S_2^{-1}S_4^{-1}))) & =& \alpha^2+\alpha^{-2}-5(\alpha+\alpha^{-1})+8\\
\tr(\rho_{8_{rr}}(S_1S_3S_2^{-1}S_4^{-1})) & = & \alpha^2+\alpha^{-2}-5(\alpha+\alpha^{-1})+7
\end{eqnarray}
It follows that $\epsilon \circ \rho_{8_r|\mathcal{A}_{H_4}}\simeq \rho_{8_r|\mathcal{A}_{H_4}}$. We then have by Proposition 4.1 of \cite{BMM} that up to conjugation in $GL_8(q)$, we have $G\leq \Omega_8^+(q^{\frac{1}{2}})$. By Proposition \ref{H3dim4}, we have that $\rho_{8_r}(\mathcal{A}_{H_3})\simeq SU_4(q^{\frac{1}{2}})$ or $\rho_{8_r}(\mathcal{A}_{H_3})\simeq SL_4(q^{\frac{1}{2}})$. The computations in the first case show here that $G$ contains a subgroup of order divisible by $\frac{1}{2}\vert SL_2(q^{\frac{1}{2}})^3\vert \frac{o(\alpha)}{\op{Gcd}(o(\alpha),4)}$. By the same arguments as for the case $\F_p(\alpha)=\F_p(\alpha+\alpha^{-1})$, we have that $G\simeq \Omega_8^+(q^{\frac{1}{2}})$.
\end{proof}

\section{Triality automorphism and the two $8$-dimensional representations}\label{sectiontriality}

The two $8$-dimensional representations being linked through triality, we first recall a few facts about triality and we will then determine the image of the Artin group inside $\rho_{8_r}\times \rho_{8_{rr}}$. We use the following notations, definitions and results of \cite{CAR}. This phenomenon was also observed in the generic case in Proposition $6.7$ of \cite{IH2}.

\begin{Def}\label{deftriality}
Let $A=\begin{pmatrix} 0 & I_4\\
I_4 & 0\end{pmatrix}$. The simple Lie algebra $\mathfrak{L}$ of type $D_4$ over $\C$ is the matrix algebra $\{T\in \mathcal{M}_8(\C), ^t\!T A+AT=0\}$. We write $(E_{i,j})_{1\leq i,j\leq 4}$ the basis of elementary matrices of $\mathcal{M}_8(\C)$.

Its Cartan decomposition is given by $\mathfrak{L}=\mathfrak{h}\oplus \underset{r\in \Phi}\bigoplus \mathfrak{L}_r$, where $\Phi=\{\pm \epsilon_i\pm \epsilon_j, 1\leq i<j\leq 4\}$, for $1\leq i<j\leq 4$,
\begin{enumerate}
\item $e_{\epsilon_i-\epsilon_j}=E_{i,j}-E_{4+j,4+i}$,
\item $e_{\epsilon_j-\epsilon_i}=E_{4+i,4+j}+E_{j,i}$,
\item $e_{\epsilon_i+\epsilon_j}=E_{i,4+j}-E_{j,4+i}$,
\item $e_{-\epsilon_i-\epsilon_j}=-E_{4+i,j}+E_{4+j,i}$.
\end{enumerate}
$L_r=\C e_r$ and $\mathfrak{h}=\left\lbrace \begin{pmatrix}
\begin{pmatrix}
\lambda_1 & 0 & 0 & 0\\
0& \lambda_2 & 0 & 0\\
0 & 0 & \lambda_3 & 0\\
0 & 0 & 0 & \lambda_4
\end{pmatrix}  & 0 \\
0 & \begin{pmatrix}
-\lambda_1 & 0 & 0 & 0\\
0& -\lambda_2 & 0 & 0\\
0 & 0 & -\lambda_3 & 0\\
0 & 0 & 0 & -\lambda_4
\end{pmatrix}
\end{pmatrix},(\lambda_1,\lambda_2,\lambda_3,\lambda_4)\in \C^4\right\rbrace$.

We let $r_1=-\epsilon_1-\epsilon_2$, $r_2=\epsilon_1-\epsilon_2$, $r_3=\epsilon_2-\epsilon_3$ and $r_4=\epsilon_3-\epsilon_4$ be the positive simple roots of $\Phi$ and $\pi=\{r_1,r_2,r_3,r_4\}$. They correspond to the following Dynkin diagram

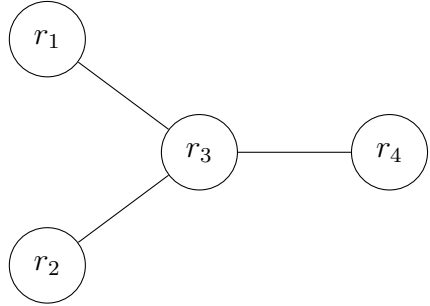
\begin{center}
\begin{tikzpicture}
[place/.style={circle,draw=black,
inner sep=1pt,minimum size=10mm}]
\node (1) at (0,3/2)[place]{$r_1$};
\node (2) at (0,-3/2)[place]{$r_2$};
\node (3) at (2,0)[place]{$r_3$};
\node (4) at (4.5,0)[place]{$r_4$};
\draw (1) to (3);
\draw (2) to (3);
\draw (3) to (4);
\end{tikzpicture}
\end{center}

The $24$ roots can then be expressed in terms of $r_1,r_2,r_3$ and $r_4$. We list here the twelve positive roots : 
$r_1$, $r_2$, $r_3$, $r_4$, $r_1+r_3$, $r_2+r_3$, $r_3+r_4$, $r_1+r_2+r_3$, $r_1+r_3+r_4$, $r_2+r_3+r_4$, $r_1+r_2+r_3+r_4$, $r_1+r_2+2r_3+r_4$.

We have then fixed the constant structures $\{N_{r,s}, r,s\in \Phi, r+s\in \Phi\}$, which are defined by the relation $[e_r,e_s]=N_{r,s}e_{r+s}$. They have the following values :\\ $N_{r_1,r_3}=N_{r_2,r_3}=N_{r_3,r_4}=N_{r_1+r_3,r_4}=N_{r_1+r_3,r_2}=N_{r_2+r_3,r_4}=N_{r_1,r_3+r_4}=N_{r_2+r_3,r_1}=N_{r_2,r_3+r_4} =N_{r_2+r_3+r_4,r_1}=N_{r_1+r_3+r_4,r_2}=N_{r_1+r_2+r_3,r_4}=N_{r_1+r_2+r_3+r_4,r_3}=1$. 

We get the remaining ones using the facts that $N_{s,r}=-N_{r,s}$ and $N_{-r,-s}=-N_{r,s}$.

We write $\mathfrak{L}_{\Z}=\op{Vect}_{\Z}(h_r,r\in \pi, e_s,s\in \Phi\}$. We then call $(h_r,r\in \pi,e_s,s\in \Phi)$ the Chevalley basis of $\mathfrak{L}$.

For an element $x\in \mathfrak{L}$, we write $\op{ad}(x): \mathfrak{L} \rightarrow \mathfrak{L}, \op{ad(x)}(y)=[x,y]=xy-yx$. For each $r\in \Phi$, since $\op{ad}(e_r)$ is a nilpotent endomorphism, we can define an element $x_r(t)$ of $\op{End}_{\Q[[t]]}(\mathfrak{L}_\Z\otimes_\Z \Q[[t]])$ by $x_r(t)=\op{exp}(t \op{ad}(e_r))$. We check on the Chevalley basis that $x_r(t)$ is actually an element of $GL_{\Z[t]}(\mathfrak{L}_\Z\otimes _\Z[t])$.

We write $\mathfrak{L}_{\F_q}=\F_q \otimes \mathfrak{L}_\Z$. For an element $u\in \F_q$, let $\theta_u$ be the morphism from $\Z[t]$ to $\F_q$ by $\theta_u(k)=\overline{k}$ and $\theta_u(t)=u$. We can extend $\theta_u$ to a morphism $\tilde{\theta_u}$ from $GL_{\Z[t]}(\mathfrak{L}_\Z\otimes _\Z[t])$ to $\op{Aut}(\mathfrak{L}_{\F_q})$. We then define for $r\in \Phi$ and $u\in \F_q$ the elements $x_r(u)=\tilde{\theta_u}(x_r(t))\in \op{Aut}(\mathfrak{L}_{\F_q})$.

The adjoint Chevalley group $\mathfrak{L}(\F_q)$ is then defined to be $<x_r(u),r\in \Phi, u\in \F_q>$.
\end{Def}

We now determine an isomorphism between the Chevalley group $\mathfrak{L}(\F_q)$ and $P\Omega_8^+(q)$, it can be found in \cite{CAR}. We have the following classical lemma

\begin{lemme}\label{Fallet}
Let $\mathcal{A}$ be a matrix algebra over $\F_q$ with Lie bracket defined the usual way. We then have that for any nilpotent matrix $y$ and any matrix $x$,
$$\op{exp}(\op{ad}(y))(x)=\op{exp}(y)x\op{exp}(y)^{-1}.$$
\end{lemme}

\begin{prop}\label{morigenisation}
We have $P\Omega_8^+(q)\simeq \mathfrak{L}(\F_q)$.
\end{prop}

\begin{proof}
The group $\Omega_8^+(q)$ is generated by long root elements which are elements of the form $I_8+ue_r$ with $r\in \Phi$ and $u\in \F_q$ (See for example \cite{W} 3.7.3). Since $e_r^2=0$, we have $\op{exp}(ue_r)=I_8+u e_r$. We define the morphism $\Psi$ from $\Omega_8^+(q)$ to $\mathfrak{L}(\F_q)$ on those generators by \\
$\Psi(I_8+ue_r)=(x\mapsto (I_8+ue_r)x(I_8+u e_r)^{-1}=\op{exp}(u e_r)x\op{exp}(u e_r)^{-1}=\op{exp}(\op{ad}(u e_r))(x)$, where the last equality follows from Lemma \ref{Fallet}.

We then have $\Psi(I_8+ue_r)=x_r(u)$, therefore the morphism is surjective. Let $y\in \op{ker}(\Psi)$, we have $yxy^{-1}=x$ for all $x\in \mathfrak{L}_{\F_q}$, therefore $y\in \F_q^{\star}I_8$, where $\F_q^\star$ is the group of invertible elements of $\F_q$. It follows that $\op{ker}(\Psi)=Z(\Omega_8^+(q))$ and we get the desired isomorphism.
\end{proof}

We now recall Proposition 12.2.3. of \cite{CAR}.

\begin{prop}
Suppose all the roots of $\mathfrak{L}$ have the same length and let $r\mapsto \tau(r)$ be a map of $\Phi$ into itself arising from a symmetry of the Dynkin diagram of $\mathfrak{L}$. Then there exists numbers $\gamma_r=\pm 1$ such that the map $x_r(u)\mapsto x_\tau(r)(\gamma_ru)$ can be extended to an automorphism of $\mathfrak{L}(K)$. The $\gamma_r$ can be chosen, therefore that $\gamma_r=1$ for all $r\in \pi \cup -\pi$. For $r=r_1+r_2$, we have $\gamma_r=\frac{\gamma_{r_1}\gamma_{r_2}N_{\tau(r_1),\tau(r_2)}}{N_{r_1,r_2}}$.
\end{prop}

We now apply this proposition with $\tau$ the following triality automorphism of the Dynkin diagram of type $D_4$ 

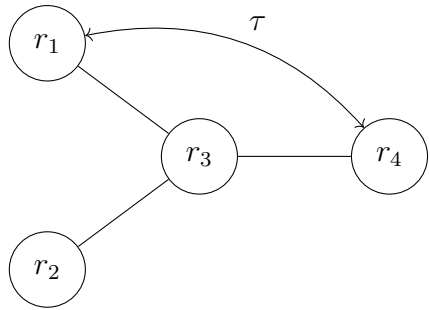
\begin{center}
\begin{tikzpicture}
[place/.style={circle,draw=black,
inner sep=1pt,minimum size=10mm}]
\node (1) at (0,3/2)[place]{$r_1$};
\node (2) at (0,-3/2)[place]{$r_2$};
\node (3) at (2,0)[place]{$r_3$};
\node (4) at (4.5,0)[place]{$r_4$};
\draw (1) to (3);
\draw (2) to (3);
\draw (3) to (4);
\draw[<->] (1) to [bend left=30]node[auto]{$\tau$} (4);
\end{tikzpicture}
\end{center}

We choose as in the proposition $\gamma_r=1$ for all $r\in \pi\cup -\pi$. We can then compute the remaining $\gamma_r$ using the induction relation in the proposition. We get

$-\gamma_{r_1+r_3}=\gamma_{r_2+r_3}=-\gamma_{r_3+r_4}=1$,\\
$\gamma_{r_1+r_2+r_3}=\gamma_{r_2+r_3+r_4}=\gamma_{r_1+r_3+r_4}=\gamma_{r_1+r_2+r_3+r_4}=\gamma_{r_1+r_2+2r_3+r_4}=1$. For the negative roots, we use the fact that $\gamma_{-r}=\gamma_r$ for all $r\in \Phi$.

We write in the following $T_r(u)=I_8+ue_r$. By Proposition \ref{morigenisation}, we get an isomorphism of $P\Omega_8^+(q)$ defined on the generators $\overline{T_r(u)}$ by $\tau(\overline{T_{\tau(r)}(\gamma_r u)})$.

\begin{prop}\label{Spin}
Assume $\F_q=\F_p(\sqrt{\alpha})=\F_p(\alpha)=\F_p(\alpha+\alpha^{-1})$. Let $G=\rho_{8_r}\times \rho_{8_{rr}}(\mathcal{A}_{H_4})$. We have that $\tau \circ \overline{\rho_{8_{rr}|\mathcal{A}_{H_4}}}\simeq \overline {\rho_{8_{r}|\mathcal{A}_{H_4}}}$.

If $\F_q=\F_p(\alpha)=\F_p(\alpha+\alpha^{-1})$, then $G\simeq Spin_8^+(q)$.

If $\F_q=\F_p(\alpha)\neq \F_p(\alpha+\alpha^{-1})$, then $G\simeq Spin_8^+(q^{\frac{1}{2}})$.
\end{prop}

\begin{proof}
Let us first show the first statement. By \cite{MR}, we have that $\mathcal{A}_{H_4}$ is generated by $u_0=T_{s_1}T_{s_2}^{-1}$, $u_1=T_{s_2}T_{s_1}T_{s_2}^{-2}$, $u_2=T_{s_2}^2T_{s_1}T_{s_2}^{-3}$, $u_3=T_{s_2}^3T_{s_1}T_{s_2}^{-4}$, $u_4=T_{s_3}T_{s_1}^{-1}$ and $u_5=T_{s_4}T_{s_1}^{-1}$. We want to consider the image of $\overline{\rho_8(u_i)}$ under the triality isomorphism for all $i\in [\![0,5]\!]$. We first have to write those elements as products of the generators $\overline{T_r(u)}$. We order lexicographically the vertices  $\{x_i,i\in [\![1,8]\!]\}$ and $\{y_i,i\in [\![1,8]\!]\}$ of the graphs of $8_r$ and $8_{rr}$  respectively, i.e. $I(x_1)=I(y_1)=\{s_1\}$, $I(x_2)=I(y_2)=\{s_2\}$, $I(x_3)=I(y_3)=\{s_1,s_3\}$, $I(x_4)=I(y_4)=\{s_2,s_3\}$, $I(x_5)=I(y_5)=\{s_1,s_4\}$, $I(x_6)=I(y_6)=\{s_2,s_4\}$, $I(x_7)=I(y_7)=\{s_1,s_3,s_4\}$ and $I(x_8)=I(y_8)=\{s_2,s_3,s_4\}$. 

By Theorem \ref{bilinwgraphs}, we have that for all $i\in [\![0,5]\!]$, $P\rho_{8_r}(u_i)P^{-1}=\!^t\rho_{8_r}(u_i)^{-1}$ and\\
 $P\rho_{8_{rr}}(u_i)P^{-1}=\!^t \rho_{8rr}(u_i)^{-1}$, where $P=\begin{pmatrix} 0 & 0 & 0 &0 & 0 & 0 & 0 & 1\\
0 & 0 & 0 & 0 & 0 & 0 & -1 & 0\\
0 & 0 & 0 & 0 &0 & 1 & 0& 0\\
0 & 0 & 0 & 0 & -1 & 0 & 0 &0\\
0 & 0 & 0 &-1 & 0 & 0 & 0 & 0\\
0 & 0 & 1 & 0 & 0 & 0 & 0 & 0\\
0 & -1 & 0 & 0 & 0 & 0 & 0 & 0\\
1 & 0 & 0 & 0 &0 &0 &0 & 0\end{pmatrix}$.

In order to write those elements in terms of the generators, we first need to work in the right basis. Let $M=\begin{pmatrix} 0& 0& 0& -1& 0 & 0& 0& 0\\
 0 & 0 & 1 & 0 & 0 & 0 &  0 & 0\\
  0 & -1 & 0 & 0 & 0 & 0 & 0 & 0\\
   1 & 0 & 0 & 0 & 0 & 0 & 0 & 0\\
   0 & 0 & 0 & 0 & 1 & 0 & 0 & 0\\
   0 & 0& 0& 0& 0& 1& 0& 0\\
    0& 0& 0 & 0& 0& 0& 1&  0\\
    0 &0 &0& 0& 0 & 0& 0& 1\end{pmatrix}$, we then have 
    $$A( M\rho_{8_r}(u_i)M^{-1} )A^{-1}=^t\!(M\rho_{8_r}(u_i)M^{-1})^{-1}$$ and $$A( M\rho_{8_{rr}}(u_i)M^{-1} )A^{-1}=^t\!(M\rho_{8_r}(u_i)M^{-1})^{-1}$$, where $A=\begin{pmatrix} 0& I_4\\
    I_4 & 0\end{pmatrix}$ as in Definition \ref{deftriality}.
    
    Let now $U_i=M\rho_{8_{rr}}(u_i)M^{-1}$  and $\tilde{U}_i=M\rho_{8_{r}}(u_i)M^{-1}$ for all $i\in [\![0,5]\!]$. We check that $\tau(\overline{U_i})=\overline{\tilde{U}_i}$ for all $i\in [\![0,5]\!]$ by explicit computations given in subsection \ref{computationstriality} of the Appendix.

This proves that $\tau\circ \overline{\rho_{8_r}|\mathcal{A}_{H_4}}\simeq\overline{\rho_{8_{rr}}|\mathcal{A}_{H_4}}$.

\smallskip

It follows that if $\F_p(\alpha)=\F_p(\alpha+\alpha^{-1})$ then $\rho_{8_{rr}}(\mathcal{A}_{H_4})\simeq \Omega_8^+(q)$ and if $\F_p(\alpha)\neq \F_p(\alpha+\alpha^{-1})$ then $\rho_{8_{rr}}(\mathcal{A}_{H_4})\simeq \Omega_8^+(q^{\frac{1}{2}})$ since it is generated by long root elements and irreducible.

\bigskip

Assume now that $\F_q=\F_p(\alpha)=\F_p(\alpha+\alpha^{-1})$.
We now want to use Goursat's Lemma (see Lemma \ref{Goursat}). We consider $G$ as a subgroup of $\Omega_8^+(q)\times \Omega_8^+(q)$. We write $\pi_1$ the projection upon the first factor and $\pi_2$ the projection upon the second factor. We write $K_1=\pi_1(G)\simeq\Omega_8^+(q)$, $K_2=\pi_2(G)\simeq \Omega_8^+(q)$, $K^1=\op{ker}(\pi_2)$ and $K^2=\op{ker}(\pi_1)$. By Goursat's Lemma, there exists an isomorphism $\varphi$ such that $G=\{(x,y)\in K_1\times K_2, \varphi(xK^1)=\varphi(yK^2)\}$. 

Let $x\in \op{ker}(\pi_1)$, we have $\pi_1(x)=I_8$. We know that projectively, we have $\tau(\overline{\pi_1(x)})=\overline{\pi_2(x)}$, therefore $\overline{\pi_2(x)}=\overline{I_8}$ and $\pi_2(x)\in\{\pm I_8\}$. 

This proves that $K^2\leq \{I_8\}\times \{\pm I_8\}$ and $K^1\leq \{\pm I_8\} \times \{I_8\}$.

We have $\rho_{8_r}((T_{s_1}T_{s_3}T_{s_2}T_{s_4})^{15}(T_{s_1}T_{s_3}T_{s_2})^{-20})=I_8=-\rho_{8_{rr}}((T_{s_1}T_{s_3}T_{s_2}T_{s_4})^{15}(T_{s_1}T_{s_3}T_{s_2})^{-20})$. It follows that $I_8\times (-I_8) \in G$, therefore $K^2=I_8\times (-I_8)$ and $K_2/K^2\simeq P\Omega_8^+(q)$. This implies that $K_1/K^1\simeq P\Omega_8^+(q)$ and $K^1=(-I_8)\times I_8$.

For every $x\in K_1$, we have exactly two elements of $K_2$ such that $\varphi(\overline{x})=\overline{y}$. This implies that $\vert G\vert =2\vert \Omega_8^+(q)\vert$. We also know that for all $x\in \Omega_8^{+}(q)$, there exists $y_x\in \Omega_8^+(q)$ such that $\tau(\overline{x})=\overline{y_x}$ and $(x,y_x)\in G$. We know there exists $h\in \mathcal{A}_{H_4}$ such that 
$$x=\rho_{8_r}(h)\rho_{8_r}((T_{s_1}T_{s_3}T_{s_2}T_{s_4})^{15}(T_{s_1}T_{s_3}T_{s_2})^{-20})$$
, therefore 
$$(x,\rho_{8_{rr}}(h)\rho_{8_{rr}}((T_{s_1}T_{s_3}T_{s_2}T_{s_4})^{15}(T_{s_1}T_{s_3}T_{s_2})^{-20}))=(x,-y_x)\in G.$$
 This proves that $G$ contains $\{(x,y)\in K_1\times K_2, \tau(\overline{x})=\overline{y}\}$. It follows using the cardinality of those two sets that $G=\{(x,y)\in K_1\times K_2,\tau(\overline{x})=\overline{y}\}$. 

We also have that the center $Z(G)$ of $G$ is equal to $\{\pm I_8\}\times \{\pm I_8\}$ since this group is included in $G$ and $Z(\Omega_8^+(q))=\{\pm I_8\}$. It follows that $G/Z(G)=\{(\overline{x},\tau(\overline{x})), x\in \Omega_8^+(q)\}\simeq P\Omega_8^+(q)$.

We also have that $G$ has two normal subgroups $Z_1$ and $Z_2$ of order $2$, $Z_1=\{\pm I_8\}\times \{I_8\}$ and $Z_2=\{ I_8\} \times \{\pm I_8\}$ such that $G/Z_1\simeq G/Z_2 \simeq \Omega_8^+(q)$. This proves that $G$ is a double cover of $\Omega_8^+(q)$. 

\smallskip

We recall here the definitions of a central extension and of a universal central extension from \cite{Aschbabook}. A central extension of a group $\Gamma$ is a pair $(H,\pi)$, where $H$ is a group and $\Pi:H\rightarrow \Gamma$ is a surjective homomorphism with $\op{ker}(\pi)\leq Z(H)$. A morphism $\alpha:(\Gamma_1,\pi_1)\rightarrow (\Gamma_2,\pi_2)$ of central extensions of $\Gamma$ is a group homomorphism $\alpha: \Gamma_1\rightarrow \Gamma_2$ with $\pi_1=\pi_2\alpha$. A central extension $(\tilde{\Gamma},\pi)$ of $\Gamma$ is universal if for each central extension $(H,\sigma)$ of $G$ there exists a unique morphism $\alpha:(\tilde{\Gamma},\pi)\rightarrow (H,\sigma)$ of central extensions. By \cite{Aschbabook} (33.1), there exists at most one universal central extension of a group $\Gamma$. By \cite{Aschbabook} (33.4), any perfect group posseses a universal central extension. The universal central extension $\tilde{\Gamma}$ is then called the universal covering group of $\Gamma$ and $\op{ker}(\pi)$ is the Schur multiplier of $\Gamma$.

\smallskip

We now show that $G$ is the universal covering group of $P\Omega_8^+(q)$. We have by \cite{GorLySol} (Theorem 6.1.4 and Table 6.1.2) that the Schur multiplier of $P\Omega_8^+(q)$ is $(\Z/2\Z)^2$. By Theorem 1.10.7 of \cite{GorLySol}, we have that the universal cover of $P\Omega_8^+(q)$ is $Spin_8^+(q)$.

By \cite{MR}, $\mathcal{A}_{H_4}$ is perfect, it follows that $G$ is perfect. We have shown above that $G/Z(G)\simeq P\Omega_8^+(q)$, this proves that $G$ is a perfect central extension of $G$. We have $Z(G)\simeq (\Z/2\Z)^2$ therefore $\vert G\vert =\vert Spin_8^+(q)\vert$. By \cite{Aschbabook} (33.8), there exists an onto morphism from $Spin_8^+(q)$ to $G$. This proves they are isomorphic since their order is equal.

\medskip

All the arguments are identical for $\F_p(\alpha)\neq \F_p(\alpha+\alpha^{-1})$.
\end{proof}

\textbf{Remark} : Note that the restrictions to $A_{H_3}$ of $\rho_{8_r}$ and $\rho_{8_{rr}}$ are identical (this can be seen on the $W$-graphs by removing $4$ from the vertices and deleting the edges from $x$ to $y$ if $I(x)=I(y)$), this proves that the projective twisted diagonal $PSL_4(q)$ inside $P\Omega_8^+(q)$ is stabilized by the triality automorphism. 

We also have that $\tau$ does not extend to a morphism from $\Omega_8^+(q)$ into itself. Assume it does, there exists $\Psi$ from $K_1\simeq\Omega_8^+(q)$ to $K_2\simeq\Omega_8^+(q)$ such that $\tau(\overline{x})=\overline{\Psi(x)}$ for all $x\in K_1$. We then have that $\overline{\rho_{8_{rr}}(u_i)}=\tau(\overline{\rho_{8_r}(u_i)})=\overline{\Psi(\rho_{8_r}(u_i)}$ for all $i\in [\![1,6]\!]$, therefore there exists $\epsilon_i\in \{-1,1\}$ such that $\rho_{8_{rr}}(u_i)=\epsilon_i\Psi(\rho_{8_r}(u_i))$. By \cite{MR}, we have that $\mathcal{A}_{H_4}$ is perfect, therefore $\epsilon_i=1$ for all $i\in [\![1,6]\!]$ and for all $h\in \mathcal{A}_{H_4}$, we have $\Psi(\rho_{8_r}(h))=\rho_{8_{rr}}(h)$. This is absurd because \\
 $\rho_{8_r}((T_{s_1}T_{s_3}T_{s_2}T_{s_4})^{15}(T_{s_1}T_{s_3}T_{s_2})^{-20})=I_8=-\rho_{8_{rr}}((T_{s_1}T_{s_3}T_{s_2}T_{s_4})^{15}(T_{s_1}T_{s_3}T_{s_2})^{-20})$.
 
 \section{Type $H_4$, high-dimensional representations}\label{H4highdim}
 
 We now give two propositions, where we determine the image of $\mathcal{A}_{H_4}$ in the remaining irreducible representations. We could use the theorem by Guralnick and Saxl \cite{GS} for all of them but we use arguments not requiring the classification of finite simple groups to determine those images when possible. The first of the two following propositions does not use the classification and the second one does.
 
 \begin{prop}
 Assume $1\sim 2$.
 \begin{enumerate}
 \item If $\F_q=\F_p(\alpha)=\F_p(\alpha+\alpha^{-1})$, then $\rho_{9_s}(\mathcal{A}_{H_4})\simeq SL_9(q^2)$ and $\rho_{\overline{9_s}}(\mathcal{A}_{H_4})\simeq SL_9(q^2)$. We also have that $\Phi_{1,2}\circ\rho_{\overline{9_s}|\mathcal{A}_{H_4}}\simeq \rho_{9_s|\mathcal{A}_{H_4}}$.
 \item If $\F_q=\F_p(\alpha)\neq \F_p(\alpha+\alpha^{-1})$, then $\rho_{9_s}(\mathcal{A}_{H_4})\simeq SL_9(q^2)$ and $\rho_{\overline{9_s}}(\mathcal{A}_{H_4})\simeq SL_9(q^2)$. We also have that $\Phi_{1,2}\circ\rho_{\overline{9_s}|\mathcal{A}_{H_4}}\simeq \rho_{9_s'|\mathcal{A}_{H_4}}$.
 \end{enumerate}

 \smallskip

If $1\nsim 2$ then we have $\rho_{9_s}(\mathcal{A}_{H_4})\simeq SU_9(q^{\frac{1}{2}})$ and $\rho_{\overline{9_s}}(\mathcal{A}_{H_4})\simeq SU_9(q^{\frac{1}{2}})$.
\smallskip

If $\F_q=\F_p(\alpha)=\F_p(\alpha+\alpha^{-1})$ then $\rho_{10_r}(\mathcal{A}_{H_4})\simeq \Omega_{10}^+(q)$.

\smallskip

If $\F_q=\F_p(\alpha)\neq \F_p(\alpha+\alpha^{-1})$ then $\rho_{10_r}(\mathcal{A}_{H_4})\simeq \Omega_{10}^+(q^{\frac{1}{2}})$.

 \end{prop}
 
 \begin{proof}
Assume $1\sim 2$ and $\F_p(\alpha)=\F_p(\alpha+\alpha^{-1})$. By Table \ref{resH4H3} and Theorem \ref{resicos}, we have $\rho_{9_s}(\mathcal{A}_{H_3})\simeq SL_3(q^2)\times SL_5(q)$ in a natural representation. By Lemma \ref{normclosiscos}, we then have that $\rho_{9_s}(\mathcal{A}_{H_4})$ is an irreducible group generated by transvections. We also have that $\rho_{9_s}$ is not self-dual, therefore $\rho_{9_s}(\mathcal{A}_{H_4})$ is included in no symplectic group. We also have that the field generated by the traces of the elements of $\rho_{9_s}(\mathcal{A}_{H_4})$ contains $\F_{q^2}$. It follows by Theorem \ref{transvections} that $\rho_{9_s}(\mathcal{A}_{H_4})\simeq SL_9(q^2)$.

\smallskip

Assume $1\sim 2$ and $\F_p(\alpha)\neq \F_p(\alpha+\alpha^{-1})$. We then have by Lemma \ref{IsomorphismH3} that $\F_p(\alpha,\xi+\xi^{-1})=\F_p(\alpha)$. We have by Theorem \ref{resicos} that $\rho_{9_s}(\mathcal{A}_{H_3})\simeq SL_3(q)\times SU_5(q^{\frac{1}{2}})$. It follows that $\rho_{9_s}(\mathcal{A}_{H_4})$ is neither unitary nor symplectic and is an irreducible group generated by transvections. We also have that the field generated by the traces of its elements contains $\F_q$. We also have by Proposition \ref{color} that it is conjugate to a subgroup of $SL_9(q)$. It follows by Theorem \ref{transvections} that $\rho_{9_s}(\mathcal{A}_{H_4})\simeq SL_9(q)$. 

\smallskip

If $1\sim 2$, we have that $\Phi_{1,2}(\alpha)\in \{\alpha,\alpha^{-1}\}$. It follows by Proposition \ref{Fieldfactorization} that $\Phi_{1,2}\circ \rho_{\overline{9_s}|\mathcal{A}_{H_4}}$ is isomorphic to the restriction of an irreducible representation of $\mathcal{H}_{H_4,\alpha}$ to $\mathcal{A}_{H_4}$. The factorization result then follows by Table \ref{resH4H3} and Proposition \ref{resrepreflicos}.

\medskip

Assume now $1\nsim 2$. We have by Theorem \ref{resicos} that $\rho_{9_s}(\mathcal{A}_{H_4})\simeq SU_3(q^{\frac{1}{2}}) \times SU_5(q^{\frac{1}{2}})$ in a natural representation. Let $\epsilon$ be the unique automorphism of order $2$ of $\F_q$. We have $\epsilon(\alpha)=\alpha^{-1}$. It follows by Proposition \ref{Fieldfactorization} that $\epsilon\circ \rho_{9_s'|\mathcal{A}_{H_4}}$ is isomorphic to the restriction of an irreducible representation of $\mathcal{H}_{H_4,\alpha}$ to $\mathcal{A}_{H_4}$. By Table \ref{resH4H3}, this implies that $\epsilon \circ\rho_{9_s'|\mathcal{A}_{H_4}}\simeq \rho_{9_s|\mathcal{A}_{H_4}}$. It follows by Lemma \ref{Harinordoquy} that up to conjugation in $GL_9(q)$, we have $\rho_{9_s}(\mathcal{A}_{H_4})\leq SU_9(q^{\frac{1}{2}})$. We have, by Lemma \ref{field}, that the field generated by the traces of the elements of $\rho_{9_s}(\mathcal{A}_{H_4})$ contains $\F_q$. It follows by Theorem \ref{transvections} that $\rho_{9_s}(\mathcal{A}_{H_4})\simeq SU_9(q^{\frac{1}{2}})$. The same arguments show that $\rho_{\overline{9_s}}(\mathcal{A}_{H_4})\simeq SU_9(q^{\frac{1}{2}})$

\bigskip

Consider now the $10$-dimensional representation $\rho_{10_r}$. Let $G=\rho_{10_r}(\mathcal{A}_{H_4})$. By Proposition \ref{color}, we can assume that it is defined over $\F_q$.

\smallskip

 Assume first that $\F_q=\F_p(\alpha)=\F_p(\alpha+\alpha^{-1})$. By Table \ref{resH4H3} and Theorem \ref{resicos}, we have that $\rho_{10_r}(\mathcal{A}_{H_3})\simeq SL_5(q)$ in a twisted diagonal representation. We have by Proposition \ref{bilinwgraphs} that up to conjugation in $GL_{10}(q)$, we have $G\leq \Omega_{10}^+(q)$. It follows by Lemma \ref{normclosiscos} that $G$ is an irreducible subgroup of $\Omega_{10}^+(q)$ generated by long root elements. It follows by Theorem \ref{theoKantor} that $G$ is isomorphic to one of the following groups
\begin{enumerate}
\item $\Omega_{10}^+(q)$
\item $\Omega_{10}^-(q^{\frac{1}{2}})$
\end{enumerate}

If $\rho_{10_r}(\mathcal{A}_{H_4})$ was conjugate to $\Omega_{10}^-(q^\frac{1}{2})$, then we would have that $\vert SL_5(q)\vert=q^{10}(q^2-1)(q^3-1)(q^4-1)(q^5-1)$ divides $\vert \Omega_{10}^{-}(q^\frac{1}{2})\vert =q^{10}(q^{\frac{5}{2}}+1)(q^4-1)(q^3-1)(q^2-1)(q-1)$. This would imply that $q^5-1$ divides $(q^{\frac{5}{2}}+1)(q-1)=q^{\frac{7}{2}}-q^{\frac{5}{2}}+q-1<q^{\frac{7}{2}}-1=q^5-1$, which is absurd. This proves that $\rho_{10_r}(\mathcal{A}_{H_4})=\Omega_{10}^+(q)$. 

\smallskip

Assume now that $\F_q=\F_p(\alpha)\neq \F_p(\alpha+\alpha^{-1})$. Let $\epsilon$ be the unique automorphism of order $2$ of $\F_q$, we have $\epsilon(\alpha)=\alpha^{-1}$. It follows by Proposition \ref{Fieldfactorization} that $\epsilon \circ\rho_{10_r}\simeq \rho_{10_r}$. We then have by Proposition $4.1$ of \cite{BMM} that up to conjugation in $GL_{10}(q)$, we have $G\leq \Omega_{10}^+(q^{\frac{1}{2}})$. We have $\rho_{10_r}(\mathcal{A}_{H_3})\simeq SU_5(q^{\frac{1}{2}})$ in a twisted diagonal representation. By the same arguments as above, we have $G\simeq \Omega_{10}^+(q^{\frac{1}{2}})$ or $G\simeq \Omega_{10}^-(q^{\frac{1}{4}})$.

 Assume by contradiction that $G\simeq \Omega_{10}^-(q^{\frac{1}{4}})$. We then have that $\vert SU_5(q^{\frac{1}{2}})\vert =q^5(q-1)(q^{\frac{3}{2}}+1)(q^2-1)(q^{\frac{5}{2}}+1)$ divides $\vert \Omega_{10}^{-}(q^{\frac{1}{4}})\vert =q^5(q^{\frac{5}{4}}+1)(q^2-1)(q^{\frac{3}{2}}-1)(q-1)(q^{\frac{1}{2}}-1)$. This implies that $q^{\frac{5}{2}}+1$ divides $(q^{\frac{5}{4}}+1)(q^{\frac{3}{2}}-1)(q^{\frac{1}{2}}-1)$. It follows that $q^{\frac{5}{2}}+1$ divides $(q^{\frac{5}{2}}+1)(q^{\frac{3}{2}}-1)(q^{\frac{1}{2}}-1)-(q^{\frac{5}{4}}+1)(q^{\frac{5}{4}}-1)(q^{\frac{3}{2}}-1)(q^{\frac{1}{2}}-1)=2(q^{\frac{3}{2}}-1)(q^{\frac{1}{2}}-1)=2q^2-2q^{\frac{3}{2}}-2q^{\frac{1}{2}}+2<q^{\frac{5}{2}}+1$ which is absurd. It follows that $G\simeq \Omega_{10}^+(q^{\frac{1}{2}})$.
\end{proof}
 
 The remaining cases will be proved as in \cite{BMM} and types $B$ and $D$ using Theorem \ref{CGFS}. We first give a Lemma which will help us prove that the groups considered are tensor-indecomposable.
 
 \begin{lemme}\label{exceptiontens2}
 Let $\ell=p^k$ with $p\notin \{2,3,5\}$ a prime and $k\geq 1$. If $r\geq 3$ then there exists no non-trivial morphism from $SL_r(\ell')$ to $SL_2(\ell^4)$ if $\ell'\in \{\ell,\ell^2,\ell^4\}$ or $SU_r(\ell')$ to $SL_2(\ell^4)$ if $\ell'\in \{\ell,\ell^2\}$.
\end{lemme}

\begin{proof}
 Let $\theta$ be a morphism from $SL_3(\ell)$ to $SL_2(\ell^4)$. We have that $\op{ker}(\theta)$ is a normal subgroup of $SL_3(\ell)$.

Assume this kernel is different from $SL_3(\ell)$. The image of $\theta$ is then non-abelian. We get an isomorphism from $PSL_3(\ell)$ to a subgroup of $SL_2(\ell^4)$. We have $\vert PSL_3(\ell)\vert =\frac{1}{\op{Gcd}(3,\ell-1)}\ell^3(\ell^2-1)(\ell^3-1)$ and $\vert SL_2(\ell^4)\vert =\ell^4(\ell^8-1)$. It follows that $\ell^3-1$ divides $3(\ell^8-1)$. This implies that $\ell^3-1$ divides $3(\ell^8-1)-3\ell^2(\ell^6-1)=3(\ell^2-1)$. We have $\ell>3$ since $p\notin \{2,3,5\}$, therefore $3\ell^2-3<\ell^3-1$ which is a contradiction. This proves that there is no non-trivial morphsim from $SL_3(\ell)$ to $SL_2(\ell^4)$. This implies in the same way that there exists no non-trivial morphism from $SL_r(\ell')$ to $SL_2(\ell^4)$ if $r\geq 3$ and $\ell'\in \{\ell,\ell^2,\ell^4\}$. 

In the same way, if there was a non trivial morphism from $SU_r(\ell')$ to $SL_2(\ell^4)$ for $r\geq 3$ and $q'\in \{\ell,\ell^2\}$ then we would have that $\ell^3+1$ divides $3(\ell^2-1)$ which is also absurd.
\end{proof}

 \begin{prop}
We write $\epsilon$ the unique automorphism of order $2$ of $\F_q$ when it exists. We then have the following results.
\begin{enumerate}
\item If $\F_q=\F_p(\sqrt{\alpha})=\F_p(\alpha)=\F_p(\alpha+\alpha^{-1})$, then $\rho_{16_r}(\mathcal{A}_{H_4})=SL_{16}(q)$ and $\rho_{16_{rr}}(\mathcal{A}_{H_4})=SL_{16}(q)$.
\item If $\F_{q^2}=\F_p(\sqrt{\alpha})\neq \F_p(\alpha)=\F_p(\alpha+\alpha^{-1})$, then $\rho_{16_r}(\mathcal{A}_{H_4})\simeq SL_{16}(q^2)$ and $\rho_{16_{rr}|\mathcal{A}_{H_4}}\simeq \varphi\circ \rho_{16_r|\mathcal{A}_{H_4}}$, where $\varphi$ is the unique automorphism of order $2$ of $\F_{q^2}$.
\item If $\F_q=\F_p(\alpha)=\F_p(\sqrt{\alpha})\neq \F_p(\alpha+\alpha^{-1})$ and $\epsilon(\sqrt{\alpha})=\sqrt{\alpha}^{-1}$, then $\rho_{16_r}(\mathcal{A}_{H_4})\simeq SU_{16}(q^{\frac{1}{2}})$ and $\rho_{16_{rr}}(\mathcal{A}_{H_4})\simeq SU_{16}(q^{\frac{1}{2}})$.
\item If $\F_q=\F_p(\alpha)=\F_p(\sqrt{\alpha})\neq \F_p(\alpha+\alpha^{-1})$ and $\epsilon(\sqrt{\alpha})=-\sqrt{\alpha}^{-1}$, then $\rho_{16_r}(\mathcal{A}_{H_4})\simeq SL_{16}(q)$ and $\rho_{16_{rr}'|\mathcal{A}_{H_4}}\simeq \epsilon \circ \rho_{16_r|\mathcal{A}_{H_4}}$.
\end{enumerate}

\medskip
 
Assume $1\sim 2$. We then have the following results
\begin{enumerate}
\item We have $\Phi_{1,2}\circ \rho_{16_t|\mathcal{A}_{H_4}}\simeq \rho_{\overline{16_t}|\mathcal{A}_{H_4}}$, $\Phi_{1,2}\circ \rho_{24_s|\mathcal{A}_{H_4}}\simeq \rho_{\overline{24_s}|\mathcal{A}_{H_4}}$, $\Phi_{1,2}\circ \rho_{24_t|\mathcal{A}_{H_4}}\simeq \rho_{\overline{24_t}|\mathcal{A}_{H_4}}$ and $\Phi_{1,2}\circ \rho_{30_s|\mathcal{A}_{H_4}}\simeq \rho_{\overline{30_s}|\mathcal{A}_{H_4}}$.
\item If $\F_q=\F_p(\alpha)=\F_p(\alpha+\alpha^{-1})$, then $\rho_{16_t}(\mathcal{A}_{H_4})\simeq \Omega_{16}^+(q)$, $\rho_{24_s}(\mathcal{A}_{H_4})\simeq \Omega_{24}^+(q)$, $\rho_{24_t}(\mathcal{A}_{H_4})\simeq \Omega_{24}^+(q)$ and $\rho_{30_s}(\mathcal{A}_{H_4})\simeq \Omega_{30}^+(q)$.
\item If $\F_q=\F_p(\alpha)\neq \F_p(\alpha+\alpha^{-1})$, then $\rho_{16_t}(\mathcal{A}_{H_4})\simeq \Omega_{16}^+(q^{\frac{1}{2}})$, $\rho_{24_s}(\mathcal{A}_{H_4})\simeq \Omega_{24}^+(q^{\frac{1}{2}})$, $\rho_{24_t}(\mathcal{A}_{H_4})\simeq \Omega_{24}^+(q^{\frac{1}{2}})$ and $\rho_{30_s}(\mathcal{A}_{H_4})\simeq \Omega_{30}^+(q^{\frac{1}{2}})$.
\end{enumerate}

\smallskip

Assume $1\nsim 2$. We then have  $\rho_{16_t}(\mathcal{A}_{H_4})\simeq \Omega_{16}^+(q^{\frac{1}{2}})$, $\rho_{24_s}(\mathcal{A}_{H_4})\simeq \Omega_{24}^+(q^{\frac{1}{2}})$, $\rho_{24_t}(\mathcal{A}_{H_4})\simeq \Omega_{24}^+(q^{\frac{1}{2}})$ and $\rho_{30_s}(\mathcal{A}_{H_4})\simeq \Omega_{30}^+(q^{\frac{1}{2}})$.

\medskip

If $\F_q=\F_p(\alpha)=\F_p(\alpha+\alpha^{-1})$, then $\rho_{18_r}(\mathcal{A}_{H_4})\simeq \Omega_{18}^+(q)$, $\rho_{25_r}(\mathcal{A}_{H_4})\simeq SL_{25}(q)$, $\rho_{36_{rr}}(\mathcal{A}_{H_4})\simeq SL_{36}(q)$ and $\rho_{40_r}(\mathcal{A}_{H_4})\simeq \Omega_{40}^+(q)$.

\smallskip

If $\F_q=\F_p(\alpha)\neq \F_p(\alpha+\alpha^{-1})$, then $\rho_{18_r}(\mathcal{A}_{H_4})\simeq \Omega_{18}^+(q^{\frac{1}{2}})$, $\rho_{25_r}(\mathcal{A}_{H_4})\simeq SU_{25}(q^{\frac{1}{2}})$, $\rho_{36_{rr}}(\mathcal{A}_{H_4})\simeq SU_{36}(q^{\frac{1}{2}})$ and $\rho_{40_r}(\mathcal{A}_{H_4})\simeq \Omega_{40}^+(q^{\frac{1}{2}})$.
 \end{prop}

\begin{proof}
We first check that in all cases, the assumptions of Theorem \ref{CGFS} are verified. We have as in the proof of primitivity in \cite{BMM} and in type $B$  that all the groups considered are primitive finite irreducible subgroup of $GL(V)$ by considering Table \ref{resH4H3}, because they contain either a natural $SL_2(q')$ or a twisted diagonal $SL_3(q')$ or $SU_3(q')$ for $q'\in \{q^{\frac{1}{2}},q,q^2\}$. Using again Table \ref{resH4H3} and Lemmas \ref{tens1}, \ref{tens2} and \ref{exceptiontens2}, all the groups considered are tensor-indecomposable. We also have $v_G(V)\leq 2\leq \max(2,\frac{\sqrt{d}}{2})$. As in page $13$ of \cite{BMM}, the restriction also shows us that we are not in case $2$ of Theorem \ref{CGFS}. This shows that they are classical groups in a natural representation.

\bigskip

Consider first the $H_4$-graphs $16_r$ and $16_{rr}$. They are not $2$-colorable, therefore we cannot apply Proposition \ref{color}.

\smallskip

Assume $\F_q=\F_p(\sqrt{\alpha})=\F_p(\alpha)=\F_p(\alpha+\alpha^{-1})$. We then have by Theorem \ref{resicos} that $\rho_{16_r}(\mathcal{A}_{H_4})$ contains a natural $SL_5(q)$, therefore the field generated by the traces of its elements contains $\F_q$. It follows that it is a classical group in a natural representation over $\F_q$. The natural $SL_5(q)$ also shows that it is preserves no non-degenerate bilinear or hermitian form over $\F_q$. It follows that $\rho_{16_r}(\mathcal{A}_{H_4})=SL_{16}(q)$. The same arguments show that $\rho_{16_{rr}}(\mathcal{A}_{H_4})=SL_{16}(q)$.

\smallskip

Assume now $\F_{q^2}=\F_p(\sqrt{\alpha})\neq \F_p(\alpha)=\F_p(\alpha+\alpha^{-1})$. By Theorem \ref{resicos}, we have that $\rho_{16_r}(\mathcal{A}_{H_4})$ contains a natural $SU_4(q)$. It follows by Lemma \ref{field} that it is a classical group in a natural representation over $\F_{q^2}$.

 Let then $\varphi$ be the unique automorphism of order $2$ of $\F_{q^2}$. We have $\varphi(\sqrt{\alpha})=-\sqrt{\alpha}$ and $\varphi(\alpha)=\alpha$. It follows by Proposition \ref{Fieldfactorization} that $\varphi\circ\rho_{16_r|\mathcal{A}_{H_4}}$ is isomorphic to the restriction of an irreducible representation of $\mathcal{H}_{H_4,\alpha}$. By Table \ref{resH4H3} and Propositions \ref{resrepreflicos} and \ref{H3dim4}, we have that $\varphi\circ\rho_{16_r|\mathcal{A}_{H_4}}\simeq \rho_{16_{rr}|\mathcal{A}_{H_4}}$. This proves by Proposition \ref{resH4derivedsubgroup} that $\varphi\circ\rho_{16_r|\mathcal{A}_{H_4}}^\star\not\simeq \rho_{16_r|\mathcal{A}_{H_4}}^\star$. We also have that $\rho_{16_r}$ is not self-dual therefore $\rho_{16_r}(\mathcal{A}_{H_4})=SL_{16}(q^2)$. The result for $16_{rr}$ follows from the factorization.

\smallskip

Assume now $\F_q=\F_p(\sqrt{\alpha})=\F_p(\alpha)\neq \F_p(\alpha+\alpha^{-1})$ and $\epsilon(\sqrt{\alpha})=\sqrt{\alpha}^{-1}$. By Lemma \ref{IsomorphismH3}, we have $1\sim 2$, $\epsilon =\Phi_{1,2}$ 	and $\Phi_{1,2}(\alpha)=\alpha^{-1}$.  It follows by Proposition \ref{Fieldfactorization} that $\Phi_{1,2} \circ\rho_{16_r|\mathcal{A}_{H_4}}$ is isomorphic to the restriction of an irreducible representation of $\mathcal{H}_{H_4,\alpha}$. By Table \ref{resH4H3}, we have $\rho_{16_r|\mathcal{A}_{H_3}}\simeq \rho_{3_s'}\times \rho_{\overline{3_s'}}\times \rho_{4_r}\times \rho_{5_r}$. It follows by Propositions \ref{resrepreflicos}, \ref{H3dim4} and \ref{H3dim5} that $\Phi_{1,2}\circ \rho_{16_r|\mathcal{A}_{H_3}}\simeq \rho_{\overline{3_s}}\times \rho_{3_s}\times \rho_{4_r'}\times \rho_{5_r'}$. It follows by Table \ref{resH4H3} that $\Phi_{1,2}\circ \rho_{16_r|\mathcal{A}_{H_3}}\simeq \rho_{16_r'}$. By Lemma \ref{Harinordoquy}, we have that $\rho_{16_r}(\mathcal{A}_{H_4})\leq SU_{16}(q^{\frac{1}{2}})$, up to conjugation in $GL_{16}(q)$. It contains a natural $SU_5(q^{\frac{1}{2}})$, therefore it is a classical group in a natural representation over $\F_q$. It follows that $\rho_{16_r}(\mathcal{A}_{H_4})\simeq SU_{16}(q^{\frac{1}{2}})$. The same arguments show that $\rho_{16_{rr}}(\mathcal{A}_{H_4})\simeq SU_{16}(q^{\frac{1}{2}})$.

\smallskip

Assume $\F_q=\F_p(\sqrt{\alpha})=\F_p(\alpha)\neq \F_p(\alpha+\alpha^{-1})$ and $\epsilon(\sqrt{\alpha})=-\sqrt{\alpha}^{-1}$. We apply here the same reasoning as before and we get that $\Phi_{1,2}\circ \rho_{16_r|\mathcal{A}_{H_3}}\simeq \rho_{\overline{3_s}}\times \rho_{3_s}\times \rho_{4_r}\times \rho_{5_r'}$. It follows by Table \ref{resH4H3} that $\Phi_{1,2}\circ \rho_{16_r|\mathcal{A}_{H_4}}\simeq \rho_{16_{rr}'|\mathcal{A}_{H_4}}$. It follows that $\rho_{16_r}(\mathcal{A}_{H_4})$ preserves no non-degenerate bilinear or hermitian form over $\F_q$. It contains a natural $SU_5(q^{\frac{1}{2}})$ therefore it is a classical group over $\F_q$. It follows that $\rho_{16_r}(\mathcal{A}_{H_4})= SL_{16}(q)$. The result for $16_{rr}$ follows from the factorization.

\bigskip

All the remaining $H_4$-graphs are $2$-colorable therefore we can assume they are defined over $\F_p(\alpha,\xi+\xi^{-1})$. We know they are classical groups in a natural representation. We will show that they are the groups given in the proposition. We first consider the $H_4$-graphs which contain weights in $\F_p(\xi+\xi^{-1})$ i.e. $\widetilde{16_t}$, $\widetilde{\overline{24_s}}$, $\widetilde{\overline{24_t}}$ and $\widetilde{30_s}$. Note first that all of those representations are self-dual, and by Proposition \ref{bilinwgraphs} and the $H_4$-graphs given in the Appendix, we have that they preserve a non-degenerate symmetric bilinear form.

\medskip

Assume first $1\sim 2$ and $\F_q=\F_p(\alpha)=\F_p(\alpha+\alpha^{-1})$. We then have by Lemma \ref{IsomorphismH3} that $\F_{q^2}=\F_p(\alpha,\xi+\xi^{-1})=\F_p(\alpha+\alpha^{-1},\xi+\xi^{-1})$, $\Phi_{1,2}$ is the unique automorphism of order $2$ of $\F_{q^2}$ and $\Phi_{1,2}(\alpha)=\alpha$.

 Since $\rho_{16_t}$ is the only $16$-dimensional self-dual representation, we have by Proposition \ref{Fieldfactorization} that $\Phi_{1,2}\circ \rho_{16_t|\mathcal{A}_{H_4}}\simeq \rho_{16_t|\mathcal{A}_{H_4}}$. The representation $\rho_{30_s}$ is the unique $30$-dimensional irreducible representation therefore we have $\Phi_{1,2}\circ \rho_{30_s|\mathcal{A}_{H_4}}\simeq \rho_{30_s|\mathcal{A}_{H_4}}$. It follows by Proposition 4.1. of \cite{BMM} that up to conjugation, we have that $\rho_{16_t}(\mathcal{A}_{H_4})\leq \Omega_{16}^+(q)$ and $\rho_{30_s}(\mathcal{A}_{H_4})\leq \Omega_{30}^+(q)$. They both contain a twisted diagonal $SL_5(q)$ therefore $\rho_{16_t}(\mathcal{A}_{H_4})$ contains an element conjugate to\\ $\diag(\diag(\alpha,\alpha^{-1},1,1,1),\diag(\alpha,\alpha^{-1},1,1,1),I_6)$ and $\rho_{30_s}(\mathcal{A}_{H_4})$ contains an element conjugate to  $\diag(\diag(\alpha,\alpha^{-1},1,1,1),\diag(\alpha,\alpha^{-1},1,1,1),I_{20}).$
  It follows that the field generated by the traces of their elements contains $2(\alpha+\alpha^{-1})$. This implies that they are classical groups defined in a natural representation over $\F_q$. It follows that $\rho_{16_t}(\mathcal{A}_{F_4})\simeq \Omega_{16}^+(q)$ and $\rho_{30_s}(\mathcal{A}_{F_4})\simeq \Omega_{30}^+(q)$.
 
 By Proposition \ref{resrepreflicos}, we have $\Phi_{1,2}\circ \rho_{3_s}\simeq \rho_{\overline{3_s}}$ and $\Phi_{1,2}\circ \rho_{3_s'}\simeq \rho_{\overline{3_s'}}$. It follows that $\Phi_{1,2}\circ \rho_{24_s|\mathcal{A}_{H_4}}\simeq \rho_{\overline{24_s}|\mathcal{A}_{H_4}}$ or $\Phi_{1,2}\circ \rho_{24_s|\mathcal{A}_{H_4}}\simeq \rho_{\overline{24_t}|\mathcal{A}_{H_4}}$. We have
 \begin{eqnarray*}
 \tr(\rho_{24_s}(S_1S_3S_2^{-1}S_4^{-1})) & = & 4(\alpha^2+\alpha^{-2})+(\xi+\xi^{-1}-15)(\alpha+\alpha^{-1})+22-3(\xi+\xi^{-1})\\
 \Phi_{1,2}(\tr(\rho_{24_s}(S_1S_3S_2^{-1}S_4^{-1}))) & = & 4(\alpha^2+\alpha^{-2})+(\xi^2+\xi^{-2}-15)(\alpha+\alpha^{-1})+22-3(\xi^2+\xi^{-2})\\
 \tr(\rho_{\overline{24_t}}(S_1S_3S_2^{-1}S_4^{-1})) & = & 4(\alpha^2+\alpha^{-2})+(\xi^2+\xi^{-2}-15)(\alpha+\alpha^{-1})+21-2(\xi^2+\xi^{-2}).
 \end{eqnarray*}
 It follows that $\Phi_{1,2}\circ \rho_{24_s|\mathcal{A}_{H_4}}\simeq \rho_{\overline{24_t}|\mathcal{A}_{H_4}}$ would imply that $0=1-\xi-\xi^{-1}=-\xi^{-1}\Phi_6(\xi)$, which is absurd. This proves that $\Phi_{1,2}\circ \rho_{24_s|\mathcal{A}_{H_4}}\simeq  \rho_{\overline{24_s}|\mathcal{A}_{H_4}}$. It follows by the same arguments as before that $\rho_{24_s}(\mathcal{A}_{H_4})\simeq \Omega_{24}^+(q)$ and $\rho_{24_t}(\mathcal{A}_{H_4})\simeq \Omega_{24}^+(q)$.
 
 \medskip
 
 Assume now $1\sim 2$ and $\F_q=\F_p(\alpha)\neq \F_p(\alpha+\alpha^{-1})$. We then have by Proposition \ref{IsomorphismH3} that $\F_p(\alpha,\xi+\xi^{-1})=\F_q$ and $\Phi_{1,2}$ is the unique automorphism of order $2$ of $\F_q$. We then have using the same arguments as before that each of the groups considered preserves a non-degenerate bilinear form over $\F_{q^{\frac{1}{2}}}$. They each contain a twised diagonal $SU_5(q^{\frac{1}{2}})$, therefore they are classical groups in a natural representation over $\F_{q^{\frac{1}{2}}}$. The result then follows.
 
\medskip

Assume now $1\nsim 2$. By Lemma \ref{IsomorphismH3}, we have $\F_p(\alpha,\xi+\xi^{-1})=\F_p(\alpha)\neq \F_p(\alpha+\alpha^{-1})$. Let $\epsilon$ be the unique automorphism of order $2$ of $\F_q$. We have $\epsilon(\alpha)=\alpha^{-1}$. It follows by Lemma \ref{Fieldfactorization} that $\epsilon \circ \rho_{16_t|\mathcal{A}_{H_4}}\simeq \rho_{16_t|\mathcal{A}_{H_4}}$ and $\epsilon \circ \rho_{30_s|\mathcal{A}_{H_4}}\simeq \rho_{30_s|\mathcal{A}_{H_4}}$.

By Proposition \ref{resrepreflicos} and Table \ref{resH4H3}, we have that $\epsilon \circ \rho_{24_s|\mathcal{A}_{H_4}}\simeq \rho_{24_s|\mathcal{A}_{H_4}}$ or $\epsilon \circ \rho_{24_s|\mathcal{A}_{H_4}}\simeq \rho_{24_t|\mathcal{A}_{H_4}}$. We have 
$$\tr(\rho_{24_t}(S_1S_3S_2^{-1}S_4^{-1}))  =  4(\alpha^2+\alpha^{-2})+(\xi+\xi^{-1}-15)(\alpha+\alpha^{-1})+21-2(\xi+\xi^{-1}).$$
Since $\xi+\xi^{-1}\in \F_{q^\frac{1}{2}}$, we have $\epsilon(\xi+\xi^{-1})=\xi+\xi^{-1}$. It follows by the same computation as in the previous case that $\epsilon \circ \rho_{24_s|\mathcal{A}_{H_4}}\simeq \rho_{24_s|\mathcal{A}_{H_4}}$. The result then follows from the same arguments as before.

\bigskip

The only remaining $H_4$-graphs to consider are $\widetilde{18_r}$, $25_r$, $36_{rr}$ and $\widetilde{40_r}$. By Proposition \ref{color}, we can assume they are defined over $\F_q=\F_p(\alpha)$. By Proposition \ref{bilinwgraphs}, we have up to conjugation in $GL_{18}(q)$ that $\rho_{18_r}(\mathcal{A}_{H_4})\leq \Omega_{18}^+(q)$ and up to conjugation in $GL_{40}(q)$ that $\rho_{40_r}(\mathcal{A}_{H_4})\leq \Omega_{40}^+(q)$.

\medskip

Assume $\F_p(\alpha)=\F_p(\alpha+\alpha^{-1})$. By Table \ref{resH4H3} and Theorem \ref{resicos}, we have that each of the groups associated to those representations contains a twisted diagonal $SL_4(q)$. It follows by the same arguments as before that the field generated by the traces of their elements contains $\F_p(\alpha+\alpha^{-1})=\F_q$. By Proposition \ref{Fieldfactorization}, if there is a field automorphism $\varphi$ and a character $\eta$ such that $(\varphi\circ \rho)\otimes \eta$ is an irreducible representation of $\mathcal{H}_{H_4,\alpha}$ then $\varphi(\alpha+\alpha^{-1})=\alpha+\alpha^{-1}$. This proves that such an automorphism must be trivial. The result then follows.

\medskip

Assume now $\F_q=\F_p(\alpha)\neq \F_p(\alpha+\alpha^{-1})$. Let $\epsilon$ be the unique automorphism of order $2$ of $\F_q$. By Proposition \ref{H3dim5}, we have that $\epsilon \circ \rho_{5_r|\mathcal{A}_{H_3}}\simeq \rho_{5_r'|\mathcal{A}_{H_3}}$. It follows by Table \ref{resH4H3} and Proposition \ref{Fieldfactorization} that all the representations considered are unitary. We have by the same arguments as before that the field generated by the traces of its elements contains $\F_{q^{\frac{1}{2}}}$. It follows that $\rho_{18_r}(\mathcal{A}_{H_4})\simeq \Omega_{18}^+(q^{\frac{1}{2}})$ and $\rho_{40_r}(\mathcal{A}_{H_4})\simeq \Omega_{40}^+(q^{\frac{1}{2}})$. We also have that $\epsilon \circ \rho_{25_r|\mathcal{A}_{H_4}}\not\simeq \rho_{25_r|\mathcal{A}_{H_4}}$ and $\epsilon \circ \rho_{36_{rr}|\mathcal{A}_{H_4}}\not\simeq \rho_{36_{rr}|\mathcal{A}_{H_4}}$. It follows that they cannot be classical groups in a natural representation over $\F_{q^{\frac{1}{2}}}$. This implies that they are classical groups in a natural representation over $\F_q$, therefore $\rho_{25_r}(\mathcal{A}_{H_4})\simeq SU_{25}(q^{\frac{1}{2}})$ and $\rho_{36_{rr}}(\mathcal{A}_{H_4})\simeq SU_{36}(q^{\frac{1}{2}})$ and the proof is concluded.
\end{proof}

\smallskip

There is now only one case remaining. We cannot apply the previous theorem to it because the restriction to $\mathcal{A}_{H_3}$ only shows that $v_G(v)\leq 4$ and $\max(2,\frac{\sqrt{48}}{2})=\sqrt{12}\leq 4$.

\begin{prop} 
If $1\nsim 2$ then we have $\rho_{48_{rr}}(\mathcal{A}_{H_4})\simeq \Omega_{48}^+(q^{\frac{1}{2}})$.
\end{prop}

\begin{proof}
By Proposition \ref{color}, we can assume that the representation is defined over $\F_p(\alpha)$. Let $G=\rho_{48_r}(\mathcal{A}_{H_4})$.

\smallskip

Assume now $1\nsim 2$. By Lemma \ref{IsomorphismH3}, we have $\F_p(\alpha,\xi+\xi^{-1})=\F_q=\F_p(\alpha)\neq \F_p(\alpha+\alpha^{-1})$. Let $\epsilon$ be the unique automorphism of $\F_q$. We have $\epsilon(\alpha)=\alpha^{-1}$. It follows by Proposition \ref{Fieldfactorization} that $\epsilon \circ \rho_{48_r|\mathcal{A}_{H_4}}\simeq \rho_{48_r|\mathcal{A}_{H_4}}$. Since we have found a symmetric non-degenerate bilinear form preserved by $G$, we have by Proposition 4.1. of \cite{BMM} that up to conjugation in $GL_{48}(q)$, we have that $G\leq \Omega_{48}^+(q^{\frac{1}{2}})$. By Theorem \ref{resicos}, we have that $G$ contains a twisted diagonal $SU_3(q^{\frac{1}{2}})$ and $\rho_{48_r}(\mathcal{A}_{H_3})\simeq SU_3(q^{\frac{1}{2}})\times SU_4(q^{\frac{1}{2}})\times SU_5(q^{\frac{1}{2}})$ or $\rho_{48_r}(\mathcal{A}_{H_3})\simeq SU_3(q^{\frac{1}{2}})\times SL_4(q^{\frac{1}{2}})\times SU_5(q^{\frac{1}{2}})$. It follows by Lemma \ref{tens2} and \ref{exceptiontens2} that $G$ is tensor-indecomposable. We have with the same arguments as in the previous proposition that $G$ is a classical group in a natural representation over $\F_{q^{\frac{1}{2}}}$ therefore $G\simeq \Omega_{48}^+(q^{\frac{1}{2}})$.
\end{proof}

\begin{conjecture}\label{fingerscrossed}
If $1\sim 2$ then we have $\rho_{48_{rr}}(\mathcal{A}_{H_4})\simeq \Omega_{48}^+(q^{\frac{1}{2}})$.
\end{conjecture}

Assuming the previous proposition is correct and has been proved. We show the following theorem which gives us the image of $\mathcal{A}_{H_4}$ inside the full Iwahori-Hecke algebra.

\begin{theo}\label{resH4}
Under our assumptions on $\alpha$ and $p$, we have the following results.
\begin{enumerate}
\item Assume $1\sim 2$ and Conjecture \ref{fingerscrossed} is true.
\begin{enumerate}
\item If $\F_q=\F_p(\sqrt{\alpha})=\F_p(\alpha)=\F_p(\alpha+\alpha^{-1})$, then the morphism from $\mathcal{A}_{H_4}$ to $\mathcal{H}_{H_4,\alpha}^\star\simeq \underset{\rho~ \mbox{ irr}}\prod GL_{n_\rho}(\F_p(\sqrt{\alpha},\xi+\xi^{-1}))$ factorizes through the surjective morphism
$$\Phi :\rightarrow SL_4(q^2)\times \Omega_6^+(q^2)\times Spin_8^+(q)\times SL_9(q^2)\times \Omega_{10}^+(q)\times SL_{16}(q)^2 \times \Omega_{16}^+(q)\times \Omega_{18}^+(q)$$
$$\times \Omega_{24}^+(q)^2\times SL_{25}(q)\times \Omega_{30}^+(q)\times SL_{36}(q)\times \Omega_{40}^+(q)\times \Omega_{48}^+(q^{\frac{1}{2}}).$$
\item If $\F_q=\F_p(\sqrt{\alpha})=\F_p(\alpha)\neq \F_p(\alpha+\alpha^{-1})$ and $\Phi_{1,2}(\sqrt{\alpha})=\sqrt{\alpha}^{-1}$, then the morphism from $\mathcal{A}_{H_4}$ to $\mathcal{H}_{H_4,\alpha}^\star\simeq \underset{\rho~ \mbox{ irr}}\prod GL_{n_\rho}(\F_p(\sqrt{\alpha},\xi+\xi^{-1}))$ factorizes through the surjective morphism
$$\Phi :\rightarrow SL_4(q)\times \Omega_6^+(q)\times Spin_8^+(q^{\frac{1}{2}})\times SL_9(q^2)\times \Omega_{10}^+(q^{\frac{1}{2}})\times SU_{16}(q^{\frac{1}{2}})^2 \times \Omega_{16}^+(q^{\frac{1}{2}})\times \Omega_{18}^+(q^{\frac{1}{2}})$$
$$\times \Omega_{24}^+(q^{\frac{1}{2}})^2\times SU_{25}(q^{\frac{1}{2}})\times \Omega_{30}^+(q^{\frac{1}{2}})\times SU_{36}(q^{\frac{1}{2}})\times \Omega_{40}^+(q^{\frac{1}{2}})\times \Omega_{48}^+(q^{\frac{1}{2}}).$$
\item If $\F_q=\F_p(\sqrt{\alpha})=\F_p(\alpha)\neq \F_p(\alpha+\alpha^{-1})$ and $\Phi_{1,2}(\sqrt{\alpha})=-\sqrt{\alpha}^{-1}$, then the morphism from $\mathcal{A}_{H_4}$ to $\mathcal{H}_{H_4,\alpha}^\star\simeq \underset{\rho~ \mbox{ irr}}\prod GL_{n_\rho}(\F_p(\sqrt{\alpha},\xi+\xi^{-1}))$ factorizes through the surjective morphism
$$\Phi :\rightarrow SL_4(q)\times \Omega_6^+(q)\times Spin_8^+(q^{\frac{1}{2}})\times SL_9(q^2)\times \Omega_{10}^+(q^{\frac{1}{2}})\times SL_{16}(q) \times \Omega_{16}^+(q^{\frac{1}{2}})\times \Omega_{18}^+(q^{\frac{1}{2}})$$
$$\times \Omega_{24}^+(q^{\frac{1}{2}})^2\times SU_{25}(q^{\frac{1}{2}})\times \Omega_{30}^+(q^{\frac{1}{2}})\times SU_{36}(q^{\frac{1}{2}})\times \Omega_{40}^+(q^{\frac{1}{2}})\times \Omega_{48}^+(q^{\frac{1}{2}}).$$
\item If $\F_{q^2}=\F_p(\sqrt{\alpha}) \neq \F_p(\alpha)=\F_p(\alpha+\alpha^{-1})$, then the morphism from $\mathcal{A}_{H_4}$ to $\mathcal{H}_{H_4,\alpha}^\star\simeq \underset{\rho~ \mbox{ irr}}\prod GL_{n_\rho}(\F_p(\sqrt{\alpha},\xi+\xi^{-1}))$ factorizes through the surjective morphism
$$\Phi :\rightarrow SL_4(q^2)\times \Omega_6^+(q^2)\times Spin_8^+(q)\times SL_9(q^2)\times \Omega_{10}^+(q)\times SL_{16}(q^2) \times \Omega_{16}^+(q)\times \Omega_{18}^+(q)$$
$$\times \Omega_{24}^+(q)^2\times SL_{25}(q)\times \Omega_{30}^+(q)\times SL_{36}(q)\times \Omega_{40}^+(q)\times \Omega_{48}^+(q^{\frac{1}{2}}).$$
\end{enumerate}
\item Assume $1\nsim 2$. We write $\epsilon$ the unique automorphism of order $2$ of 
$\F_q$.
\begin{enumerate}
\item If $\F_q=\F_p(\sqrt{\alpha})=\F_p(\alpha)\neq \F_p(\alpha+\alpha^{-1})$ and $\epsilon(\sqrt{\alpha})=\sqrt{\alpha}^{-1}$, then the morphism from $\mathcal{A}_{H_4}$ to $\mathcal{H}_{H_4,\alpha}^\star\simeq \underset{\rho~ \mbox{ irr}}\prod GL_{n_\rho}(\F_p(\sqrt{\alpha},\xi+\xi^{-1}))$ factorizes through the surjective morphism$$\Phi :\rightarrow SU_4(q^{\frac{1}{2}})^2\times \Omega_6^+(q^{\frac{1}{2}})^2\times Spin_8^+(q^{\frac{1}{2}})\times SU_9(q^{\frac{1}{2}})^2\times \Omega_{10}^+(q^{\frac{1}{2}})\times SU_{16}(q^{\frac{1}{2}})^2 \times \Omega_{16}^+(q^{\frac{1}{2}})^2$$
$$\times \Omega_{18}^+(q^{\frac{1}{2}})\times \Omega_{24}^+(q^{\frac{1}{2}})^4\times SU_{25}(q^{\frac{1}{2}})\times \Omega_{30}^+(q^{\frac{1}{2}})^2\times SU_{36}(q^{\frac{1}{2}})\times \Omega_{40}^+(q^{\frac{1}{2}})\times \Omega_{48}^+(q^{\frac{1}{2}}).$$
\item If $\F_q=\F_p(\sqrt{\alpha})=\F_p(\alpha)\neq \F_p(\alpha+\alpha^{-1})$ and $\epsilon(\sqrt{\alpha})=-\sqrt{\alpha}^{-1}$, then the morphism from $\mathcal{A}_{H_4}$ to $\mathcal{H}_{H_4,\alpha}^\star\simeq \underset{\rho~ \mbox{ irr}}\prod GL_{n_\rho}(\F_p(\sqrt{\alpha},\xi+\xi^{-1}))$ factorizes through the surjective morphism
$$\Phi :\rightarrow SU_4(q^{\frac{1}{2}})^2\times \Omega_6^+(q^{\frac{1}{2}})^2\times Spin_8^+(q^{\frac{1}{2}})\times SU_9(q^{\frac{1}{2}})^2\times \Omega_{10}^+(q^{\frac{1}{2}})\times SL_{16}(q) \times \Omega_{16}^+(q^{\frac{1}{2}})^2$$
$$\times \Omega_{18}^+(q^{\frac{1}{2}})\times \Omega_{24}^+(q^{\frac{1}{2}})^4\times SU_{25}(q^{\frac{1}{2}})\times \Omega_{30}^+(q^{\frac{1}{2}})^2\times SU_{36}(q^{\frac{1}{2}})\times \Omega_{40}^+(q^{\frac{1}{2}})\times \Omega_{48}^+(q^{\frac{1}{2}}).$$
\end{enumerate}
\end{enumerate}
\end{theo}

\begin{proof}
We have the factorizations by the previous propositions. We thus only have to show that the considered morphisms are surjective. Since $\mathcal{A}_{H_4}$ is perfect, we have by Goursat's Lemma \ref{Goursat} that the morphism is surjective unless there exists two different representations $\rho_1$ and $\rho_2$ in the decomposition such that there exists a field automorphism $\Psi$ verifying $\Psi\circ \rho_{1|\mathcal{A}_{H_4}}\simeq \rho_{2|\mathcal{A}_{H_4}}$. By Proposition \ref{Fieldfactorization}, we have that $\Psi(\alpha+\alpha^{-1})=\alpha+\alpha^{-1}$. This shows that $\Psi$ must be trivial over $\F_p(\alpha+\alpha^{-1})$. It follows by the previous propositions that there are no such representations in the decompositions and the proof is concluded.
\end{proof}

\chapter{Type $F_4$}\label{TypeF4}

In this section, we determine the image of the Artin group of type $F_4$. This group contains in a natural way two groups isomorphic to $A_{B_3}$. We will therefore use some results from the type $B$ section. The $F_4$-graphs are graphs with two parameters since there are two conjugacy classes in type $F_4$. This makes the uniqueness conditions proven in the section on $W$-graphs fail. The irreducible representations are all in low dimension. We can therefore compute easily the bilinear forms for the self-dual representations and it is not necessary to use the conjecture on $W$-graphs associated to self-dual representations. It is unclear whether the conjecture is true in type $F_4$ because of the lack of rigidity when there are two parameters. The field extensions involved in type $F_4$ are quite complicated, we have to distinguish $15$ cases which are described in the Appendix. Determining the image inside each representation uses arguments similar to the other cases. The main difference comes from the fact that $\mathcal{A}_{F_4}$ is not the normal closure of $\mathcal{A}_{B_3}$ in $A_{F_4}$, we use the fact that there are two copies of $\mathcal{A}_{B_3}$ inside $\mathcal{A}_{F_4}$ for the inductive arguments. The image inside the full Iwahori-Hecke algebra is also slightly more complicated because $\mathcal{A}_{F_4}$ is not perfect. This requires some additional computations for the proof of Theorem \ref{resultF4}.

\bigskip

Let $p$ be a prime different from $2$ and $3$. Let $\alpha,\beta\in \overline{\F_p}^\star$ such that $\alpha^4\neq 1$, $\alpha^6\neq 1$, $\alpha^{10}\neq1$, $\beta^4\neq 1$, $\beta^6\neq 1$, $\beta^{10}=1$, $(\frac{\alpha}{\beta})^6\neq 1$, $(\frac{\alpha}{\beta})^4\neq 1$, $(\alpha\beta)^6 \neq 1$, $(\alpha\beta)^4 \neq 1$, $\alpha\notin \{-\beta^2,-\beta^{-2}\}$ and $\beta\notin\{-\alpha^2,-\alpha^{-2}\}$. Write $\F_q=\F_p(\alpha,\beta)$.

There are $25$ irreducible representations of $\mathcal{H}_{F_4,\alpha,\beta}$. The highest dimensional one is of dimension $16$. Five of them are self-dual and they are represented by the $F_4$-graphs $4_1$, $6_1$, $6_2$, $12$ and $16$ given in section \ref{sectionF4graphs} of the Appendix.

\begin{Def2}
The Iwahori-Hecke algebra $\mathcal{H}_{F_4,\alpha,\beta}$ of type $F_4$ is the $\F_q$-algebra generated by the generators  $S_1,S_2,S_3,S_4$ and the following relations 
\begin{enumerate}
\item $(S_1-\alpha)(S_1+1)=(S_2-\alpha)(S_2+1)=0$.
\item $(S_3-\beta)(S_3+1)=(S_4-\beta)(S_4+1)=0$.
\item $S_1S_2S_1=S_2S_1S_2$.
\item $S_1S_3=S_3S_1$.
\item $S_1S_4=S_4S_1$.
\item $S_2S_3S_2S_3=S_3S_2S_3S_2$.
\item $S_2S_4=S_4S_2$.
\item $S_3S_4S_3=S_4S_3S_4$.
For $\sigma$ in the Coxeter group $F_4$, if $\sigma=s_{i_1}\dots s_{i_k}$ is a reduced expression we set $T_{\sigma}=S_{i_1}\dots S_{i_k}$.
\end{enumerate}
\end{Def2}

Note that we have two parameters here, therefore most of the results from the Chapter \ref{Wgraphschapter} do not hold. Moreover, we have to consider the different ways $\F_p(\sqrt{\alpha},\sqrt{\beta})$ can be a field extension of $\F_p(\alpha+\alpha^{-1},\beta+\beta^{-1})$. The Hasse diagram representing the setup of the field extensions involved is given Figure \ref{fieldsF4}.

Note that all the extensions represented by edges are of degree at most $2$ because they involve of the following polynomials : $X^2-(\sqrt{\beta}+\sqrt{\beta}^{-1})X+1$, $X^2-\alpha$, $X^2-\beta$, $X^2-(\sqrt{\alpha}+\sqrt{\alpha}^{-1})X+1$, $X^2-(\beta+\beta^{-1}+2)$, $X^2-(\beta+\beta^{-1})X+1$, $X^2-(\alpha+\alpha^{-1})X+1$, $X^2-(\alpha+\alpha^{-1}+2)$. The roles of $\alpha$ and $\beta$ are perfectly symmetric in the graph, therefore we only have to consider the cases up to permutation of $\alpha$ and $\beta$. We now try to establish what all the possibilities are. Set $\F_{q'}=\F_p(\sqrt{\alpha},\sqrt{\beta})$.

\smallskip

Assume first that $\F_p(\sqrt{\alpha},\sqrt{\beta})\neq \F_p(\sqrt{\alpha},\sqrt{\beta}+\sqrt{\beta}^{-1})$. We write $\tilde{q}=q'^{\frac{1}{2}}$.
 Then $\F_{\tilde{q}^2}=\F_p(\sqrt{\alpha},\sqrt{\beta})= \F_p(\sqrt{\alpha},\sqrt{\beta}+\sqrt{\beta}^{-1})/(X^2-(\sqrt{\beta}+\sqrt{\beta}^{-1})X+1)$. The field $\F_{\tilde{q}^2}$ has a unique subfield of degree $2$ and it is equal to $\F_{\tilde{q}}$. 
 We have that $\sqrt{\beta}$ does not belong to $\F_{\tilde{q}}$, therefore $\F_{\tilde{q}}\neq \F_p(\sqrt{\alpha}+\sqrt{\alpha}^{-1},\sqrt{\beta})$ and $\F_{\tilde{q}}\neq \F_p(\alpha,\sqrt{\beta})$. 
 Note that in this case, we also have that $\beta\notin \F_{\tilde{q}}$ because otherwise, we would have $\sqrt{\beta}=\sqrt{\beta}\frac{1+\beta^{-1}}{1+\beta^{-1}}=\frac{\sqrt{\beta}+\sqrt{\beta}^{-1}}{1+\beta^{-1}}\in \F_{\tilde{q}}$. 
 It follows that $\F_{\tilde{q}}\neq \F_p(\sqrt{\alpha},\beta)$, and that $\F_{\tilde{q}^2}=\F_p(\sqrt{\alpha},\sqrt{\beta})=\F_p'\sqrt{\alpha},\beta)=\F_p(\alpha,\sqrt{\beta})=\F_p(\sqrt{\alpha}+\sqrt{\alpha}^{-1},\sqrt{\beta})$.
  We have in the same way $\F_{\tilde{q}^2}=\F_p(\sqrt{\alpha}+\sqrt{\alpha}^{-1},\beta)=\F_p(\alpha,\beta)=\F_p(\alpha+\alpha^{-1},\sqrt{\beta})=\F_p(\alpha+\alpha^{-1},\beta)$. 
  This implies that $\F_p(\sqrt{\alpha},\sqrt{\beta})$ is an extension of degree at most $2$ of $\F_p(\alpha+\alpha^{-1},\beta+\beta^{-1})$. $\F_p(\alpha+\alpha^{-1},\beta+\beta^{-1})$ is included in $\F_p(\sqrt{\alpha},\sqrt{\beta}+\sqrt{\beta}^{-1})=\F_{\tilde{q}}$, therefore it is equal to $\F_{\tilde{q}}$.
   We can now complete all the edges in the Hasse diagram where we put dotted edges for extensions of degree $1$, full edges for extensions of degree $2$, fields equal to $\F_{\tilde{q}^2}$ in blue and fields equal to $\F_{\tilde{q}}$ in red. This can be seen in Figure \ref{fieldsF4case2}.
   
   \smallskip
   
   We get by symmetry that if $\F_p(\sqrt{\alpha},\sqrt{\beta})\neq \F_p(\sqrt{\alpha}+\sqrt{\alpha}^{-1},\sqrt{\beta})$ then we have the Hasse diagram of Figure \ref{fieldsF4case3}.
   
   \smallskip
   
   We can now assume that $\F_p(\sqrt{\alpha},\sqrt{\beta})=\F_p(\sqrt{\alpha},\sqrt{\beta}+\sqrt{\beta}^{-1})=\F_p(\sqrt{\alpha}+\sqrt{\alpha}^{-1},\sqrt{\beta})$. Assume $\F_p(\sqrt{\alpha},\sqrt{\beta})\neq \F_p(\sqrt{\alpha},\beta)$. We then write $\F_{\tilde{q}}=\F_p(\sqrt{\alpha},\beta)$ and $\F_{\tilde{q}^2}=\F_p(\sqrt{\alpha},\sqrt{\beta})$. We then have that $\sqrt{\beta}\notin \F_{\tilde{q}}$, therefore $\F_{\tilde{q}^2}=\F_p(\alpha,\sqrt{\beta})=\F_p(\alpha+\alpha^{-1},\sqrt{\beta})$. Note that $\sqrt{\beta}=\frac{\sqrt{\beta}+\sqrt{\beta}^{-1}}{1+\beta^{-1}}$, therefore $\sqrt{\beta}+\sqrt{\beta}^{-1}\notin \F_{\tilde{q}}$. It follows that $\F_{\tilde{q}^2}=\F_p(\alpha,\sqrt{\beta}+\sqrt{\beta}+\sqrt{\beta}^{-1})=\F_p(\alpha+\alpha^{-1},\sqrt{\beta}+\sqrt{\beta}^{-1})$. This proves also that $\F_{\tilde{q}}=\F_p(\alpha+\alpha^{-1},\beta+\beta^{-1})$. The only Hasse diagram possible is then given in Figure \ref{fieldsF4case4}.
   
   \smallskip
   
   We get by symmetry that if $\F_p(\sqrt{\alpha},\sqrt{\beta})\neq \F_p(\alpha,\sqrt{\beta})$ then we have the Hasse diagram in Figure \ref{fieldsF4case5}.

\smallskip

We can now assume that $\F_p(\sqrt{\alpha},\sqrt{\beta})=\F_p(\sqrt{\alpha},\sqrt{\beta}+\sqrt{\beta}+\sqrt{\beta}^{-1})=\F_p(\sqrt{\alpha},\beta)=\F_p(\alpha,\sqrt{\beta})=\F_p(\sqrt{\alpha}+\sqrt{\alpha}^{-1},\sqrt{\beta})$. Assume $\F_{\tilde{q}}=\F_p(\sqrt{\alpha},\beta+\beta^{-1})\neq \F_{q'}$. We then get that $\beta\notin \F_{\tilde{q}}$ and $\sqrt{\beta}+\sqrt{\beta}^{-1}\notin \F_{\tilde{q}}$. The only Hasse diagram possible is then given in Figure \ref{fieldsF4case6}.

\smallskip

We get by symmetry that if $\F_p(\alpha+\alpha^{-1},\sqrt{\beta})\neq \F_{q'}$ then we get the Hasse diagram given in Figure \ref{fieldsF4case7}.

\smallskip

Assume now that $\F_{\tilde{q}}=\F_p(\sqrt{\alpha}+\sqrt{\alpha}^{-1},\beta)\neq \F_{q'}$. We then have that $\sqrt{\alpha}\notin \F_{\tilde{q}}$ and $\sqrt{\beta}\notin \F_{\tilde{q}}$. It follows that $\alpha\notin \F_{\tilde{q}}$ since $\sqrt{\alpha}=\frac{\sqrt{\alpha}+\sqrt{\alpha}^{-1}}{1+\alpha^{-1}}$. We also have that $\sqrt{\beta}+\sqrt{\beta}^{-1}\notin \F_{\tilde{q}}$ since $\sqrt{\beta}=\frac{\sqrt{\beta}+\sqrt{\beta}^{-1}}{1+\beta^{-1}}$. The only Hasse diagram possible is then given in Figure \ref{fieldsF4case8}.

\smallskip

By symmetry, if $\F_p(\alpha,\sqrt{\beta}+\sqrt{\beta}^{-1})\neq \F_{q'}$ then we get the Hasse diagram given in Figure \ref{fieldsF4case9}.

\smallskip

Assume now $\F_{\tilde{q}}=\F_p(\alpha,\beta)\neq \F_{q'}$. We then get that $\sqrt{\alpha}+\sqrt{\alpha}^{-1}\notin \F_{\tilde{q}}$ and $\sqrt{\beta}+\sqrt{\beta}^{-1}\notin \F_{\tilde{q}}$. We then get the Hasse diagram given in Figure \ref{fieldsF4case10}.

\smallskip

Assume now $\F_{\tilde{q}}=\F_p(\sqrt{\alpha}+\sqrt{\alpha}^{-1},\sqrt{\beta}+\sqrt{\beta}^{-1})\neq \F_{q'}$. We have $\alpha\notin \F_{\tilde{q}}$ and $\beta\notin \F_{\tilde{q}}$. We then get the Hasse diagram given in Figure \ref{fieldsF4case11}.

\smallskip

We can now assume that $\F_{q'}=\F_p(\sqrt{\alpha},\beta+\beta{-1})=\F_p(\sqrt{\alpha}+\sqrt{\alpha}^{-1},\beta)=\F_p(\sqrt{\alpha}+\sqrt{\alpha}^{-1},\sqrt{\beta}+\sqrt{\beta}^{-1})=\F_p(\alpha,\beta)=\F_p(\alpha,\sqrt{\beta}+\sqrt{\beta}^{-1})=\F_p(\alpha+\alpha^{-1})$. Assume $\F_{\tilde{q}}=\F_p(\sqrt{\alpha}+\sqrt{\alpha}^{-1},\beta+\beta^{-1})\neq \F_q'$. We then have $\alpha \notin \F_{\tilde{q}}$, $\beta\notin \F_{\tilde{q}}$ and $\sqrt{\beta}+\sqrt{\beta}^{-1}\notin \F_{\tilde{q}}$. This gives us the Hasse diagram in Figure \ref{fieldsF4case12}.

\smallskip

By symmetry, if $\F_p(\alpha+\alpha^{-1},\sqrt{\beta}+\sqrt{\beta}^{-1})\neq \F_{q'}$ then we get in Hasse diagram given in Figure \ref{fieldsF4case13}.

\smallskip

Assume now $\F_{\tilde{q}}=\F_p(\alpha+\alpha^{-1},\beta)\neq \F_{q'}$. We have $\alpha\notin \F_{\tilde{q}}$, $\sqrt{\alpha}+\sqrt{\alpha}^{-1}\notin \F_{\tilde{q}}$ and $\sqrt{\beta}+\sqrt{\beta}^{-1} \notin \F_{\tilde{q}}$. We then get the Hasse diagram given in Figure \ref{fieldsF4case14}.

\smallskip

By symmetry, if $\F_p(\alpha,\beta+\beta^{-1})\neq \F_{q'}$ we get the Hasse diagram given in Figure \ref{fieldsF4case15}.

\smallskip

We can now assume $\F_{q'}=\F_p(\sqrt{\alpha}+\sqrt{\alpha}^{-1},\beta+\beta^{-1})=\F_p(\alpha+\alpha^{-1},\beta)=\F_p(\alpha,\beta+\beta^{-1})=\F_p(\alpha+\alpha^{-1},\sqrt{\beta}+\sqrt{\beta}^{-1})$. We then either have case $1$ which is $\F_p(\sqrt{\alpha},\sqrt{\beta})=\F_p(\alpha,\beta)$ or the Hasse diagram given in Figure \ref{fieldsF4case16}.

The $F_4$-graphs we will be considering are given in Figures \ref{F4graphs1}, \ref{F4graphs2} and \ref{F4graphs3}. They are taken from \cite{G-P} (11.3.2). They are slightly different because we consider left-actions instead of right actions. We here have two parameters, therefore the rules to read the $F_4$-graphs are more complicated. When there is an edge between $x$ and $y$ in the graph then as in \cite{G-P}, we use the following conventions
\begin{enumerate}
\item $\mu_{x,y}^s=0$ if one of the following is satisfied
\begin{enumerate}
\item $I(x)$ consists of three elements, $I(y)$ consists of one element and $s=s_1$.
\item $s=s_2$, $\{s_2,s_4\}\subset I(x)$ and $I(y)$ contains $s_3$ but not $s_4$.
\end{enumerate}
\item If the above conditions are not verified then $\mu_{x,y}^s=1$ except if
\begin{enumerate}
\item The is an integer labeling the edge between $x$ and $y$, then $\mu_{x,y}^s$ is equal to that integer.
\item If $s\in \{s_3,s_4\}$ and there is a label on the edge between $x$ and $y$ which is not an integer, then $\mu_{x,y}^s$ is equal to that label, where
$$a=\sqrt{u}\sqrt{v}^{-1}+\sqrt{v}\sqrt{u}^{-1}, b=u\sqrt{v}^{-1}+\sqrt{v}u^{-1}$$
$$f=\sqrt{u}+\sqrt{u}^{-1}, g=u+u^{-1}, h=\sqrt{u}^3+\sqrt{u}^{-3}$$.
\end{enumerate}
\end{enumerate}

We now prove that the algebra is split semisimple as we did with the other types using the Schur elements.

\begin{prop2}
Under our assumptions on $p$, $\alpha$ and $\beta$, $\mathcal{H}_{F_4,\alpha,\beta}$ is split semisimple over $\F_p(\sqrt{\alpha},\sqrt{\beta})$, the representations afforded by the $W$-graphs are irreducible and pairwise non-isomorphic over $\F_q$. The restrictions of the irreducible representations of $\mathcal{H}_{F_4,\alpha,\beta}$ to $\mathcal{H}_{B_3,\alpha,\beta}$ are the same as in the generic case.
\end{prop2}

\begin{proof}
We will apply Proposition \ref{Tits}. Let $A=\Z[\sqrt{u}^{\pm 1},\sqrt{v}^{\pm 1}]$ and $F=\Q(\sqrt{u},\sqrt{v})$. We have a symetrizing trace defined by $\tau(T_0)=1$ and $\tau(T_{\sigma})=0$ for all $\sigma\in F_4\setminus \{1_{F_4}\}$. $\mathcal{H}_{F_4,u,v}$ is then a free symmetric $F$-algebra of rank $1152$. $A$ is integrally closed. Let $\theta$ be the ring homomorphism from $A$ to $L=\F_q$ defined by $\theta(\sqrt{u})=\sqrt{\alpha}$, $\theta(\sqrt{v})=\sqrt{\beta}$ and $\theta(k)=\overline{k}$. We know $FH$ is split. The basis formed by the elements $T_\sigma$, $\sigma\in F_4$ verifies the conditions of the Proposition \ref{Tits}. The $F_4$-graphs considered remain connected since the weights don't vanish after specialization. In order to verify this, we only need to check that $\theta(a)\neq 0$, $\theta(b)\neq 0$, $\theta(f)\neq 0$, $\theta(g+1)\neq 0$, $\theta(2-g)\neq 0$ and $\theta(h)\neq 0$. We have $\theta(a)=\sqrt{\alpha}\sqrt{\beta}^{-1}+\sqrt{\beta}\sqrt{\alpha}^{-1}=\sqrt{\alpha}^{-1}\sqrt{\beta}\Phi_2(\frac{\alpha}{\beta})\neq 0$, $\theta(b)=\alpha\sqrt{\beta}^{-1}+\sqrt{\beta}\alpha^{-1}=\alpha\sqrt{\beta}^{-1}(1+\beta\alpha^{-2})\neq 0$, $\theta(f)=\sqrt{\alpha}+\sqrt{\alpha}^{-1}=\sqrt{\alpha}^{-1}\Phi_2(\alpha)\neq 0$, $\theta(g+1)=\alpha+\alpha^{-1}+1=\alpha^{-1}\Phi_3(\alpha)\neq 0$, $\theta(2-g)=2-\alpha-\alpha^{-1}=(1-\alpha)(1-\alpha^{-1})$ and $\theta(\sqrt{u}^3+\sqrt{u}^{-3})=\sqrt{\alpha}^3+\sqrt{\alpha}^{-3}=\sqrt{\alpha}^{-3}(\alpha^3+1)\neq 0$.

We now only need to check that the Schur elements can be specialized and do not vanish under the specialization. This is clear from Table \ref{SchurelmtsF4}. (This table is taken from Table 11.1. \cite{G-P})

\begin{table}
\begin{enumerate}
\item $1_1$ : $\Phi_2(u)\Phi_2(v)\Phi_3(u)\Phi_3(v)(uv^2+1)(u^2v+1)\Phi_4(uv)\Phi_6(uv)\Phi_2(uv)^2$.
\item $1_3$ : $u^{-6}\Phi_2(u)\Phi_2(v)\Phi_3(u)\Phi_3(v)(u+v^2)(u^2+v)\Phi_4(\frac{v}{u})\Phi_6(\frac{v}{u})\Phi_2(\frac{v}{u})^2$.
\item $2_1$ : $v^{-2}\Phi_2(u)^2\Phi_3(u)\Phi_6(u)\Phi_3(v)\Phi_2(uv)\Phi_2(\frac{u}{v})(1+u^2v)(u^2+v)$.
\item $2_3$ : $u^{-2}\Phi_2(v)^2\Phi_3(v)\Phi_6(v)\Phi_3(u)\Phi_2(uv)\Phi_2(\frac{v}{u})(u+v^2)(uv^2+1)$.
\item $4_1$ : $2u^{-1}v^{-3}\Phi_3(u)\Phi_3(v)\Phi_2(\frac{v}{u})^2\Phi_2(uv)^2$.
\item $4_2$ : $v^{-1}\Phi_2(u)\Phi_2(v)\Phi_3(u)\Phi_3(v)\Phi_2(\frac{v}{u})\Phi_6(uv)\Phi_2(uv)^2$.
\item $4_4$ : $u^{-3}v^{-1}\Phi_2(u)\Phi_2(v)\Phi_3(u)\Phi_3(v)\Phi_2(uv)\Phi_6(\frac{v}{u})\Phi_2(\frac{v}{u})^2$.
\item $6_1$ : $3u^{-1}v^{-3}\Phi_2(u)^2\Phi_2(v)^2\Phi_6(\frac{v}{u})\Phi_2(uv)^2$.
\item $6_2$ : $3u^{-1}v^{-3}\Phi_2(u)^2\Phi_2(v)^2\Phi_6(uv)\Phi_2(\frac{v}{u})^2$.
\item $8_1$ : $u^{-1}v^{-3}\Phi_2(u)^2\Phi_3(u)\Phi_6(u)\Phi_3(v)(u+v^2)(uv^2+1)$.
\item $8_3$ : $u^{-3}v^{-1}\Phi_3(u)\Phi_2(v)^2\Phi_3(v)\Phi_6(v)(u^2+v)(u^2v+1)$.
\item $9_1$ : $(uv)^{-2}\Phi_2(u)\Phi_2(v)(u+v^2)(u^2+v)\Phi_4(uv)\Phi_2(uv)^2$.
\item $9_2$ : $(uv)^{-2}\Phi_2(u)\Phi_2(v)(1+uv^2)(1+u^2v)\Phi_4(\frac{u}{v})\Phi_2(\frac{v}{u})^2$.
\item $12$ : $6u^{-1}v^{-3}\Phi_6(u)\Phi_6(v)\Phi_2(\frac{v}{u})^2\Phi_2(uv)^2$.
\item $16$ : $2u^{-1}v^{-3}\Phi_3(u)\Phi_3(v)\Phi_4(uv)\Phi_4(\frac{v}{u})$.
\end{enumerate}
\caption{Schur elements in type F4}\label{SchurelmtsF4}
\end{table}
\end{proof}

\begin{prop2}\label{resF4derivedsubgroup}
The restrictions to $\mathcal{A}_{F_4}$ of the representations afforded by those $W$-graphs are absolutely irreducible and the representations of dimension greater than one are pairwise non-isomorphic.
\end{prop2}

\begin{proof}
As in \cite{BMM} Lemma $3.4$, we only need to prove that $A_{F_4}=<s_1,s_2,s_3,s_4>$ is generated by $A_{B_3}=<s_1,s_2,s_3>$ and $\mathcal{A}_{F_4}$. This true because $s_4=s_4s_3^{-1}s_3$, $s_4s_3^{-1}\in \mathcal{A}_{F_4}$ and $s_3\in A_{B_3}$.

Let now $\rho$ and $\rho'$ be two representations of $\mathcal{H}_{F_4,\alpha,\beta}$ such that $\rho_{|\mathcal{A}_{F_4}}\simeq \rho'_{|\mathcal{A}_{F_4}}$. We know by Lemma \ref{abel} that there exists a character $\chi;A_{F_4}\rightarrow \F_q$ such that $\rho\simeq \rho'\otimes \chi$. Since $F_4/\mathcal{A}_{F_4}\simeq <\overline{s_1},\overline{s_3}> \simeq \Z^2$, there exists $(x,y)\in \F_q^2$ such that $\rho(S_1)=x\rho'(S_1)$ and $\rho(S_3)=y\rho'(S_3)$. The eigenvalues of $\rho(S_1)$ are $-1$ and $\alpha$, therefore $\{-1,\alpha\}=\{-x,x\alpha\}$. It follows that $x=1$ since $\alpha^2\neq 1$. In the same way $y=1$ since $\beta^2\neq 1$. It follows that $\rho\simeq \rho'$.
\end{proof}

Before determining the image of the Artin groups inside this Iwahori-Hecke algebra, we need as in the other cases a Lemma on Artin groups which will allow us to use the restriction from $E_6$ to $D_5$.

\begin{lemme2}\label{genB3F4}
$\mathcal{A}_{F_4}$ is generated by $\mathcal{A}_{B_3,1}$ and $\mathcal{A}_{B_3,2}$, where we identify $A_{B_3,1}$ (resp $A_{B_3,2})$ as a subgroup of $A_{F_4}$ using the natural isomorphism from $A_{B_3}$ to $<s_1,s_2,s_3>$ (resp $<s_2,s_3,s_4>$).
\end{lemme2}

\begin{proof}
By \cite{MR}, we have that $\mathcal{A}_{F_4}$ is generated by $s_2s_1^{-1}$, $s_1s_2s_1^{-2}$, $s_3s_4^{-1}$ and $s_4s_3s_4^{-2}$. We have that $s_2s_1^{-1}$ and $s_1s_2s_1^{-2}$ belong to $\mathcal{A}_{B_3,1}$ and $s_3s_4^{-1}$ and $s_4s_3s_4^{-2}$ are in $\mathcal{A}_{B_3,2}$, therefore the proof is complete.
\end{proof}

We now give a proposition where we determine the image of $\mathcal{A}_{F_4}$ with respect to the representations $\rho_{2_3}$, $\rho_{4_1}$ and $\rho_{4_2}$. The special phenomenon appearing for the representation $\rho_{4_1}$ was also observed in the generic case, see Lemma $2.22$ of \cite{IH2} for an analogous result.

\begin{prop2}
If $\F_{\tilde{q}}=\F_p(\alpha)=\F_p(\alpha+\alpha^{-1})$ then $\rho_{2_3}(\mathcal{A}_{F_4})\simeq SL_2(\tilde{q})$.

If $\F_{\tilde{q}}=\F_p(\alpha)\neq \F_p(\alpha+\alpha^{-1})$ then $\rho_{2_3}(\mathcal{A}_{F_4})\simeq SU_2(\tilde{q}^{\frac{1}{2}})$.

If $\F_{\tilde{q}}=\F_p(\beta)=\F_p(\beta+\beta^{-1})$ then $\rho_{2_3}(\mathcal{A}_{F_4})\simeq SL_2(\tilde{q})$.

If $\F_{\tilde{q}}=\F_p(\beta)\neq \F_p(\beta+\beta^{-1})$ then $\rho_{2_3}(\mathcal{A}_{F_4})\simeq SU_2(\tilde{q}^{\frac{1}{2}})$.

If $\F_{\tilde{q}_1}=\F_p(\alpha)=\F_p(\alpha+\alpha^{-1})$ and $\F_{\tilde{q}_2}=\F_p(\beta)=\F_p(\beta+\beta^{-1})$ then $\rho_{4_1}(\mathcal{A}_{F_4})\simeq SL_2(\tilde{q}_1)\circ SL_2(\tilde{q}_2)$.

If $\F_{\tilde{q}_1}=\F_p(\alpha)\neq \F_p(\alpha+\alpha^{-1})$ and $\F_{\tilde{q}_2}=\F_p(\beta)=\F_p(\beta+\beta^{-1})$ then $\rho_{4_1}(\mathcal{A}_{F_4})=SU_2(\tilde{q}_1^{\frac{1}{2}})\circ SL_2(\tilde{q}_2)$.

If $\F_{\tilde{q}_1}=\F_p(\alpha)\neq \F_p(\alpha+\alpha^{-1})$ and $\F_{\tilde{q}_2}=\F_p(\beta)\neq\F_p(\beta+\beta^{-1})$ then $\rho_{4_1}(\mathcal{A}_{F_4})=SU_2(\tilde{q}_1^{\frac{1}{2}})\circ SU_2(\tilde{q}_2^{\frac{1}{2}})$.

In cases $11$, $12$, $13$ and $16$, we have $\rho_{4_2}(\mathcal{A}_{F_4})\simeq SU_4(q^{\frac{1}{2}})$ and $\rho_{4_4}(\mathcal{A}_{F_4})\simeq SU_4(q^{\frac{1}{2}})$. In all the remaining cases, we have $\rho_{4_2}(\mathcal{A}_{F_4})\simeq SL_4(q)$ and $\rho_{4_4}(\mathcal{A}_{F_4})\simeq SL_4(q)$.

\end{prop2}

\begin{proof}

We have that $\rho_{2_3|\mathcal{H}_{A_2,1}}$ is isomorphic to the representation $\rho_{[2,1]}$ of $\mathcal{H}_{A_2,\alpha}$ and $\rho_{2_1|\mathcal{H}_{A_2,2}}$ is isomorphic to the representation $\rho_{[2,1]}$ of $\mathcal{H}_{A_2,\beta}$, where $A_{A_2,1}=<S_1,S_2>$ and $A_{A_2,2}=<S_3,S_4>$. We have $\rho_{2_3}(S_3)=\rho_{2_3}(S_4)=\begin{pmatrix}
\beta & 0\\
0 & \beta
\end{pmatrix}$ and $\rho_{2_1}(S_1)=\rho_{2_1}(S_2)=\begin{pmatrix}
\alpha & 0\\
0 & \alpha
\end{pmatrix}$. It then follows that $\rho_{2_3}(\mathcal{A}_{F_4})=\rho_{2_3}(\mathcal{A}_{A_2,1})$ and $\rho_{2_3}(\mathcal{A}_{F_4})=\rho_{2_3}(\mathcal{A}_{A_2,2})$. The result then follows from Lemma 3.5 of \cite{BM}.

\smallskip

We now consider the representation $\rho_{4_1}$. We order the basis indexed by the vertices in the lexicographic way, i.e. $I(x_1)=\{s_1,s_3\}$, $I(x_2)=\{s_1,s_4\}$, $I(x_3)=\{s_2,s_3\}$ and $I(x_4)=\{s_2,s_4\}$. Let $X=\begin{pmatrix}
0 & 0 & 0 & 1\\
0 & 1 & 0 & 0\\
0 & 0 & 1 & 0\\
1 & 0 & 0 & 0
\end{pmatrix}$. We then have $\rho_{4_1|\mathcal{H}_{A_2,1}}=X\begin{pmatrix}
\rho_{2_3|\mathcal{H}_{A_2,1}} & 0\\
0 & \rho_{2_3|\mathcal{H}_{A_2,1}}
\end{pmatrix}X^{-1}$ and $\rho_{4_1|\mathcal{H}_{A_2,2}}=\begin{pmatrix}
\rho_{2_1|\mathcal{H}_{A_2,2}} & 0\\
0 & \rho_{2_1|\mathcal{H}_{A_2,2}}
\end{pmatrix}$. Note now that for any matrices $M,N\in GL_2(\overline{\F_p})$, we have that
$[X\begin{pmatrix}
M & 0\\
0 & M
\end{pmatrix}X^{-1},\begin{pmatrix}
N & 0\\
0 & N
\end{pmatrix}]=I_4$. We also have that if $\{X\begin{pmatrix}
N_1 & 0\\
0 & N_1
\end{pmatrix}X^{-1},N_1\in SL_2(\overline{\F_p})\}\cap \{\begin{pmatrix}
N_2 & 0\\
0 & N_2
\end{pmatrix}, N_2\in SL_2(\overline{\F_p})\}=\{\pm I_4\}$. We have by \cite{MR} that $\mathcal{A}_{F_4}$ is generated by $s_2s_1^{-1}$, $s_1s_2s_1^{-2}$, $s_3s_4^{-1}$ and $s_4s_3s_4^{-2}$. The result then follows from what was proven above for $2_1$ and $2_3$.

\smallskip

Consider now the representation $\rho_{4_2}$. We know, using Table \ref{resF44B31} and Table \ref{resF44B32}, that $\rho_{4_2|\mathcal{H}_{B_{3,1}}}$ is isomorphic to the representation $\rho_{[2],[1]}\otimes \rho_{[3],\emptyset}$ of $\mathcal{H}_{B_3,\alpha,\beta}$ and $\rho_{4_2|\mathcal{H}_{B_{3,2}}}$ is isomorphic to the representation $\rho_{[2],[1]}\otimes \rho_{[3],\emptyset}$ of $\mathcal{H}_{B_3,\beta,\alpha}$.

By Lemmas \ref{platypus} and \ref{platypus456}, $\rho_{4_2}(\mathcal{A}_{B_3,1})\in \{SL_3(q), SU_3(q^{\frac{1}{2}})\}$, $\rho_{4_2}(\mathcal{A}_{B_3,2})\in \{SL_3(q), SU_3(q^{\frac{1}{2}})\}$ and $\rho_{4_2}(\mathcal{A}_{B_2}) \in \{SL_2(q),SL_2(q^{\frac{1}{2}}), SU_2(q^{\frac{1}{2}})\}$. Take now a transvection $t$ of $\rho_{4_2}(\mathcal{A}_{B_2})$. We know that for $i\in \{1,2\}$, $\rho_{4_2}(\mathcal{A}_{B_3,i})$ is normally generated by that transvection in $\rho_{4_2}(\mathcal{A}_{B_3,i}$. Since $\mathcal{A}_{F_4}$ is generated by $\mathcal{A}_{B_3,1}$ and $\mathcal{A}_{B_3,2}$, we have that the normal closure in $\rho_{4_2}(\mathcal{A}_{F_4})$ of $t$ is equal to $\rho_{4_2}(\mathcal{A}_{F_4})$. This proves that $\rho(\mathcal{A}_{F_4})$ is an irreducible subgroup of $SL_4(\F_p(\sqrt{\alpha},\sqrt{\beta}))$ generated by transvections. We can now use Theorem \ref{transvections}. Let $\F_{q_0}=\F_p(\sqrt{\alpha},\sqrt{\beta})$. There exists $q'$ dividing $q_0$ such that $\rho_{4_2}(\mathcal{A}_{F_4})$ is conjugate in $GL_4(q_0)$ to $SL_4(q')$, $SU_4(q'^{\frac{1}{2}})$ or $SP_4(q')$. We also have that $\rho_{4_2}\not\simeq \rho_{4_2}^\star$, therefore the symplectic case is excluded. Note that $\rho_{4_2}(\mathcal{A}_{F_4})$ contains either a natural $SL_3(q)$ or a natural $SU_3(q^{\frac{1}{2}})$, therefore we have that $q$ divides $q'$ and $q'\in \{q,q_0\}$.

\begin{figure}
\begin{tiny}
\begin{tabular}{ |p{0.4cm}|p{1.1cm}|p{1.1cm}|p{1.1cm}|p{1.1cm}|p{1.1cm}|p{1.1cm}|p{1.1cm}|p{1.1cm}|p{1.1cm}|p{1.1cm}|  }
 \hline
 & $(\emptyset,[1^3])$ & $([3],\emptyset)$ & $([1^3],\emptyset)$ & $(\emptyset,[3])$ & $([2,1],\emptyset)$ & $(\emptyset,[2,1])$ & $([2],[1])$ & $([1],[1^2])$ & $([1],[2])$ &  $([1^2],[1])$  \\
 \hline

 $4_1$   &   &  &    &  & $1$  &  $1$&  &  &  &  \\
 \hline
 $4_2$    &   & $1$ &    &  &  &  & $1$  &  &  &  \\
 \hline
 $4_2'$   & $1$  &  &    &  &  &  &  & $1$ &  &  \\
 \hline
 $4_4$  &   &  &    & $1$ &  &  &  &  & $1$ &  \\
 \hline
 $4_4'$ &   &  & $1$   &  &  &  &  &  &  & $1$ \\
 \hline
\end{tabular}
\end{tiny}
\caption{Restriction of the $4$-dimensional representations to $\mathcal{H}_{B_3,1}\simeq \mathcal{H}_{B_3,\alpha,\beta}$}\label{resF44B31}
\end{figure}

\begin{figure}
\begin{tiny}
\begin{tabular}{ |p{0.4cm}|p{1.1cm}|p{1.1cm}|p{1.1cm}|p{1.1cm}|p{1.1cm}|p{1.1cm}|p{1.1cm}|p{1.1cm}|p{1.1cm}|p{1.1cm}|  }
 \hline
 & $(\emptyset,[1^3])$ & $([3],\emptyset)$ & $([1^3],\emptyset)$ & $(\emptyset,[3])$ & $([2,1],\emptyset)$ & $(\emptyset,[2,1])$ & $([2],[1])$ & $([1],[1^2])$ & $([1],[2])$ &  $([1^2],[1])$  \\
 \hline

 $4_1$   &   &  &    &  & $1$  &  $1$&  &  &  &  \\
 \hline
 $4_2$    &   & $1$ &    &  &  &  & $1$  &  &  &  \\
 \hline
 $4_2'$   & $1$  &  &    &  &  &  &  & $1$ &  &  \\
 \hline
 $4_4$  &   &  &   $1$ &  &  &  &  &  &  & $1$ \\
 \hline
 $4_4'$ &   &  &   & $1$ &  &  &  &  & $1$ &  \\
 \hline
\end{tabular}
\end{tiny}
\caption{Restriction of the $4$-dimensional representations to $\mathcal{H}_{B_3,2}\simeq \mathcal{H}_{B_3,\beta,\alpha}$}\label{resF44B32}
\end{figure}

It now only remains to show that we have the correct groups. Assume $\F_{q^2}=\F_p(\sqrt{\alpha},\sqrt{\beta})\neq \F_p(\alpha,\beta)=\F_q$. We then have $\F_q=\F_p(\alpha+\alpha^{-1},\beta+\beta^{-1})$, therefore $\rho_{4_2}(\mathcal{A}_{F_4})$ contains a natural $SL_3(q)$. There exists then a unique automorphism $\epsilon$ of order $2$ of $\F_{q^2}$. Since $\F_q$ is fixed by $\epsilon$, we have that $\epsilon\circ \rho_{4_2}\simeq \rho_{4_2}$ by considering the restriction to $\mathcal{A}_{B_3,1}$. It follows by Lemma \ref{Ngwenya} that up to conjugation, we have $\rho_{4_2}(\mathcal{A}_{F_4})\leq SL_4(q)$. This proves that $q'=q$ and $\rho_{4_2}(\mathcal{A}_{F_4})\simeq SL_4(q)$ in cases $4$, $5$ and $10$. We can now assume $\F_p(\sqrt{\alpha},\sqrt{\beta})=\F_p(\alpha,\beta)$, this implies that $q_0=q$, therefore $q'=q$. This implies that $\rho_{4_2}(\mathcal{A}_{F_4})$ is conjugate in $GL_4(q)$ to $SL_4(q)$ or $SU_4(q^{\frac{1}{2}})$. If $\F_p(\alpha,\beta)=\F_p(\alpha+\alpha^{-1},\beta+\beta^{-1})$ then we have again $q'=q$ and $\rho_{4_2}(\mathcal{A}_{F_4})\simeq SL_4(q)$. If $\F_p(\alpha,\beta)=\F_p(\alpha+\alpha^{-1},\beta)\neq\F_p(\alpha+\alpha^{-1},\beta+\beta^{-1})$ then we have that there exists a unique automorphism $\epsilon$ of order $2$ of $\F_q$ and that for all representation $\varphi$ of $\mathcal{H}_{B_3,\alpha,\beta}$, we have $\epsilon\circ \varphi^\star\simeq\varphi$. This implies that $\epsilon \circ \rho_{4_2}^\star \simeq \rho_{4_2}$. It follows that $\rho_{4_2}(\mathcal{A}_{F_4})\simeq SU_4(q^{\frac{1}{2}})$. Finally, if $\F_p(\alpha,\beta)=\F_p(\alpha,\beta+\beta^{-1})\neq \F_p(\alpha+\alpha^{-1},\beta)=\F_p(\alpha+\alpha^{-1},\beta+\beta^{-1})$ then $\rho_{4_2}(\mathcal{A}_{F_4})$ contains a natural $SL_3(q)$, therefore $\rho_{4_2}(\mathcal{A}_{F_4})\simeq SL_4(q)$. The only remaining possibility is $\F_p(\alpha,\beta)=\F_p(\alpha+\alpha^{-1},\beta)\neq \F_p(\alpha,\beta+\beta^{-1})=\F_p(\alpha+\alpha^{-1},\beta+\beta^{-1})$. The arguments are identical to the previous case. To conclude, we have $\rho_{4_2}(\mathcal{A}_{F_4})\simeq SU_4(q^{\frac{1}{2}})$ in cases $11$, $12$, $13$ and $16$ and $\rho_{4_2}(\mathcal{A}_{F_4})\simeq SL_4(q)$ otherwise. By the restriction table, we can  apply the same reasoning to $\rho_{4_4}$ and we get the same result.
\end{proof}

\begin{prop2}
In cases $1$, $4$, $5$ and $10$, we have $\rho_{8_2}(\mathcal{A}_{F_4})\simeq SL_8(q)$ and $\rho_{8_3}(\mathcal{A}_{F_4})\simeq SL_8(q)$.

In cases $11$, $12$, $13$ and $16$, we have $\rho_{8_2}(\mathcal{A}_{F_4})\simeq SU_8(q^{\frac{1}{2}})$ and $\rho_{8_3}(\mathcal{A}_{F_4})\simeq SU_8(q^{\frac{1}{2}})$.

In cases $2$, $6$, $9$ and $15$, we have $\rho_{8_2}(\mathcal{A}_{F_4})\simeq SU_8(q^{\frac{1}{2}})$ and $\rho_{8_3}(\mathcal{A}_{F_4})\simeq SL_8(q^{\frac{1}{2}})$.

In cases $3$, $7$, $8$ and $14$, we have $\rho_{8_2}(\mathcal{A}_{F_4})\simeq SL_8(q^{\frac{1}{2}})$ and $\rho_{8_3}(\mathcal{A}_{F_4})\simeq SU_8(q^{\frac{1}{2}})$.
\end{prop2}

\begin{proof}

\begin{table}
\begin{tiny}
\centering
\begin{tabular}{ |p{0.4cm}|p{1.1cm}|p{1.1cm}|p{1.1cm}|p{1.1cm}|p{1.1cm}|p{1.1cm}|p{1.1cm}|p{1.1cm}|p{1.1cm}|p{1.1cm}|  }
 \hline
 & $(\emptyset,[1^3])$ & $([3],\emptyset)$ & $([1^3],\emptyset)$ & $(\emptyset,[3])$ & $([2,1],\emptyset)$ & $(\emptyset,[2,1])$ & $([2],[1])$ & $([1],[1^2])$ & $([1],[2])$ &  $([1^2],[1])$  \\
 \hline

 $8_2$   &   &  &    &  &  $1$ &  &$1$  &  &  & $1$ \\
 \hline
 $8_2'$    &   & &    &  &  & $1$ &  & $1$ & $1$ &  \\
 \hline
 $8_3$   &   & $1$ &    & $1$ &  &  & $1$ & & $1$ &  \\
 \hline
 $8_3'$  &  $1$ &  & $1$   &  &  &  &  & $1$ &  & $1$ \\
 \hline
\end{tabular}
\end{tiny}
\smallskip
\caption{Restriction of the $8$-dimensional representations to $\mathcal{H}_{B_3,1}\simeq \mathcal{H}_{B_3,\alpha,\beta}$.}\label{resF48B31}
\end{table}

\begin{table}
\begin{tiny}
\centering
\begin{tabular}{ |p{0.4cm}|p{1.1cm}|p{1.1cm}|p{1.1cm}|p{1.1cm}|p{1.1cm}|p{1.1cm}|p{1.1cm}|p{1.1cm}|p{1.1cm}|p{1.1cm}|  }
 \hline
 & $(\emptyset,[1^3])$ & $([3],\emptyset)$ & $([1^3],\emptyset)$ & $(\emptyset,[3])$ & $([2,1],\emptyset)$ & $(\emptyset,[2,1])$ & $([2],[1])$ & $([1],[1^2])$ & $([1],[2])$ &  $([1^2],[1])$  \\
 \hline

 $8_2$   &   & $1$ &    & $1$ &  &  & $1$ & & $1$ &  \\
 \hline
 $8_2'$    &  $1$ &  & $1$   &  &  &  &  & $1$ &  & $1$ \\
 \hline
 $8_3$   &   &  &    &  &  $1$ &  &$1$  &  &  & $1$ \\
 \hline
 $8_3'$  &   & &    &  &  & $1$ &  & $1$ & $1$ &  \\
 \hline
\end{tabular}
\end{tiny}
\smallskip
\caption{Restriction of the $8$-dimensional representations to $\mathcal{H}_{B_3,2}\simeq \mathcal{H}_{B_3,\beta,\alpha}$.}\label{resF48B32}
\end{table}

Let $\F_{\tilde{q}}=\F_p(\sqrt{\alpha},\sqrt{\beta})$ and $G_{8_2}=\rho_{8_2}(\mathcal{A}_{F_4})$. Consider now the representation $\rho_{8_2}$. We have in cases $1$, $4$, $5$, $10$, $11$, $12$, $13$ and $16$ by Tables \ref{resF48B31} and \ref{resF48B32} and Lemmas \ref{platypus} and \ref{platypus456} that $\rho_{8_2}(\mathcal{A}_{B_3,1})$ and $\rho_{8_2}(\mathcal{A}_{B_3,2})$ are generated by transvections. Since $\mathcal{A}_{F_4}$ is generated by $\mathcal{A}_{B_3,1}$ and $\mathcal{A}_{B_3,2}$, we have that $G_{8_2}$ is  generated by transvections. We can then apply Theorem \ref{transvections}. We get that there exists $q'$ dividing $\tilde{q}$ such that $\rho_{8_2}(\mathcal{A}_{F_4})$ is conjugate in $GL_8(\tilde{q})$ to $SL_8(q')$, $SU_8(q'^{\frac{1}{2}})$ or $SP_8(q')$. The symplectic case is exluded because $\rho_{8_2}\not\simeq \rho_{8_2}^\star$. In cases $4$, $5$ and $10$, we have $\tilde{q}=q^2$ and $\rho_{8_2}(\mathcal{A}_{B_3,1})\simeq SL_3(q)^2\times SL_2(r)$ if $\F_r=\F_p(\alpha)=\F_p(\alpha+\alpha^{-1})$ and $\rho_{8_2}(\mathcal{A}_{B_3,1})\simeq SL_3(q)^2\times SU_2(r^{\frac{1}{2}})$ if $\F_r=\F_p(\alpha)\neq \F_p(\alpha+\alpha^{-1})$. It follows that $q$ divides $q'$ since $G_{8_2}$ contains a natural $SL_3(q)$. Since $\tilde{q}=q^2$, there exists a unique automorphism $\Phi$ of order $2$ of $\F_{q^2}$. For any representation $\varphi$ of $\mathcal{H}_{B_3,1}$, we have $\Phi \circ \varphi \simeq \varphi$, therefore by Table \ref{resF48B31}, we have that $\Phi\circ \rho_{8_2}\simeq \rho_{8_2}$. It follows by Lemma \ref{Ngwenya} that up to conjugation in $GL_8(q^2)$, $G_{8_2}$ is a subgroup of $SL_8(q)$. Furthermore, $G_{8_2}\simeq SL_8(q)$ because $SU_8(q^{\frac{1}{2}})$ is not conjugate in $GL_8(q^2)$ to a subgroup of $SL_8(q)$.

 In cases $1$, $11$, $12$, $13$ and $16$, we have $\tilde{q}=q$ and $G_{8_2}$ contains either a natural $SL_3(q)$ or a natural $SU_3(q^{\frac{1}{2}})$. We then have by Lemma \ref{field} that $q'=q$. Hence $G_{8_2}=SL_8(q)$ in case $1$ and $G_{8_2}\simeq SU_8(q^{\frac{1}{2}})$ in cases $11$, $12$, $13$ and $16$.

\smallskip

Assume now $\F_p(\sqrt{\alpha},\sqrt{\beta})=\F_p(\alpha,\beta)=\F_q=\F_p(\alpha+\alpha^{-1},\beta)\neq \F_p(\alpha,\beta+\beta^{-1})=\F_p(\alpha+\alpha^{-1},\beta+\beta^{-1})$. Then $\rho_{8_2}(\mathcal{A}_{B_3,1})\simeq SL_3(q)\times SU_2(r^{\frac{1}{2}})$ if $\F_r=\F_p(\alpha)\neq \F_p(\alpha+\alpha^{-1})$ and $\rho_{8_2}(\mathcal{A}_{B_3,1})\simeq SL_3(q)\times SL_2(r)$ if $\F_r=\F_p(\alpha)=\F_p(\alpha+\alpha^{-1})$. There exists a unique automorphism $\epsilon$ of order $2$ of $\F_q$ and if we consider the representations appearing in the restriction of $\rho_{8_2}$ to $\mathcal{H}_{B_3,1}$,   then $\epsilon \circ \rho_{[2],[1]}\simeq \rho_{[1],[2]}$ and $\epsilon \circ \rho_{[1^2],[1]}\simeq \rho_{[1],[1^2]}$ by Proposition \ref{isomorphisme}. Hence Table \ref{resF48B31} gives $\epsilon \circ \rho_{8_2}\simeq \rho_{8_2'}\simeq \rho_{8_2}^\star$, and it follows that $\epsilon \circ \rho_{8_2}^\star \simeq \rho_{8_2}$. By Lemma \ref{Harinordoquy}, we have that up to conjugation in $GL_8(q)$, $G_{8_2}$ is a subgroup of $SU_8(q^{\frac{1}{2}})$. By the exact same arguments as in the proof of Lemma \ref{platypus456}, we conclude that $G_{8_2}\simeq SU_8(q^{\frac{1}{2}})$. This corresponds to cases $2$, $6$, $9$ and $15$.

\smallskip

Assume now $\F_p(\sqrt{\alpha},\sqrt{\beta})=\F_p(\alpha,\beta)=\F_q=\F_p(\alpha,\beta+\beta^{-1})\neq \F_p(\alpha+\alpha^{-1},\beta)=\F_p(\alpha+\alpha^{-1},\beta+\beta^{-1})$. Then $\rho_{8_2}(\mathcal{A}_{B_3,1})\simeq SL_3(q)\times SU_2(r^{\frac{1}{2}})$, where $\F_r=\F_p(\alpha)\neq \F_p(\alpha+\alpha^{-1})$. There exists a unique automorphism $\epsilon$ of order $2$ of $\F_q$ and if we consider the representations appearing in the restriction of $\rho_{8_2}$ to $\mathcal{H}_{B_3,1}$, $\epsilon \circ \rho_{[2,1]} \simeq \rho_{[1^2],[1]}$ and $\epsilon \circ \rho_{[1^2],[1]}\simeq \rho_{[2],[1]}$ by Proposition \ref{isomorphisme}. It follows by Table \ref{resF48B31} that $\epsilon \circ \rho_{8_2}\simeq \rho_{8_2}$. We then have by Lemma \ref{Ngwenya} that up to conjugation in $GL_8(q)$, $G_{8_2}$ is a subgroup of $SL_8(q^{\frac{1}{2}})$. We can then again use the arguments of the proof of Lemma \ref{platypus456} to get that $G_{8_2}\simeq SL_8(q^{\frac{1}{2}})$. This corresponds to cases $3$, $7$, $8$ and $14$.

\medskip

The results for the representation $\rho_{8_3}$ are symmetric with regards to $\alpha$ and $\beta$ to the results for $\rho_{8_2}$ since Tables \ref{resF48B31} and \ref{resF48B32} are identical after permutation of $8_2$ and $8_3$.
\end{proof}

\begin{prop2}
In cases $11$, $12$, $13$ and $16$, we have $\rho_{9_1}(\mathcal{A}_{F_4})\simeq SU_9(q^{\frac{1}{2}})$ and $\rho_{9_2}(\mathcal{A}_{F_4})\simeq SU_9(q^{\frac{1}{2}})$.

In all the remaining case, we have $\rho_{9_1}(\mathcal{A}_{F_4})\simeq SL_9(q)$ and $\rho_{9_2}(\mathcal{A}_{F_4})\simeq SL_9(q)$.

In cases $2$, $6$, $9$ and $15$, if $\epsilon$ is the unique automorphism of order $2$ of $\F_q$ then $\epsilon \circ \rho_{9_1|\mathcal{A}_{F_4}} \simeq \rho_{9_2|\mathcal{A}_{F_4}}$.

In cases $3$, $7$, $8$ and $14$, if $\epsilon$ is the unique automorphism of order $2$ of $\F_q$ then $\epsilon \circ \rho_{9_1|\mathcal{A}_{F_4}} \simeq \rho_{9_2'|\mathcal{A}_{F_4}}$.
\end{prop2}

\begin{proof}

\begin{table}
\begin{tiny}
\centering
\begin{tabular}{ |p{0.4cm}|p{1.1cm}|p{1.1cm}|p{1.1cm}|p{1.1cm}|p{1.1cm}|p{1.1cm}|p{1.1cm}|p{1.1cm}|p{1.1cm}|p{1.1cm}|  }
 \hline
 & $(\emptyset,[1^3])$ & $([3],\emptyset)$ & $([1^3],\emptyset)$ & $(\emptyset,[3])$ & $([2,1],\emptyset)$ & $(\emptyset,[2,1])$ & $([2],[1])$ & $([1],[1^2])$ & $([1],[2])$ &  $([1^2],[1])$  \\
 \hline

 $9_1$   &   & $1$ &    &  & $1$  &  & $1$&  & $1$ &  \\
 \hline
 $9_1'$    &  $1$ & &    &  &  & $1$ &  & $1$ &  & $1$ \\
 \hline
 $9_2$   &   &  &    & $1$ &  & $1$ & $1$ & & $1$ &  \\
 \hline
 $9_2'$  &   &  & $1$ &  & $1$ &  &  & $1$ &  & $1$ \\
 \hline
\end{tabular}
\end{tiny}
\smallskip

\caption{Restriction of the $9$-dimensional representations to $\mathcal{H}_{B_3,1}\simeq \mathcal{H}_{B_3,\alpha,\beta}$.}\label{resF49B31}
\end{table}

\begin{table}
\begin{tiny}
\centering
\begin{tabular}{ |p{0.4cm}|p{1.1cm}|p{1.1cm}|p{1.1cm}|p{1.1cm}|p{1.1cm}|p{1.1cm}|p{1.1cm}|p{1.1cm}|p{1.1cm}|p{1.1cm}|  }
 \hline
 & $(\emptyset,[1^3])$ & $([3],\emptyset)$ & $([1^3],\emptyset)$ & $(\emptyset,[3])$ & $([2,1],\emptyset)$ & $(\emptyset,[2,1])$ & $([2],[1])$ & $([1],[1^2])$ & $([1],[2])$ &  $([1^2],[1])$  \\
 \hline

 $9_1$   &   &  &    & $1$ &  & $1$ & $1$ & & $1$ &  \\
 \hline
 $9_1'$     &   &  & $1$ &  & $1$ &  &  & $1$ &  & $1$ \\
 \hline
 $9_2$   &   & $1$ &    &  & $1$  &  & $1$&  & $1$ &  \\
 \hline
 $9_2'$   &  $1$ & &    &  &  & $1$ &  & $1$ &  & $1$ \\
 \hline
\end{tabular}
\end{tiny}
\smallskip

\caption{Restriction of the $9$-dimensional representations to $\mathcal{H}_{B_3,2}\simeq \mathcal{H}_{B_3,\beta,\alpha}$.}\label{resF49B32}
\end{table}

Let $\F_{\tilde{q}}=\F_p(\sqrt{\alpha},\sqrt{\beta})$ and $G_{9_1}=\rho_{9_1}(\mathcal{A}_{F_4})$. By Tables \ref{resF49B31} and \ref{resF49B32}, in cases $1$, $4$, $5$, $10$, $11$, $12$, $13$ and $16$, $\rho_{9_1}(\mathcal{A}_{B_3,1})$ and $\rho_{9_2}(\mathcal{A}_{B_3,2})$ are generated by transvections. It follows that $\rho_{9_1}(\mathcal{A}_{F_4})$ is generated by transvections. We get by Theorem \ref{transvections} that there exists $q'$ dividing $\tilde{q}$ such that $G_{9_1}$ is conjugate in $GL_9(\tilde{q})$ to $SL_9(q')$, $SU_9(q'^{\frac{1}{2}})$ or $SP_9(q')$. The symplectic case is excluded since we have $\rho_{9_1}\not\simeq \rho_{9_1}^\star$. In case $1$, we have $\tilde{q}=q$ and $\F_p(\alpha,\beta)=\F_p(\alpha+\alpha^{-1},\beta+\beta^{-1})$. By Theorem \ref{result1}, $G_{9_1}$ contains a natural $SL_3(q)$, therefore $q'=q$ and $G_{9_1}=SL_9(q)$. In cases $4$, $5$ and $10$, we have $\tilde{q}=q^2$ and $\F_q=\F_p(\alpha,\beta)=\F_p(\alpha+\alpha^{-1},\beta+\beta^{-1})$. We then have again by Theorem \ref{result1} that $G_{9_1}$ contains a natural $SL_3(q)$, therefore $q'\in \{q,q^2\}$. There exists a unique automorphism $\Phi$ of order $2$ of $\F_{q^2}$ and if we consider the restrictions of $\rho_{9_1}$ to $\mathcal{H}_{B_3,1}$, then $\Phi\circ \rho_{[2],[1]|\mathcal{A}_{F_4}} \simeq \rho_{[2],[1]|\mathcal{A}_{F_4}}$ and $\Phi\circ \rho_{[1],[2]|\mathcal{A}_{F_4}}\simeq \rho_{[1],[2]|\mathcal{A}_{F_4}}$. By Table \ref{resF49B31}, we have either $\Phi\circ  \rho_{9_1|\mathcal{A}_{F_4}}\simeq \rho_{9_1|\mathcal{A}_{F_4}}$ or $\Phi\circ \rho_{9_1|\mathcal{A}_{F_4}}\simeq \rho_{9_2|\mathcal{A}_{F_4}}$. Assume now by contradiction that $\Phi\circ \rho_{9_1|\mathcal{A}_{F_4}}\simeq \rho_{9_2|\mathcal{A}_{F_4}}$. There exists a character $\chi$ from $\mathcal{A}_{F_4}$ to $\F_{q^2}$ such that $\Phi\circ\rho_{9_1}\simeq \rho_{9_2}\otimes \chi$ by Lemma \ref{abel}. Since $A_{F_4}/\mathcal{A}_{F_4}\simeq <\overline{S_2},\overline{S_3}>$, it follows that there exists $y\in \F_{q^2}$ such that $\Phi(\rho_{9_1}(S_3))$ is conjugate in $GL_9(q)$ to $y\rho_{9_2}(S_3)$.  The eigenvalues of $\rho_{9_1}(S_3)$ are $\beta$ with multiplicity $3$ and $-1$ with multiplicity $6$, therefore the eigenvalues of $\Phi(\rho_{9_1}(S_3))$ are $\beta$ with multiplicity $6$ and $-1$ with multiplicity $3$. On the other hand, the eigenvalues of $y\rho_{9_2}(S_3)$ are $y\beta$ with multiplicity $3$ and $-y$ with multiplicity $6$. This implies that $y\beta=-1$ and $-y=\beta$. It follows that $y=-\beta=-\beta^{-1}$, therefore $\beta^2=1$ which is absurd. This proves that $\Phi\circ \rho_{9_1}\simeq \rho_{9_1}$. Hence, up to conjugation in $GL_9(q)$, $G_{9_1}$ is a subgroup of $SL_9(q)$ by Lemma \ref{Harinordoquy}. It follows that $G_{9_1}\simeq SL_9(q)$ in cases $4$, $5$ and $10$.

In cases $11$, $12$, $13$ and $16$ we have $\tilde{q}=q$ and $G_{9_1}$ contains a natural $SU_3(q^{\frac{1}{2}})$. By Lemma \ref{field}, we have $q'=q$. There exists a unique automorphism $\epsilon$ of order $2$ of $\F_q$ and for any representation $\varphi$ appearing in the decomposition of the restriction of $\rho_{9_1}$ from $\mathcal{H}_{F_4}$ to $\mathcal{H}_{B_3,1}$, we have $\epsilon \circ \varphi^{\star}\simeq \varphi$. It follows that $\epsilon \circ \rho_{9_1}^\star \simeq \rho_{9_1}$ or $\epsilon \circ \rho_{9_1}^\star \simeq \rho_{9_2}$. The eigenvalues of $\epsilon(\rho_{9_1}^\star(S_3))$ are $\beta^{-1}$ with multiplicity $3$ and $-1$ with multiplicity $6$, whereas the eigenvalues of $v\rho_{9_2}(S_3)$ are $v\beta$ with multiplicity $3$ and $-v$ with multiplicity $6$. We get in the same way as in cases $4$, $5$ and $10$ that $\epsilon\circ\rho_{9_1}^\star\simeq \rho_{9_1}$. By Lemma \ref{Ngwenya}, we get that $G_{9_2}$ is conjugate to a subgroup of $SU_9(q^{\frac{1}{2}})$, therefore $G_{9_2}\simeq SU_9(q^{\frac{1}{2}})$.

\medskip

Assume now $\F_p(\sqrt{\alpha},\sqrt{\beta})=\F_p(\alpha,\beta)=\F_q=\F_p(\alpha+\alpha^{-1},\beta)\neq \F_p(\alpha,\beta+\beta^{-1})=\F_p(\alpha+\alpha^{-1},\beta+\beta^{-1})$. We then have $\rho_{9_1}(\mathcal{A}_{B_3,1})\simeq SL_3(q)\times SU_2(r^{\frac{1}{2}})$ if $\F_r=\F_p(\alpha)\neq \F_p(\alpha+\alpha^{-1})$ and $\rho_{9_1}(\mathcal{A}_{B_3,1})\simeq SL_3(q)\times SL_2(r)$ if $\F_r=\F_p(\alpha)=\F_p(\alpha+\alpha^{-1})$. There exists a unique automorphism $\epsilon$ of order $2$ of $\F_q$ and if we consider the representations appearing in the restriction of $\rho_{9_1}$ to $\mathcal{H}_{B_3,1}$, we have by Proposition \ref{isomorphisme} that $\epsilon \circ \rho_{[2],[1]}\simeq \rho_{[1],[2]}$ and $\epsilon \circ \rho_{[1],[2]}\simeq \rho_{[2],[1]}$. It follows that $\epsilon \circ \rho_{9_1|\mathcal{A}_{F_4}}\simeq \rho_{9_1|\mathcal{A}_{F_4}}$ or $\epsilon \circ \rho_{9_1|\mathcal{A}_{F_4}}\simeq \rho_{9_2|\mathcal{A}_{F_4}}$. We have here $\epsilon(\alpha)=\alpha$ and $\epsilon(\beta)=\beta^{-1}$. Assume now by contradiction that $\epsilon \circ \rho_{9_1|\mathcal{A}_{F_4}}\simeq \rho_{9_1|\mathcal{A}_{F_4}}$. There exists $v\in \F_{q}$ such that $\epsilon(\rho_{9_1})(S_3)$ conjugate to $v\rho_{9_1}(S_3)$. The eigenvalues of $\epsilon(\rho_{9_1}(S_3))$ are $\beta^{-1}$ with multiplicity $6$ and $-1$ with multiplicity $3$ and the eigenvalues of $v\rho_{9_1}(S_1)$ are $v\beta$ with multiplicity $6$ and $-v$ with multiplicity $3$. This implies that $v=1=\beta^{-2}$, therefore $\beta^2=1$ which is absurd. This implies that $\epsilon \circ \rho_{9_1|\mathcal{A}_{F_4}}\simeq \rho_{9_2|\mathcal{A}_{F_4}}$. We will now use Theorem \ref{CGFS}. We have that $G_{9_1}$ contains a natural $SL_2(r)$ or a natural $SU_2(r^{\frac{1}{2}})$. This implies that $v_{G_{9_1}}(V)\leq 2=\max(2,\frac{\sqrt{9}}{2})$ and that we are not in the second case of Theorem \ref{CGFS}. It also implies by Lemma \ref{tens1} that $G_{9_1}$ is tensor-indecomposable. However, $\mathcal{A}_{F_4}$ is not normally generated by $\mathcal{A}_{B_3}$, therefore we cannot use the same arguments as before to show that $G_{9_1}$ is primitive in the non-monomial case. Nevertheless, we can show in the same way that if $G_{9_1}\subset (GL_3(q)\times GL_3(q)\times GL_3(q))\rtimes \mathfrak{S}_3$ then $\rho_{9_1}(\mathcal{A}_{B_3,1})\subset GL_3(q)\times GL_3(q)\times GL_3(q)$ and $\rho_{9_1}(\mathcal{A}_{B_3,2})\subset GL_3(q)\times GL_3(q)\times GL_3(q)$. Since $G_{9_1}$ is generated by $\rho_{9_1}(\mathcal{A}_{B_3,1})$ and $\rho_{9_1}(\mathcal{A}_{B_3,2})$, this contradicts the irreducibility, and $G_{9_1}$ is therefore primitive. It follows that $G_{9_1}$ is a classical group over $\F_{q'}$ for some $q'$ dividing $q$ in a natural representation. We have that $\rho_{9_1|\mathcal{A}_{F_4}}\not\simeq \rho_{9_1|\mathcal{A}_{F_4}}^\star$, therefore $G_{9_1}$ preserves no non-degenerate bilinear form. We have that $G_{9_1}$ contains $\rho_{9_1}(\mathcal{A}_{B_3,1})$ which is equal up to conjugation to $\left\lbrace\begin{pmatrix}
M & 0 & 0 & 0\\
0 & \epsilon(M)  & 0& 0\\
 0  & 0 & N& 0\\
 0 & 0 & 0  & 1
\end{pmatrix}, M\in SL_3(q), N\in SL_2(\tilde{r})\right\rbrace$, where $\tilde{r}=r$ if $\F_r=\F_p(\alpha)=\F_p(\alpha+\alpha^{-1})$ and $\tilde{r}=r^{\frac{1}{2}}$ otherwise. It follows that $G_{9_1}$ contains up to conjugation $\op{diag}(\alpha,\alpha^{-1},1,\alpha,\alpha^{-1},1,1,1,1)$ and $\op{diag}(\beta,\beta^{-1},1,\beta,\beta^{-1},1,1,1,1)$. Thus, the field generated by the traces of the elements of $G_{9_1}$ contains $5+2(\alpha+\alpha^{-1})$ and $5+2(\beta+\beta^{-1})$. Since $p\neq 2$, this field contains $\alpha+\alpha^{-1}$ and $\beta+\beta^{-1}$, therefore $q^{\frac{1}{2}}$ divides $q'$. This implies that $G_{9_1}$ is conjugate to $SU_9(q^{\frac{1}{4}})$, $SL_9(q^{\frac{1}{2}})$, $SU_9(q^{\frac{1}{2}})$ or $SL_9(q)$. We have $\epsilon \circ \rho_{9_1|\mathcal{A}_{F_4}}\simeq \rho_{9_2}$, therefore $\epsilon \circ \rho_{9_1|\mathcal{A}_{F_4}}^{\star}\simeq \rho_{9_2}^{\star}\simeq \rho_{9_2'}$. We have by Proposition \ref{resF4derivedsubgroup} that $\rho_{9_2|\mathcal{A}_{F_4}}\not\simeq \rho_{9_1|\mathcal{A}_{F_4}}$ and $\rho_{9_2'|\mathcal{A}_{F_4}}\not\simeq \rho_{9_1|\mathcal{A}_{F_4}}$. It follows that $G_{9_1}\simeq SL_9(q)$. This corresponds to cases $2$, $6$, $9$ and $15$.

\medskip

Assume now $\F_p(\sqrt{\alpha},\sqrt{\beta})=\F_p(\alpha,\beta)=\F_q=\F_p(\alpha,\beta+\beta^{-1})\neq \F_p(\alpha+\alpha^{-1},\beta)=\F_p(\alpha+\alpha^{-1},\beta+\beta^{-1})$. We then have $\rho_{9_1}(\mathcal{A}_{B_3,1})\simeq SL_3(q)\times SU_2(r^{\frac{1}{2}})$, where $\F_r=\F_p(\alpha)\neq \F_p(\alpha+\alpha^{-1})$. There exists a unique automorphism $\epsilon$ of order $2$ of $\F_q$ and if we consider the representations appearing in the restriction of $\rho_{9_1}$ to $\mathcal{H}_{B_3,1}$, $\epsilon \circ \rho_{[2,1]} \simeq \rho_{[1^2],[1]}$ and $\epsilon \circ \rho_{[1],[2]}\simeq \rho_{[1],[1^2]}$ by Proposition \ref{isomorphisme}. It follows that $\epsilon \circ \rho_{9_1|\mathcal{A}_{F_4}}\simeq \rho_{9_1'|\mathcal{A}_{F_4}}$ or $\epsilon \circ \rho_{9_1|\mathcal{A}_{F_4}}\simeq \rho_{9_2'|\mathcal{A}_{F_4}}$. We have here $\epsilon(\alpha)=\alpha^{-1}$ and $\epsilon(\beta)=\beta$. We have 

\begin{eqnarray*}
\tr(\rho_{9_1}(S_1S_2^{-1}S_3S_4^{-1})) & = & \frac{\alpha^2\beta^2-3\alpha^2\beta-2\alpha\beta^2+2\alpha^2+5\alpha\beta+\beta^2-3\alpha-2\beta+1}{\alpha\beta}\\
\epsilon(\tr(\rho_{9_1}(S_1S_2^{-1}S_3S_4^{-1}))) & = & \frac{\alpha^2-2\alpha^2\beta-2\alpha\beta^2+\alpha^2+5\alpha\beta+\beta^2-3\alpha-3\beta+2}{\alpha\beta}\\
\tr(\rho_{9_1'}(S_1S_2^{-1}S_3S_4^{-1})) & = & \frac{\alpha^2\beta^2-2\alpha^2\beta-3\alpha\beta^2+\alpha^2+5\alpha\beta+2\beta-2\alpha-3\beta+1}{\alpha\beta}
\end{eqnarray*}

It follows that

$$\epsilon(\tr(\rho_{9_1}(S_1S_2^{-1}S_3S_4^{-1})))-\tr(\rho_{9_1'}(S_1S_2^{-1}S_3S_4^{-1}))=\frac{(\alpha-1)(\beta^2-1)}{\alpha\beta}\neq 0.$$

This implies that $\epsilon \circ \rho_{9_1|\mathcal{A}_{F_4}}\not\simeq \rho_{9_1'|\mathcal{A}_{F_4}}$, therefore  $\epsilon \circ \rho_{9_1|\mathcal{A}_{F_4}}\simeq \rho_{9_2'|\mathcal{A}_{F_4}}$. We then get in the same ways as in cases $2$, $6$, $9$ and $15$ that we can apply Theorem \ref{CGFS} to get that $G_{9_1}$ is a classical group over $\F_{q'}$ in a natural representation and that $q^{\frac{1}{2}}$ divides $q'$. This implies that $G_{9_1}$ is conjugate to $SU_9(q^{\frac{1}{4}})$, $SL_9(q^{\frac{1}{2}})$, $SU_9(q^{\frac{1}{2}})$ or $SL_9(q)$. We have $\epsilon \circ \rho_{9_1|\mathcal{A}_{F_4}}\simeq \rho_{9_2'|\mathcal{A}_{F_4}}\not\simeq \rho_{9_1}$ and $\epsilon \circ \rho_{9_1|\mathcal{A}_{F_4}}^\star\simeq \rho_{9_2|\mathcal{A}_{F_4}}\not\simeq \rho_{9_1}$. It follows that $G_{9_1}$ is conjugate to $SL_9(q)$. This corresponds to cases $3$, $7$, $8$ and $14$.

\bigskip

By what was proven in cases $2$, $3$, $6$, $7$, $8$, $9$, $14$ and $15$ and Table \ref{resF49B32}, we have the second part of the proposition.
\end{proof}

\begin{prop2}
In cases $11$, $12$, $13$ and $16$, we have $\rho_{6_1}(\mathcal{A}_{F_4})\simeq \Omega_6^+(q^{\frac{1}{2}})$ and $\rho_{6_2}(\mathcal{A}_{F_4})\simeq \Omega_6^+(q^{\frac{1}{2}})$.

In all the remaining cases, we have $\rho_{6_1}(\mathcal{A}_{F_4})\simeq \Omega_6^+(q)$ and $\rho_{6_2}(\mathcal{A}_{F_4})\simeq \Omega_6^+(q)$.
\end{prop2}

 \begin{proof}
 \begin{table}
\begin{tiny}
\centering
\begin{tabular}{ |p{0.4cm}|p{1.1cm}|p{1.1cm}|p{1.1cm}|p{1.1cm}|p{1.1cm}|p{1.1cm}|p{1.1cm}|p{1.1cm}|p{1.1cm}|p{1.1cm}|  }
 \hline
 & $(\emptyset,[1^3])$ & $([3],\emptyset)$ & $([1^3],\emptyset)$ & $(\emptyset,[3])$ & $([2,1],\emptyset)$ & $(\emptyset,[2,1])$ & $([2],[1])$ & $([1],[1^2])$ & $([1],[2])$ &  $([1^2],[1])$  \\
 \hline

 $6_1$  &   &  &    &  &  &  & & & $1$ & $1$ \\
 \hline
 $6_2$    &   &  &  &  &  &  & $1$ & $1$ &  &  \\
 \hline
\end{tabular}
\end{tiny}
\smallskip

\caption{Restriction of the $6$-dimensional representations to $\mathcal{H}_{B_3,1}\simeq \mathcal{H}_{B_3,\alpha,\beta}$.}\label{resF46B31}
\end{table}

\begin{table}
\begin{tiny}
\centering
\begin{tabular}{ |p{0.4cm}|p{1.1cm}|p{1.1cm}|p{1.1cm}|p{1.1cm}|p{1.1cm}|p{1.1cm}|p{1.1cm}|p{1.1cm}|p{1.1cm}|p{1.1cm}|  }
 \hline
 & $(\emptyset,[1^3])$ & $([3],\emptyset)$ & $([1^3],\emptyset)$ & $(\emptyset,[3])$ & $([2,1],\emptyset)$ & $(\emptyset,[2,1])$ & $([2],[1])$ & $([1],[1^2])$ & $([1],[2])$ &  $([1^2],[1])$  \\
 \hline

 $6_1$  &   &  &    &  &  &  & & & $1$ & $1$ \\
 \hline
 $6_2$     &   &  &  &  &  &  & $1$ & $1$ &  &  \\ 
 \hline
\end{tabular}
\end{tiny}
 
 \smallskip
 \caption{Restriction of the $6$-dimensional representations to $\mathcal{H}_{B_3,1}\simeq \mathcal{H}_{B_3,\alpha,\beta}$.}\label{resF46B32}
 \end{table}
 
Let $\F_{\tilde{q}}=\F_p(\sqrt{\alpha},\sqrt{\beta})$, $G_{6_1}=\rho_{6_1}(\mathcal{A}_{F_4})$ and $G_{6_2}=\rho_{6_2}(\mathcal{A}_{F_4})$.

Let $P_1$ be the anti-diagonal matrix with coefficients $(1,-1,1,1,-1,1)$. We write the matrices of $\rho_{6_1}$ with respect to the basis $(e_{x_1},e_{x_2},e_{x_3},e_{x_4},e_{x_5},e_{x_6})$, where $I(x_1)=\{s_1,s_3,s_4\}$, $I(x_2)=\{s_2,s_4\}$, $I(x_3)=\{s_3\}$ and $I(x_j)=\{s_1,s_2,s_3,s_4\}\setminus I(x_{6-j})$ for $j\in \{4,5,6\}$. We then have for all $g\in G_{6_1}$, $P_1gP_1^{-1}=^t\!g^{-1}$. This proves that up to conjugation in $GL_6(\tilde{q})$, $G_{6_1}\leq \Omega_6^+(\tilde{q})$. By Proposition \ref{bilinwgraphs}, we have that up to conjugation, $G_{6_2}\leq \Omega_6^+(\tilde{q})$.

 By Tables \ref{resF46B31} and \ref{resF46B32} and Theorems \ref{result1} to \ref{result6}, we have that ($\rho_{6_1}(\mathcal{A}_{B_3,1})\simeq SL_3(q)$ or $\rho_{6_1}(\mathcal{A}_{B_3,1})\simeq SU_3(q^{\frac{1}{2}})$), ($\rho_{6_1}(\mathcal{A}_{B_3,2})\simeq SL_3(q)$ or $\rho_{6_1}(\mathcal{A}_{B_3,1})\simeq SU_3(q^{\frac{1}{2}})$) and ($\rho_{6_1}(\mathcal{A}_{B_2})\simeq SL_2(q)$ or $\rho_{6_1}(\mathcal{A}_{B_2})\simeq SL_2(q^{\frac{1}{2}})$. Note that those isomorphisms are given by twisted diagonal embeddings, therefore a transvection in $SL_3(q)$ is mapped to a long root element of $\Omega_6^+(q)$ by the isomorphism $\rho_{6_1}(\mathcal{A}_{B_3,1})\simeq SL_3(q)$. This proves that $\rho_{6_1}(\mathcal{A}_{B_3,1})$ and $\rho_{6_1}(\mathcal{A}_{B_3,2})$ are generated by long root elements. If $t$ is a long root element of $\rho_{6_1}(\mathcal{A}_{B_2})$ then its normal closure in $\rho_{6_1}(\mathcal{A}_{B_3,1})$ is equal to $\rho_{6_1}(\mathcal{A}_{B_3,1})$ and its normal closure in $\rho_{6_1}(\mathcal{A}_{B_3,2})$ is equal to $\rho_{6_1}(\mathcal{A}_{B_3,2})$. It follows that the normal closure of $t$ in $G_{6_1}$ is equal to $G_{6_1}$ since by Lemma \ref{genB3F4}, $\mathcal{A}_{B_3,1}$ and $\mathcal{A}_{B_3,2}$ generate $\mathcal{A}_{F_4}$. This proves that $G_{6_1}$ is an irreducible subgroup of $\Omega_6^+(\tilde{q})$ generated by a conjugacy class of long root elements. Since $O_p(G_{6_1})$ is normal in $G_{6_1}$ and $V=\F_q^6$ is an irreducible $\F_qG_{6_1}$-module, we apply Clifford's Theorem \cite[Theorem 11.1]{C-R} and get that $Res^{G_{6_1}}_{O_p(G_{6_1})}(V)$ is semisimple. Since $O_p(G_{6_1})$ is a $p$-group, the unique irreducible $\F_qO_p(G_{6_1})$-module is the trivial module, therefore $O_p(G_{6_1})$ acts trivially on $V$. It follows that $O_p(G_{6_1})$ is trivial. By Theorem \ref{theoKantor}, $G_{6_1}$ is isomorphic to one of the following groups for some $q'$ dividing $\tilde{q}$
 \begin{enumerate}
 \item $\Omega_6^+(q')$ in a natural representation,
 \item $\Omega_6^-(q'^{\frac{1}{2}})$ as a subgroup of $\Omega_6^+(q')$, where $\Omega_6^+(q')$ is in a natural representation,
 \item $SU_3(q')$ as a subgroup of $\Omega_6^+(q')$, where $\Omega_6^+(q')$ is in a natural representation.
 \end{enumerate}
 
 Assume now that we are in cases $4$, $5$ or $10$. We have $\tilde{q}=q^2$ and $\F_q=\F_p(\alpha,\beta)=\F_p(\alpha+\alpha^{-1},\beta+\beta^{-1})$. There exists a unique automorphism $\Phi$ of order $2$ of $\F_q$. We have $\Phi\circ \rho_{6_1|\mathcal{A}_{F_4}}\simeq \rho_{6_1|\mathcal{A}_{F_4}}$ or $\Phi\circ \rho_{6_1|\mathcal{A}_{F_4}}\simeq \rho_{6_2|\mathcal{A}_{F_4}}$. By Table \ref{resF46B31}, we cannot have $\Phi\circ \rho_{6_1|\mathcal{A}_{F_4}}\simeq \rho_{6_2|\mathcal{A}_{F_4}}$. It follows that $\Phi\circ \rho_{6_1|\mathcal{A}_{F_4}}\simeq \rho_{6_1|\mathcal{A}_{F_4}}$. Hence, by Lemma \ref{Harinordoquy}, $G_{6_1}$ is conjugate in $GL_6(q^2)$ to a subgroup of $\Omega_6^+(q)$. Since it contains up to conjuugation a twisted diagonal $SL_3(q)$, it contains up to conjugation $\diag(\alpha,\alpha^{-1},1,\alpha^{-1},\alpha,1)$ and  $\diag(\beta,\beta^{-1},1,\beta^{-1},\beta,1)$. This implies the field generated by the traces of its elements contains $2(\alpha+\alpha^{-1})+2$ and $2(\beta+\beta^{-1})+2$. It follows that $q'=q$.
 
  If $G_{6_1}\simeq SU_3(q)$ then we would have that $SL_3(q)$ is isomorphic to a subgroup of $SU_3(q)$, therefore $q^3(q^2-1)(q^3-1)$ divides $q^3(q^2-1)(q^3+1)$. This would imply that $q^3-1$ divides $q^3+1$, therefore $q^3-1$ divides $q^3+1-(q^3-1)=2$. This is absurd, therefore $G_{6_1}\not\simeq SU_3(q)$.
  
  Assume now by contradiction that $G_{6_1}\simeq \Omega_6^-(q^{\frac{1}{2}})$.  We have $\vert \Omega_6^-(q^{\frac{1}{2}})\vert =\frac{1}{2} q^3
(q^{\frac{3}{2}}+1)(q-1)(q^2-1)$ and $\vert SL_3(q)\vert =q^3(q^2-1)(q^3-1)$. Since $\op{Gcd}(q^3-1,q^3)=1$, $q^3-1$ divides $(q^{\frac{3}{2}}+1)(q-1)=q^{\frac{5}{2}}-q^{\frac{3}{2}}+q-1<q^{\frac{5}{2}}-1<q^3-1$. This is absurd, therefore $G_{6_1}\not\simeq \Omega_6^-(q^{\frac{1}{2}})$.

It follows that in cases $4$, $5$ and $10$
, we have $G_{6_1}\simeq \Omega_6^+(q)$.

\medskip

We can now assume $\tilde{q}=q$, therefore $q'$ divides $q$. By the same reasoning as above using Lemma \ref{field}, we have that $\F_{q'}$ contains $F_p(\alpha+\alpha^{-1},\beta+\beta^{-1})$ in all the remaining cases. In case $1$, we can apply the same reasoning as above to get that $G_{6_1}\simeq \Omega_6^+(q)$. 

\smallskip

In cases $11$, $12$, $13$ and $16$, there exists a unique automorphism $\epsilon$ of order $2$ of $\F_q$ and for any representation $\varphi$ of $\mathcal{H}_{B_3,1}$, we have $\epsilon \circ \varphi^{\star} \simeq \varphi$. It then follows by Table \ref{resF46B31} that $\epsilon \circ \rho_{6_1|\mathcal{A}_{F_4}}^{\star} \simeq \rho_{6_1|\mathcal{A}_{F_4}}$. Then, by \cite{BMM} (Proposition 4.1.), $G_{6_1}$ is conjugate in $GL_6(q)$ to a subgroup of $\Omega_6^+(q^{\frac{1}{2}})$. It follows that $q'$ divides $q^{\frac{1}{2}}$, therefore $q'=q^{\frac{1}{2}}$ since $\F_{q'}$ contains $\F_p(\alpha+\alpha^{-1},\beta+\beta^{-1})=\F_{q^{\frac{1}{2}}}$. This proves that $G_{6_1}$ is isomorphic to one of the following groups

 \begin{enumerate}
 \item $\Omega_6^+(q^{\frac{1}{2}})$ in a natural representation,
 \item $\Omega_6^-(q^{\frac{1}{4}})$ as a subgroup of $\Omega_6^+(q^{\frac{1}{2}})$, where $\Omega_6^+(q^{\frac{1}{2}})$ is in a natural representation,
 \item $SU_3(q^{\frac{1}{2}})$ as a subgroup of $\Omega_6^+(q^{\frac{1}{2}})$, where $\Omega_6^+(q^{\frac{1}{2}})$ in a natural representation.
 \end{enumerate}

Assume now by contradiction that $G_{6_1}\simeq SU_3(q^{\frac{1}{2}})$. We know that $\rho_{6_1}(\mathcal{A}_{B_3,1})\simeq  SU_3(q^{\frac{1}{2}})$. We only need to show that $\rho_{6_1}(\mathcal{A}_{B_3,1})\neq G_{6_1}$. 

Let $X=\begin{pmatrix}
0 & 1 & 0 & 0 & 0 & \frac{\alpha^2+\beta}{\alpha\sqrt{\beta}}\\
\frac{\sqrt{\alpha}(\beta+1)}{\alpha\beta+1} & -\frac{\alpha+\beta}{\alpha\beta+1}& 0 & 0 & \frac{(\beta+1)(\alpha^2+\beta)}{\sqrt{\alpha}\sqrt{\beta}(\alpha\beta+1)} & -\frac{(\alpha+\beta)(\alpha^2+\beta)}{\alpha\sqrt{\beta}(\alpha\beta+1)}\\
-\frac{\alpha^2+\beta}{\sqrt{\alpha}(\alpha\beta+1)}& 0 & \frac{(\alpha^2+\beta)(\alpha+\beta)}{\alpha\sqrt{\beta}(\alpha\beta+1)}& 0 &  -\frac{(\alpha^4+\beta)^2}{\alpha\sqrt{\alpha}\sqrt{\beta}(\alpha\beta+1)}& 0\\
1 & -\frac{\alpha+\beta}{\sqrt{\alpha}(\beta+1)} & 0 & \frac{\alpha^2+\beta}{\alpha(\beta+1)}& 0 & 0\\
0 & \frac{\alpha\beta+1}{\sqrt{\alpha}(\beta+1)} & 0 & -\frac{(\alpha\beta+1)(\alpha^2+\beta)}{\alpha(\alpha+\beta)(\beta+1)}& 0 & 0\\
0 & 0 & 0 & \frac{\alpha\beta+1}{\alpha+\beta} & 0 & 0
\end{pmatrix}$. We have with respect to the previous basis that for all $g\in  \mathcal{A}_{B_3,1}$, 
$$X\rho_{6_1}(g)X^{-1}=\begin{pmatrix}\rho_{[1],[2]}(g) & 0\\
0 & \rho_{[1^2],[1]}(g)\end{pmatrix}.$$
Note that here, the matrix $\rho_{6_1}(q)$ with respect to the previous basis is again denoted by the same symbol in order to simplify the notations.

 We have $(X\rho_{6_1}(S_3S_4^{-1})X^{-1})_{2,6}=-\frac{\alpha^3+\beta^3}{\sqrt{\alpha}(\alpha\beta+1)^2\beta}\neq 0$. This proves that $\rho_{6_1}(S_3S_4^{-1})\notin \rho_{6_1}(\mathcal{A}_{B_3,1})$, therefore $G_{6_1}\not\simeq SU_3(q^{\frac{1}{2}})$.

Assume by contradiction that $G_{6_1}\simeq \Omega_6^-(q^{\frac{1}{4}})$. Then $\vert \Omega_6^-(q^{\frac{1}{4}})\vert =\frac{1}{2}q^{\frac{3}{2}}(q^{\frac{3}{4}}+1)(q^{\frac{1}{2}}-1)(q-1)$, $\vert SU_3(q^{\frac{1}{2}})\vert =q^{\frac{3}{2}}(q-1)(q^{\frac{3}{2}}+1)\vert $, and $q^{\frac{3}{2}}+1$ divides $(q^{\frac{3}{4}}+1)(q-1)\frac{1}{2}$. Therefore, $q^{\frac{3}{2}}+1$ divides $(q^{\frac{3}{4}}+1)(q-1)\frac{1}{2}-(q^{\frac{3}{4}}-1)(q-1)\frac{1}{2}=q-1<q^{\frac{3}{2}}+1$. This is absurd. Thus, we have $G_{6_1}\simeq \Omega_6^+(q^{\frac{1}{2}})$.

\medskip

Assume now $\F_p(\sqrt{\alpha},\sqrt{\beta})=\F_p(\alpha,\beta)=\F_q=\F_p(\alpha+\alpha^{-1},\beta)\neq \F_p(\alpha,\beta+\beta^{-1})=\F_p(\alpha+\alpha^{-1},\beta+\beta^{-1})$. There exists then a unique automorphism $\epsilon$ of order $2$ of $\F_q$.
We have by Proposition \ref{isomorphisme} that $\epsilon\circ \rho_{[1],[2]|\mathcal{A}_{B_3,1}}\simeq \rho_{[2],[1]|\mathcal{A}_{B_3,1}}$. This proves that $\epsilon \circ \rho_{6_1|\mathcal{A}_{F_4}}\simeq \rho_{6_2|\mathcal{A}_{F_4}}$ by Table \ref{resF46B31}. We know that $q'\in \{q^{\frac{1}{2}},q\}$. We have by Proposition \ref{resF4derivedsubgroup} that $\epsilon \circ \rho_{6_1|\mathcal{A}_{F_4}}\not\simeq \rho_{6_1|\mathcal{A}_{F_4}}$, therefore we cannot have $q'=q^{\frac{1}{2}}$. It follows that $q'=q$. The proof then uses the same arguments as in case $1$, therefore we get that $G_{6_1}\simeq \Omega_6^+(q)$ in cases $2$, $6$, $9$ and $15$.

\medskip

Finally, assume $\F_p(\sqrt{\alpha},\sqrt{\beta})=\F_p(\alpha,\beta)=\F_q=\F_p(\alpha,\beta+\beta^{-1})\neq \F_p(\alpha+\alpha^{-1},\beta)=\F_p(\alpha+\alpha^{-1},\beta+\beta^{-1})$. There exists then a unique automorphism $\epsilon$ of order $2$ of $\F_q$. We have by Proposition \ref{isomorphisme} that $\epsilon\circ \rho_{[1],[2]|\mathcal{A}_{B_3,1}}\simeq \rho_{[1],[1^2]|\mathcal{A}_{B_3,1}}$. This proves that $\epsilon \circ \rho_{6_1|\mathcal{A}_{F_4}}\simeq \rho_{6_2|\mathcal{A}_{F_4}}$ by Table \ref{resF46B31}. We know that $q'\in \{q^{\frac{1}{2}},q\}$. We have by Proposition \ref{resF4derivedsubgroup} that $\epsilon \circ \rho_{6_1|\mathcal{A}_{F_4}}\not\simeq \rho_{6_1|\mathcal{A}_{F_4}}$, therefore we cannot have $q'=q^{\frac{1}{2}}$. It follows that $q'=q$. The proof then uses the same arguments as in case $1$, therefore we get that $G_{6_1}\simeq \Omega_6^+(q)$ in cases $3$, $7$, $8$ and $14$.

\bigskip

By Tables \ref{resF46B31} and \ref{resF46B32}, we only need to check that $\rho_{6_2}(\mathcal{A}_{F_4})\neq \rho_{6_2}(\mathcal{A}_{B_3,1})$ to get the same results for the representation $\rho_{6_2}$. Let

\begin{small} $$Y=\begin{pmatrix}
 0 & 0 & 0 & 0 & \alpha+\beta & \sqrt{\alpha}(\beta+1)\\
 0 & 0 & \sqrt{\alpha}(\beta+1) & 0 & -\alpha\beta-1 & 0\\
 0 & 0 & -\frac{\alpha^2\beta+1}{\sqrt{\alpha}}& 0 & 0 & 0\\
 1 & 0 & -\frac{\sqrt{\alpha}\sqrt{\beta}(\alpha+1)}{\alpha^2\beta+1} & 0 & \frac{\sqrt{\beta}\alpha}{\alpha^2\beta+1}& 0\\
 -\frac{\alpha+\beta}{\alpha\beta+1} & \frac{\sqrt{\alpha}(\beta+1)}{\alpha\beta+1} & -\frac{\sqrt{\alpha}\sqrt{\beta}}{\alpha\beta+1}& 0 & -\frac{\sqrt{\beta}\alpha(\alpha+\beta)}{(\alpha^2\beta+1)(\alpha\beta+1)}& -\frac{\alpha\sqrt{\alpha}\sqrt{\beta}(\beta+1)}{(\alpha^2\beta+1)(\alpha\beta+1)}\\
 \frac{\alpha(\beta+1)(\alpha+\beta)}{(\alpha^2\beta+1)(\alpha\beta+1)} & -\frac{\sqrt{\alpha}\alpha(\beta+1)^2}{(\alpha\beta+1)(\alpha^2\beta+1)} & \frac{\sqrt{\alpha}\alpha\sqrt{\beta}(\beta+1)}{(\alpha\beta+1)(\alpha^2\beta+1)} & \frac{\alpha(\beta+1)}{\alpha^2\beta+1} & -\frac{\alpha\sqrt{\beta}(\alpha+1)(\beta+1)}{(\alpha^2\beta+1)(\alpha\beta+1)} & -\frac{\sqrt{\alpha}\sqrt{\beta}\alpha(\beta+1)}{(\alpha^2\beta+1)(\alpha\beta+1)}
 \end{pmatrix}.$$\end{small}
 
 We consider the matrices corresponding to the representation $\rho_{6_2}$ with respect to the basis $(e_{x_1},e_{x_2},e_{x_3},e_{x_4},e_{x_5},e_{x_6})$ with $I(x_1)=\{s_1,s_2\}$, $I(x_2)=\{s_1,s_3\}$, $I(x_3)=\{s_1,s_4\}$, $I(x_4)=\{s_2,s_3\}$, $I(x_5)=\{s_2,s_4\}$, $I(x_6)=\{s_3,s_4\}$. We then have that for all $g\in \mathcal{H}_{B_3,1}$, 
$$Y\rho_{6_2}(g)Y^{-1}=\begin{pmatrix}
\rho_{[1],[2]}(g) & 0\\
0 & \rho_{[2],[1^2]}(g)
\end{pmatrix}.$$
 We have that $(Y\rho_{6_2}(S_3S_4^{-1})Y^{-1})_{1,5}=-\frac{\alpha+\beta}{\sqrt{\beta}}\neq 0$, therefore $G_{6_2}\neq \rho_{6_2}(\mathcal{A}_{B_3,1})$ and the proof is concluded. \end{proof}
 
 \begin{prop2}
 In cases $1$, $4$, $5$ and $10$, we have $\rho_{12}(\mathcal{A}_{F_4})\simeq \Omega_{12}^+(q)$ and $\rho_{16}(\mathcal{A}_{F_4})\simeq \Omega_{16}^+(q)$.
 In all the remaining cases, we have $\rho_{12}(\mathcal{A}_{F_4})\simeq \Omega_{12}^+(q^{\frac{1}{2}})$ and $\rho_{16}(\mathcal{A}_{F_4})\simeq \Omega_{16}^+(q)$.
 \end{prop2}
 
 \begin{proof}
 
\begin{table}
\begin{tiny}
\centering
\begin{tabular}{ |p{0.4cm}|p{1.1cm}|p{1.1cm}|p{1.1cm}|p{1.1cm}|p{1.1cm}|p{1.1cm}|p{1.1cm}|p{1.1cm}|p{1.1cm}|p{1.1cm}|  }
 \hline
 & $(\emptyset,[1^3])$ & $([3],\emptyset)$ & $([1^3],\emptyset)$ & $(\emptyset,[3])$ & $([2,1],\emptyset)$ & $(\emptyset,[2,1])$ & $([2],[1])$ & $([1],[1^2])$ & $([1],[2])$ &  $([1^2],[1])$  \\
 \hline

 $12$  &   &  &    &  &  &  & $1$ & $1$ & $1$ & $1$ \\
 \hline
 $16$     &   &  &  &  & $1$ & $1$ & $1$ & $1$ & $1$ & $1$ \\ 
 \hline
\end{tabular}
\end{tiny} 

\smallskip
 \caption{Restriction of the high dimensional representations to $\mathcal{H}_{B_3,1}\simeq \mathcal{H}_{B_3,\alpha,\beta}$}\label{resF4B31}
 \end{table}
 
 \begin{table}
\begin{tiny}
\centering
\begin{tabular}{ |p{0.4cm}|p{1.1cm}|p{1.1cm}|p{1.1cm}|p{1.1cm}|p{1.1cm}|p{1.1cm}|p{1.1cm}|p{1.1cm}|p{1.1cm}|p{1.1cm}|  }
 \hline
 & $(\emptyset,[1^3])$ & $([3],\emptyset)$ & $([1^3],\emptyset)$ & $(\emptyset,[3])$ & $([2,1],\emptyset)$ & $(\emptyset,[2,1])$ & $([2],[1])$ & $([1],[1^2])$ & $([1],[2])$ &  $([1^2],[1])$  \\
 \hline

 $12$  &   &  &    &  &  &  & $1$ & $1$ & $1$ & $1$ \\
 \hline
 $16$     &   &  &  &  & $1$ & $1$ & $1$ & $1$ & $1$ & $1$ \\ 
 \hline
\end{tabular}
\end{tiny} 

\smallskip
 \caption{Restriction of the high dimensional representations to $\mathcal{H}_{B_3,2}\simeq \mathcal{H}_{B_3,\beta,\alpha}$}\label{resF4B32}
 \end{table}
 
Let $\F_{\tilde{q}}=\F_p(\sqrt{\alpha},\sqrt{\beta})$, $G_{12}=\rho_{12}(\mathcal{A}_{F_4})$ and $G_{16}=\rho_{16}(\mathcal{A}_{F_4})$. In cases $1$, $4$, $5$, $10$, $11$, $12$, $13$ and $16$, $G_{12}$ contains a twisted diagonal $SL_3(q)$ or a twisted diagonal $SU_3(q^{\frac{1}{2}})$ by Table \ref{resF4B31} and Theorems \ref{result1} and \ref{result3}. Note that in cases $4$, $5$ and $10$ we have that $\Phi\circ\rho_{12|\mathcal{A}_{F_4}}\simeq \rho_{12|\mathcal{A}_{F_4}}$, where $\Phi$ is the unique automorphism of order $2$ of $\F_{\tilde{q}}$. It follows  by Lemma \ref{Ngwenya} that in cases $4$, $5$ and $10$, $G_{12}$ is conjugate in $GL_{12}(\tilde{q})$ to a subgroup of $GL_{12}(q)$. We then have that $v_{G_{12}}(v)\leq 2=\max(2,\frac{\sqrt{12}}{2})$, $G_{12}$ is primitive and we are not in case $(2)$ of Theorem \ref{CGFS}. In order to apply the theorem we have to prove that $G_{12}$ is tensor-indecomposable but we cannot use Lemma \ref{tens2} because $12\leq 16$. The arguments in the proof of Lemma \ref{tens2} \cite{BMM} only require $d\geq 16$ in order to have $a\geq 3$ and $b\geq 4$  if $a,b\geq 3$. This is still true for $d\leq 16$ unless $d=9$, therefore we can still apply Lemma \ref{tens2}. We still have to prove that we do not have $G_{12}\leq GL_2(q)\otimes GL_6(q)$. Assume by contradiction that $G_{12}\leq GL_2(q)\otimes GL_6(q)$. We then have a morphim from $G_{12}$ to $GL_2(q)$. Consider the restriction of this morphism to $\rho_{12}(\mathcal{A}_{B_3,1})$. Table \ref{resF4B31} and Propositions \ref{platypus} and \ref{platypus456} give us a morphism from $SL_3(q)\times SL_3(q)$ or $SU_3(q^{\frac{1}{2}})\times SU_3(q^{\frac{1}{2}})$ to $GL_2(q)$. We can consider the restriction to each factor and we get a morphism from $SL_3(q)$ to $GL_2(q)$ or a morphism from $SU_3(q^{\frac{1}{2}})$ to $GL_2(q)$. If this morphism is non-trivial then we get an isomorphism from $PSL_3(q)$ to a subgroup of $GL_2(q)$ or an isomorphism form $PSU_3(q^{\frac{1}{2}})$ to a subgroup of $GL_2(q)$. This would imply by considering the orders of those groups that $\frac{1}{(3,q-1)}q^3(q^2-1)(q^3-1)$ or $\frac{1}{(3,q-1)}q^3(q^2-1)(q^3+1)$ divides $q^2(q^2-1)(q-1)$ which is absurd since $q^3$ cannot divide $q^2$. This proves that the restriction of this morphism to each factor is trivial, therefore the restriction of this morphism to $\rho_{12}(\mathcal{A}_{B_3,1})$ is trivial. The restriction to $\rho_{12}(\mathcal{A}_{B_3,2})$ is also trivial by Table \ref{resF4B32} and the same arguments as above. This would imply by Lemma \ref{genB3F4} that this morphism is trivial, which contradicts the irreducibility of $G_{12}$. Therefore $G_{12}$ is tensor-indecomposable. So, we can apply Theorem \ref{CGFS} and we get that $G_{12}$ is a classical group in a natural representation. Consider now this representation with respect to the basis $(e_{x_i})_{i\in [\![1,12]\!]}$, where $I(x_1)=\{s_2\}$, $I(x_2)=\{s_3\}$, $I(x_3)=\{s_1,s_2\}$, $I(x_4)=\{s_1,s_3\}$, $I(x_5)=\{s_1,s_3\}$, $I(x_6)=\{s_1,s_4\}$ and for $j\in [\![7,12]\!]$, $I(x_j)=\{s_1,s_2,s_3,s_4\}\setminus\{I(x_{13-j})\}$. Let $P$ be the anti-diagonal matrix with coefficients $(3,3,-1,1,-3,-1,-1,-3,1,-1,3,3)$. We then have for all $i\in \{1,2\}$ that $P\rho_{12}(S_i)P^{-1}=-\alpha^t\!\rho_{12}(S_i)^{-1}$ and for all $j\in \{3,4\}$, $P\rho_{12}(S_j)P^{-1}=-\beta ^t\!\rho_{12}(S_j)$. This proves that $G_{12}$ is conjugate in $GL_{12}(\tilde{q})$ to a subgroup of $\Omega_{12}^+(q)$. Since it contains a twisted diagonal $SL_3(q)$ or a twisted diagonal $SU_3(q^{\frac{1}{2}})$, the field generated by its traces contains $\F_p(\alpha+\alpha^{-1},\beta+\beta^{-1})$. In case $1$, $4$, $5$ and $10$ this proves that $G_{12}$ is a classical group over $\F_q$ and since $G_{12}$ is conjugate to a subgroup of $\Omega_{12}^+(q)$, we get $G_{12}\simeq \Omega_{12}^+(q)$. 

Assume now that we are in case $11$, $12$, $13$ or $16$. Then there exists an automorphism $\epsilon$ of order $2$ of $\F_q$, and $\epsilon \circ \rho_{12|\mathcal{A}_{F_4}}\simeq \rho_{12|\mathcal{A}_{F_4}}$. It follows by Lemma \ref{Ngwenya} that $G_{12}$ is conjugate to a subgroup of $\Omega_{12}^+(q^{\frac{1}{2}})$. Hence, $G_{12}$ is a classical group in a natural representation over $\F_{q^{\frac{1}{2}}}$, therefore $G_{12}\simeq \Omega_{12}^+(q^{\frac{1}{2}})$.

\smallskip

Assume now $\F_p(\sqrt{\alpha},\sqrt{\beta})=\F_p(\alpha,\beta)=\F_q=\F_p(\alpha+\alpha^{-1},\beta)\neq \F_p(\alpha,\beta+\beta^{-1})=\F_p(\alpha+\alpha^{-1},\beta+\beta^{-1})$. There exists then a unique automorphism $\epsilon$ of order $2$ of $\F_q$. We have $\epsilon \circ \rho_{12|\mathcal{A}_{F_4}}\simeq \rho_{12|\mathcal{A}_{F_4}}$. This implies by Lemma \ref{Ngwenya} that $G_{12}$ is conjugate to a subgroup of $\Omega_{12}^+(q^{\frac{1}{2}})$. We know by Theorem \ref{result5} and Tables \ref{resF4B31} and \ref{resF4B32} that $\rho_{12}(\mathcal{A}_{B_3,1})\simeq SL_3(q)$ and $\rho_{12}(\mathcal{A}_{B_3,2})\simeq SL_3(q)$. We have that $G_{12}$ is primitive, irreducible and that the field generated by its traces contains $\F_{q^{\frac{1}{2}}}$. Let us show that $G_{12}$ is also tensor-indecomposable. If it was tensor-decomposable then we would have $G_{12}\leq GL_2(q^{\frac{1}{2}})\otimes GL_6(q^{\frac{1}{2}})$ or $G_{12}\leq GL_3(q^{\frac{1}{2}})\otimes GL_4(q^{\frac{1}{2}})$. By the same argument as in the above cases, we have that $G_{12}\leq GL_2(q^{\frac{1}{2}})\otimes GL_6(q^{\frac{1}{2}})$ is impossible.

 Assume now by contradiction that $G_{12}\leq GL_3(q^{\frac{1}{2}})\otimes GL_4(q^{\frac{1}{2}})$. We would then have a morphism from $G_{12}$ to $GL_3(q^{\frac{1}{2}})$. The restriction of this morphism from $G_{12}$ to $\rho_{12}(\mathcal{A}_{B_3,1})$ or to $\rho_{12}(\mathcal{A}_{B_3,2})$ gives us a morphism from $SL_3(q)$ to $GL_3(q^{\frac{1}{2}})$. If this morphism is non-trivial, we get an isomorphism from $PSL_3(q)$ to a subgroup of $GL_3(q^{\frac{1}{2}})$. This is absurd since $\vert PSL_3(q)\vert =\frac{1}{(3,q-1)}q^3(q^2-1)(q^3-1)$, $\vert GL_3(q^{\frac{1}{2}})\vert =q^{\frac{3}{2}}(q^{\frac{1}{2}}-1)(q-1)(q^{\frac{3}{2}}-1)$ and $q^3$ does not divide $q^{\frac{3}{2}}$. This proves that the restriction of this morphism to $\rho_{12}(\mathcal{A}_{B_3,1})$ and $\rho_{12}(\mathcal{A}_{B_3,2})$ is trivial. It follows by Lemma \ref{genB3F4} that this morphism is trivial and it contradicts the irreducibility of $G_{12}$. This proves that $G_{12}$ is tensor-indecomposable. We will now explicitely find a matrix $g\in G_{12}$ such that $\dim(g-I_{12})=2$. We consider our matrices in the same basis as before. 
 
 Let $X=(X_1 X_2 X_3)$ and $Y=(Y_1 Y_2 Y_3 Y_4)$, where $X_1$, $X_2$, $X_3$ and $Y_1$, $Y_2$, $Y_3$ and $Y_4$ are given in section \ref{sectionmatricesF412} of the Appendix. 

We then have that $X\rho_{12|\mathcal{A}_{B_3,1}}X^{-1}=\begin{pmatrix}
\rho_{[2],[1]} & 0 & 0 & 0\\
0 & \rho_{[2],[1]}^{\star} & 0 & 0\\
0 & 0 & \epsilon \circ \rho_{[2],[1]}& 0\\
0 & 0 & 0 & \epsilon \circ \rho_{[2],[1]}^\star
\end{pmatrix}$ and 
$$Y^{-1}\rho_{12|\mathcal{A}_{B_3,2}}Y=\begin{pmatrix}
\rho_{[2],[1]} & 0 & 0 & 0\\
0 & \rho_{[2],[1]}^{\star} & 0 & 0\\
0 & 0 & \epsilon \circ \rho_{[2],[1]}& 0\\
0 & 0 & 0 & \epsilon \circ \rho_{[2],[1]}^\star
\end{pmatrix}.$$
It follows that 
$$L_1=(XY)(I_{12}+E_{2,3}-E_{6,5}+E_{8,9}-E_{12,11})(XY)^{-1}\in XG_{12}X^{-1}$$ and $$L_2=I_{12}+E_{2,3}-E_{6,5}+E_{8,9}-E_{12,11}\in XG_{12}X^{-1}.$$
 We then have $L_3=[L_2,L_1]\in XG_{12}X^{-1}$. We have 
 $$L_3=I_{12}+2\frac{\Phi_6(\alpha)\Phi_2(\alpha\beta)}{(\alpha^2\beta+1)(\beta+1)(\alpha+1)^2}E_{6,9}-2\frac{\sqrt{\alpha}\Phi_6(\alpha)\Phi_2(\alpha\beta)^2}{\Phi_2(\alpha)^3(\alpha+\beta)\Phi_2(\alpha^2\beta)\Phi_2(\beta)}E_{12,3}.$$
  Hence, $\dim((L_3-I_{12})V)=2$. We can then apply Theorem \ref{CGFS} to get that $G_{12}$ is a classical group in a natural representation over $\F_{q'}$ for some $q'$ dividing $q$. We have that $q'$ divides $q^{\frac{1}{2}}$ since $G_{12}$ is conjugate to a subgroup of $\Omega_{12}^+(q^{\frac{1}{2}})$. We have
 $$\diag(\alpha,\alpha^{-1},1,\alpha^{-1},\alpha,1,\alpha,\alpha^{-1},1,\alpha^{-1},\alpha,1)\in XG_{12}X^{-1},$$
  $$\diag(\beta,\beta^{-1},1,\beta^{-1},\beta,1,\beta^{-1},\beta,1,\beta,\beta^{-1},1)\in XG_{12}X^{-1}$$
   therefore the field $\F_{q'}$ generated by the traces of the elements of $G_{12}$ contains $4(\alpha+\alpha^{-1}+1)$ and $4(\beta+\beta^{-1}+1)$. It follows that $\F_{q'}$ contains $\F_p(\alpha+\alpha^{-1},\beta+\beta^{-1})=F_{q^{\frac{1}{2}}}$, therefore $q'=q^{\frac{1}{2}}$. We can then conclude that $G_{12}\simeq \Omega_{12}^+(q^{\frac{1}{2}})$.

\smallskip

By complete symmetry of $\alpha$ and $\beta$ for the $12$-dimensional representation, we get that in cases $3$, $7$, $8$ and $14$, we have $G_{12}\simeq \Omega_{12}^+(q^{\frac{1}{2}})$.

\bigskip

Consider now the $16$-dimensional representation. We consider the matrices corresponding to the representation $\rho_{16}$ with respect to the basis $(e_{x_i})_{i\in [\![1,16]\!]}$, where $(I(x_i))_{i\in [\![1,16]\!]}$ is ordered lexicographically i.e. $I(x_1)=\{s_2\}$, $I(x_2)=\{s_3\}$, $I(x_3)=\{s_1,s_2\}$, $I(x_4)=I(x_5)=I(x_6)=\{s_1,s_3\}$, $I(x_7)=I(x_8)=\{s_2,s_4\}$ and $I(x_i)=\{s_1,s_2,s_3,s_4\}\setminus I(x_{17-i})$ for all $i \in [\![9,16]\!]$ and such that $\mu^{s_1}_{x_4,x_2}=1$, $\mu^{s_3}_{x_5,x_3}=-1$, $\mu_{x_{10},x_2}^{s_2}=1$, $\mu_{x_{15},x_7}^{s_2}=1$, $\mu_{x_{14},x_{11}}^{s_3}=2$ and $\mu_{x_{15},x_{13}}^{s_1}=1$. Let $P_{16}$ be the symmetric antidiagonal matrix with coefficients with coefficients in the first $8$ rows respectively equal to $1$, $1$, $-1$, $-1$, $1$, $-2$, $-1$ and $2$. We then have for all $g\in G_{16}$, $PgP^{-1}=^t\!g^{-1}$. Since $P$ is symmetric, we get that $G_{16}$ is up to conjugation a subgroup of $\Omega_{16}^+(\tilde{q})$. Note that by Table \ref{resF4B31} and Theorem \ref{empty}, $G_{16}$ contains up to conjugation a diagonal $SL_2(r)$ or a diagonal $SL_2(r^{\frac{1}{2}})$, where $\F_r=\F_p(\alpha)$ in all cases. It follows that $G_{16}$ is irreducible, primitive, $v_{G_{16}}(V)\leq 2$ and we are not in the second case of Theorem \ref{CGFS}. It only remains to show that $G_{16}$ is tensor-indecomposable in order to apply Theorem \ref{CGFS}. Note that in cases $4$, $5$ and $10$, we have $\Phi\circ \rho_{16|\mathcal{A}_{F_4}} \simeq \rho_{16|\mathcal{A}_{F_4}}$, therefore $G_{16}$ is conjugate to a subgroup of $\Omega_{16}^+(q)$. In all the remaining cases, we have $\tilde{q}=q$. By Lemma \ref{tens2}, we have $G_{16}$ tensor-indecomposable or $G_{16}\leq GL_2(q)\otimes GL_8(q)$. Depending on the cases, we have that $G_{16}$ contains up to conjugation either 
$$H_1=\left\lbrace\begin{pmatrix}
M & 0 & 0\\
0 & ^t\!M^{-1} & 0 \\
0 & 0 & I_{10}
\end{pmatrix}, M\in SL_3(q)\right\rbrace,$$
 $$H_2=\left\lbrace\begin{pmatrix}
M & 0 & 0\\
0 & ^t\!M^{-1} & 0 \\
0 & 0 & I_{10}
\end{pmatrix}, M\in SL_3(q)\right\rbrace$$
 or
  $$H_3=\left\lbrace\begin{pmatrix}M &0 & 0 & 0 & 0\\
0 & ^t\!M^{-1} & 0 & 0 & 0\\
0 & 0 & \epsilon(M) &0 & 0\\
0 & 0 & 0 & \epsilon(^t\!M^{-1}) & 0\\
0 & 0 & 0& 0 & I_4\end{pmatrix}, M\in SL_3(q)\right\rbrace$$
 where $\epsilon$ is the unique automorphism of $\F_q$ ($H_2$ and $H_3$ only appear when $\F_p(\alpha,\beta)\neq \F_p(\alpha+\alpha^{-1},\beta+\beta^{-1})$). Assume by contradiction that $G_{16} \leq GL_2(q)\otimes GL_8(q)$. The restriction of the morphism from $G_{16}$ to $GL_2(q)$ to $H_1$, $H_2$ or $H_3$ is then trivial since $\vert PSL_3(q)\vert$ does not divide $\vert GL_2(q)\vert$. This implies that the eigenvalues of semisimple elements of $H_1$, $H_2$ or $H_3$ all have multiplicity divisible by $8$. This is absurd because there exists $\xi\in \F_q^{\star}$ such that $\xi^2\neq 1$, $H_1$ and $H_2$ contain elements conjugate to $\diag(\xi,\xi,\xi^{-1},\xi^{-1},1,1,1,1,1,1,1,1,1,1,1,1)$ or $H_3$ contains an element conjugate to $\diag(\xi,\xi,\xi,\xi,\xi^{-1},\xi^{-1},\xi^{-1},\xi^{-1},1,1,1,1,1,1,1,1)$. It follows that $G_{16}$ is tensor-indecomposable in all cases. We also have that in all cases $G_{16}$ contains up to conjugation $\diag(\alpha,\alpha,\alpha,\alpha,\alpha^{-1},\alpha^{-1},\alpha^{-1},\alpha^{-1},1,1,1,1,1,1,1,1)$ and \\
 $\diag(\beta,\beta,\beta,\beta,\beta^{-1},\beta^{-1},\beta^{-1},\beta^{-1},1,1,1,1,1,1,1,1)$, therefore the field $\F_{q'}$ generated by the traces of the elements of $G_{16}$ contains $\F_p(\alpha+\alpha^{-1},\beta+\beta^{-1})$. Hence, by Theorem \ref{CGFS}, $G_{16}$ is a classical group in a natural representation over $\F_{q'}$.

\smallskip

Assume now we are in case $1$, $4$, $5$ or $10$. We have $\F_q=\F_p(\alpha+\alpha^{-1},\beta+\beta^{-1})$, therefore $q'=q$ and $G_{16}$ is a classical group in a natural representation over $\F_q$. Since $G_{16}$ is conjugate to a subgroup of $\Omega_{16}^+(q)$, we get that $G_{16}\simeq \Omega_{16}^+(q)$. 

In all the remaining cases, there exists a unique automorphism $\epsilon$ of order $2$ of $\F_q$ and $\epsilon \circ \rho_{16|\mathcal{A}_{F_4}}\simeq \rho_{16|\mathcal{A}_{F_4}}$. It follows that $G_{16}$ is conjugate to a subgroup of $\Omega_{16}^+(q^{\frac{1}{2}})$. It then follows that $q'=q^{\frac{1}{2}}$, therefore $G_{16}\simeq \Omega_{16}^+(q^{\frac{1}{2}})$.
 \end{proof}
 
 \begin{theo2}\label{resultF4}
We write $\F_{\tilde{q}}=\F_p(\sqrt{\alpha},\sqrt{\beta})$, $\F_{r_{\alpha}}=\F_p(\alpha+\alpha^{-1})$ and $\F_{r_{\beta}}=\F_p(\beta+\beta^{-1})$.

In cases $1$, $4$, $5$ and $10$, the morphism from $\mathcal{A}_{F_4}$ to $\mathcal{H}_{F_4,\alpha,\beta}^\star\simeq \underset{\rho \mbox{ irr}}\prod GL_{n_\rho}(\tilde{q})$ factorizes through the surjective morphism
$$\Phi: \mathcal{A}_{F_4} \rightarrow (SL_2(r_{\alpha})\circ SL_2(r_{\beta}))\times SL_4(q)^2\times \Omega_6^+(q)^2$$
$$\times SL_8(q)^2\times SL_9(q)^2\times \Omega_{12}^+(q)\times \Omega_{16}^+(q).$$
In cases $11$, $12$, $13$ and $16$, the morphism from $\mathcal{A}_{F_4}$ to $\mathcal{H}_{F_4,\alpha,\beta}^\star\simeq \underset{\rho \mbox{ irr}}\prod GL_{n_\rho}(\tilde{q})$ factorizes through the surjective morphism
$$\Phi: \mathcal{A}_{F_4} \rightarrow  \times (SL_2(r_{\alpha})\circ SL_2(r_{\beta}))\times SU_4(q^{\frac{1}{2}})^2\times \times \Omega_6^+(q^{\frac{1}{2}})^2$$
$$\times SU_8(q^{\frac{1}{2}})^2\times SU_9(q^{\frac{1}{2}})^2\times \Omega_{12}^+(q^{\frac{1}{2}})\times \Omega_{16}^+(q^{\frac{1}{2}}).$$
In cases $2$, $6$, $9$ and $15$, the morphism from $\mathcal{A}_{F_4}$ to $\mathcal{H}_{F_4,\alpha,\beta}^\star\simeq \underset{\rho \mbox{ irr}}\prod GL_{n_\rho}(\tilde{q})$ factorizes through the surjective morphism
$$\Phi: \mathcal{A}_{F_4} \rightarrow  (SL_2(r_{\alpha})\circ SL_2(r_{\beta}))\times SL_4(q)\times \times \Omega_6^+(q)$$
$$\times SL_8(q^{\frac{1}{2}})\times SU_8(q^{\frac{1}{2}})\times SL_9(q)\times \Omega_{12}^+(q^{\frac{1}{2}})\times \Omega_{16}^+(q^{\frac{1}{2}}).$$
In cases $3$, $7$, $8$ and $14$, the morphism from $\mathcal{A}_{F_4}$ to $\mathcal{H}_{F_4,\alpha,\beta}^\star\simeq \underset{\rho \mbox{ irr}}\prod GL_{n_\rho}(\tilde{q})$ factorizes through the surjective morphism
$$\Phi: \mathcal{A}_{F_4} \rightarrow  (SL_2(r_{\alpha})\circ SL_2(r_{\beta}))\times SL_4(q) \times \Omega_6^+(q)$$
$$\times SL_8(q)\times SU_8(q^{\frac{1}{2}})\times SL_9(q)\times \Omega_{12}^+(q^{\frac{1}{2}})\times \Omega_{16}^+(q^{\frac{1}{2}}).$$
 \end{theo2}
 
 \begin{proof}
Note the by \cite{MR}, $\mathcal{A}_{F_4}$ is not perfect, this makes the proof of this theorem more complex than for the previous types. The result follows easily from Goursat's Lemma in cases $2$, $3$, $6$, $7$, $8$, $9$, $14$ and $15$.

\medskip

Assume now we are in case $1$, $4$ $5$ or $10$. We must show that there exists no field automorphism $\Psi$ and no character $z:\mathcal{A}_{F_4}\rightarrow \F_q^\star$ such that $\rho_{4_2|\mathcal{A}_{F_4}}\simeq (\Psi \circ \rho_{4_4|\mathcal{A}_{F_4}})\otimes z$ or $\rho_{4_2|\mathcal{A}_{F_4}}\simeq (\Psi \circ \rho_{4_4'|\mathcal{A}_{F_4}})\otimes z$ or $\rho_{6_1|\mathcal{A}_{F_4}}\simeq (\Psi \circ \rho_{6_2|\mathcal{A}_{F_4}})\otimes z$ or $\rho_{8_1|\mathcal{A}_{F_4}}\simeq (\Psi \circ \rho_{8_2|\mathcal{A}_{F_4}})\otimes z$ or $\rho_{8_1|\mathcal{A}_{F_4}}\simeq (\Psi \circ \rho_{8_2'|\mathcal{A}_{F_4}})\otimes z$ or $\rho_{9_1|\mathcal{A}_{F_4}}\simeq (\Psi \circ \rho_{9_2|\mathcal{A}_{F_4}})\otimes z$ or $\rho_{9_1|\mathcal{A}_{F_4}}\simeq (\Psi \circ \rho_{9_2'|\mathcal{A}_{F_4}})\otimes z$. Assume by contradiction that there exists such a field automorphism $\Psi$ and such a character $z$. By \cite{MR}, $\mathcal{A}_{F_4}$ is generated by $p_0=S_2S_1^{-1}$, $p_1=S_1S_2S_1^{-2}$, $q_0=S_3S_4^{-1}$ and $q_1=S_4S_3S_4^{-2}$.

\medskip

Assume first that $\rho_{4_2|\mathcal{A}_{F_4}}\simeq (\Psi \circ \rho_{4_4|\mathcal{A}_{F_4}})\otimes z$.
We have $\tr(\rho_{4_2}(p_0))= \tr(\rho_{4_4}(p_0))=3-(\alpha+\alpha^{-1})$ and $\tr(\rho_{4_2}(q_0))= \tr(\rho_{4_4}(q_0))=3-(\beta+\beta^{-1})$

\smallskip

\textbf{1.} Assume $3-\alpha-\alpha^{-1}\neq 0$. We have 
$$\tr(\rho_{4_2}(p_1))=\tr(\rho_{4_4}(p_1))=\tr(\rho_{4_2}(p_0p_1^{-1}))=\tr(\rho_{4_4}(p_0p_1^{-1}))=3-\alpha-\alpha^{-1}\neq 0.$$
 It follows that $z(p_0)=z(p_1)=z(p_0p_1^{-1})=\frac{3-\alpha-\alpha^{-1}}{\Psi(3-\alpha-\alpha^{-1})}$ therefore $z(p_1)^{-1}=1$ and $z(p_0)=z(p_1)=1$.

\smallskip

\textbf{1.1.} Assume $3-(\beta+\beta^{-1})\neq 0$. We have
$$\tr(\rho_{4_2}(q_1))=\tr(\rho_{4_4}(q_1))=\tr(\rho_{4_2}(q_0q_1^{-1}))=\tr(\rho_{4_4}(q_0q_1^{-1}))=3-\beta-\beta^{-1}\neq 0.$$
 It follows that $z(q_0)=z(q_1)=z(q_0q_1^{-1})=\frac{3-\beta-\beta^{-1}}{\Psi(3-\beta-\beta^{-1})}$ therefore $z(q_1)^{-1}=1$ and $z(q_0)=z(q_1)=1$. This proves that $z$ is the trivial character. It follows that $\Psi(\alpha+\alpha^{-1})=\alpha+\alpha^{-1}$ and $\Psi(\beta+\beta^{-1})=\beta+\beta^{-1}$. This implies that $\Psi$ is the trivial automorphism. This would imply that $\rho_{4_2|\mathcal{A}_{F_4}}\simeq \rho_{4_4|\mathcal{A}_{F_4}}$ which is absurd by Proposition \ref{resH4derivedsubgroup}.

\smallskip

\textbf{1.2.} Assume $A=3-\beta-\beta^{-1}=0$. We have 
$$\tr(\rho_{4_2}(q_0^2))=\tr(\rho_{4_4}(q_0^2))=\tr(\rho_{4_2}(q_1^2))=\tr(\rho_{4_4}(q_1^2))=(A-2)^2=4\neq 0.$$
 It follows that $z(q_0^2)=z(q_1^2)=\frac{4}{\Psi(4)}=1$. This implies that  $z(q_0)\in \{\pm 1\}$. We have 
 $$\tr(\rho_{4_2}(q_0q_1^2))=\tr(\rho_{4_4}(q_0q_1^2))=A^3-7A^2+15A-8=-8.$$
  It follows that $z(q_0q_1^2)=1$ therefore $z(q_0)=1$. We also have $$\tr(\rho_{4_2}(q_1q_0^2))=\tr(\rho_{4_4}(q_1q_0^2))=A^3-7A^2+15A-8=-8$$
   therefore $z(q_1)=1$. We get as before that $z$ is trivial. This implies that $\rho_{4_2|\mathcal{A}_{F_4}}\simeq \Psi\circ\rho_{4_4|\mathcal{A}_{F_4}}$. It then follows from Proposition \ref{Fieldfactorization} that $\Psi$ is trivial which is absurd by Proposition \ref{resH4derivedsubgroup}.

\smallskip

\textbf{2.} Assume $B=3-\alpha-\alpha^{-1}=0$. We have 
$$\tr(\rho_{4_2}(p_0^2))=\tr(\rho_{4_4}(p_0^2))=\tr(\rho_{4_2}(p_1^2))=\tr(\rho_{4_4}(p_1^2))=(B-2)^2=4\neq 0.$$
 It follows that $z(p_0^2)=z(p_1^2)=\frac{4}{\Psi(4)}=1$. This implies that  $z(p_0)\in \{\pm 1\}$. We have 
 $$\tr(\rho_{4_2}(p_0p_1^2))=\tr(\rho_{4_4}(p_0p_1^2))=B^3-7B^2+15B-8=-8.$$
  It follows that $z(p_0p_1^2)=1$ therefore $z(p_0)=1$. We also have $$\tr(\rho_{4_2}(p_1p_0^2))=\tr(\rho_{4_4}(p_1p_0^2))=B^3-7B^2+15B-8=-8$$
   therefore $z(p_1)=1$.

\smallskip

\textbf{2.1.} Assume $3-\beta-\beta^{-1}\neq 0$. By the same computations as in \textbf{1.1.} we have that $z(q_0)=z(q_1)=1$. It follows that $z$ is trivial. This implies that $\rho_{4_2|\mathcal{A}_{F_4}}\simeq \Psi\circ\rho_{4_4|\mathcal{A}_{F_4}}$. It then follows from Proposition \ref{Fieldfactorization} that $\Psi$ is trivial which is absurd by Proposition \ref{resH4derivedsubgroup}.

\smallskip

\textbf{2.2.} Assume $3-\beta-\beta^{-1}=0$. We then have $\alpha+\alpha^{-1}=3=\beta+\beta^{-1}$. This implies that $\alpha\in \{\beta,\beta^{-1}\}$ and contradicts our assumptions on $\alpha$ and $\beta$.

\medskip

Note that all the computations above were made in $\mathcal{A}_{A_2,1}$  or $\mathcal{A}_{A_2,2}$. Since the restrictions to those subgroups are stable by the transposed inverse operation, we have that $\rho_{4_2|\mathcal{A}_{F_4}}\simeq (\Psi \circ \rho_{4_4'|\mathcal{A}_{F_4}})\otimes z$.

\medskip

Assume now $\rho_{6_1|\mathcal{A}_{F_4}}\simeq (\Psi \circ \rho_{6_2|\mathcal{A}_{F_4}})\otimes z$. We have 
\begin{tiny}
$$\tr(\rho_{6_1}(p_0))=\tr(\rho_{6_1}(p_1))=\tr(\rho_{6_1}(p_0p_1^{-1}))=\tr(\rho_{6_2}(p_0))=\tr(\rho_{6_2}(p_1))=\tr(\rho_{6_2}(p_0p_1^{-1}))=-2\alpha^{-1}(\alpha-1)^2\neq 0.$$
$$\tr(\rho_{6_1}(q_0))=\tr(\rho_{6_1}(q_1))=\tr(\rho_{6_1}(q_0q_1^{-1}))=\tr(\rho_{6_2}(q_0))=\tr(\rho_{6_2}(q_1))=\tr(\rho_{6_2}(q_0q_1^{-1}))=-2\beta^{-1}(\beta-1)^2\neq 0.$$
\end{tiny}
It follows that $z(p_0)=z(p_1)=z(q_0)=z(q_1)=1$. We then have by Proposition \ref{Fieldfactorization} that $\Psi$ is trivial therefore $\rho_{6_1|\mathcal{A}_{F_4}}\simeq \rho_{6_2|\mathcal{A}_{F_4}}$.

\medskip

Assume now $\rho_{8_2|\mathcal{A}_{F_4}}\simeq (\Psi \circ \rho_{8_3|\mathcal{A}_{F_4}})\otimes z$. We have 
$$\tr(\rho_{8_2}(p_0))=\tr(\rho_{8_2}(p_1))=\tr(\rho_{8_2}(p_0p_1^{-1}))=6-2(\alpha+\alpha^{-1}),$$
$$\tr(\rho_{8_3}(p_0))=\tr(\rho_{8_3}(p_1))=\tr(\rho_{8_3}(p_0p_1^{-1}))=5-3(\alpha+\alpha^{-1}).$$
 Note that those quantities are non-zero because we have 
 $$\tr(\rho_{8_2}(p_0))=z(p_0)\Psi(\tr(\rho_{8_3}(p_0))),$$
 $$3\tr(\rho_{8_2}(p_0))-2\tr(\rho_{8_2}(p_0))=18-10=8\neq 0.$$
  It follows that $z(p_0)=z(p_1)=z(p_0p_1^{-1})=1$. We also have
  $$\tr(\rho_{8_2}(q_0))=\tr(\rho_{8_2}(q_1))=\tr(\rho_{8_2}(q_0q_1^{-1}))=6-2(\beta+\beta^{-1}),$$
  $$\tr(\rho_{8_3}(q_0))=\tr(\rho_{8_3}(q_1))=\tr(\rho_{8_3}(q_0q_1^{-1}))=5-3(\beta+\beta^{-1}).$$
   It follows that $z(q_0)=z(q_1)=1$. We then have that $z$ is trivial and $\Psi$ is trivial therefore $\rho_{8_2|\mathcal{A}_{F_4}}\simeq \rho_{8_3|\mathcal{A}_{F_4}}$ which contradicts Proposition \ref{resH4derivedsubgroup}.

\medskip

All the computations above give the same results if we substitute $\rho_{8_3}$ by $\rho_{8_3'|\star}$ therefore $\rho_{8_2|\mathcal{A}_{F_4}}\not\simeq (\Psi \circ \rho_{8_3'|\mathcal{A}_{F_4}})\otimes z$.

\medskip

Assume now $\rho_{9_1|\mathcal{A}_{F_4}}\simeq (\Psi \circ \rho_{9_2|\mathcal{A}_{F_4}})\otimes z$. We have 
\begin{footnotesize}
$$\tr(\rho_{9_1}(p_0))=\tr(\rho_{9_1}(p_1))=\tr(\rho_{9_1}(p_0p_1^{-1}))=\tr(\rho_{9_2}(p_0))=\tr(\rho_{9_2}(p_1))=\tr(\rho_{9_2}(p_0p_1^{-1}))=-3\alpha^{-1}(\alpha-1)^2\neq 0.$$
$$\tr(\rho_{9_1}(q_0))=\tr(\rho_{9_1}(q_1))=\tr(\rho_{9_1}(q_0q_1^{-1}))=\tr(\rho_{9_2}(q_0))=\tr(\rho_{9_2}(q_1))=\tr(\rho_{9_2}(q_0q_1^{-1}))=-3\beta^{-1}(\beta-1)^2\neq 0.$$
\end{footnotesize}
It follows that $z(p_0)=z(p_1)=z(q_0)=z(q_1)=1$. We then have by Proposition \ref{Fieldfactorization} that $\Psi$ is trivial therefore $\rho_{9_1|\mathcal{A}_{F_4}}\simeq \rho_{9_2|\mathcal{A}_{F_4}}$. All the computations give the same results if we substitute $\rho_{9_3}$ by $\rho_{9_3'}$. This concludes the proof in cases $1$, $4$, $5$ and $10$.

\bigskip

In cases $11$, $12$, $13$ and $16$, we have the result using the same computations since $\Psi$ is either the trivial automorphism or the unique automorphism $\epsilon$ of order $2$ of $\F_q$ if $z$ is trivial. \end{proof}
 
 \newpage
 
 $ $ 
 
 \newpage

\chapter{Appendix}

\section{Erratum of the papers \cite{BM} and \cite{BMM} in type $A$}\label{ErratBMBMM}

This Ph.D. thesis is mainly based on \cite{BM} and \cite{BMM}. It uses many results from those papers. All the results are correct, however some proofs are incomplete. We give here a few corrections of those papers.

\medskip

In the proof of Proposition 3.1 of \cite{BM}, $B_{n-1}$ should be replaced by $B_n$.

\smallskip

In section 3.2. of \cite{BM}, $R(B)$ should be replaced by $-R(B_n)$. The determinant of a matrix $M\in GU_n(q^{\frac{1}{2}})$ does not verify $\epsilon (\det(M))=\det(M)$ but $\epsilon(\det(M))=\det(M)^{-1}$. We have $\det(-R(\sigma_i))=-\alpha $ therefore $\epsilon(-\alpha)=-\alpha^{-1}$. It follows that $\epsilon(\alpha+\alpha^{-1})=\alpha+\alpha^{-1}$ and the conclusion remains true since $\F_q=\F_p(\alpha+\alpha^{-1})$.

\smallskip

In the proof of Lemma $3.5$ of \cite{BM}, in order to prove that $[\overline{G},\overline{G}]$ is not abelian by abelian, it is necessary to check that $(s_1s_2^{-1})(s_1^{-1}s_2)\pm(s_1^{-1}s_2)(s_1s_2^{-1})\neq 0$. In the paper, it was only checked that $(s_1s_2^{-1})(s_1^{-1}s_2)-(s_1^{-1}s_2)(s_1s_2^{-1})\neq 0$ and the matrix given is wrong. We have 
$$(s_1s_2^{-1})(s_1^{-1}s_2)-(s_1^{-1}s_2)(s_1s_2^{-1})=\begin{pmatrix}
\frac{\alpha^4+\alpha^3-\alpha-1}{\alpha^2} & \frac{\alpha^3-1}{\alpha^2}\\
-\frac{\alpha^5+\alpha^4+\alpha^3-\alpha^2-\alpha-1}{\alpha^2} & -\frac{\alpha^4+\alpha^3-\alpha-1}{\alpha^2}
\end{pmatrix}$$
This matrix is non-zero because it is assumed in the Lemma that the order of $\alpha$ is not in $\{1,2,3,4,5,6,10\}$. We have
$$(s_1s_2^{-1})(s_1^{-1}s_2)+(s_1^{-1}s_2)(s_1s_2^{-1})=\begin{pmatrix}
\frac{\alpha^4-\alpha^3+2\alpha^2-\alpha+1}{\alpha^2} & \frac{\alpha^3+1}{\alpha^2}\\
-\frac{\alpha^5+\alpha^4+\alpha^3+\alpha^2+\alpha+1}{\alpha^2} & -\frac{\alpha^4+\alpha^3+\alpha+1}{\alpha^2}
\end{pmatrix}$$
This matrix is also non-zero by the conditions on the order of $\alpha$.

\smallskip

In the proof of Proposition $3.7$, it is not proved that $\epsilon \circ \rho^\star$ is a representation of the Hecke algebra. However, it is true by Proposition \ref{Fieldfactorization} of this Ph.D. thesis since $\epsilon(\alpha)=\alpha^{-1}$.

\smallskip

In the proof of section $4$, it should be $k=\ulcorner \frac{N}{2}\urcorner$ instead of $k=\llcorner \frac{N}{2} \lrcorner$ in order to prove that $G_r\not\simeq A\rtimes \mathfrak{S}_N$ has less transvection than the natural $SL_{a}{q}$ with $a\geq \frac{N}{2}$. The proof for $N=5$ then becomes a subcase of the proof for $N\geq 6$ because we have that $k=3$ when $N=5$.

\medskip

In Lemma $3.2$ of \cite{BMM}, the action of $\mathcal{L}$ satisfies $\frac{\mathcal{L}s_r\mathcal{L}^{-1}}{-\alpha}=^t\!s_r^{-1}$ but not $\frac{\mathcal{L}s_r\mathcal{L}^{-1}}{-\alpha\nu(\lambda)}=^t\!s_r^{-1}$.

\smallskip

After Theorem $5.1$ and Lemma $5.5$ of \cite{BM} it should say that we need to prove that $G$ is primitive not imprimitive. 

\smallskip

In \cite{BM}, after Lemma $5.3$, Proposition $2.4.$ is used in order to show that a morphism from $R(\mathcal{A}_{A_n})$ to $PSL_2(q)$ must be trivial. It is used for $n\geq 6$ eventhough Proposition $2.4$ is only valid for $n\geq 7$. The proof is therefore not correct for $n=6$. The only partition to consider is $[3,2,1]$ which is of dimension $16$. By the induction assumption, we have $R_{[3,2,1]}(\mathcal{A}_{A_5})\simeq SL_5(q)\times SL_6(q)$. It follows by Lemma \ref{exceptiontens2} of this Ph.D. thesis that the morphism to $SL_2(q)$ is indeed trivial and concludes the proof.

\section{$H_4$-graphs}\label{sectionnewH4graphs}
Below are the $W$-graphs we use for type $H_4$, they are taken from \cite{G-P} except for the new ones we found verifying the properties of Theorem \ref{bilinwgraphs} ($\widetilde{16_t}$, $\widetilde{18_r}$, $\widetilde{\overline{24_s}}$, $\widetilde{\overline{24_t}}$, $\widetilde{30_s}$, $\widetilde{40_r}$ and $\widetilde{48_{rr}}$. We rearranged the vertices in the ones taken from \cite{G-P} in order to male the connected components of the restrictions to $\mathcal{A}_{H_3}$ appear. In what follows $\beta=\xi+\xi^{-1}$ and if the graph is $2$-colorable then we give a possible coloring.

\bigskip

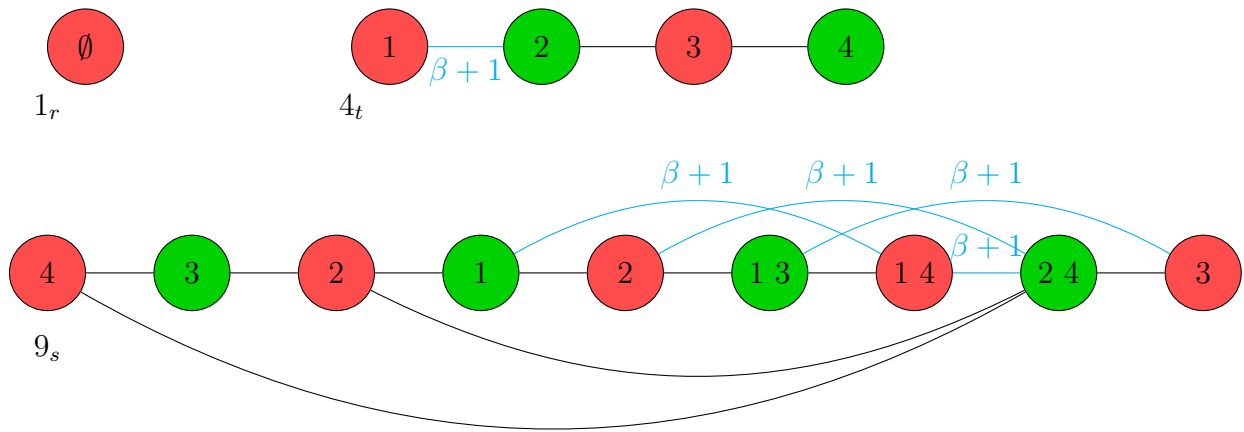
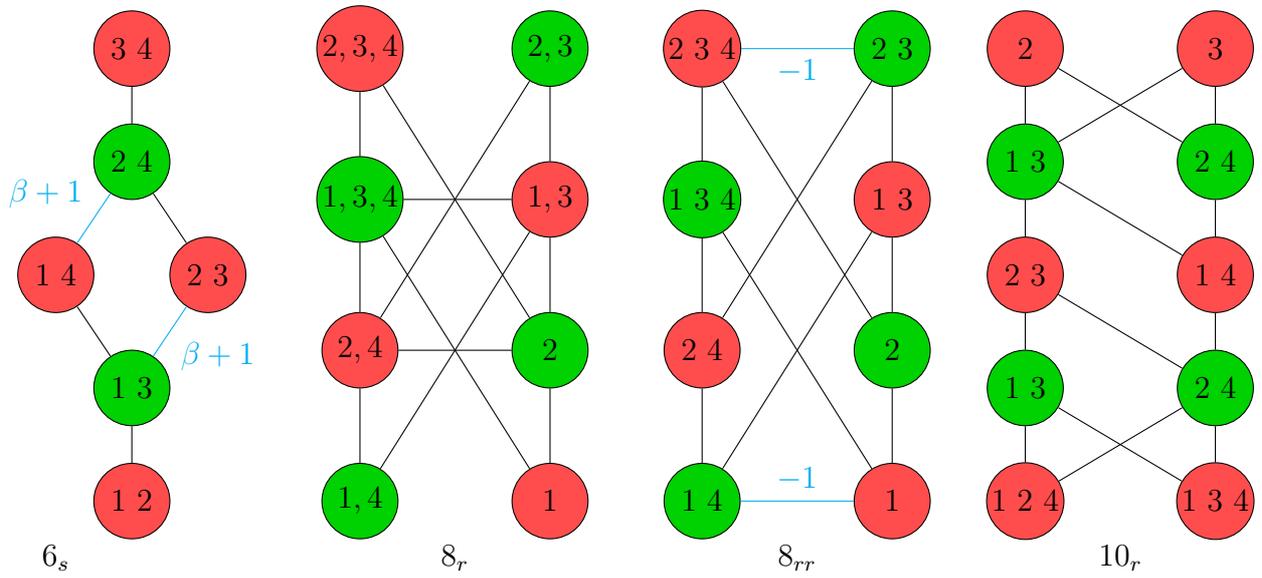
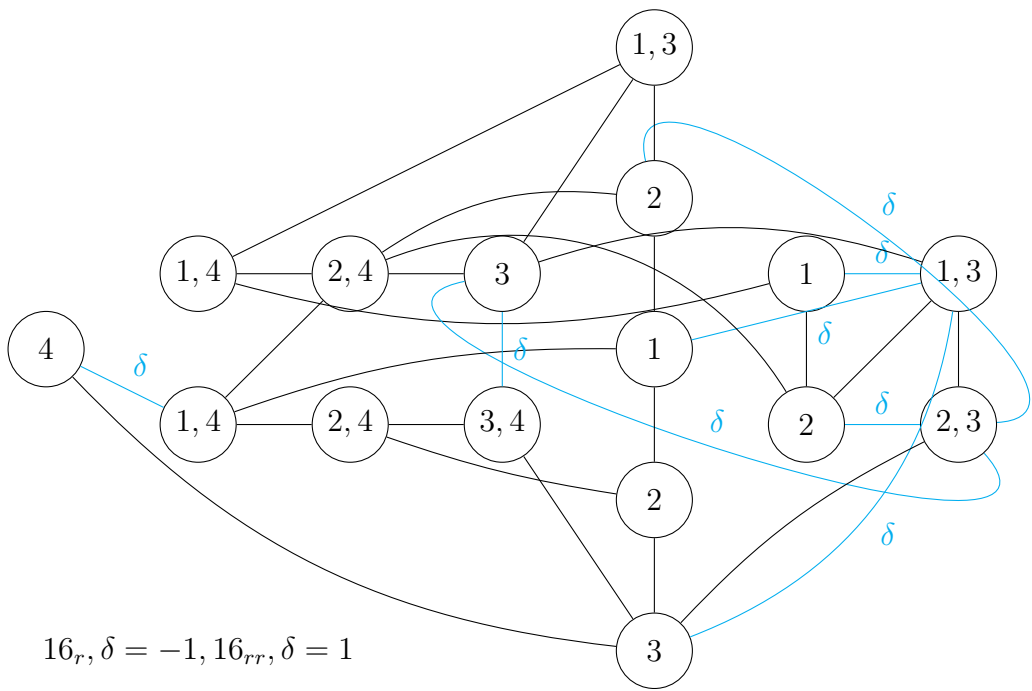
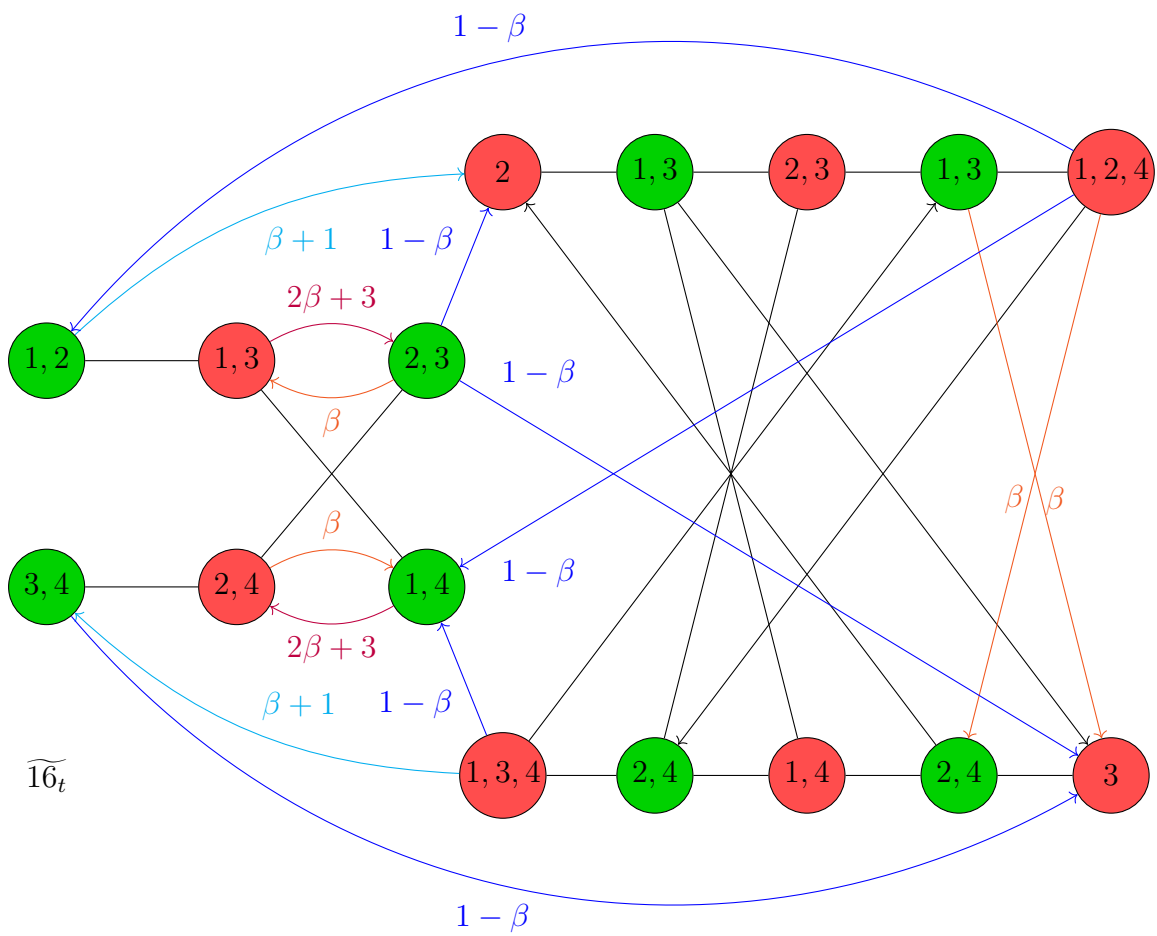
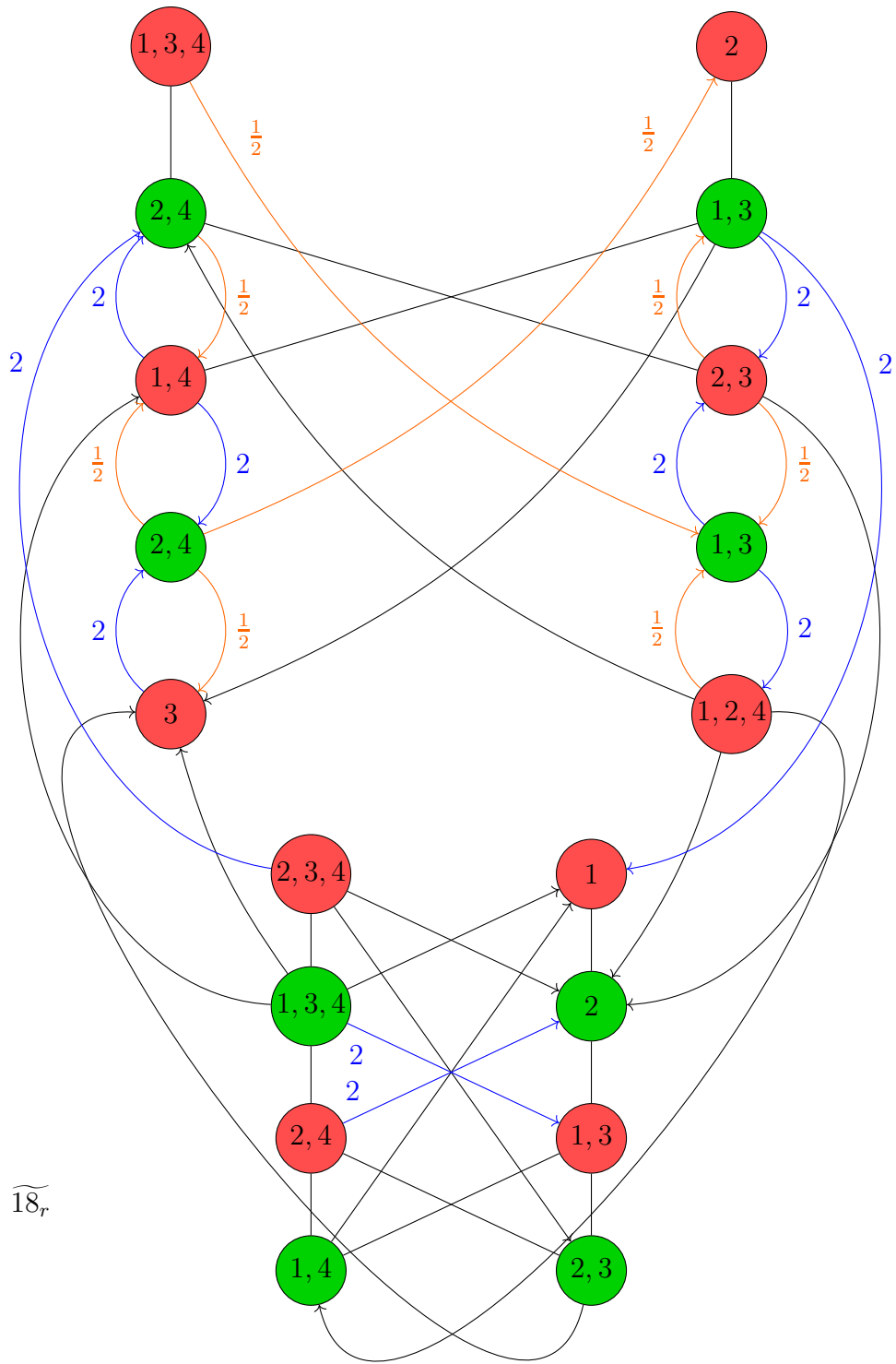
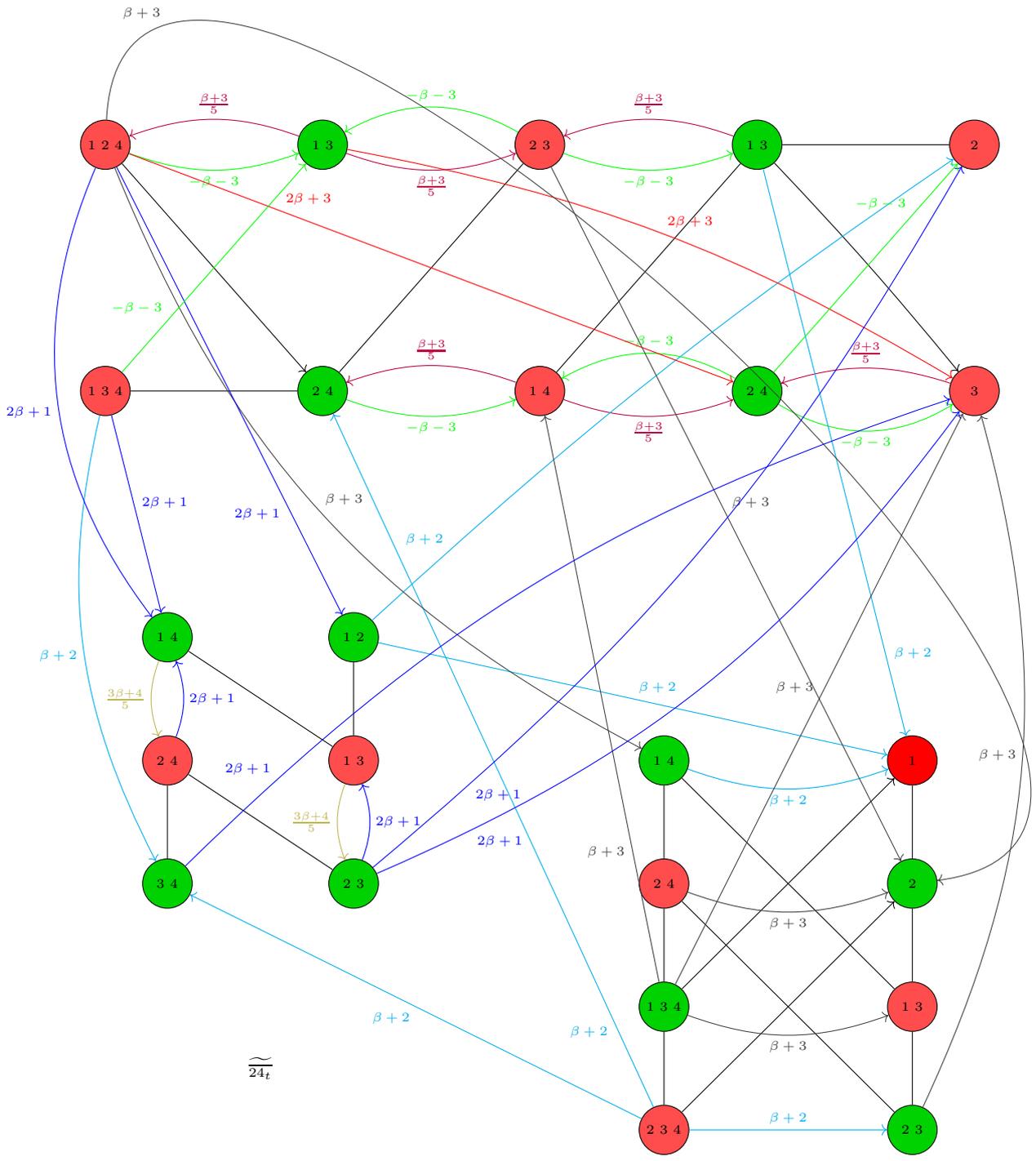
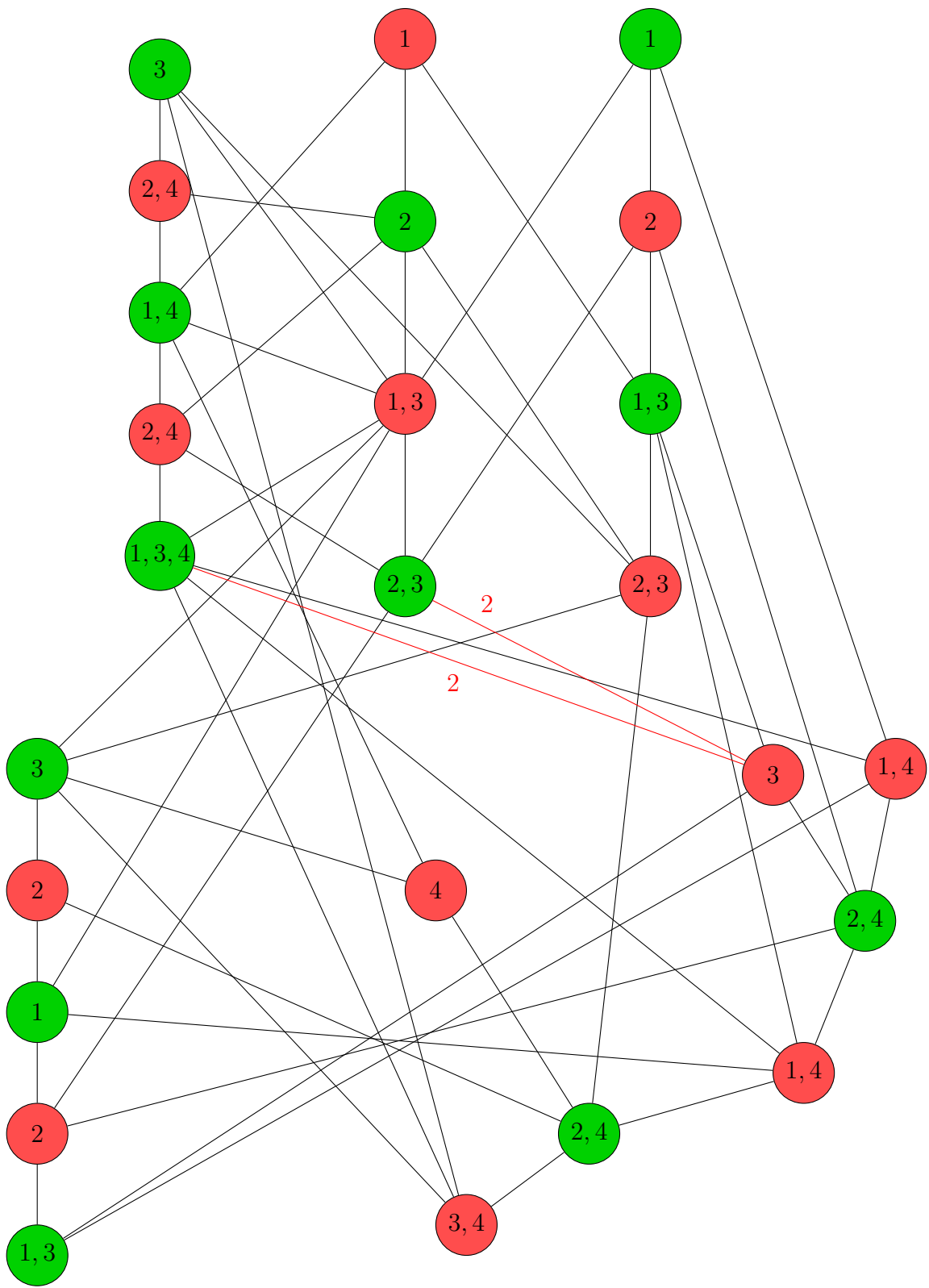
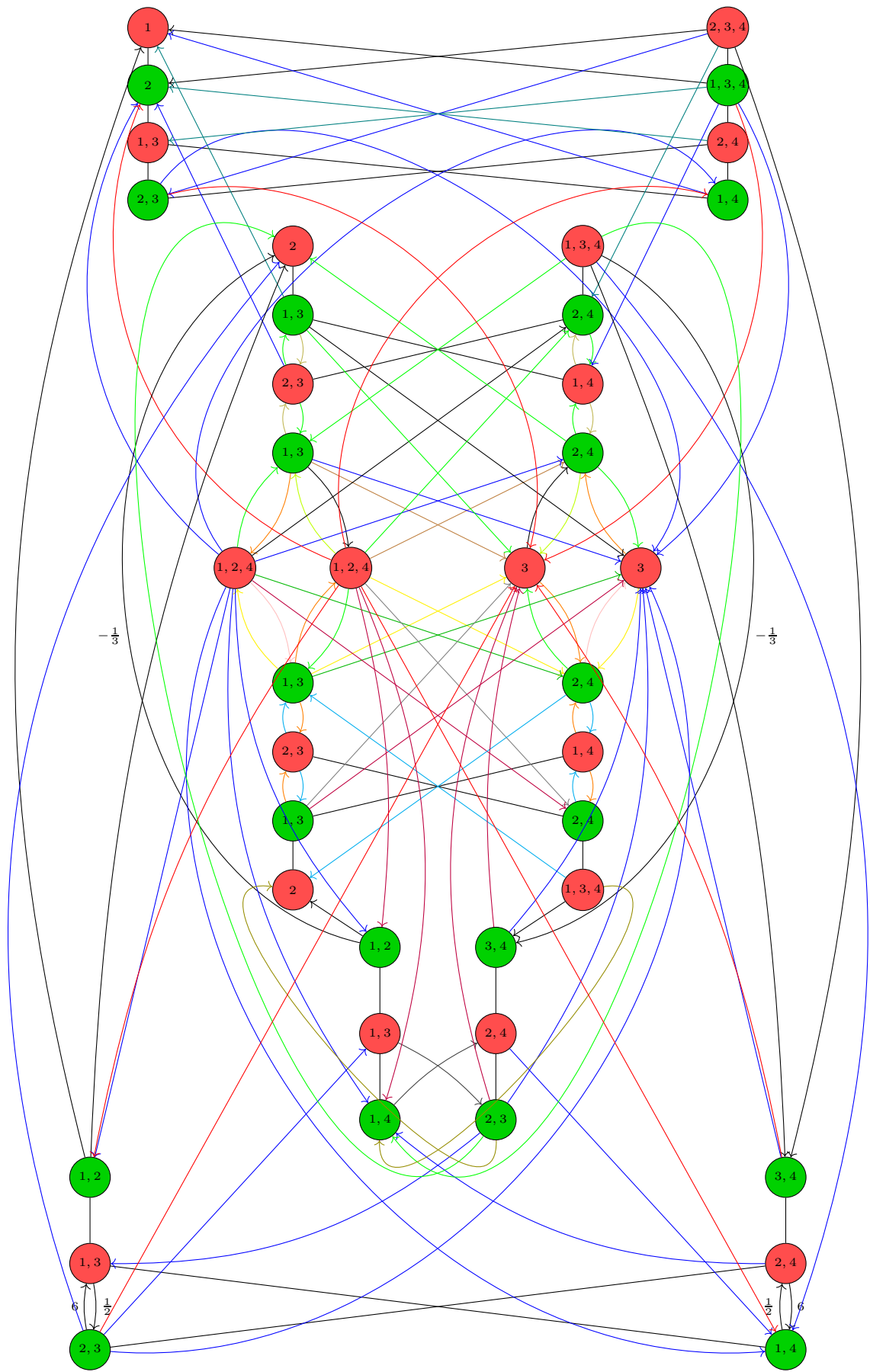
\begin{figure}[h]
\centering

\end{tiny}
\end{figure}

\clearpage
\newpage

In the 40-dimensional $H_4$-graph, the blue edges are of weight $2$, the red edges are of weight $\frac{4}{3}$, the cyan edges are of weight $\frac{6}{5}$, the orange edges are of weight $\frac{5}{6}$, the green edges are of weight $\frac{2}{3}$, the dark yellow edges are of weight $\frac{3}{2}$, the purple edges are of weight $-\frac{2}{3}$, the dark gray edges are of weight $-\frac{1}{2}$, the yellow edges are of werght $-\frac{5}{9}$, the dark green edges are of weight $\frac{5}{3}$, the teal edges are of weight $3$, the lime edges are of weight $\frac{4}{9}$, the pink edges are of weight $-\frac{4}{5}$, the brown edges are of weight $\frac{7}{3}$ and the olive edges are of weight $-2$.

In the 48-dimensional $H_4$-graph, we have omitted the weights on the edges for clarity. One can still observe the symmetry appearing in this self-dual $H_4$-graph.

\section{Computations in type $H_4$}\label{ComputationsH4}

\subsection{Proof that $E_3$ and $D_3$ cannot vanish simultaneously}\label{computationD3E3}

\begin{small}
\begin{eqnarray*}-\frac{v^{22}}{\Phi_6(\alpha)}D_3& = &v^{40}-2v^{38}+4v^{36}-8v^{34}+11v^{32}-16v^{30}+22v^{28}-26v^{26}
+30v^{24}-34v^{22}+34v^{20}-34v^{18}\\
& + & 30v^{16}-26v^{14}+22v^{12}-16v^{10}
+11v^8-8v^6+4v^4-2v^2+1\\
-\frac{v^{18}}{\Phi_6(\alpha)}E_3 & = & v^{32}-v^{30}-v^{28}-v^{26}+4v^{24}-v^{22}-v^{20}-3v^{18}+4v^{16}-3v^{14}-v^{12}-v^{10}+4v^8-v^6-v^4-v^2+1.
\end{eqnarray*}
\end{small}

\begin{tiny}
\begin{eqnarray*}
v_5 & = & \op{Rem}(-\frac{v^{22}}{\Phi_6(\alpha)}D_3,-\frac{v^{18}}{\Phi_6(\alpha)}E_3)\\
& = & -5v^{30}+8v^{28}+2v^{26}-v^{24}-13v^{22}+10v^{20}+12v^{18}-7v^{16}-3v^{14}+8v^{12}+10v^{10}-15v^8-v^6+2v^4+8v^2-5\\
v_6 & = & 5^2\op{Rem}(-\frac{v^{18}}{\Phi_6(\alpha)}E_3)\\
&  = & 9v^{28}-24v^{26}+32v^{24}-14v^{22}+65v^{20}-74v^{18}+64v^{16}-44v^{14}+49v^{12}
-70v^{10}+50v^8-18v^6+21v^4-26v^2+10\\
v_7 & = & \frac{3^3}{5^2}\op{Rem}(v_5,v_6)\\
 &= & 6v^{26}+11v^{24}+16v^{22}+8v^{20}+4v^{18}+7v^{16}-2v^{14}-2v^{12}-4v^{10}+5v^8
 -2v^2+1\\
 v_8 & = & 2^2\op{Rem}(v_6,v_7)\\
 & = & 329v^{24}+328v^{22}+452v^{20}-230v^{18}+457v^{16}-218v^{14}+166v^{12}-418v^{10}
 +335v^8-72v^6+96v^4-164v^2+67\\
 v_9 & = & \frac{7^2 47^2}{2^3 3^3}\op{Rem}(v_7,v_8)\\
 &= & 1380v^{22}+2656v^{20}-414v^{18}+2007v^{16}-853v^{14}+1549v^{12}-1871v^{10}
 +603v^8-327v^6+765v^4-361v^2-11\\
 v_{10} & = & \frac{2^2 3^2 5^2 23^2}{7^2 47^2}\op{Rem}(v_8,v_9)\\
  & = & 5006v^{20}-3519v^{18}+4857v^{16}-3413v^{14}+4199v^{12}-4291v^{10}+2403v^8
  -1437v^6+1545v^4-1061v^2+284\\
  v_{11} & =& \frac{2503^2}{3^2 5^2 23^2}\op{Rem}(v_9,v_{10})\\
  & =& 41901v^{18}-30018v^{16}+24300v^{14}-16298v^{12}+30252v^{10}-39028v^8+15158v^6
  -3244v^4+17330v^2-11407\\
  v_{12} & = & \frac{3^1 13967^2}{2503^2}\op{Rem}(v_{10},v_{11})\\
  & = & 187017v^{16}-140575v^{14}+57067v^{12}+30187v^{10}+61161v^8-100305v^6
  -48597v^4+25593v^2+28241\\
  v_{13} & = & \frac{3^1 17^2 19^2 193^2}{2^2 13967^2}\op{Rem}(v_{11},v_{12})\\
  & = & 188627v^{14}-351293v^{12}+243691v^{10}-254562v^8+400992v^6-128403v^4
  +161367v^2-173764\\
  v_{14} & =& \frac{53^2 3559^2}{3^1 17^2 19^2 193^2}\op{Rem}(v_{12},v_{13})\\
  & = & 617414v^{12}+43395v^{10}-171149v^8-1265240v^6-205046v^4+61567v^2+670171\\
  v_{15} & =& \frac{2^2 7^2 44101^2}{3^1 53^2 3559^2}\op{Rem}(v_{13},v_{14})\\
  & =& 1148529v^{10}+110455v^8-1012180v^6-958108v^4-25091v^2+792599\\
  v_{16} & =& \frac{3^2 382843^2}{2^2 7^2 44101^2}\op{Rem}(v_{14},v_{15})\\
  & = & 1295951v^8-2644727v^6-709010v^4-1262572v^2+2357251\\
  v_{17} & =& \frac{881^2 1471^2}{3^2 382843^2}\op{Rem}(v_{15},v_{16})\\
  & = & 5888346v^6+1914360v^4+352576v^2-4674729\\
  v_{18} & = & \frac{2^1 3^1 981391^2}{881^2 1471^2}\op{Rem}(v_{16},v_{17})\\
  & = & 723244v^4-172511v^2-264507\\
  v_{19} & = &\frac{2^3 180811^2}{599^1 981391^2}\op{Rem}(v_{17},v_{18})\\
  &= & 1495v^2-1569\\
  v_{20} & = & \frac{5^2 13^2 23^2}{2^3 3^1 180811^2}\op{Rem}(v_{18},v_{19})\\
  &= & 1.  
\end{eqnarray*}
\end{tiny}

\subsection{Computations in case 2.1.}\label{computationscase21}

\begin{tiny}
\begin{eqnarray*}
E_9 & = & -\frac{10-6(4-E_1)+E_3-E_1^2+(4-E_1)^2}{2-E_1}\\
& = & v^8+v^{-8}\\
D_9 & = &  -\frac{10-6(4+D_1)+D_3-D_1^2+(4+D_1)^2}{2+D_1}\\
& = & \Phi(E_9)\\
E_{10} & = &  \tr([n_2,u_2][u_2^{-1},n_2])-E_1E_9^2-2E_1^2+2E_9^2\\
D_{10} & = & \tr([n_1,u_1][u_1^{-1},n_1])+D_1D_9^2-2D_1^2+2D_9^2\\
& = & \Phi(E_{10})\\
E_{11} & = & E_9^3+4E_1^2+E_1E_{10}-14E_1E_9-4E_{10}E_9-2E_9^2-10E_{10}\\
D_{11} & = & D_9^3+4D_1^2-D_1D_{10}+14D_1D_9-4 D_{10} D_9-2 D_9^2-10 D_{10}\\
& = & \Phi(E_{11})\\
E_{12} & = &E_1 E_{10}+E_9^3+4 E_1 E_{11}-180 E_1^2+18 E_1 E_9-2 E_9^2-18 E_{10}-9 E_{11}+(8 56) E_1-40 E_9-208\\
&= & 0\\
D_{12} & = & \frac{(D_1+2) v^2 (-D_1 D_{10}+D_9^3-4 D_1 D_{11}-180 D_1^2-18 D_1 D_9-2 D_9^2-18 D_{10}-9 D_{11}-(8 56) D_1-40 D_9-208)}{4\Phi_6(\alpha)}\\
& = & 0\\
E_{13} & = &-E_1^2 E_9+34 E_1^2+12 E_1 E_{10}+4 E_1 E_9+E_{10}^2-64 E_1-12 E_{10}-4 E_9+28\\
& = & 0\\
D_{13} & = & \frac{v^2 (D_1+2)^2 (-D_1^2 D_9+34 D_1^2-12 D_1 D_{10}-4 D_1 D_9+D_{10}^2+64 D_1-12 D_{10}-4 D_9+28)}{\Phi_6(\alpha)}\\
D_{14} & = & -\frac{1}{2}(D_{13}-(2 (v^{16}+v^{-16})-v^{12}-v^{-12}-4 (v^{10}+v^{-10})-26 (v^8+v^{-8})\\
& & +4 (v^6+v^{-6})+13 (v^4+v^{-4})+152 (v^2+v^{-2})+184) D_{12})\\
D_{15} & = & 76 ((2 19) D_{12}-(v^2+v^{-2}) D_{14})+787 D_{14}\\
D_{16} & = & -\frac{1}{2^3 19^2}(453385 (453385 D_{14}-(76 (v^2+v^{-2})) D_{15})-250464003 D_{15})\\
D_{17} & = & \frac{1}{2^1 5^2 90677^2}(805808162 (805808162 D_{15}-(453385 (v^2+v^{-2})) D_{16})-245248536209874 D_{16})\\
D_{18} & = & -\frac{1}{402904081^2}(109342364049/2 (437369456196 D_{16}-(805808162 (v^2+v^{-2})) D_{17})-57979709384117131616 D_{17})\\
D_{19} & = & -\frac{1}{3^8 1349905729^2}(1835659768232 (1835659768232 D_{17}-(109342364049 (v^2+v^{-2})) D_{18})\\
& & -470719001514611845893347 D_{18})\\
D_{20} & = & \frac{1}{2^5 233^2 984796013^2}(26716821821819 (26716821821819 D_{18}-(1835659768232 (v^2+v^{-2})) D_{19})\\
& & +246084500040920924910400141 D_{19})\\
D_{21} &= &\frac{1}{1291^2 20694672209^2}(7431837830116907 (37159189150584535 D_{19}-(106867287287276 (v^2+v^{-2})) D_{20})\\
& & -161710748575151990511450485922 D_{20})\\
D_{22} & = &\frac{1}{2^1 64271^2 362953^2 318589^2}(207561548775276615471 (207561548775276615471 D_{20}-(37159189150584535 (v^2+v^{-2})) D_{21})\\
& & -(2 12489945002283355770270866923659197255) D_{21})\\
D_{23} & = & \frac{1}{(2^1 3^2 5^1 229^2 302127436354114433^2}(32987388636176124480235619 (32987388636176124480235619 D_{21}\\
& & 
-(207561548775276615471 (v^2+v^{-2})) D_{22})+(2 12206056324453264983259133859103049535987217108) D_{22})\\
D_{24} & = & \frac{1}{2^4 7^2 1523^2 31068281^2 99593906868959^2}(185954429283234131709345979553 (557863287849702395128037938659 D_{22}\\
& & -(32987388636176124480235619 (v^2+v^{-2})) D_{23})-2739412859329752149612296806686633985128432534267734101 D_{23})\\
D_{25} & = & -\frac{1}{73^2 25284900439^2 21330668861813^2 4723^2}(202259863244780558381860453499585 (202259863244780558381860453499585 D_{23}\\
& & -(557863287849702395128037938659 (v^2+v^{-2})) D_{24})\\
& & -92984479634762534856210727781070979304376139845640812904182336 D_{24}),
\end{eqnarray*}
\begin{eqnarray*}
D_{26} & = & -\frac{1}{3^1 5^2 677^2 224699^2 10708193576771347^2 24833257^2}(66321637240913770895312913339793214\\
& &  (66321637240913770895312913339793214 D_{24}-(202259863244780558381860453499585 (v^2+v^{-2})) D_{25})\\
& & -14465919591689307128573983437824887265672029343238691854081861559749 D_{25})\\
D_{27} & = & \frac{1}{2^2 223^2 597482859539^2 248882825907331565531^2}\\
& & (2019335401386293786207166724075859341 (2019335401386293786207166724075859341 D_{25}-\\
& & (66321637240913770895312913339793214 (v^2+v^{-2})) D_{26})
\\
& & +297580655462944983601072219405031933105913989851705016725521711922510649 D_{26})\\
D_{28} & = & \frac{1}{569^2 1439^2 2466240348741368415391921411051^2}\\
& & (4412343134447284474117530515743496149 (61772803882261982637645427220408946086 D_{26}\\
& & -(2019335401386293786207166724075859341 (v^2+v^{-2})) D_{27})\\
& & +13337596486249616633416128365715003750495014302489045160099715949334689029 D_{27})\\
 D_{29} & = & \frac{1}{31^2 218069171^2 28721196055172805209^2 22725361^2 }\\
& &  (763361856645910935772663964451373952681 (763361856645910935772663964451373952681 D_{27}\\
& &  -(61772803882261982637645427220408946086 (v^2+v^{-2})) D_{28})\\
& &  -149724800808298139234895716448044923226445865165412568988595381125439091482374 D_{28})\\
 D_{30} & = & -\frac{1}{2^1 7^1 419^2 5598301^2 11550103^2 2551712030421623^2 3631^2 3041^2}\\
& &  (1126366500715763371124174330938368977431333 (1126366500715763371124174330938368977431333 D_{28}\\
& & -((v^2+v^{-2}) 763361856645910935772663964451373952681) D_{29})\\
& &  +156243595081359809872606121440348819964297847254205462786652267503931502270227271 D_{29})\\
 D_{31} & = & \frac{1}{2^2 19^2 326983565485534051144425519451185707^2 181301^2}\\
& &  (88867285707427504552815003873633980867964653 
 (88867285707427504552815003873633980867964653 D_{29}\\
& &  -(1126366500715763371124174330938368977431333 (v^2+v^{-2})) D_{30})\\
& & -123184813670701842459356619811841515374511032937050012802961098943874095248072007393959 D_{30})\\
 D_{32} & = & \frac{1}{2^2 29^2 449^2 1091^2 18414023^2 339722369745704256525298189901^2}\\
& &  (1713951505154451849278412459322972408517919889 
 (1713951505154451849278412459322972408517919889 D_{30}\\
& &  -(88867285707427504552815003873633980867964653 (v^2+v^{-2})) D_{31})\\
& &  +(2^1 213353666340007623638526192809579151292119580681379269874979045556633724290006464855951269)
  D_{31})\\
 D_{33} & = & -\frac{1}{43^2 97273570621^2 26396954918953349^2 130471771^2 4157^2 28621^2}\\
& &  (29676624191545618144648111624944621294859696166 
 (29676624191545618144648111624944621294859696166 D_{31}\\
 & & - ((v^2+v^{-2}) 1713951505154451849278412459322972408517919889) D_{32})\\
& &  +8239748870740490166734825221497528333002908462542913649313150812445359919137914031034615093 
 D_{32})\\
 D_{34} & = &\frac{1}{2^2 53^2 107073391^2 3080296628471719923564479600753^2 848857^2}\\
& &  (1202360930403560931030222837534182401486698364611 
 (1202360930403560931030222837534182401486698364611 D_{32}\\
 & & -(29676624191545618144648111624944621294859696166 
 (v^2+v^{-2})) D_{33})\\
& &  +101467727823064281002362989326473766812185962223175981805992928804837837463110204469450768695497
  D_{33})\\
 D_{35} & =&\frac{1}{3^2 107^2 576391299578767^2 274076093^2 290942059^2 13931^2 28657^2 204137^2}\\
& &  (50592899583888726481372641670551635013829487634578 
 (50592899583888726481372641670551635013829487634578 D_{33}\\
& &  -(1202360930403560931030222837534182401486698364611 (v^2+v^{-2})) D_{34})\\
& &  -18716563042659124012459968685711353673244482249337370871121262458352499400472107136941540852623
  D_{34})\\
 D_{36} & = & \frac{1}{2^4 29^2 1187^2 522300221197^2 7694404573578921876239719^2 182858701^2}\\
& &  (857145950896171400671963795556989894042194131289487 
 (857145950896171400671963795556989894042194131289487 D_{34}\\
& &  -(50592899583888726481372641670551635013829487634578 (v^2+v^{-2})) D_{35})\\
& & +4411923185656392383365159760236306142723974489779682676837077861812027157604126881190104182451178043
D_{35}),
\end{eqnarray*}
\begin{eqnarray*}
 D_{37} & = & \frac{1}{7^2 1311258103^2 93383157195708795083578685538920598221647^2}\\
& & ( 189712945006349928111821235829016019978456271676888 
 (189712945006349928111821235829016019978456271676888 D_{35}\\
 & & -(857145950896171400671963795556989894042194131289487 (v^2+v^{-2})) D_{36})\\
& &  -3888871381671794196458802315710490482078226712051900740103091728501684497594902265929767805231194888739
  D_{36})\\
 D_{38} & = & -\frac{1}{2^12 61^2 128478995994071^2 3025833381100035119025676342890881^2}\\
& &  (24773656778414901634803870348969176884565517861269127 
 (24773656778414901634803870348969176884565517861269127 D_{36}\\
& &  -(189712945006349928111821235829016019978456271676888 (v^2+v^{-2})) D_{37})\\
& &  +116014108817226025636250906637610633082682989655401422927845421603685431302028558553419600166912131824443
  D_{37})\\
 D_{39} & = & \frac{1}{3^2 23^2 1493^2 2830547^2 544917713^2 31598314303^2 24535887822332807^2 201101^2}\\
& & ( 1462774104358166447320328081145651681708102851750470
 (1462774104358166447320328081145651681708102851750470 D_{37}\\
 & & -(24773656778414901634803870348969176884565517861269127 (v^2+v^{-2})) D_{38})\\
& & -107114048131428196602036042282893363776335322672552405610021367908555000955590515740746720015779761857821 
D_{38})\\
D_{40} & = &  \frac{1}{2^2 5^2 472233067^2 335188641021288661^2 924126824062838315412481^2}\\
& & (447126463542130736420808986261322507372482025744507 
(447126463542130736420808986261322507372482025744507 D_{38}\\
& & -(1462774104358166447320328081145651681708102851750470 (v^2+v^{-2})) D_{39})\\
& &+1362236736392092735711595816371027107045738556864046551337719758325370857704119596126985689847376416801
 D_{39})\\
D_{41} & = &\frac{1}{3^2 8278065226744470551^2 4566184955538414377682528233^2 3943^2}\\
& & (7693200087216050876804839304036908158990939786094
 (7693200087216050876804839304036908158990939786094 D_{39}\\
& & -(447126463542130736420808986261322507372482025744507 (v^2+v^{-2})) D_{40})\\
& & +12380421501975427887574839221558916249368633135750618849551513659803018520966646313776000175920433209 
D_{40})\\
D_{42} & = & \frac{1}{2^4 23^2 28775297939^2 45354556404323956764088141663033^2 128147^2}\\
& & (444262371078411275036530891280549760044894125887
 (444262371078411275036530891280549760044894125887 D_{40}\\
& & -(7693200087216050876804839304036908158990939786094 (v^2+v^{-2})) D_{41})\\
& & +1359689975955312021939066413608950103421830908431132145007295475564434519858469029413502565524083 
D_{41})\\
D_{43} & = & \frac{1}{2^2 3^2 53^2 49609439671346607175679^2 56322002047749796817567^2}\\
& & (2598201402519308463837536678307826538637637361 
(12991007012596542319187683391539132693188186805 D_{41}\\
& & -(444262371078411275036530891280549760044894125887 (v^2+v^{-2})) D_{42})\\
 & &-1423853873843032305071517991227836398944647254463545764863856507677353245106560774424962103775
 D_{42})\\
D_{44} & = & \frac{1}{7 19^2 61^2 185233449385388602061^2 12102357082912717651939^2}\\
 & &(5724569388018125395560846055677725147411014
 (5724569388018125395560846055677725147411014 D_{42}\\
 & &-(12991007012596542319187683391539132693188186805 (v^2+v^{-2})) D_{43})\\
& & -(7^1 16858333738933659887870586331090479589286797477594035643445468458518447355054012844780395)
 D_{43})\\
D_{45} & = & \frac{1}{5^1 2^2 53^2 54005371585076654675102321279978539126519^2}\\
& & (173537434709794154672446152614045596469
 (173537434709794154672446152614045596469 D_{43}\\
& & -(5724569388018125395560846055677725147411014 (v^2+v^{-2})) D_{44})\\
& & +(5 1729710475616727676629256114810220610438295967049693913920779869327786123824237459) D_{44})\\
D_{46} & = & \frac{1}{7^1 11836998751^2 14660594155687780275955370219^2}\\
& & (472466105843823212817581554083188 (58585797124634078389380112706315312 D_{44}\\
& & -(173537434709794154672446152614045596469 (v^2+v^{-2})) D_{45})\\
& & -486773091634646981124420181550725259305095623574103826313324090433645617 D_{45})\\
D_{47} & = & \frac{1}{2^3 118116526460955803204395388520797^2}\\
& & (7735916375013449867448090407 (7735916375013449867448090407 D_{45}\\
& & -(14646449281158519597345028176578828 (v^2+v^{-2})) D_{46})\\
& & -297284847843573291454908916261083417941832630990370255741557001 D_{46}),
\end{eqnarray*}
\begin{eqnarray*}
D_{48} & = & \frac{1}{7^4 17^2 29^2 31^1 43 79^1 1013^1 36537661^2 368487522735331^2 19813^2 1747^2}\\
& & (2976484557725606944157223 (2976484557725606944157223 D_{46}
-(30943665500053799469792361628 (v^2+v^{-2})) D_{47})\\
& & +99497340316827556722995396776286275093114733864673142 D_{47})\\
& = & 1.
\end{eqnarray*}
\end{tiny}

\subsection{Computations on triality}\label{computationstriality}

\begin{tiny}
    \begin{eqnarray*}
    U_0 & = & T_{r_1}(-v-v^{-1}) T_{r_1+r_3}(-1-v^{-2}) T_{-r_1}(v) T_{r_3+r_4}(1) T_{-r_3-r_4}(1) T_{-r_3}(-v)  T_{r_4}(-v-v^{-1}) T_{-r_4}(v) T_{r_4}(-2v^{-3}) \\
     &  &T_{r_3+r_4}(v^{-2}) T_{r_1}(-v^{-3})\\
 \tau(\overline{U_0})& = &  \overline{T_{r_4}(-v-v^{-1})} ~\overline{T_{r_3+r_4}(1+v^{-2})}~\overline{ T_{-r_4}(v) }~\overline{T_{r_1+r_3}(-1)}~\overline{ T_{-r_1-r_3}(-1)}~\overline{ T_{-r_3}(-v)} ~ \overline{T_{r_1}(-v-v^{-1})}~\overline{ T_{-r_1}(v) }~\overline{T_{r_1}(-2v^{-3})}\\
     &  & \overline{T_{r_1+r_3}(-v^{-2})}~\overline{ T_{r_4}(-v^{-3})}\\
      & = & \overline{\tilde{U}_0}.
      \end{eqnarray*}
      \begin{eqnarray*}
      U_1 & = & T_{-r_1-r_3}(v^{-2}-1) T_{r_1+r_3}(-v^2+1) T_{r_1}(-\frac{v^4-v^2+1}{v^3}) T_{r_1+r_3+r_4}(v^{-1}-v-1) T_{-r_1-r_3-r_4}(2)  T_{r_1+r_3+r_4}(-1) \\
      & &T_{r_3+r_4}(v^{-2}-2) T_{-r_3-r_4}(-v^2) T_{r_3+r_4}(2) T_{-r_3-r_4}(-1) T_{-r_3}(\frac{3v^6-4v^4+3v^2-1}{v^5})  T_{r_4}(-\frac{6v^{10}-11v^8+15v^6-11v^4+4v^2-1}{v^4(v^5-v^3+v)})\\
      & & T_{-r_4}(-v^5+v^3-v) T_{r_4}(\frac{v^8-2v^6+5v^4-3v^2+1}{v^5(v^4-v^2+1)}) T_{r_3+r_4}(-\frac{2v^8-3v^6+6v^4-3v^2+1}{v^4}) T_{r_3}(-2v^5+v^3-v) \\
      & & T_{-r_1-r_3-r_4}(-2v^4-v^3+4v^2+v-2) T_{-r_1-r_3}(\frac{2v^8+v^7-6v^6-3v^5+10v^4+3v^3-6v^2-v+2}{v^5}) T_{-r_1}(v-2)\\
      \tau(\overline{U_1})& = & \overline{T_{-r_3-r_4}(1-v^{-2})}~\overline{ T_{r_3+r_4}(v^2-1)} ~\overline{T_{r_4}(-\frac{v^4-v^2+1}{v^3})} ~\overline{T_{r_1+r_3+r_4}(v^{-1}-v-1)}~ \overline{T_{-r_1-r_3-r_4}(2)} ~\overline{ T_{r_1+r_3+r_4}(-1)} \\
      & &\overline{T_{r_1+r_3}(2-v^{-2})}~\overline{ T_{-r_1-r_3}(v^2)}~\overline{ T_{r_1+r_3}(-2)}~\overline{ T_{-r_1-r_3}(1)}~\overline{ T_{-r_3}(\frac{3v^6-4v^4+3v^2-1}{v^5}) }~\overline{ T_{r_1}(-\frac{6v^{10}-11v^8+15v^6-11v^4+4v^2-1}{v^4(v^5-v^3+v)})}\\
      & & \overline{T_{-r_1}(-v^5+v^3-v)}~\overline{ T_{r_1}(\frac{v^8-2v^6+5v^4-3v^2+1}{v^5(v^4-v^2+1)})}~\overline{ T_{r_1+r_3}(\frac{2v^8-3v^6+6v^4-3v^2+1}{v^4})}~\overline{ T_{r_3}(-2v^5+v^3-v)} \\
      & & \overline{T_{-r_1-r_3-r_4}(-2v^4-v^3+4v^2+v-2)}~\overline{ T_{-r_3-r_4}(-\frac{2v^8+v^7-6v^6-3v^5+10v^4+3v^3-6v^2-v+2}{v^5})}~\overline{ T_{-r_4}(v-2)}\\
      & = &\overline{\tilde{U}_1}.
      \end{eqnarray*}
      \begin{eqnarray*}
      U_2 & = & T_{-r_1-r_3}(-\frac{v^6-v^4-1}{v^2(v^4-v^2+1)}) T_{r_1+r_3}(-v^4+v^2-1) T_{r_1}(-\frac{2v^8-3v^6+3v^4-2v^2+1}{v^5}) T_{r_1+r_3+r_4}(-\frac{(v^2-1)(v^4-v^2+1)}{v^3})\\
      & &  T_{r_3+r_4}(\frac{2v^6-2v^4+v^2-1}{v^2(v^4-v^2+1)}) T_{-r_3-r_4}(v^4-v^2+1) T_{-r_3}(-
  \frac{(5v^{14}-13v^{12}+19v^{10}-20v^8+15v^6-8v^4+3v^2-1}{v^7(v^4-v^2+1)})\\
  & &  T_{r_4}(-\frac{5v^{16}-18v^{14}+30v^{12}-34v^{10}+31v^8-21v^6+10v^4-3v^2+1)}{(v^4-v^2+1)^2v^7(v^2-1)}) T_{-r_4}(v^3(v^6-v^4+v^2-1))\\
  & &  T_{r_4}(\frac{v^8-2v^6+2v^4-3v^2+1}{v^7(v^2-1)})  T_{r_3+r_4}(-\frac{v^{12}-2v^{10}+2v^8-4v^6+3v^4-2v^2+1}{v^4(v^4-v^2+1}) T_{r_3}(-\frac{(v^6-v^4-1)v^3}{v^4-v^2+1})\\
  & & T_{-r_1-r_3-r_4}(-\frac{(v^6-v^4-1)v^3}{v^4-v^2+1}) T_{-r_1-r_3}(\frac{v^{12}-2v^{10}+v^8-3v^6+2v^4-2v^2+1}{v^4(v^4-v^2+1)}) T_{-r_1}(\frac{v^3(v^2+1)}{v^4-v^2+1})\\
  \tau(\overline{U_2}) & = & \overline{T_{-r_3-r_4}(\frac{v^6-v^4-1}{v^2(v^4-v^2+1)})}~ \overline{T_{r_3+r_4}(v^4-v^2+1)} ~\overline{T_{r_4}(-\frac{2v^8-3v^6+3v^4-2v^2+1}{v^5})}~ \overline{T_{r_1+r_3+r_4}(-\frac{(v^2-1)(v^4-v^2+1)}{v^3})}\\
      & &  \overline{T_{r_1+r_3}(-\frac{2v^6-2v^4+v^2-1}{v^2(v^4-v^2+1)})}~\overline{ T_{-r_1-r_3}(-v^4+v^2-1)} \overline{T_{-r_3}(-
  \frac{(5v^{14}-13v^{12}+19v^{10}-20v^8+15v^6-8v^4+3v^2-1}{v^7(v^4-v^2+1)})}\\
  & & \overline{ T_{r_1}(-\frac{5v^{16}-18v^{14}+30v^{12}-34v^{10}+31v^8-21v^6+10v^4-3v^2+1)}{(v^4-v^2+1)^2v^7(v^2-1)})} ~\overline{T_{-r_1}(v^3(v^6-v^4+v^2-1))}\\
  & &  \overline{T_{r_1}(\frac{v^8-2v^6+2v^4-3v^2+1}{v^7(v^2-1)})} ~\overline{ T_{r_1+r_3}(\frac{v^{12}-2v^{10}+2v^8-4v^6+3v^4-2v^2+1}{v^4(v^4-v^2+1)})}~\overline{ T_{r_3}(-\frac{(v^6-v^4-1)v^3}{v^4-v^2+1})}\\
  & & \overline{T_{-r_1-r_3-r_4}(-\frac{(v^6-v^4-1)v^3}{v^4-v^2+1})}~\overline{ T_{-r_3-r_4}(-\frac{v^{12}-2v^{10}+v^8-3v^6+2v^4-2v^2+1}{v^4(v^4-v^2+1)})}~\overline{ T_{-r_4}(\frac{v^3(v^2+1)}{v^4-v^2+1})}\\
  & = & \overline{\tilde{U}_2}.
      \end{eqnarray*}
      \begin{eqnarray*}
      U_3 & = & T_{-r_1-r_3}(-\frac{v^8-v^6+v^4+1}{v^2(v^6-v^4+v^2-1)}) T_{r_1}(\frac{2v^{12}-4v^{10}+5v^8-5v^6+4v^4-2v^2+1}{v^7}) T_{r_1+r_3}(v^6-v^4+v^2-1) \\
&&      T_{r_1+r_3+r_4}(\frac{(v^4-v^2+1)(v^6-v^4+v^2-1)}{v^5})T_{r_3+r_4}(\frac{2v^8-2v^6+2v^4-v^2+1}{v^2(v^6-v^4+v^2-1)}) \\
& & 
      T_{-r_3}(\frac{5v^{20}-15v^{18}+28v^{16}-39v^{14}+44v^{12}-40v^{10}+30v^8-18v^6+9v^4-3v^2+1}{v^9(v^6-v^4+v^2-1)}) T_{-r_3-r_4}(-v^6+v^4-v^2+1)\\
      & &        T_{r_4}(-\frac{5v^{28}-20v^{26}+48v^{24}-85v^{22}+125v^{20}-155v^{18}+167v^{16}-155v^{14}+126v^{12}
       -89v^{10}+54v^8-27v^6+12v^4-3v^2+1}{v^9(v^8-v^6+v^4-v^2+1)(v^6-v^4+v^2-1)^2})\\
       & & T_{-r_4}(v^{13}-v^{11}+v^9-v^7+v^5) T_{r_4}(\frac{v^{16}-2v^{14}+3v^{12}-4v^{10}+6v^8-5v^6+3v^4-3v^2+1}{v^9(v^8-v^6+v^4-v^2+1)})\\ 
       & & T_{r_3+r_4}(\frac{v^{16}-2v^{14}+3v^{12}-3v^{10}+5v^8-4v^6+3v^4-2v^2+1}{v^4(v^6-v^4+v^2-1)}) T_{r_3}(\frac{(v^8-v^6+v^4+1)v^5}{v^6-v^4+v^2-1}) T_{-r_1-r_3-r_4}(\frac{(v^8-v^6+v^4+1)v^5}{v^6-v^4+v^2-1})\\
       & &  T_{-r_1-r_3}(-\frac{v^{16}-2v^{14}+3v^{12}-2v^{10}+4v^8-3v^6+2v^4-2v^2+1}{v^4(v^6-v^4+v^2-1)}) T_{-r_1}(\frac{v^5(v^2+1)}{v^6-v^4+v^2-1})\\
       \tau(\overline{U_3}) & = & \overline{T_{-r_3-r_4}(\frac{v^8-v^6+v^4+1}{v^2(v^6-v^4+v^2-1)}) }~ \overline{T_{r_4}(\frac{2v^{12}-4v^{10}+5v^8-5v^6+4v^4-2v^2+1}{v^7})} ~ \overline{T_{r_3+r_4}(-v^6+v^4-v^2+1)} \\
&&     \overline{ T_{r_1+r_3+r_4}(\frac{(v^4-v^2+1)(v^6-v^4+v^2-1)}{v^5})}~ \overline{T_{r_1+r_3}(-\frac{2v^8-2v^6+2v^4-v^2+1}{v^2(v^6-v^4+v^2-1)}) }\\
& & 
     \overline{ T_{-r_3}(\frac{5v^{20}-15v^{18}+28v^{16}-39v^{14}+44v^{12}-40v^{10}+30v^8-18v^6+9v^4-3v^2+1}{v^9(v^6-v^4+v^2-1)})} \overline{ T_{-r_1-r_3}(v^6-v^4+v^2-1)}\\
      & &       \overline{ T_{r_1}(-\frac{5v^{28}-20v^{26}+48v^{24}-85v^{22}+125v^{20}-155v^{18}+167v^{16}-155v^{14}+126v^{12}
       -89v^{10}+54v^8-27v^6+12v^4-3v^2+1}{v^9(v^8-v^6+v^4-v^2+1)(v^6-v^4+v^2-1)^2})}\\
       & & \overline{T_{-r_1}(v^{13}-v^{11}+v^9-v^7+v^5)} \overline{T_{r_1}(\frac{v^{16}-2v^{14}+3v^{12}-4v^{10}+6v^8-5v^6+3v^4-3v^2+1}{v^9(v^8-v^6+v^4-v^2+1)})}\\ 
       & & \overline{T_{r_1+r_3}(-\frac{v^{16}-2v^{14}+3v^{12}-3v^{10}+5v^8-4v^6+3v^4-2v^2+1}{v^4(v^6-v^4+v^2-1)})} ~\overline{ T_{r_3}(\frac{(v^8-v^6+v^4+1)v^5}{v^6-v^4+v^2-1})} ~\overline{ T_{-r_1-r_3-r_4}(\frac{(v^8-v^6+v^4+1)v^5}{v^6-v^4+v^2-1})}\\
       & &  \overline{T_{-r_3-r_4}(\frac{v^{16}-2v^{14}+3v^{12}-2v^{10}+4v^8-3v^6+2v^4-2v^2+1}{v^4(v^6-v^4+v^2-1)})} ~\overline{T_{-r_4}(\frac{v^5(v^2+1)}{v^6-v^4+v^2-1})}\\
       & = & \overline{\tilde{U}_3}.
      \end{eqnarray*}
      \begin{eqnarray*}
      U_4 & = & T_{r_1+r_3+r_4}(v^{-2}) T_{-r_1-r_3-r_4}(-v^2-1) T_{-r_2}(v) T_{r_1+r_3+r_4}(1) T_{-r_4}(v^{-1}) T_{r_1+r_2+2r_3+r_4}(v) T_{-r_1-r_2-r_3-r_4}(v+v^{-1}) \\
      & & T_{r_3}(v^2+1) T_{-r_1}(v^{-1}) T_{r_2+r_3}(-v-v^{-1}) T_{-r_1-r_3-r_4}(-1-v^{-2})\\
      \tau(\overline{U_4})& = & \overline{T_{r_1+r_3+r_4}(v^{-2})}~\overline{ T_{-r_1-r_3-r_4}(-v^2-1)}~ \overline{T_{-r_2}(v)}~\overline{ T_{r_1+r_3+r_4}(1)}~ \overline{T_{-r_1}(v^{-1})}~ \overline{T_{r_1+r_2+2r_3+r_4}(v)}~ \overline{ T_{-r_1-r_2-r_3-r_4}(v+v^{-1})} \\
      & & \overline{T_{r_3}(v^2+1)} ~ \overline{T_{-r_4}(v^{-1})}~ \overline{T_{r_2+r_3}(-v-v^{-1})}~ \overline{T_{-r_1-r_3-r_4}(-1-v^{-2})}\\
      & = & \overline{\tilde{U}_4}.
      \end{eqnarray*}
      \begin{eqnarray*}
      U_5 & = & T_{-r_1-r_2-r_3-r_4}(-1) T_{r_1+r_2+r_3+r_4}(v^2+1) T_{r_1+r_2+r_3}(v^3+v) T_{-r_4}(-v) T_{r_1+r_2+2r_3+r_4}(v(v^2+2)) T_{-r_1-r_2-r_3-r_4}(-v^{-2}) \\
      & & T_{r_3}(-v-v^{-1}) T_{-r_1}(v^{-1}) T_{r_2+r_3+r_4}(v^3) T_{r_2}(v) T_{r_1+r_2+r_3+r_4}(v^4+v^2)\\
      \tau(\overline{U_5}) & = & \overline{T_{-r_1-r_2-r_3-r_4}(-1)}~\overline{ T_{r_1+r_2+r_3+r_4}(v^2+1)} ~\overline{ T_{r_1+r_2+r_3}(v^3+v)} ~ \overline{T_{-r_1}(-v)}~\overline{ T_{r_1+r_2+2r_3+r_4}(v(v^2+2))} ~ \overline{T_{-r_1-r_2-r_3-r_4}(-v^{-2})} \\
      & & \overline{T_{r_3}(-v-v^{-1})}~ \overline{T_{-r_4}(v^{-1})} ~ \overline{T_{r_1+r_2+r_3}(v^3)}~ \overline{T_{r_2}(v)} ~ \overline{T_{r_1+r_2+r_3+r_4}(v^4+v^2)}\\
      & = & \overline{\tilde{U}_5}.
    \end{eqnarray*}
    \end{tiny}

\section{$E_6$-graphs}\label{E6graphssection}

We give in this section the $E_6$-graphs $\widetilde{10_s}$ and $\widetilde{20_s}$ which verify the properties of Theorem \ref{bilinwgraphs}. We do not give the other new $E_6$-graphs and $E_8$-graphs which were obtained because of the large number of vertices and edges. They can be downloaded at \cite{newgraphsEsterle}.

\begin{center}

\end{center}

\section{$F_4$-graphs and field extensions}\label{sectionF4graphs}

\clearpage

\begin{figure}
\centering
%
\caption{$F_4$-graphs}\label{F4graphs1}
\end{figure}

\begin{figure}
\centering
%
\caption{$F_4$-graphs}\label{F4graphs2}
\end{figure}

\begin{figure}
\centering
\begin{tiny}
%
\end{tiny}
\caption{$F_4$-graphs}\label{F4graphs3}
\end{figure}

\begin{figure}
\centering
\begin{tiny}
%
\end{tiny}
\caption{Field extensions in type $F_4$}\label{fieldsF4}
\end{figure}

 \begin{figure}
\centering
\begin{tiny}
%
\end{tiny}
\caption{Field extensions in type $F_4$ in case $2$}\label{fieldsF4case2}
\end{figure}

\begin{figure}
\centering
\begin{tiny}
%
\end{tiny}
\caption{Field extensions in type $F_4$ in case $3$}\label{fieldsF4case3}
\end{figure}

\begin{figure}
\centering
\begin{tiny}
%
\end{tiny}
\caption{Field extensions in type $F_4$ in case $4$}\label{fieldsF4case4}
\end{figure}

\begin{figure}
\centering
\begin{tiny}
%
\end{tiny}
\caption{Field extensions in type $F_4$ in case $12$}\label{fieldsF4case12}
\end{figure}

\clearpage

\begin{figure}
\centering
\begin{tiny}
%
\end{tiny}
\caption{Field extensions in type $F_4$ in case $16$}\label{fieldsF4case16}
\end{figure}

\section{Connectedness of the $E_8$-graphs}

\subsection{$\widetilde{4480_y}$}\label{conec4480}
 If $p=7$ then the edges disappearing are $  493\rightarrow 50 $,  $570\rightarrow 26 $,  $598\rightarrow 27 $,  $1313\rightarrow 26 $,  $1313\rightarrow 293 $, 
          $ 1372\rightarrow 27 $, $ 1372\rightarrow 328 $, $ 1398\rightarrow 389 $, $ 1477\rightarrow 424 $, 
          $ 1495\rightarrow 419 $, $ 1526\rightarrow 441 $, $ 2026\rightarrow 357 $, $ 2102\rightarrow 371 $,  $ 2283\rightarrow 779 $, $ 2323\rightarrow 26 $, $ 2323\rightarrow 843 $, $ 2361\rightarrow 874 $, 
          $ 2416\rightarrow 646 $, $ 2416\rightarrow 891 $, $ 2421\rightarrow 27 $, $ 2421\rightarrow 886 $, 
          $ 2474\rightarrow 926 $, $ 2504\rightarrow 680 $, $ 2504\rightarrow 962 $, $ 2530\rightarrow 977 $, 
          $ 2618\rightarrow 996 $, $ 2667\rightarrow 646 $, $ 2688\rightarrow 680 $, $ 2734\rightarrow 691 $, $ 2933\rightarrow 57 $, $ 2933\rightarrow 119 $, 
          $ 3011\rightarrow 63 $, $ 3011\rightarrow 136 $,  $ 3475\rightarrow 27 $, $ 3485\rightarrow 646 $, $ 3485\rightarrow 1863 $, $ 3504\rightarrow 680 $, 
          $ 3504\rightarrow 1951 $, $ 3519\rightarrow 1977 $, $ 3555\rightarrow 27 $, $ 3555\rightarrow 2007 $, 
          $ 3559\rightarrow 27 $, $ 3590\rightarrow 2065 $, $ 3595\rightarrow 154 $, $ 3595\rightarrow 691 $, 
          $ 3595\rightarrow 2060 $, $ 3601\rightarrow 691 $, $ 3607\rightarrow 691 $, $ 3607\rightarrow 2120 $, 
          $ 3638\rightarrow 2158 $, $ 3702\rightarrow 2198 $, $ 3790\rightarrow 61 $, $ 3790\rightarrow 328 $, 
          $ 3790\rightarrow 424 $, $ 3790\rightarrow 874 $, $ 3790\rightarrow 880 $, $ 3790\rightarrow 886 $, 
          $ 3790\rightarrow 1747 $, $ 3801\rightarrow 47 $, $ 3801\rightarrow 50 $, $ 3801\rightarrow 977 $, 
          $ 3801\rightarrow 1793 $, $ 3801\rightarrow 1977 $, $ 3835\rightarrow 14 $, $ 3835\rightarrow 996 $, 
          $ 3835\rightarrow 1814 $, $ 3835\rightarrow 2065 $, $ 3864\rightarrow 424 $, $ 3870\rightarrow 441 $, $ 4040\rightarrow 26 $, $ 4040\rightarrow 27 $, $ 4040\rightarrow 611 $, $ 4040\rightarrow 2955 $,  $ 4057\rightarrow 617 $, 
          $ 4057\rightarrow 691 $, $ 4057\rightarrow 3004 $,   $ 4062\rightarrow 15 $,
          $ 4062\rightarrow 2986 $,  $ 4078\rightarrow 26 $, $ 4092\rightarrow 27 $,  $ 4092\rightarrow 3083 $, $ 4110\rightarrow 2379 $, $ 4124\rightarrow 2455 $,   $ 4153\rightarrow 691 $, $ 4153\rightarrow 3109 $, $ 4188\rightarrow 3168 $,  $ 4245\rightarrow 27 $, $ 4312\rightarrow 112 $, $ 4327\rightarrow 886 $, 
          $ 4345\rightarrow 57 $, $ 4345\rightarrow 1470 $, $ 4362\rightarrow 63 $, $ 4362\rightarrow 1548 $,  $ 4369\rightarrow 15 $, $ 4369\rightarrow 169 $,
          $ 4418\rightarrow 3 $, $ 4418\rightarrow 119 $, $ 4418\rightarrow 1470 $,  $ 4421\rightarrow 3 $, $ 4420\rightarrow 691 $, $ 4424\rightarrow 136 $, 
          $ 4424\rightarrow 1548 $,  $ 4431\rightarrow 27 $, 
          $ 4431\rightarrow 680 $, $ 4431\rightarrow 3988 $,  $ 4434\rightarrow 27 $, $ 4434\rightarrow 680 $, $ 4450\rightarrow 27 $, 
          $ 4454\rightarrow 31 $, $ 4454\rightarrow 47 $, $ 4454\rightarrow 50 $, $ 4454\rightarrow 236 $, 
          $ 4454\rightarrow 389 $, $ 4454\rightarrow 441 $, $ 4454\rightarrow 922 $, $ 4454\rightarrow 926 $, 
          $ 4454\rightarrow 1006 $, $ 4454\rightarrow 2060 $, $ 4454\rightarrow 3109 $, $ 4454\rightarrow 3883 $,  $ 4455\rightarrow 403 $, $ 4455\rightarrow 441 $, 
          $ 4455\rightarrow 2158 $,
          $ 4455\rightarrow 3168 $, $ 4455\rightarrow 3911 $, $ 4466\rightarrow 112 $, $ 4466\rightarrow 419 $, $ 4467\rightarrow 646 $,  $  4478\rightarrow 60 $ and $ 4478\rightarrow 63 $.
          
           They can be replaced by the paths $ 493\rightarrow 739\rightarrow 80\rightarrow 2323\rightarrow 515\rightarrow 50 $, $ 570\rightarrow 80\rightarrow 2323\rightarrow 2416\rightarrow 4040\rightarrow 2247\rightarrow 504\rightarrow 26 $, 
  $ 598\rightarrow 124\rightarrow 2421\rightarrow 2530\rightarrow 523\rightarrow 27 $, $ 1313\rightarrow 80\rightarrow 2323\rightarrow 2416\rightarrow 4040\rightarrow 2247\rightarrow 504\rightarrow 26 $, $ 1313\rightarrow 1398\rightarrow 1843\rightarrow 293 $, 
  $ 1372\rightarrow 124\rightarrow 2421\rightarrow 2530\rightarrow 523\rightarrow 27 $, $ 1372\rightarrow 3315\rightarrow 3790\rightarrow 3801\rightarrow 1668\rightarrow 328 $, $ 1398\rightarrow 178\rightarrow 293\rightarrow 389 $,
  $ 1477\rightarrow 1853\rightarrow 1976\rightarrow 293\rightarrow 1863\rightarrow 424 $, $ 1495\rightarrow 1875\rightarrow 2009\rightarrow 3666\rightarrow 2196\rightarrow 419 $,
  $ 1526\rightarrow 1939\rightarrow 2125\rightarrow 4378\rightarrow 3870\rightarrow 3913\rightarrow 2234\rightarrow 441 $, $ 2026\rightarrow 255\rightarrow 1888\rightarrow 357 $, 
  $ 2102\rightarrow 261\rightarrow 357\rightarrow 371 $, $ 2283\rightarrow 2361\rightarrow 3485\rightarrow 1465\rightarrow 1885\rightarrow 2858\rightarrow 874\rightarrow 1863\rightarrow 328\rightarrow 779 $, 
  $ 2323\rightarrow 2416\rightarrow 4040\rightarrow 2247\rightarrow 504\rightarrow 26 $, $ 2323\rightarrow 2416\rightarrow 4040\rightarrow 2247\rightarrow 504\rightarrow 26\rightarrow 90\rightarrow 843 $, 
  $ 2361\rightarrow 2474\rightarrow 4062\rightarrow 493\rightarrow 739\rightarrow 80\rightarrow 2323\rightarrow 515\rightarrow 50\rightarrow 874 $, $ 2416\rightarrow 4040\rightarrow 2247\rightarrow 504\rightarrow 26\rightarrow 646 $, 
  $ 2416\rightarrow 4040\rightarrow 2247\rightarrow 504\rightarrow 26\rightarrow 90\rightarrow 843\rightarrow 891 $, $ 2421\rightarrow 2530\rightarrow 523\rightarrow 27 $, 
  $ 2421\rightarrow 2530\rightarrow 523\rightarrow 27\rightarrow 112\rightarrow 886 $, $ 2474\rightarrow 2911\rightarrow 2980\rightarrow 4144\rightarrow 3119\rightarrow 926 $, 
  $ 2504\rightarrow 4057\rightarrow 2286\rightarrow 550\rightarrow 617\rightarrow 680 $, $ 2504\rightarrow 4057\rightarrow 2618\rightarrow 4188\rightarrow 2505\rightarrow 962 $, 
  $ 2530\rightarrow 2618\rightarrow 4188\rightarrow 2638\rightarrow 3083\rightarrow 3168\rightarrow 886\rightarrow 977 $, $ 2618\rightarrow 3000\rightarrow 3139\rightarrow 962\rightarrow 3004\rightarrow 996 $, 
  $ 2667\rightarrow 603\rightarrow 2424\rightarrow 646 $, $ 2688\rightarrow 2734\rightarrow 617\rightarrow 680 $, 
  $ 2734\rightarrow 2813\rightarrow 680\rightarrow 691 $, $ 2933\rightarrow 3011\rightarrow 56\rightarrow 57 $, $ 2933\rightarrow 3011\rightarrow 92\rightarrow 119 $, 
  $ 3011\rightarrow 3066\rightarrow 57\rightarrow 63 $, $ 3011\rightarrow 3066\rightarrow 1201\rightarrow 136 $, $ 3475\rightarrow 1499\rightarrow 197\rightarrow 27 $, 
  $ 3485\rightarrow 1465\rightarrow 2374\rightarrow 646 $, $ 3485\rightarrow 1465\rightarrow 1885\rightarrow 2858\rightarrow 874\rightarrow 1863 $, 
  $ 3504\rightarrow 1558\rightarrow 2491\rightarrow 680 $, $ 3504\rightarrow 1558\rightarrow 1941\rightarrow 2897\rightarrow 926\rightarrow 1951 $, 
  $ 3519\rightarrow 178\rightarrow 293\rightarrow 1863\rightarrow 424\rightarrow 1977 $, $ 3555\rightarrow 1584\rightarrow 2540\rightarrow 27 $, 
  $ 3555\rightarrow 1584\rightarrow 1990\rightarrow 2923\rightarrow 977\rightarrow 2007 $, $ 3559\rightarrow 2708\rightarrow 686\rightarrow 27 $, 
  $ 3590\rightarrow 3638\rightarrow 4391\rightarrow  1526\rightarrow 1939\rightarrow 2125\rightarrow 4378\rightarrow 3870\rightarrow 3913\rightarrow 2234\rightarrow 441\rightarrow 2065 $, $ 3595\rightarrow 4062\rightarrow 455\rightarrow 475\rightarrow 15\rightarrow 154 $, 
  $ 3595\rightarrow 4062\rightarrow 2570\rightarrow 691 $, $ 3595\rightarrow 4369\rightarrow 4454\rightarrow 3795\rightarrow 1951\rightarrow 2060 $, 
  $ 3601\rightarrow 2743\rightarrow 2806\rightarrow 691 $, $ 3607\rightarrow 1623\rightarrow 2596\rightarrow 691 $, 
  $ 3607\rightarrow 1623\rightarrow 2107\rightarrow 3016\rightarrow 996\rightarrow 2120 $, $ 3638\rightarrow 4391\rightarrow 4455\rightarrow 3850\rightarrow 3904\rightarrow 3168\rightarrow 3988\rightarrow 2158 $, 
  $ 3702\rightarrow 4153\rightarrow 2618\rightarrow 3000\rightarrow 3139\rightarrow 962\rightarrow 3004\rightarrow 996\rightarrow 2120\rightarrow 2198 $, $ 3790\rightarrow 1749\rightarrow 229\rightarrow 61 $, 
  $ 3790\rightarrow 3801\rightarrow 1668\rightarrow 328 $, $ 3790\rightarrow 3801\rightarrow 3864\rightarrow 3870\rightarrow 389\rightarrow 424 $, 
  $ 3790\rightarrow 1885\rightarrow 2858\rightarrow 874 $, $ 3790\rightarrow 1675\rightarrow 1738\rightarrow 880 $, 
  $ 3790\rightarrow 310\rightarrow 112\rightarrow 886 $, $ 3790\rightarrow 3801\rightarrow 1668\rightarrow 1747 $, $ 3801\rightarrow 285\rightarrow 743\rightarrow 47 $,
  $ 3801\rightarrow 285\rightarrow 743\rightarrow 50 $, $ 3801\rightarrow 1990\rightarrow 2923\rightarrow 977 $, 
  $ 3801\rightarrow 3864\rightarrow 1747\rightarrow 1793 $, $ 3801\rightarrow 3864\rightarrow 3870\rightarrow 389\rightarrow 424\rightarrow 1977 $, 
  $ 3835\rightarrow 1800\rightarrow 450\rightarrow 14 $, $ 3835\rightarrow 2107\rightarrow 3016\rightarrow 996 $, 
  $ 3835\rightarrow 2057\rightarrow 3878\rightarrow 1814 $, $ 3835\rightarrow 3870\rightarrow 3913\rightarrow 2234\rightarrow 441\rightarrow 2065 $, 
  $ 3864\rightarrow 3870\rightarrow 389\rightarrow 424 $, $ 3870\rightarrow 3913\rightarrow 2234\rightarrow 441 $, 
  $ 4040\rightarrow 2247\rightarrow 504\rightarrow 26 $, $ 4040\rightarrow 2281\rightarrow 523\rightarrow 27 $, $ 4040\rightarrow 2247\rightarrow 568\rightarrow 611 $, 
  $ 4040\rightarrow 2247\rightarrow 504\rightarrow 26\rightarrow 90\rightarrow 843\rightarrow 891\rightarrow 2955 $, $ 4057\rightarrow 2286\rightarrow 550\rightarrow 617 $, 
  $ 4057\rightarrow 2586\rightarrow 2734\rightarrow 2813\rightarrow 680\rightarrow 691 $, $ 4057\rightarrow 2618\rightarrow 4188\rightarrow 2505\rightarrow 2628\rightarrow 3004 $, 
  $ 4062\rightarrow 455\rightarrow 475\rightarrow 15 $, $ 4062\rightarrow 2285\rightarrow 815\rightarrow 2472\rightarrow 2606\rightarrow 2986 $, 
  $ 4078\rightarrow 93\rightarrow 2326\rightarrow 26 $, $ 4092\rightarrow 2281\rightarrow 523\rightarrow 27 $, $ 4092\rightarrow 4188\rightarrow 4303\rightarrow 3083 $,
  $ 4110\rightarrow 4124\rightarrow 4220\rightarrow 2379 $, $ 4124\rightarrow 2310\rightarrow 2379\rightarrow 2455 $, 
  $ 4153\rightarrow 2504\rightarrow 4057\rightarrow 2286\rightarrow 550\rightarrow 617\rightarrow 680\rightarrow 691 $, $ 4153\rightarrow 2504\rightarrow 4057\rightarrow 2286\rightarrow 550\rightarrow 617\rightarrow 680\rightarrow 691\rightarrow 1166\rightarrow 3109 $, 
  $ 4188\rightarrow 2638\rightarrow 3083\rightarrow 3168 $, $ 4245\rightarrow 2549\rightarrow 686\rightarrow 27 $, 
  $ 4312\rightarrow 2944\rightarrow 1137\rightarrow 112 $, $ 4327\rightarrow 4466\rightarrow 3197\rightarrow 1035\rightarrow 112\rightarrow 886 $, $ 4345\rightarrow 3280\rightarrow 56\rightarrow 57 $,
  $ 4345\rightarrow 3280\rightarrow 1415\rightarrow 1470 $, $ 4362\rightarrow 3351\rightarrow 1754\rightarrow 63 $, 
  $ 4362\rightarrow 3439\rightarrow 4398\rightarrow 1548 $, $ 4369\rightarrow 473\rightarrow 480\rightarrow 15 $, 
  $ 4369\rightarrow 2474\rightarrow 3504\rightarrow 169 $, $ 4418\rightarrow 2675\rightarrow 22\rightarrow 3 $, $ 4418\rightarrow 2662\rightarrow 92\rightarrow 119 $, 
  $ 4418\rightarrow 3477\rightarrow 357\rightarrow 1470 $, $ 4420\rightarrow 2813\rightarrow 680\rightarrow 691 $, $ 4421\rightarrow 139\rightarrow 142\rightarrow 3 $, 
  $ 4424\rightarrow 2701\rightarrow 622\rightarrow 136 $, $ 4424\rightarrow 3514\rightarrow 371\rightarrow 1548 $, 
  $ 4431\rightarrow 3485\rightarrow 1465\rightarrow 2374\rightarrow 646\rightarrow 27 $, $ 4431\rightarrow 3485\rightarrow 1465\rightarrow 2374\rightarrow 646\rightarrow 680 $, 
  $ 4431\rightarrow 3628\rightarrow 2158\rightarrow 4401\rightarrow 3168\rightarrow 3988 $, $ 4434\rightarrow 2678\rightarrow 646\rightarrow 27 $, 
  $ 4434\rightarrow 2678\rightarrow 646\rightarrow 680 $, $ 4450\rightarrow 2798\rightarrow 646\rightarrow 27 $, $ 4454\rightarrow 39\rightarrow 749\rightarrow 31 $, 
  $ 4454\rightarrow 285\rightarrow 743\rightarrow 47 $, $ 4454\rightarrow 285\rightarrow 743\rightarrow 50 $, $ 4454\rightarrow 1691\rightarrow 207\rightarrow 236 $, 
  $ 4454\rightarrow 3783\rightarrow 1692\rightarrow 389 $, $ 4454\rightarrow 3793\rightarrow 1763\rightarrow 441 $, 
  $ 4454\rightarrow 1941\rightarrow 2897\rightarrow 922 $, $ 4454\rightarrow 1941\rightarrow 2897\rightarrow 926 $, 
  $ 4454\rightarrow 3051\rightarrow 3133\rightarrow 1006 $, $ 4454\rightarrow 3795\rightarrow 1951\rightarrow 2060 $, 
  $ 4454\rightarrow 2071\rightarrow 3892\rightarrow 3109 $, $ 4454\rightarrow 3795\rightarrow 3817\rightarrow 3883 $, 
  $ 4455\rightarrow 3155\rightarrow 3986\rightarrow 403 $, $ 4455\rightarrow 3790\rightarrow 3801\rightarrow 3864\rightarrow 3870\rightarrow 389\rightarrow 424\rightarrow 441 $, 
  $ 4455\rightarrow 3850\rightarrow 3904\rightarrow 3168\rightarrow 3988\rightarrow 2158 $, $ 4455\rightarrow 3850\rightarrow 3904\rightarrow 3168 $, 
  $ 4455\rightarrow 3850\rightarrow 3904\rightarrow 3168\rightarrow 3988\rightarrow 2158\rightarrow 4401\rightarrow 3911 $, $ 4466\rightarrow 3197\rightarrow 1035\rightarrow 112 $, 
  $ 4466\rightarrow 3976\rightarrow 2228\rightarrow 419 $, $ 4467\rightarrow 103\rightarrow 2356\rightarrow 646 $, 
  $ 4478\rightarrow 2540\rightarrow 2923\rightarrow 60 $ and $ 4478\rightarrow 2506\rightarrow 57\rightarrow 63 $. This proves that this $E_8$-graph remains connected after specialization for $p=7$.
  
 If $p=11$ then the edges disappearing are $ 3790\rightarrow 57 $, $ 4418\rightarrow 26 $, $ 4418\rightarrow 27 $, $ 4424\rightarrow 691 $, $ 4425\rightarrow 27 $, 
  $ 4454\rightarrow 56 $, $ 4454\rightarrow 63 $, $ 4455\rightarrow 63 $, $ 4449\rightarrow 27 $ and $ 4454\rightarrow 32 $. They can be replaced by the paths $ 3790\rightarrow 351\rightarrow 800\rightarrow 57 $, $ 4418\rightarrow 2715\rightarrow 568\rightarrow 26 $, $ 4418\rightarrow 2675\rightarrow 646\rightarrow 27 $, 
  $ 4424\rightarrow 2701\rightarrow 680\rightarrow 691 $, $ 4425\rightarrow 2776\rightarrow 698\rightarrow 27 $, $ 4454\rightarrow 39\rightarrow 749\rightarrow 56 $, 
  $ 4454\rightarrow 369\rightarrow 820\rightarrow 63 $, $ 4455\rightarrow 369\rightarrow 820\rightarrow 63 $, $ 4449\rightarrow 2708\rightarrow 686\rightarrow 27 $ and
  $ 4454\rightarrow 1691\rightarrow 207\rightarrow 32 $. This proves that this $E_8$-graph remains connected after specialization for $p=11$.
  
  \subsection{$\widetilde{5670_y}$}\label{connec5670}
  
  If $p=7$ then the edges disappearing are $ 257\rightarrow 66 $, $ 1474\rightarrow 156 $, $ 1481\rightarrow 141 $, $ 1485\rightarrow 156 $, $ 1485\rightarrow 538 $, 
  $ 1540\rightarrow 243 $, $ 1545\rightarrow 244 $, $ 1556\rightarrow 191 $, $ 1565\rightarrow 1724 $, $ 1571\rightarrow 247 $, 
  $ 1610\rightarrow 171 $, $ 1610\rightarrow 280 $, $ 1634\rightarrow 171 $, $ 1686\rightarrow 191 $, $ 1724\rightarrow 1642 $, 
  $ 1766\rightarrow 171 $, $ 1790\rightarrow 473 $, $ 1923\rightarrow 801 $, $ 1968\rightarrow 505 $, $ 1984\rightarrow 822 $, 
  $ 2068\rightarrow 894 $, $ 2090\rightarrow 927 $, $ 2112\rightarrow 2227 $, $ 2133\rightarrow 2056 $, 
  $ 2133\rightarrow 2324 $, $ 2141\rightarrow 1014 $, $ 2166\rightarrow 473 $, $ 2195\rightarrow 244 $, 
  $ 2369\rightarrow 473 $, $ 2409\rightarrow 382 $, $ 2410\rightarrow 1835 $, $ 2477\rightarrow 425 $, $ 2765\rightarrow 505 $, 
  $ 2833\rightarrow 516 $, $ 2843\rightarrow 856 $, $ 2873\rightarrow 616 $, $ 2873\rightarrow 714 $, $ 2931\rightarrow 1020 $, 
  $ 2955\rightarrow 746 $, $ 2955\rightarrow 1019 $, $ 3040\rightarrow 631 $, $ 3063\rightarrow 1075 $, 
  $ 3347\rightarrow 3538 $, $ 3444\rightarrow 3559 $, $ 3479\rightarrow 818 $, $ 3497\rightarrow 1417 $, 
  $ 3570\rightarrow 793 $, $ 3591\rightarrow 818 $, $ 3615\rightarrow 3538 $, $ 3642\rightarrow 890 $, 
  $ 3642\rightarrow 1019 $, $ 3706\rightarrow 984 $, $ 3706\rightarrow 1103 $, $ 3836\rightarrow 3261 $, 
  $ 3947\rightarrow 4106 $, $ 3967\rightarrow 1083 $, $ 3978\rightarrow 1103 $, $ 4029\rightarrow 1192 $, 
  $ 4029\rightarrow 3947 $, $ 4210\rightarrow 141 $, $ 4210\rightarrow 191 $, $ 4254\rightarrow 191 $, 
  $ 4254\rightarrow 2174 $, $ 4479\rightarrow 1642 $, $ 4568\rightarrow 1693 $, $ 4568\rightarrow 1965 $, 
  $ 4588\rightarrow 1704 $, $ 4596\rightarrow 206 $, $ 4596\rightarrow 2608 $, $ 4651\rightarrow 2740 $, 
  $ 4652\rightarrow 2029 $, $ 4652\rightarrow 2716 $, $ 4657\rightarrow 3530 $, $ 4687\rightarrow 1965 $, 
  $ 4744\rightarrow 3581 $, $ 4777\rightarrow 3603 $, $ 4781\rightarrow 2029 $, $ 4815\rightarrow 2828 $, 
  $ 4849\rightarrow 3687 $, $ 4853\rightarrow 2080 $, $ 4853\rightarrow 2192 $, $ 4870\rightarrow 3748 $, 
  $ 4878\rightarrow 243 $, $ 4878\rightarrow 2101 $, $ 4925\rightarrow 243 $, $ 4925\rightarrow 2716 $, 
  $ 4957\rightarrow 2798 $, $ 5040\rightarrow 2631 $, $ 5055\rightarrow 2798 $, $ 5133\rightarrow 4186 $, 
  $ 5155\rightarrow 2838 $, $ 5166\rightarrow 2906 $, $ 5166\rightarrow 3703 $, $ 5198\rightarrow 3302 $, 
  $ 5198\rightarrow 3505 $, $ 5198\rightarrow 3881 $, $ 5246\rightarrow 3194 $, $ 5289\rightarrow 3262 $, 
  $ 5391\rightarrow 4061 $, $ 5424\rightarrow 4100 $, $ 5427\rightarrow 3476 $, $ 5427\rightarrow 4126 $, 
  $ 5428\rightarrow 746 $, $ 5428\rightarrow 793 $, $ 5428\rightarrow 4131 $, $ 5465\rightarrow 1075 $, 
  $ 5480\rightarrow 1417 $, $ 5480\rightarrow 1461 $, $ 5480\rightarrow 3985 $, $ 5480\rightarrow 4115 $, 
  $ 5500\rightarrow 3905 $, $ 5500\rightarrow 4037 $, $ 5500\rightarrow 4061 $, $ 5515\rightarrow 4186 $, 
  $ 5515\rightarrow 4197 $, $ 5530\rightarrow 1461 $, $ 5530\rightarrow 4190 $ and $ 5605\rightarrow 5414 $.
  
  They can be replaced by the paths $ 257\rightarrow 58\rightarrow 77\rightarrow 66 $, $ 1474\rightarrow 551\rightarrow 1585\rightarrow 156 $, $ 1481\rightarrow 1580\rightarrow 1479\rightarrow 141 $, 
  $ 1485\rightarrow 215\rightarrow 250\rightarrow 156 $, $ 1485\rightarrow 4204\rightarrow 1478\rightarrow 538 $, 
  $ 1540\rightarrow 2097\rightarrow 1534\rightarrow 243 $, $ 1545\rightarrow 1693\rightarrow 377\rightarrow 244 $, 
  $ 1556\rightarrow 1692\rightarrow 1538\rightarrow 191 $, $ 1565\rightarrow 1532\rightarrow 1563\rightarrow 1724 $, $  $, 
  $ 1610\rightarrow 335\rightarrow 4241\rightarrow 171 $, $ 1610\rightarrow 1760\rightarrow 440\rightarrow 280 $, 
  $ 1634\rightarrow 335\rightarrow 4241\rightarrow 171 $, $ 1686\rightarrow 399\rightarrow 4254\rightarrow 4379\rightarrow 4542\rightarrow 191 $, 
  $ 1724\rightarrow 1667\rightarrow 1636\rightarrow 1642 $, $ 1766\rightarrow 2332\rightarrow 279\rightarrow 171 $, 
  $ 1790\rightarrow 1921\rightarrow 485\rightarrow 473 $, $ 1923\rightarrow 2441\rightarrow 1302\rightarrow 801 $, 
  $ 1968\rightarrow 3220\rightarrow 1141\rightarrow 505 $, $ 1984\rightarrow 1959\rightarrow 807\rightarrow 822 $, 
  $ 2068\rightarrow 870\rightarrow 562\rightarrow 894 $, $ 2090\rightarrow 4904\rightarrow 2076\rightarrow 927 $, 
  $ 2112\rightarrow 421\rightarrow 4346\rightarrow 2227 $, $ 2133\rightarrow 2197\rightarrow 2085\rightarrow 2056 $, 
  $ 2133\rightarrow 2197\rightarrow 2411\rightarrow 2324 $, $ 2141\rightarrow 2313\rightarrow 1219\rightarrow 1014 $, 
  $ 2166\rightarrow 4878\rightarrow 2537\rightarrow 473 $, $ 2195\rightarrow 2409\rightarrow 377\rightarrow 244 $, 
  $ 2369\rightarrow 2534\rightarrow 5546\rightarrow 473 $, $ 2409\rightarrow 2477\rightarrow 409\rightarrow 382 $, 
  $ 2410\rightarrow 2470\rightarrow 2395\rightarrow 1835 $, $ 2477\rightarrow 405\rightarrow 962\rightarrow 425 $, 
  $ 2765\rightarrow 3967\rightarrow 1141\rightarrow 505 $, $ 2833\rightarrow 4097\rightarrow 1179\rightarrow 516 $, 
  $ 2843\rightarrow 850\rightarrow 587\rightarrow 856 $, $ 2873\rightarrow 5189\rightarrow 3328\rightarrow 616 $, 
  $ 2873\rightarrow 5189\rightarrow 3524\rightarrow 714 $, $ 2931\rightarrow 3376\rightarrow 1223\rightarrow 1020 $, 
  $ 2955\rightarrow 3059\rightarrow 2975\rightarrow 746 $, $ 2955\rightarrow 3059\rightarrow 2975\rightarrow 1019 $, 
  $ 3040\rightarrow 619\rightarrow 584\rightarrow 631 $, $ 3063\rightarrow 3161\rightarrow 3082\rightarrow 1075 $, 
  $ 3347\rightarrow 3260\rightarrow 3474\rightarrow 3538 $, $ 3444\rightarrow 1325\rightarrow 5250\rightarrow 3559 $, 
  $ 3479\rightarrow 5447\rightarrow 3620\rightarrow 818 $, $ 3497\rightarrow 3435\rightarrow 1292\rightarrow 1417 $, 
  $ 3570\rightarrow 3561\rightarrow 765\rightarrow 793 $, $ 3591\rightarrow 3606\rightarrow 3630\rightarrow 818 $, 
  $ 3615\rightarrow 3630\rightarrow 3616\rightarrow 3538 $, $ 3642\rightarrow 3700\rightarrow 916\rightarrow 890 $, 
  $ 3642\rightarrow 3700\rightarrow 3644\rightarrow 1019 $, $ 3706\rightarrow 959\rightarrow 409\rightarrow 984 $, 
  $ 3706\rightarrow 1102\rightarrow 496\rightarrow 1103 $, $ 3836\rightarrow 2605\rightarrow 2020\rightarrow 3261 $, 
  $ 3947\rightarrow 3767\rightarrow 3943\rightarrow 4106 $, $ 3967\rightarrow 1031\rightarrow 1241\rightarrow 1083 $, 
  $ 3978\rightarrow 4126\rightarrow 1288\rightarrow 1103 $, $ 4029\rightarrow 4035\rightarrow 4004\rightarrow 1192 $, 
  $ 4029\rightarrow 4035\rightarrow 4004\rightarrow 3947 $, $ 4210\rightarrow 4222\rightarrow 1568\rightarrow 141 $, 
  $ 4210\rightarrow 4220\rightarrow 4217\rightarrow 191 $, $ 4254\rightarrow 4379\rightarrow 4542\rightarrow 191 $, 
  $ 4254\rightarrow 4379\rightarrow 2236\rightarrow 2174 $, $ 4479\rightarrow 5157\rightarrow 4501\rightarrow 1642 $, 
  $ 4568\rightarrow 5177\rightarrow 2872\rightarrow 1693 $, $ 4568\rightarrow 4588\rightarrow 1994\rightarrow 1965 $, 
  $ 4588\rightarrow 4430\rightarrow 4640\rightarrow 1704 $, $ 4596\rightarrow 4437\rightarrow 4624\rightarrow 206 $, 
  $ 4596\rightarrow 4437\rightarrow 2283\rightarrow 2608 $, $ 4651\rightarrow 4448\rightarrow 2295\rightarrow 2740 $, 
  $ 4652\rightarrow 5231\rightarrow 4667\rightarrow 2029 $, $ 4652\rightarrow 5231\rightarrow 4667\rightarrow 2716 $, 
  $ 4657\rightarrow 4452\rightarrow 3358\rightarrow 3530 $, $ 4687\rightarrow 5035\rightarrow 4712\rightarrow 1965 $, 
  $ 4744\rightarrow 5332\rightarrow 4745\rightarrow 3581 $, $ 4777\rightarrow 5109\rightarrow 4801\rightarrow 3603 $, 
  $ 4781\rightarrow 4755\rightarrow 1971\rightarrow 2029 $, $ 4815\rightarrow 5084\rightarrow 4821\rightarrow 2828 $, 
  $ 4849\rightarrow 4864\rightarrow 3712\rightarrow 3687 $, $ 4853\rightarrow 4859\rightarrow 2117\rightarrow 2080 $, 
  $ 4853\rightarrow 4859\rightarrow 2226\rightarrow 2192 $, $ 4870\rightarrow 4369\rightarrow 3230\rightarrow 3748 $, 
  $ 4878\rightarrow 4906\rightarrow 2290\rightarrow 243 $, $ 4878\rightarrow 4906\rightarrow 2110\rightarrow 2101 $, 
  $ 4925\rightarrow 4988\rightarrow 2147\rightarrow 243 $, $ 4925\rightarrow 4986\rightarrow 4934\rightarrow 2716 $, 
  $ 4957\rightarrow 3980\rightarrow 5055\rightarrow 2343\rightarrow 482\rightarrow 2798 $, $ 5040\rightarrow 4755\rightarrow 1994\rightarrow 2631 $, 
  $ 5055\rightarrow 2343\rightarrow 482\rightarrow 2798 $, $ 5133\rightarrow 5514\rightarrow 5146\rightarrow 4186 $, 
  $ 5155\rightarrow 4492\rightarrow 1574\rightarrow 2838 $, $ 5166\rightarrow 5217\rightarrow 5168\rightarrow 2906 $, 
  $ 5166\rightarrow 5217\rightarrow 5168\rightarrow 3703 $, $ 5198\rightarrow 3134\rightarrow 5223\rightarrow 3302 $, 
  $ 5198\rightarrow 3134\rightarrow 793\rightarrow 3505 $, $ 5198\rightarrow 3134\rightarrow 1103\rightarrow 3881 $, 
  $ 5246\rightarrow 3063\rightarrow 5177\rightarrow 3194 $, $ 5289\rightarrow 3161\rightarrow 5203\rightarrow 3262 $, 
  $ 5391\rightarrow 3297\rightarrow 1185\rightarrow 4061 $, $  $, $ 5427\rightarrow 5294\rightarrow 3262\rightarrow 3476 $, 
  $ 5427\rightarrow 5442\rightarrow 4165\rightarrow 4126 $, $ 5428\rightarrow 4137\rightarrow 3574\rightarrow 746 $, 
  $ 5428\rightarrow 3479\rightarrow 5447\rightarrow 793 $, $ 5428\rightarrow 4137\rightarrow 3574\rightarrow 4131 $, 
  $ 5465\rightarrow 5035\rightarrow 473\rightarrow 1075 $, $ 5480\rightarrow 3873\rightarrow 2686\rightarrow 1417 $, 
  $ 5480\rightarrow 4007\rightarrow 2772\rightarrow 1461 $, $ 5480\rightarrow 4015\rightarrow 3837\rightarrow 3985 $, 
  $ 5480\rightarrow 4133\rightarrow 3578\rightarrow 4115 $, $ 5500\rightarrow 5116\rightarrow 2745\rightarrow 3905 $, 
  $ 5500\rightarrow 5515\rightarrow 4086\rightarrow 4037 $, $ 5500\rightarrow 4074\rightarrow 4166\rightarrow 4061 $, 
  $ 5515\rightarrow 5421\rightarrow 5456\rightarrow 4186 $, $ 5515\rightarrow 4199\rightarrow 4150\rightarrow 4197 $, 
  $ 5530\rightarrow 4103\rightarrow 1449\rightarrow 1461 $, $ 5530\rightarrow 5527\rightarrow 4097\rightarrow 4190 $ and
  $ 5605\rightarrow 5594\rightarrow 5613\rightarrow 5414 $. This proves this $E_8$-graph remains connected after specialization when $p=7$.
  
  \subsection{$\widetilde{7168}$}\label{connec7168}
  For $p=7$, the edges disappearing are $ 3397\rightarrow 403 $, $ 3995\rightarrow 4 $, $ 4205\rightarrow 772 $, $ 4795\rightarrow 86 $, $ 4795\rightarrow 157 $, 
  $ 4993\rightarrow 4 $, $ 5755\rightarrow 157 $, $ 5755\rightarrow 300 $, $ 5938\rightarrow 4 $, $ 6397\rightarrow 2964 $, 
  $ 6498\rightarrow 157 $, $ 6538\rightarrow 181 $, $ 6766\rightarrow 3772 $, $ 6869\rightarrow 4 $, $ 6869\rightarrow 181 $, 
  $ 6869\rightarrow 1414 $, $ 6988\rightarrow 86 $, $ 6988\rightarrow 300 $, $ 6988\rightarrow 631 $, $ 7012\rightarrow 4 $, 
  $ 7012\rightarrow 671 $, $ 7012\rightarrow 1414 $, $ 7012\rightarrow 2374 $, $ 7071\rightarrow 4 $, $ 7083\rightarrow 4 $, 
  $ 7083\rightarrow 181 $, $ 7083\rightarrow 2374 $, $ 7165\rightarrow 86 $, $ 7165\rightarrow 98 $, $ 7165\rightarrow 157 $, 
  $ 7165\rightarrow 300 $, $ 7165\rightarrow 1231 $, $ 7165\rightarrow 2176 $ and $ 7165\rightarrow 3174 $.
  
  They can be replaced by the paths $ 3397\rightarrow 4795\rightarrow 1156\rightarrow 1362\rightarrow 86\rightarrow 403 $, $ 3995\rightarrow 30\rightarrow 106\rightarrow 4 $, $ 4205\rightarrow 4795\rightarrow 4993\rightarrow 48\rightarrow 106\rightarrow 4\rightarrow 157\rightarrow 772 $, 
  $ 4795\rightarrow 1156\rightarrow 1362\rightarrow 86 $, $ 4795\rightarrow 4993\rightarrow 48\rightarrow 106\rightarrow 4\rightarrow 157 $, $ 4993\rightarrow 48\rightarrow 106\rightarrow 4 $, 
  $ 5755\rightarrow 2109\rightarrow 453\rightarrow 157 $, $ 5755\rightarrow 470\rightarrow 2264\rightarrow 300 $, $ 5938\rightarrow 30\rightarrow 106\rightarrow 4 $, 
  $ 6397\rightarrow 6869\rightarrow 4877\rightarrow 5085\rightarrow 1414\rightarrow 2964 $, $ 6498\rightarrow 3593\rightarrow 916\rightarrow 157 $, 
  $ 6538\rightarrow 3627\rightarrow 972\rightarrow 181 $, $ 6766\rightarrow 6869\rightarrow 4592\rightarrow 1234\rightarrow 181\rightarrow 3772 $, $ 6869\rightarrow 48\rightarrow 106\rightarrow 4 $, 
  $ 6869\rightarrow 4592\rightarrow 1234\rightarrow 181 $, $ 6869\rightarrow 4877\rightarrow 5085\rightarrow 1414 $, 
  $ 6988\rightarrow 2009\rightarrow 344\rightarrow 86 $, $ 6988\rightarrow 3719\rightarrow 2278\rightarrow 300 $, 
  $ 6988\rightarrow 4969\rightarrow 3533\rightarrow 631 $, $ 7012\rightarrow 18\rightarrow 49\rightarrow 4 $, $ 7012\rightarrow 1439\rightarrow 3015\rightarrow 671 $, 
  $ 7012\rightarrow 5264\rightarrow 3772\rightarrow 1414 $, $ 7012\rightarrow 6766\rightarrow 7083\rightarrow 5796\rightarrow 2176\rightarrow 2374 $, $ 7071\rightarrow 30\rightarrow 106\rightarrow 4 $, 
  $ 7083\rightarrow 30\rightarrow 106\rightarrow 4 $, $ 7083\rightarrow 2166\rightarrow 510\rightarrow 181 $, $ 7083\rightarrow 5796\rightarrow 2176\rightarrow 2374 $, 
  $ 7165\rightarrow 2009\rightarrow 344\rightarrow 86 $, $ 7165\rightarrow 81\rightarrow 91\rightarrow 98 $, $ 7165\rightarrow 418\rightarrow 1148\rightarrow 157 $, $ 7165\rightarrow 2522\rightarrow 2577\rightarrow 300 $, 
  $ 7165\rightarrow 2022\rightarrow 2667\rightarrow 1231 $, $ 7165\rightarrow 2022\rightarrow 854\rightarrow 2176 $ and
  $ 7165\rightarrow 5406\rightarrow 5623\rightarrow 3174 $. This proves this $E_8$-graph remains connected after specialization when $p=7$.
  
  \section{Matrices for the $12$-dimensional representation of $\mathcal{H}_{F_4,\alpha,\beta}$}\label{sectionmatricesF412}
  
  $X_1$, $X_2$ and $X_3$  are respectively
\begin{tiny} $$\setcounter{MaxMatrixCols}{20} \begin{pmatrix}
  0&0&0&0\\
 0&0&0&0\\
 0&0&0&0\\
 1&0&-\frac {\alpha+1}{3\sqrt{\alpha}}&0\\
 \frac {\alpha+\beta}{\alpha\beta+1}&0&-\frac {\alpha^2+\alpha\beta+\alpha+\beta}{ 3\left( \alpha\beta+1 \right) \sqrt{\alpha}}&\frac {\alpha\beta+\alpha+\beta+1}{3(\alpha\beta+1)}\\
 \frac {\alpha \left( \alpha\beta^2+2\,\alpha\beta+\beta^2+\alpha+2\,\beta+1 \right) }{\alpha^3\beta^2+\alpha^2\beta+\alpha\beta+1}&\frac {\alpha \left( \beta+1 \right) }{\alpha^2\beta+1}&-\frac { \left( \alpha^2+\alpha\beta+\alpha+\beta \right) \sqrt{\alpha} \left( \beta+1 \right) }{ 3\left( \alpha\beta+1 \right)  \left( \alpha^2\beta+1 \right) }&\frac { \left( \alpha\beta+\alpha+\beta+1 \right) \alpha \left( \beta+1 \right) }{ 3\left( \alpha\beta+1 \right)  \left( \alpha^2\beta+1 \right) }\\
 \alpha\beta+1&0&0&0\\
 -\alpha-\beta&0&0&0\\
 0&\alpha+\beta&0&0\\
 0&0&0&0\\
 0&0&0&0\\
 0&0&0&0
         \end{pmatrix}
 $$ 
  $$\setcounter{MaxMatrixCols}{20} \begin{pmatrix}
  0&0& 0&0\\
  0&\sqrt{\alpha}\beta+\sqrt{\alpha}& 0&0 \\
  0&-\frac {\alpha^2\beta+1}{\sqrt{\alpha}} &0&0 \\
  0&-\frac { \left( 2\,\alpha^2\beta-\alpha^2+\alpha\beta+\alpha-\beta+2 \right) \sqrt{\beta}}{3(\alpha^2\beta^2+\alpha^2\beta+\beta+1)} & 0&-\frac{ \sqrt{\beta}}{\beta+1} \\
 -\frac {\sqrt{\alpha} \left( \beta+1 \right) }{\alpha\beta+1} &-\frac {\sqrt{\beta} \left( \alpha\beta+4\,\alpha+4\,\beta+1 \right) }{3(\alpha\beta^2+\alpha\beta+\beta+1)}& 0&-\frac { \left( \alpha+\beta \right) \sqrt{\beta}}{ \left( \alpha\beta+1 \right)  \left( \beta+1 \right) }\\
 -\frac { \left( \beta+1 \right) ^{2}\sqrt{\alpha}\alpha}{ \left( \alpha\beta+1 \right)  \left( \alpha^2\beta+1 \right) }& -\frac { \left( \alpha\beta+4\,\alpha+4\,\beta+1 \right) \alpha\sqrt{\beta}}{3 \left( \alpha\beta+1 \right)  \left( \alpha^2\beta+1 \right) }&-\frac {\sqrt{\alpha} \left( \alpha\beta+\alpha+\beta+1 \right) }{3(\alpha^2\beta+1)}& -\frac { \left( \alpha\beta+\alpha+\beta+1 \right) \alpha\sqrt{\beta}}{\alpha^3\beta^2+\alpha^2\beta+\alpha\beta+1}\\
 0&0 & 0&\frac {\alpha\sqrt{\beta} \left( \alpha\beta+1 \right) }{\alpha^2+\beta}\\
 \sqrt{\alpha} \left( \beta+1 \right) &-\frac { \left( \alpha+\beta \right) \alpha\sqrt{\beta}}{\alpha^2+\beta} &0&-\frac { \left( \alpha+\beta \right) \alpha\sqrt{\beta}}{\alpha^2+\beta}\\
 -\frac {\alpha^2+\beta}{\sqrt{\alpha}}&\frac { \left( \alpha+\beta \right) \sqrt{\beta}}{\beta+1}&0&0\\
 0&1 &0&1\\
 0&\frac {\alpha^3\beta+\alpha\beta^2+\alpha^2+\beta}{\alpha^3+2\alpha^2\beta+\alpha\beta^2+\alpha^2+2\alpha\beta+\beta^2} & 0&\frac {\alpha\beta+1}{\alpha+\beta}\\ 
 0&\frac {\alpha \left( \beta+1 \right)  \left( \alpha\beta+1 \right) }{\alpha^3+2\,\alpha^2\beta+\alpha\beta^2+\alpha^2+2\,\alpha\beta+\beta^2} & 0&0
         \end{pmatrix}
 $$ 
$$
\begin{pmatrix}
\alpha+\beta&-\sqrt{\alpha}\beta-\sqrt{\alpha}&0&0\\
-\alpha\beta-1&0&0&0\\ 
0&0&0&0\\
\frac {\sqrt{\beta}\sqrt{\alpha} \left( 2\,\alpha\beta-\alpha-\beta+2 \right) }{3(\alpha^2\beta^2+\alpha^2\beta+\beta+1)}&0&\frac {\sqrt{\beta} \left( \alpha+1 \right) }{\sqrt{\alpha} \left( \beta+1 \right) }&0\\ 
\frac { \left( \alpha+\beta \right) \sqrt{\beta}\sqrt{\alpha} \left( 2\,\alpha\beta-\alpha-\beta+2 \right) }{3 \left( \alpha^2\beta^2+\alpha^2\beta+\beta+1 \right)  \left( \alpha\beta+1 \right) }&-\frac {\alpha\sqrt{\beta} \left( 2\,\alpha\beta-\alpha-\beta+2 \right) }{ 3\left( \alpha\beta+1 \right)  \left( \alpha^2\beta+1 \right) }&\frac { \left( \alpha^2+\alpha\beta+\alpha+\beta \right) \sqrt{\beta}}{\sqrt{\alpha} \left( \beta+1 \right)  \left( \alpha\beta+1 \right) }&\frac {\sqrt{\beta}\sqrt{\alpha}}{\alpha\beta+1}\\ 
\frac {\sqrt{\beta}\sqrt{\alpha} \left( \alpha^2\beta+\alpha^2+2\,\alpha\beta+2\,\alpha+\beta+1 \right) }{3(\alpha^3\beta^2+\alpha^2\beta+\alpha\beta+1)}&-\frac { \left( \alpha\beta+\alpha+\beta+1 \right) \alpha\sqrt{\beta}}{3(\alpha^3\beta^2+\alpha^2\beta+\alpha\beta+1)}&\frac { \left( \alpha^2+\alpha\beta+\alpha+\beta \right) \sqrt{\beta}\sqrt{\alpha}}{ \left( \alpha\beta+1 \right)  \left( \alpha^2\beta+1 \right) }&\frac { \left( \beta+1 \right) \sqrt{\beta}\sqrt{\alpha}\alpha}{ \left( \alpha\beta+1 \right)  \left( \alpha^2\beta+1 \right) }\\ 
-\frac { \left( \alpha+\beta \right) \sqrt{\beta}\sqrt{\alpha} \left( \alpha\beta+1 \right) }{ \left( \beta+1 \right)  \left( \alpha^2+\beta \right) }&0&0&0\\ 
\frac {\sqrt{\beta}\sqrt{\alpha} \left( \alpha+\beta \right) ^{2}}{ \left( \beta+1 \right)  \left( \alpha^2+\beta \right) }&-\frac { \left( \alpha+\beta \right) \alpha\sqrt{\beta}}{\alpha^2+\beta}&0&\frac { \left( \beta+1 \right) \sqrt{\beta}\sqrt{\alpha}\alpha}{\alpha^2+\beta}\\ 
0&0&0&-\sqrt{\beta}\sqrt{\alpha}\\ 
-\frac {\sqrt{\alpha}}{\alpha+1}&\frac {\alpha \left( \beta+1 \right) }{\alpha^2+\alpha\beta+\alpha+\beta}&-\frac {\alpha^2+\beta}{ \left( \alpha+\beta \right) \sqrt{\alpha}}&-\frac {\sqrt{\alpha} \left( \beta+1 \right) }{\alpha+\beta}\\ 
-\frac { \left( \alpha\beta+1 \right) \sqrt{\alpha}}{\alpha^2+\alpha\beta+\alpha+\beta}&0&-\frac { \left( \alpha\beta+1 \right)  \left( \alpha^2+\beta \right) }{ \left( \alpha+\beta \right) ^{2}\sqrt{\alpha}}&0\\
0&0&-\frac {\sqrt{\alpha} \left( \beta+1 \right)  \left( \alpha\beta+1 \right) }{ \left( \alpha+\beta \right) ^{2}}&0
\end{pmatrix}
 $$
 \end{tiny}
and $Y_1$, $Y_2$, $Y_3$ and $Y_4$ are respectively

\begin{tiny}
$$
\begin{pmatrix}
\frac {\sqrt{\alpha}}{\alpha+1}&\frac { \left( \alpha+\beta \right) \sqrt{\alpha}}{ \left( \alpha\beta+1 \right)  \left( \alpha+1 \right) }\\
 0&\frac { \left( \alpha^3+2\,\alpha^2\beta+\alpha\beta^2+\alpha^2+2\,\alpha\beta+\beta^2 \right) \sqrt{\alpha}}{ \left( \alpha^2\beta^2+\alpha^2\beta+\alpha\beta+\alpha \right)  \left( \alpha+1 \right) }\\
 \frac { \left( \alpha^3+\alpha^2\beta+2\,\alpha^2+2\,\alpha\beta+\alpha+\beta \right) \sqrt{\alpha}}{ \left( \sqrt{\alpha}^{5}+\sqrt{\alpha}\alpha\beta+\sqrt{\alpha}\alpha+\sqrt{\alpha}\beta \right)  \left( \alpha+1 \right) }&\frac { \left( \alpha+\beta \right)  \left( \alpha^3+\alpha^2\beta+2\,\alpha^2+2\,\alpha\beta+\alpha+\beta \right) \sqrt{\alpha}}{ \left( \alpha\beta+1 \right)  \left( \sqrt{\alpha}^{5}+\sqrt{\alpha}\alpha\beta+\sqrt{\alpha}\alpha+\sqrt{\alpha}\beta \right)  \left( \alpha+1 \right) }\\
 0&\frac { \left( \alpha+\beta \right)  \left( \alpha^3+\alpha^2\beta+2\,\alpha^2+2\,\alpha\beta+\alpha+\beta \right) }{ \left( \alpha\beta+1 \right)  \left( \sqrt{\alpha}^{5}+\sqrt{\alpha}\alpha\beta+\sqrt{\alpha}\alpha+\sqrt{\alpha}\beta \right) }\\
 0&\frac {\alpha+\beta}{\alpha\beta+1}\\
 0&0\\
 \frac { \left( 2\,\alpha^2\beta+2\,\alpha\beta^2-\alpha^2-2\,\alpha\beta-\beta^2+2\,\alpha+2\,\beta \right) \sqrt{\alpha}}{ \left( \sqrt{\alpha}\alpha\beta^2+\sqrt{\alpha} \right)  \left( \alpha+1 \right) }&\frac { \left( \alpha^3\beta+\alpha^2\beta^2+4\,\alpha^3+9\,\alpha^2\beta+5\,\alpha\beta^2+5\,\alpha^2+9\,\alpha\beta+4\,\beta^2+\alpha+\beta \right) \sqrt{\alpha}}{ \left( \sqrt{\alpha}^{5}\beta^2+\sqrt{\alpha}^{5}\beta+\sqrt{\alpha}\alpha\beta^2+2\,\sqrt{\alpha}\alpha\beta+\sqrt{\alpha}\alpha+\sqrt{\alpha}\beta+\sqrt{\alpha} \right)  \left( \alpha+1 \right) }\\
 -\frac { \left( \alpha^2\sqrt{\beta}-2\,\alpha\sqrt{\beta}+\sqrt{\beta} \right) \sqrt{\alpha}}{ \left( \alpha\beta+\alpha \right)  \left( \alpha+1 \right) }&-\frac { \left( \alpha+\beta \right)  \left( \alpha^2\sqrt{\beta}-2\,\alpha\sqrt{\beta}+\sqrt{\beta} \right) \sqrt{\alpha}}{ \left( \alpha\beta+1 \right)  \left( \alpha\beta+\alpha \right)  \left( \alpha+1 \right) }\\
 \frac {\sqrt{\beta} \left( 2\,\alpha^2\beta+2\,\alpha\beta^2-\alpha^2-2\,\alpha\beta-\beta^2+2\,\alpha+2\,\beta \right) \sqrt{\alpha}}{ \left( \beta+1 \right)  \left( \sqrt{\alpha}\alpha\beta^2+\sqrt{\alpha} \right)  \left( \alpha+1 \right) }&\frac { \left( \alpha+\beta \right) \sqrt{\beta} \left( 2\,\alpha^2\beta+2\,\alpha\beta^2-\alpha^2-2\,\alpha\beta-\beta^2+2\,\alpha+2\,\beta \right) \sqrt{\alpha}}{ \left( \alpha\beta+1 \right)  \left( \beta+1 \right)  \left( \sqrt{\alpha}\alpha\beta^2+\sqrt{\alpha} \right)  \left( \alpha+1 \right) }\\
 0&\frac { \left( \alpha^2+\alpha\beta+\alpha+\beta \right) \sqrt{\beta} \left( 2\,\alpha^2\beta+2\,\alpha\beta^2-\alpha^2-2\,\alpha\beta-\beta^2+2\,\alpha+2\,\beta \right) \sqrt{\alpha}}{ \left( \sqrt{\alpha}\alpha\beta+\sqrt{\alpha} \right)  \left( \beta+1 \right)  \left( \sqrt{\alpha}\alpha\beta^2+\sqrt{\alpha} \right)  \left( \alpha+1 \right) }\\
 0&0\\
 0&-\frac { \left( \alpha+\beta \right)  \left( \alpha^2\sqrt{\beta}-2\,\alpha\sqrt{\beta}+\sqrt{\beta} \right) }{ \left( \alpha\beta+1 \right)  \left( \alpha\beta+\alpha \right) }
\end{pmatrix}
$$

$$
\begin{pmatrix}
\frac { \left( \alpha^2\beta+\alpha\beta^2+\alpha\beta+\beta^2 \right) \sqrt{\alpha}}{ \left( \alpha^2\beta^3+\alpha\beta^2+\alpha\beta+1 \right)  \left( \alpha+1 \right) }&0&0\\
\frac { \left( \alpha\beta+\beta \right)  \left( \alpha^3+2\,\alpha^2\beta+\alpha\beta^2+\alpha^2+2\,\alpha\beta+\beta^2 \right) \sqrt{\alpha}}{ \left( \alpha\beta^2+1 \right)  \left( \alpha^2\beta^2+\alpha^2\beta+\alpha\beta+\alpha \right)  \left( \alpha+1 \right) }&0& 0\\
 \frac { \left( \alpha^2\beta+\alpha\beta+\alpha+1 \right)  \left( \alpha^2\beta+\alpha\beta^2+\alpha\beta+\beta^2 \right) \sqrt{\alpha}}{ \left( \sqrt{\alpha}\alpha\beta+\sqrt{\alpha} \right)  \left( \alpha^2\beta^3+\alpha\beta^2+\alpha\beta+1 \right)  \left( \alpha+1 \right) }&0& 0\\
\frac { \left( \alpha^2\beta+\alpha\beta+\alpha+1 \right)  \left( \alpha^2\beta+\alpha\beta^2+\alpha\beta+\beta^2 \right) }{ \left( \sqrt{\alpha}\alpha\beta+\sqrt{\alpha} \right)  \left( \alpha^2\beta^3+\alpha\beta^2+\alpha\beta+1 \right) }&0&  0\\
 \frac {\alpha^2\beta+\alpha\beta^2+\alpha\beta+\beta^2}{\alpha^2\beta^3+\alpha\beta^2+\alpha\beta+1}&0& 0\\
\frac { \left( \alpha\beta+1 \right)  \left( \alpha^2\beta+\alpha\beta+\alpha+1 \right)  \left( \alpha^2\beta+\alpha\beta^2+\alpha\beta+\beta^2 \right) }{\sqrt{\beta} \left( \sqrt{\alpha}\alpha\beta+\sqrt{\alpha} \right)  \left( \alpha^2\beta^3+\alpha\beta^2+\alpha\beta+1 \right)  \left( \alpha+1 \right) }&0& 0\\
\frac { \left( \alpha\beta+\beta \right)  \left( \alpha^3\beta+\alpha^2\beta^2+4\,\alpha^3+9\,\alpha^2\beta+5\,\alpha\beta^2+5\,\alpha^2+9\,\alpha\beta+4\,\beta^2+\alpha+\beta \right) \sqrt{\alpha}}{ \left( \alpha\beta^2+1 \right)  \left( \sqrt{\alpha}^{5}\beta^2+\sqrt{\alpha}^{5}\beta+\sqrt{\alpha}\alpha\beta^2+2\,\sqrt{\alpha}\alpha\beta+\sqrt{\alpha}\alpha+\sqrt{\alpha}\beta+\sqrt{\alpha} \right)  \left( \alpha+1 \right) } & 0& \frac {\sqrt{\alpha}}{\alpha+\beta}\\
\frac { \left( \alpha\beta^2+\alpha\beta+\beta+1 \right)  \left( \alpha^2\beta+\alpha\beta^2+\alpha\beta+\beta^2 \right) \sqrt{\alpha}}{ \left( \alpha\sqrt{\beta}\beta+\sqrt{\beta} \right)  \left( \alpha^2\beta^3+\alpha\beta^2+\alpha\beta+1 \right)  \left( \alpha+1 \right) }&0& 0\\
\frac { \left( \alpha^3\sqrt{\beta}\beta+\alpha^2\sqrt{\beta}^{5}+\alpha^3\sqrt{\beta}+3\,\alpha^2\sqrt{\beta}\beta+2\,\alpha\sqrt{\beta}^{5}+2\,\alpha^2\sqrt{\beta}+3\,\alpha\sqrt{\beta}\beta+\sqrt{\beta}^{5}+\alpha\sqrt{\beta}+\sqrt{\beta}\beta \right) \sqrt{\alpha}}{ \left( \sqrt{\alpha}^{5}\beta^3+\sqrt{\alpha}\alpha\beta^2+\sqrt{\alpha}\alpha\beta+\sqrt{\alpha} \right)  \left( \alpha+1 \right) }&-\frac {\sqrt{\alpha}}{\alpha\sqrt{\beta}+\sqrt{\beta}}& \frac { \left( \alpha\beta+1 \right) \sqrt{\alpha}}{ \left( \alpha+\beta \right)  \left( \alpha\sqrt{\beta}+\sqrt{\beta} \right) }\\
\frac { \left( \alpha+\beta \right)  \left( \alpha^3\sqrt{\beta}\beta+\alpha^2\sqrt{\beta}^{5}+\alpha^3\sqrt{\beta}+3\,\alpha^2\sqrt{\beta}\beta+2\,\alpha\sqrt{\beta}^{5}+2\,\alpha^2\sqrt{\beta}+3\,\alpha\sqrt{\beta}\beta+\sqrt{\beta}^{5}+\alpha\sqrt{\beta}+\sqrt{\beta}\beta \right) \sqrt{\alpha}}{ \left( \sqrt{\alpha}\beta+\sqrt{\alpha} \right)  \left( \sqrt{\alpha}^{5}\beta^3+\sqrt{\alpha}\alpha\beta^2+\sqrt{\alpha}\alpha\beta+\sqrt{\alpha} \right)  \left( \alpha+1 \right) }&-\frac { \left( \alpha^2+\alpha\beta+\alpha+\beta \right) \sqrt{\alpha}}{ \left( \sqrt{\alpha}\alpha+\sqrt{\alpha}\beta \right)  \left( \alpha\sqrt{\beta}+\sqrt{\beta} \right) }&0\\
\frac { \left( \alpha\beta+1 \right) \sqrt{\alpha} \left( \alpha^2\beta+\alpha\beta+\alpha+1 \right)  \left( \alpha^2\beta+\alpha\beta^2+\alpha\beta+\beta^2 \right) }{ \left( \alpha+1 \right) ^{2}\sqrt{\beta} \left( \sqrt{\alpha}\alpha\beta+\sqrt{\alpha} \right)  \left( \alpha^2\beta^3+\alpha\beta^2+\alpha\beta+1 \right) }&0& 0\\
 \frac { \left( \alpha\beta^2+\alpha\beta+\beta+1 \right)  \left( \alpha^2\beta+\alpha\beta^2+\alpha\beta+\beta^2 \right) }{ \left( \alpha\sqrt{\beta}\beta+\sqrt{\beta} \right)  \left( \alpha^2\beta^3+\alpha\beta^2+\alpha\beta+1 \right) }&0&0\\
\end{pmatrix}
$$

$$
\begin{pmatrix}
0&0&0&0\\
0&0&\frac {\sqrt{\alpha}}{\alpha+\beta}& -\frac {\sqrt{\alpha}}{\alpha+\beta}\\
0&0&0& 0\\
0&0&0&0\\
0&0&0& 0\\
0&0&0&0\\
-\frac { \left( \alpha\beta^2+1 \right) \sqrt{\alpha}}{ \left( \alpha\beta+\beta \right)  \left( \alpha+\beta \right) }&0&3\,\frac {\alpha}{ \left( \alpha+1 \right)  \left( \alpha+\beta \right) }& -3\,\frac {\alpha}{ \left( \alpha+1 \right)  \left( \alpha+\beta \right) }\\
0&-\frac {\sqrt{\alpha}}{\alpha\sqrt{\beta}+\sqrt{\beta}}&\frac {\sqrt{\alpha}}{\alpha\sqrt{\beta}+\sqrt{\beta}}& 0\\
0&-3\,\frac {\alpha}{ \left( \alpha\sqrt{\beta}+\sqrt{\beta} \right)  \left( \alpha+1 \right) }&3\,\frac {\alpha}{ \left( \alpha\sqrt{\beta}+\sqrt{\beta} \right)  \left( \alpha+1 \right) }& 0\\
0&0&3\,\frac {\sqrt{\alpha}}{\alpha\sqrt{\beta}+\sqrt{\beta}}&0\\ 
0&0&0&0\\
0&0&\frac {\alpha+1}{\alpha\sqrt{\beta}+\sqrt{\beta}}& 0
\end{pmatrix}
$$

$$
\begin{pmatrix}
-\beta&-\beta&-\beta\\ 
-\frac { \left( \alpha^3+2\,\alpha^2\beta+\alpha\beta^2+\alpha^2+2\,\alpha\beta+\beta^2 \right) \beta}{\alpha\beta^2+\alpha^2}&0&0\\ 
0&0&0\\ 
 0&0&0\\ 
-\frac { \left( \alpha+1 \right) \beta}{\sqrt{\alpha}}&0&0\\ 
 0&0&0\\
-3\,\frac { \left( \alpha^3+2\,\alpha^2\beta+\alpha\beta^2+\alpha^2+2\,\alpha\beta+\beta^2 \right) \sqrt{\alpha}\beta}{ \left( \alpha+1 \right)  \left( \alpha\beta^2+\alpha^2 \right) }&-\frac { \left( 3\,\beta^2+3\,\alpha \right) \sqrt{\alpha} \left( \alpha^3+2\,\alpha^2\beta+\alpha\beta^2+\alpha^2+2\,\alpha\beta+\beta^2 \right) \beta}{ \left( \alpha\beta+\beta^2+\alpha+\beta \right)  \left( \alpha+1 \right)  \left( \alpha\beta^2+\alpha^2 \right) }&-\frac { \left( 3\,\beta^2+3\,\alpha \right) \sqrt{\alpha} \left( \alpha^3+2\,\alpha^2\beta+\alpha\beta^2+\alpha^2+2\,\alpha\beta+\beta^2 \right) \beta}{ \left( \alpha\beta+\beta^2+\alpha+\beta \right)  \left( \alpha+1 \right)  \left( \alpha\beta^2+\alpha^2 \right) }\\ 
\frac { \left( \alpha\sqrt{\beta}\beta-3\,\alpha\sqrt{\beta}-3\,\sqrt{\beta}\beta+\sqrt{\beta} \right) \beta}{\beta^2+\alpha}&\frac { \left( \alpha\sqrt{\beta}\beta-3\,\alpha\sqrt{\beta}-3\,\sqrt{\beta}\beta+\sqrt{\beta} \right) \beta}{\beta^2+\alpha}&\frac { \left( \alpha^2+\beta \right) \sqrt{\beta}}{\alpha}\\ 
-3\,\frac {\sqrt{\beta}\beta\sqrt{\alpha} \left( \alpha^3+2\,\alpha^2\beta+\alpha\beta^2+\alpha^2+2\,\alpha\beta+\beta^2 \right) }{ \left( \beta+1 \right)  \left( \alpha+1 \right)  \left( \alpha\beta^2+\alpha^2 \right) }&-3\,\frac {\sqrt{\beta}\beta\sqrt{\alpha} \left( \alpha^3+2\,\alpha^2\beta+\alpha\beta^2+\alpha^2+2\,\alpha\beta+\beta^2 \right) }{ \left( \beta+1 \right)  \left( \alpha+1 \right)  \left( \alpha\beta^2+\alpha^2 \right) }&0\\ 
 -3\,\frac {\sqrt{\beta}\beta \left( \alpha^3+2\,\alpha^2\beta+\alpha\beta^2+\alpha^2+2\,\alpha\beta+\beta^2 \right) }{ \left( \beta+1 \right)  \left( \alpha\beta^2+\alpha^2 \right) }&0&0\\ 
 0&0&\frac { \left( \alpha+\beta \right) \sqrt{\beta}}{\sqrt{\alpha}}\\ 
\frac { \left( \alpha+1 \right)  \left( \alpha\sqrt{\beta}\beta-3\,\alpha\sqrt{\beta}-3\,\sqrt{\beta}\beta+\sqrt{\beta} \right) \beta}{\sqrt{\alpha} \left( \beta^2+\alpha \right) }&0&0
\end{pmatrix}
$$
\end{tiny}

 \clearpage

\bibliographystyle{plain}
\bibliography{mybib}

\begin{thebibliography}{10}

\bibitem{Matrics}
Plateforme de calcul {M}atri{CS} de l'université de {P}icardie {J}ules
  {V}erne.
\newblock
  \url{https://www.u-picardie.fr/recherche/presentation/plateformes/plateforme-matrics-382844.kjsp}.

\bibitem{A-L}
Dean Alvis and George Lusztig.
\newblock The representations and generic degrees of the {H}ecke algebra of
  type {$H_{4}$}.
\newblock {\em J. Reine Angew. Math.}, 336:201--212, 1982.

\bibitem{A-K}
S.~Ariki and K.~Koike.
\newblock A {H}ecke algebra of {$(\Z/r\Z)\wr S_n$} and construction of its
  irreducible representations.
\newblock {\em Adv. Math.}, 106(2):216--243, 1994.

\bibitem{ASCHMAXSUBGRPS}
M.~Aschbacher.
\newblock On the maximal subgroups of the finite classical groups.
\newblock {\em Invent. Math.}, 76(3):469--514, 1984.

\bibitem{Aschbabook}
M.~Aschbacher.
\newblock {\em Finite group theory}, volume~10 of {\em Cambridge Studies in
  Advanced Mathematics}.
\newblock Cambridge University Press, Cambridge, second edition, 2000.

\bibitem{Bouc}
Serge Bouc.
\newblock {\em Biset functors for finite groups}, volume 1990 of {\em Lecture
  Notes in Mathematics}.
\newblock Springer-Verlag, Berlin, 2010.

\bibitem{Bourb}
Nicolas Bourbaki.
\newblock {\em Commutative algebra. {C}hapters 1--7}.
\newblock Elements of Mathematics (Berlin). Springer-Verlag, Berlin, 1998.
\newblock Translated from the French, Reprint of the 1989 English translation.

\bibitem{BN}
R.~Brauer and C.~Nesbitt.
\newblock On the modular characters of groups.
\newblock {\em Ann. of Math. (2)}, 42:556--590, 1941.

\bibitem{BHRC}
J.~N. Bray, D.~F. Holt, and C.~M. Roney-Dougal.
\newblock {\em The maximal subgroups of the low-dimensional finite classical
  groups}, volume 407 of {\em London Mathematical Society Lecture Note Series}.
\newblock Cambridge University Press, Cambridge, 2013.
\newblock With a foreword by Martin Liebeck.

\bibitem{BMR}
Michel Brou\'e, Gunter Malle, and Rapha\"el Rouquier.
\newblock Complex reflection groups, braid groups, {H}ecke algebras.
\newblock {\em J. Reine Angew. Math.}, 500:127--190, 1998.

\bibitem{BM}
O.~Brunat and I.~Marin.
\newblock Image of the braid groups inside the finite {T}emperley-{L}ieb
  algebras.
\newblock {\em Math. Z.}, 277(3-4):651--664, 2014.

\bibitem{BMM}
Olivier Brunat, Kay Magaard, and Ivan Marin.
\newblock Image of the braid groups inside the finite {I}wahori-{H}ecke
  algebras.
\newblock {\em J. Reine Angew. Math.}, 733:161--182, 2017.

\bibitem{CAR}
Roger~W. Carter.
\newblock {\em Simple groups of {L}ie type}.
\newblock Wiley Classics Library. John Wiley \& Sons, Inc., New York, 1989.
\newblock Reprint of the 1972 original, A Wiley-Interscience Publication.

\bibitem{CCNPW}
J.~H. Conway, R.~T. Curtis, S.~P. Norton, R.~A. Parker, and R.~A. Wilson.
\newblock {\em Atlas of finite groups}.
\newblock Oxford University Press, Eynsham, 1985.
\newblock Maximal subgroups and ordinary characters for simple groups, With
  computational assistance from J. G. Thackray.

\bibitem{Cox}
H.~S.~M. Coxeter.
\newblock The complete enumeration of finite groups of the form
  $r_i^2=(r_ir_j)^{k_{i,j}}=1$.
\newblock {\em Journal of the London Mathematical Society}, s1-10(1):21--25,
  1935.

\bibitem{C-R}
C.~W. Curtis and I.~Reiner.
\newblock {\em Methods of representation theory. {V}ol. {I}}.
\newblock John Wiley \& Sons, Inc., New York, 1981.
\newblock Pure and Applied Mathematics, A Wiley-Interscience Publication.

\bibitem{newgraphsEsterle}
A.~Esterle.
\newblock Self-dual {W}-graphs and bilinear form, 2018.
\newblock \url{http://www.lamfa.u-picardie.fr/marin/EsterleData-en.html}.

\bibitem{FriedVolklein}
Michael~D. Fried and Helmut V\"olklein.
\newblock The inverse {G}alois problem and rational points on moduli spaces.
\newblock {\em Math. Ann.}, 290(4):771--800, 1991.

\bibitem{CHEVIE}
M.~Geck, G.~Hiss, F.~L{\"u}beck, G.~Malle, and G.~Pfeiffer.
\newblock {\sf CHEVIE} -- {A} system for computing and processing generic
  character tables for finite groups of {L}ie type, {W}eyl groups and {H}ecke
  algebras.
\newblock {\em Appl. Algebra Engrg. Comm. Comput.}, 7:175--210, 1996.

\bibitem{G-P}
M.~Geck and G.~Pfeiffer.
\newblock {\em Characters of finite {C}oxeter groups and {I}wahori-{H}ecke
  algebras}, volume~21 of {\em London Mathematical Society Monographs. New
  Series}.
\newblock The Clarendon Press, Oxford University Press, New York, 2000.

\bibitem{GorLySol}
Daniel Gorenstein, Richard Lyons, and Ronald Solomon.
\newblock {\em The classification of the finite simple groups. {N}umber 3.
  {P}art {I}. {C}hapter {A}}, volume~40 of {\em Mathematical Surveys and
  Monographs}.
\newblock American Mathematical Society, Providence, RI, 1998.
\newblock Almost simple $K$-groups.

\bibitem{Gour}
E.~Goursat.
\newblock Sur les substitutions orthogonales et les divisions régulières de
  l'espace.
\newblock {\em Annales Scientifiques de l'ENS 3ème série tome 6}, pages
  9--102, 1889.

\bibitem{GS}
R.~M. Guralnick and J.~Saxl.
\newblock Generation of finite almost simple groups by conjugates.
\newblock {\em J. Algebra}, 268(2):519--571, 2003.

\bibitem{Gyo}
Akihiko Gyoja.
\newblock On the existence of a {$W$}-graph for an irreducible representation
  of a {C}oxeter group.
\newblock {\em J. Algebra}, 86(2):422--438, 1984.

\bibitem{HAT}
A.~Hatcher.
\newblock {\em Algebraic topology}.
\newblock Cambridge University Press, Cambridge, 2002.

\bibitem{HOEF}
P.~N. Hoefsmit.
\newblock {\em Representations of {H}ecke algebras of finite groups with
  {BN}-pairs of classical type}.
\newblock ProQuest LLC, Ann Arbor, MI, 1974.
\newblock Thesis (Ph.D.)--The University of British Columbia (Canada).

\bibitem{HUP}
B.~Huppert.
\newblock {\em Endliche {G}ruppen. {I}}.
\newblock Die Grundlehren der Mathematischen Wissenschaften, Band 134.
  Springer-Verlag, Berlin-New York, 1967.

\bibitem{K}
W.~M. Kantor.
\newblock Subgroups of classical groups generated by long root elements.
\newblock {\em Trans. Amer. Math. Soc.}, 248(2):347--379, 1979.

\bibitem{KAR}
G.~Karpilovsky.
\newblock {\em The {S}chur multiplier}, volume~2 of {\em London Mathematical
  Society Monographs. New Series}.
\newblock The Clarendon Press, Oxford University Press, New York, 1987.

\bibitem{KL}
David Kazhdan and George Lusztig.
\newblock Representations of {C}oxeter groups and {H}ecke algebras.
\newblock {\em Invent. Math.}, 53(2):165--184, 1979.

\bibitem{MM}
G.~Malle and B.~H. Matzat.
\newblock {\em Inverse {G}alois theory}.
\newblock Springer Monographs in Mathematics. Springer-Verlag, Berlin, 1999.

\bibitem{M}
I.~Marin.
\newblock Irréductibilité générique des produits tensoriels de monodromies.
\newblock {\em Bull. Soc. Math. Fr.}, 132:201–232, 2004.

\bibitem{M2}
I.~Marin.
\newblock L'alg\`ebre de {L}ie des transpositions.
\newblock {\em J. Algebra}, 310(2):742--774, 2007.

\bibitem{IH2}
I.~Marin.
\newblock Infinitesimal {H}ecke algebras {II}.
\newblock arXiv:0911.1879, 2009.

\bibitem{Branch}
I.~Marin.
\newblock Branching properties for the groups {$G(de,e,r)$}.
\newblock {\em J. Algebra}, 323(4):966--982, 2010.

\bibitem{MR}
J.~Mulholland and D.~Rolfsen.
\newblock Local indicability and commutator subgroups of {A}rtin groups.
\newblock {\em arXiv:math/0606116}, 2008.

\bibitem{NAR2}
Hiroshi Naruse.
\newblock Representation theory of {W}eyl group of type {$C_n$}.
\newblock {\em Tokyo J. Math.}, 8(1):177--190, 1985.

\bibitem{NAR1}
Hiroshi Naruse.
\newblock On an isomorphism between {S}pecht module and left cell of
  {$\mathcal{S}_n$}.
\newblock {\em Tokyo J. Math.}, 12(2):247--267, 1989.

\bibitem{NivZucMont}
Ivan Niven, Herbert~S. Zuckerman, and Hugh~L. Montgomery.
\newblock {\em An introduction to the theory of numbers}.
\newblock John Wiley \& Sons, Inc., New York, fifth edition, 1991.

\bibitem{P}
H.~Pollatsek.
\newblock Irreducible groups generated by transvections over finite fields of
  characteristic two.
\newblock {\em J. Algebra}, 39(1):328--333, 1976.

\bibitem{RIB}
Kenneth~A. Ribet.
\newblock Galois action on division points of {A}belian varieties with real
  multiplications.
\newblock {\em Amer. J. Math.}, 98(3):751--804, 1976.

\bibitem{SZ}
V.~N. Sere{\v{z}}kin and A.~E. Zalesski{\u\i}.
\newblock Linear groups generated by transvections.
\newblock {\em Izv. Akad. Nauk SSSR Ser. Mat.}, 40(1):26--49, 221, 1976.

\bibitem{SV}
K.~Strambach and H.~V\"olklein.
\newblock Generalized braid groups and rigidity.
\newblock {\em J. Algebra}, 175(2):604--615, 1995.

\bibitem{V}
Helmut V\"olklein.
\newblock {\em Groups as {G}alois groups}, volume~53 of {\em Cambridge Studies
  in Advanced Mathematics}.
\newblock Cambridge University Press, Cambridge, 1996.
\newblock An introduction.

\bibitem{Wa}
A.~Wagner.
\newblock Collineation groups generated by homologies of order greater than
  {$2$}.
\newblock {\em Geom. Dedicata}, 7(4):387--398, 1978.

\bibitem{W}
R.~A. Wilson.
\newblock {\em The finite simple groups}, volume 251 of {\em Graduate Texts in
  Mathematics}.
\newblock Springer-Verlag London, Ltd., London, 2009.

\end{thebibliography}
\end{document}